\documentclass[paper=a4,
bibliography=totoc,
BCOR=15mm,             
DIV=16,
numbers=noenddot, 
footinclude=false, 
footlines=1.1, 
headinclude=true, 
]{scrbook}
\addtokomafont{chapter}{\fontsize{15}{2}\selectfont}
\addtokomafont{section}{\fontsize{14}{2}\selectfont}
\addtokomafont{subsection}{\fontsize{13}{2}\selectfont}
\addtokomafont{subsubsection}{\fontsize{12}{2}\selectfont}

\RedeclareSectionCommand[
  afterindent=false,
  beforeskip=\baselineskip,
  afterskip=0.1\baselineskip]{section}
\RedeclareSectionCommand[
  afterindent=false,
  beforeskip=\baselineskip,
  afterskip=0.1\baselineskip]{subsection}
\RedeclareSectionCommand[
  afterindent=false,
  beforeskip=\baselineskip,
  afterskip=0.1\baselineskip]{subsubsection}
\RedeclareSectionCommand[
  beforeskip=0.6\baselineskip]{paragraph}

\usepackage[british]{babel}                         
\usepackage[utf8]{inputenc}
\usepackage{amsmath,amssymb,amsthm,amsfonts}        
\usepackage{bbm}                                    
\usepackage{aligned-overset}
\usepackage{mathtools}
\usepackage{float}
\usepackage{scrhack}
\usepackage{graphicx}
\usepackage[dvipsnames]{xcolor}
\usepackage{microtype}                              
\usepackage{lmodern}                                
\usepackage[color=sixclassRdYlBu5, textsize=small]{todonotes}
\definecolor{sixclassRdYlBu5}{rgb}{0.57,0.75,0.86}
\usepackage{etoolbox}
\usepackage[headsepline,automark]{scrlayer-scrpage}                               
\usepackage{relsize}                                
\usepackage{upgreek}
\usepackage{pdfpages}
\usepackage{sidecap} 
\usepackage[
algo2e,
lined,
boxed,
commentsnumbered,
linesnumbered,
ruled
]{algorithm2e}
\usepackage{xurl}                                   
\usepackage{csquotes}                               
\usepackage{enumitem}                              

\usepackage{booktabs}                               
\usepackage{verbatim}                               
\usepackage{subcaption}                             
\usepackage[format=plain,labelformat=simple,font=small]{caption}                  
\deffootnote[1em]{1em}{1em}{\textsuperscript{\thefootnotemark}}
\DeclareTOCStyleEntry[numwidth=4em,beforeskip=-0.2cm]{tocline}{figure} 
\DeclareTOCStyleEntry[numwidth=4em,beforeskip=-0.2cm]{tocline}{table} 

\usepackage{pgfplots}

\usepackage{pgfkeys}
\newenvironment{customlegend}[1][]{%
	\begingroup
	\csname pgfplots@init@cleared@structures\endcsname
	\pgfplotsset{#1}%
}{%
	\csname pgfplots@createlegend\endcsname
	\endgroup
}%
\def\addlegendimage{\csname pgfplots@addlegendimage\endcsname}

\pgfplotsset{compat=newest}
\usepgfplotslibrary{groupplots}
\usepgfplotslibrary{dateplot}
\usepgfplotslibrary{units}
\usetikzlibrary{spy,backgrounds}
\usepackage{pgfplotstable}

\usepackage{tikz}
\usetikzlibrary{external}
\usetikzlibrary{arrows,positioning,calc,chains}
\usepackage{xintexpr}
\tikzset{
	>=stealth',
	punkt/.style={
		rectangle,
		draw=black, thick,
		text width=6.5em,
		minimum height=15em,
		minimum width=15em,
		fill = {rgb:red,108;green,152;blue,221},
	},
	blue-rec/.style={
		rectangle,
		draw=black, thick,
		text width=6.5em,
		minimum height=11em,
		minimum width=11em,
		fill = {rgb:red,108;green,152;blue,221},
	},
	Smatrixgreen/.style={
		rectangle,
		rounded corners,
		draw=black, thick,
		text width=1.5em,
		minimum height=4em,
		minimum width=4em,
		fill = green!80,
	},
	Smatrixgreen2/.style={
		rectangle,
		rounded corners,
		draw=black, thick,
		text width=1.5em,
		minimum height=2em,
		minimum width=2em,
		fill = green!80,
	},
	Smatrixblue/.style={
		rectangle,
		rounded corners,
		draw=black, thick,
		text width=1.5em,
		minimum height=4em,
		minimum width=4em,
		fill = blue!20,
	},
	Smatrixhuge/.style={
		rectangle,
		rounded corners,
		draw=black, thick,
		text width=1.5em,
		minimum height=16em,
		minimum width=16em,
		fill = black!5
	},
	Smatrixblack/.style={
		rectangle,
		rounded corners,
		draw=black, thick,
		text width=0em,
		minimum height=3.75em,
		minimum width=1.75em,
		fill = black!10
	},
	Smatrixblack2/.style={
		rectangle,
		rounded corners,
		draw=black, thick,
		text width=1.5em,
		minimum height=1.75em,
		minimum width=3.75em,
		fill = black!10
	},
	Smatrixred/.style={
		rectangle,
		rounded corners,
		draw=black, thick,
		text width=1.5em,
		minimum height=4em,
		minimum width=4em,
		fill = red!80,
	},
	pil/.style={
		->,
		shorten <=5pt,
		shorten >=5pt, anchor=north}
}

\setkomafont{sectioning}{\normalfont\normalcolor\bfseries}
\colorlet{chaptercolor}{Gray}

\renewcommand\raggedchapter{\raggedleft}
\setkomafont{chapter}{\normalfont\Huge}
\setkomafont{chapterprefix}{\normalfont\large}
\newkomafont{chapternumber}{\mdseries\fontsize{50pt}{60pt}\selectfont}

  \tikzset{
      headings/base/.style = {
        outer sep = 0pt,
        inner sep = 5pt,
      },
      headings/chapterbackground/.style = {
        headings/base,
        shade,
        left color = white,
        right color = chaptercolor,
      },
      headings/chapapp/.style = {
        headings/base,
        text = chaptercolor,
        font = \usekomafont{chapterprefix}
      },
      headings/chapternumber/.style= {
        headings/base,
        text = chaptercolor,
        font = \usekomafont{chapternumber}
      },
      headings/chapterline/.style = {
        chaptercolor,
        line width = 2pt
      }
  }

  \RedeclareSectionCommand[
  beforeskip=2\baselineskip,
  ]{chapter}

  \renewcommand\raggedchapter{\raggedleft}
  \makeatletter
  \renewcommand*\chapterlinesformat[3]{%
    \Ifstr{#1}{chapter}{%
      \begin{tikzpicture}[baseline=(title.base)]
        \node(title){%
          \parbox[t]
            {\dimexpr\textwidth-2\pgfkeysvalueof{/pgf/inner xsep}\relax}
            {\raggedchapter #3}%
        };
        \node[headings/chapapp,anchor=south east]
          at (title.north east){\Ifstr{#2}{}{}{\chapapp}\strut};
        \useasboundingbox
          (current bounding box.north west)
          rectangle
          ([yshift=-10pt]current bounding box.south east);
        \draw[headings/chapterline]
          (current bounding box.south east)++(+.5\pgflinewidth,0)--+(0,0.15\paperheight);
        \node[anchor=base west,headings/chapternumber]
          at([xshift=10pt]title.base-|current bounding box.east){#2};
      \end{tikzpicture}
      \par
    }{%
      \@hangfrom{#2}{#3}
    }%
  }
  \makeatother

\renewcaptionname{british}{\appendixname}{Chapter} 

\usepackage[hidelinks]{hyperref}
 \usepackage[capitalize,sort&compress,noabbrev]{cleveref}
\crefname{subsubsubappendix}{Appendix}{Appendix} 

\usepackage{nameref}
\usepackage[
  acronym,            
  toc,                
  xindy,
  section,
  nogroupskip,
  numberedsection=nameref,
]{glossaries}
\usepackage{glossary-longbooktabs}
\setglossarystyle{long3col-booktabs}


\makeglossaries 

\newcommand{\mytitle}{Adaptive Reduced Basis Methods for Multiscale Problems and Large-scale PDE-constrained Optimization}
	
\title{\mytitle}
\author{Tim Keil}
\date{\today}

\theoremstyle{plain}
\newtheorem{lemma}{Lemma}[section]
\newtheorem{theorem}[lemma]{Theorem}
\newtheorem{remark}[lemma]{Remark}
\newtheorem{corollary}[lemma]{Corollary}
\newtheorem{definition}[lemma]{Definition}
\newtheorem{problem}[lemma]{Problem}
\newtheorem{proposition}[lemma]{Proposition}
\newtheorem{assumption}{Assumption}

\newtheorem*{proposition_red_grad}{Repetition of \cref{prop:true_corrected_reduced_gradient_adj}}
\newtheorem*{proposition_red_Hessian}{Repetition of \cref{prop:true_corrected_reduced_Hessian}}

\DeclarePairedDelimiterX\Set[1]\{\}{%

\,#1\,
}

\DeclareMathOperator{\esssup}{ess \, sup}

\DeclareMathOperator{\supp}{supp}

\DeclareMathOperator{\essinf}{ess \, inf}
\DeclareMathOperator*{\argmin}{argmin}
\DeclareMathOperator*{\argmax}{arg\,max}
\DeclareMathOperator{\Span}{span}

\DeclarePairedDelimiterX\norm[1]{\|}{\|}{#1}
\DeclarePairedDelimiterX\abs[1]{|}{|}{#1}

\DeclarePairedDelimiterX\mengenA[1]{\lbrace}{\rbrace}{#1}
\DeclarePairedDelimiterX\mengenB[2]{\lbrace}{\rbrace}{#1\, \delimsize\vert \, #2}
\makeatletter   
\newcommand{\set}[2][\relax]{
	\ifx#1\relax \ensuremath{
		\mengenA*{#2}}
	\else \ensuremath{%
		\mengenB*{#1}{#2}}
	\fi}
\makeatother

\newcommand{\hS}[1]{\hspace{#1pt}}

\newcommand{\Params}{\mathcal{P}}
\newcommand{\Paramsad}{\mathcal{P}} 
\newcommand{\ParamsTrain}{\mathcal{P}_\textnormal{train}}
\newcommand{\ParamsVer}{\mathcal{P}_\textnormal{val}}
\newcommand{\muBar}{{\overline{\mu}}}

\newcommand{\J}{\mathcal{J}}
\newcommand{\LL}{\mathcal{L}}
\newcommand{\HH}{\mathcal{H}}
\newcommand{\Jhat}{\hat{\mathcal{J}}}
\newcommand{\uclassic}{u_c}
\newcommand{\R}{\mathbb{R}}
\newcommand{\N}{\mathbb{N}}
\newcommand{\pr}{\textnormal{pr}}
\newcommand{\du}{\textnormal{du}}
\newcommand{\D}{d}
\newcommand{\red}{r}

\newcommand{\Jnoncor}{\hat{J}}
\newcommand{\cJhatn}{{{\Jhat_\red}}}
\newcommand{\cHhatn}{\HH_{\red,\mu}}
\newcommand{\cHhath}{\HH_{h,\mu}}
\newcommand{\cHhat}{\HH_{\mu}}

\newcommand{\dred}[1]{\tilde d_{#1}}
\newcommand{\gradred}[1]{\widetilde\nabla_{#1}}
\newcommand{\noncorgrad}{\widetilde\nabla}
\newcommand{\cont}[1]{\gamma_{#1}}
\newcommand{\infsup}[1]{\gamma^{\textnormal{pg}}_{#1}}
\newcommand{\Proj}{\mathrm{P}}

\DeclareMathOperator{\integralend}{d}
\newcommand{\intend}[1]{\integralend \hS{-1.25} #1}
\newcommand{\dx}{\intend{x}}
\newcommand{\ds}{\intend{s}}

\newcommand{\bformd}{a_{\mu}}
\newcommand{\lformd}{l_{\mu}}
\newcommand{\kformd}{k_{\mu}}
\newcommand{\jformd}{j_{\mu}}
\newcommand{\resd}{r_{\mu}}

\DeclareRobustCommand{\rchi}{{\mathpalette\irchi\relax}}
\newcommand{\irchi}[2]{\raisebox{\depth}{$#1\chi$}}
\DeclarePairedDelimiterX\skal[2]{(}{)}{#1\,,\,#2}

\newcommand{\dimVH}{{N_H}}
\newcommand{\dimVh}{{N_h}}

\newcommand{\grid}{\mathcal{T}_h}

\newcommand{\Gridh}{\mathcal{T}_H}

\newcommand{\K}{\mathbb{K}}
\newcommand{\Kmu}{\mathbb{K}_\mu}
\newcommand{\Kmurb}{\mathbb{K}_\mu^{rb}}
\newcommand{\Ktmu}{\mathbb{K}_{T,\mu}}
\newcommand{\Ktmuref}{\mathbb{K}_{T,\mu_{\textnormal{ref}}}}
\newcommand{\Ktmurb}{\mathbb{K}_{T,\mu}^{rb}}
\newcommand{\Ktqo}{\mathbb{K}_{T,\xi}^{0}}
\newcommand{\Ktqrb}{\mathbb{K}_{T,\xi}^{rb}}
\newcommand{\KT}{J_T}

\newcommand{\kl}{\ell}
\newcommand{\el}{l}
\newcommand{\uhm}[1][\mu]{u_{h,#1}}
\newcommand{\phm}[1][\mu]{u_{h,#1}}
\newcommand{\uhkm}[1][\mu]{u_{H,\kl,#1}}
\newcommand{\phkm}[1][\mu]{p_{H,\kl,#1}}
\newcommand{\uhkmcoeff}{\underline{u}_{H,\kl,\mu}}

\newcommand{\uhkmsm}{u_{H,\kl,\mu}^{\textnormal{ms}}}
\newcommand{\phkmsm}{p_{H,\kl,\mu}^{\textnormal{ms}}}
\newcommand{\uhkmsmcoeff}{\underline{u}_{H,\kl,\mu}^{\textnormal{ms}}}
\newcommand{\uhkmsmrb}{u_{H,\kl,\mu}^{\textnormal{ms,rb}}}
\newcommand{\uhkmsmrblod}{u_{H,\kl,\mu}^{\textnormal{ms,rblod}}}
\newcommand{\phkmsmrb}{p_{H,\kl,\mu}^{\textnormal{ms,rb}}}

\newcommand{\uhkmrblod}{u_{H,\kl,\mu}^{\textnormal{rblod}}}
\newcommand{\phkmrblod}{p_{H,\kl,\mu}^{\textnormal{rblod}}}
\newcommand{\uhkmrb}{u_{H,\kl,\mu}^{\textnormal{rb}}}
\newcommand{\phkmrb}{p_{H,\kl,\mu}^{\textnormal{rb}}}

\newcommand{\uhkmrblodcoeff}{\underline{u}_{H,\kl,\mu}^{\textnormal{rblod}}}
\newcommand{\phkmrblodcoeff}{\underline{p}_{H,\kl,\mu}^{\textnormal{rblod}}}

\newcommand{\ehms}{e^{\textnormal{ms}}_{H,\kl}}
\newcommand{\eh}{e_{H,\kl}}

\newcommand{\vf}{v^{\textnormal{f}}}
\newcommand{\vft}[1][T]{v^\textnormal{f}_{#1}}
\newcommand{\uft}[1][T]{u^\textnormal{f}_{#1}}
\newcommand{\ufti}[1]{u^\textnormal{f}_{T_{#1}}}
\newcommand{\vfti}[1]{v^\textnormal{f}_{T_{#1}}}

\newcommand{\IH}{\mathcal{I}_H}
\newcommand{\CI}{C_{\mathcal{I}_H}}
\newcommand{\CO}{C_{\kl,\textnormal{ovl}}}
\newcommand{\CPG}{\infsup{\kl}}
\newcommand{\QQ}{{\mathcal{Q}}}
\newcommand{\QQm}{{\QQ_{\mu}}}
\newcommand{\QQkm}{{\QQ_{\kl,\mu}}}
\newcommand{\QQkmpr}{{\QQ_{\kl,\mu}^{\pr}}}
\newcommand{\QQktmpr}{{\QQ_{\kl,\mu}^{T,\pr}}}
\newcommand{\QQkmrb}{{\QQ_{\kl,\mu}^\textnormal{rb}}}
\NewDocumentCommand{\QQktm}{O{T} O{\mu}}{\QQ^{#1}_{\kl,#2}}
\NewDocumentCommand{\QQktmrb}{O{T} O{\mu}}{\QQ^{#1,rb}_{\kl,#2}}
\NewDocumentCommand{\QQktmts}{O{T} O{\mu}}{\QQ^{#1,TS}_{\kl,#2}}
\newcommand{\Vfh}{V^{\textnormal{f}}_{h}}
\newcommand{\Vfhkt}[1][T]{V^{\textnormal{f}}_{h,\kl,#1}}
\newcommand{\Vfrbkt}[1][T]{V^{\textnormal{f},\textnormal{rb}}_{\kl,#1}}
\newcommand{\Wfrbkt}[1][T]{W^{\textnormal{f},\textnormal{rb}}_{\kl,#1}}

\newcommand{\Tlast}{T_{\abs{\mathcal{T}_H}}}
\newcommand{\Vts}{\mathfrak{V}}
\newcommand{\Wtsrb}{\mathfrak{W}^{rb}}
\newcommand{\Vtsrbpr}{\mathfrak{V}^{\textnormal{rb,pr}}}
\newcommand{\Vtsrb}{\mathfrak{V}^{\textnormal{rb}}}
\newcommand{\Vtsrblod}{\mathfrak{V}^{\textnormal{rblod}}}
\newcommand{\Vtsrbdu}{\mathfrak{V}^{\textnormal{rb,du}}}
\newcommand{\pts}{\mathfrak{p}}
\newcommand{\uts}{\mathfrak{u}}
\newcommand{\utsm}{\mathfrak{u}_{\mu}}
\newcommand{\ptsm}{\mathfrak{p}_{\mu}}
\newcommand{\utsmrb}{\mathfrak{u}_{\mu}^{\textnormal{rb}}}
\newcommand{\utsmrblod}{\mathfrak{u}_{\mu}^{rblod}}
\newcommand{\ptsmrb}{\mathfrak{p}_{\mu}^{\textnormal{rb}}}
\newcommand{\ptsmrblod}{\mathfrak{p}_{\mu}^{\textnormal{rblod}}}
\newcommand{\vts}{\mathfrak{v}}
\newcommand{\Btsm}{\mathfrak{B}_{\mu}}
\newcommand{\Btsq}[1][q]{\mathfrak{B}_{\xi}}
\newcommand{\Fts}{\mathfrak{F}}
\newcommand{\sumT}{\sum_{T\in\Gridh}}
\newcommand{\Qa}{\Xi_a}

\newcommand{\dimVfrbkt}{N_T}
\newcommand{\rbt}[2][T]{\varphi_{#1,#2}}
\newcommand{\rbts}[1]{\mathfrak{b}_{#1}}
\newcommand{\rbtscoeff}[1]{\underline{\mathfrak{b}}_{#1}}
\newcommand{\rbestt}[2][T]{\psi_{#1,#2}}
\newcommand{\rbestts}[1]{\mathfrak{c}_{#1}}
\newcommand{\rbesttscoeff}[1]{\underline{\mathfrak{c}}_{#1}}

\newcommand{\diff}[2]{\frac{{\partial #1}}{{\partial #2}} }
\DeclareMathOperator{\id}{Id} 
\newcommand{\dd}{{\text d}}

\newcommand{\Jac}{J_\textnormal{ac}} 
\newcommand{\stepsizek}{s^{(k)}}

\DeclareFontFamily{U}{matha}{\hyphenchar\font45}
\DeclareFontShape{U}{matha}{m}{n}{
	<5> <6> <7> <8> <9> <10> gen * matha
	<10.95> matha10 <12> <14.4> <17.28> <20.74> <24.88> matha12
}{}
\DeclareSymbolFont{matha}{U}{matha}{m}{n}
\DeclareFontSubstitution{U}{matha}{m}{n}

\DeclareFontFamily{U}{mathx}{\hyphenchar\font45}
\DeclareFontShape{U}{mathx}{m}{n}{
	<5> <6> <7> <8> <9> <10>
	<10.95> <12> <14.4> <17.28> <20.74> <24.88>
	mathx10
}{}
\DeclareSymbolFont{mathx}{U}{mathx}{m}{n}
\DeclareFontSubstitution{U}{mathx}{m}{n}

\DeclareMathDelimiter{\vvvert}{0}{matha}{"7E}{mathx}{"17}

\DeclarePairedDelimiterXPP{\snorm}[1]{}{\lVert}{\rVert}{_{1}}{\ifblank{#1}{\:\cdot\:}{#1}}
\DeclarePairedDelimiterXPP{\anorm}[1]{}{\lVert}{\rVert}{_{a_\mu}}{\ifblank{#1}{\:\cdot\:}{#1}}
\DeclarePairedDelimiterXPP{\Snorm}[1]{}{\vvvert}{\vvvert}{_{1}}{\ifblank{#1}{\:\cdot\:}{#1}}
\DeclarePairedDelimiterXPP{\Anorm}[1]{}{\vvvert}{\vvvert}{_{a_\mu}}{\ifblank{#1}{\:\cdot\:}{#1}}
\DeclarePairedDelimiterXPP{\SMnorm}[1]{}{\vvvert}{\vvvert}{_{1,\mu}}{\ifblank{#1}{\:\cdot\:}{#1}}

\allowdisplaybreaks[1]

\newglossaryentry{params}{
  name={$\Params$},
  sort={P},
  description={Parameter space with dimension $P$}}

\newglossaryentry{R}{
	name={$\R$},
	sort={R},
	description={Real numbers}}

\newacronym{pde}{PDE}{Partial Differential Equation}
\newacronym{pdes}{PDEs}{Partial Differential Equations}

\hypersetup{
	pdftitle={\mytitle},
	pdfauthor={Tim Keil}
}

\BeforeTOCHead[toc]{%
  \KOMAoptions{parskip=false}
  \RedeclareSectionCommand[afterskip=1sp minus 1sp]{chapter}
  \RedeclareSectionCommand[afterskip=1sp minus 1sp]{section}
  \RedeclareSectionCommand[afterskip=1sp minus 1sp]{subsection}
  \RedeclareSectionCommand[afterskip=1sp minus 1sp]{subsubsection}
}
\DeclareTOCStyleEntry[beforeskip=.1cm]{chapter}{chapter}
\DeclareTOCStyleEntry[beforeskip=0.05cm]{section}{section}
\DeclareTOCStyleEntry[beforeskip=-0.05cm]{default}{subsection}
\DeclareTOCStyleEntry[beforeskip=-0.05cm]{default}{subsubsection}
\usepackage[titles]{tocloft}

\numberwithin{figure}{chapter}
\numberwithin{table}{chapter}

\begin{document}

\setlist[description]{font=\normalfont\bfseries\space,labelindent=10pt}
\frontmatter
\begin{titlepage}
\tikzexternaldisable
\begin{tikzpicture}[
	overlay,remember picture,
	every node/.append style={inner sep=0pt}
	]
	\coordinate (O) at (current page.center);
	\node[xshift=-17cm,yshift=-.7cm] (glocke) at (O) {\includegraphics[height=1.1\paperheight]{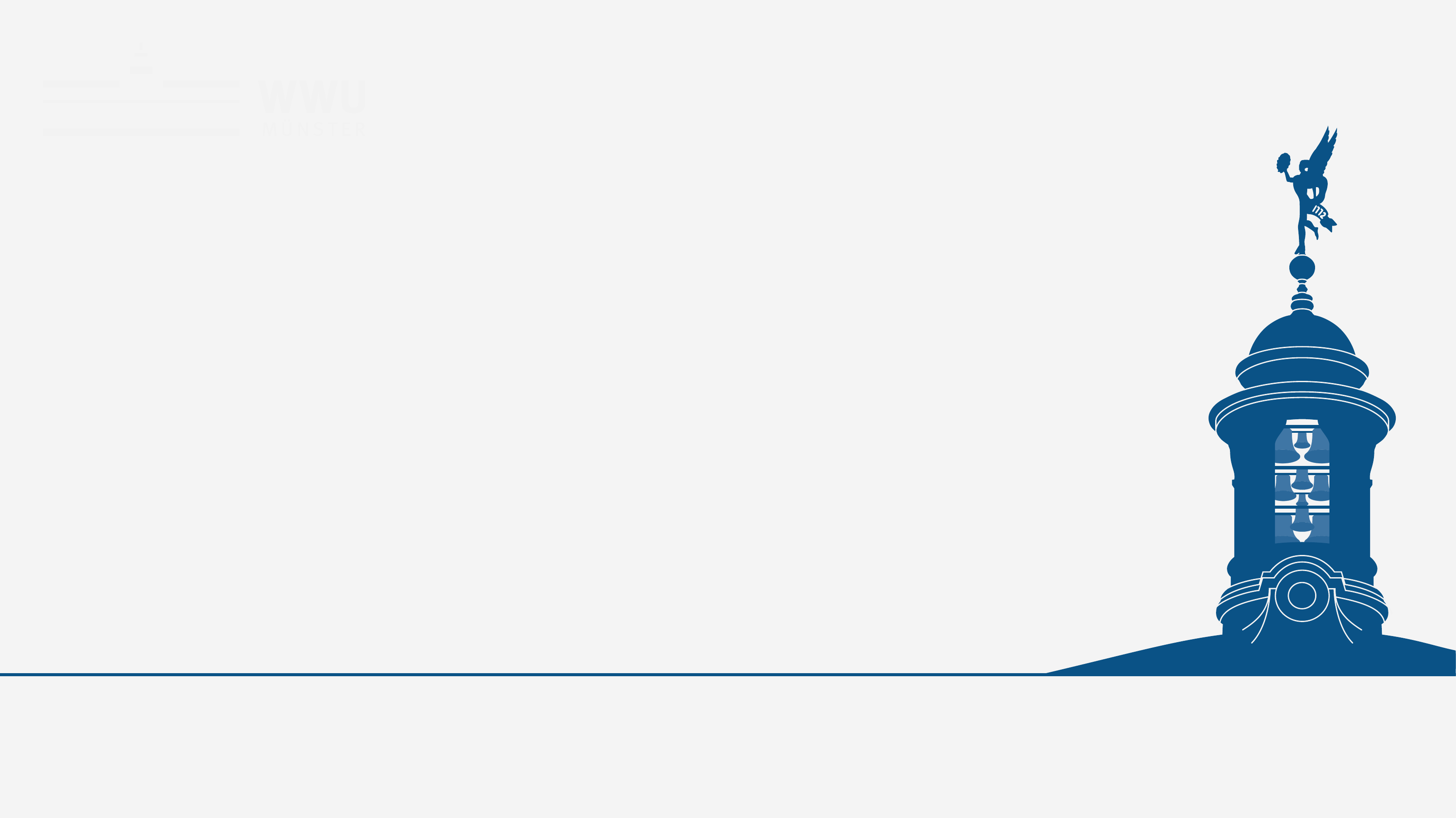}};
	\node[anchor=north west,yshift=-14mm,xshift=16mm] (wwu) at  (current page.north west){\includegraphics[width=75mm, keepaspectratio]{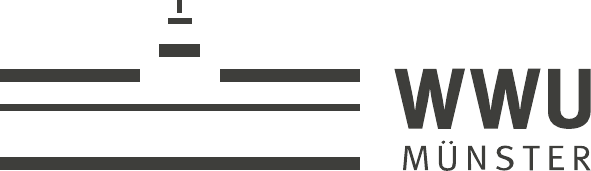}};
	\node[anchor=south west,yshift=12mm,xshift=16mm] (fb10) at  (current page.south west){\includegraphics[height=2.3cm, keepaspectratio]{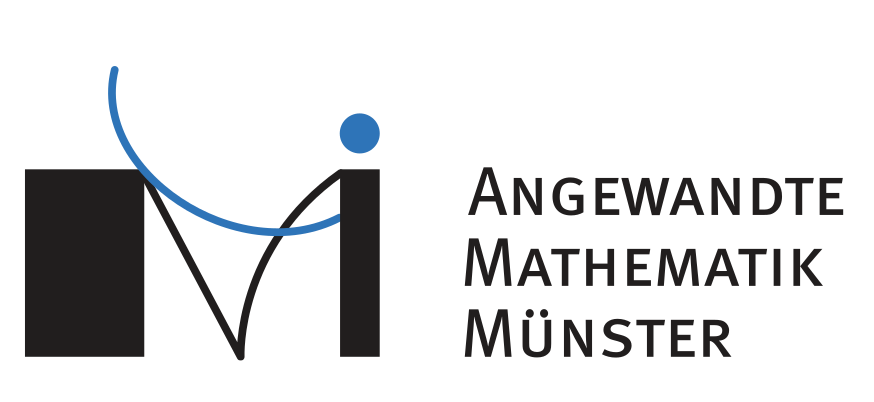}};
	\node[anchor=south east,yshift=9mm,xshift=-16mm] (mm) at  (current page.south east){\includegraphics[height=2.2cm, keepaspectratio]{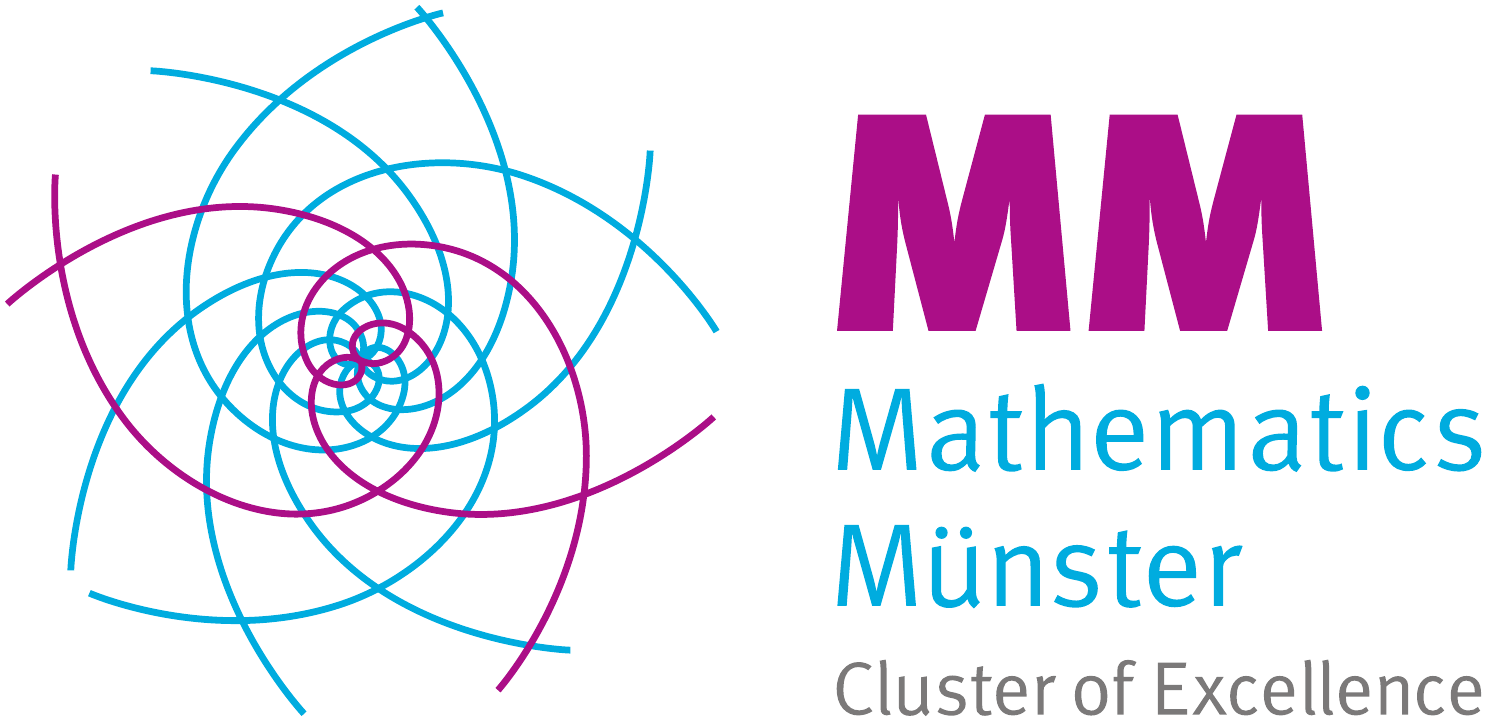}};
\end{tikzpicture}%

\def\einzugText{16mm}
\colorlet{ColorGrey}{black!60}

\begin{tikzpicture}[
	overlay,
	remember picture,
	every node/.append style={inner sep=0pt,align=left,anchor=north west}
	]
	\node[
		xshift=\einzugText,
		yshift=-65mm,
		font=\sffamily\Large
	] (text1) at (current page.north west){
		Dissertation
	};
	\node[
		xshift=\einzugText,
		yshift=-73mm,
		font=\rmfamily\Huge\bfseries,
		text width=13cm
	] (text2) at (current page.north west){
		\sffamily Adaptive Reduced Basis \\Methods for Multiscale \\Problems and Large-scale \\PDE-constrained Optimization
	};
	\node[
		xshift=138,
		yshift=-220mm,
		font=\Huge\scshape
	] (text3) at (current page.north west){
		\sffamily{Tim Keil}
	};
	\node[
	xshift=150,
	yshift=-230mm,
	font=\huge\scshape,
	ColorGrey
	] (text4) at (current page.north west){
	\sffamily{- 2022 -}
	};
	
\end{tikzpicture}

\tikzexternalenable
\end{titlepage}

\begin{titlepage}
\centering

~\vspace{10\baselineskip}

{\Large
MATHEMATIK\par}
\vspace{2\baselineskip}

{\Huge \bfseries
	\mytitle\par}
\vfill

{\large
Inauguraldissertation \\
zur Erlangung des Doktorgrades der Naturwissenschaften \\
im Fachbereich Mathematik und Informatik \\
der Mathematisch-Naturwissenschaftlichen Fakultät \\
der Westfälischen Wilhelms-Universität Münster
\par}
\vspace{4\baselineskip}

{\large
vorgelegt von \\
Tim Keil \\
aus Wickede/Ruhr
\par
}


\vspace{2\baselineskip}

{\large
- 2022 - 
}

\end{titlepage}

\thispagestyle{empty}
\par\vspace*{\fill}
{\textcolor{Gray}\hrule}
~\\[\baselineskip]
\textbf{Dekan:}  \hS{52} Prof.\ Dr.\ Xiaoyi Jiang \\
\textbf{Erstgutachter:} \hS{12} Prof.\ Dr.\ Mario Ohlberger \\
\textbf{Zweitgutachter:} \hS{5} Prof.\ Dr.\ Stefan Volkwein\\
~\\
\textbf{Tag der mündlichen Prüfung:}\hspace{1cm}
\begin{minipage}[t]{4cm}
       22.06.2022 \\
	\end{minipage} \\
\textbf{Tag der Promotion:}
\hspace{2.721cm} \begin{minipage}[t]{4cm}
      22.06.2022 \\
	\end{minipage}

\chapter*{Abstract/Zusammenfassung}
\vspace{2cm}
\section*{Abstract}
\vspace{0.1cm}
\pagenumbering{roman}
\noindent
Model order reduction is an enormously growing field that is particularly suitable for numerical simulations in real-life applications such as engineering and various natural science disciplines.
Here, partial differential equations are often parameterized towards, e.g., a physical parameter.
Furthermore, it is likely to happen that the repeated utilization of standard numerical methods like the finite element method (FEM) is considered too costly or even inaccessible.

This thesis presents recent advances in model order reduction methods with the primary aim to construct online-efficient reduced surrogate models for parameterized multiscale phenomena and accelerate large-scale PDE-constrained parameter optimization methods.
In particular, we present several different adaptive RB approaches that can be used in an error-aware trust-region framework for progressive construction of a surrogate model used during a certified outer optimization loop.
In addition, we elaborate on several different enhancements for the trust-region reduced basis (TR-RB) algorithm and generalize it for parameter constraints.
Thanks to the a posteriori error estimation of the reduced model, the resulting algorithm can be considered certified with respect to the high-fidelity model.
Moreover, we use the first-optimize-then-discretize approach in order to take maximum advantage of the underlying optimality system of the problem.

In the first part of this thesis, the theory is based on global RB techniques that use an accurate FEM discretization as the high-fidelity model.
In the second part, we focus on localized model order reduction methods and develop a novel online efficient reduced model for the localized orthogonal decomposition (LOD) multiscale method.
The reduced model is internally based on a two-scale formulation of the LOD and, in particular, is independent of the coarse and fine discretization of the LOD.

The last part of this thesis is devoted to combining both results on TR-RB methods and localized RB approaches for the LOD.
To this end, we present an algorithm that uses adaptive localized reduced basis methods in the framework of a trust-region localized reduced basis (TR-LRB) algorithm.
The basic ideas from the TR-RB are followed, but FEM evaluations of the involved systems are entirely avoided.

Throughout this thesis, numerical experiments of well-defined benchmark problems are used to analyze the proposed methods thoroughly and to show their respective strength compared to approaches from the literature.

\newpage

\hphantom{text}
\vspace{2cm}
\section*{Zusammenfassung}
\vspace{0.1cm}

\noindent Modellreduktion ist ein enorm wachsendes Gebiet, das sich besonders für numerische Simulationen in realen Anwendungen wie Ingenieurwissenschaften und verschiedenen naturwissenschaftlichen Disziplinen eignet.
Hier werden PDEs oft auf z.B. eine physikalische Größe parametrisiert.
Darüber hinaus ist es wahrscheinlich, dass das wiederholte Lösen mit numerischen Standardmethoden wie die Finite-Elemente-Methode (FEM) als zu kostspielig oder sogar unzugänglich angesehen wird.

In dieser Dissertation präsentieren wir die jüngsten Fortschritte bei Methoden zur Modellreduction mit dem primären Ziel, online-effiziente reduzierte Ersatzmodelle für parametrisierte Mehrskalenproblemen zu konstruieren und Methoden zur Optimierung von PDE-bedingten Parameter Problemen zu beschleunigen.
Insbesondere stellen wir mehrere verschiedene adaptive RB-Ansätze vor, die in einem error-aware Trust-Region-Framework für eine progressive Konstruktion eines Ersatzmodells verwendet werden können, das in einer zertifizierten äußeren Optimierungsschleife verwendet wird.
Darüber hinaus erarbeiten wir verschiedene Erweiterungen für den Trust-Region-Reduced-Basis (TR-RB)-Algorithmus und verallgemeinern ihn für Parameterbedingungen.
Dank der a-posteriori Fehlertherorie des reduzierten Modells kann der resultierende Algorithmus in Bezug auf das High-Fidelity Modell als zertifiziert angesehen werden.
Darüber hinaus verwenden wir den first-optimize-then-discretize Ansatz, um das zugrunde liegende Optimalitätssystem des Problems maximal auszunutzen.

Im ersten Teil dieser Arbeit basiert die Theorie auf globalen RB-Techniken, die eine FEM-Diskretisierung als High-Fidelity-Modell verwenden.
Im zweiten Teil konzentrieren wir uns auf lokalisierte Modellordnungsreduktionsmethoden und entwickeln ein neuartiges online effizientes reduziertes Modell für die lokalisierte orthogonale Zerlegung (LOD) Mehrkalenmethode.
Das reduzierte Modell basiert intern auf einer zweiskalen-Formulierung der LOD und ist insbesondere unabhängig von der Grob- und Feindiskretisierung der LOD.

Der letzte Teil dieser Arbeit widmet sich der Kombination beider Ergebnisse zu TR-RB-Methoden und lokalisierten RB-Ansätzen für die LOD.
Wir  stellen einen Algorithmus vor, der adaptive Localized-Reduced-Basis-Methoden im Rahmen eines Trust-Region-Localized-Reduced-Basis (TR-LRB)-Algorithmus verwendet.
Dabei wird den Grundgedanken des TR-RB gefolgt, auf FEM-Auswertungen der beteiligten Systeme jedoch gänzlich verzichtet.

Wir werden numerische Experimente mit wohldefinierten Benchmark-Problemen verwenden, um die vorgeschlagenen Methoden gründlich zu analysieren und ihre jeweilige Stärke im Vergleich zu Ansätzen aus der Literatur aufzuzeigen.

\vfill

\pagebreak
\clearpage
 \chapter*{Acknowledgements/Danksagung}
\vspace{2cm}
\section*{Acknowledgements}
\vspace{0.1cm}
\noindent
First of all, I gratefully acknowledge the unconditional support of my Ph.D. supervisor Prof. Mario Ohlberger, who provided me with an excellent working atmosphere and always supported or helped whenever necessary.
Moreover, I would like to thank Prof. Stefan Volkwein and his workgroup in Konstanz for a very fruitful collaboration.
I am particularly grateful to my co-authors directly related to this thesis: Luca Mechelli, Mario Ohlberger, Stephan Rave, Felix Schindler, and Stefan Volkwein.
I am also profoundly thankful to Fredrik Hellman and Prof. Axel M{\aa}lqvist for their support before and at the beginning of my Ph.D. and for the early research experience they shared with me.
In this context, I also express my thanks to the whole workgroup Ohlberger.

During my time at Mathematics Münster, I met excellent people who simplified the daily working day and were available for many discussions or lunch- and coffee breaks.
Particular thanks go to Chiara D'Onofrio, Julia Schleu{\ss}, Marie Tacke, Jannes Bantje, Fjedor Gaede, and Hendrik Kleikamp, who I by now can not just call colleagues but good friends.

For proofreading (at least parts) of this thesis and constructive comments, I would like to thank Lena Pillkahn, Julia Schleu\ss, Hendrik Kleikamp, and Luca Mechelli.
I want to express additional gratitude to Luca, who was immensely supportive during the whole Ph.D. project, especially in stressful times when one could always expect an answer, no matter what time or day it was.
Concerning the last 1-2 years, I also want to thank Hendrik, who supported me more than he would think he did, not just because he made the lonely Corona times in the office much more bearable.

I would also like to thank my friends who always believed in me and supported me unconditionally.
Of course, nothing would have been possible without my family.
Last but not least, my dearest thanks are dedicated to Kathrin for always having patience, supporting me with everything, accepting that I sometimes got lost in the math world, and for spending the last four years with me.

\vspace{0.2cm}
\noindent \textbf{Thank you!}

\vspace{1cm}

\noindent This work is supported by the Deutsche Forschungsgemeinschaft (DFG, German Research Foundation) under Germany's Excellence Strategy -- EXC 2044 -- 390685587, Mathematics Münster and by the DFG under contract OH 98/11-1.
\newpage

\makeatletter
\renewcommand{\@pnumwidth}{3em}
\renewcommand{\@tocrmarg}{4em}
\makeatother

\makeatletter
\renewcommand\tableofcontents{%
    \@starttoc{toc}%
}
\makeatother
\makeatletter
\renewcommand\listoftables{%
        \@starttoc{lot}%
}
\makeatother
\makeatletter
\renewcommand\listoffigures{%
        \@starttoc{lof}%
}
\makeatother

\chapter*{Contents}
\pagestyle{plain}        
\currentpdfbookmark{Contents}{Contents}
\markboth{Contents}{Contents}
\vspace{-0.85cm}
\tableofcontents
\chapter*{List of Figures}
\markboth{List of Figures}{List of Figures}
\listoffigures
\chapter*{List of Tables}
\listoftables

\mainmatter 
\chapter{Introduction}\label{chap:introduction}

\section{Motivation}

Partial differential equations (PDEs) with large- or multiscale and parameterized data functions have seen tremendous research activities in the last decades.
Efficient and stable numerical solution methods are relevant for countless real-life applications in geology, engineering, physics, chemistry, biomedicine, and other natural science disciplines.
While nowadays' mathematical models, especially those informed by real-life measurements, become more and more complex, the computational burden and memory requirements of standard numerical simulation tools, such as the finite element method (FEM), quickly reach the limits of computer technology.
Thus, the interest in more flexible, localizable, and efficient approximation tools is vast.

The computational efficiency is particularly of interest when not only a single deterministic problem is to be solved.
A famous example is the class of parameterized PDEs, where multiple (physical) parameters are varied for several different simulations, meaning that the involved data functions, boundary values, or other parameters may change.
Parameterized PDEs occur in all of the aforementioned disciplines. They are of significant interest in, for instance, optimization problems, uncertainty quantification, inverse problems, parameter identification, or other types of decision-making problems.
Since, depending on the application, the parameterized system needs to be solved for thousands of parameter samples, model order reduction (MOR) techniques are designed to decrease the computational effort for such many-query scenarios.
These methods replace the high-fidelity model with a suitable surrogate model that approximates the solution but is cheap to evaluate or store.

A real-world example is the optimal design of composite materials.
Given a suitable objective, parameters that account for different material components or different distributions of the fibers in the material need to be optimized.
Therefore, the task results in a parameter optimization problem, where multiple evaluations of the parameterized PDE are required.
Another example consists of materials whose fiber distribution can be modeled by stochastic perturbations of reference material.
Solution methods to solve these stochastic problems, e.g., the Monte Carlo method, also require thousands of realizations of the parameterized problem.
Focusing on the benchmark problem in this thesis, we also mention the PDE-constrained optimization problem of the optimal material choice for parameterized walls, doors, windows, and heaters to attain a desired stationary temperature in a specific room of a building.

Furthermore, reduced models for a parameterized PDE can, under suitable assumptions, be stored with low memory capacity and, most importantly, be evaluated quickly, which constitutes another significant advantage of model order reduction.
This makes reduced models particularly interesting for real-time decision-making, where an approximate solution needs to be found in a minimal amount of time.

A descriptive instance of such a real-time scenario is a motorsport race in Formula~1\textsuperscript{\texttrademark}.
Nowadays, designing a Formula~1 car is a challenging task carried out by engineering experts with many numerical simulations involved.
For instance, to achieve the best aerodynamics of the car, the engineer is interested in generating the most achievable downforce within the regulations and the best possible friction concerning the weight.
Advanced simulation techniques can estimate the car's lap time, given that it does not have any damage.
Let us now assume that the driver damages the car during a race.
Then, the simulation must also be capable of computing the approximated lap time, given the respective damage to the car.
The team must decide whether the car needs to be repaired (which probably costs the driver a good position in the race) or whether the car can still perform relatively well to finish the race with the damage.
Obviously, the team can not just send an extensive computation to a server and wait for the simulation to be finished but needs an immediate decision since otherwise, the race could be ruined.
This is where model order reduction methods for real-time scenarios are extremely useful.
If the reduced model has been parameterized towards the specific damage and can be evaluated efficiently within seconds, the team can immediately decide.

All the discussed scenarios can easily encounter the problem that standard approximation methods like FEM can not be used and localized methods are required instead.
Famous examples of such a case are parameterized large- or multiscale problems, for instance, groundwater flow or the aforementioned heterogeneous composite materials, porous media, or metamaterials.
Another particularly challenging task is associated with examples in which the underlying parameterization is high-dimensional.

While the real-world examples given are for motivation purposes only, this thesis is devoted to deriving mathematical methods for efficiently solving parameterized large- and multiscale elliptic problems. 
Moreover, we intend to apply global and local reduced schemes to linear-quadratic PDE-constrained optimization problems that involve elliptic PDEs.
First, we aim to efficiently obtain a numerical approximation of PDE-constrained parameter optimization problems using MOR techniques where large-dimensional parameter spaces are possible.
Second, we develop a MOR approach for a localized multiscale method suited for fast real-time evaluations of large systems.
Lastly, we combine both concepts to deduce an algorithm for solving optimization problems constrained by large- or multiscale PDEs with large parameterizations.

In the following, we provide a literature overview of the fields relevant to this thesis.
Subsequently, we clarify the goal and contribution of this work.
We emphasize that more specific references and explanations of the related mathematical background are given in \cref{chap:background}.

\section{Literature overview}

One of the main mathematical concepts of this thesis are grid-based numerical approximation schemes for PDEs. Among others, finite volume (FV), discontinuous Galerkin (DG), and, in particular, continuous finite element (FE) methods have an extensive theoretical history and have been summarized in many books.
For standard works on FV and DG methods, we refer to \cite{BHO2017,cockburn2012discontinuous,eymard2000finite}, and for the finite element method (FEM) we recommend \cite{Braess, ciarlet2002finite}.
A major ingredient of FEM is the construction of the respective finite element mesh with the aim to enable accurate finite-dimensional approximations of the true solution.
In many applications, especially those directly linked to real-life applications (see above), the underlying data functions require a particularly small resolution of the FE-mesh, which implies a high computational cost of the corresponding approximation scheme with increasing degrees of freedoms (DoFs).
Such scenarios like large- or multiscale problems have been of increasing interest for many decades.

Mesh adaptivity has been advised to minimize the number of DoFs, for instance, based on a posteriori error estimation \cite{verfurth1994posteriori}.
A wide spread method is the adaptive finite element method (AFEM) \cite{MNS2002}, see also \cite{dorfler1996convergent}.
Other methods to resolve the high computational effort of global discretizations are domain decomposition, multigrid, or multiscale methods -- They will be detailed further below.

An important problem class directly linked to real-life applications, which has already been presented in the motivation, are parameterized systems.
Usually, when applying model order reduction methods to such parameterized problems, the FE-mesh is assumed to be accurate enough such that discretization errors can be neglected.

\vspace{0.2cm}
\paragraph{Model order reduction for parameterized problems} 
As motivated above, MOR methods are particularly used for many-query or real-time scenarios, where the idea is to replace the high-fidelity model with a surrogate model.
Many MOR methods for parameterized problems are by now well-established. For a recent overview, we refer to \cite{BBH+2015}.

In this thesis, we solely concentrate on the projection-based reduced basis (RB) method that builds a reduced-order model (ROM) by projecting the underlying equations onto a problem-adapted subspace spanned by carefully chosen solutions of the full-order model (FOM).
The RB method is further based on a computational splitting into an offline- and an online phase.
In the expensive offline phase, the parameterized full-order model (for example, a standard FEM approach) is used to construct a reduced basis containing full-order solution information.
Such a basis can be found, e.g., with the help of a goal-oriented greedy-search algorithm, using an a posteriori error estimator to measure the accuracy of the surrogate model and to find quasi-optimal basis functions that minimize the reduction error \cite{binev2011convergence,Haasdonk2013}.
Alternatively, the construction of reduced bases using proper orthogonal decomposition (POD) may be used as a purely data-driven approach \cite{GV17}.
The reduced basis is used to perform a Galerkin projection of the full-order model.
Consequently, only a system of the reduced basis' size must be solved in the online phase, resulting in a speedup of multiple magnitudes.
Meanwhile, RB methods have been extended to many different fields, such as parabolic problems, inf-sup stable problems, and nonlinear problems.
For an introduction and overview of the recent development, we refer to the monographs and collections \cite{MR3672144,MR3701994,Hesthaven2016,Quarteroni2016}.
Moreover, we point to the tutorial introduction in~\cite{HAA2017} and to \cite{successandlimits} for a recent overview of open questions and challenges.
More details are further given in \cref{sec:RB_challenges}.

\vspace{0.2cm}
\paragraph{Model order reduction for PDE-constrained optimization problems} 
The class of PDE-constrained optimization problems is of great interest for many applications (see above).
Usually, the underlying PDE is solved numerically by a grid-based approach and is thus subject to the same computational issues as discussed above.
To remedy this, mesh-adaptivity for optimal control problems has been advised in \cite{BKR2000,MR2556843,Clever2012,MR2905006,MR1887737,MR2892950}.
The solution method of the optimization problem can either be followed by the \emph{first-discretize-then-optimize} or by the \emph{first-optimize-then-discretize} approach.
The main difference lies in whether the optimization problem is entirely considered in the discrete setting or whether infinite-dimensional considerations help to minimize the effect of the discretization of the system.
In this thesis, we solely concentrate on the \emph{first-optimize-then-discretize} approach that is often based on the Lagrangian functional and corresponding optimality conditions.
So-called \emph{all-at-once} approaches concentrate on solving the optimality system directly.
We refer to the class of sequential quadratic programming (SQP) methods \cite{Clever2012,ganzler2006sqp,HV06} or Gau{\ss}-Newton-type methods \cite[Chapter 18]{Nocedal}.
A particularly suitable approach for uniquely solvable PDE systems is based on the reduced formulation of a PDE-constrained optimization problem and is particularly targeted in this thesis.
For an overview of PDE-constrained optimization and general optimization methods, we further refer to \cite{HPUU2009,Nocedal}, for instance.

Since the control variables can often be interpreted as parameters in the PDE constraints, MOR methods can help to accelerate the solution process.
There exists a large amount of literature using reduced-order surrogate models for optimization methods.
A posteriori error estimates for the reduced approximation of linear-quadratic optimization problems and parameterized optimal control problems with control constraints were studied, e.g., in \cite{Dede2012,GK2011,KTV13,NRMQ2013,OP2007}.
In \cite{dihl15}, an RB approach is proposed, which also enables an estimation of the actual error on the control variable and not only on the gradient of the output functional.
Certified reduced basis methods for parametrized elliptic optimal control problems with distributed controls were studied in \cite{KTGV2018}.
With the help of an a posteriori error estimator, ROMs can be constructed with respect to the desired accuracy but also with respect to a local area in the parameter set \cite{EPR2010,HDO2011}.
For very high-dimensional parameter sets, simultaneous parameter and state reduction have been advised \cite{HO2014,HO2015,LWG2010}.

These references have in common that the offline phase for constructing a reduced-order surrogate can be ruled out of the computational efficiency of the method, implicitly assuming that the optimization problems are not aimed to be solved once, but rather with changing parameters, such as different output functionals.
Using MOR for solving a single parameter optimization problem has raised the challenge that the offline time can not be considered negligible, which can easily be significant for a prohibitively expensive forward problem.
To remedy this bottleneck, it is beneficial to use optimization methods that optimize on a local level of the control variable, assuming the surrogate only to be accurate enough in a respective parameter region.
Hence, we require an approach that goes beyond the classical offline-online decomposition.
RB methods have recently been advised with a progressive construction of ROMs in \cite{BMV2020,GHH2016,ZF2015}. 
Also, localized RB methods that are based on efficient localized a posteriori error control and online enrichment \cite{BEOR2017,SO15} overcome traditional offline-online splitting and are thus particularly well suited for applications in optimization or inverse problems \cite{OSS2018,OS2017}.

\vspace{0.2cm}
\paragraph{Trust-region reduced-order methods} 

Trust-region (TR) approaches are a class of optimization methods that are particularly tailored towards the aim to use locally accurate surrogate models for the expensive objective function without the necessity to construct a globally accurate surrogate (see \cite{TRbook,Nocedal}).
The (nonlinear) objective is replaced by a model function that can be evaluated with much less effort and whose are used to solve an optimization sub-problem in a local area of the parameter set.
Subsequently, the respective model function is updated or replaced by means of the current outer iterate.
With this strategy, TR methods ensure global convergence for locally convergent methods.
Constraints on the control and the metric for the trust-region radius can affect the convergence of the method.
The choice of the sub-problem solver can also significantly influence the convergence speed.
The TR strategy can be combined with second-order methods for nonlinear optimization: with the Newton method to solve the reduced optimization problem, and with the Gau{\ss}-Newton or sequential quadratic programming (SQP) method, cf.~\cite[Chapter 18]{Nocedal} or \cite{Clever2012,ganzler2006sqp,HV06} for the all-at-once approach, where the entire optimality system is solved at once.

If employed with reduced-order models for parameterized systems, the choice of the model function consists of the construction of a locally accurate reduced model.
During the TR optimization process, one usually moves away from the original parameters for which the initial reduced-order model was built, and the quality of the surrogate can not be guaranteed any more.
Therefore, a priori and a posteriori error analysis are required to ensure accurate reduced-order approximations for the optimization problem; cf.~\cite{GV17,HV08,KTV13}.
One suitable choice for the local model is a reduced-order discretization of the objective (e.g., by utilizing a second-order Taylor approximation).
To ensure convergence to stationary points, the accuracy of the model function and its gradient has to be monitored.
In \cite{RoggTV17}, a posteriori error bounds are utilized to control the approximation quality of the gradient.
We also refer to \cite{GGMRV17}, where the authors apply basis update strategies to improve the reduced-order approximation scheme concerning the optimization goal.
Another TR approach was proposed in \cite{AFS00,SS13} to control the quality of the (POD) reduced-order model, referred to as TR-POD, a meanwhile well-established method in applications; cf.~\cite{BC08,CAN12}.
TR algorithms with successive construction of surrogate models have also been conducted in the context of topology optimization \cite{webster2021enriched, yano2021globally}.

In an error-aware TR method, the TR radius is directly characterized by the a posteriori error estimator for the cost functional of the surrogate model.
In this way, the offline phase of the RB method can entirely be omitted since the RB model can be adaptively enriched during the outer optimization loop.
By this procedure, the surrogate model eventually will have a high accuracy around the optimum of the optimization problem, ignoring the accuracy of the part of the parameter set where the outer (and inner) optimization loop is not active.
Error-aware TR-RB methods can be utilized in many different ways. 
In \cite{YM2013}, the TR framework is combined with an efficient RB error bound for defining the trust-region in the design optimization of vibrating structures using frequency-domain formulations.
One possible TR-RB approach has been extensively studied in \cite{QGVW2017} for linear parametric elliptic equations, which ensures convergence of the non-local TR-RB.
Note that the experiments in \cite{QGVW2017} are for up to six-dimensional parameter sets without inequality constraints.

\vspace{0.2cm}
\paragraph{Numerical methods for large- and multiscale methods}
The above approaches assume that a global discretization scheme and solution method such as FEM is (at least) accessible. However, this assumption is very restrictive, especially for real-life examples.
Therefore, further research is devoted to localized methods that rely on a spatial splitting of the computational domain and instead consider coarse systems that are only locally informed by the high-fidelity resolution.

In this context, we refer to the class of domain decomposition methods \cite{DD,TW2005} that are based on such a splitting and instead solve small local problems on (overlapping or non-overlapping) sub-domains.
The resulting global model then only inherits a few degrees of freedom from the local problems.

Numerical multiscale methods are based on a similar idea of a spatial domain decomposition but are more motivated by classical homogenization theory for PDEs \cite{allaire,nguetseng1989general}.
In the last two decades, there has been a tremendous development of such methods with the primary intention to resolve the finest required scale only locally and collect the gathered fine-scale information in an effective coarse-scale global system.

Some of these multiscale methods, such as the heterogeneous multiscale method (HMM)~\cite{EE2005,Ohlberger2005,EEnq}, are based on ideas from mathematical homogenization \cite{allaire} and aim at computing effective coefficients for an appropriate coarse-scale equation.
In contrast, approaches, such as the multiscale finite element method (MsFEM)~\cite{Efendiev2009,HOS2014,MsFEM}, its generalized variants (GMsFEM) \cite{CEL2014,Efendiev2013}, or the generalized finite element method (GFEM)~\cite{LiptonBabuska}, construct coarse-scale elements that incorporate the local fine-scale features of the solution and then approximate the solution in the space spanned by these multiscale elements.
Many multiscale methods, especially the ones directly derived from homogenization theory, are designed to solve specific multiscale problems that exhibit strong assumptions on the structure of the problem (for instance, local periodicity).

An approach that does not require structural assumptions on the fine-scale data is the localized orthogonal decomposition (LOD)~\cite{HMP2014,HP2013,MP14}, which itself is based on the variational multiscale method (VMM)~\cite{Hughes1995387,HUGHES19983,LM2005}.
This method instead uses splitting the full fine-scale approximation space into a negligible fine-scale component and an orthogonal multiscale space in which the solution is sought.
This multiscale space is then approximated by computing localized auxiliary approximations of the orthogonal projection onto the fine-scale space (corrector problems).
For an introduction of the LOD, applications, and extensions, we point to \cite{LODbook} and the references therein.
A crucial detail of multiscale methods is the possibility to treat the sub-problems entirely in parallel, allowing for efficient parallel implementations on multiple CPUs or HPC clusters.
We particularly mention the Petrov--Galerkin formulation of the LOD (PG--LOD), introduced in~\cite{elf}.
This variant has lower storage and communication requirements for the computed fine-scale data in comparison to the original Galerkin formulation.
Since this variant is particularly well suited for the case of substantial storage consumptions where local fine-scale correctors can not be stored, it will extensively be used throughout this thesis.

We also note the concept of gamblets \cite{Gamblets}, which has independently been developed but is based on similar ideas as the LOD.
Here, a multi-level strategy reconstructs the fine-scale information based on a hierarchy.
The approach is also applicable to more general energy-minimizing problems.

\vspace{0.2cm}
\paragraph{Model order reduction for parameterized large- and multiscale problems}
For large parameterized multiscale problems, RB methods need to be combined with numerical multiscale methods. Otherwise, the FOM solution for a single snapshot parameter might already be computationally infeasible. 
One possible approach is to speed up the solution of the individual cell problems using RB techniques, as is done in~\cite{Bo08} in the context of numerical homogenization, in~\cite{AB12,AB13,AB14,AB14_2,A19} for the HMM, in~\cite{Efendiev2013,HesthavenZhangEtAl2015,Nguyen2008} for the MsFEM, or in~\cite{RBLOD} for the LOD.
This approach is applicable both for problems where the fine-scale data variation over the domain is parameterized, resulting in a single ROM for all cell problems~\cite{AB12,AB13,AB14,AB14_2,A19,Bo08,HesthavenZhangEtAl2015}, or for general parameterized problems, where for each cell problem, a dedicated ROM is built~\cite{RBLOD,Nguyen2008}.

Some recent works have considered the case where the multiscale data varies only in some cells, e.g.~caused by perturbations, by a sequence over a period of time, or by other parameterizations, such that individual fine-scale solutions can be reused~\cite{HKM20,HM19,MV20} or a ROM can be constructed in the case of parameterized perturbations~\cite{MV21}.

These approaches have in common that, while the constructed surrogates are independent of the resolution of the fine-scale mesh, the effort for their evaluation still scales with the size of the used coarse mesh.
Hence, if the coarse mesh itself is large, the repeated assembly and solution of the coarse system can still become a computational bottleneck.
Further, the influence of the approximation error of the cell-problem ROMs on the coarse system is not rigorously controlled.
Let us also mention that the work on the localized reduced basis multiscale method (LRBMS) also incorporates RB techniques in a domain decomposition setting by building separated RB spaces on every coarse element \cite{AHKO2012,SO15}.
For an overview of localized model order reduction and application to parameterized multiscale problems, we refer to the review article \cite{BIORSS21}.

An online-efficient RB-ROM for locally periodic homogenization problems was developed for the HMM in~\cite{OS12,schaefer13,OSS2018} based on the two-scale formulation in~\cite{Ohlberger2005}.
An online-efficient RB in the context of periodic homogenization was derived in~\cite{OS12}.
Apart from these works, we are unaware of other approaches that yield fully online-efficient multiscale ROMs.

\section{Goal of this thesis}

This thesis studies and presents several advances in adaptive trust-region reduced basis methods for PDE-constrained parameter optimization problems.
On top of that, we focus on an efficient reduced-order model for the localized orthogonal decomposition method for multiscale problems and combine the two approaches to a trust-region localized RB (TR-LRB) method.

At first, we elaborate on recent improvements of the adaptive certified error-aware TR-RB algorithm that has been introduced in \cite{QGVW2017} and is based on \cite{YM2013}.
The method aims to accelerate the solution process of a single PDE-constrained parameter optimization problem with the first-optimize-then-discretize approach.
As the first contribution of this thesis, we discuss several different reduced basis strategies for the model function in the TR-RB algorithm.
In particular, a non-conforming dual (NCD) approach is introduced. The respective reduced objective functional is based on the Lagrangian approach associated with the optimization problem and permits more accurate approximations of the objective functional.
Furthermore, we discuss alternative choices of the reduced model for the optimality system, some of which proved to be infeasible for a computationally efficient model.
At last, a Petrov--Galerkin choice of the reduced model to allow for a more straightforward computation of the required derivatives is presented.

For controlling the error of the reduced models, we provide efficiently computable a posteriori error estimates for all reduced quantities, such as the objective functional, its (approximate) gradient and hessian, and the optimal parameter.
Different from \cite{QGVW2017}, we discuss higher-order TR-RB methods using the projected Newton method to solve the TR sub-problems.
Furthermore, we rigorously prove the convergence of the TR-RB method with bilateral inequality constraints on the parameters with two different results for either RB-based approximations as well as for the (more general) finite-dimensional case, including a new optional enrichment.
With regard to the different reduced models and several new features of the TR-RB algorithm, we devise new adaptive enrichment strategies for the progressive construction of the RB spaces.
We demonstrate the main concepts of our new TR-RB methods in numerical experiments.
We show that our new techniques outperform existing model reduction approaches for large-scale optimization problems in well-defined benchmark problems under given circumstances.
On top of that, we further discuss additional improvements to the TR-RB method, including a relaxed TR-RB method (R-TR-RB) where the strong certification of the algorithm is left out on purpose with the aim of a faster (but still certified) convergence.

The second part of this thesis is devoted to the localized orthogonal decomposition method.
In particular, we investigate the PG--LOD concerning its computational requirements for many-query scenarios of parameterized PDEs.
In this context, we discuss a non-RB-based adaptive method for reusing fine-scale corrector information from former solutions.
Subsequently, we introduce a new two-scale reduced basis method for the LOD (TSRBLOD) that takes both the local corrector problems and the coarse-scale problem into account to produce a single small-size ROM that no longer requires explicit solutions of local sub-problems in the online phase.
The model is based on a new two-scale formulation of the LOD in a single variational problem, inspired by the theory of two-scale convergence in mathematical homogenization~\cite{allaire}.
As the computation of solution snapshots for this model would still be computationally expensive, we combine this formulation with a preceding reduction of the fine-scale corrector problems similar to the RBLOD approach in~\cite{RBLOD}.
Rigorous and efficient a posteriori bounds derived from the two-scale formulation control the error over this entire two-stage reduction process.

Another main goal of this thesis is the combination of TR-RB methods for PDE-constrained parameter optimization problems and localized model order reduction techniques for parameterized multiscale problems.
In this context, we introduce trust-region localized reduced basis methods (TR-LRB) that progressively construct a localized reduced model following a trust-region approach based on the TR-RB method for global approximation schemes.
To this end, we combine the RBLOD and TSRBLOD to a specific instance of a TR-LRB method.
For this, we devise a suitable a posteriori error estimation and demonstrate the method in numerical examples.

\section{Overview of the author's contribution}

The content of this thesis has partly been published in \cite{paper2,HKM20,paper1,paper3,paper4,tsrblod}, where \cite{paper2,paper1,paper3} are contained in \cref{chap:TR_RB}, the main results of \cite{HKM20} are briefly presented in \cref{sec:adaptive_PGLOD_from_hekema}, the primary content of \cref{chap:TSRBLOD} is based on \cite{tsrblod}, and \cref{chap:TR_TSRBLOD} is an extended version of \cite{paper4}.
In what follows, we chronologically itemize this thesis' author's research items, including a short description and how it is related to the work at hand:

{\setlist[description]{font=\normalfont\bfseries\space,labelindent=0pt}
\begin{description}
	\item[\cite{HKM20}]  F. Hellman, T. Keil, and A. M{\aa}lqvist. Numerical upscaling of perturbed diffusion problems. SIAM
	Journal on Scientific Computing, 42(4):A2014–A2036, 2020.
	
	\begin{addmargin}[1em]{0em}
		While the main idea of the adaptive PG--LOD method for successively solving perturbed problems was already analyzed in the author's master's thesis, relevant additional results and details have been introduced in the referenced paper.
		Since they also fit into the setting of many-query scenarios for the LOD, we briefly refactor the results in \cref{sec:adaptive_PGLOD_from_hekema}.
	\end{addmargin}
	\item[\cite{paper1}] T. Keil, L. Mechelli, M. Ohlberger, F. Schindler, and S. Volkwein. A non-conforming dual approach
	for adaptive trust-region reduced basis approximation of PDE-constrained parameter optimization.
	ESAIM. Mathematical Modelling and Numerical Analysis, 55(3):1239, 2021.
	
	\begin{addmargin}[1em]{0em}
	In the first work concerning global TR-RB methods for solving PDE-constrained parameter optimization problems, we introduced the NCD-corrected functional and analyzed BFGS-based methods with special emphasis on the features that have been added to the TR-RB approach proposed in \cite{QGVW2017}.
	\end{addmargin}
	\item[\cite{paper2}] S. Banholzer, T. Keil, L. Mechelli, M. Ohlberger, F. Schindler, and S. Volkwein. An adaptive
	projected newton non-conforming dual approach for trust-region reduced basis approximation of
	PDE-constrained parameter optimization. Pure and Applied Functional Analysis, 7(5):1561–1596,
	2022.
	\begin{addmargin}[1em]{0em}
	In this follow-up paper of \cite{paper1}, we further elaborated on Newton-type TR-RB methods with optional RB basis enrichment and introduced a way to validate the approximation quality of the optimal parameter.
	\end{addmargin}
	\item[\cite{paper3}] T. Keil and M. Ohlberger. ‘Model Reduction for Large-Scale Systems’. In: Large-Scale
	Scientific Computing. Ed. by I. Lirkov and S. Margenov. Cham: Springer International
	Publishing, 2022, pp. 16–28.
	
	\begin{addmargin}[1em]{0em}
	In these proceedings, we deviated from the NCD-corrected approach and discussed a Petrov--Galerkin variant of the TR-RB algorithm.
	\end{addmargin}
	\item[\cite{tsrblod}]  T. Keil and S. Rave. An online efficient two-scale reduced basis approach for the localized orthogonal
	decomposition. arXiv preprint arXiv:2111.08643, 2021.

	\begin{addmargin}[1em]{0em}
	This paper was mainly motivated by weaknesses of the RBLOD proposed in \cite{RBLOD} concerning the online efficiency of the reduced scheme.
	The RBLOD constitutes the first work regarding localized MOR for the LOD on parameterized multiscale problems.
	In our work, we devised a two-scale-based TSRBLOD surrogate model for the PG--LOD that is well-suited for many-query simulations and remarkably efficient in real-time applications of large-scale problems.
	\end{addmargin}

	\item[\cite{AMLenopt}]  T. Keil, H. Kleikamp, R. J. Lorentzen, M. B. Oguntola, and M. Ohlberger. Adaptive machine
	learning based surrogate modeling to accelerate PDE-constrained optimization in enhanced oil
	recovery. Advances in Computational Mathematics, 48(73), 2022.
	
	\begin{addmargin}[1em]{0em}
	This collaboration was initiated by our work on \cite{paper1}.
	In the referenced paper, we aimed to accelerate the solution process of PDE-constrained optimization problems in enhanced oil recovery.
	Since the high-fidelity solver of the three-phase-flow polymer model was considered a black-box model, we did not apply RB-based methods to reduce the system.
	Instead, we used non-intrusive machine-learning models.
	In particular, we devised deep neural networks (DNNs) to learn the input-output map of the control space to the objective functional.
	Just as in the TR-RB method, the surrogate model is updated progressively without the aim to be accurate for the whole parameter space.
	Although using DNNs prevents using certified error estimates from assessing the accuracy of the surrogate, we used a FOM-based adaptive algorithm that still converges with respect to a FOM criterion.
	This work is not included in this thesis and is only shorty mentioned in \cref{sec:non-certified_TR}.
	\end{addmargin}
		
	\item[\cite{paper4}]  T. Keil and M. Ohlberger. A relaxed localized trust-region reduced basis approach for optimization
	of multiscale problems. arXiv preprint arXiv:2203.09964, 2022.
	
	\begin{addmargin}[1em]{0em}
	A natural contribution that is directly related to this thesis is the combination of the TR-RB method, proposed in \cite{paper2,paper1,paper3}, and the two-scale RB approach for the PG--LOD, introduced in \cite{tsrblod}.
	The work in \cite{paper4} consists of a first instance towards localized TR-RB methods that do not incorporate FEM solves and global RB methods but online adaptive localized RB enrichment.
	\end{addmargin}
\end{description}
}
\noindent Furthermore, the author has contributed to the following software packages to carry out the numerical experiments for the listed publications.

\begin{description}
	\item[\texttt{pyMOR} \cite{pymor}:] The model order reduction library \texttt{pyMOR} is a pure \texttt{Python} software that has been developed for many years and is by now a convenient and widely used software for implementing MOR-related methods.
	The entire RB-related numerical experiments of this thesis have been implemented with \texttt{pyMOR}, for which the author of this thesis has contributed to \texttt{pyMOR}'s development.
	For instance, \texttt{pyMOR}'s usability for parameter derivatives has been extended, and an elaborate tutorial on PDE-constrained parameter optimization methods has been added in order to introduce users to the topic and to simplify first points of contact with \texttt{pyMOR}.
	More implementational details on how \texttt{pyMOR} is used and extended are given in \cref{sec:software_TRRB}, \cref{sec:software_TSRBLOD}, and \cref{sec:software_TR_TSRBLOD}.
	\item[\texttt{gridlod} \cite{gridlod}:] The pure \texttt{Python} library \texttt{gridlod} was mainly developed by Fredrik Hellman and mainly consists of a implementing the PG--LOD with specific focus on a parallelizable code design.
	Thus, \texttt{gridlod} is the groundwork for the LOD-related numerical experiments of this thesis.
	For the numerical experiments in \cite{HKM20} and \cite{paper4,tsrblod}, the author added relevant contributions to the development of the software package and implemented \texttt{pyMOR}-bindings for using \texttt{gridlod} together with MOR methods.
	Again, more details on the implementation are given in \cref{sec:software_TSRBLOD} and \cref{sec:software_TR_TSRBLOD}.
\end{description}

\newpage
\section{Outline}

This thesis is organized as follows:
In \cref{chap:background}, we give a brief introduction to the basic concepts of this thesis.
Among others, we introduce (parameterized) elliptic PDEs, the theory of PDE-constrained optimization methods, iterative optimization methods, and the basic principles of the numerical approximation of (parameterized) PDEs.
Moreover, we introduce model order reduction techniques and describe their relevance to this thesis.

In \cref{chap:TR_RB}, we carefully describe the TR-RB approach for the spatially global FEM method.
In particular, we discuss several strategies for the reduced model and device the corresponding a posteriori error analysis.
Below, we specify the TR-RB in greater depth and provide the corresponding convergence result.
During the numerical experiments, we are going to set up a well-defined benchmark problem and to demonstrate the TR-RB method on several examples with different emphases.
Lastly, we end the chapter with an overview of further variants of the TR-RB method.

\cref{chap:TSRBLOD} is devoted to localized model order reduction with special emphasis on the LOD method for parameterized systems.
We start with a general view on localized model order reduction and subsequently introduce the Petrov--Galerkin variant of the LOD method while focusing on its computational complexity for a many-query simulation.
In \cref{sec:two_scale} and the sections thereafter, we introduce the TSRBLOD, including corresponding numerical experiments.

In \cref{chap:TR_TSRBLOD}, we merge the ideas from \cref{chap:TR_RB} and \cref{chap:TSRBLOD} to TR-LRB methods.
While we first formulate these methods in a general way, we describe the particular instance of a TR-LRB algorithm based on the RB-based reduced formulations of the LOD and show numerical experiments.

Concluding remarks and an outlook to further research perspectives are given in \cref{chap:outlook}.

\chapter{Mathematical background}\label{chap:background}
\noindent
In this chapter, we introduce the mathematical background of this thesis with the primary intention to provide a more mathematical view of the concepts discussed in \cref{chap:introduction}.
In \cref{sec:preliminaries}, we provide a brief introduction to the weak formulation of partial differential equations and discuss parameterized PDEs.
Afterwards, these problem classes are used for PDE-constrained optimization problems, presented in \cref{sec:background_PDE-constr}, followed by a brief overview of general iterative optimization methods in \cref{sec:iterative_optimization}.
For solving PDE-constrained optimization problems numerically, in \cref{sec:background_npdgl}, we review the basic concepts of the numerical approximation of PDEs.
We finalize the chapter by introducing model order reduction techniques for the efficient approximation of parameterized PDEs in \cref{sec:background_model_order_reduction}.
To neglect ambiguities, this chapter is also devoted to specifying the necessary notation used throughout this work.
We emphasize that the sections are not aimed at theoretical completeness in the respective field but rather serve as a brief introduction to the relevant theory.
We provide literature for further reading of the discussed topic in every section.

\section{Preliminaries}
\label{sec:preliminaries}

We start with an introduction to the basic mathematical theory of elliptic PDEs, their weak formulation, and existence results.
Moreover, we formulate parameterized elliptic PDEs and introduce suitable derivatives of the involved components.

\subsection{A prototypical model problem in weak formulation}
The \emph{stationary heat equation} or \emph{stationary diffusion equation} can be considered as a prototypical model problem for the entire thesis.
Let $\Omega \subset \mathbb{R}^d$, $d=1,2,3$, be a bounded polygonal or polyhedral domain with boundary $\Gamma$.
The classical formulation of stationary heat is a second-order elliptic PDE: Find $\uclassic$ such that
\begin{equation}
\begin{aligned}
\label{eq:classical_PDE}
- \nabla \cdot A  \nabla \uclassic &= f, \qquad \text{ in }\Omega, \\
\uclassic & = 0, \qquad \text{ on }\Gamma, 
\end{aligned}
\end{equation}
where we assume to be given a source term $f$ and a positive definite diffusion coefficient~$A$.
To keep the theory for this section as simple as possible, we consider homogeneous Dirichlet boundary conditions and refer, for instance, to \cref{sec:TRRB_num_experiments} for other boundary conditions.
For such an arbitrary complex domain $\Omega$ and data functions $A$ and $f$, solving~\eqref{eq:classical_PDE} can be a challenging task.
Most importantly, in lots of engineering models, it is not even given that a classical solution exists in $C^2(\Omega)$, for instance, due to discontinuities in the diffusion coefficient~$A$.
To circumvent this, linear functional analysis (e.g.~\cite{alt}) and, in particular, the theory of Sobolev and Hilbert spaces, enable a weak formulation of~\eqref{eq:classical_PDE}, allowing for less regularity of the solution and (under sufficient assumptions) admits a unique (weak) solution.
For the weak (or variational) formulation, let $V$ be a real-valued Hilbert space, equipped with an inner product $(\cdot,\cdot)$ and a corresponding norm $\norm{\cdot}$.
We define an \emph{elliptic problem} in the following way:

\begin{problem}[Elliptic problem]
	\label{def:elliptic_problem}
	Let $a: V \times V \to \R$ denote a continuous and coercive bilinear form and $l \in V'$ a continuous linear functional (definitions of the terms are given below).
	As an elliptic problem, we consider the task to find $u \in V$ such that
	\begin{equation}
	\label{eq:weak_problem}
	a(u,v) =  l(v) \qquad\qquad \text{for all } \, v \in V.
	\end{equation} 
\end{problem}

\noindent Elliptic problems occur in many applications. As far as this thesis is concerned, we omit a further discussion on non-elliptic problems.
However, we indicate that the spectrum of problems that can be written as \eqref{eq:weak_problem} goes far beyond the examples given in this thesis.

Under suitable assumptions, we can show uniqueness and existence of a solution $u \in V$ of \cref{def:elliptic_problem}, which can be proven by the \emph{Lax-Milgram Theorem} (see, e.g.,~\cite[Theorem 4.3.16]{hilbertspaces}).
For the simple example of the stationary diffusion~\eqref{eq:classical_PDE}, we have $V = H^1_0(\Omega)=\{v\in H^1(\Omega)\ |\ \text{tr}(v)=0\}$ as the Hilbert space of weakly differentiable functions with vanishing boundary values.
Standard norms that are based on inner products on $V$ are the standard $L^2$-, $H^1$-semi, and $H^1$-norm, i.e.~
\begin{equation} \label{eq:norms}
\begin{split}
\norm{v}^2_{L^2(\Omega)}  &\coloneqq \int_\Omega \abs{v(x)}^2\dx, \\
|v|^2_{H^1(\Omega)}  &\coloneqq \int_\Omega \abs{\nabla v(x)}^2\dx, \\
\norm{v}_{H^1(\Omega)}^2 &\coloneqq  |v|_{H^1(\Omega)}^2 + \norm{v}_{L^2(\Omega)}^2.
\end{split}
\end{equation}
In the following, we also use the abbreviation $\snorm{v} \coloneqq |v|_{H^1(\Omega)}$ for the $H^1$-semi-norm.
In the case where $V = H^1_0(\Omega)$, the semi-norm is a natural choice for the norm of $V$ since it is a norm on $V$ due to Friedrich's inequality (see \cite[Chapter II.1]{Braess}).
Assuming $f \in L^2(\Omega)$ and $A \in L^\infty(\Omega, \R^{d \times d})$, a multiplication of \eqref{eq:classical_PDE} by so-called test functions in $V$ and a simple use of integration by parts gives
\begin{equation} \label{eq:weak_components}
a(v,w) \coloneqq \int_{\Omega}^{} \left( A \nabla v \right) \cdot \nabla w \dx, \qquad\qquad l(v)\coloneqq 
\int_\Omega {f} v \dx.
\end{equation}
For verifying that $a$ and $l$ fulfill the assumptions in \cref{def:elliptic_problem}, we note that a bilinear form $a: V \times V \to \R$ is called \emph{continuous} if there exists a constant $\gamma_a>0$ such that
$$
\abs*{a(v,w)} \leq \gamma_a \, \norm{v} \cdot \norm{w} \qquad\qquad \text{for all }\, v,w \in V
$$
and \emph{coercive} if there exists a constant $\alpha_a > 0$ such that
$$
a(v,v) \geq \alpha_a \norm{v}^2 \qquad\qquad \text{for all }\, v \in V.
$$
To prove that $a$ and $l$ from \eqref{eq:weak_components} are continuous, Cauchy-Schwarz's inequality can be used.
Moreover, for the coercivity of $a$, we require additional elliptic assumptions on $A$ and Friedrich's inequality; more details are given in \cref{chap:TSRBLOD}.

There exists a large amount of literature for solving the variational equation~\eqref{eq:weak_problem}.
In this work, we are solely focusing on grid-based methods, meaning to obtain an approximation of $u \in V$ from a numerical point of view, i.e., to fully discretize the domain $\Omega$ and use a finite-dimensional subspace $V_h \subset V$.
We discuss a standard discretization approach for~\eqref{eq:weak_problem} and related techniques in \cref{sec:background_npdgl}.
Localized discretization techniques are introduced in \cref{chap:TSRBLOD}.
For now, we follow the infinite-dimensional case and introduce \emph{parameterized elliptic problems}.

\subsection{Parameterized elliptic problems}
\label{sec:param_problems}

\cref{def:elliptic_problem} is formulated for a fixed domain $\Omega$, bilinear form $a$, and right-hand side $l$.
Many physical phenomena and applications can additionally be interpreted as a \emph{parameterized problem}, meaning that at least one of the components of \cref{def:elliptic_problem} is dependent on a physical parameter $\mu \in \Params$\glsadd{params} as an element of a parameter space (or parameter set) $\Params \subset \R^P$, where $P \in \mathbb{N}$.
In this case, one is not interested in a single solution but solutions for many parameter samples from the parameter space.
Parametric problems have many (real-world) applications, e.g., shape optimization for parametric $\Omega_\mu$ (see~\cite{shapeoptimization}). However, we solely concentrate on parametric problems where only the bilinear form and right-hand side are assumed to be parametric, and hence, the function space $V$ does not change.

Moreover, we introduce the notion of output functionals associated with the parameterized problem, often used in applications, where not the solution state of the problem is the quantity of interest but instead a single number, e.g., an average over a specific part of the domain.
With particular note to the primary goal of this thesis, output functionals can also be interpreted as objective functionals for optimal control problems with PDE constraints; cf.~\cref{sec:background_PDE-constr}. 

\begin{problem}[Parameterized elliptic problem]
	\label{def:param_elliptic_problem}
	For every $\mu \in \Params$, let $a_\mu: V \times V \to \R$ denote a continuous and coercive bilinear form with continuity and coercivity constants $\cont{a_\mu} > 0$ and $\alpha_{a_\mu} > 0$, and $l_\mu \in V'$ a continuous linear functional. 
	We consider a parameterized problem as the task to find $u_\mu \in V$, for a fixed parameter $\mu \in \Params$, as the solution of
	\begin{equation}
	\label{eq:param_elliptic_problem}
	a_\mu(u_\mu,v) =  l_\mu(v) \qquad\qquad \text{for all } \, v \in V.
	\end{equation} 
	We also define the bounded map $\mathcal{S}:\Params \to V$, $\mu \mapsto u_\mu\coloneqq \mathcal S(\mu)$ as the solution map (or parameter-to-state map) of~\eqref{eq:param_elliptic_problem}.
	Furthermore, we assume to be given an output functional $\Jhat(\mu) \coloneqq \J(u_\mu, \mu)$, for an arbitrary functional $\J: V \times \Params \to \R$.
\end{problem}

We again note that the sub-index $\mu$ can, for instance, inherit from a parametric diffusion coefficient $A_\mu$, right-hand side $f_\mu$, as opposed to $A$ and $f$ given in~\eqref{eq:weak_components}, or from parametric boundary conditions that are hidden in the weak formulation.
Importantly, however, the function space $V$ is fixed.
Since this thesis intersects with concepts from optimal control theory, we emphasize that the sub-index $\mu$ does not stand for an abbreviation of partial derivatives (as is often the case in optimal control theory).
Given the assumptions in \cref{def:param_elliptic_problem}, the existence of the solution map $\mathcal{S}$ again follows from the Lax-Milgram Theorem.
This thesis does not generally assume linearity for the output functional $\J$.
Moreover, we use a slightly different notation for $\J$ concerning the parameter dependency, meaning that we do not use the abbreviation $\J_\mu$.
Furthermore, in the definition of the PDE-constrained optimization problem~\eqref{Phat} from below, the functional $\Jhat$ is introduced as the \emph{reduced objective functional}, which aligns with the notation of \cref{def:param_elliptic_problem} since the definition is the same.
However, output functionals of parameterized systems are not only restricted to PDE-constrained optimization, which is why we already include them in the general definition of \cref{def:param_elliptic_problem}.

As already mentioned, parametric problems can be tackled by model order reduction techniques.
To motivate these more precisely, one can think about finite-dimensional subspaces $V_h \subset V$, e.g., based on a grid-based discretization, cf.~\cref{sec:background_npdgl}.
These subspaces can be high-dimensional, making it costly to compute solutions numerically.
The main idea of model order reduction, particularly of RB methods, is to find carefully chosen \emph{snapshots} $u_{h,\mu} \in V_h$ as solutions of the finite-dimensional version of \cref{def:param_elliptic_problem} (cf.~\cref{def:discr_param_elliptic_problem}) and to construct low dimensional subspaces $V_{\text{red}} \subset V_h$ from those snapshots.
Consequently, solutions on these subspaces can be computed quickly for a parameter $\mu \in \Params$.
While we discuss and present these techniques in \cref{sec:background_model_order_reduction}, we already introduce a critical assumption on the parameter dependence of $a_\mu$ and~$l_\mu$ which is crucial for an \emph{online efficient} surrogate model.
~
\begin{definition}[Affine parameter separability]
	\label{def:parameter_separable}
	A form like $a_\mu$ (or $l_\mu$ and $\J$) is called parameter separable with $\Xi^a \in \N$ non-parametric components $a_\xi: V \times V \to \R$ for $1 \leq \xi \leq \Xi^a$, and respective parameter functionals $\theta_\xi^a$, if
	\begin{align}
	a_\mu(u, v) = \sum_{\xi = 1}^{\Xi^a} \theta_{\xi}^a(\mu)\, a_{\xi}(u, v).
	\end{align}
\end{definition}

\noindent Thus, parameter separability means separating the Hilbert space affinely from the parameter space.
The following assumption looks quite restrictive but is naturally satisfied in many applications.
If the assumption is not fulfilled, empirical interpolation~\cite{BarraultMadayEtAl2004} can be used to approximate the separability by suitable components.
For further examples, we refer to \cite{DHO2012, FHRS2016, CEGG2014}, for instance.

\begin{assumption}[Parameter separability of $a_\mu$, $l_\mu$ and $\J$]
	\label{asmpt:parameter_separable}
	We assume that all parts of \cref{def:param_elliptic_problem}, i.e. $a_\mu$, $l_\mu$ and $\J$, fulfill \cref{def:parameter_separable} with respective components and parameter functionals.
\end{assumption}

We emphasize that assumptions in this thesis are meant to be valid throughout the entire work.
We always explicitly clarify if assumptions are left out on purpose.

A famous example of a parameterized elliptic problem with affinely decomposed data functions is the so-called elliptic \emph{thermal block problem}, where the domain $\Omega$ is decomposed into $N_a \times N_b$ conductivity fields $\Omega_i$ yielding a parameter space $\Params = \R_{>0}^{(N_a\cdot N_b)}$.
Each parameter component $\mu_i$ for each block scales the conductivity field linearly in $A_\mu$, i.e. $\theta_i(\mu) = \mu_i$ and $A_q = 1_{\Omega_i}$.
Suited to the thermal block problem, an exemplary output functional can be defined as an $L^2$-misfit with an additional Tikhonov-regularization term, e.g.
\begin{equation}
	\label{eq:example_J}
	\mathcal{J}_{L^2}(u, \mu) \coloneqq \frac{1}{2} \int_{D}^{} (u - u^{\text{d}})^2 \dx + \frac{1}{2} 
\sum^{P}_{i=1} \mu_i^2,
\end{equation}
on a domain of interest $D \subseteq \Omega$ and with a desired state $u^{\text{d}}$.
The output functional $\mathcal{J}_{L^2}$ plays a crucial role throughout this thesis and justifies the fact that output functionals of parameterized systems, in general, can not be assumed linear.
The thermal block problem, combined with an $L^2$-misfit, is further used in the numerical experiments in \cref{sec:loc_experiments}.

For the \emph{a posteriori} 
error control of the reduced-order models, as well as for the general definition of PDE-constrained optimal control problems, it is crucial to introduce the \emph{primal residual} of \cref{def:param_elliptic_problem}.
The name \emph{primal} 
refers to the context of PDE-constrained optimization problems, where \cref{def:param_elliptic_problem} is also called the \emph{state-} or \emph{primal equation}.

\begin{definition}[Primal residual]\label{def:primal_residual}
	For given $u \in V$, $\mu \in \Params$, we introduce the primal residual $r_\mu^\pr(u) \in V'$ associated with~\eqref{eq:param_elliptic_problem} by
	\begin{align}
	r_\mu^\pr(u)[v] \coloneqq l_\mu(v) - a_\mu(u, v) \qquad\qquad \text{for all } \, v \in V.
	\label{eq:primal_residual}
	\end{align}	
\end{definition}
\noindent For a general formulation of PDE-constrained parameter optimization problems (see Problem~\eqref{P}), it is also helpful to use an abstract formulation of \cref{def:param_elliptic_problem}, which is mathematically equivalent to the residual vanishing on $V$.

\begin{problem}[Abstract formulation of a parameterized elliptic problem] 
	\label{def:abstract_param_equation}
	With the operators $\mathcal{A}_\mu: V\to V'$ defined as $\mathcal{A}_\mu(u)\coloneqq a_\mu(u,\cdot)$ and $\mathcal{L}_\mu \in V'$ defined as $\mathcal{L}_\mu\coloneqq l_\mu$, we define $e: V \times \Params \to V'$ by
	\begin{equation}
	\label{eq:abstract_state_eq}
	e(u,\mu) \coloneqq \mathcal{L}_\mu-\mathcal{A}_\mu(u).
	\end{equation}
	Then, the task: Find $u_\mu \in V$ such that
	\begin{equation}
		e(u_\mu,\mu) = 0
	\end{equation} 
	is equivalent to \cref{def:param_elliptic_problem} and, for all $u,v\in V$, we have $e(u,\mu)[v] = r^\pr_\mu(u)[v]$.
\end{problem}

For PDE-constrained parameter optimization problems, we aim to consider (directional) derivatives of operators like $e(\cdot\,,\,\cdot)$.
For this purpose, we define a formal way of computing functional derivatives on Hilbert spaces in the following section.

\subsection{Fr\'echet derivatives}

This section briefly introduces the notion of Fr\'echet- and G\^ateaux derivatives, which deduce respective derivatives of the output functional in PDE-constrained optimization problems.
The definitions are mainly based on \cite[Section 1.4]{HPUU2009}.
Fr\'echet derivatives are defined for normed spaces and serve as a generalization of the derivative of vector-valued functions to the infinite-dimensional case.
The definition also strongly connects to the G\^ateaux derivative for generalizing directional derivatives on Banach spaces.
While both terms are not equivalent (which we do not discuss in detail), we use the G\^ateaux derivative to compute the Fr\'echet derivative.
The section is also devoted to clarifying the notation of derivatives used throughout this work.
Let us start with a definition of directional derivative in Banach spaces.

\begin{definition}[Directional derivative] \label{def:directional_derivative}
	Let $X$ and $Y$ be Banach spaces and let $F: U \subset X \to Y$ be an operator with an open subset $U \neq \emptyset$.
	The directional derivative of $F$ at $u \in U$ in direction $w \in X$ is defined as
	\begin{equation}
		\D F(u)[w] = \lim_{t\to 0} \frac{F(u + tw) - F(u)}{t} \in Y.
	\end{equation}
	If the limit exists for all $w \in X$, we call $F$ directionally differentiable.
\end{definition}

\noindent Given the definition of directional differentiability, we can define the G\^ateaux- and Fr\'echet differentiability.

\begin{definition}[G\^ateaux and Fr\'echet differentiability] 
	Let $F$ be given as in \cref{def:directional_derivative}. We call $F$
	\begin{itemize}
		\item G\^ateaux differentiable at $u \in U$ if $F$ is directionally differentiable at $u$ and the directional derivative $\D F(u) \in Y$ is a bounded and linear operator.
		\item Fr\'echet differentiable at $u \in U$ if $F$ is G\^ateaux differentiable at $u$ and if the following approximation condition holds:
		\begin{equation} \label{eq:frechet_condition}
			\norm{F(u+h) - F(u) - \D F(u)[h]}_Y = o(\norm{h}_X) \qquad \text{ for }\norm{h}_X \to 0.
		\end{equation}
	\end{itemize}
\end{definition}

\noindent Note that with~\eqref{eq:frechet_condition}, Fr\'echet differentiability can also be interpreted as the existence of a bounded linear functional $\D F(u) \in Y$ such that
\begin{equation}
	\lim_{\norm{h}_X \to 0} \frac{\norm{F(u+h) - F(u) - \D F(u)[h]}_Y}{\norm{h}_X} = 0.
\end{equation}
Moreover, standard derivative properties such as linearity and the chain rule apply to the Fr\'echet derivative.
The remainder of this subsection is devoted to an exemplary computation of derivatives of bilinear forms $a_\mu$, linear functionals $l_\mu$, and solution maps of \cref{def:param_elliptic_problem}.
Thus, we require the Fr\'echet derivatives on the Hilbert space $V$ and the parameter space $\Params$.
We start with a suitable assumption that is supposed to be valid for the remainder of this thesis.

\begin{assumption}[Differentiability of $a_\mu$, $l_\mu$ and $\J$]
	\label{asmpt:differentiability}
	We assume $a_{\mu}$, $l_{\mu}$, and $\J$ from \cref{def:param_elliptic_problem} to be Fr\'echet differentiable w.r.t~each argument $u$, $v$, and $\mu$.
	Given \cref{asmpt:parameter_separable}, this follows from assuming that also the separable parts of $a_{\mu}$, $l_{\mu}$, and $\J$ are Fr\'echet differentiable.
\end{assumption}

Regarding notation, for instance, interpreting the bilinear form $a_\mu$ as a map $a_\mu: V \times V \times \Params \to \R$, $\left(u, v, \mu\right) \mapsto a_\mu(u, v)$, we denote the Fr\'echet derivatives of $a_\mu$ w.r.t.~the first, second and third argument of said map in the direction of $w\in V$, $\nu \in \R^P$ by $\partial_u a_\mu(u, v)[w] \in \R$, $\partial_v a_\mu(u, v)[w] \in \R$ and $\partial_\mu a_\mu(u, v)\cdot \nu \in \R$, respectively (noting that the dual pairing of the latter is simply the Euclidean product).
Similarly, interpreting the linear functional $l_\mu$ as a map $l_\mu: V \times \Params \to \R$, $(v, \mu) \mapsto l_\mu(v)$, we denote the Fr\'echet derivatives of $l$ w.r.t.~the first and second argument of said map in the direction of $w \in V$, $\nu \in \R^P$ by $\partial_v l_\mu(v)[w] \in \R$ and $\partial_\mu l_\mu(v) \cdot \nu \in \R$, respectively.
The same notation for (bi-)linear forms is also used for the cost functional $\J$.
We omit the word Fr\'{e}chet when referring to the derivatives of $a_\mu$, $l_\mu$, and $\J$ to simplify the notation.
For the canonical partial derivatives we may also use a shorthand, for instance $\partial_\mu l_\mu(v) \cdot e_i \coloneqq \partial_{\mu_i} l_\mu(v)$, where $e_i$ denotes the $i$-th canonical direction. 
Moreover, as usual, the vector of all canonical parameter derivatives of, e.g., $l_\mu$ is defined as the gradient $\nabla_\mu l_\mu$, where the sub-index is left out if the respective derivative component is clear.

Note that we denote the derivatives w.r.t.~the symbol of the argument in the original definition of the functional or bilinear form, not w.r.t.~the symbol of the actual argument, i.e.~we use $\partial_u a_\mu(u_\mu,v)$ for the derivative w.r.t.~the first argument, not $\partial_{u_\mu} a_\mu(u_\mu, v)$ or $\partial_v a_\mu(u, p)$ for the derivative w.r.t.~the second argument, not $\partial_p a_\mu(u, p)$.

Given \cref{asmpt:differentiability}, we can compute the directional derivatives of (bi-)linear forms $a_\mu$ and~$l_\mu$ in the following way.

\begin{proposition}[Derivatives w.r.t.~$V$]
	\label{rem:gateaux_wrt_V}
    For $u, v \in V$, $\mu \in \Params$, the derivatives of $a$ and $l$ w.r.t.~arguments in $V$ in the direction of $w \in V$ are given by
	\begin{align*}
	\partial_u a_\mu(u, v)[w] = a_\mu(w, v), \qquad\qquad 
	\partial_v a_\mu(u, v)[w] = a_\mu(u, w), \qquad\text{ and }\qquad
	\partial_v l_\mu(v)[w] 	  = l_\mu(w), 
	\end{align*}
	respectively.
\end{proposition}
 \begin{proof}
   For $\mu \in \Params$ and $u, v \in V$, we obtain
   \begin{align*}
     \partial_u a_\mu(u, v)[w] = \lim_{t \to 0} \frac{a_\mu(u + tw, v) - a_\mu(u, v)}{t}
                               = \lim_{t \to 0} \frac{a_\mu(tw, v)}{t}
                               = a_\mu(w, v),
   \end{align*}
   using the definition of the directional derivative in the first equality, and the linearity of $a_\mu$ w.r.t.~$u$ and $v$.
   The other derivatives of $a_\mu$ and $l_\mu$ can be obtained similarly.
 \end{proof}

\noindent Further, we can compute the partial derivatives of $a_\mu$ and $l_\mu$ w.r.t.~the parameter by using their separable decomposition from \cref{asmpt:parameter_separable}.

\begin{proposition}[Derivatives w.r.t.~$\Params$]
	\label{rem:gateaux_wrt_P}
	For $\mu \in \Params$, $u, v \in V$, the derivatives of $a$ and $l$ w.r.t.~$\mu$ in the direction of $\nu \in \R^P$ are given by
	\begin{align*}
	\partial_\mu a_\mu(u, v) \cdot \nu &= \sum_{\xi = 1}^{\Xi^a} \big(\partial_\mu 
\theta_\xi^a(\mu) \cdot \nu\big)\, a_\xi(u, v)&&\text{and}&&&
	\partial_\mu l_\mu(v) \cdot \nu &= \sum_{\xi = 1}^{\Xi^l} \big(\partial_\mu 
\theta_\xi^l(\mu) \cdot \nu\big)\, l_\xi(v),
	\end{align*}
	respectively, if $u$ and $v$ do not depend on $\mu$.
\end{proposition}
\begin{proof}
	This follows directly from \cref{asmpt:parameter_separable} and the linearity of the derivative.
\end{proof}

\noindent We also introduce the following shorthand notation for the derivative of functionals and bilinear forms w.r.t.~the parameter in the direction of $\nu \in \R^P$, e.g.~for $\mu \in \Params$, we introduce
\begin{align*}
\partial_\mu l_\mu &\cdot \nu \in V'&&& v \mapsto \big(\partial_\mu l_\mu \cdot \nu\big)(v) &\coloneqq 
\partial_\mu l_\mu(v)\cdot \nu&&\text{and}\\
\partial_\mu a_\mu &\cdot \nu \in (V \times V \to \R)&&& u, v \mapsto \big(\partial_\mu a_\mu \cdot 
\nu\big)(u, v) &\coloneqq \partial_\mu a_\mu(u, v)\cdot \nu,
\end{align*}
and note that $\partial_\mu l_\mu$ and $\partial_\mu a_\mu$ are continuous and separable w.r.t.~the parameter, owing to \cref{asmpt:parameter_separable}.


Fr\'echet derivatives of the solution map $\mathcal S(\mu)\coloneqq u_\mu$ of \cref{def:param_elliptic_problem} are commonly used for reduced basis methods and optimization, e.g., for constructing Taylor RB spaces that consist of the primal solution as well as their respective sensitivities (see~\cite{HAA2017}) or for computing first or second-order derivatives of output functionals, as we see later.

\begin{proposition}[Fr\'{e}chet derivative of the solution map in any direction]
	\label{prop:solution_dmu_eta}
	Considering the solution map $\mathcal S:\Params \to V$, $\mu \mapsto u_\mu$ of \cref{def:param_elliptic_problem}, we denote its Fr\'{e}chet derivative w.r.t. a direction $\eta \in \R^P$ by $d_{\eta} u_\mu \in V$, which is given as the solution of
	\begin{align}\label{eq:primal_sens}
	a_\mu(d_{\eta}u_\mu, v) = \partial_\mu r_\mu^\pr(u_\mu)[v] \cdot \eta &&\text{for all } v \in V.
	\end{align}
\end{proposition}
\begin{proof}
	Given the solution $u_\mu \in V$ of \cref{def:param_elliptic_problem} for $\mu \in \Params$, we obtain
	\begin{align*}
	d_{\eta} a_\mu(u_\mu, v) &= \partial_u a_\mu(u_\mu, v)[d_{\eta}u_\mu] 
	+ \partial_\mu a_\mu(u_\mu,v) \cdot \eta\\
	&= a_\mu(d_{\eta}u_\mu, v) + \partial_\mu a_\mu(u_\mu,v) \cdot \eta,
	\end{align*}
	using the chain rule in the first equality and \cref{rem:gateaux_wrt_V} in the second equality.
	Since~$u_\mu$ solves~\eqref{eq:param_elliptic_problem}, we obtain the desired result by differentiating~\eqref{eq:param_elliptic_problem} w.r.t.~$\mu$ in the direction~$\eta$, yielding $d_{\eta} a_\mu(u_\mu, v) = d_{\eta} l_\mu(v) = \partial_\mu l_\mu(v)\cdot \eta$ (since $v$ does not depend on $\mu$) and using the definition of the residual~\eqref{eq:primal_residual}.
	The existence and uniqueness of a solution of \eqref{eq:primal_sens} again follow from Lax-Milgram.
\end{proof}

In conclusion, apart from the different right-hand side, \eqref{eq:primal_sens} has the same form as \cref{def:param_elliptic_problem}.
The previous proposition motivates the following problem definition.

\begin{problem}[Partial derivatives of the solution map] \label{def:primal_sensitivities}
	For the $i$-th component of $\Params$, the canonical directions $\eta=e_i$ lead to the partial derivatives of the solution map given by the solution $d_{\mu_i} u_\mu \in V$, such that
	  \begin{align} \label{eq:primal_sensitivities}
	    a_\mu(d_{\mu_i}u_\mu, v) = \partial_{\mu_i} r_{\mu}^\pr(u_\mu)[v] &&\text{for all } v \in V.
	  \end{align}
	Partial derivatives of the solution map are also called sensitivities.
\end{problem}

\noindent We are now prepared to introduce PDE-constrained parameter optimization problems.

\section{Elliptic PDE-constrained parameter optimization}\label{sec:background_PDE-constr}

This section introduces PDE-constrained parameter optimization problems and discusses the existence results of locally optimal solutions.
As the primary reference for this section, we used \cite{HPUU2009}. We also recommend \cite{troltzsch2010optimal} for more details on PDE-constrained parameter optimization problems.
We also note that, although the theory in \cref{sec:background_problem_formulation_and_existence} is quite general, in the sequel of the section, we solely consider optimization problems where the underlying PDE is a parameterized elliptic problem in the sense of \cref{def:param_elliptic_problem}.
Moreover, in this section, we stay in the infinite-dimensional setting.

\subsection{Problem formulation and existence result}
\label{sec:background_problem_formulation_and_existence}
This work is solely concerned with the particular case of a PDE-constrained parameter optimization problem, the definition of which is detailed in \cref{def:parameter_optimal_control}.
The theory is based on general PDE-constrained optimization (or optimal control) problems.
Hence, we start with a general formulation of an optimal control problem with PDE constraints, defined on arbitrary Banach spaces $X$, $Y$, and $Z$.

\begin{definition}[General PDE-constrained optimization problem, cf.~{\cite[Section 1.5.2]{HPUU2009}}]
	\label{def:general_optimal_control}
	Let $X$ and $Y$ be reflexive Banach spaces and let $Z$ be a Banach space.
	Let $J: X \times Y \to \R$ and $E: X \times Y \to Z$ be continuous.
	Then, we consider the PDE-constrained optimization problem as finding a locally optimal solution $(\bar{x}, \bar{y})$ of
	\begin{equation} 
		\min_{(x,y)\in X\times Y} J(x, y), \qquad\qquad \text{ subject to } \qquad\qquad 
E(x, y) = 0, \qquad x \in X_{ad}, y \in Y_{ad}, \label{P_gen} \tag{$\textnormal{P}_{\textnormal{gen}}$} 
	\end{equation}
	where $X_{ad}$ and $Y_{ad}$ are "admissible" subsets of $X$ and $Y$, respectively.
\end{definition}
\noindent Note that, due to the mild assumptions on $E$, this problem is considered a nonlinear optimization problem, which is the typical case for PDE-constrained optimization problems.
Furthermore, the objective functional is generally non-convex and can have many local minima.
Thus, the existence of a uniquely defined solution of \eqref{P_gen} can not be expected.
For a corresponding existence result, we follow~\cite{HPUU2009} to introduce the following assumptions:

\begin{assumption}[Admissibility of Problem~\eqref{P_gen}, cf.~{\cite[Assumption 1.44]{HPUU2009}}]
	\label{asmpt:admissibility_of_control_problem} \hfill \\
\hphantom{space} 1. The control space $Y_{ad} \subset Y$ is convex, bounded and closed. \\
\hphantom{space} 2. The state space $X_{ad} \subset X$ is convex and closed, such that 
\eqref{P_gen} has an admissible (or \\ \hphantom{space 1. } feasible) point. \\
\hphantom{space} 3. The state equation $E(x, y) = 0$ has a bounded solution operator $y \in Y_{ad} 
\mapsto x(y) \in X_{ad}$. \\
\hphantom{space} 4. The map $(x, y) \in X \times Y \mapsto E(x, y) \in Z$ is continuous under weak 
convergence. \\
\hphantom{space} 5. $J$ is sequentially weakly lower semi-continuous.
\end{assumption}
\noindent The existence result for solutions of Problem~\eqref{P_gen} can then be stated as follows.

\begin{theorem}[Existence of optimal solutions, cf.~{\cite[Theorem 1.45]{HPUU2009}}] 
\label{thm:existence_optimal_control}
	Let \cref{asmpt:admissibility_of_control_problem} be true.
	Then Problem~\eqref{P_gen} has an optimal solution $(\bar{x}, \bar{y}) \in X \times Y$.
\end{theorem}

\noindent Note that the (nonlinear) formulation $E(x, y) = 0$ is set to be a general PDE-based equality constraint with arbitrary control dependence.
We refer to~\cite[Chapter 1]{HPUU2009} for examples.
Following the preliminaries in \cref{sec:preliminaries}, we now state what we call a PDE-constrained parameter optimization problem.
Recall that $V$ denotes a real-valued Hilbert space.

\begin{definition}[PDE-constrained parameter optimization problem]
	\label{def:parameter_optimal_control}
	Let $\J: V \times \Params \to \R$ be a continuous objective functional and let $e: V \times \Params \to V'$ represent the parametric state equation as in \cref{def:abstract_param_equation}.
	Then, we consider the PDE-constrained parameter optimization problem as to find a locally optimal solution $(\bar{u}, \bar{\mu})$ of
	\begin{equation}
	\min_{(u,\mu)\in V\times \Params} \J(u, \mu), \qquad\qquad \text{ subject to }
	\qquad\qquad e(u, \mu) = 0, \qquad u \in V, \mu \in \Paramsad, \tag{$\textnormal{P}$} 
\label{P}
	\end{equation}
	with a convex, bounded and closed admissible parameter space $\Params \subset \R^P$.
\end{definition} 

The phrase "parameter" refers to the control variable being hidden in the state equation as a parameter $\mu$.
Note that our notation substantially differs from the source~\cite{HPUU2009}, which is clarified in the following.
Despite the tendency of optimal control theory to name the control variable as $u$, we emphasize that, in the context of RB methods, the name $u$ is always reserved for the solution $u_\mu$ of the parameterized state equation~\eqref{eq:param_elliptic_problem}.
Moreover, in PDE-constrained optimal control problems, the control is typically denoted as a function in space and/or time. In our application, the control space is the parameter space of the parameterized elliptic problem and, thus, a subset of $\R^P$.
To conclude, in what follows, $\mu \in \Params$ is used as the control variable with control space $\Params$.
In the remainder of this thesis, the term "control space" and "parameter space" can thus be used equivalently for $\Params$.

Typical choices for constraints on the parameter space are simple bound constraints, i.e. $\Paramsad = [\mu_a,\mu_b]$, with bounds $\mu_a, \mu_b \in \R^P$, cf.~\cref{sec:constr_optimization}.
We formulate the corresponding existence result as a consequence of \cref{thm:existence_optimal_control}.

\begin{corollary}[Existence of optimal solutions] \label{cor:existence_of_param_opt_control}
	Problem~\eqref{P} admits an optimal solution $(\bar{u}, \bar{\mu}) \in V\times\Params$, where $ \bar{u} = u_{\bar{\mu}} = \mathcal{S}(\bar{\mu})$.
\end{corollary}
\begin{proof}
	To use \cref{thm:existence_optimal_control}, Problem~\eqref{P} has to fulfill \cref{asmpt:admissibility_of_control_problem}, where $Y_{ad} \coloneqq \Paramsad$, $X_{ad} \coloneqq V$, $E \coloneqq e$ and $J\coloneqq\J$.
	While 1. is given by definition, 2. is automatically given for a Hilbert space $V$.
	For 3., we refer to the solution map $\mathcal{S}$ from \cref{def:param_elliptic_problem}, and for an elaborated discussion on 4. and 5., we refer to~\cite[Theorem 1.45]{HPUU2009}.
\end{proof}

\noindent We note that the assumptions in Problem~\eqref{P} are not fulfilled if no constraints on the control space $\Params$ are enforced, i.e., $\Params = \R^P$.
In that case, additional assumptions on $\J$ are enforced for an existence result, which we do not consider further in this thesis (see~\cite{HPUU2009}).

Since \cref{def:param_elliptic_problem} admits a unique solution and we have $\bar{u} = u_{\bar{\mu}}$ for the optimal solution, Problem~\eqref{P} can be reformulated as a so-called \emph{reduced} optimization problem, cf.~\cite[Chapter 1.6]{HPUU2009}.
\begin{definition}[Reduced PDE-constrained parameter optimization problem]
	\label{def:red_parameter_optimal_control}
	For $\J$, $e$ and $\Paramsad$ from Problem~\eqref{P}, let the reduced cost functional be defined as
	$$\Jhat: \Params\mapsto \mathbb{R},\,\mu\mapsto \Jhat(\mu)
	\coloneqq \J(u_\mu, \mu)= \J( \mathcal S(\mu),\mu).$$
	Then, the reduced problem
	\begin{align}
	\min_{\mu \in \Paramsad} \Jhat(\mu)
	\tag{$\hat{\textnormal{P}}$}\label{Phat}
	\end{align}
	 is equivalent to~\eqref{P}.
\end{definition}
\noindent In this formulation case, the term "reduced" is not related to reduced models in the MOR sense but is solely referring to the fact that the optimization functional is only considered on $\Params$.
The important difference to \eqref{P} is that the PDE constraints are no longer explicitly present but hidden in the functional.
Therefore, derivatives are more involved and no separate consideration of $u$ and $\mu$ is possible.

The contrary approach of solving \eqref{P} without a reduced formulation are so-called \emph{all-at-once} approaches, which approximate a solution of \eqref{P} without explicit use of the solution map $\mathcal{S}$.
These methods are usually entirely based on optimality conditions, cf.~\cref{sec:background_optimality_conditions}, which make them particularly efficient.
Example of all-at-once methods are, for instance, Gau{\ss}-Newton-type or sequential quadratic programming (SQP) approaches, cf.~\cite[Chapter 18]{Nocedal} or \cite{Clever2012,ganzler2006sqp,HV06}.

To conclude,~\eqref{Phat} and \eqref{P} define an elliptic PDE-constrained parameter optimization problem that we aim to solve numerically.
If not specified differently, we refer to~\eqref{Phat} as the PDE-constrained parameter optimization problem, which, in this thesis, always exists and is equivalent to~\eqref{P}.

While \cref{cor:existence_of_param_opt_control} ensures the existence of a solution for the constrained case, it does not provide a method to compute it.
It is not even clear how we can evaluate the objective functional $\Jhat$ as this involves solving the parametric elliptic state equation~\eqref{eq:param_elliptic_problem} for which we did not provide a solution method yet.
However, we want to anticipate that computing an approximation of the solution of~\eqref{eq:param_elliptic_problem} can be very demanding, and thus, the computational effort for evaluations of $\Jhat$ can not be neglected.
On the contrary, this thesis is particularly concerned about problems where evaluating $\Jhat$ with standard methods can even be prohibitively costly, aiming for alternative solution methods; cf.~\cref{chap:TR_TSRBLOD}.
We focus on the discretization of the state equation and the respective solution method in \cref{sec:background_npdgl}. First, we proceed with the infinite-dimensional case and note that, as far as the thesis is concerned, the results can easily be transferred to the finite-dimensional case.
We now discuss derivatives of $\Jhat$ and optimality conditions as the basis for deriving iterative optimization methods.

\subsection{Derivatives of the objective functional}
\label{sec:background_derivatives_of_J}

Standard solution methods for arbitrary optimization problems use derivative information of~$\Jhat$.
In the simplest case, the gradient, whose negative value is the steepest descent direction, can be approximated numerically by, e.g., finite differences or stochastic approaches.
In such cases, only evaluations of~$\Jhat$ are required without direct access to $\nabla_\mu \Jhat$.
However, firstly, the accuracy of these approaches can be inferior, and secondly, these approaches require evaluating~$\Jhat$ unnecessarily often, which can become a computational bottleneck (as shortly discussed before).
Due to that, we do not consider approaches that only require evaluations of $\Jhat$ and instead compute the gradient information explicitly (or approximate it in other ways).
There are several ways to derive the exact gradient of $\Jhat$.
One possibility is to use partial derivatives of the state equation as given in \cref{def:primal_sensitivities}.
To be precise, by using the chain rule, we have:
\begin{equation}
\begin{split}\label{eq:sensitivity_gradient} 
(\nabla_\mu \Jhat(\mu))_i &= \partial_{\mu_i} \J(u_{\mu}, \mu) + \partial_u \J(u_{\mu}, 
\mu)[d_{\mu_i} u_{\mu}].
\end{split}	
\end{equation}

\noindent Consequently, computing the full gradient of $\Jhat$ with the sensitivity approach means to solve \cref{def:primal_sensitivities} for every canonical direction in $\R^P$.
The computational complexity of \cref{def:primal_sensitivities} scales with the complexity of \cref{def:param_elliptic_problem} and thus, the overall complexity of $\nabla_\mu \Jhat$ heavily scales with the dimension of $\Params$.
The so-called \emph{adjoint approach} is a more suitable way to obtain~$\nabla_\mu \Jhat$.
Here, we use the dual equation, associated with \eqref{Phat}, as an auxiliary problem:

\begin{problem}[Dual or adjoint equation] \label{def:dual_solution}
	 For a fixed $\mu \in \Params$, given the solution $u_\mu \in V$ of \cref{def:param_elliptic_problem}, we define the dual solution $p_{\mu} \in V$ by	
	\begin{align} \label{eq:dual_solution}
	a_\mu(q, p_\mu) = \partial_u \J(u_\mu, \mu)[q]
	&&\text{for all } q \in V.
	\end{align}
\end{problem}
\noindent A solution of \cref{def:dual_solution} is also called \emph{Lagrangian multiplier}, which is detailed in \cref{sec:background_optimality_conditions}.
For now, we can use the dual equation to rewrite $\nabla_\mu \Jhat$ in the following way:
\begin{equation}
\begin{split}\label{eq:adjoint_gradient_derivation} 
(\nabla_\mu \Jhat(\mu))_i &= \partial_{\mu_i} \J(u_{\mu}, \mu) + \partial_u \J(u_{\mu}, 
\mu)[d_{\mu_i} u_{\mu}] \\
&= \partial_{\mu_i} \J(u_{\mu}, \mu) + a_\mu(d_{\mu_i} u_{\mu}, p_\mu) \\
&= \partial_{\mu_i} \J(u_{\mu}, \mu) + \partial_{\mu_i} r_{\mu}^\pr(u_\mu)[p_{\mu}],
\end{split}	
\end{equation}
where we have used~\eqref{eq:sensitivity_gradient} in the first,~\eqref{eq:dual_solution} in the second, and~\eqref{eq:primal_sensitivities} in the third equality.
In conclusion, to obtain $\nabla_\mu \Jhat$, we no longer need to solve $\dim{\Params}$ many sensitivity equations, but only one single, dual equation~\eqref{eq:dual_solution}.
The evaluation of the partial derivatives of the primal residual $\partial_{\mu_i} r_{\mu}^\pr$, instead, are cheap.
Thus, we have:

\begin{proposition}[Gradient of $\Jhat$] \label{prop:gradient_Jhat}
	The gradient of $\Jhat$ is given as
\begin{equation}\label{eq:adjoint_gradient} 
	\nabla_\mu \Jhat(\mu) = \nabla_\mu \J(u_{\mu}, \mu) + \nabla_\mu r_{\mu}^\pr(u_\mu)[p_{\mu}].
	\end{equation}
\end{proposition}
\begin{proof} See~\eqref{eq:adjoint_gradient_derivation}.
\end{proof}

We remark that~\eqref{eq:adjoint_gradient_derivation} relies on the fact that both $u_\mu$ and $p_\mu$ belong to the same space $V$; cf.~\cite{HPUU2009}.
For a non-conforming choice of the dual and primal space,~\eqref{eq:adjoint_gradient} does not hold, which plays a significant role in the construction of the surrogate models in \cref{chap:TR_RB}.
To derive the dual equation and to prove~\eqref{eq:adjoint_gradient}, we can also use the Lagrangian point of view of optimization problems, which is closely related to the first-order necessary optimality conditions for~\eqref{Phat}; cf.~\cref{sec:background_optimality_conditions}.
We emphasize that, while the formulation in~\eqref{eq:adjoint_gradient} seems advantageous, it is not always possible to solve \cref{def:dual_solution} since the dual of $a_\mu$ may not be accessible, for instance, if the primal system is solved with a black-box solver.

For sufficient optimality conditions, we also require information on the second derivative of~$\Jhat$, i.e., the Hessian of~$\Jhat$.
However, from an optimization point of view, we do not require the Hessian component-wise, but only the application to a particular direction $\nu \in \R^P$.
Since the Hessian requires sensitivities of both $u_\mu$ and $p_\mu$, we also define the sensitivities for the dual solutions.

\begin{proposition}[Directional derivative of the dual solution map at any direction]
	\label{prop:dual_solution_dmu_eta}
	Considering the dual solution map $\Params \to V$, $\mu \mapsto p_\mu$, where $p_\mu$ is the solution of \cref{def:dual_solution}, we denote its directional derivative w.r.t. a direction $\eta$ by $d_{\eta} p_\mu \in V$, which is given as the solution of 
	\begin{equation} \label{eq:dual_sens}
	a_\mu(q, d_{\eta} p_\mu) = -\partial_\mu a_\mu(q, p_\mu)\cdot \eta 
	+ \partial_\mu \partial_u \J(u_\mu, \mu)[q] \cdot \eta \qquad\qquad
	\text{for all }\, q \in V.
	\end{equation}
\end{proposition}
 \begin{proof}
 	Using \cref{rem:gateaux_wrt_V}, the result follows along the same lines as in the proof of \cref{prop:solution_dmu_eta}.
 \end{proof} 

\noindent Corresponding dual sensitivities with the canonical directions $\eta:=e_i$ can be defined analogously to \cref{def:primal_sensitivities}.
Note that the right-hand side of \cref{prop:dual_solution_dmu_eta} is not just the derivative of the residual of the dual equation, like it is the case for the primal sensitivity in \cref{def:primal_sensitivities}, since the right-hand side of \eqref{eq:dual_solution} depends on $u_\mu$.

\begin{proposition}[Application of the Hessian of $\Jhat$ to a direction $\eta$]
	\label{prop:Hessian_Jhat}
	Given a direction $\eta \in \R^P$, we have
	\begin{align}
	\cHhat(\mu) \cdot \eta
	= \nabla_\mu\Big(&\partial_u \J(u_\mu, \mu)[d_{\eta}u_\mu]
	+ l_\mu(d_{\eta}p_\mu) - a_\mu(d_{\eta}u_\mu, p_\mu) - a_\mu(u_\mu, d_{\eta}p_\mu)
	\notag\\
	&+\partial_\mu\big(\J(u_\mu, \mu)+ l_\mu(p_\mu)- a_\mu(u_\mu, p_\mu)\big)\cdot \eta\Big),
	\notag
	\end{align}
	where $u_\mu, p_\mu \in V$ denote the primal and dual solutions, respectively. 
\end{proposition}
\begin{proof}
	Given two directions $\eta,\nu\in \R^P$, we have
	\begin{align*}
	\big(\cHhat(\mu) \cdot \eta\big) \cdot \nu &= d_{\eta}\, \big(\partial_\mu \Jhat(\mu) 
\cdot \nu\big) \\
	&= \underbrace{d_{\eta} \partial_\mu \J(u_\mu, \mu)\cdot \nu}_{=: (i)} + 
\underbrace{d_{\eta}\partial_\mu l_\mu(p_\mu)\cdot \nu,}_{=: (ii)} - \underbrace{d_{\eta} 
\partial_\mu a_\mu(u_\mu, p_\mu)\cdot \nu}_{=: (iii)}
	\end{align*}
	using the definition of the Hessian, leaving us with three terms.
	For all terms, we can use the chain rule, the fact that we can exchange differentiation w.r.t.~$V$ and $\Params$, and \cref{rem:gateaux_wrt_V}.
	Regarding $(i)$, we obtain
	\begin{equation}
	(i) = \partial_u\big(\partial_\mu \J(u_\mu, \mu)\cdot \nu\big)[d_{\eta}u_\mu] + 
\partial_\mu\big(\partial_\mu \J(u_\mu, \mu)\cdot \nu\big)\cdot \eta.
	\end{equation}
	For $(ii)$, we have
	\begin{align*}
	(ii) &= \underbrace{\partial_v \big(\partial_\mu l_\mu(p_\mu)\cdot 
\nu\big)[d_{\eta}p_\mu]}_{= \partial_\mu \big(\partial_v l_\mu(p_\mu)[d_{\eta}p_\mu]\big)\cdot \nu} 
+ \partial_\mu\big(\partial_\mu l_\mu(p_\mu)\cdot \nu\big)\cdot \eta \\
	&= \partial_\mu\big( l_\mu(d_{\eta}p_\mu) + \partial_\mu l_\mu(p_\mu)\cdot \eta\big)\cdot 
\nu
	\end{align*}
	and concerning $(iii)$, it holds
	\begin{align*}
	(iii) &= \underbrace{\partial_u\big(\partial_\mu a_\mu(u_\mu, p_\mu) \cdot 
\nu\big)[d_{\eta}u_\mu]}_{\partial_\mu a_\mu(d_{\eta}u_\mu, p_\mu)\cdot \nu} \\
	&+ \underbrace{\partial_v\big(\partial_\mu a_\mu(u_\mu, p_\mu) \cdot 
\nu\big)[d_{\eta}p_\mu]}_{\partial_\mu a_\mu(u_\mu, d_{\eta}p_\mu)\cdot \nu} + 
\partial_\mu\big(\partial_\mu a_\mu(u_\mu, p_\mu) \cdot \nu\big)\cdot \eta, \\
	&= \partial_\mu\Big(a_\mu(d_{\eta}u_\mu, p_\mu) + a_\mu(u_\mu, d_{\eta}p_\mu) + 
\partial_\mu a_\mu(u_\mu, p_\mu) \cdot \eta\Big)\cdot \nu.
	\end{align*}
	Now, since the direction $\nu$ is arbitrary, we obtain the claim, considering the Riesz-representative of the directional derivative.
\end{proof}

\noindent To conclude, computing the Hessian of $\Jhat$ also requires solving for sensitivities of $u_\mu$ and $p_\mu$, which makes the computational effort proportional to the dimension of $\Params$.

\subsection{Optimality conditions}
\label{sec:background_optimality_conditions}

The optimality conditions that we use for the problem class of PDE-constrained optimization problems are strongly related to the Karush–Kuhn–Tucker (KKT) conditions, see \cite{HPUU2009}, for instance.
In the general case of a nonlinear and non-convex optimization problem~\eqref{P}, first-order necessary conditions for a (locally) optimal solution can be derived by the Lagrangian functional.

\begin{definition}[Lagrangian functional, see~\cite{HPUU2009}] \label{def:lagrangian_functional}
	The Lagrangian functional of~\eqref{P} is defined by
	\begin{equation}\LL(u,\mu,p) = \J(u,\mu) + r_\mu^\pr(u)[p] \label{eq:lagrangian_functional}
	\end{equation}
	for $(u,\mu)\in V\times\Params$ and for $p\in V$.
\end{definition}

\noindent In what follows, we formulate first-order necessary optimality conditions for \eqref{Phat}.

\begin{proposition}[First-order necessary optimality conditions, cf.~{\cite[Cor. 1.3]{HPUU2009}}]
	\label{prop:first_order_opt_cond}
	Let $(\bar u, \bar \mu) \in V \times \Paramsad$ be a local optimal solution to~\eqref{P}.
	Moreover, let \cref{asmpt:differentiability} hold.
	Then, there exists an associated unique Lagrange multiplier $\bar p\in V$, such that the following first-order necessary optimality conditions hold:
	\begin{subequations}
		\label{eq:optimality_conditions}
		\begin{align}
		r_{\bar \mu}^\pr(\bar u)[v] &= 0 &&\text{for all } v \in V,
		\label{eq:optimality_conditions:u}\\
		\partial_u \J(\bar u,\bar \mu)[v] - a_\mu(v,\bar p) &= 0 &&\text{for all } v \in V,
		\label{eq:optimality_conditions:p}\\
		(\partial_\mu \J(\bar u,\bar \mu)+\partial_{\mu} r^\pr_{\bar\mu}(\bar u)[\bar p]) 
\cdot (\nu-\bar \mu) &\geq 0 &&\text{for all } \nu \in \Params. 
		\label{eq:optimality_conditions:mu}
		\end{align}	
	\end{subequations}
	The tuple $(\bar u, \bar \mu) \in V \times \Paramsad$ is called a first-order stationary (or critical) point.
\end{proposition}
\begin{proof} For a complete proof, we refer to~\cite[Cor. 1.3]{HPUU2009}.
	In order to obtain conditions~\eqref{eq:optimality_conditions:u} -~\eqref{eq:optimality_conditions:mu} from the Lagrangian functional $\LL$ from \cref{def:lagrangian_functional}, we reformulate the optimization problem~\eqref{P} as a first-order problem:
	\begin{equation}
			0 = \nabla \LL = \begin{pmatrix}  
		 	\nabla_u \LL,\nabla_\mu \LL ,  \nabla_q \LL
		 \end{pmatrix}^T =: F.
	\end{equation}
	Therefore, we seek the solution $(\bar{u},\bar{\mu},\bar{\lambda})$ of the first-order equality 
	\begin{equation}
		F(\bar{u},\bar{\mu},\bar{\lambda}) = 0. 
	\end{equation} 
	For $u, p, q \in V$, $\mu \in \Params$, we obtain
	\begin{equation} \label{eq:L_du}
	\begin{split} 
	  \nabla_u \LL(u, \mu, p)[q] &= \partial_u \J(u, \mu)[q] - \partial_u a_\mu(q, p)[v] \\
	                             &= \partial_u \J(u, \mu)[q] - a_\mu(q, p),
	\end{split}
	\end{equation}	
	using the definition of $\LL$, the primal residual~\eqref{eq:primal_residual} and the fact that $l$ does not depend on $u$ in the first equality, and \cref{rem:gateaux_wrt_V} in the second equality. We also see that
	\begin{equation} \label{eq:L_dp}
	\begin{split}
	\nabla_p \LL(u, \mu, p)[w] &= \partial_v l_\mu(p)[w] - \partial_v a_\mu(u, p)[w] \\
	&= l_\mu(w) - a_\mu(u, w) \\
	&= r_\mu^\pr(u)[w],
	\end{split}
	\end{equation}	
	using the definition of $\LL$, the primal residual~\eqref{eq:primal_residual} and the fact that $\J$ does not depend on $p$ in the first equality, \cref{rem:gateaux_wrt_V} in the second equality, and the definition of the primal residual~\eqref{eq:primal_residual} in the third equality.
	For the parameter derivative of $\LL$, we have
	\begin{equation} \label{eq:L_dmu}
	     \nabla_\mu \LL(u, \mu, p)\cdot \nu = \partial_\mu \J(u, \mu)\cdot \nu
	     + \partial_{\mu} r^\pr_\mu (u)[p] \cdot \nu,
	\end{equation}
	which is simply a consequence of the definition of $\LL$ and the primal residual~\eqref{eq:primal_residual}.
	From ~\eqref{eq:L_du} and~\eqref{eq:L_dp}, we deduce~\eqref{eq:optimality_conditions:p} and~\eqref{eq:optimality_conditions:u}, respectively.
	For~\eqref{eq:optimality_conditions:mu}, it is essential to account for the fact that we have restricted the control space to a closed and convex admissible control set $\Paramsad$.
	By \cite[Theorem 1.48]{HPUU2009}, we know that $\bar{\mu}$ fulfills the variational equality
	\begin{equation} \label{eq:variational_equality}
			(\nabla \Jhat(\bar{\mu})) \cdot (\nu-\bar \mu) \geq 0 \qquad\qquad \text{for all }\, \nu \in \Params.
	\end{equation}
	By definition of $\LL$, evaluating the Lagrangian functional at the tuple $(u_\mu, \mu, p)$, with $\mu\in\Params, p \in V$ and $u_\mu$ solving~\eqref{eq:param_elliptic_problem} for $\mu \in \Params$, is equivalent to evaluating the reduced objective functional $\Jhat$ at $\mu$, which we see by
	\begin{align} \label{eq:lag_equals_red}
	\LL(u_\mu, \mu, p) = \J(u_\mu,\mu) + 0 = \Jhat(\mu).
	\end{align} 
	Further, since $u_\mu$ is dependent on $\mu$, we have
	\begin{equation}\label{eq:alternative_gradient_derivation}
		(\nabla_\mu \Jhat(\mu))_i = (\nabla_\mu \LL(u_\mu, \mu, p))_i = \partial_u \mathcal{L}(u_{\mu},\mu,p)[d_{\mu_i}u_{\mu}]
		+ \partial_{\mu_i} \mathcal{L}(u_{\mu},\mu,p).
	\end{equation}
	From $\nabla_p \LL = 0$ and~\eqref{eq:L_dp}, we observe that $\bar{u} = \mathcal S(\bar{\mu}) = u_{\bar{\mu}}$ and from $\nabla_u \LL = 0$ we conclude that $\partial_u \mathcal{L}(\bar{u},\bar{\mu},\bar{p})=0$.
	Thus, the first term of the right-hand side in~\eqref{eq:alternative_gradient_derivation} vanishes for the optimal tuple $(\bar{u},\bar{\mu},\bar{p})$.
	Hence, we obtain~\eqref{eq:optimality_conditions:mu} by combining~\eqref{eq:alternative_gradient_derivation},~\eqref{eq:L_dmu} and~\eqref{eq:variational_equality}.
\end{proof}

\noindent In conclusion, a stationary point of~\eqref{P} fulfills the first-order optimality conditions~\eqref{eq:optimality_conditions}.
We also see that~\eqref{eq:optimality_conditions:u} corresponds to the primal equation~\eqref{eq:param_elliptic_problem}, whereas from~\eqref{eq:optimality_conditions:p}, we obtain the dual equation~\eqref{eq:dual_solution}.
Equation~\eqref{eq:optimality_conditions:mu} further characterizes the stationary point in case of additional parameter constraints and includes the gradient of $\Jhat$.
In the proof, we also demonstrated an alternate way of deriving the gradient of $\Jhat$ as in~\eqref{eq:adjoint_gradient}, namely to combine~\eqref{eq:L_dmu} and~\eqref{eq:alternative_gradient_derivation}.
 
The stated necessary optimality conditions provide a strong basis for developing iterative optimization methods to find stationary points of~\eqref{P}.
However, to verify that a stationary point $(u, \mu)$ that satisfies~\eqref{eq:optimality_conditions} is a locally optimal solution of~\eqref{P}, we also require sufficient conditions for excluding saddle points or local maxima.
These are formulated in the following:

\begin{proposition}[Second-order sufficient optimality conditions]\label{prop:second_order}
	Let \cref{asmpt:differentiability} hold true.
	Suppose that $\bar \mu\in \Params$ satisfies the first-order necessary optimality conditions~\eqref{eq:optimality_conditions}.
	If $\cHhat(\bar \mu)$ is positive definite on the \emph{critical cone} $\mathcal C(\bar\mu)$ at $\bar\mu\in\Paramsad$, i.e., if $\nu \cdot(\cHhat(\bar \mu)\cdot \nu) > 0$ for all $\nu\in\mathcal C(\bar\mu)\setminus\{0\}$, with
	\begin{align*}
	\mathcal C(\bar\mu)\coloneqq \big\{\nu\in\Params\,\big|\, \exists \mu\in\Paramsad,\,c_1>0:
	\nu = c_1(\mu-\bar \mu),\, \nabla_\mu \Jhat(\bar\mu)\cdot \nu = 0  \big\},
	\end{align*}
	then $\bar \mu$ is a strict local minimum of~\eqref{Phat}.
\end{proposition} 
\begin{proof}
	For this result we refer to~\cite{CasTr15,Nocedal}, for instance.
\end{proof}
\noindent For so-called small residual problems (i.e, $\|\partial_u \J(\bar u,\bar \mu)\|_V$ is small), one can ensure that the second-order sufficient optimality conditions hold.
This can be proven analogously to~\cite[Section 3.3]{Vol2001}.

\subsection{First-optimize vs. first-discretize}
\label{sec:first_opt_first_dis}

This section tackled PDE-constrained parameter optimization problems from the infinite-dimensional perspective and highlighted the coherence to optimality conditions, which were derived using the Lagrangian approach.
Certainly, for a numerical solution method of Problem \eqref{P}, we require a suitable discretization scheme for the involved PDEs.
There still exist many different strategies to proceed with the optimization procedure.
Roughly, we can divide these into the \emph{first-optimize-then-discretize} and \emph{first-discretize-then-optimize} perspectives.
The latter means that the first step after formulating \eqref{P} is to fully discretize all involved quantities and utilize methods that are not necessarily connected to an infinite-dimensional theory like the optimality conditions in \cref{sec:background_optimality_conditions}.
Usually, in such a case, Lagrangian approaches are not considered, no dual problems are used, and the gradient information is solely based on the discrete point of view.

Instead, as also used in the sequel, the \emph{first-optimize-then-discretize} approach tries to minimize the effect of the discrete setting on the optimization problem, and, for instance, with the Lagrangian approach, most of the considerations are first done in the infinite-dimensional setting.
Such a strategy is also done in the \emph{all-at-once} approach, where the whole optimality system is discretized and solved simultaneously.
Since we are interested in an iterative procedure of solving \eqref{P}, we now introduce iterative optimization methods for general objective functions.
Finally, in \cref{sec:background_npdgl}, we proceed with the numerical approximation of PDEs.

\section{Iterative optimization methods}
\label{sec:iterative_optimization}

The entirety of theory and algorithms for numerically solving optimization problems like~\eqref{Phat} can not be covered in this thesis.
We focus on the essential techniques and point to~\cite{Nocedal} for further reading.
In the sequel, we review iterative optimization methods for finding a stationary point of a general
optimization task, i.e., finding an optimal point in the input- or control space of a functional $J$ with arbitrary constraints and complexity:
\begin{align}
	\min_{\mu \in \Paramsad} J(\mu)
	\tag{$\textnormal{P}_\textnormal{simple}$}\label{P_simple},
\end{align}
where the constraints may be hidden in the functional or in the input space $\Params$, cf.~\eqref{Phat}.

Many optimization tasks can't deduce an optimal solution analytically.
In such a case, an intuitive way of finding an optimum is to evaluate $J$ at a sufficient amount of points that fulfill the constraints, trying to replicate the landscape of the objective functional as well as possible.
For smooth and sufficiently convex examples, this procedure may be helpful.
However, for arbitrary objective functionals, this would mean evaluating $J$ almost infinitely often since it is never sure that the optimum has been found.
On the contrary, the optimum may be only visible within a finer resolution in the input space.
Moreover, in cases where the \emph{curse of dimensionality} is present, meaning that the control space is high-dimensional, such procedures are very inefficient.
The bottleneck even amplifies if single evaluations of $J$ are expensive (as it, e.g., is the case with PDE-constraints for $\Jhat$ in~\eqref{Phat}).

Thus, there is a need for using iterative algorithms to evaluate the functional $J$ only at "useful" points.
For this, we assume to be given an initial guess $\mu^0$ in the control space of $J$ and iteratively update this point once we find a sufficiently lower value of $J$.
There exists a vast literature concerning iterative optimization methods and their fundamentals.
For a classical book, we refer to~\cite{kelley}.

We emphasize that, just as Problem \eqref{Phat}, Problem~\eqref{P_simple}, is not generally convex and does not have a unique optimal point.
Instead, we aim to search local minima and cannot expect a global minimum.
Techniques for tackling global minima need to be more sophisticated and may have immense computational complexity.
Therefore, in this thesis, we are only concerned with finding a local minimum of an optimization task.

In the sequel, we briefly explain iterative optimization methods for solving \eqref{P_simple} in the unconstrained and constrained case of the control space.
Moreover, we shortly introduce trust-region methods.

\subsection{Optimization methods for unconstrained optimization}
\label{sec:unconstrained_optimization}
For this subsection, we assume, in contrast to \cref{sec:background_PDE-constr}, that we are solving an unconstrained optimization problem \eqref{P_simple} without any constraints on the control, i.e., $\Paramsad \coloneqq \R^P$.
In the easiest case, searching for an optimum of the differentiable objective functional $J$ is equal to finding a first-order stationary point $\bar{\mu}$ that fulfills the necessary ´first-order condition
\begin{equation}\label{simple_foc}
	\norm{\nabla J(\bar{\mu})} = 0,
\end{equation}
and an appropriate sufficient second-order condition (similar to \cref{prop:second_order}).
In this section, we briefly review the steepest-descent, Newton and BFGS methods for tackling unconstrained minimization tasks and combine them with line-search methods.
We then transfer these methods to the constrained case in \cref{sec:constr_optimization}.

A descending optimization method can be described by an iterative sequence $(\mu^{(k)})_k$, where, in each iteration $k$, a descent direction $d^{(k)} \in \R^P$ is obtained and used as a search direction, combined with an appropriate line-search technique.
For the given direction $d^{(k)}$, an appropriate step size $\stepsizek > 0$ to update $\mu^{(k)}$ needs to be found following the associated line, i.e.
\begin{equation}\label{eq:iteration_procedure}
	\mu^{(k+1)} = \mu^{(k)} + \stepsizek d^{(k)}.
\end{equation}
The descent direction $d^{(k)}$ as well as the step size $\stepsizek$ can be determined in several ways.
While extensive literature exists considering choices of $d^{(k)}$ and $\stepsizek$ that enable differently robust, convergent, and fast algorithms, we limit this section to three popular ways of computing $d^{(k)}$ (Gradient, Newton, and BFGS), which are targeted in the subsequent subsections.

Beforehand, we explain details of the accompanying line-search method, where we assume to be given a descendent direction $d^{(k)}$.
Concerning \eqref{eq:iteration_procedure}, we desire to find $\stepsizek$, such that, on the line determined by $d^{(k)}$, the updated iterate admits a decrease in $J(\mu^{(k+1)})$.
In general, another minimization problem is required to find the optimal step size.
In practice, we do not necessarily need the perfect step size but rather a sufficiently good one.
Indeed, in some cases, an appropriate step size $\stepsizek$ can be computed explicitly, and convergence results can be shown.
A more computationally feasible approach is the so-called Armijo- or (more strongly) Wolfe condition \cite{Nocedal} for approximately computing $\stepsizek$ while still fulfilling a suitable sufficient decrease.

The Armijo (or Armijo-Goldstein) condition for determining an appropriate step size $\stepsizek$ is given by
\begin{equation}\label{eq:armijo}
	J\left(\mu^{(k)} + \stepsizek d^{(k)}\right) - J\left(\mu^{(k)}\right) \leq \beta \stepsizek \nabla J\left(\mu^{(k)}\right)^T d^{(k)},		
\end{equation} 
with a number $\beta \in (0,1)$ that, according to~\cite{Nocedal}, is typically chosen as $\beta = 10^{-4}$.
The Armijo condition~\eqref{eq:armijo} prevents the line-search from overshooting in a specific direction, i.e., to keep~$\stepsizek$ small enough and enables a sufficient decrease in each iteration.
For more details on the decrease, we refer to~\cite{Nocedal,kelley}.
Moreover, we do not want to use only a tiny step size since this would require more outer algorithm iterations.
In practice, it is not feasible to check~\eqref{eq:armijo} for every number $\stepsizek \in \R^+$, since this again does not circumvent the fact that evaluations of $J$ might be costly.
Instead, we can decrease the step size iteratively, i.e. $\stepsizek \coloneqq \kappa^j$ for $j=0, 1, \dots$ and $\kappa \in (0, 1)$.

To terminate the outer optimization loop, we also require a termination criterion verifying that the algorithm iterated to a stationary point $\bar{\mu}$.
With regard to the standard first-order critical condition stated in~\eqref{simple_foc}, we use a tolerance $\tau_{\text{FOC}}$ and stop the algorithm if $$\norm{\nabla J(\mu^{(k)})}< \tau_{\text{FOC}}.$$
Next, we are concerned with computing the descent direction $d^{(k)} \in \R^P$. 

\subsubsection{Steepest-descent method}

The so-called steepest-descent (or gradient-descent) method is one of the most popular standard optimization methods.
We compute the gradient of $J$ at each iteration point $\mu^{(k)}$ and obtain the search direction by $d^{(k)} \coloneqq - \nabla J(\mu^{(k)})$.
By using the Armijo condition as explained above, the standard Armijo-based steepest-descent method can be summarized by Algorithm~\ref{alg:gradient_descent}.

\begin{algorithm2e}
	\KwData{$\mu^0$, $\kappa$, $\beta$, $\tau_{\text{FOC}}$}
	\KwResult{Approximate stationary point $\bar{\mu}$}
	$k \leftarrow 0$\;
	\While{$\norm{\nabla J(\mu^{(k)})} < \tau_{FOC}$}{
		compute descent direction $d^{(k)} = - \nabla J(\mu^{(k)})$\;
		search $j$ such that~\eqref{eq:armijo} is fulfilled with $\stepsizek=\kappa^j$\;
		$\mu^{(k+1)} \leftarrow \mu^{(k)} + \stepsizek d^{(k)}$\;
		$k \leftarrow k + 1$\;
	}
	\caption{Steepest-descent method}\label{alg:gradient_descent}
\end{algorithm2e}

In general, the steepest-descent method is known to be a first-order optimization method, where "first-order" refers to the convergence order of the method.
It is relatively easy to find optimization problems where the gradient-descent method requires a lot of iterations and especially struggles to find an approximate stationary point if the region near the optimum is particularly flat.
A well-known example where the gradient-descent method has a weak convergence behavior is the Rosenbrock function, introduced in~\cite{rosenbrock}.
Because of the poor convergence of the steepest descent method, we instead consider higher-order optimization methods.

\subsubsection{Newton's method}
\label{sec:newton_method}
The most popular higher-order optimization method is Newton's method, where we assume to have explicit access to the Hessian of the objective functional $\HH$, or at least to the application of the Hessian to an arbitrary direction.
Then, we obtain the descent direction by solving 
\begin{equation}\label{eq:newtons_descent}
\HH(\mu^{(k)}) d^{(k)} = - \nabla J(\mu^{(k)}).
\end{equation}
Multiple techniques can be used to solve~\eqref{eq:newtons_descent}.
An intuitive way would be to compute the inverse of $\HH$ explicitly.
Admittedly, this may result in vast computational requirements, especially for large control spaces.
Instead, we can solve the linear system by direct linear solvers such as a Cholesky factorization or iterative linear solvers such as GMRES or the conjugate gradient (CG) method.
However, in the first place, if the Hessian is symmetric positive definite (spd), these methods either do not converge or may result in a non-descending direction $d^{(k)}$.
For references on the discussed solvers and further details, we again point to \cite{kelley} and the references therein.
If close enough to the optimum, and if \eqref{eq:newtons_descent} is solved exactly, a step size of $\stepsizek=1$ and thus $j=0$ is enough, and therefore, a line-search method is not needed.
However, far away from the local minima, Newton's method can not be considered convergent, which is why we combine it with a line-search method.
Further, in the case where the Hessian is not spd, for instance, the truncated CG-Newton method can be used, where the CG iterations are followed as long as the algorithm fails due to the missing assumptions on the Hessian.
In this case, we only approximate the solution of \eqref{eq:newtons_descent}, referring to inexact-Newton methods.
Suppose the descent direction is indeed only an approximate solution of~\eqref{eq:newtons_descent}.
In that case, we call the method an inexact Newton method, and it is required to verify that $d^{(k)}$ is, indeed, a descendent direction.
This can be checked by $\nabla J(\mu^{(k)})^T d^{(k)} < 0$ and is trivially the case for the gradient-descent method.
For a more elaborated discussion on how to approximate the solution of~\eqref{eq:newtons_descent}, we refer to~\cite{kelley}.

Given the (approximate) solution of~\eqref{eq:newtons_descent}, we again apply a line-search and follow the Armijo rule to verify the sufficient decrease, which is referred to as the damped- or relaxed Newton method.
We shall also mention that Newton's method is commonly combined with a so-called Wolfe condition, see \cite[p.~38]{Nocedal}.
The above-explained damped Newton method can be summarized by \cref{alg:newtons_method}.

\begin{algorithm2e}
	\KwData{$\mu^0$, $\kappa$, $\beta$, $\tau_{\text{FOC}}$}
	\KwResult{Approximate stationary point $\bar{\mu}$}
	$k \leftarrow 0$\;
	\While{$\norm{\nabla J(\mu^{(k)})} < \tau_{FOC}$}{
		compute descent direction $d^{(k)}$ by solving~\eqref{eq:newtons_descent} with an 
		(inexact) linear solver\;
		search $j$ such that~\eqref{eq:armijo} is fulfilled with $\stepsizek=\kappa^j$ \label{newton:ls}\;
		$\mu^{(k+1)} \leftarrow \mu^{(k)} + \stepsizek d^{(k)}$\;
		$k \leftarrow k + 1$\;
	}
	\caption{Damped Newton method}\label{alg:newtons_method}
\end{algorithm2e}

While the (damped) Newton method defines the gold standard for higher optimization methods, there also exists a class of quasi-Newton methods which is another type of inexact-Newton methods.
A prevalent quasi-Newton method which was introduced in~\cite{BFGS} is the Broyden–Fletcher–Goldfarb–Shanno (BFGS) algorithm.
In the following, we shortly revise the basic algorithm since, for the constrained case, it is often used throughout this thesis.

\hspace{0.5cm}
\subsubsection{BFGS method}

Let $B^{(k)}$ be an approximation of the Hessian $\HH$ at the current iterate $\mu^{(k)}$.
Then, analogously to \eqref{eq:newtons_descent}, we determine the descent direction by solving the equation
\begin{equation}\label{eq:descent_with_approx_Hessian}
	B^{(k)} d^{(k)} = - \nabla J(\mu^{(k)}).
\end{equation}
After each iteration, a new update of the Hessian approximation $B^{(k)}$ can be performed by
\begin{align*}
	x &\coloneqq \mu^{(k+1)} - \mu^{(k)}, \\
	y &\coloneqq \nabla J(\mu^{(k+1)}) - \nabla J(\mu^{(k)}), \\
	B^{(k)} &\coloneqq B^{(k)} + \frac{yy^T}{y^Tx} - \frac{(B^{(k)}x)(B^{(k)}x)^T}{x^TB^{(k)}x}.
\end{align*}
We are again interested in an efficient way of solving~\eqref{eq:descent_with_approx_Hessian}.
In the particular case of the BFGS method, it can be shown (cf.~\cite[Lemma 4.1.1]{kelley}) that there exists an explicit iterative formula for $(B^{(k)})^{-1}$, which reads
\begin{equation}\label{eq:update_Bk}
	(B^{(k)})^{-1} = \left(I - \frac{xy^T}{y^Tx} \right) B^{(k-1)} \left(I - \frac{yx^T}{y^Tx} 
\right) + \frac{xx^T}{y^Tx},
\end{equation}
commonly known as the Sherman-Morrison-Woodbury formula.
Consequently, we can compute the descent direction as $d^{(k)} = - (B^{(k)})^{-1} \nabla J(\mu^{(k)})$ efficiently.
We use the identity as initial value~$B^{(0)}$.
In an iterative procedure, $B^{(k)}$ may have a deplorable approximation behavior.
The so-called curvature condition commonly tracks this by $x^Ty > 0$. Whenever the condition is violated in an update of $B^{(k+1)}$, the matrix is instead reverted to the identity.
Just as in the inexact versions of Newton's method, a potential descent direction $d^{(k)}$ is only accepted if $\nabla J(\mu^{(k)})^T d^{(k)} < 0$.
We can again follow an Armijo-rule analogously to the above-described procedure with the computed descent direction and report the BFGS algorithm in Algorithm~\ref{alg:BFGS}.
For more theoretical details on the BFGS algorithm, we recommend~\cite[Chapter 4]{kelley}.

\begin{algorithm2e}
	\KwData{$\mu^0$, $\kappa$, $\beta$, $\tau_{\text{FOC}}$}
	\KwResult{Approximate stationary point $\bar{\mu}$}
	$k \leftarrow 0$\;
	$B^0 \leftarrow I$\;
	\While{$\norm{\nabla J(\mu^{(k)})} < \tau_{FOC}$}{
		compute descent direction $d^{(k)} = - (B^{(k)})^{-1} \nabla J(\mu^{(k)})$\;
		search $j$ such that~\eqref{eq:armijo} is fulfilled with $\stepsizek=\kappa^j$\;
		$\mu^{(k+1)} \leftarrow \mu^{(k)} + \stepsizek d^{(k)}$\;
		update $(B^{(k)})^{-1}$ with~\eqref{eq:update_Bk}. If $x^Ty > 0$, set $(B^{(k)})^{-1} 
\leftarrow I$\;
		$k \leftarrow k + 1$\;
	}
	\caption{BFGS method}\label{alg:BFGS}
\end{algorithm2e}

We also emphasize that, in practice, it can become a storage problem that the BFGS method stores the Hessian approximation as an $n \times n$ matrix.
For these cases, there exist \emph{limited-memory} variants of the BFGS method (L-BFGS)~\cite{LBFGS}, where the idea is to solely store vectors that are used for an implicit representation of the Hessian.
\hspace{0.5cm}

\subsubsection{More inexact optimization methods}
The amount of variants of the above-mentioned basic optimization methods is vast, mainly depending on the specific application and problem.
Just as the mentioned simplifications of Newton's method, many optimization methods can be considered \emph{inexact optimization methods}.
For instance, gradient information might not be explicitly accessible in the model.
For these cases, it is possible to use inaccurate derivative information, for instance, automatic differentiation, finite differences, or stochastic approaches to approximate the gradient.
Importantly, it can be expected that the numerically approximated derivative information can suffer significant numerical errors due to instabilities, as, for instance, is known for finite differences.
On top of that, the numerical optimization method can be susceptible to approximation errors in the gradient.
Similar techniques can be used to approximate the Hessian.
We refer to~\cite[Section 2.3.1]{kelley} for an elaborated discussion on errors in functions, gradients, and Hessians.
As mentioned, the truncated variants for non-spd Hessians are also inexact methods.

\subsection{Constrained optimization}
\label{sec:constr_optimization}

The approaches in \cref{sec:unconstrained_optimization} are discussed for an optimization problem without constraints on the parameter space.
Typically, optimization tasks are very likely to be subject to constraints of several types.
In \cref{sec:background_PDE-constr}, we already discussed an example of equality constraints given by the PDE and we introduced the control space $\Params$ as a bounded subset of $\R^P$.
For instance, simple bound constraints $\Paramsad\coloneqq[\mu_\mathsf{a},\mu_\mathsf{b}]$, with vectors $\mu_\mathsf{a}, \mu_\mathsf{b} \in \R^P$ may be enforced on \eqref{P_simple}.
Consequently, local minima may lie outside the constraints and (global) stationary points may not be found at all.
To tackle this, necessary and sufficient optimality conditions need to be adjusted.
Condition~\eqref{simple_foc} can simply be modified to
\begin{equation}\label{simple_constr_foc}
\nabla J(\bar{\mu})^T \cdot (\mu - \bar{\mu}) \geq 0 \qquad\qquad \text{for all }\,\mu \in \Paramsad.
\end{equation}
This type of condition is well aligned with the first-order optimality condition for the PDE-constrained case in \cref{prop:first_order_opt_cond}.

To replace the termination criterion, we use a projection operator $P_{\Paramsad}: \R^P \to \Paramsad$ that maps an inadmissible control point to the respective (unique) point inside the constraints.
It can be shown (see~\cite[Theorem 5.2.4]{kelley})) that a $\bar{\mu}$ is a stationary point if and only if 
\begin{equation}\label{simple_proj_foc}
	\bar{\mu} = P_{\Paramsad}(\bar{\mu} - s \nabla J(\bar{\mu})) \qquad\qquad 
	\text{for all }\,s \geq 0,
\end{equation} 
or in other words
\begin{equation}\label{simple_proj_foc_norm}
\norm{\bar{\mu} - P_{\Paramsad}(\bar{\mu} - \nabla J(\bar{\mu}))} = 0.
\end{equation} 
Hence, we are left with finding an appropriate projection $P_{\Params}$ for the respective constraints.
For simple bound constraints this can be defined by
\begin{align*}
(P_{\Paramsad}(\mu))_i\coloneqq \left\{ \begin{array}{ll}
(\mu_\mathsf{a})_i & \text{if } (\mu)_i\leq (\mu_\mathsf{a})_i, \\
(\mu_\mathsf{b})_i & \text{if } (\mu)_i\geq (\mu_\mathsf{b})_i, \\
(\mu)_i & \text{otherwise,}
\end{array} \right. 
&& \text{for } i=1,\ldots,P.
\end{align*}

\noindent With the help of $P_{\Paramsad}$ and~\eqref{simple_proj_foc_norm}, we can thus detect a stationary point of the respective constrained optimization problem.
It is also clear that every iterate $\mu^{(k+1)}$ should fulfill the constraints, i.e., we change the iteration update to 
$$
\mu^{(k+1)} = P_{\Paramsad}(\mu^{(k)} + \stepsizek d^{(k)}).
$$
For the constrained case, the Armijo-type rule~\eqref{eq:armijo} has to be changed as well.
It can be replaced by the following Armijo-type sufficient decrease condition:
\begin{equation}
\label{eq:Armijo}J(\mu^{(k+1)}) - J(\mu^{(k)}) \leq  -\frac{\beta}{\stepsizek} \big\| \mu^{(k)} - 
\mu^{(k+1)}\big\|^2_2.
\end{equation}

From the above-discussed details, we can derive projected versions of the steepest-descent, BFGS, and Newton method, see \cite[Chapter 5]{kelley}.
It is important to note that these methods also incorporate active set strategies, which means that a component of $\mu^{(k)}$ that lies on the boundary of the constraints (and is thus called an active point) is treated differently in the optimization routine.
For more insights on active set strategies, see \cite[Chapter 5]{kelley}.

Notably, the convergence analysis of projected methods can not immediately be transferred from the unconstrained case and requires more involved considerations.
We do not detail the projected methods and their convergence results.
Whenever we use them in the sequel of this thesis, we instead point to the respective sections in the references.

\subsection{Trust-region methods}
\label{sec:background_TR_methods}

We present another class of iterative optimization algorithms called trust-region methods.
The main idea of these algorithms is to obtain a robust globally convergent algorithm while keeping computational effort low.

At each iteration $k\geq 0$, we consider a so-called model function $m^{(k)}$, which is a cheaply computable approximation of the cost functional $J$ in a neighborhood of the control $\mu^{(k)}$.
As the definition suggests, this neighborhood is called the "trust-region".
Therefore, for $k\geq 0$, given a TR radius $\delta^{(k)}$, we consider the TR minimization sub-problem
\begin{equation}
\label{TRsubprob}
\min_{s\in \mathbb{R}^P} m^{(k)}(s) \, \text{ subject to } \|s\|_k \leq \delta^{(k)},\, 
\widetilde{\mu}\coloneqq \mu^{(k)}+s \in\Paramsad 
\, \text{ for all }  v\in V.
\end{equation}
Under suitable assumptions on $m^{(k)}$, problem~\eqref{TRsubprob} admits a unique solution $\bar s^{(k)}$, which is used to compute the next iterate $\mu^{(k+1)} = \mu^{(k)} + \bar s^{(k)}$.
Trust-region methods are thus very flexible, and many variants can be derived. To discuss this flexibility, in Algorithm~\ref{alg:BTR}, we follow~\cite[Algorithm 6.1.1]{TRbook} to define what is called a basic trust-region (BTR) algorithm with a standard first-order critical point termination criterion for an unconstrained optimization problem.

\begin{algorithm2e}
	\KwData{initial point $\mu^{(0)}$, TR-radius $\delta^{(0)}$, and condition $0<\eta<1$}
	\KwResult{Approximate stationary point $\bar{\mu}$}
	$k \leftarrow 0$\;
	\While{$\norm{\nabla J(\mu^{(k)})} < \tau_{FOC}$}{
		\textbf{Step 1: Model definition.} Choose $\norm{\cdot}_k$ and define model 
$m^{(k)}$\;
		\textbf{Step 2: Step calculation.} Compute a step $s^{(k)}$ as a solution of 
\eqref{TRsubprob}\;
		\textbf{Step 3: Acceptance of trial point.} Compute $J(\mu^{(k)} + s^{(k)})$ and
		$$
		\rho^{(k)} = \frac{J(\mu^{(k)})-J(\mu^{(k)} + s^{(k)})}{m^{(k)}(0)-m^{(k)}(s^{(k)})}.
		$$
		If $\rho^{(k)} \geq \eta$, then accept $\mu^{(k+1)} = \mu^{(k)} + s^{(k)}$, otherwise reject 
step\;
		\textbf{Step 4: Trust-region radius update.} Shrink or enlarge TR radius if 
required or possible\;
		$k \leftarrow k + 1$\;
	}
	\caption{Basic trust-region algorithm for unconstrained problems}\label{alg:BTR}
\end{algorithm2e}

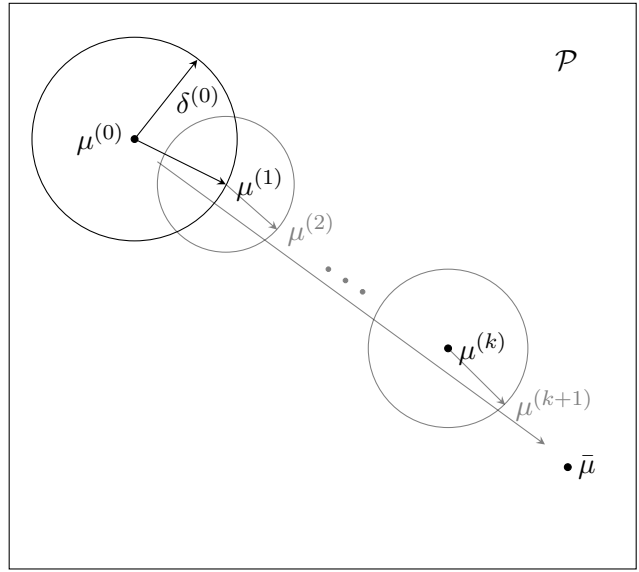
\begin{figure} \centering
	\begin{tikzpicture}[scale=1.5, >=stealth]
	\draw[fill=white] (-2.5,-2) rectangle (3,3);
	\draw[black](2.2,2.5) node[right] (scrip) {$\Params$};
	\draw[fill=black] (-1.4,1.8) circle [radius=0.03, color=black];
	\draw[black](-1.4,1.8) node[left] (scrip) {$\mu^{(0)}$};
	\draw[fill=black] (2.4,-1.1) circle [radius=0.03, color=black];
	\draw[black](2.4,-1.1) node[right] (scrip) {$\bar{\mu}$};
	\draw[->, opacity=0.5] (-1.2,1.6) -- (2.2,-0.9) node {}; 
	\draw[->,black] (-1.4,1.8) -- (-0.85,2.5) node [midway,right] {$\delta^{(0)}$};
	\draw[black] (-1.4,1.8) circle [radius=0.9, color=black];
	\draw[->,black] (-1.4,1.8) -- (-0.6,1.4) node [right] {$\mu^{(1)}$};
	\draw[black, opacity=0.5] (-0.6,1.4) circle [radius=0.6, color=black];
	\draw[->,black, opacity=0.5] (-0.6,1.4) -- (-0.15,1) node [right] {$\mu^{(2)}$};
	\draw[fill=black, opacity=0.5] (0.3,0.65) circle [radius=0.02, color=black];
	\draw[fill=black, opacity=0.5] (0.45,0.55) circle [radius=0.02, color=black];
	\draw[fill=black, opacity=0.5] (0.6,0.45) circle [radius=0.02, color=black];
	\draw[fill=black] (1.35,-0.05) circle [radius=0.03, color=black];
	\draw[black](1.35,-0.05) node[right] (scrip) {$\mu^{(k)}$};
	\draw[black, opacity=0.5] (1.35,-0.05) circle [radius=0.7, color=black];
	\draw[->,black, opacity=0.5] (1.35,-0.05) -- (1.85,-0.55) node [right] {$\mu^{(k+1)}$};
	\end{tikzpicture}	
	\caption[Visualization of TR methods.]{Visualization of TR methods. The method starts from an initial guess $\mu^{(0)} \in \Params$ and iteratively updates the model function while converging to $\bar{\mu}$.}
	\label{fig:TR_visualization}
\end{figure}

We refer to \cref{fig:TR_visualization} for a visualization of Algorithm~\ref{alg:BTR}.
This basic description does not cover all TR variants and does not detail assumptions that have to be made for the individual components of the algorithm.
A significant part of the algorithm is to choose the model in Step~1.
In general, only a few assumptions need to be fulfilled for the model function $m^{(k)}$, concerning regularity and approximation properties of $J$.
We refer to~\cite[Section 6.2.2]{TRbook} for an elaborated discussion on suitable assumptions.
As a commonly used example, we introduce the following example of a model function:
$$
m^{(k)}(s) \coloneqq J(\mu^{(k)}) + \nabla J(\mu^{(k)})^T s + \frac{1}{2} s^T \HH(\mu^{(k)}) s. 
$$
Note that, for computing $m^{(k)}(s)$, we do not need to evaluate (the possibly expensive) functional~$J$, but instead use a quadratic approximation around the current iterate $\mu^{(k)}$.
One of these variants is the TR-Newton-CG Steihaug method, presented in~\cite[Algorithm 7.2]{Nocedal}, where CG approximations are used for the Hessian.
The variant will be used in the sequel of this thesis to mimic a truncated Newton line-search algorithm; cf. \cref{sec:exp_paper2}.

The solution of the sub-problem~\eqref{TRsubprob} requires iterative optimization algorithms like the ones discussed in \cref{sec:iterative_optimization} with the significant difference that the model function can cheaply be evaluated.
Importantly, these sub-problems are usually solved concerning a sufficient decrease condition.
If such a sufficient point can not be found, the so-called Cauchy point (that can always be found with the gradient-descent method) can be used instead.
For more details, we refer to \cite[Section 6.3]{TRbook}.

If a point is accepted in Step~3 of \cref{alg:BTR}, the model function shows a reasonably good accuracy w.r.t.~$J$.
Then, we may have a chance to increase the TR radius further (in Step~4).
If the point is not accepted, the model is assumed to have poor approximation quality in this region, and we need to decrease the TR-radius in Step~4 and recompute the sub-problem.
We also emphasize that the choice of the first-order critical termination criterion has to be adjusted to the respective application, e.g., with a projected version if constraints on the parameter set are enforced.
For an elaborated convergence study of unconstrained and constrained TR-Algorithms, we refer to~\cite{TRbook}, for instance.
Concerning the TR method that we develop in \cref{chap:TR_RB}, we present the convergence study in detail; cf.~\cref{sec:TRRB_convergence_analysis}.

Trust-region methods have been used in many fields (see the discussion in~\cite[Section 1.3]{TRbook}, for instance) and the number of variants is remarkable.
In particular, their robustness for complex optimization tasks has shown advantageous in many regards.
For the work at hand, the idea of a TR method is crucial since it enables optimization methods with certified but adaptive surrogate modeling.
In \cref{chap:TR_RB}, we change Step~1 to an error-aware version, where we use a model function that is a surrogate model for the high-fidelity function $J$.
Let us assume that the error of this surrogate model can be efficiently bounded by an a posteriori error estimator $\Delta_{\Jhat}$.
Then, we can replace the norm $\norm{\cdot}_k$ by this estimator, which results in a more meaningful trust-region since it is no longer a relatively unrelated metric object (such as a circle in \cref{fig:TR_visualization}) but is instead associated to the actual error characteristic that the surrogate model produces.
Moreover, since the surrogate model is built progressively based on all iterates $\mu^{(k)}$, we note that, different from \cref{fig:TR_visualization}, the TR contains all iterates (and their respective local regions).
We precisely describe the algorithm in \cref{chap:TR_RB} and present elaborated illustrations of the resulting error-aware trust-regions; cf.~\cref{sec:proof_of_concept}.

For solving PDE-constrained parameter optimization problems such as \eqref{Phat}, it remains to evaluate~$\Jhat$ for an arbitrary point in $\Params$, and thus, to solve the parameterized PDE that is formulated in \cref{def:param_elliptic_problem}.
For this purpose, the following section provides an overview of the numerical approximation of PDEs.

\section{Numerical approximation of PDEs}\label{sec:background_npdgl}

This section is devoted to the numerical solution method for parameterized elliptic problems, defined in \cref{def:param_elliptic_problem}.
Recall that the involved equation is called \emph{the primal equation} in the application of PDE-constrained parameter optimization in \cref{sec:background_PDE-constr}.
Every time we want to evaluate the objective functional for a new parameter $\mu \in \Params$, we require a solution method for deducing the state $u_\mu \in V$ as a solution of \cref{def:param_elliptic_problem}.
As a first step toward a numerical approximation of $u_\mu$, we require a finite-dimensional subspace of $V$.
In this thesis, we solely use grid-based approximation methods that take advantage of a spatial discretization of the underlying PDE to obtain a finite-dimensional subspace $V_h \subset V$.
The approximations $u_{h,\mu} \in V_h$ can then be found by the solution of a system of (linear) equations.
In what follows, we introduce the basic idea of the finite element method (FEM) and 
discuss challenges that arise with specific problem classes such as multiscale methods.
 
\subsection{The finite element method}
\label{sec:fem}

The finite element method has been tremendously developed and used by both engineers and academics to numerically approximate PDEs based on the so-called Ritz-Galerkin approach.
For standard literature, we refer to \cite{Braess, ciarlet2002finite}.
The main idea is to discretize the computational domain $\Omega$ in an appropriate way, meaning to define a grid (also called mesh) $\mathcal{T}_h$ which consists of non-overlapping elements $t \in \mathcal{T}_h$ with a simple shape.
This family of elements is chosen such that the whole domain $\Omega$ can be replicated by it, i.e. $\overline{\bigcup_{t \in \mathcal{T}_h}t} = \Omega$.
The corresponding mesh is often called finite element mesh (FE mesh) as it consists of (finitely many) elements, nodes, and edges.
For instance, these elements are intervals in one dimension, triangles or quadrilaterals in two dimensions, and tetrahedrons in three dimensions.
FE meshes provide the basis of the finite element space (FE space) $V_h \subset V$, approximating the infinite-dimensional function space $V$.
A typical choice of the FE space is based on $\mathcal{T}_h$- piecewise linear functions that are continuous on~$\Omega$:
$$
\mathcal{P}_1(\mathcal{T}_h) \coloneqq \set[v \in C^0(\overline{\Omega})]{v \big|_t \text{ is a linear polynomial of 
degree $\leq 1$, for every $t \in \mathcal{T}_h$}}.
$$
Then, for any FE mesh $\mathcal{T}_h$, we define the corresponding FE space by
$$
V_h = V \cap \mathcal{P}_1(\mathcal{T}_h).
$$
We note that there exists a large variety of FE spaces, e.g., for higher polynomial order; see also \cite{Braess,ciarlet2002finite} for more finite elements.
With the help of the discretized FE space, we formulate the standard FEM for parameterized problems; cf. \cref{def:param_elliptic_problem}.

\begin{problem}[Discrete parametrized elliptic problem]
	\label{def:discr_param_elliptic_problem}
	Let the assumptions of \cref{def:param_elliptic_problem} be fulfilled.
	For a fixed parameter $\mu \in \Params$, we seek the discrete version of $u_\mu \in V$, denoted by $u_{h,\mu} \in V_h$, such that
	\begin{equation}
	\label{eq:discr_param_elliptic_problem}
	a_\mu(u_{h,\mu},v_h) =  l_\mu(v_h) \qquad\qquad \text{for all } \, v_h \in V_h.
	\end{equation} 
	We also define the discrete solution map by $\mathcal{S}_h:\Params \to V_h$, $\mu \mapsto u_{\mu, h}\coloneqq \mathcal S_h(\mu)$.
	Furthermore, the discrete output functional is defined as $\Jhat_h(\mu) \coloneqq \J(u_{h,\mu}, 
\mu)$.
\end{problem}

The choice of the FE space $V_h$, as well as the variational formulation in \cref{def:discr_param_elliptic_problem}, is a so-called continuous Galerkin (CG) approximation.
It is important to note that, for the forms~$a_\mu$,~$l_\mu$ and~$\J$, we do not introduce an approximate version $a_{h,\mu}$, $l_{h,\mu}$ or $\J_h$, although inserting functions of $V_h$ changes the way we evaluate these.
In our case, however, the approximate versions only mean replacing integral expressions with numerical quadrature rules (which are exact in the piece-wise linear case).
This does not necessarily hold in a more general case of $V_h$.
For instance, higher-order choices of $V_h$ impede exact integral approximations.
Moreover, in the case of discontinuous Galerkin (DG) methods (which is a non-conforming method, i.e., $V_h \not\subset V$), the approximate bilinear form $a_{h,\mu}$ additionally consists of penalty terms of the non-conforming discretization; see \cite{cockburn2012discontinuous}.
Since this thesis solely deals with $\mathcal{P}_1$-conforming approaches, we omit these cases of approximate forms and stick to the notation in \cref{def:discr_param_elliptic_problem}.
In our case, coercivity and continuity assumptions on the forms $a_\mu$ and $l_\mu$ also translate to $V_h$.
Moreover, we know that $V_h \subset V$ is also a Hilbert space (as a closed subspace of $V$), inheriting the scalar product of $V$.
Hence, we can still apply Lax-Milgram's theorem to obtain a unique solution for every parameter.

Recall that \cref{def:discr_param_elliptic_problem} is a Ritz-Galerkin approach with the particular case of a FE space $V_h$ (which then makes it a finite element method).
In other words, Equation~\eqref{eq:discr_param_elliptic_problem} is a Ritz-Galerkin projection of the infinite-dimensional formulation of Equation~\eqref{eq:param_elliptic_problem}, meaning that the \emph{test space} for test functions $v_h$ is given by $V_h$ and the \emph{ansatz space} for the solutions $u_{h,\mu} \in V_h$ is the same.
Choices where test- and ansatz spaces are not equal are called Petrov--Galerkin approaches and can have advantages in terms of storage consumptions and approximability; cf. \cref{sec:TRRB_pg_approach} and \cref{sec:PG--LOD}.
Petrov--Galerkin formulations are also often used for stabilizing purposes, e.g., for transport-dominated problems \cite{brooks1982streamline}.
We also refer to \cite{babuvska1994special, griffiths1978analysis} and the references therein.
Apart from these works, we note that Petrov--Galerkin formulations have become very helpful in various applications, two of which will further be discussed in \cref{sec:TRRB_pg_approach} and \cref{sec:LOD}.

In practice, \cref{def:discr_param_elliptic_problem} results in a system of linear equations which can be solved with a linear solver.
As we shall see in the next section, the number of unknowns scales with the number of nodes in the FE grid.
This is of significant importance for the computational effort of FEM.

\subsection{Algebraic formulation}
\label{sec:algebraic_fem}

In this section, we focus on the computational procedure for obtaining an approximate solution of Problem~ \ref{def:discr_param_elliptic_problem}.
The FE space $V_h$ is finite-dimensional by construction, and therefore, it also has a finite-dimensional basis.
For the simple construction of $\mathcal{P}_1$, these basis functions are called shape- or hat-functions, mainly because of their geometric interpretation for $d=1,2$, cf. \cref{correctorplots}(b) for a 2d-plot.
Let $(\phi_i)_{i=1}^{N_h}$, $N_h \coloneqq \dim{V_h}$, denote the $\mathcal{P}_1$ finite basis of $V_h$. Then, this basis represents a partition of unity, i.e., only one basis function is $1$ at a single node, whereas the others vanish.
Thus, all basis functions have meager support, which proves advantageous in terms of storage and the choice of the linear solver since the small sparsity pattern enables a quick assembly of the system matrix.

Given the basis $(\phi_i)_{i=1}^{N_h}$, we can write a solution $u_{h, \mu} \in V_h$ of~\eqref{eq:discr_param_elliptic_problem} as
$$
	u_{h, \mu} \coloneqq \sum_{i=1}^{N_h} \underline{u_{h, \mu}}^{\mkern-15mu i \mkern+15mu} \phi_i,
$$
where $\underline{u_{h, \mu}} \in \R^{N_h}$ denotes the coefficient vector of $u_{h, \mu}$.
Following this representation, we can use the linearity of $l_\mu$ and linearity in the second argument of $a_\mu$ to write~\eqref{eq:discr_param_elliptic_problem} as to finding  $\underline{u_{h, \mu}} \in \R^{N_h}$, such that
$$
	a_\mu( \sum_{i=1}^{N_h} \underline{u_{h, \mu}}^{\mkern-15mu i \mkern+15mu} \phi_i, \phi_j ) 
= l_\mu(\phi_j)
	\qquad\qquad \text{for all }\,j=1,\dots,N_h.
$$
Using linearity of $a_\mu$ in the first argument results in
$$
\sum_{i=1}^{N_h} a_\mu( \phi_i, \phi_j )\underline{u_{h, \mu}}^{\mkern-15mu i \mkern+15mu} = 
l_\mu(\phi_j)
\qquad\qquad \text{for all }\,j=1,\dots,N_h.
$$
Hence, by defining the matrix $(\mathbb{A}_\mu)_{ji} \coloneqq a_\mu(\phi_i, \phi_j)$ and the vector $(\mathbb{F}_\mu)_i \coloneqq l_\mu(\phi_i)$, a solution $u_{h, \mu}$ solves~\eqref{eq:discr_param_elliptic_problem} if and only if its component vector $\underline{u_{h, \mu}} \in \R^{N_h}$ solves the linear system of equations:
\begin{equation}\label{eq:fem_LSE}
	\mathbb{A}_\mu \underline{u_{h, \mu}} = \mathbb{F}_\mu.
\end{equation}
The matrix $\mathbb{A}_\mu \in \R^{N_h \times N_h}$ is called the stiffness matrix and the vector $\mathbb{F}_\mu \in \R^{N_h}$ is called the load vector.
In fact, to obtain a solution of~\eqref{eq:discr_param_elliptic_problem}, we need to solve a system of linear equations of size $N_h$.
As known from linear algebra, using a direct solver can scale up to $\mathcal{O}(N_h^3)$.
Luckily, the stiffness matrix has a minimal sparsity pattern due to the small support of the basis functions in $V_h$, which simplifies the solution procedure.
Furthermore, the sparse stiffness matrix is relatively easy to store.
In conclusion, the computational efficiency is restricted to linear solvers of~\eqref{eq:fem_LSE} and scales with the number of basis functions in $V_h$.

Lastly, also the output functional $\Jhat_h$ of \cref{def:discr_param_elliptic_problem} can be written in an algebraic form.
In the easiest case where the output functional is linear, we have that $(\mathbb{J}_\mu)_i = J(\phi_i, \mu)$ is a vector. Thus, $\Jhat_h$ can be computed by 
\begin{equation}
	\Jhat_h = \underline{u_{h, \mu}}^T \mathbb{J}_\mu.
\end{equation}
Note that nonlinear output functionals, such as the linear-quadratic one $\J_{L^2}$ from \eqref{eq:example_J} that is used in \cref{chap:TR_RB}, can also be written within an algebraic form.

\subsection{Error bounds}
\label{sec:error_bounds_for_fem}

We now investigate the approximation error of FEM and discuss the choice of the mesh size of $\mathcal{T}_h$.
The associated mesh size $h$ of a FE mesh $\mathcal{T}_h$ is defined as the maximum diameter of all elements in $\mathcal{T}_h$.
If $h$ is fixed, the C\'ea Lemma provides an abstract a-priori result on the accuracy of the corresponding Ritz-Galerkin method, which we formulate in the sequel.

\begin{lemma}[C\'ea's Lemma, see {\cite[Lemma 4.2]{Braess}}]\label{lem:ceas_lemma}
	For a fixed parameter $\mu \in \Params$, the exact solution $u_\mu \in V$ of \cref{def:param_elliptic_problem}, and the discrete solution $u_{h, \mu} \in V_h$ of \cref{def:discr_param_elliptic_problem}, we have
	\begin{equation}
		\norm{u_\mu - u_{h, \mu}} \leq \frac{\cont{a_\mu}}{\alpha_{a_\mu}} \inf_{v_h \in 
V_h} \norm{u_\mu - v_h},
	\end{equation}
	where $\cont{a_\mu}$ and $\alpha_{a_\mu}$ denote the continuity and coercivity constants of $a_\mu$.
	In addition, Galerkin orthogonality is fulfilled, i.e.
	\begin{equation}
		a_\mu(u_\mu - u_{h, \mu}, v_h) = 0 \qquad\qquad \text{for all }\, v_h \in V_h.
	\end{equation}
\end{lemma}

\noindent From C\'ea's Lemma, we conclude that the FEM approximation yield a quasi best-approximation $u_{h,\mu}$ in the respective space $V_h$.
\cref{lem:ceas_lemma} also provides the groundwork for a-priori results of the FEM approximation in the $L^2$- and $H^1$- norm as defined in \eqref{eq:norms}.
The following result can be proven with an interpolation estimate of an interpolation operator that maps functions from $V$ to $V_h$ as well as with the Lemma of Aubin-Nietsche.

\begin{theorem}[$L^2$- and $H^1$- a-priori error bound for FEM, see {\cite[Chapter II.7]{Braess}}]\label{thm:apriori_fem}
	For a fixed parameter $\mu \in \Params$, let $u_{h, \mu} \in V_h$ be the solution of \cref{def:discr_param_elliptic_problem}.
	If the exact solution $u_\mu \in V$ of \cref{def:param_elliptic_problem} is $H^2(\Omega)$-regular, we have
	\begin{equation}
		\norm{u_\mu - u_{h, \mu}}_{L^2(\Omega)} \leq C_{L^2} h^2\qquad\qquad \text{and} 
\qquad\qquad \norm{u_\mu - u_{h, \mu}}_{H^1(\Omega)} \leq C_{H^1} h,
	\end{equation}
	with constants $C_{L^2}, C_{H^1} > 0$, independent of $h$.
\end{theorem}

For more information on the regularity of $u_\mu$, we also refer to \cite[Chapter II.7]{Braess}.
From \cref{thm:apriori_fem}, we see that, depending on the norm, the approximation accuracy scales linearly or quadratically with the mesh size $h$.
We note that this result can be extended to higher-order polynomials in the FE space and obtains higher-order convergence results if further regularity properties are given.
\cref{thm:apriori_fem} can be understood as a justification for using FEM since the desired accuracy can always be reached for a sufficiently small $h$.
As highlighted in \cref{sec:algebraic_fem}, the computational effort of FEM heavily scales with the size of $h$, and thus, we desire a choice of $h$ where we do not waste unnecessary resources.
Moreover, \cref{thm:apriori_fem} does not give rise to a rule for how to choose $h$ for given specific data functions.

Choosing an appropriate FE mesh has been extensively studied, and many related approaches have been developed.
One possible way is to use a multi-grid approach, i.e., refining the mesh at specific parts of the computational domain, for instance, w.r.t. the regularity of the solution.
We particularly mention the family of methods that are based on a posteriori error analysis~\cite{verfurth1994posteriori,Ver2013}.
Here, the idea is to obtain an estimator that can be computed after solving the approximate solution $u_{h, \mu} \in V_h$ and estimates the error w.r.t. the unknown real solution~$u_\mu$.
Hence, it can posteriorly be verified whether the mesh size is sufficiently small.
There exist multiple methods that result from this theory.
For instance, the adaptive finite element method (AFEM)  uses a localized a posteriori error estimator for indicating regions to refine the mesh.
For details, we refer to \cite{MNS2002} and the references therein.
For other related approaches, see the discussion in \cite[Chapter II.8ff]{Braess}.

This thesis mainly neglects theory and ideas on finding an appropriate mesh size $h$ for the discretization.
We always assume to priorly know a sufficient mesh size for our problem, which is further formalized in \cref{asmpt:truth} and \cref{asmpt:loc_truth}.
In the subsequent section, we instead motivate the case where we already know that the given data functions require a mesh size $h$ that results in a prohibitively large system.

\subsection{Large scale and multiscale problems}\label{sec:multiscale_problems}

The purpose of this section is to briefly motivate the need for alternative solution methods for \cref{def:discr_param_elliptic_problem} that can handle the case of a prohibitively small mesh size $h$.
These phenomena are often called large- or multiscale problems.
They have seen tremendous development in the last decades, as they are particularly motivated by real-world simulations where the associated data structure becomes more and more complex.
An extensive literature overview of large- and multiscale problems and some real-world applications have been given in \cref{chap:introduction}.
For an elaborated motivation of multiscale problems, we additionally refer to \cite{LODbook}.

From a numerical point of view, it is obvious that a complex data structure requires an appropriately small mesh size that accounts for the resolution of the associated data functions.
In a case where the data functions are already defined w.r.t. a prescribed data mesh, for instance, given by a pixel-based material description, it is easy to find such a mesh size.
Moreover, it is not feasible to expect a sufficiently accurate approximation if the mesh size does not entirely capture the data.

In a more analytical setting where the data functions are instead given by mathematical expressions or objects, an appropriate mesh size is not directly apparent.
These kinds of problems are closely related to homogenization theory, where rapid oscillations are present.
In what follows, we present a numerical example that aims at presenting the issue of multiscale problems and is depicted from \cite{Pe15}.

\begin{figure}
	\centering
	\includegraphics[scale=0.8]{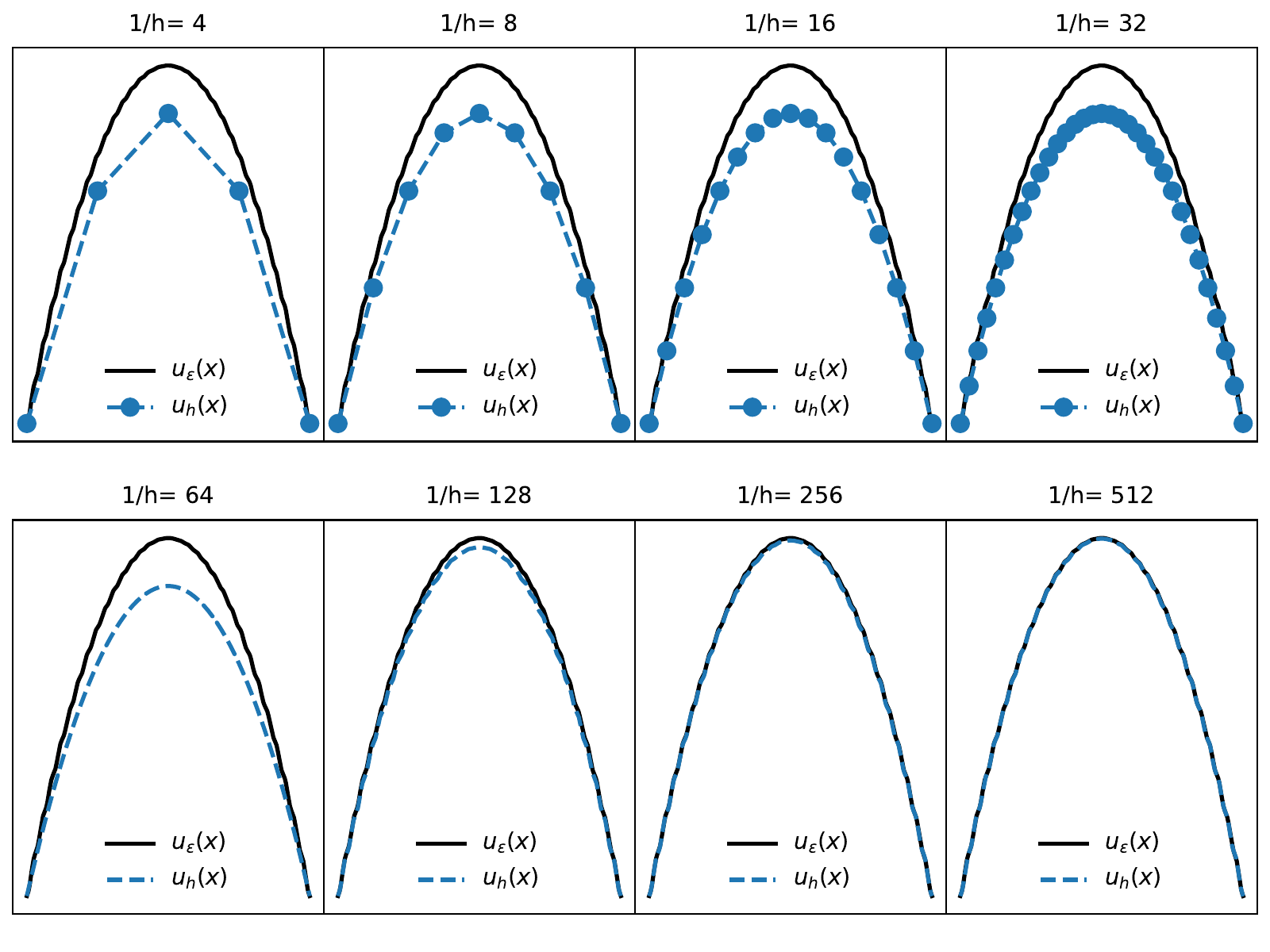}	
	\caption[FEM approximations of the 1-dimensional multiscale problem.]{FEM approximations of the 1-dimensional multiscale problem with varying mesh size $h$ and $\varepsilon=1/128$.}
	\label{fig:1d-ms_example}
\end{figure}

Let $A_{\varepsilon}$ be a rapidly oscillating diffusion coefficient with frequency $\varepsilon >0$, such that
\[
A_{\varepsilon}(x) := \frac{1}{4}\Big( 2 - \cos(2 \pi x / \varepsilon) \Big)^{-1}.
\]
With $f \equiv 1$, the exact solution of \eqref{eq:classical_PDE} can be stated explicitly by
\begin{equation} \label{eq:assmpt_exp_example}
u_{\varepsilon}(x) = 4 (x-x^2)- 4 \varepsilon \left( \frac{1}{4 \pi} \sin(2 \pi \frac{x}{\varepsilon}) - \frac{1}{2 \pi}x \sin(2 \pi \frac{x}{\varepsilon}) - 
\frac{\varepsilon}{4 \pi^2} \cos(2 \pi \frac{x}{\varepsilon}) + \frac{\varepsilon}{4 \pi^2} \right),
\end{equation}
where the oscillations of $A_{\varepsilon}$ are only present in $u_\varepsilon$ by a factor $\varepsilon$.
On the other hand, $u_{\varepsilon}$ admits a coarse-scale behavior that is independent of $\varepsilon$, cf.~\cref{fig:1d-ms_example}.
From a numerical perspective, one could expect that FEM can already capture the coarse behavior with a relatively coarse mesh size $h > \varepsilon$ and that the approximation error that occurs due to the non-captured oscillations is of size $\varepsilon$.
However, as we see in \cref{fig:1d-ms_example}, where we approximate the problem with $\varepsilon=10^{-5}$, this expectation is wrong.
On the contrary, the approximation error for $h > \varepsilon$ is large and the coarse behavior is not captured at all.

A theoretical explanation of this behavior can be given by \cref{thm:apriori_fem}.
A closer look at the proof shows that the error bound is also dependent on the term $\norm{\nabla^2 u_\varepsilon}_{L^2(\Omega)}$; cf.~\cite[Theorem 7.3]{Braess}.
Hence, the approximability of $u_\varepsilon$ also depends on its second-order derivative.
Although the coarse behavior of $u_\varepsilon$ is mainly unaffected by the rapid changes, the term $\norm{\nabla^2 u_\varepsilon}_{L^2(\Omega)}$ is also approximately the size of the frequency $\varepsilon^{-1}$.
Thus, following \cref{thm:apriori_fem}, we can only expect a sufficient numerical approximation if $h \ll \varepsilon$. 
We conclude that the numerical approximation of such problems results in an arbitrarily large system for $\varepsilon \to 0$.

General grid-free homogenization theory is devoted to finding a suitable convergence concept for classifying a suitable homogenized limit of $\varepsilon \to 0$.
In the case of PDEs, for instance, the two-scale convergence, which goes back to \cite{allaire} and \cite{nguetseng1989general} can be used.
The convergence is internally based on a so-called asymptotic expansion:
\begin{equation}\label{eq:ass_exp}
u_\varepsilon(x) = u^0(x) + \varepsilon u^1(x, x/\varepsilon) + \dots,
\end{equation}
with small $\varepsilon > 0$, and where $u^1(x, x/\varepsilon)$ is $1$-periodic in the second argument.
The term $u^0$ is called the coarse part of $u^\varepsilon$, and the remaining scales are considered micro-scales.
In particular, the solution in \eqref{eq:assmpt_exp_example} can be written as \eqref{eq:ass_exp}.
However, in general, the fulfillment of an asymptotic expansion \eqref{eq:ass_exp} is a pretty strong assumption and is only given for a specific class of multiscale problems.

Many fields that study multiscale problems are inspired by homogenization theory and often build on several similar assumptions (for instance, periodicity in the data).
Let us also mention that multiscale problems are by far not only of interest for grid-based numerical approximation methods but are also of great interest in other grid-free fields, for instance in general energy-minimizing approaches in the calculus of variation.
Famous concepts for convergence theory in homogenization are the so-called \emph{H-convergence} \cite{murat2018h} or \emph{$\Gamma$-convergence} \cite{dal2012introduction}, see also \cite{murat1997calculus}.

Homogenization theory is primarily devoted to understanding the case where $\varepsilon \to 0$, whereas, instead, numerical multiscale methods are designed to approximate the solution $u_\varepsilon$ for a small but fixed~$\varepsilon$.
In a nutshell, numerical multiscale methods follow the idea of coarse-scale approximations, where the coarse-scale system incorporates reconstructions of fine scales that can be computed locally.
At the same time, as indicated before, numerical multiscale methods can also be used for large-scale problems.
In \cref{sec:LOD}, we introduce the well-established localized orthogonal decomposition (LOD) method as a specific instance of numerical multiscale methods that can handle arbitrarily rough coefficients.

In this thesis, parameterized large- and multiscale problems mainly serve as a particular example of problems where standard methods such as FEM are infeasible, given the limited computational resources.
To tackle this, we discuss localized approaches and extensions to parameterized problems by using localized model order reduction.
This is the main contribution of \cref{chap:TSRBLOD}, where we present how the LOD can be extended to solve parameterized multiscale problems efficiently.

For generally accelerating the solution process of parameterized problems, we now introduce model order reduction methods, one of the main theoretical foundations of this thesis.

\section{Model order reduction for parameterized PDEs}\label{sec:background_model_order_reduction}

In the former section, we provided an overview of how parameterized PDEs like \cref{def:param_elliptic_problem} can be approximated numerically.
With the concepts explained, the overall goal of numerically solving PDE-constrained parameter optimization problems like \eqref{Phat}, based on a discretization of the involved equations, is now finally possible.
This is further formalized as \eqref{Phat_h} in \cref{chap:TR_RB}.
Regardless of the optimization routine and under suitable assumptions, model order reduction (MOR) techniques can be used to speed up the solution process for multiple parameter samples of \cref{def:discr_param_elliptic_problem}.

MOR conducted for the numerical approximation of PDEs is widespread and has been used in uncountable applications \cite{MR3701994}.
In what follows, we introduce the basic idea of MOR methods for the numerical approximation of PDEs.
In particular, we introduce the idea of the reduced basis method (RBM).
Subsequently, in \cref{sec:background_RB_methods_for_PDEopt}, we explain the coherence to PDE-constrained optimization problems, introduce existing approaches from the literature, and motivate the main results of this thesis.
For a tutorial introduction to RBM, we refer to \cite{HAA2017}.


To recall from \cref{sec:param_problems}, we are interested in finding the solution $u_\mu \in V$ of a parameterized elliptic \cref{def:param_elliptic_problem}, for a fixed parameter instance $\mu \in \Params$ of a parameter space $\Params \subset \R^P$.
In the discrete setting of \cref{sec:background_npdgl}, this translates to finding the discrete solution $u_{h, \mu} \in V_h$ of \cref{def:discr_param_elliptic_problem}.
We can employ the same FEM techniques for every new parameter $\mu \in \Params$.
Although the solution method that we use for solving \cref{def:discr_param_elliptic_problem} can be very efficient in itself, standard techniques such as FEM do not use information that has been obtained from former samples of $\Params$.
In contrast, the idea of MOR methods is to use specific former solutions of \cref{def:discr_param_elliptic_problem} for constructing a reduced basis space and a corresponding reduced model.

MOR methods usually use the phrases \emph{full-order model} (FOM) and \emph{reduced-order model} (ROM).
In the standard case, we typically consider \cref{def:discr_param_elliptic_problem} as the FOM since it has full complexity and achieves maximum accuracy (for fixed and reasonably small $h$).
Moreover, we consider the FOM as the expensive model and the ROM as a low dimensional and desirably efficient model yet to be constructed.
More generally, the ROM is often called a \emph{surrogate model}.
The construction of such a surrogate model requires evaluations of the FOM and can thus be very expensive.
In MOR methods, the phase for constructing a ROM is called the \emph{offline phase}, whereas the evaluation of the ROM is called the \emph{online phase}.
Having this in mind, as also illustrated in the motivation of this thesis, we are interested in two scenarios, where MOR methods show their great advantage:

\begin{figure}[t] \centering
	\begin{tikzpicture}
	\definecolor{color0}{rgb}{0.65,0,0.15}
	\definecolor{color1}{rgb}{0.84,0.19,0.15}
	\definecolor{color2}{rgb}{0.96,0.43,0.26}
	\definecolor{color3}{rgb}{0.99,0.68,0.38}
	\definecolor{color4}{rgb}{1,0.88,0.56}
	\definecolor{color5}{rgb}{0.67,0.85,0.91}
	\definecolor{color6}{rgb}{0.27,0.46,0.71}
	\begin{axis}[
	name=left,
	anchor=west,
	width=8cm,
	height=5cm,
	tick align=outside,
	tick pos=left,
	legend style={nodes={scale=0.7}, fill opacity=0.8, draw opacity=1, text opacity=1, 
		anchor=north west, xshift=-4.7cm, yshift=0cm, draw=white!80!black
	},
	x grid style={white!69.0196078431373!black},
	xlabel={number of samples $\mu$},
	ylabel={run time [s]},
	xmajorgrids,
	xtick style={color=black},
	y grid style={white!69.0196078431373!black},
	ymajorgrids,
	ytick style={color=black},
	ytick pos=left,
	]
	\addplot [thick, color6, mark size=0]
	table {%
		0 0
		10 15
		20 30
		30 45
		40 60
		50 75
	};
	\addlegendentry{FOM}
	\addplot [thick, color2, mark size=0]
	table {%
		0 0
		0.01 35
		50 35.7 
	};
	\addlegendentry{ROM}
	\end{axis}
	\end{tikzpicture}
	\caption[Example of overall wall time comparison in a many-query scenario.]{Example of overall wall time comparison in a many-query scenario. The offline time to construct the fast ROM requires $35$ seconds. After a certain point of parameter samples, this offline time eventually pays off.}
	\label{fig:many-query}
\end{figure}
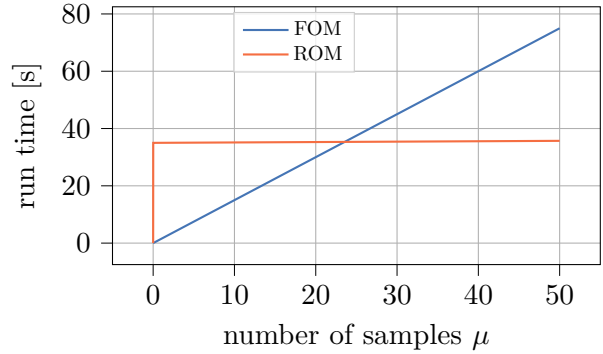

\hspace{0.2cm}
\begin{description}
	\item[Many-query scenarios:] In many applications, parametric problems need to be solved for multiple, if not thousands, of parameter samples $\mu \in \Params$.
	In this context, one is interested in reducing the \emph{overall effort} (or overall-efficiency) of the numerical simulation.
	To construct a suitable ROM in advance may take a lot of offline time but eventually pays off after a certain number of solutions are required.
	For a simple illustration of such an efficiency comparison, we refer to \cref{fig:many-query}.
	A particular instance of many-query scenarios is PDE-constrained parameter optimization problems as presented in \cref{sec:background_PDE-constr}.
	Other examples are Monte Carlo methods.
	\item[Real-time scenarios:] Instead of aiming at overall-efficiency in a many-query simulation, MOR methods can also be extremely effective even when the offline time is comparably large.
	In real-time scenarios, it is desired to approximate the solution $u_{h,\mu} \in V_h$ for an unknown parameter $\mu \in \Params$ in a minimal amount of time and with minimal storage requirements.
	In this context, (potentially large) offline times can be assumed to be negligible.
	Examples of real-time scenarios are applications where a decision needs to be taken quickly (think about the example of a Formula 1 race) or applications where only a ROM can be installed on a small device with low computational power (for example, a smartphone or raspberry-pi).
\end{description}
\hspace{0.2cm}

We emphasize that the efficient utilization of MOR methods for both scenarios highly depends on the problem class.
In this thesis, we are only concerned with the elliptic case stated in \cref{def:discr_param_elliptic_problem}, for which many MOR methods have been demonstrated and developed.
As a particularly well-suited instance of these methods, we now introduce the RBM for the efficient numerical approximation of parameterized PDEs.

\subsection{The reduced basis method}
\label{sec:background_RB_method}
The main idea of the RB method is to construct a low dimensional space from full-order solutions of \cref{def:discr_param_elliptic_problem}.
Hence, let $V_\red \subset V_h$ be a low dimensional subspace of the FE space $V_h$, with dimension $N_\red \coloneqq \dim(V_\red) \ll N_h$ and a corresponding basis $\Psi_\red \coloneqq \set{\psi_1, \dots, \psi_{N_\red}}$.
Then, we call~$V_\red$ a reduced basis space with a reduced basis $\Psi_\red$ and formulate the ROM of \cref{def:discr_param_elliptic_problem} as follows:

\begin{problem}[Reduced parameterized elliptic problem]
	\label{def:red_param_problem}
	Let the assumptions of \cref{def:param_elliptic_problem} be fulfilled.
	For a fixed parameter $\mu \in \Params$, we seek the reduced solution $u_{\red, \mu} \in V_\red$, such that
	\begin{equation}
	\label{eq:red_param_problem}
	a_\mu(u_{\red,\mu},v_\red) =  l_\mu(v_\red) \qquad\qquad \text{for all } \, v_\red \in V_\red.
	\end{equation} 
	We also define the reduced solution map by $\mathcal{S}_\red:\Params \to V_\red$, $\mu \mapsto u_{\red, \mu}\coloneqq \mathcal S_\red(\mu)$ and the reduced output functional as $\Jhat_\red(\mu) \coloneqq \J(u_{\red,\mu}, \mu)$.
\end{problem}

Since $V_\red$ is a Hilbert space (with the inherited product), this problem also admits a unique solution by Lax Milgram's theorem.
In a view of the discussion in \cref{sec:fem}, Equation~\eqref{eq:red_param_problem} can be considered as a Galerkin projection of \cref{def:discr_param_elliptic_problem} on the reduced space $V_\red$.
A simple modification of the a priori error bound in \cref{lem:ceas_lemma} suggests that we can consider the reduced solution $u_{\red,\mu}$ as the quasi best-approximation result in $V_\red$.

As motivated in the former subsection, we are concerned with the computational effort that is required for solving~\eqref{eq:red_param_problem}.
Moreover, we need to construct the reduced basis $\Psi_\red$.
Before answering these questions in detail, we first introduce a posteriori error estimation for certifying the approximation quality of the ROM.

\subsection{A posteriori error control}
\label{sec:background_a_post_errors}

While \cref{lem:ceas_lemma} gives rise that the formulation of \cref{def:red_param_problem} is feasible, it does not yet tell us whether the approximation quality of $V_\red$ is sufficient.
The overall approximation error can be divided into two parts, the discretization error that is due to the choice of $V_h$ and the reduction error in $V_\red$, i.e. using the triangle inequality, we have
\begin{equation}
	\norm{u_{\mu} - u_{\red, \mu}} \leq \underbrace{\norm{u_{\mu} - u_{h, \mu}}}_{\text{discr. error}} + \underbrace{\norm{u_{h,\mu} - u_{\red, \mu}}}_{\text{red. error}}.
\end{equation}
There exist multiple approaches concerning the error control of the RB method that either bound the overall approximation error or solely concentrate on the reduction error.
As stated in \cref{sec:background_npdgl}, for the remainder of this thesis, we consider the discretization error to be negligible, meaning that an appropriate mesh size $h$ is known, cf. \cref{asmpt:truth}.

In \cref{sec:error_bounds_for_fem}, we discussed residual-based a posteriori error estimation for finding a sufficiently small choice for the mesh size $h$ for FEM.
To resume, the idea is to value the accuracy of an already computed approximation $u_{\red,\mu} \in V_\red$ with neither accessing the FEM solution $u_{h,\mu} \in V_h$ nor the true solution $u_{\mu} \in V$.
The situation in RB methods is even more comfortable since reduced approximations are considered cheap, and a good approximation quality is very important since we want to evaluate the ROM for many parameter samples.
Moreover, we know that we can quickly increase the approximation quality of $V_\red$ by adding another FOM basis function to $\Psi_\red$.
The standard a posteriori error estimator for the model reduction error is based on the residual of~\eqref{eq:discr_param_elliptic_problem}, evaluated with the reduced solution $u_{\red, \mu}$. 
Analogously to \cref{def:primal_residual}, for every $u_{\red, \mu} \in V_\red$, we define the corresponding residual $r_\mu(u_{\red, \mu}) \in V'_h$ by
\begin{align}
r_\mu(u_{\red, \mu})[v_h] \coloneqq l_\mu(v_h) - a_\mu(u_{\red, \mu}, v_h) \qquad\qquad \text{for all } \, v_h \in 
V_h.
\label{eq:RB_residual}
\end{align}	
Now, we state the standard a posteriori error result for RB methods. We refer to \cite{HAA2017}, for instance.

\begin{proposition}[A posteriori estimate for the model order reduction error]
	\label{prop:standard_RB_estimator}
	For $\mu \in \Params$ let $u_{h, \mu} \in V_h$ be the solution of \cref{def:discr_param_elliptic_problem} and $u_{\red, \mu} \in V_\red$ the solution of \cref{def:red_param_problem}.
	Then, it holds
\begin{equation} \label{eq:a_post_RB}
	\norm{u_{h, \mu} - u_{\red, \mu}} \leq \Delta(\mu) \coloneqq \alpha_{a_\mu}^{-1} \, 
\|r_\mu(u_{\red, \mu})\|_{V'_h}.
\end{equation}
\end{proposition}
\begin{proof}
	With the shorthand $e_{h, \mu} \coloneqq u_{h, \mu} - u_{\red, \mu}$, we have
	\begin{align*}
	\alpha_{a_\mu}\, \|e_{h, \mu}\|^2 &\leq a_\mu(e_{h, \mu}, e_{h, \mu})
	= \underbrace{a_\mu(u_{h, \mu}, e_{h, \mu})}_{= l_\mu(e_{h, \mu})} - a_\mu(u_{\red, \mu}, 
	e_{h, \mu})\\
	&= r_\mu(u_{\red, \mu})[e_{h, \mu}]
	\leq \|r_\mu(u_{\red, \mu})\|_{V'_h}\, \|e_{h, \mu}\|
	\end{align*}
	using the coercivity of $a_\mu$ in the first inequality, the definition of $e_{h, \mu}$ in the first equality, the fact that $u_{h, \mu}$ solves~\eqref{eq:discr_param_elliptic_problem} and the definition of the residual~\eqref{eq:RB_residual} in the second equality, and the continuity of the residual in the second inequality. Dividing by $\|e_{h, \mu}\|$ gives the desired result.
\end{proof}

It is clear that the resulting estimate does not necessarily need to be sharp and other choices of estimators are possible (for example, scaled with different constants).
From a MOR point of view, we are interested in a sharp estimator that is computationally affordable.
To characterize the estimator in its approximation quality, we can use the notion of effectivity.
In a nutshell, an estimator can be called effective if the overestimation factor can be bounded by a small constant.
It can be shown that $\Delta$ is an effective estimator, whereas, for other instances in this thesis, such a result does not exist.
For more information about the effectiveness and reliability of estimators, we again refer to \cite{HAA2017}.

Reduced models often inherit multiple reduced quantities of interest, such as the objective functional, derivative information, or other auxiliary equations.
A posteriori error analysis can also be used to control the reduction error for these quantities.
In the formulation of \cref{def:param_elliptic_problem} we introduced the objective functional $\Jhat$ and, in \cref{def:red_param_problem}, we defined a reduced version $\Jhat_\red$ of the discrete functional $\Jhat_h$.
For RB-based PDE-constrained parameter optimization methods, we are also interested in an estimator $\Delta_{\Jhat_\red}$, such that
$$|\Jhat_h(\mu) - \Jhat_\red(\mu)| \leq \Delta_{\Jhat}(\mu),$$
which will further be motivated in \cref{sec:background_RB_methods_for_PDEopt} and \cref{chap:TR_RB}.
Suitable error estimation for linear output functionals has extensively been targeted in \cite{HAA2017}, where also a primal-dual approach for a more accurate reduced functional is discussed.
We refer to \cref{chap:TR_RB}, in particular to \cref{sec:TRRB_estimation_for_J}, for an elaborated discussion on nonlinear and corrected reduced output functionals, and their error estimation.

\subsection{Offline-online decomposition}
\label{sec:background_off_on_RB}

So far, we have stated the reduced model and introduced a way to certificate it with an a posteriori error estimator.
Certainly, RB methods only show their remarkable power when the solution method for \cref{def:red_param_problem}, as well as the evaluation of the error estimator $\Delta$, is independent of the dimension of $V_h$, and hence, independent of the fine mesh size $h$.
This is the case when all computations requiring the fine resolution can be done in the offline phase, and only cheap computations are left for the online phase.

To achieve such an \emph{offline-online decomposition}, the so-called parameter separability (or affine decomposition) from \cref{asmpt:parameter_separable} is key.
To repeat, all involved forms (like $a_\mu$, $l_\mu$, and $\J$) can be written in the following form:
$$
a_\mu(u, v) = \sum_{\xi = 1}^{\Xi^a} \theta_{\xi}^a(\mu)\, a_{\xi}(u, v).
$$
We now give a detailed explanation of how the offline-online decomposition can be built from \cref{asmpt:parameter_separable}.

\subsubsection{Offline-online decomposition of the ROM}
Let us start with an algebraic view on \cref{def:red_param_problem}.
Similar to \cref{sec:algebraic_fem}, we use the reduced basis $\Psi_\red = \set{\psi_1, \dots, \psi_{N_\red}}$ to write the solution $u_{\red, \mu} \in V_\red$ of \cref{def:red_param_problem} as
$$
u_{\red, \mu} \coloneqq \sum_{i=1}^{N_\red} \underline{u_{\red, \mu}}^{\mkern-15mu i \mkern+15mu} \psi_i,
$$
where $\underline{u_{\red, \mu}} \in \R^{N_\red}$ again denotes the component vector of $u_{\red, \mu}$.
By defining the matrix $(\mathbb{A}_{\red,\mu})_{ji} \coloneqq a_\mu(\psi_i, \psi_j)$ and the vector $(\mathbb{L}_{\red,\mu})_i \coloneqq l_\mu(\psi_i)$, a function $u_{\red, \mu}$ solves~\eqref{eq:red_param_problem} if and only if its component vector $\underline{u_{\red, \mu}} \in \R^{N_h}$ solves
\begin{equation}\label{eq:RB_LSE}
	\mathbb{A}_{\red,\mu} \, \underline{u_{\red, \mu}} = \mathbb{L}_{\red,\mu}.
\end{equation}
In contrast to the sparse stiffness matrix $\mathbb{A}_\mu$ for solving the high-fidelity equation \eqref{eq:fem_LSE}, the reduced stiffness matrix $\mathbb{A}_{\red,\mu}$ is dense, which does certainly not harm the speed of a solution method for~\eqref{eq:RB_LSE} as long as $N_\red$ is small.
Hence, the solution method for the ROM is already independent of $h$.
However, the assembly of $\mathbb{A}_{\red,\mu}$ and $ \mathbb{L}_{\red,\mu}$ is not.
For every new parameter $\mu$, the basis functions of $\Psi_\red$ still need to be inserted into $a_\mu$ and $l_\mu$.
At this point, the affine decomposition helps.
To be precise, in the offline phase, we can pre-compute $a_\xi(\phi_i, \phi_j)$ and $l_\xi(\phi_i)$ for every affine component $\xi \in \set{1, \dots, \Xi}$ and basis function $\psi_i \in \Psi_\red$ to obtain matrices and vectors $\mathbb{A}_\xi$ and $\mathbb{L}_\xi$.
Consequently, for a new parameter $\mu$, the matrix $\mathbb{A}_{\red,\mu}$ and the vector~$\mathbb{L}_{\red,\mu}$ can be quickly assembled by
\begin{equation*}
	\mathbb{A}_{\red,\mu} = \sum_{\xi = 1}^{\Xi_a} 
	\theta_{\xi}^a(\mu)\,\mathbb{A}_\xi
	\qquad\qquad\text{and}\qquad\qquad
	\mathbb{L}_{\red,\mu} = \sum_{\xi = 1}^{\Xi_l} \theta_{\xi}^l(\mu)\,\mathbb{L}_\xi.
\end{equation*}
Moreover, the same trick can be used for evaluating the output functional $\Jhat_\red$ with an offline-online decomposition.
For this, we require to pre-compute $\mathbb{J}_\xi$ for all $\xi$ and all basis functions.
In the linear case, we then have
\begin{equation*}
	\mathbb{J} = \sum_{\xi = 1}^{\Xi_\J} \theta_{\xi}^{\J}(\mu) \mathbb{J}_\xi
	\qquad\qquad\text{and}\qquad\qquad
	\Jhat_\red(\mu) =  \underline{u_{\red, \mu}}^T \mathbb{J}, 
\end{equation*}
which again can be generalized to the nonlinear case of $\J$.
In conclusion, we showed that all parts of the ROM can be offline-online decomposed given \cref{asmpt:parameter_separable}.
In terms of numerical effort, we observe that the assembly of \eqref{eq:RB_LSE} is of order $\mathcal{O}(\Qa N_\red^2 + \Xi_l N_\red)$ and solving the system is of order $\mathcal{O}(N_\red^3)$. Therefore, the computational effort is independent of $N_h$.

\subsubsection{Offline-online decomposition of the coercivity constant}
\label{sec:min_theta}

For the efficient computation of the error bound in \cref{prop:standard_RB_estimator}, we require a quickly computable function $\alpha_{a,\text{LB}}$ serving as a lower bound for $\alpha_{a_\mu}$ for any $\mu \in \Params$, such that
$$0 < \alpha_{a,\text{LB}}(\mu) \leq \alpha_{a_\mu}.$$

\noindent If \cref{asmpt:parameter_separable} is fulfilled, such a function can be obtained by the \emph{min-theta approach} (see \cite[Prop.~2.35]{HAA2017}), where we have
\begin{equation}
	\alpha_{a,\text{LB}}(\mu) \coloneqq  \alpha_{a_{\check{\mu}}} \cdot \min_{\xi = 1,\dots,\Xi_a} \frac{\theta_\xi^a(\mu)}{\theta_\xi^a(\check{\mu})}.
\end{equation}
Here, $\alpha_{a_{\check{\mu}}}$ is the coercivity constant of a fixed parameter $\check{\mu}$.
The proof that $\alpha_{a,\text{LB}}(\mu)$, indeed, defines a lower bound of $\alpha_{a_\mu}$ is a simple computation.

What remains is to compute $\alpha_{a_{\check{\mu}}}$, which, in particular, is dependent on the inner product of $V_h$.
If the inner product of $V_h$ is defined by the mesh-independent energy product $(u, v) \coloneqq a_{\check{\mu}}(u, v)$ for a fixed parameter $\check{\mu} \in \Params$, the coercivity constant $\alpha_{a_{\check{\mu}}}$ is one.
Since the bilinear form $a_{\mu}$ is symmetric, continuous, and coercive, this energy-product is a product on $V$.
This strategy is used in \cref{chap:TR_RB}.

Lastly, we mention that, in case \cref{asmpt:parameter_separable} is not fulfilled, the more general successive constraint method~\cite{pat07} yields less sharp estimates and is computationally demanding, both offline and online.

\subsubsection{Offline-online decomposition of the error bound}
\label{sec:background_off_on_error}
Next, we elaborate on the offline-online decomposition of the error estimator $\Delta$ in~\eqref{eq:a_post_RB}.
Evaluating $\Delta$ at $\mu \in \Params$ involves the computation of the coercivity constant $\alpha_{a_\mu}$ and the dual norm of the residual on $V_h$, which are both computationally demanding.
The efficient computation of the coercivity constant has already been targeted in the previous section.
In what follows, we explain the offline-online decomposition for the residual, using \cref{asmpt:parameter_separable}.

The dual norm of the residual can be computed by the Riesz-representative, i.e.
$$
\|r_\mu(u_{\red, \mu})\|_{V'_h} = \|v_{r, \mu}\|_{V_h},
$$
where $v_{r, \mu} \in V_h$ is the unique solution of $r_\mu(u_{\red, \mu})[v] =  (v_{r, \mu}, v)_{V_h}$, for all $v \in V_h$.
For computing the Riesz-representative, we require the product matrix of $V_h$ and the residual $r_\mu(u_{\red, \mu})$, which then leaves us with a linear system that still scales with $N_h$.
Luckily, the computation of the Riesz-representative can also be decomposed w.r.t. the affine decomposition in \cref{asmpt:parameter_separable}.
To see this, let $\Xi_r \coloneqq \Xi_l + N_\red \Xi_a$ and the residual components $r_\xi \in V'_h$, $\xi= 1, \dots, \Xi_r$, be defined by
\begin{align*}
		(r_1, \dots, r_{\Xi_r}) \coloneqq (l_1, \dots, l_{\Xi_l}&, a_1(\psi_1, \cdot), \dots, a_{\Xi_a}(\psi_1, 
		\cdot), \\
		\dots&, a_1(\psi_{N_\red}, \cdot), \dots, a_{\Xi_a}(\psi_{N_\red}, \cdot)).
\end{align*}
Then, for the component vector $\underline{u_{\red, \mu}} \in \R^{N_\red}$ of $u_{\red, \mu}$, we can define the parameter functionals
\begin{align*}
	(\theta^r_1, \dots, \theta^r_{\Xi_r}) \coloneqq (\theta^f_1, \dots, \theta^f_{\Xi_l}&, 
	- \theta^a_1 \cdot  \underline{u_{\red, \mu}}^{\mkern-15mu 1 \mkern+15mu} , \dots, 
	- \theta^a_{\Xi_a} \cdot  \underline{u_{\red, \mu}}^{\mkern-15mu 1 \mkern+15mu}, \\
	\dots&, - \theta^a_1 \cdot  \underline{u_{\red, \mu}}^{\mkern-15mu N_\red \mkern+15mu},
	\dots, - \theta^a_{\Xi_a} \cdot  \underline{u_{\red, \mu}}^{\mkern-15mu N_\red 
\mkern+15mu}).
\end{align*}
With the respective Riesz-representatives, $r_\xi[v] =  (v_{\xi}, v)_{V_h}$, we can thus write
\begin{equation}
	r_\mu(u_{\red, \mu})[v] = \sum_{\xi = 1}^{\Xi_r} \theta^r_\xi(\mu)r_\xi[v]
	\qquad\qquad\text{and}\qquad\qquad
	v_{r,\mu} = \sum_{\xi = 1}^{\Xi_r} \theta^r_\xi(\mu) v_\xi,
\end{equation}
for all $\mu \in \Params$ and $v \in V_h$. In total, we have
\begin{equation} \label{eq:classic_riesz_offon}
	\|v_{r, \mu}\|^2_{V_h} =  \sum_{\xi_1 = 1}^{\Xi_r} \sum_{\xi_2 = 1}^{\Xi_r} 
	\theta^r_{\xi_1}(\mu) \theta^r_{\xi_2}(\mu) \, (v_{\xi_1}, v_{\xi_2})_{V_h}.
\end{equation}
With this representation, we can pre-compute all $(v_{\xi_1}, v_{\xi_2})_{V_h}$ terms in the offline phase and consequently yield the offline-online decomposition.

As shown in~\cite{buhr14}, computing~\eqref{eq:classic_riesz_offon} is not numerically stable and can only reach half of the machine precision.
We shortly review the approach that has been proposed in~\cite{buhr14}
since it plays a vital role in \cref{chap:TSRBLOD}.
To remedy the instability of~\eqref{eq:classic_riesz_offon}, we use an orthonormal basis of $(v_\xi)_{\xi=1}^{\Xi_r}$, which we call $\Psi_\red^r = \set{\psi_1^r, \dots, \psi_{\Xi_r}^{r}}$.
Therefore, for every $v_\xi$, we have the representation $v_\xi = \sum_{\xi = 1}^{\Xi_r} \bar{v_\xi} \psi_\xi^r$, where $\bar{v_\xi} = (v_\xi, \psi_\xi^r)_{V_h}$ and we observe the alternative computation of~\eqref{eq:classic_riesz_offon} by
\begin{equation} \label{eq:stable_riesz_offon}
\|v_{r, \mu}\|^2_{V_h} =  \sum_{\xi_1 = 1}^{\Xi_r} \left( \sum_{\xi_2 = 1}^{\Xi_r} 
\theta^r_{\xi_1}(\mu) \bar{v_\xi} \right)^2,
\end{equation}
which is again offline-online decomposable.
While this approach is numerically stable, it requires the computation of the respective orthonormal basis (e.g., with Gram-Schmidt), which computationally affects the offline phase.

In fact, the computable offline-online decomposed version of $\Delta$ from \cref{prop:standard_RB_estimator} is
\begin{equation}\label{eq:computable_delta}
	\Delta(\mu) \leq \alpha_{a, \text{LB}}(\mu)^{-1} \, 
	\|v_{r,\mu}\|_{V_h},	
\end{equation}
where $\|v_{r,\mu}\|_{V_h}$ can be computed by~\eqref{eq:classic_riesz_offon} or~\eqref{eq:stable_riesz_offon}.
The numerical effort for the online evaluation of the error bound is of order $\mathcal{O}((\Xi_l +\Qa N_\red)^2 + \Qa)$, where the second term corresponds to the min-theta approach.

Finally, we conclude that, given \cref{asmpt:parameter_separable}, it is possible to derive an offline-online decomposed version of \cref{def:red_param_problem} and its certification via \cref{prop:standard_RB_estimator}.

\subsection{Basis generation}
\label{sec:background_basis_gen_RB}

In the previous sections, we have assumed the existence of a suitable reduced basis space $V_\red$.
It remains to explain how we efficiently construct such a space, i.e., how we find appropriate reduced basis functions.
The approximation quality of these basis functions is crucial for a good reduced model.
Furthermore, we aim at a small number of basic functions, such that~\eqref{eq:red_param_problem} can be solved quickly.
In the literature, there have been plenty of proposals on the construction of RB spaces; cf.~\cite{HAA2017}.
The most commonly used approaches are the proper orthogonal decomposition (POD) method \cite{KV2001,GV17} and the so-called greedy-search algorithm~\cite{binev2011convergence}.

In a POD, one assumes to be given the FOM solutions $(u_{h, \mu})_{\mu \in \Params_{\text{train}}}$ for a training set $\Params_{\text{train}} \subset \Params$.
The main idea of the POD is to extract the (in the $L^2$-sense) most valuable information from this set, given a prescribed tolerance or basis size.
The POD approach has proven advantageous in terms of the approximation quality of the RB space since the accuracy in the $L^2$-sense can be estimated by the respective next eigenvalue whose modes have not entered the RB space.
The POD is particularly handy for the basis construction of non-stationary problems since it can also be understood as a principal component analysis of full time trajectories.
In fact, the POD is especially helpful for data-based approaches, which, however, can also be understood as a computational disadvantage since it requires to solve \cref{def:discr_param_elliptic_problem} for the whole training set.
Since this thesis is mainly devoted to stationary problems with expensive forward problems and, at the same time, large parameter sets, we instead focus on a greedy-search algorithm to construct $V_\red$.

In contrast to POD, the idea of a greedy-search algorithm~\cite{binev2011convergence} is to iteratively find the current parameter $\mu \in \Params_{\text{train}}$, which, given the current ROM, admits the largest approximation error w.r.t. the FOM (in the $L^\infty$-sense).
In the strong version of the greedy-search algorithm case, the worst parameter is found by computing the actual model reduction or orthogonal projection error, which, just as in the POD, also means to solve \cref{def:discr_param_elliptic_problem} for the entire training set.
Once a prescribed tolerance $\varepsilon > 0$ is reached, the algorithm is stopped.
In conclusion, while POD constructs a quasi-optimal basis in the $L^2$-sense, the progressive procedure of the greedy-search constructs a quasi-optimal basis in the $L^\infty$-sense.
To reduce the computational cost required to compute the actual error, the so-called weak greedy-search algorithm uses the reliable a posteriori error estimate $\Delta$ for finding the worst parameter.
Thus, only solutions $u_{h, \mu}$ (so-called \emph{snapshots}) for parameters that, at some point in the iteration, were identified as the worst approximated parameters need to be computed.
Due to this, the training set $\Params_{\text{train}}$ can be chosen much larger than for the POD or the strong greedy-search algorithm.
We resume the weak greedy-search algorithm in \cref{alg:weak_greedy}.

\begin{algorithm2e}
	\KwData{Training set $\Params_{\text{train}}$, tolerance $\varepsilon>0$}
	\KwResult{$V_\red$}
	$V_\red \leftarrow \{0\}$\;
	\While{$\max_{\mu\in\Params_{\textnormal{train}}} \Delta(\mu) > \varepsilon$}{
		$\mu^* \leftarrow \argmax_{\mu \in \Params_{\text{train}}} \Delta(\mu)$\;
		$V_\red \leftarrow \Span(V_\red \cup \{u_{h, \mu^*}\})$\;
	}
	\caption{Weak greedy algorithm for the generation of $V_\red$.}\label{alg:weak_greedy}
\end{algorithm2e}

Such a greedy-search algorithm was first considered a heuristic, proving very advantageous in practice.
Theoretical results regarding the convergence behavior of the constructed space were first given in \cite{buffa2012priori} and have further been developed in \cite{binev2011convergence}.
The approximation result can be formalized with the so-called \emph{Kolmogorov n-width}, which, transferred to our notation, can be described by the maximum distance of an $n$-dimensional subspace of $V$ to a compact and closed subset $\mathcal{M}$ of $V$, i.e.
\begin{equation}\label{eq:kolmogorov}
d_n(\mathcal{M}) = \inf_{\substack{V_\red \subset V\\ \dim(V_\red)=n}} d(V_\red, \mathcal{M}),
\end{equation}
where the distance of $V_\red$ to $\mathcal{M}$ is defined as
\begin{equation}\label{eq:distance}
	d(V_\red, \mathcal{M}) = \sup_{u \in \mathcal{M}} \inf_{u_\red \in V_\red} \norm{u - u_\red}.
\end{equation}
For theoretical foundations, we refer to \cite{devore1993constructive,Pin1985}.
Furthermore, we note that an example of $\mathcal{M}$ is the training set $\ParamsTrain$.
In general, the infimum in \eqref{eq:kolmogorov} is considered impossible to find.
However, in \cite{buffa2012priori,binev2011convergence}, it was shown that, under mild assumptions, the greedy-search algorithm produces a reduced space where the Kolmogorov n-width can be bounded by a factor with an exponential (or algebraic) convergence rate.
Of course, fast convergence of the Kolmogorov n-width can only be obtained when the training set is chosen large enough and is heavily dependent on the general complexity of the problem.

By now, several variants for the greed-search algorithm have been considered; cf~\cite{HDO2011}.
We particularly mention so-called goal-oriented greedy-search algorithms, where other error estimators than $\Delta(\mu)$, tailored towards the respective quantity of interest, are used to find the goal-oriented worst parameter; cf.~\cref{sec:background_RB_methods_for_PDEopt}.
Furthermore, adaptive construction of the training set for significantly large parameter spaces was advised in~\cite{HDO2011}.

\subsection{Extensions and challenges}
\label{sec:RB_challenges}

RB methods have tremendous success for various problem classes and, at the same time, revealed lots of open questions left for future research.
For a recent overview, we refer to~\cite{successandlimits} and \cite{HAA2017}.

This thesis is mainly concerned with problems that can be considered as part of the "ideal world"\footnote{Phrase used in \cite{successandlimits}} of coercive and affinely decomposed problems.
We recommend \cite{MR3672144,MR3701994,Hesthaven2016,Quarteroni2016} and the references therein for literature concerning extensions of RB methods, such as for time-dependent problems, inf-sup stable problems, nonlinear problems, and other related approaches.

In \cref{sec:background_off_on_RB}, the online efficiency for real-time scenarios was observed based on the affine decomposition from \cref{asmpt:parameter_separable}.
As aforementioned, this restrictive assumption has successfully been addressed by the empirical interpolation \cite{BarraultMadayEtAl2004}, which has later been called the discrete empirical interpolation method (DEIM) in~\cite{CS2010} and was utilized in, e.g., \cite{CEGG2014,DHO2012,FHRS2016}.

Several concerns that arise in RB methods are still present.
The most prominent example is the large class of nonlinear problems and those where the Kolmogorov n-width is decaying significantly slowly.
While tackling these problems is ongoing work, we again mention that, in this thesis, we only consider coercive and affinely decomposed problems.
Two main challenges that are already present in such an "ideal world" are described in the following:

\begin{description}
\item[Curse of dimensionality] In the context of RB methods, this famous concern refers to the fact that the dimension of~$\Params$ or the parameter dependency (for example, the number of affine coefficients $\Xi$) of the involved functions can be large.
Then, again, the decay of the Kolmogorov N-width may be substantially slow.
Consequently, a global efficient surrogate model is difficult (and time-consuming) to find, and the overall (and online) efficiency suffers significantly.
Concerning \cref{fig:many-query}, this means that the construction of the reduced model may not pay off at all.
Hence, given the model reduction's specific application, there is a critical point where globally (w.r.t.~$\Params$) accurate surrogate models can be too expensive to build since many basis functions are required for obtaining a sufficiently rich reduced space.

\item[Inaccessible global discretization]
As another essential challenge of RB methods, by construction, the FOM of standard RB methods incorporates high-fidelity solutions obtained by a traditional solution method such as FEM on a fine mesh $\mathcal{T}_h$.
As we discussed in \cref{sec:multiscale_problems}, future research is also devoted to the case where computations on the full fine mesh are considered prohibitively costly.
For RB methods, this means that no FEM snapshots are available and localized model order reduction methods, instead, need to be followed.
\end{description}

\noindent In \cref{sec:background_outlook_rb}, we explain how these challenges are tackled in this thesis.
Beforehand, we shall discuss the additional characteristics that occur when RB methods accelerate the solution method of a PDE-constrained parameter optimization problem.

\subsection{RB methods for accelerating PDE-constrained optimization}
\label{sec:background_RB_methods_for_PDEopt}

As indicated in \cref{chap:introduction}, reduced basis methods have been enormously used for solving PDE-constrained optimization problems.
The initial idea for using RB methods for PDE-constrained optimization problems was to construct the surrogate such that it can be reused for several different optimization problems, either in a many-query optimization context or for real-time applications, where, for instance, the objective functional changes w.r.t. a parameter.
Thus, these strategies were not designed for only one instance of an optimization problem, such that the (potentially large) offline time to construct the reduced model can be considered negligible \cite{Dede2012,dihl15,GK2011,KTV13,NRMQ2013,OP2007}.

As discussed in \cref{sec:first_opt_first_dis}, these methods also differ in terms of their discretization strategy.
It can generally be said that the \emph{first-optimize-then-discretize} approach is very efficient for a single optimization problem.
On the other hand, the \emph{first-discretize-then-optimize} approach shows advantages for many problems of the same type because, for instance, objective functionals can be replaced more easily.

Regardless of the solution approach, as discussed in \cref{sec:background_basis_gen_RB}, a suitable greedy-search algorithm for PDE-constrained parameter optimization can follow a goal-oriented strategy.
In this case, the greedy-search terminates if the a posteriori error estimator for, e.g., the reduced functional $\Delta_{\Jhat}$ or its derivative, Hessian, or parameter information is small enough.
Suitable choices for tolerances for constructing an appropriately accurate reduced model have been discussed in the given references.
It has undoubtedly been shown that in cases where the offline phase is ruled out, RB methods have a large impact on the computational time in the online phase, especially if the utilized optimization methods require many iteration steps.

However, it was unclear how RB methods can help to solve only a single optimization problem.
Then, overall-efficiency is required, and thus, the offline time can not be ignored.
In particular, concerning \cref{fig:many-query}, it can be assumed that the many-query application of a single optimization process is likely to become a many-query instance where the overall required run time of the fully FOM-based method is not as large as the offline time for constructing the ROM.
Thus, we face a many-query scenario with only a small number of samples.
Clearly, this highly depends on the optimization problem and the required iterations of the (higher-order) optimization method.

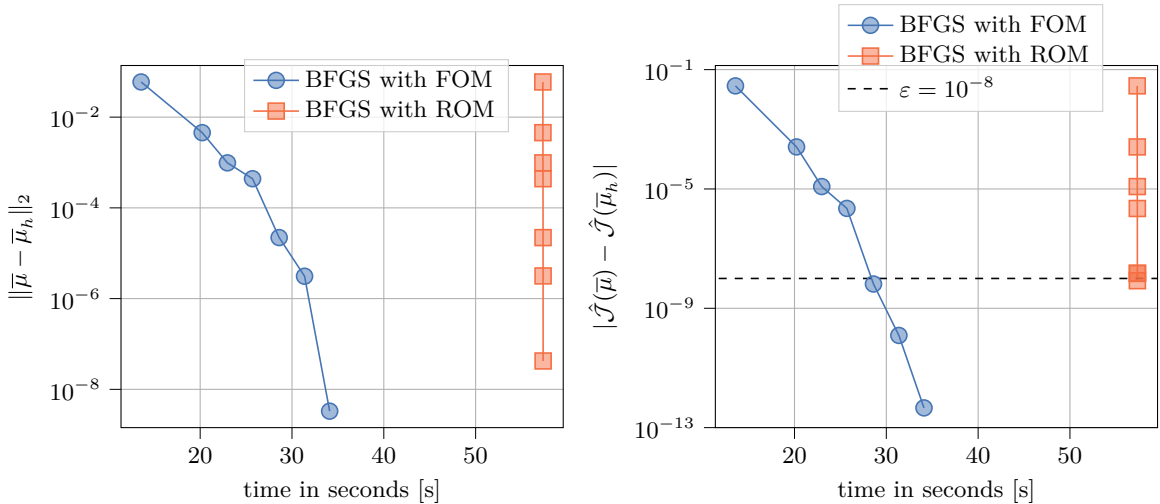
\begin{figure}
	\centering \footnotesize
\begin{tikzpicture}

\definecolor{color0}{rgb}{0.65,0,0.15}
\definecolor{color1}{rgb}{0.84,0.19,0.15}
\definecolor{color2}{rgb}{0.96,0.43,0.26}
\definecolor{color3}{rgb}{0.99,0.68,0.38}
\definecolor{color4}{rgb}{1,0.88,0.56}
\definecolor{color5}{rgb}{0.67,0.85,0.91}
\definecolor{color6}{rgb}{0.27,0.46,0.71}

\begin{axis}[
legend cell align={left},
legend style={fill opacity=0.8, draw opacity=1, text opacity=1, at={(1,0.5)}, anchor=west, draw=white!80!black, xshift=-4.2cm, yshift=2cm},
log basis y={10},
width=7.4cm,
tick align=outside,
tick pos=left,
x grid style={white!69.0196078431373!black},
xlabel={time in seconds [s]},
xmajorgrids,
xmin=11.3789805087261, xmax=59.5417068631388,
xtick style={color=black},
y grid style={white!69.0196078431373!black},
ylabel={$\|\muBar - \muBar_h \|_2$},
ymajorgrids,
ymin=1.43358704972122e-09, ymax=0.137346959465612,
ymode=log,
ytick style={color=black}
]
\addplot [semithick, color6, mark=*, mark size=3, mark options={solid, fill opacity=0.5}]
table {%
13.5681953430176 0.0595701054235534
20.2206008434296 0.00454837491910481
22.9678139686584 0.000981823017957724
25.7047007083893 0.000438749415106937
28.6020307540894 2.21513967971858e-05
31.3699910640717 3.11480015640991e-06
34.1037094593048 3.30532942670661e-09
};
\addlegendentry{BFGS with FOM}
\addplot [semithick, color2, mark=square*, mark size=3, mark options={solid, fill opacity=0.5}]
table {%
57.3219746891409 0.0595701054235534
57.3296977821738 0.00454834446743468
57.3342732731253 0.000981925795483599
57.3393246475607 0.000438770217454726
57.3440281692892 2.21103172356924e-05
57.3483583275229 3.16105762533912e-06
57.3524920288473 4.22237955834022e-08
};
\addlegendentry{BFGS with ROM}
\end{axis}

\end{tikzpicture}
\begin{tikzpicture}

\definecolor{color0}{rgb}{0.65,0,0.15}
\definecolor{color1}{rgb}{0.84,0.19,0.15}
\definecolor{color2}{rgb}{0.96,0.43,0.26}
\definecolor{color3}{rgb}{0.99,0.68,0.38}
\definecolor{color4}{rgb}{1,0.88,0.56}
\definecolor{color5}{rgb}{0.67,0.85,0.91}
\definecolor{color6}{rgb}{0.27,0.46,0.71}

\begin{axis}[
legend cell align={left},
legend style={fill opacity=0.8, draw opacity=1, text opacity=1, at={(1,0.5)}, anchor=west, draw=white!80!black, xshift=-4.2cm, yshift=2.5cm},
log basis y={10},
width=7.4cm,
tick align=outside,
tick pos=left,
x grid style={white!69.0196078431373!black},
xlabel={time in seconds [s]},
xmajorgrids,
xmin=11.3789805087261, xmax=59.5417068631388,
xtick style={color=black},
y grid style={white!69.0196078431373!black},
ylabel={$|\Jhat(\muBar) - \Jhat(\muBar_h)|$},
ymajorgrids,
ymin=1e-13, ymax=0.137346959465612,
ymode=log,
ytick style={color=black}
]
\addplot [semithick, color6, mark=*, mark size=3, mark options={solid, fill opacity=0.5}]
table {%
	13.5681953430176 0.0285193569089879
	20.2206008434296 0.000256583741397698
	22.9678139686584 1.20366392426519e-05
	25.7047007083893 2.2304496809511e-06
	28.6020307540894 6.51699316733811e-09
	31.3699910640717 1.22980736705358e-10
	34.1037094593048 4.58300064565265e-13
};
\addlegendentry{BFGS with FOM}
\addplot [semithick, color2, mark=square*, mark size=3, mark options={solid, fill opacity=0.5}]
table {%
	57.3219746891409 0.0285193569089879
	57.3296977821738 0.000256570483964325
	57.3342732731253 1.20238219647462e-05
	57.3393246475607 2.21576590053729e-06
	57.3440281692892 8.40050873307518e-09
	57.3483583275229 1.47985401710571e-08
	57.3524920288473 1.49275125593817e-08
};
\addlegendentry{BFGS with ROM}
\addplot [semithick, black, mark size=0, dashed]
table {%
	10 1e-08
	60 1e-08
};
\addlegendentry{$\varepsilon = 10^{-8}$}
\end{axis}

\end{tikzpicture}
	\caption[Total run time comparison FOM vs. global ROM.]{Total run time comparison. Depicted is the optimization error w.r.t. the parameter (left) and its objective value (right).}
	\hspace{0.5cm}
	\label{fig:timings_Greedy}		
\end{figure}

In \cref{fig:timings_Greedy}, we anticipate exemplary timings for a specific instance of a PDE-constrained optimization problem with a 2-dimensional parameter space; see \cref{sec:proof_of_concept} for a detailed introduction.
\cref{fig:timings_Greedy} illustrates the optimization steps needed for finding an approximate local optimum $\muBar_h$.
In particular, it is shown how the optimization error evolves depending on the algorithm's required computational run time.
While the computational time for the FOM method scales linearly with the number of iteration counts, the online time for the ROM is almost zero, which results in an almost vertical plot on the time scale.
Admittedly, the offline time required beforehand exceeds the overall time needed for the entire FOM algorithm.
Therefore, \cref{fig:timings_Greedy} is certainly an example of \cref{fig:many-query}, where the usage of a reduced model does not pay off at all.

We remark that, for the goal-oriented greedy-search algorithm, we prescribed a tolerance of $\varepsilon = 10^{-8}$ in the objective functional, which resulted in an RB space of dimension $12$.
The choice of $\varepsilon$ is well visible in \cref{fig:timings_Greedy} since, while the first steps of FOM and ROM match perfectly during the algorithm, the lack of accuracy beyond $\varepsilon = 10^{-8}$ prevents the ROM algorithm from stepping further.
The accuracy in the optimal parameter seems to be sufficient in this experiment.
However, it is not generally given that the error in the parameter is of the same size as the prescribed tolerance for the functional.
On the contrary, in \cref{num_test:12params}, we present an experiment where it is much more challenging to find the optimal parameter although the approximate optimal function value is already found.

Recall that the presented example has a $2$-dimensional parameter space with an affine decomposition with $\Xi=2$ components.
RB methods are known to be very efficient for many-query simulations of such problem.
Nevertheless, in this experiment, the global construction of a ROM is already too costly.
The situation further abbreviates if the size of the parameter space is even beyond the standard dimensions of RB, which we introduced as the curse of dimensionality.
For this reason, in the aforementioned references, global RB methods were only considered helpful if not only one optimization problem is to be solved.

Adaptive reduced models can resolve the issue of a too large offline time for a single PDE-constrained parameter optimization problem.
As the primary motivation for such an adaptive construction, we recall that we are solely interested in finding an optimal point somewhere in the parameter space~$\Params$.
From an optimization perspective, it does not matter how we find this point as long as we know that our specific algorithm converges to a critical point w.r.t. the FOM.
To accomplish this, it is only essential that the ROM is sufficiently accurate in a neighborhood of a local optimum.
The rest of the parameter space $\Params$ is only relevant for pushing the optimizer towards this region.
In contrast, a greedy-search algorithm assumes that a specific tolerance shall be reached in the whole parameter space.
Adaptive ROM strategies may instead only update the model along the optimization path and ignore the rest of the parameter space accordingly.
Such methods have, for instance, been introduced in  \cite{BMV2020, GHH2016, ZF2015}.
Recently, based on \cite{YM2013}, a certified adaptive algorithm for RB methods was presented in \cite{QGVW2017}, where an error-aware trust-region algorithm is used.

\subsection{Outlook to this thesis}
\label{sec:background_outlook_rb}

The goal of this thesis with respect to model order reduction can now further be specified.
As mentioned at the beginning of this section, MOR methods are of major interest for many-query and real-time scenarios, and there still remain challenges that also occur in the "ideal world".
In this thesis, we are concerned with both many-query and real-time scenarios and aim to address the challenges of the curse of dimensionality as well as the inaccessible global discretization; cf.~\cref{sec:RB_challenges}.

In \cref{chap:TR_RB}, we are interested in further developing the TR-RB algorithm as a fast algorithm for finding the solution of a single PDE-constrained parameter optimization problem following the first-optimize-then-discretize approach.
Therefore, we consider the overall effort within an adaptive RB approach.
We show that even in cases where the curse of dimensionality is present, overall-efficiency can be reached since it is not required to construct a globally sufficient surrogate model but only a locally trustable model.
In the TR-RB algorithm, we adapt the local surrogate model following the optimization path.
We also present that a posteriori error estimation is vital for such an approach.

In \cref{chap:TSRBLOD}, we focus on model order reduction for localized approximation techniques that mainly account for the challenge of inaccessible global discretizations.
In particular, we deduce an entirely online efficient reduced model based on a multiscale method with local problems and thus, do not require global high-fidelity methods.
Importantly, in the presented approach, the local problems of the multiscale method are entirely hidden in the offline phase, reducing the multiscale method to a single (low dimensional) system which we achieve by a two-scale formulation of the reduced basis localized orthogonal decomposition method (TSRBLOD).
The resulting online-efficiency makes the approach especially suitable for real-time scenarios.

Finally, in \cref{chap:TR_TSRBLOD}, we combine both approaches from \cref{chap:TR_RB} and \cref{chap:TSRBLOD} to solve PDE-constrained parameter optimization problems that require localized methods for the primal state equation.
For this, the online efficiency of the TSRBLOD reduced model in \cref{chap:TSRBLOD} and the TR-RB algorithm for deducing overall-efficiency also for significant parameter dependencies are combined to TR-LRB methods.
We also show that the localized method lowers the curse of dimensionality in the system if the parameterization can spatially be localized.
\chapter{Trust-region reduced basis methods}\label{chap:TR_RB}
\noindent
In this chapter, we review our contribution to the development of trust-region reduced basis (TR-RB) methods for solving PDE-constrained parameter optimization problems like~\eqref{Phat}.
We focus on a computationally overall-efficient algorithm and do not consider approaches where the offline time is negligible; cf. the discussion in \cref{sec:background_RB_methods_for_PDEopt}.
The presented results are partly published in~\cite{paper2, paper1, paper3}.
While we consider these papers sequentially, we provide supplementary explanations, discuss further challenges, future research perspectives, and unpublished insights into the algorithms.

The main idea of TR-RB methods is to use an RB-based surrogate model as the model function in the TR sub-problem and to construct this surrogate model adaptively along the optimization path.
To yield a certified and robust algorithm, the TR methodology is enhanced by an error-aware version, where the trust-region is characterized by the error estimator of the RB surrogate model.
As we shall see in this chapter, while in \cite{QGVW2017} a basic TR-RB algorithm is presented, many improvements concerning the reduced model, the basis construction, and the algorithm itself are possible.

Having in mind the discussed challenges of global approximations (see \cref{sec:multiscale_problems} and \cref{sec:RB_challenges}), we emphasize that, for the remainder of the chapter, we solely deal with global RB methods, meaning that we do not use localized approaches for superseding FEM as the FOM method.
Instead, we refer to \cref{chap:TR_TSRBLOD}, where the presented TR algorithm is pursued to localized approaches.

Given $\mu_\mathsf{a},\mu_\mathsf{b}\in\mathbb{R}^P$ with $P \in \N$, we consider the compact, bounded, and convex admissible parameter set
\[ 
\Paramsad:= \left\{\mu\in\mathbb{R}^P\,|\,\mu_\mathsf{a} \leq \mu \leq \mu_\mathsf{b} \right\} \subseteq \R^P,
\] 
where $\leq$ is to be understood component-wise.
Analogously to \cref{chap:background}, let $V$ be a real-valued Hilbert space with inner product $(\cdot \,,\cdot)$ and induced norm $\|\cdot\|$.
As motivated in \cref{sec:background_PDE-constr}, we are interested in approximating PDE-constrained parameter optimization problems.
In this chapter, we assume the cost functional $\J$ to satisfy the following linear-quadratic continuous structure:
\begin{assumption}[Linear-quadratic continuous objective functional] \label{asmpt:quadratic_J}
	The cost functional $\J$ can be written in the form
	\begin{equation}\label{eq:quadratic_J}
	\J: V \times \Params \to \R, \quad  (u,\mu) \mapsto \J(u, \mu)
	\coloneqq \Theta(\mu) + j_\mu(u) + k_\mu(u, u),
	\end{equation}
	where $\Theta: \Params \to \mathbb{R}$ denotes a parameter function and, for each $\mu \in \Params$, $j_\mu \in V'$ is a parameter-dependent continuous linear functional, and $k_\mu: V \times V \to \R$ a continuous symmetric bilinear form.	
\end{assumption}
We emphasize that this formulation can represent many objective functions, for instance, the already stated $L^2$-misfit from~\eqref{eq:example_J}, also used in \cref{sec:TRRB_benchmark_problem}.
Furthermore, it enables a direct algebraic formulation and thus simplifies the assembly of the right-hand side of the dual problem.
The theory in this chapter is restricted to the symmetric linear-quadratic case but can be generalized to other objective functionals, simply by computing the Fr\`echet derivatives accordingly.
Instead, the involved parameter functionals $\Theta$ and $\theta_\xi$ are arbitrary as long as they are differentiable.
Resuming \cref{def:red_parameter_optimal_control}, we consider the following reduced optimization problem, where the term "reduced" is associated with the optimization problem, not to reduced-order models:
\begin{align}
	\min_{\mu \in \Paramsad} \Jhat(\mu),
	\tag{$\hat{\textnormal{P}}$}\label{Phat_second}
\end{align}
with $\Jhat(\mu) \coloneqq \J(u_{\mu}, \mu)$.
Moreover, to deduce derivatives of the involved terms, we remind the reader to the differentiability assumption, stated in \cref{asmpt:differentiability}.

The chapter is organized as follows:
First, in \cref{sec:TRRB_MOR_for_PDEconstr}, we combine the theory of \cref{sec:background_PDE-constr} and \cref{sec:background_npdgl} to formulate the finite-dimensional FOM for the optimization problem.
Second, we explain several advances for a corresponding ROM and focus on the necessity of a new correction term in the reduced functional.
Subsequently, in \cref{sec:TRRB_a_post_error_estimates}, we elaborate on a posteriori error bounds for the derived reduced quantities. 
\cref{sec:TR_RB_algorithm} is devoted to presenting the TR-RB algorithm, followed by an extensive convergence study.
Furthermore, a list of derived TR-RB variants is given.
After a brief introduction to the employed software (cf.~\cref{sec:software_TRRB}), we numerically analyze the method and compare it to existing approaches from the literature in \cref{sec:TRRB_num_experiments}.
Lastly, in \cref{sec:TRRB_related_approaches}, we discuss related techniques that incorporate further ideas for enhancing the TR-RB method. In particular, we introduce the relaxed TR-RB method which plays a significant role in \cref{chap:TR_TSRBLOD}.

\section{MOR for PDE-constrained parameter optimization}
\label{sec:TRRB_MOR_for_PDEconstr}

In \cref{sec:background_PDE-constr}, the theory for PDE-constrained parameter optimization problems has been carefully introduced for the case of the infinite-dimensional real-valued Hilbert space $V$.
We now intend to define the FOM for the optimization problem~\eqref{Phat} based on the grid-based numerical approximation of PDEs as discussed in \cref{sec:background_npdgl}.
Further, we introduce different ways of constructing a corresponding ROM for the optimization problem.
The resulting optimality system is, in general, not equivalent to a Ritz-Galerkin projection of the FOM onto a reduced space.
For this reason, we also introduce a non-conforming dual-corrected (NCD-corrected) approach in \cref{sec:TRRB_ncd_approach}, and additionally introduce a Petrov--Galerkin reduced model; cf. \cref{sec:TRRB_pg_approach}. 

\subsection{Full order model}
\label{sec:TRRB_fom}
In order to discretize the optimization problem analogously to \cref{sec:background_npdgl}, we assume that a finite-dimensional FE space $V_h \subset V$ is given.
The FOM version of the finite-dimensional optimality system of~\eqref{P} (and~\eqref{Phat}) in \cref{prop:first_order_opt_cond} can be obtained by a Ritz-Galerkin projection of equations~\eqref{eq:optimality_conditions} onto $V_h$.
To be precise, we seek for each $\mu \in \Params$ the solution $u_{h, \mu} \in V_h$ of \cref{def:discr_param_elliptic_problem}, i.e.
\begin{align}
\bformd(u_{h, \mu}, v_h) = \lformd(v_h) &&\textnormal{for all } \, v_h \in V_h,
\label{eq:state_h}
\end{align}
which, in the context of optimization, is called the \emph{discrete primal equation}.
Hence, with the primal residual from \cref{def:primal_residual}, we have $\resd^\pr(u_{h, \mu})[v_h] = 0$ for all $v_h \in V_h$, $\mu \in \Params$.
For each $\mu \in \Params$, the discrete approximation of $p_{\mu} \in V$ from \cref{def:dual_solution} can be stated as the solution $p_{h, \mu} \in V_h$ of the \emph{discrete dual equation}
\begin{align}
\bformd(v_h, p_{h, \mu}) = \partial_u \J(u_{h, \mu}, \mu)[v_h] = \jformd(v_h) + 2 \kformd(v_h, u_{h, \mu}) &&\textnormal{for all } \, v_h \in V_h,
\label{eq:dual_solution_h}
\end{align}
where the second equality directly follows from~\eqref{eq:quadratic_J} and \cref{rem:gateaux_wrt_V} applied on each argument of the symmetric bilinear form $\kformd(u_{h, \mu}, u_{h, \mu})$.
For given $u_h, p_h \in V_h$, we also define the dual residual $r_\mu^\du(u_h, p_h) \in V'_h$ associated with~\eqref{eq:dual_solution_h} by
\begin{align}\label{eq:dual_residual}
	r_\mu^\du(u_h, p_h)[q_h] := j_\mu(q_h) + 2k_\mu(q_h, u_h) - a_\mu(q_h, p_h)&&\textnormal{for all } \,q_h \in V_h.
\end{align}
Hence, we have $\resd^\du(u_{h, \mu}, p_{h, \mu})[q_h] = 0$ for all $q_h \in V_h$, $\mu \in \Params$.
Similarly, the \emph{discrete primal derivative equations} for solving for $d_{\nu} u_{h, \mu} \in V_h$ as well as \emph{discrete dual derivative equations} for solving for $d_{\nu} p_{h, \mu} \in V_h$ at any direction $\nu \in \mathbb{R}^P$ follow analogously to Propositions \ref{prop:solution_dmu_eta} and \ref{prop:dual_solution_dmu_eta}.
Furthermore, $\Jhat$ is approximated by the \emph{discrete reduced functional}, cf. the definition in \cref{def:discr_param_elliptic_problem}, i.e.
\begin{align}
\Jhat_h(\mu) := \J(u_{h, \mu}, \mu) = \LL(u_{h,\mu},\mu,p_h) && \textnormal{for all } \, p_h\in V_h,
\label{eq:Jhat_h}
\end{align}
with the Lagrangian functional $\LL$ from \cref{def:lagrangian_functional}, and where $u_{h, \mu} \in V_h$ is the solution of~\eqref{eq:state_h}.
The second equation of \eqref{eq:Jhat_h} follows trivially because of the conforming choice for the test- and ansatz space of the primal- and dual equation, i.e. $\resd^\pr(u_{h, \mu})[p_h] = 0$, for any $p_h \in V_h$.
Finally, we formulate the discrete version of~\eqref{Phat} as the optimization problem 
\begin{align}
\min_{\mu \in \Paramsad} \Jhat_h(\mu).
\tag{$\hat{\textnormal{P}}_h$}\label{Phat_h}
\end{align}
By $\bar\mu_h \in \Paramsad$ we denote a locally optimal solution of~\eqref{Phat_h}.
We emphasize again that the above-explained FOM follows a conforming choice, i.e., $u_{h,\mu}$ and $p_{h,\mu}$ belong to the same space~$V_h$.
Consequently, the first- and second-order optimality conditions in \cref{prop:first_order_opt_cond} and \cref{prop:second_order}, the first-order solution map derivatives in \cref{def:primal_sensitivities} and \cref{prop:gradient_Jhat}, and the computation of the Hessian in \cref{prop:Hessian_Jhat} from \cref{sec:background_PDE-constr}, can be derived analogously for $V_h$, with all quantities replaced by their discrete counterparts.
We do not restate these discrete versions for brevity.

As usual in the context of RB methods and as discussed in \cref{sec:error_bounds_for_fem}, we eliminate the issue of choosing an appropriate mesh size $h$ by assuming that the high-dimensional space $V_h$ is accurate enough to approximate the actual solution.

\begin{assumption}[The FOM is the ``truth'']
	\label{asmpt:truth}
	We assume that the primal discretization error $\|u_\mu - u_{h, \mu}\|$, the dual error $\|p_\mu - p_{h, \mu}\|$, the primal sensitivity errors $\|d_{\mu_i} u_\mu - d_{\mu_i} u_{h, \mu}\|$ and the dual sensitivity errors $\|d_{\mu_i} p_\mu - d_{\mu_i} p_{h, \mu}\|$ are negligible for all $\mu \in \Params$, $1 \leq i \leq P$.
	This implies that also $|\Jhat(\mu) - \Jhat_h(\mu)|$ can be considered negligible.
\end{assumption}
\noindent As a consequence, and also since this chapter is solely concerned about the discrete optimization problem~\eqref{Phat_h}, we often write $\J$ instead of $\J_h$.

Recall that the resulting FEM approximation may be computationally demanding, which is why, in what follows, we define suitable reduced versions for all involved quantities.
For the following subsections, let $V_\red^\pr, V_\red^\du \subset V_h$ be problem adapted primal and dual RB spaces of low dimensions $n := \dim V_\red^\pr$ and $m := \dim V_\red^\du$, the construction of which is detailed in \cref{sec:TRRB_construct_RB}.
We stress here that, for various reasons, $V_\red^\pr$ and $V_\red^\du$ might not coincide, which implies further discussions of the RB approximation of the optimality system~\eqref{eq:optimality_conditions}.

\subsection{Reduced model -- Standard approach}
\label{sec:TRRB_standard_approach}

We start with what we call the \emph{standard RB approach} for the optimality conditions, initially proposed in~\cite{QGVW2017}.
Analogously to \cref{sec:background_RB_method}, we obtain the RB approximation of the primal and dual equations as follows:

\begin{subequations}
	\label{eq:optimality_conditionsRB}
	\begin{itemize}
		\item RB approximation for~\eqref{eq:optimality_conditions:u}: For each $\mu \in \Params$, the primal RB variable $u_{\red, \mu} \in V_\red^\pr$ is defined through
		\begin{align}
		\bformd(u_{\red, \mu}, v_\red) = \lformd(v_\red) &\qquad \textnormal{for all } \, v_\red \in V_\red^\pr.
		\label{eq:state_red}
		\end{align}
		\item RB approximation for~\eqref{eq:optimality_conditions:p}: For each $\mu \in \Params$ and with $u_{\red, \mu} \in V_\red^\pr$ from above, the approximate dual/adjoint RB variable $p_{\red, \mu} \in V_\red^\du$ satisfies the approximate dual equation, defined by
		\begin{align}
		\bformd(q_\red, p_{\red, \mu}) = \partial_u \J(u_{\red, \mu}, \mu)[q_\red] = \jformd(q_\red) + 2 \kformd(q_\red, u_{\red, \mu}) &&\textnormal{for all } \, q_\red \in V_\red^\du.
		\label{eq:dual_solution_red}
		\end{align}
	\end{itemize}
\end{subequations}
Analogously to \cref{def:red_param_problem}, we define the \emph{RB solution map} $\mathcal{S}_\red : \Params \to V_\red^\pr$ by $\mu \mapsto u_{\red, \mu}$ and similarly the \emph{RB dual solution map} $\mathcal{A}_\red:\Params \to V_\red^\du$ by $\mu \mapsto p_{\red, \mu}$, where $u_{\red, \mu} \in  V_\red^\pr$ and $p_{\red, \mu} \in  V_\red^\du$ denote the primal and dual reduced solutions of~\eqref{eq:state_red} and~\eqref{eq:dual_solution_red}, respectively.
Equation~\eqref{eq:state_red} is a standard Galerkin projection of~\eqref{eq:optimality_conditions:u} onto the reduced space $V_\red^\pr$. Importantly, the reduced dual equation~\eqref{eq:dual_solution_red} is not a Galerkin projection of~\eqref{eq:optimality_conditions:p} onto $V_\red^\du$ since the right-hand side incorporates the primal RB variable $u_{\red, \mu}$.
Thus, we consider $p_{\red, \mu}$ only as an approximate dual RB variable, which, as we see later, plays a crucial role for the reduced gradient information.

To approximate $\Jhat_h$, the standard choice of an \emph{RB reduced functional} simply means to replace the discrete FOM solution $u_{h,\mu}$ by their RB approximation $u_{\red,\mu} \in  V_\red^\pr$, i.e.
\begin{align}
\Jnoncor_\red(\mu) := \J(u_{\red, \mu}, \mu)=\J(\mathcal S_r(\mu), \mu),
\label{eq:Jhat_red}
\end{align}
where $u_{\red, \mu} \in V_\red^\pr$ is the solution of~\eqref{eq:state_red}.
The resulting reduced optimization problem reads as finding a locally optimal solution $\bar \mu_\red \in \Params$ of
\begin{align}
\min_{\mu \in \Paramsad} \Jnoncor_\red(\mu)
\label{eq:Phat_red_uncorrected}\tag{$\hat{\textnormal{P}}^J_\red$}
\end{align}
and a solution of the optimality system~\eqref{eq:optimality_conditions} is approximated by the RB triple $(u_{\red,\bar\mu_\red},\bar\mu_\red,p_{\red,\bar\mu_\red})$.

As proposed in~\cite{QGVW2017}, for computing an approximation of the gradient of $\Jnoncor_\red$, the gradient from \cref{prop:gradient_Jhat} can be reduced by replacing $u_{\mu} \in V$ and $p_{\mu} \in V$ with their RB counterparts $u_{\red,\mu} \in  V_\red^\pr$ and $p_{\red,\mu} \in  V_\red^\du$.
However, it can not be guaranteed in general that the computed gradient is the actual gradient of $\Jnoncor_\red$, as long as $V^\pr_\red$ and $V^\du_\red$ are not equal.
To see this, we consider again the Lagrangian and note that, for $1\leq i\leq P$ and all $p\in V^\pr_\red$, it holds
\begin{equation}
\Jnoncor_\red(\mu)  = \mathcal{L}(u_{\red,\mu},\mu,p), \qquad
\label{eq:Gradient_lagrangian}
\big(\nabla_\mu \Jnoncor_\red(\mu)\big)_i  =  \partial_u \mathcal{L}(u_{\red,\mu},\mu,p)[d_{\mu_i}u_{\red,\mu}]
+ \partial_{\mu_i} \mathcal{L}(u_{\red,\mu},\mu,p),	
\end{equation}
with special emphasis on the fact that~\eqref{eq:Gradient_lagrangian} does not hold for $p \in V^\du_\red$ since we would not necessarily have $r_\mu^\pr(u_{\red, \mu})[p]=0$.
Now, following~\cite{QGVW2017}, we define the \emph{inexact gradient} $\noncorgrad_\mu \Jnoncor_\red:\Params \to \R^P$ by 
\begin{equation}
\label{naive:red_grad}
\big(\noncorgrad_\mu \Jnoncor_\red(\mu)\big)_i := \partial_{\mu_i}\J(u_{\red, \mu}, \mu) + \partial_{\mu_i} r_\mu^\pr(u_{\red, \mu})[p_{\red, \mu}] = \partial_{\mu_i} \mathcal{L}(u_{\red,\mu},\mu,p_{\red,\mu})
\end{equation}
for all $1 \leq i \leq P$ and $\mu \in \Params$, where $u_{\red, \mu} \in V_\red^\pr$ and $p_{\red, \mu} \in V_\red^\du$ denote the primal and approximate dual reduced solutions of~\eqref{eq:state_red} and~\eqref{eq:dual_solution_red}, respectively.
With the superscript~$\sim$, we stress that $\noncorgrad_\mu \Jnoncor_\red(\mu)$ is not the actual gradient of $\Jnoncor_\red$, but its approximation.
In the FOM case, due to the conforming case,~\eqref{naive:red_grad} is already the exact formulation of $\nabla_\mu\Jhat_h$, which we showed in \cref{prop:gradient_Jhat} and in the proof of \cref{prop:first_order_opt_cond}, where the $\partial_u \mathcal{L}$ term in~\eqref{eq:Gradient_lagrangian} vanishes.
Now that $p_{\red,\mu}$ is an element of $V_\red^\pr$, we can neither use~\eqref{eq:Gradient_lagrangian} nor do we have $\partial_u\mathcal L(u_{\red,\mu},\mu, p_{\red,\mu})[d_{\mu_i}u_{\red,\mu}]= 0$ since~\eqref{eq:dual_solution_red} is not the dual equation with respect to the optimization problem~\eqref{eq:Phat_red_uncorrected}, cf.~\cite[Section~1.6.4]{HPUU2009}.
This would only be true if $V^\du_\red\subseteq V^\pr_\red$.
Thus, with the choice made in~\cite{QGVW2017}, Equation~\eqref{naive:red_grad} defines only an approximation of the actual gradient of $\Jnoncor_\red$.

This inexact choice introduces an additional approximation error in reconstructing the solution of the optimality system~\eqref{eq:optimality_conditions}, which is well visible in our numerical experiments; cf. \cref{sec:TRRB_estimator_study}.
For the TR-RB algorithm, the standard RB approach leads to a significant lack of accuracy, requiring additional outer optimization steps to enrich the RB space further.
Based on the previous remarks, we propose to add a correction term to $\Jnoncor_\red$ that resolved the approximability issue of the non-conforming approach.

\subsection{Reduced model -- NCD-corrected approach} 
\label{sec:TRRB_ncd_approach}
As introduced in~\cite{paper1} and pursued in~\cite{paper2,banholzer2022trust}, we discuss a novel approach for the ROM of the optimization problem.
The presented strategy has been used in the context of adaptive finite elements in~\cite{BKR2000, Rannacher2006} and the primal-dual RB (or dual-weighted residual (DWR)) approach for linear output functionals~\cite[Section~2.4]{HAA2017}.
The main motivation for the DWR approach was the fact that the standard reduced (linear) output functional (not necessarily connected to an optimization problem) only admitted an error estimator that was linearly decreasing.
Instead, for the compliant case (the output functional is equal to the right-hand side of the PDE), it was possible to show an improved error estimate with quadratic convergence behavior.
It was shown that the use of a dedicated dual problem like \cref{def:dual_solution} enables to show such an improved estimate for a corrected functional that is equivalent to what is stated in \eqref{eq:Jhat_red_corected} below.
Although it is mathematically the same as the primal-dual approach, we emphasize that the following NCD-correction approach was found by proving \cref{prop:error_reduced_quantities} and, in the context of optimization problems, can primarily be justified by the use of the Lagrangian functional for the optimization problem.

In fact, closely related to the optimality conditions in \cref{sec:background_optimality_conditions}, we seek to minimize the Lagrangian corresponding to problem~\eqref{Phat_h}.
In the previous section, we showed that~\eqref{eq:Gradient_lagrangian} only holds if $p \in V^\pr_\red$, which is not aligned with the reduced model in Equations~\eqref{eq:optimality_conditionsRB}.
From an optimization point of view, we still aim to use the Lagrangian functional and thus, we define the \emph{NCD-corrected RB-reduced functional} by
\begin{align}
\cJhatn(\mu) := \mathcal{L}(u_{\red,\mu},\mu,p_{\red,\mu}) = \Jnoncor_\red(\mu) + r_\mu^\pr(u_{\red,\mu})[p_{\red,\mu}],
\label{eq:Jhat_red_corected}
\end{align} 
with $u_{\red, \mu} \in V_\red^\pr$ and $p_{\red, \mu} \in V_\red^\du$ the solutions of~\eqref{eq:state_red} and~\eqref{eq:dual_solution_red}, respectively.
In other words, we augment the standard RB approach from the previous section by using a true Lagrangian approach for the optimization problem.
Note that $\cJhatn$ coincides with the functional $\Jnoncor_\red$ in~\eqref{eq:Jhat_red} if $V_\red^\du = V_\red^\pr$.
To conclude, we consider the \emph{NCD-corrected RB-reduced optimization problem} of finding a locally optimal solution $\bar \mu_\red$ of
\begin{align}
\min_{\mu \in \Paramsad} \cJhatn(\mu).
\tag{$\hat{\textnormal{P}}_\red$}\label{Phat_\red}
\end{align}
It remains to show how the exact gradient of $\Jhat_\red$ can be computed efficiently.
Starting with the standard chain rule, we obtain the following result:
\begin{proposition}[Gradient of $\Jhat_\red$ -- Sensitivity Approach]
	\label{prop:true_corrected_reduced_gradient_sens}
	The $i$-th component of the true gradient $\nabla_\mu\cJhatn : \Params \to \R^P$ of $\cJhatn$ is given by
	\begin{align*}
		\big(\nabla_\mu \cJhatn(\mu)\big)_i & = \partial_{\mu_i}\J(u_{\red, \mu}, \mu) 
		+ \partial_{\mu_i} r_\mu^\pr(u_{\red, \mu})[p_{\red, \mu}] 
		+ r_\mu^\pr(u_{\red, \mu})[d_{\mu_i}p_{\red, \mu}]
		+ r_\mu^\du(u_{\red, \mu},p_{\red, \mu})[d_{\mu_i} u_{\red, \mu}] 
	\end{align*}
	for all $1 \leq i \leq P$ and $\mu \in \Params$, where $u_{\red, \mu} \in V_\red^\pr$ and $p_{\red, \mu} \in V_\red^\du$ solve~\eqref{eq:optimality_conditionsRB}, $d_{\mu_i}u_{\red, \mu} \in V_\red^\pr$ and $d_{\mu_i}p_{\red, \mu} \in V_\red^\du$ denote the derivatives of the RB primal and dual solution maps.
\end{proposition}
\begin{proof} 
	The assertion follows from the chain rule and \cref{rem:gateaux_wrt_V}, see also~\cite[Section~1.6.1]{HPUU2009}.
	Using the definition of the gradient and~\eqref{eq:Jhat_red_corected} we obtain
	\begin{align*}
		\big(\nabla_\mu \cJhatn(\mu)\big)_i = d_{\mu_i}\J(u_{\red, \mu}, \mu) + d_{\mu_i} r_\mu^\pr(u_{\red, \mu})[p_{\red, \mu}],
	\end{align*}
	where we recall that $d_{\mu_i}$ denotes the total derivative.
	For the first term, we have
	\begin{align*}
		d_{\mu_i}\J(u_{\red, \mu}, \mu) = \partial_u \J(u_{\red, \mu}, \mu)[d_{\mu_i} u_{\red, \mu}]
		+ \partial_{\mu_i} \J(u_{\red, \mu}, \mu)
	\end{align*}
	using the chain rule. 
	On the other hand, we obtain for the residual term
	\begin{align*}
		d_{\mu_i} r_\mu^\pr(u_{\red, \mu})[p_{\red, \mu}] 
		&= d_{\mu_i} l_\mu(p_{\red, \mu}) - d_{\mu_i} a_\mu(u_{\red, \mu}, p_{\red, \mu})\\
		&= \partial_{\mu_i} r_\mu^\pr(u_{\red, \mu})[p_{\red, \mu}]
		 + r_\mu^\pr(u_{\red, \mu})[d_{\mu_i}p_{\red, \mu}] 
		 - a(d_{\mu_i}u_{\red, \mu}, p_{\red, \mu}),												  
	\end{align*}
	using the definition of the residual~\eqref{eq:primal_residual} in the first equality and the chain rule and \cref{rem:gateaux_wrt_V} in the second equality.
	We obtain the desired result by combining the previous two calculations and using \eqref{eq:dual_solution_h} and \eqref{eq:dual_residual}.
\end{proof}

While this expression for the gradient is helpful to prove the a posteriori estimate for $\nabla_\mu \cJhatn(\mu)$ (see~\cref{prop:grad_Jhat_error_sens}, it is not ideal from a computational perspective.
Analogously to the explanations in \cref{sec:background_derivatives_of_J}, the reason is the presence of sensitivities, requiring additional solutions of $(2\cdot P)$ equations.
Luckily, we can again use auxiliary problems to diminish the additional complexity.
\begin{proposition}[Gradient of $\Jhat_\red$ -- Adjoint Approach]
	\label{prop:true_corrected_reduced_gradient_adj}
	The $i$-th component of the true gradient of $\cJhatn$ is given by
	\begin{align*}
	\big(\nabla_\mu \cJhatn(\mu)\big)_i & = \partial_{\mu_i}\J(u_{\red,\mu},\mu) 
	+ \partial_{\mu_i}r_\mu^\pr(u_{\red,\mu})[p_{\red,\mu}+w_{\red,\mu}] 
	- \partial_{\mu_i}r^\du_\mu (u_{\red,\mu},p_{\red,\mu})[z_{\red,\mu}] 
	\end{align*}
	where $u_{\red, \mu} \in V_\red^\pr$ and $p_{\red, \mu} \in V_\red^\du$ denote the RB approximate primal and dual solutions of~\eqref{eq:state_red} and~\eqref{eq:dual_solution_red}, $z_{\red,\mu} \in V_\red^\du$ solves
	\begin{equation}
	\label{A2:z_eq_}
	a_\mu(z_{\red,\mu},q) = -r_\mu^\pr(u_{\red,\mu})[q] \quad \textnormal{for all } q\in V_\red^\du,
	\end{equation}
	and $w_{\red,\mu} \in V_\red^\pr$ solves
	\begin{equation}
	\label{A2:w_eq_}
	a_\mu(v,w_{\red,\mu}) = r_\mu^\du(u_{\red,\mu},p_{\red,\mu})[v]-2k_\mu(z_{\red,\mu},v), \quad \textnormal{for all } v\in V^\pr_r.
	\end{equation}
\end{proposition}
The proof of this proposition is displaced to \cref{sec:derivation_of_ncd_grad_and_hess} and follows the idea of auxiliary variables as shown in~\cite[Section~1.6.2]{HPUU2009}. This reference is, in itself, an alternative derivation of the dual variable as stated in \cref{def:dual_solution}. 

In conclusion, we require four equations to evaluate the true gradient of the NCD-corrected functional.
Furthermore, it can easily be seen that, for $V^\pr_\red = V^\du_\red$, the right-hand side of~\eqref{A2:w_eq_} and~\eqref{A2:z_eq_} vanish, and thus, the formulation of the corrected gradient agrees with the theory of the conforming case.

For using Newton's method in our algorithm, we require a reduced Hessian, for which the reduced sensitivities are needed.

\begin{definition}[Partial derivatives of the RB primal and dual solution maps]
	\label{prop:red_sens_}
	The derivatives of the reduced primal and dual solution maps associated with~\eqref{eq:optimality_conditionsRB} in direction $\nu \in \mathbb{R}^P$ as the solutions $d_{\nu} u_{\red, \mu} \in V_\red^\pr$ and $d_{\nu} p_{\red, \mu} \in V_\red^\du$ are defined as
	\begin{align}
	a_\mu(d_{\nu} u_{\red, \mu}, v_\red) &= \partial_{\mu} r_\mu^\pr(u_{\red, \mu})[v_\red]\cdot\nu &&\textnormal{for all } \, v_\red \in V_\red^\pr \text{ and}\label{eq:true_primal_reduced_sensitivity}\\
	a_\mu(q_\red, d_{\nu} p_{\red, \mu}) &= -\partial_\mu a_\mu(q_\red, p_{r, \mu})\cdot\nu + d_\mu\partial_u \J(u_{\red, \mu}, \mu)[q_\red]\cdot\nu\notag\\
	&= \partial_\mu r_\mu^\du(u_{\red, \mu}, p_{\red, \mu})[q_\red]\cdot\nu + 2k_\mu(q_\red, d_{\nu} u_{\red, \mu}) &&\textnormal{for all } \, q_\red \in V_\red^\du,
	\label{eq:true_dual_reduced_sensitivity}
	\end{align}
	respectively, analogously to Propositions \ref{prop:solution_dmu_eta} and \ref{prop:dual_solution_dmu_eta}, where the last equality holds for quadratic functionals as in~\eqref{eq:quadratic_J}.
\end{definition}

Many choices for a reduced Hessian are possible.
As a standard RB version in the sense of \cref{sec:TRRB_standard_approach}, it is again feasible to consider the FOM Hessian from Proposition \ref{prop:Hessian_Jhat} and reduce it by replacing all FOM quantities with their respective reduced counterpart.
While this approach may be a better approximation of the FOM Hessian, it is not the true Hessian of the NCD-corrected functional, resulting in a quasi-Newton method.
To prevent an overload of the work at hand, we omit a further discussion of this approach.

The true Hessian of the NCD-corrected RB reduced functional can, again, be computed by following the approach in~\cite[Section~1.6.4]{HPUU2009}.
The derivation is deferred to \cref{sec:derivation_of_ncd_grad_and_hess}.
\begin{proposition}[Hessian of the NCD-corrected RB reduced functional]
	\label{prop:true_corrected_reduced_Hessian}
	Given a direction $\nu\in\Params$, the evaluation of the Hessian $\cHhatn$ of $\cJhatn$ is
	\[
	\begin{aligned}
		\cHhatn(\mu)\cdot\nu = \nabla_\mu \Big(&j_\mu(d_\nu u_{\red,\mu})+2k_\mu(u_{\red,\mu},d_\nu u_{\red,\mu})  - a_\mu(d_\nu u_{\red,\mu},p_{\red,\mu}+w_{\red,\mu}).\\
		& + r^\pr_\mu(u_{\red,\mu})[d_\nu p_{\red,\mu}+d_\nu w_{\red,\mu}] - 2 k_\mu (z_{\red,\mu},d_\nu u_{\mu,\red}) \\
		&  + a_\mu(z_{\red,\mu},d_\nu p_{\red,\mu}) - r^\du_\mu(u_{\red,\mu},p_{\red,\mu})[d_\nu z_{\red,\mu}]   \\
		& +\partial_\mu (\J(u_{\red,\mu},\mu) + r^\pr_\mu(u_{\red,\mu})[p_{\red,\mu}+w_{\red,\mu}]- r_\mu^\du(u_{\red,\mu},p_{\red,\mu})[z_{\red,\mu}])\cdot \nu \Big),
	\end{aligned}
	\]
	where $d_\nu u_{\red,\mu}, w_{\red,\mu} \in V^\pr_\red$ and $d_\nu p_{\red,\mu}, z_{\red,\mu} \in V^\du_r$ solve~\eqref{eq:true_primal_reduced_sensitivity},~\eqref{A2:w_eq},~\eqref{eq:true_dual_reduced_sensitivity} and~\eqref{A2:z_eq}, respectively.
	Furthermore, $d_\nu z_{\red,\mu} \in V^\du_r$ solves
	\begin{equation}
		\label{z_eq_sens_}
		a_\mu(d_\nu z_{\red,\mu},q) = -\partial_\mu(r_\mu^\pr(u_{\red,\mu})[q] + a_\mu(z_{\red,\mu},q))\cdot \nu + a_\mu(d_\nu u_{\red,\mu},q)
	\end{equation}
	for all $q\in V^\du_r$ and $w_{\red,\mu} \in V^\pr_r$ solves
	\begin{equation}
		\label{w_eq_sens_}
		\begin{aligned}
			a_\mu(v,d_\nu w_{\red,\mu}) = \partial_\mu (&r_\mu^\du(u_{\red,\mu},p_{\red,\mu})[v]-2k_\mu(z_{\red,\mu},v)-a_\mu(v,w_{\red,\mu}))\cdot\nu \\
			 &+2k_\mu(v,d_\nu u_{\red,\mu}-d_\nu z_{\red,\mu})-a_\mu(v,d_\nu p_{\red,\mu}), \quad \textnormal{for all }\, v\in V^\pr_r.
		\end{aligned}
	\end{equation}
\end{proposition}
We remark that the true Hessian $\cHhatn(\mu)$ can also be computed without the use of auxiliary functions $z_{\red, \mu}$ and $w_{\red, \mu}$ and their derivatives.
However, we would have to compute second-order derivatives of $u_{\red, \mu}$ and $p_{\red, \mu}$ which aggravates the computations and makes the Hessian inefficiently callable because the second direction can not be pulled out.

\subsection{Reduced model -- Sensitivity snapshots} 
\label{sec:TRRB_sensitivity_approach}

In terms of numerical accuracy w.r.t.~the FOM functional, a solution $d_{\mu_i} u_{\red, \mu} \in V_\red^\pr$ of~\eqref{eq:true_primal_reduced_sensitivity} does not necessarily need to be a good approximation of the FOM version $d_{\mu_i} u_{h, \mu} \in V_h$ even though $u_{h, \mu}$ is contained in $V_\red^\pr$.
The reason is that the high-dimensional sensitivities are not generally contained in the respective reduced space (cf.~Proposition~\ref{prop:primal_solution_dmui_error}).

To remedy this, it can be an idea to compute the FOM sensitivities for all canonical directions and either include them to the respective primal and dual space $V_\red^\pr$ and $V_\red^\du$ (thus forming primal and dual Taylor RB spaces; cf. \cite{HAA2017}) or distribute all directional sensitivities to problem adapted RB spaces for the primal and dual sensitivities w.r.t.~all canonical directions, i.e. $V_\red^{\pr,d_{\mu_i}}, V_\red^{\du,d_{\mu_i}} \subset V_h$ for all $i$.
Just as the initial non-conforming, this choice for the primal and dual reduced spaces, is considered a variational crime and is formulated next.

\begin{definition}[Approximate partial derivatives of the RB primal and dual solution maps]
	\label{prop:red_sens}
	Considering the reduced primal and dual solution maps $\mathcal{S}_\red$ and $\mathcal{A}_\red$ of~\eqref{eq:state_red} and~\eqref{eq:dual_solution_red}, respectively, we define their approximate partial derivatives w.r.t.~the $i$-th component of $\mu$ by $\dred{\mu_i} u_{\red, \mu} \in V_\red^{\pr,d_{\mu_i}}$ and $\dred{\mu_i} p_{\red, \mu} \in V_\red^{\du,d_{\mu_i}}$, respectively, as solutions of the sensitivity equations
	\begin{align}
		\bformd(\dred{\mu_i} u_{\red, \mu}, v_\red) &= \partial_\mu \resd^\pr(u_{\red, \mu})[v_\red] \cdot e_i&&\textnormal{for all } \, v_\red \in V_\red^{\pr,d_{\mu_i}},
		\label{eq:red_sens_pr} \\
		\bformd(q_\red, \dred{\mu_i} p_{\red, \mu}) &= \partial_\mu \resd^\du(u_{\red, \mu}, p_{\red, \mu})[q_\red] \cdot e_i 
		+ 2\kformd(q_\red, \dred{\mu_i} u_{\red, \mu})&&\textnormal{for all } \,q_\red \in V_\red^{\du,d_{\mu_i}}.
		\label{eq:red_sens_du}
	\end{align}
	Similarly, we denote the approximate partial derivatives in direction $\nu \in \mathbb{R}^P$ by $\dred{\nu} u_{\red, \mu}$ and $\dred{\nu} p_{\red, \mu}$, respectively, defined by substituting $e_i$ with $\nu$ above.
\end{definition}
Following Propositions \ref{prop:solution_dmu_eta} and \ref{prop:dual_solution_dmu_eta}, we obtain  $\dred{\mu_i} u_{\red, \mu} = d_{\mu_i} u_{\red, \mu}$, if $V_\red^{\pr,d_{\mu_i}} =V_\red^{\pr}$, and $\dred{\mu_i} p_{\red, \mu} = d_{\mu_i} p_{\red, \mu}$, if $V_\red^{\du,d_{\mu_i}} =V_\red^{\du}$.
Moreover, the approximate partial derivatives depend on the choice of the corresponding reduced approximation spaces.
Now, we can also formulate the reduced gradient as a standard approach, but this time, the approximated sensitivities are used.

\begin{definition}[Approximate gradient of the NCD-corrected RB reduced functional]
	\label{prop:red_approx_grad}
	We define the approximate gradient $\gradred{\mu}\cJhatn:\Params \to \R^P$ of $\cJhatn$ by
	\begin{align}
		\big(\gradred{\mu} \cJhatn(\mu)\big)_i :=  \partial_{\mu_i}\J(u_{\red,\mu}, \mu) 
		&+ \partial_{\mu_i}r_\mu^\pr(u_{\red,\mu})[p_{\red,\mu}] 
		+ r_\mu^\pr(u_{\red, \mu})[\dred{\mu_i}p_{\red, \mu}] + r_\mu^\du(u_{\red, \mu},p_{\red, \mu})[\dred{\mu_i} u_{\red, \mu}], 
	\end{align}
	for $1 \leq i \leq P$, where $u_{\red,\mu} \in V_\red^\pr$, $p_{\red,\mu} \in V_\red^\du$ denote the reduced primal and dual solutions and $\dred{\mu_i} u_{\red, \mu} \in V_\red^{\pr,d_{\mu_i}}$ and $\dred{\mu_i} p_{\red, \mu} \in V_\red^{\du,d_{\mu_i}}$ denote the solutions of~\eqref{eq:red_sens_pr} and~\eqref{eq:red_sens_du}. 
\end{definition}
\noindent
Analogously, Hessian information can be derived similarly by replacing the respective solutions of \cref{prop:Hessian_Jhat} in a standard way.

As can be seen in our numerical experiments in \cref{sec:TRRB_estimator_study}, if FOM sensitivities are used for constructing the respective reduced spaces, the resulting gradient information has the strongest accuracy w.r.t. the FOM gradient.
This is also expected since we add the most FOM information to the reduced spaces.
However, for every enrichment step, the amount of FOM computations scales with the dimension of the parameter space.
This computational burden for constructing the sensitivity-based RB spaces is why we entirely omit such a strategy in the TR-RB strategies, which is further discussed in \cref{sec:TRRB_construct_RB}.

\subsection{Reduced model -- Petrov--Galerkin approach} 
\label{sec:TRRB_pg_approach}

\cref{sec:TRRB_ncd_approach} showed how the NCD-correction term attains a more suitable reduced model.
As shown in~\cite{paper3}, it can also be an idea to use a Petrov--Galerkin (PG) projection of the involved FOM or ROM systems.
Then, no correction term is required, i.e.~\eqref{eq:Gradient_lagrangian} is not violated, although we plug a dual solution $p^\textnormal{pg}_{\red, \mu} \in V_\red^\du$ into the Lagrangian.
Furthermore, with the PG projection, we can extract actual gradient and Hessian information from the formulas that are used for the FOM, which was not possible in \cref{sec:TRRB_standard_approach} or \cref{sec:TRRB_ncd_approach}.
To be precise, in the Petrov--Galerkin approach for the optimality system, we interchange the test spaces for the primal and dual equations.
In the FOM formulation, we use a FE space $V_h^\du$ as test space for~\eqref{eq:state_h} and a FE space $V_h^\pr$, potentially different from $V_h^\du$, as a test space for~\eqref{eq:dual_solution_h}.
Existence of solutions of PG approaches can be obtained by the inf-sup stability, which reads 
\begin{equation} \label{eq:inf_sup_trivial}
0 < \infsup{\mu} := \inf_{0\neq w_h\in V^\pr_h}\sup_{0\neq v_h\in V^\du_h}
\frac{a_\mu(w_h, v_h)}{\norm{w_h}\norm{v_h}}.
\end{equation}
We refer to \cite{bathe2001inf} and the references therein for more information on inf-sup stable problems.
A suitable reduced version of the PG-based FOM optimality system with PG reduced quantities can be obtained as follows:
\begin{subequations}
	\label{eq:optimality_conditionsRBPG}
	\begin{itemize}
		\item PG-RB approximation for~\eqref{eq:optimality_conditions:u}: 
		For each $\mu \in \Params$ the primal PG-RB variable $u^\textnormal{pg}_{\red, \mu} \in V_\red^\pr$ is defined through
		\begin{align}
		\bformd(u^\textnormal{pg}_{\red, \mu}, v_\red) = \lformd(v_\red) &\qquad \textnormal{for all } v_\red \in V_\red^\du.
		\label{eq:state_red_pg}
		\end{align}
		\item PG-RB approximation for~\eqref{eq:optimality_conditions:p}: For each $\mu \in \Params$, and with  $u^\textnormal{pg}_{\red, \mu} \in V_\red^\pr$ from above, the dual/adjoint PG-RB variable $p^\textnormal{pg}_{\red, \mu} \in V_\red^\du$ satisfies the approximate dual equation, defined by 
		\begin{align}
		\bformd(q_\red, p^\textnormal{pg}_{\red, \mu}) = \partial_u \J(u^\textnormal{pg}_{\red, \mu}, \mu)[q_\red] = \jformd(q_\red) + 2 \kformd(q_\red, u^\textnormal{pg}_{\red, \mu}) &&\textnormal{for all } q_\red \in V_\red^\pr.
		\label{eq:dual_solution_red_pg}
		\end{align}
	\end{itemize}
\end{subequations}
We define the PG-RB reduced optimization functional by 
\begin{align}
\cJhatn^\textnormal{pg}(\mu) := \J(u^\textnormal{pg}_{\red, \mu}, \mu) 
\label{eq:Jhat_red_pg}
\end{align} 
with $u^\textnormal{pg}_{\red, \mu} \in V_\red^\pr$ being the solution of~\eqref{eq:state_red_pg}.
We then consider the \emph{Petrov--Galerkin RB reduced optimization problem} by finding a locally optimal solution $\bar \mu_\red$ of
\begin{align}
\min_{\mu \in \Paramsad} \cJhatn^\textnormal{pg}(\mu).
\tag{$\hat{\textnormal{P}}^\textnormal{pg}_\red$}\label{Phat_red_pg}
\end{align}
In contrast to the NCD approach, we do not need to correct the reduced functional, as the primal and dual solutions automatically satisfy $r_\mu^\pr(u^\textnormal{pg}_{\red,\mu})[p^\textnormal{pg}_{\red,\mu}] = 0$.
In particular, we compute the gradient of the reduced functional with respect to the parameters solely based on the primal and dual solution of the PG-RB approximation, i.e.
\begin{align} \label{eq:red_gradient}
\big(\nabla_\mu \cJhatn^\textnormal{pg}(\mu)\big)_i & = \partial_{\mu_i}\J(u^\textnormal{pg}_{\red,\mu},\mu) + \partial_{\mu_i}r_\mu^\pr(u^\textnormal{pg}_{\red,\mu})[p_{\red,\mu}]
\end{align}
for all $1 \leq i \leq P$ and $\mu \in \Params$, where $u^\textnormal{pg}_{\red, \mu} \in V_\red^\pr$ and $p^\textnormal{pg}_{\red, \mu} \in V_\red^\du$ denote the PG-RB primal and dual reduced solutions of~\eqref{eq:state_red_pg} and~\eqref{eq:dual_solution_red_pg}, respectively.
For the Hessian, also the FOM formula can directly be used, which we do not discuss further.

\subsection{Derivation of the NCD-corrected gradient and Hessian}
\label{sec:derivation_of_ncd_grad_and_hess}
In what follows, we prove \cref{prop:true_corrected_reduced_gradient_adj} and \cref{prop:true_corrected_reduced_Hessian}, which are repeated in the following:
\begin{proposition_red_grad}[Gradient of $\Jhat_\red$ -- Adjoint Approach]
	The $i$-th component of the true gradient of $\cJhatn$ is given by
	\begin{align*}
		\big(\nabla_\mu \cJhatn(\mu)\big)_i & = \partial_{\mu_i}\J(u_{\red,\mu},\mu) 
		+ \partial_{\mu_i}r_\mu^\pr(u_{\red,\mu})[p_{\red,\mu}+w_{\red,\mu}] 
		- \partial_{\mu_i}r^\du_\mu (u_{\red,\mu},p_{\red,\mu})[z_{\red,\mu}] 
	\end{align*}
	where $u_{\red, \mu} \in V_\red^\pr$ and $p_{\red, \mu} \in V_\red^\du$ denote the RB approximate primal and dual solutions of \eqref{eq:state_red} and \eqref{eq:dual_solution_red}, $z_{\red,\mu} \in V_\red^\du$ solves
	\begin{equation}
		\label{A2:z_eq}
		a_\mu(z_{\red,\mu},q) = -r_\mu^\pr(u_{\red,\mu})[q] \quad \textnormal{for all } q\in V_\red^\du,
	\end{equation}
	and $w_{\red,\mu} \in V_\red^\pr$ solves
	\begin{equation}
		\label{A2:w_eq}
		a_\mu(v,w_{\red,\mu}) = r_\mu^\du(u_{\red,\mu},p_{\red,\mu})[v]-2k_\mu(z_{\red,\mu},v), \quad \textnormal{for all } v\in V^\pr_r.
	\end{equation}
\end{proposition_red_grad}
\begin{proof}
	Recall, that the corrected functional is equal to the Lagrangian functional evaluated at the tuple $(u_{\red,\mu},\mu,p_{\red,\mu})$.  
	Following the observations in \cite[Section~1.6.2]{HPUU2009}, we require the dual pairing $\langle\cdot,\cdot\rangle_{V',V}$, which, in our case, can be interpreted as $\langle L, v\rangle_{V',V} = L[v]$ for $L \in V'$ and $v \in V$.
	We have
	\[
	\begin{aligned}
		\partial_\mu \cJhatn(\mu)\cdot\nu 
		&= \langle \partial_u \LL(u_{\red,\mu},\mu,p_{\red,\mu}),\partial_\mu
		u_{\red,\mu}\cdot\nu\rangle_{V',V} \\
		& \qquad + \partial_\mu \LL(u_{\red,\mu},\mu,p_{\red,\mu})\cdot\nu \\
		& \qquad +  \langle \partial_p \LL(u_{\red,\mu},\mu,p_{\red,\mu}),\partial_\mu 
		p_{\red,\mu}\cdot\nu\rangle_{V',V} \\
		& = \left((\partial_\mu u_{\red,\mu})^*\partial_u \LL(u_{\red,\mu},\mu,p_{\red,\mu})\right)\cdot\nu \\
		& \qquad + \partial_\mu \LL(u_{\red,\mu},\mu,p_{\red,\mu})\cdot\nu \\
		& \qquad +  \left((\partial_\mu p_{\red,\mu})^*\partial_p 
		\LL(u_{\red,\mu},\mu,p_{\red,\mu})\right)\cdot\nu,
	\end{aligned}
	\]
	where ${}^*$ denotes the adjoint of an operator. We now use the operator $e: V \times \Params \to V'$, which we defined in \cref{def:abstract_param_equation} for an abstract and equivalent definition for \cref{def:param_elliptic_problem}, such that solving \eqref{eq:param_elliptic_problem} is equivalent to
	\begin{equation}
		\label{A2:state}
		e(u_{\red,\mu},\mu):= \mathcal{L}_\mu-\mathcal{A}_\mu(u_{\red,\mu}) = 0
	\end{equation}
	with appropriate operators $\mathcal{A}_\mu(u):= a_\mu(u,\cdot)$ and $\mathcal{L}_\mu:= l_\mu(\cdot)$, such that $e(u,\mu)[v]\equiv r^\pr_\mu(u)[v]$ for all $v\in V^\pr_\red$.
	By computing the derivative of \eqref{A2:state}, we then have
	\[
	\langle\partial_u e(u_{\red,\mu},\mu),\partial_\mu u_{\red,\mu}\cdot\nu\rangle_{V',V}+\partial_\mu e(u_{\mu,r},\mu)\cdot\nu = 0, \qquad \textnormal{for all } \nu\in\Params.
	\]
	Thus, we obtain
	\[
	(\partial_\mu u_{\red,\mu})^* = - (\partial_\mu e(u_{\red,\mu},\mu))^*(\partial_u e(u_{\red,\mu},\mu))^{-*},
	\]
	and hence,
	\[
	\left((\partial_\mu u_{\red,\mu})^*\partial_u \LL(u_{\red,\mu},\mu,p_{\red,\mu})\right)\cdot\nu 
	= - \left((\partial_\mu e(u_{\red,\mu},\mu))^*(\partial_u e(u_{\red,\mu},\mu))^{-*}\partial_u
	\LL(u_{\red,\mu},\mu,p_{\red,\mu})\right)\cdot\nu.
	\]
	Now, let $w^1_{\red,\mu}$ be the solution of
	\[
	(\partial_u e(u_{\red,\mu},\mu))^{*} w^1_{\red,\mu} = \partial_u \LL(u_{\red,\mu},\mu,p_{\red,\mu}),
	\]
	i.e. for all $v\in V^\pr_r$, we have
	\[
	\begin{aligned}
		& \qquad \qquad \langle (\partial_u e(u_{\red,\mu},\mu))^{*} w^1_{\red,\mu}, v\rangle_{V',V} &&= \langle \partial_u \LL(u_{\red,\mu},\mu,p_{\red,\mu}), v\rangle_{V',V} \\
		&\Leftrightarrow \quad \qquad
		 \langle (\partial_u e(u_{\red,\mu},\mu)) v , w^1_{\red,\mu}\rangle_{V',V} && = \langle \partial_u \LL(u_{\red,\mu},\mu,p_{\red,\mu}), v\rangle_{V',V} \\
		&\Leftrightarrow \quad \qquad \quad \qquad \langle -A_\mu(v) , w^1_{\red,\mu}\rangle_{V',V} && = \langle \partial_u \J(u_{\red,\mu},\mu)+\partial_u r_\mu^\pr(u_{\red,\mu})[p_{\red,\mu}], v \rangle_{V',V} \\
		&\Leftrightarrow \quad \qquad \qquad \qquad \qquad -a_\mu(v,w^1_{\red,\mu}) && = \partial_u \J(u_{\red,\mu},\mu)[v]-a_\mu(v,p_{\red,\mu}).
	\end{aligned}
	\]
	So $w^1_{\red,\mu}$ solves
	\begin{equation}
		\label{A2:w1_eq}
		a_\mu(v,w^1_{\red,\mu}) = -r_\mu^\du(u_{\red,\mu},p_{\red,\mu})[v], \quad \textnormal{for all } v\in V^\pr_r.
	\end{equation}
	Moreover, we have
	\[
	\begin{aligned}
		((\partial_\mu e(u_{\red,\mu},\mu))^* w^1_{\red,\mu})\cdot \nu & = \langle \partial_\mu e(u_{\red,\mu},\mu) \cdot \nu, w^1_{\red,\mu} \rangle_{V',V} \\
		& = \partial_\mu l_\mu(w^1_{\red,\mu})\cdot \nu - \partial_\mu a_\mu(u_{\red,\mu},w^1_{\red,\mu}) \cdot \nu,
	\end{aligned}
	\]
	and therefore,
	\[
	\left((\partial_\mu u_{\red,\mu})^*\partial_u \LL(u_{\red,\mu},\mu,p_{\red,\mu})\right)\cdot\nu =
	\partial_\mu a_\mu(u_{\red,\mu},w^1_{\red,\mu}) \cdot \nu -\partial_\mu l_\mu(w^1_{\red,\mu})\cdot
	\nu = -\partial_\mu r_\mu^\pr(u_{\red,\mu})[w^1_{\red,\mu}]\cdot\nu,
	\]
	where $w^1_{\red,\mu}$ solves \eqref{A2:w1_eq}.
	Note that since the right-hand side of \eqref{A2:w1_eq} is the dual residual, the FOM counterpart of $w^1_{\red,\mu}$ is zero. Thus, using the RB space of the primal variable for $w^1_{\red,\mu}$ allows to compute it cheaply without additional effort for performing an enrichment. The third piece of the gradient can be computed similarly.
	Let us introduce first $A^\du_\mu(p):= a_\mu(\cdot\,,p)$ and $L^\du_\mu:= \partial_u J(u_{\red,\mu},\mu)[\cdot]= j_\mu(\cdot)+2k_\mu(\cdot,u_{\red,\mu})$.
	Moreover, let 
	\[
	e^\du(p,\mu):= L_\mu^\du-A_\mu^\du(p).
	\]
	As before, we have
	\[
	\langle \partial_p e^\du(p_{\red,\mu},\mu),\partial_\mu p_{\red,\mu}\cdot\nu\rangle_{V',V} + \partial_\mu e^\du(p_{\red,\mu},\mu)\cdot\nu = 0 \quad \textnormal{for all } \nu\in\Params,
	\]
	and thus
	\[
	(\partial_\mu p_{\red,\mu})^* = -( \partial_\mu e^\du(p_{\red,\mu},\mu))^*(\partial_p e^\du(p_{\red,\mu},\mu))^{-*}.
	\]
	Hence, let $z_{\red,\mu}$ be the solution of
	\[
	(\partial_p e^\du(p_{\red,\mu},\mu))^{*}z_{\red,\mu} = \partial_p \LL(u_{\red,\mu},\mu,p_{\red,\mu}),
	\]
	and we obtain for all $q\in V^\du_r$:
	\[
	\begin{aligned}
		&\qquad \quad \langle(\partial_p e^\du(p_{\red,\mu},\mu))^{*} z_{\red,\mu}, q\rangle_{V',V} &&= \langle \partial_p \LL(u_{\red,\mu},\mu,p_{\red,\mu}), q\rangle_{V',V} \\
		&\Leftrightarrow \qquad \langle (\partial_p e^\du(p_{\red,\mu},\mu)) q , z_{\red,\mu}\rangle_{V',V} && = \langle \partial_p \LL(u_{\red,\mu},\mu,p_{\red,\mu}), q\rangle_{V',V} \\
		&\Leftrightarrow \quad\qquad \qquad \langle -A^\du_\mu(q) , z_{\red,\mu}\rangle_{V',V} && = \langle \partial_p r_\mu^\pr(u_{\red,\mu})[p_{\red,\mu}], q \rangle_{V',V} \\
		&\Leftrightarrow \qquad\qquad\qquad \qquad -a_\mu(z_{\red,\mu},q) && = r_\mu^\pr(u_{\red,\mu})[q].
	\end{aligned}
	\]
	So $z_{\red,\mu}$ solves \eqref{A2:z_eq}.
	Moreover
	\[
	\begin{aligned}
		((\partial_\mu e^\du(p_{\red,\mu},\mu))^* z_{\red,\mu})\cdot \nu & = \langle \partial_\mu e^\du(p_{\red,\mu},\mu) \cdot \nu, z_{\red,\mu} \rangle_{V',V} \\
		& = \langle\partial_\mu (\partial_u \J(u_{\red,\mu},\mu)[\cdot]-a_\mu(\cdot,p_{\red,\mu}))\cdot \nu, z_{\red,\mu}\rangle_{V',V} \\
		& = \langle\partial_\mu(j_\mu(\cdot)+2k_\mu(\cdot,u_{\red,\mu})-a_\mu(\cdot,p_{\red,\mu}))\cdot\nu, z_{\red,\mu}\rangle_{V',V} \\
		& = \partial_\mu j_{\mu}(z_{\red,\mu})\cdot\nu + 2\partial_\mu k_\mu(z_{\red,\mu},u_{\red,\mu})\cdot\nu \\
		& \quad  + \langle 2\partial_v k_\mu(z_{\red,\mu},u_{\red,\mu}),\partial_\mu u_{\red,\mu}\cdot\nu\rangle_{V',V} -\partial_\mu a_\mu(z_{\red,\mu},p_{\red,\mu})\cdot\nu \\
		& = \partial_\mu j_{\mu}(z_{\red,\mu})\cdot\nu + 2\partial_\mu k_\mu(z_{\red,\mu},u_{\red,\mu})\cdot\nu \\
		& \quad + \left((\partial_\mu u_{\red,\mu})^*2\partial_v k_\mu(z_{\red,\mu},u_{\red,\mu})\right)\cdot\nu -\partial_\mu a_\mu(z_{\red,\mu},p_{\red,\mu})\cdot\nu.
	\end{aligned}
	\]
	Note that $\partial_v$ here refers to the second argument of $k_\mu$.
	Proceeding as before, we have
	\[
	(\partial_\mu u_{\red,\mu})^*2\partial_v k_\mu(z_{\red,\mu},u_{\red,\mu}) = -(\partial_\mu e(u_{\red,\mu},\mu))^* w^2_{\red,\mu},
	\]
	where $w^2_{\red,\mu}$ solves
	\[
	\begin{aligned}
		\langle \partial_u e(u_{\red,\mu},\mu)v, w^2_{\red,\mu}\rangle_{V',V} &= \langle2\partial_v k_\mu(z_{\red,\mu},u_{\red,\mu}),v\rangle_{V',V} \quad \textnormal{for all } v\in V^\pr_r 
	\end{aligned}
	\]
	So, $w^2_{\red,\mu}$ solves 
	\begin{equation}
		\label{A2:w2_eq}
		a_\mu(v,w^2_{\red,\mu}) = -2k_\mu(z_{\red,\mu},v) \quad \textnormal{for all } v\in V^\pr_r.
	\end{equation}
	Moreover
	\[
	\begin{aligned}
		((\partial_\mu e(u_{\red,\mu},\mu))^* w^2_{\red,\mu})\cdot \nu & = \langle \partial_\mu e(u_{\red,\mu},\mu) \cdot \nu, w^2_{\red,\mu} \rangle_{V',V} \\
		& = \partial_\mu l_\mu(w^2_{\red,\mu})\cdot \nu - \partial_\mu a_\mu(u_{\red,\mu},w^2_{\red,\mu}) \cdot \nu.
	\end{aligned}
	\]
	Therefore 
	\[
	\begin{aligned}
		\left((\partial_\mu p_{\red,\mu})^*\partial_p \LL(u_{\red,\mu},\mu,p_{\red,\mu})\right)\cdot\nu & = \partial_\mu a_\mu(z_{\red,\mu},p_{\red,\mu})\cdot\nu - \partial_\mu j_{\mu}(z_{\red,\mu})\cdot\nu - 2\partial_\mu k_\mu(z_{\red,\mu},u_{\red,\mu})\cdot\nu \\
		& \qquad - \partial_\mu a_\mu(u_{\red,\mu},w^2_{\red,\mu}) \cdot \nu +\partial_\mu l_\mu(w^2_{\red,\mu})\cdot \nu \\
		& = \partial_\mu \left( - r_\mu^\du(u_{\red,\mu},p_{\red,\mu})[z_{\red,\mu}] + r_\mu^\pr(u_{\red,\mu})[w^2_{\red,\mu}]\right)\cdot \nu,
	\end{aligned}
	\]
	where $z_{\red,\mu}$ and $w^2_{\red,\mu}$ solve \eqref{A2:z_eq} and \eqref{A2:w2_eq}, respectively. 
	At last, 
	\[
	\partial_\mu \LL(u_{\red,\mu},\mu,p_{\red,\mu})\cdot\nu = \partial_\mu \J(u_{\red,\mu},\mu)\cdot\nu + \partial_\mu r_\mu^\pr(u_{\red,\mu})[p_{\red,\mu}]\cdot\nu.
	\]
	Summing everything together and defining $w_{\red,\mu}= w^2_{\red,\mu}-w^1_{\red,\mu}$, by linearity, we obtain the claim.
\end{proof}

\noindent
With the same strategy as before, we can also prove \cref{prop:true_corrected_reduced_Hessian}. 

\begin{proposition_red_Hessian}
	Given a direction $\nu\in\Params$, the evaluation of the Hessian $\cHhatn$ of $\cJhatn$ is
	\[
	\begin{aligned}
		\cHhatn(\mu)\cdot\nu = \nabla_\mu \Big( &j_\mu(d_\nu u_{\red,\mu})+2k_\mu(u_{\red,\mu},d_\nu u_{\red,\mu})  - a_\mu(d_\nu u_{\red,\mu},p_{\red,\mu}+w_{\red,\mu}) \\
		& + r^\pr_\mu(u_{\red,\mu})[d_\nu p_{\red,\mu}+d_\nu w_{\red,\mu}] - 2 k_\mu (z_{\red,\mu},d_\nu u_{\mu,\red}) \\
		&  + a_\mu(z_{\red,\mu},d_\nu p_{\red,\mu}) - r^\du_\mu(u_{\red,\mu},p_{\red,\mu})[d_\nu z_{\red,\mu}]   \\
		& +\partial_\mu (\J(u_{\red,\mu},\mu) + r^\pr_\mu(u_{\red,\mu})[p_{\red,\mu}+w_{\red,\mu}]- r_\mu^\du(u_{\red,\mu},p_{\red,\mu})[z_{\red,\mu}])\cdot \nu \Big)
	\end{aligned}
	\]
	where $d_\nu u_{\red,\mu}, w_{\red,\mu} \in V^\pr_\red$ and $d_\nu p_{\red,\mu}, z_{\red,\mu} \in V^\du_r$ solve \eqref{eq:true_primal_reduced_sensitivity}, \eqref{A2:w_eq}, \eqref{eq:true_dual_reduced_sensitivity} and \eqref{A2:z_eq}, respectively.
	Furthermore, $d_\nu z_{\red,\mu} \in V^\du_r$ solves
	\begin{equation}
		\label{z_eq_sens}
		a_\mu(d_\nu z_{\red,\mu},q) = -\partial_\mu(r_\mu^\pr(u_{\red,\mu})[q] + a_\mu(z_{\red,\mu},q))\cdot \nu + a_\mu(d_\nu u_{\red,\mu},q)
	\end{equation}
	for all $q\in V^\du_r$ and $w_{\red,\mu} \in V^\pr_r$ solves
	\begin{equation}
		\label{w_eq_sens}
		\begin{aligned}
			a_\mu(v,d_\nu w_{\red,\mu}) = \partial_\mu (&r_\mu^\du(u_{\red,\mu},p_{\red,\mu})[v]-2k_\mu(z_{\red,\mu},v)-a_\mu(v,w_{\red,\mu}))\cdot\nu \\
			& +2k_\mu(v,d_\nu u_{\red,\mu}-d_\nu z_{\red,\mu})-a_\mu(v,d_\nu p_{\red,\mu}), \quad \textnormal{for all }\, v\in V^\pr_r.
		\end{aligned}
	\end{equation}
\end{proposition_red_Hessian}
\begin{proof}
	\newcommand{\fdir}{\nu}
	\newcommand{\sdir}{\eta}
	Given two directions $\sdir,\fdir\in \R^P$, we have
	\[
	(\cHhatn(\mu)\cdot\fdir)\cdot\sdir = d_\fdir (\partial_\mu \cJhatn(\mu)\cdot \sdir).
	\]
	We now compute the derivative $d_\fdir$ for each piece of the gradient in Theorem~\ref{prop:true_corrected_reduced_gradient_adj}.
	First part:
	\[
	\begin{aligned}
		d_\fdir(\partial_\mu \J(u_{\red,\mu},\mu)\cdot \sdir) & = \partial_u(\partial_\mu \J(u_{\red,\mu},\mu)\cdot\sdir)[d_\fdir u_{\red,\mu}]+ \partial_\mu(\partial_\mu\J(u_{\red,\mu},\mu)\cdot\sdir)\cdot \fdir \\
		& = \partial_\mu (\partial_u \J(u_{\red,\mu},\mu)[d_\fdir u_{\red,\mu}]) + \partial_\mu(\partial_\mu\J(u_{\red,\mu},\mu)\cdot\fdir)\cdot \sdir\\
		& = \partial_\mu (j_\mu(d_\fdir u_{\red,\mu})+2k_\mu(u_{\red,\mu},d_\fdir u_{\red,\mu})+\partial_\mu \J(u_{\red,\mu},\mu)\cdot\fdir)\cdot \sdir,
	\end{aligned}
	\]
	where $d_\fdir u_{\red,\mu}$ solves \eqref{eq:primal_sens}.
	Second part:
	\[
	\begin{aligned}
		d_\fdir (\partial_\mu r^\pr_\mu(u_{\red,\mu})[p_{\red,\mu}+w_{\red,\mu}]\cdot \sdir) & = \partial_\mu(\partial_\mu r^\pr_\mu(u_{\red,\mu})[p_{\red,\mu}+w_{\red,\mu}]\cdot \sdir)\cdot \fdir \\ &\qquad + \partial_u (\partial_\mu r^\pr_\mu(u_{\red,\mu})[p_{\red,\mu}+w_{\red,\mu}]\cdot \sdir)[d_\fdir u_{\red,\mu}] \\&\qquad + \partial_v (\partial_\mu r^\pr_\mu(u_{\red,\mu})[p_{\red,\mu}+w_{\red,\mu}]\cdot \sdir) [d_\fdir p_{\red,\mu}+d_\fdir w_{\red,\mu}] \\
		& = \partial_\mu\left(\partial_\mu r^\pr_\mu(u_{\red,\mu})[p_{\red,\mu}+w_{\red,\mu}]\cdot \fdir - a_\mu(d_\fdir u_{\red,\mu},p_{\red,\mu}+w_{\red,\mu}) \right. \\
		& \quad \left. + r^\pr_\mu(u_{\red,\mu})[d_\fdir p_{\red,\mu}+d_\fdir w_{\red,\mu}]\right) \cdot\sdir,
	\end{aligned}
	\]
	where $d_\fdir p_{\red,\mu}$ solves \eqref{eq:dual_sens} and $d_\fdir w_{\red,\mu}$ solves \eqref{w_eq_sens}.
	Third part:
	\[
	\begin{aligned}
		d_\fdir (\partial_\mu r_\mu^\du(u_{\red,\mu},p_{\red,\mu})[z_{\red,\mu}]\cdot \sdir) = \partial_\mu (& \partial_\mu r_\mu^\du(u_{\red,\mu},p_{\red,\mu})[z_{\red,\mu}]\cdot \sdir)\cdot \fdir \\
		& + \partial_u (\partial_\mu r_\mu^\du(u_{\red,\mu},p_{\red,\mu})[z_{\red,\mu}]\cdot \sdir) [d_\fdir u_{\red,\mu}] \\
		& +  \partial_p (\partial_\mu r_\mu^\du(u_{\red,\mu},p_{\red,\mu})[z_{\red,\mu}]\cdot \sdir) [d_\fdir p_{\red,\mu}] \\
		& +  \partial_q (\partial_\mu r_\mu^\du(u_{\red,\mu},p_{\red,\mu})[z_{\red,\mu}]\cdot \sdir) [d_\fdir z_{\red,\mu}] \\
		& = \partial_\mu \left(\partial_\mu r_\mu^\du(u_{\red,\mu},p_{\red,\mu})[z_{\red,\mu}]\cdot \fdir + 2 k_\mu (z_{\red,\mu},d_\fdir u_{\mu,\red}) \right. \\
		&\quad \left. -a_\mu(z_{\red,\mu},d_\fdir p_{\red,\mu}) + r^\du_\mu(u_{\red,\mu},p_{\red,\mu})[d_\fdir z_{\red,\mu}]\right) \cdot\sdir,
	\end{aligned}
	\]
	where $d_\fdir z_{\red,\mu}$ solves \eqref{z_eq_sens}. Summing everything together, we obtain the claim.
\end{proof}

\subsection{Summary}

In this section, we derived several different models for the optimality system of \eqref{Phat}.
The full-order model from \cref{sec:TRRB_fom} can be used as a benchmark model with exact finite-dimensional approximations of the optimality system, which is ensured by \cref{asmpt:truth}.
To be precise, the primal and dual equations are simply replaced by their discrete counterparts and the objective functional $\Jhat$ as well as its first- and second-order derivatives can directly be computed analog to the finite-dimensional setting.

For the resulting finite-dimensional optimization problem \cref{Phat_h} and the corresponding optimality system, we showed that we can use several different RB schemes.
All of these approaches have in common that the RB spaces for the primal and dual spaces are chosen differently which makes the resulting reduced optimality system non-conforming.

As a naive and standard approach, in \cref{sec:TRRB_standard_approach}, we discussed to simply perform the respective Galerkin projections on the primal and dual systems with the respective RB spaces and to compute the objective functional and its derivatives simply by replacing the FOM functions by their ROM counterparts.
However, as shown, this results in an inexact ROM.

To remedy this, in \cref{sec:TRRB_ncd_approach}, we presented the non-conforming dual corrected approach.
We further presented how the exact gradient and Hessian information can efficiently be computed.

As an alternative ROM, in \cref{sec:TRRB_sensitivity_approach}, we further suggested adding more FOM information to the respective spaces to resolve the inexactness of the standard approach.
In \cref{sec:TRRB_pg_approach}, we additionally discussed to follow a Petrov-Galerkin ansatz where the primal and dual equations are reduced with the opposite RB space.
This approach showed advantages in the exact computation of the derivative information.

As the next step towards an efficient algorithm to solve the underlying PDE-constrained optimization problem, we require a rigorous error control of all discussed reduced methods, which is devised in the next section.

\section{A posteriori error analysis for the reduced models} 
\label{sec:TRRB_a_post_error_estimates}
In order to control the respective models' accuracy, we are now concerned about the a posteriori error analysis for the presented ROMs.
The resulting estimates for the objective functional is significant for the error-aware TR method, whereas an estimate for the control variable can be used as a post-processing of the algorithm; cf. \cref{sec:TR_RB_algorithm}.
If used inside the TR algorithm, the error estimates must be efficiently computable.
Recall that \cref{asmpt:parameter_separable} is crucial for offline-online efficient models as it allows to precompute most of the required FOM terms; cf. \cref{sec:background_off_on_RB}.

As before, for any functional $l \in V_h'$ or bilinear form $a: V_h \times V_h \to \R$, we denote their dual or operator norms $\|l\|$ and $\|a\|$, given by the continuity constants $\cont{l}$ and $\cont{a}$, respectively.
Moreover, for $\mu \in \Params$, we denote the coercivity constant of $\bformd$ w.r.t.~the $V_h$-norm by $\alpha_{\bformd} > 0$.

\subsection{The reduced primal and dual solutions}
\label{sec:estimation_for_primal_and_dual}

We start with the residual-based a posteriori error estimation for the primal variable.
Despite a slightly adapted notation, the estimate is completely equivalent to the standard result that we derived in \cref{prop:standard_RB_estimator}.
\begin{proposition}[Upper bound on the primal model reduction error] \label{prop:primal_rom_error}
	For $\mu \in \Params$, let $u_{h, \mu} \in V_h$ be the solution of~\eqref{eq:state_h} and $u_{\red, \mu} \in V_\red^\pr$ the solution of~\eqref{eq:state_red}.
	Then it holds
	\begin{align}
	\|u_{h, \mu} - u_{\red, \mu}\| \leq \Delta_\pr(\mu)  := \alpha_{\bformd}^{-1}\, \|\resd^\pr(u_{\red, \mu})\|_{V'}.
	\notag
	\end{align}
\end{proposition}
For the reduced dual problem, a similar idea can be used to derive the following estimation, accounting for the fact that $p_{\red, \mu}$ is not a Galerkin projection of $p_{h, \mu}$.
The stated proof can also be found in~\cite[Lemma 3]{QGVW2017}.
\begin{proposition}[Upper bound on the dual model reduction error]
	\label{prop:dual_rom_error}
	For $\mu \in \Params$, let $p_{h, \mu} \in V_h$ be the solution of~\eqref{eq:dual_solution_h} and $p_{\red, \mu} \in V_\red^\du$ the solution of~\eqref{eq:dual_solution_red}.
	Then it holds
	\begin{align}
	\|p_{h, \mu} - p_{\red, \mu}\| &\leq \Delta_\du(\mu) := \alpha_{\bformd}^{-1}\Big(2 \cont{\kformd}\;\Delta_\pr(\mu) + \|\resd^\du(u_{\red, \mu}, p_{\red, \mu})\|_{V'}\Big).
	\notag
	\end{align}
\end{proposition}
\begin{proof}
	Using the shorthand $e_{h, \mu}^\du := p_{h, \mu} - p_{\red, \mu}$, we have
	\begin{align*}
		\alpha_{\bformd}\, \|e_{h, \mu}^\du\|^2 &\leq a_\mu(e_{h, \mu}^\du, e_{h, \mu}^\du) = \underbrace{a_\mu(e_{h, \mu}^\du, p_{h, \mu})}_{=\partial_u \J(u_{h, \mu}, \mu)[e_{h, \mu}^\du]}
		- \,\, a_\mu(e_{h, \mu}^\du, p_{\red, \mu}) \\
		&= \underbrace{\partial_u \J(u_{h, \mu}, \mu)[e_{h, \mu}^\du]- \partial_u \J(u_{\red, \mu}, \mu)[e_{h, \mu}^\du]}_{= 2 k_\mu(e_{h, \mu}^\du, e_{h, \mu}^\pr)} 
		+ \big(\underbrace{\partial_u \J(u_{\red, \mu}, \mu)[e_{h, \mu}^\du] - a_\mu(e_{h, \mu}^\du, p_{\red, \mu})}_{ = r_{h, \mu}^\du(u_{\red, \mu}, p_{\red, \mu})[e_{h, \mu}^\du]}\big)\\
		&\leq 2 \|k_\mu\|\; \|e_{h, \mu}^\du\|\; \|e_{h, \mu}^\pr\| + \|r_{h, \mu}^\du(u_{\red, \mu}, p_{\red, \mu})\|\; \|e_{h, \mu}^\du\|,
	\end{align*}
	using the coercivity of $a_\mu$ in the first inequality, the definition of $e_{h, \mu}^\du$ in the first equality, the fact that $p_{h, \mu}$ solves~\eqref{eq:dual_solution_h} in the second equality, adding zero in the third equality, using the definition of $\J$ and \cref{rem:gateaux_wrt_V} and the definition of the discrete dual residual~\eqref{eq:dual_residual} in the intermediate equations, and the continuity of $k_\mu$ and the residual in the last inequality.
	We obtain the desired result utilizing Proposition \ref{prop:primal_rom_error}.
\end{proof}

\subsection{The reduced sensitivities}
\label{sec:RB_estimates_for_sens}
We also derive error estimates for the reduction error of the reduced sensitivities from~\eqref{eq:true_primal_reduced_sensitivity} and~\eqref{eq:true_dual_reduced_sensitivity} since we need them for the error estimation of the derivatives of the reduced objective functional.
For $v_h \in V_h$, a direction $e_i \in \R^P$, and quadratic objective functional~\eqref{eq:quadratic_J}, the corresponding FOM residuals of the above-mentioned equations are given by
\begin{align}
	\resd^{\pr,d_{\mu_i}}(u_{h, \mu}, d_{\mu_i} u_{h, \mu})[v_h] &:= \partial_{\mu_i} \resd^\pr(u_{h, \mu})[v_h] - \bformd(d_{\mu_i} u_{h, \mu}, v_h),
	\label{sens_res_pr}\\
	\resd^{\du,d_{\mu_i}}(u_{h, \mu}, p_{h, \mu},\mkern-2mu d_{\mu_i} u_{h, \mu},\mkern-2mu d_{\mu_i} p_{h, \mu})[v_h] 
	&:=  \partial_{\mu_i} \resd^\du(u_{h, \mu}, p_{h, \mu})[v_h]\mkern-3mu +\mkern-3mu 2\kformd(v_h, \mkern-2mu d_{\mu_i} u_{h, \mu})\mkern-3mu - \mkern-3mu \bformd(v_h,\mkern-2mu d_{\mu_i} p_{h, \mu}).
	\label{sens_res_du}
\end{align}
We note that the residuals for the approximated sensitivities in~\eqref{eq:red_sens_pr} and~\eqref{eq:red_sens_du} can be defined analogously. 
\begin{proposition}[Residual based upper bound on the model reduction error of the sensitivity of the primal solution map]
	\label{prop:primal_solution_dmui_error}
	For $\mu \in \Params$ and $1 \leq i \leq P$, let $d_{\mu_i}u_{h, \mu} \in V_h$ be the solution of the discrete version of~\eqref{eq:primal_sensitivities} and $d_{\mu_i}u_{\red, \mu} \in V_\red^{\pr,d_{\mu_i}}$ be the solution of~\eqref{eq:true_primal_reduced_sensitivity}.
	We then have
	\begin{align}
		\|d_{\mu_i}u_{h, \mu} - d_{\mu_i}u_{\red, \mu}\| &\leq \Delta_{d_{\mu_i}\pr}(\mu) := \alpha_{\bformd}^{-1}\Big(\cont{\partial_{\mu_i} \bformd} \Delta_\pr(\mu) + \|\resd^{\pr,d_{\mu_i}}(u_{\red, \mu}, d_{\mu_i}u_{\red, \mu})\| \Big).
		\notag
	\end{align}
\end{proposition}
\begin{proof}
	Using the shorthand $d_{\mu_i} e_{h, \mu}^\pr := d_{\mu_i}u_{h, \mu} - d_{\mu_i}u_{\red, \mu}$, we obtain
	\begin{align*}
		\alpha_{\bformd}\, &\|d_{\mu_i} e_{h, \mu}^\pr\|^2 \leq \bformd(d_{\mu_i} e_{h, \mu}^\pr, d_{\mu_i} e_{h, \mu}^\pr) \\
		&= \bformd(d_{\mu_i}u_{h, \mu}, d_{\mu_i} e_{h, \mu}^\pr) -\, \bformd(d_{\mu_i}u_{\red, \mu}, d_{\mu_i} e_{h, \mu}^\pr) \\
		&= \partial_{\mu_i} \resd^\pr(u_{h, \mu})[d_{\mu_i} e_{h, \mu}^\pr] - \bformd(d_{\mu_i}u_{\red, \mu}, d_{\mu_i} e_{h, \mu}^\pr) \\
		&= \partial_{\mu_i} \resd^\pr(u_{h, \mu})[d_{\mu_i} e_{h, \mu}^\pr] - \partial_{\mu_i} \resd^\pr(u_{\red, \mu})[d_{\mu_i} e_{h, \mu}^\pr] + \partial_{\mu_i} \resd^\pr(u_{\red, \mu})[d_{\mu_i} e_{h, \mu}^\pr] - \bformd(d_{\mu_i}u_{\red, \mu}, d_{\mu_i} e_{h, \mu}^\pr) \\	
		&= - \partial_{\mu_i} \bformd(e_{h, \mu}^\pr, d_{\mu_i} e_{h, \mu}^\pr) + \resd^{\pr,d_{\mu_i}}(u_{\red, \mu}, d_{\mu_i}u_{\red, \mu})[d_{\mu_i} e_{\red, \mu}^\pr] \\
		&\leq \cont{\partial_{\mu_i} \bformd} \;\|e_{h, \mu}^\pr\| \;\|d_{\mu_i} e_{h, \mu}^\pr\|  + \|\resd^{\pr,d_{\mu_i}}(u_{\red, \mu}, d_{\mu_i}u_{\red, \mu})\|\; \|d_{\mu_i} e_{h, \mu}^\pr\|
	\end{align*}
	using the coercivity of $\bformd$ in the first inequality, the definition of $d_{\mu_i} e_{h, \mu}^\pr$ in the first equality, Proposition~\ref{prop:solution_dmu_eta} applied to $u_{h, \mu}$ in the second equality, the definition of the discrete sensitivity primal residual~\eqref{sens_res_pr} in the fourth equality and the continuity of $\partial_{\mu_i} \bformd$ in the last inequality.
\end{proof}
We emphasize that the same result can be shown for the approximated sensitivity $\dred{\mu_i} u_{\red, \mu}$ by replacing $d_{\mu_i} u_{\red, \mu}$ and using the equation~\eqref{eq:red_sens_pr} instead of~\eqref{eq:true_primal_reduced_sensitivity}.
We call the resulting error estimator $\Delta_{\dred{\mu_i}\pr}(\mu)$.
Next, the dual sensitivity model reduction error is deduced.

\begin{proposition}[Residual based upper bound on the model reduction error of the sensitivity of the dual solution map]
	\label{prop:dual_solution_dmui_error}
	For $\mu \in \Params$ and $1 \leq i \leq P$, let $d_{\mu_i}p_{h, \mu} \in V_h$ be the solution of the discrete version of~\eqref{eq:dual_sens} and $d_{\mu_i}p_{\red, \mu} \in V_\red^{\pr,d_{\mu_i}}$ be the solution of~\eqref{eq:true_dual_reduced_sensitivity}.
	We then obtain
	\begin{align*}
		\|d_{\mu_i}p_{h, \mu} - &d_{\mu_i}p_{\red, \mu}\| \leq \Delta_{ d_{\mu_i}\du}(\mu)\hspace*{5cm}\text{with}\\
		\Delta_{ d_{\mu_i}\du}(\mu) := \alpha_{\bformd}^{-1}\Big(
		&2 \cont{\partial_{\mu_i} \kformd} \; \Delta_\pr(\mu) +  \cont{\partial_{\mu_i} \bformd} \; \Delta_\du(\mu) \\ & \qquad+ 2 \cont{\kformd}  \; \Delta_{d_{\mu_i}\pr}(\mu) 
		+ \| \resd^{\du,d_{\mu_i}}(u_{\red, \mu}, p_{\red, \mu}, d_{\mu_i}u_{\red, \mu}, d_{\mu_i}p_{\red, \mu}) \| \Big).
	\end{align*}
\end{proposition}
\begin{proof}
	Using the shorthand $d_{\mu_i} e_{h, \mu}^\du := d_{\mu_i}p_{h, \mu} - d_{\mu_i}p_{\red, \mu}$ and $e_{h, \mu}^\du := p_{h, \mu} - p_{\red, \mu}$, we obtain
	\begin{align*}
		\alpha_{\bformd}\;&\|d_{\mu_i} e_{h, \mu}^\du\|^2 \leq \bformd(d_{\mu_i} e_{h, \mu}^\du, d_{\mu_i} e_{h, \mu}^\du) \\
		&= \underbrace{\bformd(d_{\mu_i} e_{h, \mu}^\du, d_{\mu_i}p_{h, \mu})}_{= \partial_{\mu_i} \resd^\du(u_{h, \mu}, p_{h, \mu})[d_{\mu_i} e_{h, \mu}^\du] 
			+ 2\kformd(d_{\mu_i} e_{h, \mu}^\du, d_{\mu_i}u_{h, \mu})} - \bformd(d_{\mu_i} e_{h, \mu}^\du, d_{\mu_i}p_{\red, \mu})
		\\
		&=\partial_{\mu_i} \resd^\du(u_{h, \mu}, p_{h, \mu})[d_{\mu_i} e_{h, \mu}^\du] + 2\kformd(d_{\mu_i} e_{h, \mu}^\du, d_{\mu_i}u_{h, \mu})
		- \partial_{\mu_i} \resd^\du(u_{\red, \mu}, p_{\red, \mu})[d_{\mu_i} e_{h, \mu}^\du] \\
		&\qquad - 2\kformd(d_{\mu_i} e_{h, \mu}^\du, d_{\mu_i}u_{\red, \mu})
		+ \partial_{\mu_i} \resd^\du(u_{\red, \mu}, p_{\red, \mu})[d_{\mu_i} e_{h, \mu}^\du] \\
		&\qquad + 2\kformd(d_{\mu_i} e_{h, \mu}^\du, d_{\mu_i}u_{\red, \mu}) 
		- \bformd(d_{\mu_i} e_{h, \mu}^\du, d_{\mu_i}p_{\red, \mu}) \\
		&= \partial_{\mu_i} \jformd(d_{\mu_i} e_{h, \mu}^\du) + 2 \partial_{\mu_i}\kformd(d_{\mu_i} e_{h, \mu}^\du, u_{h, \mu}) - \partial_{\mu_i}\bformd(d_{\mu_i} e_{h, \mu}^\du, p_h)
		- \partial_{\mu_i}\jformd(d_{\mu_i} e_{h, \mu}^\du) \\
		&\qquad + 2 \partial_{\mu_i}\kformd(d_{\mu_i} e_{h, \mu}^\du, u_{\red, \mu})
		-\partial_{\mu_i}\bformd(d_{\mu_i} e_{h, \mu}^\du, p_\red) +\! 2\kformd(d_{\mu_i} e_{h, \mu}^\du, d_{\mu_i}u_{h, \mu}) 
		\\ &\qquad
		- 2\kformd(d_{\mu_i} e_{h, \mu}^\du, d_{\mu_i}u_{\red, \mu}) + \resd^{\du,d_{\mu_i}}(u_{\red, \mu}, p_{\red, \mu}, d_{\mu_i}u_{\red, \mu}, d_{\mu_i}p_{\red, \mu})[d_{\mu_i} e_{h, \mu}^\du]	\\
		&= 2 \partial_{\mu_i}\kformd(d_{\mu_i} e_{h, \mu}^\du, e_{h, \mu}^\pr) - \partial_{\mu_i}\bformd(d_{\mu_i} e_{h, \mu}^\du, e_{h, \mu}^\du)
		+ 2\kformd(d_{\mu_i} e_{h, \mu}^\du, d_{\mu_i} e_{h, \mu}^\pr) \\
		&\qquad
		+ \resd^{\du,d_{\mu_i}}(u_{\red, \mu}, p_{\red, \mu}, d_{\mu_i}u_{\red, \mu}, d_{\mu_i}p_{\red, \mu})[d_{\mu_i} e_{h, \mu}^\du]\\
		&\leq \big(2 \cont{\partial_{\mu_i} \kformd} \; \|e_{h, \mu}^\pr\| +  \cont{\partial_{\mu_i} \bformd} \; \|e_{h, \mu}^\du\|  \big)\|d_{\mu_i} e_{h, \mu}^\du\| \\
		&\qquad + 2 \cont{\kformd}  \; \|d_{\mu_i} e_{h, \mu}^\pr\| \; \|d_{\mu_i} e_{h, \mu}^\du\|
		+ \| \resd^{\du,d_{\mu_i}}(u_{\red, \mu}, p_{\red, \mu}, d_{\mu_i}u_{\red, \mu}, d_{\mu_i}p_{\red, \mu}) \| \; \|d_{\mu_i} e_{h, \mu}^\du\|,
	\end{align*}
	using the coercivity of $\bformd$ in the first inequality, the definition of $d_{\mu_i} e_{h, \mu}^\du$ in the first equality, \cref{prop:dual_solution_dmu_eta} applied to $p_{h, \mu}$ in the second equality, the definition of the dual residual in~\eqref{eq:dual_residual} in the third equality and continuity of all parts in the last inequality.
\end{proof}
Again, the same result holds for $\dred{\mu_i} p_{\red, \mu}$ if we replace $d_{\mu_i} p_{\red, \mu}$ and use~\eqref{eq:red_sens_du} instead of~\eqref{eq:true_dual_reduced_sensitivity}.
The resulting error estimator is called $\Delta_{\dred{\mu_i}\du}(\mu)$.

\subsection{The reduced functional}
\label{sec:TRRB_estimation_for_J}

In the following, we state the error estimation result of the standard reduced objective functional $\Jnoncor_\red$ in~\eqref{eq:Jhat_red}, which was presented in~\cite[Theorem~4]{QGVW2017}.
We furthermore show an improved version by using, in contrast to~\cite{QGVW2017}, the NCD-corrected reduced functional $\Jhat_\red$, which results in an optimal higher-order a posteriori upper bound without lower order terms.
It is important to emphasize that \cref{prop:Jhat_error} is essential for the error-aware TR-RB method in
\cref{sec:TR_RB_algorithm} since the trust-region is characterized by the accuracy of the reduced objective functional.

\begin{proposition}[Upper bound on the model reduction error of the reduced output]
	\label{prop:Jhat_error} 
	\hfill
	\begin{enumerate}
		\item [\emph{(i)}] With the notation from above, we have for the standard RB reduced cost functional
		\begin{equation*}
		|\Jhat_h(\mu) - \Jnoncor_\red(\mu)| \leq \Delta_{\Jnoncor_\red}(\mu) :=  \Delta_\pr(\mu) \|\resd^\du(u_{\red, \mu}, p_{\red,\mu})\| + \Delta_\pr(\mu)^2 \cont{\kformd} + \big|r_{\mu}^\pr(u_{\red, \mu})[p_{\red,\mu}]\big|.
		\end{equation*}	
		\item [\emph{(ii)}] Furthermore, we have for the NCD-corrected RB reduced cost functional (or equivalently for the Lagrangian for any $p\in V_h$)
		\begin{align*}
		|\Jhat_h(\mu) - \cJhatn(\mu)| &= |\mathcal{L}(u_{h,\mu},\mu,p)- \mathcal{L}(u_{\red,\mu},\mu,p)| \\ &\leq \Delta_{\cJhatn}(\mu) :=  \Delta_\pr(\mu) \|\resd^\du(u_{\red, \mu}, p_{\red,\mu})\| + \Delta_\pr(\mu)^2 \cont{\kformd}.
		\end{align*}
	\end{enumerate}
\end{proposition}
\begin{proof}
	We prove $(i)$ analogously to~\cite[Theorem~4]{QGVW2017}.
	Using the shorthand $e_{h, \mu}^\pr := u_{h, \mu} - u_{\red, \mu}$ and $a_\mu(e_{h, \mu}^\pr,p_{\red,\mu}) = r_{\mu}^\pr(u_{\red, \mu})[p_{\red,\mu}]$, we have 
	\begin{align*}
		|\Jhat_h(\mu)& - \Jnoncor_\red(\mu)| 
		= |\J(u_{h, \mu}, \mu) - \J(u_{\red, \mu}, \mu)|
		\\
		&= |j_\mu(e_{h, \mu}^\pr) + \kformd(u_{h, \mu}, u_{h, \mu}) - \kformd(u_{\red, \mu}, u_{\red, \mu}) \underbrace{- \,a_\mu(e_{h, \mu}^\pr,p_{\red,\mu}) + r_{\mu}^\pr(u_{\red, \mu})[p_{\red,\mu}]}_{=0}| \\
		&= |r_{\mu}^\du(u_{\red, \mu}, p_{\red,\mu})[e_{h, \mu}^\pr] - 2\kformd(u_{\red, \mu},e_{h, \mu}^\pr)
		+ \kformd(u_{h, \mu}, u_{h, \mu}) - \kformd(u_{\red, \mu}, u_{\red, \mu}) + r_{\mu}^\pr(u_{\red, \mu})[p_{\red,\mu}]| \\
		&= |r_{\mu}^\du(u_{\red, \mu}, p_{\red,\mu})[e_{h, \mu}^\pr] + \kformd(e_{h, \mu}^\pr, e_{h, \mu}^\pr) + r_{\mu}^\pr(u_{\red, \mu})[p_{\red,\mu}]| \\
		&\leq \|r_{\mu}^\du(u_{\red, \mu}, p_{\red,\mu})\|\; \|e_{h, \mu}^\pr\| +\cont{\kformd} \; \|e_{h, \mu}^\pr\|^2 + \big|r_{\mu}^\pr(u_{\red, \mu})[p_{\red,\mu}]\big|,
	\end{align*}
	where we used the definition of the dual residual in the second equality and Cauchy-Schwarz for the inequality. 
	The assertion follows by using Proposition~\ref{prop:primal_rom_error}. 
	
	Since $\cJhatn(\mu) = \Jnoncor_\red(\mu) + r_\mu^\pr(u_{\red,\mu})[p_{\red,\mu}]$, the proof for  $(ii)$ follows along the same lines without the necessity to add $0$ in the second equality.
\end{proof}

It can be seen in the proof of \cref{prop:Jhat_error} that the NCD-correction is a natural enhancement for the standard approach since the residual term actually needs to be added for proving~$(i)$.
Furthermore, the enhanced accuracy is easily visible: The estimator $\Delta_{\cJhatn}$ only consists of terms that converge twice towards zero (if $V^\pr_\red$ and $V^\pr_\red$ grow) and the uncorrected version $\Delta_{\Jnoncor_\red}$ contains an additional term that decreases only of order one; cf. \cref{sec:TRRB_estimator_study}.

For the convergence study of the TR-RB algorithm, we state the following remark, which is a simple property of residual-based error estimation and inherits from the fact that the Riesz-representative of the residual is continuous in $\Params$.
\begin{remark}\label{continuity_of_estimator}
	The estimator $\Delta_{\cJhatn}(\mu)$ is continuous w.r.t.~$\mu \in \Params$.
\end{remark}

While \cref{prop:Jhat_error} already constitutes an enhanced error behavior for the NCD-corrected functional, the same can be proven for the reduced gradient information.

\subsection{The gradient of the reduced functional}
\label{sec:TRRB_RB_estimates_for_grad}

We continue with an error analysis of the reduced derivative information.
For the inexact gradient, we follow the proof from~\cite[Theorem~4]{QGVW2017}.
In the following proposition, we also show a first error estimate for the NCD-corrected gradient, which will be enhanced in \cref{prop:grad_Jhat_error_sens}.
\begin{proposition}[Upper bound on the model reduction error of the gradient of the reduced output]
	\label{prop:grad_Jhat_error}
	\hfill
	\begin{enumerate}
		\item [\emph{(i)}] For the inexact gradient $\noncorgrad_\mu \Jnoncor_\red(\mu)$ from the standard-RB approach~\eqref{naive:red_grad}, we have  
		\begin{align*}
		&\big\|\nabla_\mu \Jhat_h(\mu) - \noncorgrad_\mu \Jnoncor_\red(\mu)\big\|_2 \leq \Delta_{\noncorgrad \Jnoncor_\red}(\mu) = \big\|\underline{\Delta_{\noncorgrad \Jnoncor_\red}(\mu)}\big\|_2 \quad\quad\text{with} \\
		\big(\underline{\Delta_{\noncorgrad \Jnoncor_\red}(\mu)}\big)_i  
		&:= 2\Delta_\pr(\mu) \|u_{\red, \mu}\|\; \cont{\partial_{\mu_i} \kformd}
		+ \Delta_\pr(\mu)\big(\cont{\partial_{\mu_i} \jformd} + \cont{\partial_{\mu_i} \bformd}\; \|p_{\red, \mu}\|\big)\\
		&\,\,+ \Delta_\du(\mu)\big(\cont{\partial_{\mu_i} \lformd} + \cont{\partial_{\mu_i} \bformd}\; \|u_{\red, \mu}\|\big)
		+ \Delta_\pr(\mu)\; \Delta_\du(\mu)\; \cont{\partial_{\mu_i} \bformd}
		+ (\Delta_\pr)^2(\mu)\; \cont{\partial_{\mu_i} \kformd}.
		\end{align*}
		\item [\emph{(ii)}] For the gradient $\nabla_{\mu} \cJhatn(\mu)$ of the NCD-corrected reduced functional, expressed by the adjoint approach from \cref{prop:true_corrected_reduced_gradient_adj}, we have 
		\begin{align*}
		&\big\|\nabla_\mu \Jhat_h(\mu) - \nabla_{\mu} \cJhatn(\mu)\big\|_2 \leq \Delta^*_{\nabla \Jhat_\red}(\mu) 
		= \big\|\underline{\Delta^*_{\nabla \Jhat_\red}(\mu)}\big\|_2 \quad\quad\text{with} \\
		\big(\underline{\Delta^*_{\nabla \Jhat_\red}(\mu)}\big)_i &:=
		2\Delta_\pr(\mu) \|u_{\red, \mu}\|\; \cont{\partial_{\mu_i} \kformd}
		+ \Delta_\pr(\mu)\big(\cont{\partial_{\mu_i} \jformd} + \cont{\partial_{\mu_i} \bformd}\; \|p_{\red, \mu}\|\big)\\
		&+ \Delta_\du(\mu)\big(\cont{\partial_{\mu_i} \lformd} + \cont{\partial_{\mu_i} \bformd}\; \|u_{\red, \mu}\|\big)
		+ \Delta_\pr(\mu)\; \Delta_\du(\mu)\; \cont{\partial_{\mu_i} \bformd}
		+ (\Delta_\pr)^2(\mu)\; \cont{\partial_{\mu_i} \kformd} \\
		&+ (\cont{\partial_{\mu_i}l_{\mu}} + \cont{\partial_{\mu_i}a_{\mu}} \| u_{\red,\mu} \|) \alpha_{\bformd}^{-1} \big( \|r_\mu^\du(u_{\red,\mu},p_{\red,\mu})\| 
		+ 2\cont{k_\mu} \alpha_{\bformd}^{-1} \|r_\mu^\pr(u_{\red,\mu})\| \big) \\
		&+ \alpha_{\bformd}^{-1} \|r_\mu^\pr(u_{\red,\mu})\| \big( \cont{\partial_{\mu_i}j}+ 2\cont{\partial_{\mu_i}k} \|u_{\red,\mu}\|
		+ \cont{\partial_{\mu_i}a} \|p_{\red,\mu}\| \big).
		\end{align*}
	\end{enumerate}
\end{proposition}
\begin{proof}
	(i) For $\Delta_{\noncorgrad \Jnoncor_\red}(\mu)$, we have 
	\begin{align*}
	\big(\nabla_\mu \Jhat_h(\mu) \!-\! \noncorgrad_{\mu} \cJhatn(\mu)\big)_i
	&=\partial_{\mu_i}\J(u_{h, \mu}, \mu) - \partial_{\mu_i}\J(u_{\red, \mu}, \mu) 
	+ \partial_{\mu_i}r_\mu^\pr(u_{h,\mu})[p_{h,\mu}] 
	-\partial_{\mu_i}r_\mu^\pr(u_{\red,\mu})[p_{\red,\mu}].
	\end{align*}
	Regarding the first contribution, we obtain, using $\|u_{h, \mu}\| \leq \|e_{h, \mu}^\pr\| + \|u_{\red, \mu}\|$, the estimate
	\begin{align*}
	\big|\partial_{\mu_i}\J(u_{h, \mu}, \mu) - \partial_{\mu_i}\J(u_{\red, \mu}, \mu)\big| &=
	|\partial_{\mu_i}j_\mu(e_{h, \mu}^\pr) 
	+ \partial_{\mu_i}k_{\mu}(e_{h, \mu}^\pr, u_{\red, \mu}) 
	+ \partial_{\mu_i}k_{\mu}(u_{h, \mu}, e_{h, \mu}^\pr)|
	\\
	&\leq \Delta_\pr(\mu)\Big( \cont{\partial_{\mu_i} \jformd} + \cont{\partial_{\mu_i} \kformd}\big(2 \|u_{\red, \mu}\| + \Delta_\pr(\mu)\big)\Big).
	\end{align*} 
	For the other contributions, we refer to~\cite[Theorem~5]{QGVW2017}.
	
	(ii) For the adjoint estimator $\Delta^*_{\nabla_\mu \Jhat_\red}$, we have 
	\begin{align*}
	\big(\nabla_\mu \Jhat_h(\mu) - \nabla_{\mu} \cJhatn(\mu)\big)_i
	&=\partial_{\mu_i}\J(u_{h, \mu}, \mu) - \partial_{\mu_i}\J(u_{\red, \mu}, \mu) 
	+\partial_{\mu_i}r_\mu^\pr(u_{h,\mu})[p_{h,\mu}] 
	- \partial_{\mu_i}r_\mu^\pr(u_{\red,\mu})[p_{\red,\mu}] \\
	& \qquad \qquad - \partial_{\mu_i} r_\mu^\pr(u_{\red,\mu})[w_{\red,\mu}] 
	+ \partial_{\mu_i} r_\mu^\du(u_{\red,\mu},p_{\red,\mu})[z_{\red,\mu}]. 
	\end{align*}
	The first line is equal to the estimator $\Delta_{\noncorgrad \Jnoncor_\red}(\mu)$, the first term of the second line can be estimated by
	\begin{align*}
	\partial_{\mu_i} r_\mu^\pr(u_{\red,\mu})[w_{\red,\mu}] 
	&\leq \cont{\partial_{\mu_i}l_{\mu}}\|w_{\red, \mu}\| + \cont{\partial_{\mu_i}a_{\mu}} \| u_{\red,\mu} \| \|w_{\red, \mu}\|. 
	\end{align*}
	The second term can analogously be estimated by
	\begin{align*}
	\partial_{\mu_i} r_\mu^\du(u_{\red,\mu},p_{\red,\mu})[z_{\red, \mu}] \leq 
	\cont{\partial_{\mu_i}j} \|z_{\red, \mu}\| + 2\cont{\partial_{\mu_i}k} \| z_{\red, \mu} \| \|u_{\red,\mu}\|
	+ \cont{\partial_{\mu_i}a} \|z_{\red, \mu}\| \|p_{\red,\mu}\|.
	\end{align*}
	We also have 
	\begin{align*}
	\alpha_{\bformd} \|w_{\red, \mu}\|^2 &\leq \bformd(w_{\red, \mu},w_{\red, \mu}) = r_\mu^\du(u_{\red,\mu},p_{\red,\mu})[w_{\red, \mu}]-2k_\mu(z_{\red, \mu},w_{\red, \mu}) \\
	&\leq \|r_\mu^\du(u_{\red,\mu},p_{\red,\mu})\| \|w_{\red, \mu}\| + 2\cont{k_\mu} \|z_{\red, \mu}\| \|w_{\red, \mu}\|,
	\end{align*}
	which gives
	\begin{equation*}
	\|w_{\red, \mu}\| \leq \alpha_{\bformd}^{-1} \left( \|r_\mu^\du(u_{\red,\mu},p_{\red,\mu})\| + 2\cont{k_\mu} \|z_{\red, \mu}\| \right).
	\end{equation*}
	For $z_\mu$, we estimate
	\begin{align*}
	\alpha_{\bformd} \|z_{\red, \mu}\|^2 &\leq \bformd(z_{\red, \mu},z_{\red, \mu}) = -r_\mu^\pr(u_{\red,\mu})[z_\mu]\leq \|r_\mu^\pr(u_{\red,\mu})\| \|z_{\red, \mu}\|.
	\end{align*}
	Summation attains the assertion.
\end{proof}
We note that the estimator $\Delta^*_{\nabla \Jhat_\red}(\mu)$ for the NCD-corrected gradient does not suggest a better approximation of the FOM gradient since more terms are added to the standard estimate.
In order to derive a more effective estimator with higher-order convergence, we take advantage of the alternative sensitivity-based expression of $\nabla_{\mu} \cJhatn$ in \cref{prop:true_corrected_reduced_gradient_sens}, although it was declared computationally inefficient.
Using the residual-based a posteriori error estimates for the sensitivities presented in \cref{sec:RB_estimates_for_sens}, we can state two a posteriori error bounds on the model reduction error of the true gradient $\nabla_{\mu} \cJhatn$ and the approximated gradient $\gradred{\mu} \cJhatn$ of the NCD-corrected functional.

\begin{proposition}[Upper bound on the model reduction error of the gradient of the reduced output -- sensitivity approach] 
	\label{prop:grad_Jhat_error_sens}
	\hfill
	\begin{enumerate}
		\item [\emph{(i)}] For the gradient $\nabla_{\mu} \cJhatn(\mu)$ of the NCD-corrected RB reduced functional, expressed with sensitivities according to  \cref{prop:true_corrected_reduced_gradient_sens}, we have
		\begin{align*}
		&\big\|\nabla_\mu \Jhat_h(\mu) - \nabla_{\mu} \cJhatn(\mu)\big\|_2 \leq \Delta_{\nabla\Jhat_\red}(\mu) 
		= \big\|\underline{\Delta_{\nabla\Jhat_\red}(\mu)}\big\|_2 \quad\quad\text{with} \\
		\big(\underline{\Delta_{\nabla\Jhat_\red}(\mu)}\big)_i &:= \cont{\partial_{\mu_i} \kformd} \, \left(\Delta_\pr(\mu)\right)^2 + \cont{\bformd} \, 
		\Delta_{d_{\mu_i}pr}(\mu) \,\Delta_\du(\mu) \\ & \qquad\qquad + \| \resd^{\du,d_{\mu_i}}(u_{\red, \mu}, p_{\red, \mu}, d_{\mu_i}u_{\red, \mu}, d_{\mu_i}p_{\red, \mu}) \| \,
		\Delta_\pr(\mu).
		\end{align*}
		\item [\emph{(ii)}] Furthermore, we have for the approximate gradient $\gradred{\mu} \cJhatn$ from \cref{prop:red_approx_grad}
		\begin{align*}
		&\big\|\nabla_\mu \Jhat_h(\mu) - \gradred{\mu} \cJhatn(\mu)\big\|_2 \leq \Delta_{\widetilde\nabla\Jhat_\red}(\mu) 
		= \big\|\underline{\Delta_{\widetilde\nabla\Jhat_\red}(\mu)}\big\|_2 \quad\quad\text{with} \\
		\big(\underline{\Delta_{\widetilde\nabla\Jhat_\red}(\mu)}\big)_i &:= \cont{\partial_{\mu_i} \kformd} \, \left(\Delta_\pr(\mu)\right)^2 + \cont{\bformd} \, 
		\Delta_{\dred{\mu_i}\pr}(\mu) \,\Delta_\du(\mu) \\ & \qquad\qquad+ \| \resd^{\du,d_{\mu_i}}(u_{\red, \mu}, p_{\red, \mu}, \dred{\mu_i}u_{\red, \mu}, \dred{\mu_i}p_{\red, \mu}) \| \,
		\Delta_\pr(\mu).
		\end{align*}
	\end{enumerate}
\end{proposition}
\begin{proof}
	(i) To prove the first assertion, we use $r_\mu^\pr(u_{h, \mu})[d_{\mu_i}p_{\red, \mu}] = 0$ and $r_\mu^\du(u_{h, \mu},p_{h, \mu})[d_{\mu_i} u_{\red, \mu}]=0$ to obtain
	\begin{align*}
	\big(\nabla_\mu \Jhat_h(\mu) &- \nabla_\mu \cJhatn(\mu)\big)_i \\ 
	&=\partial_{\mu_i}\J(u_{h, \mu}, \mu) - \partial_{\mu_i}\J(u_{\red, \mu}, \mu) + \partial_{\mu_i}r_\mu^\pr(u_{h,\mu})[p_{h,\mu}] - \partial_{\mu_i}r_\mu^\pr(u_{\red,\mu})[p_{\red,\mu}]  \\
	&\qquad - r_\mu^\pr(u_{\red, \mu})[d_{\mu_i}p_{\red, \mu}] - r_\mu^\du(u_{\red, \mu},p_{\red, \mu})[d_{\mu_i} u_{\red, \mu}] \\
	&= \partial_{\mu_i} \jformd(e_{h, \mu}^\pr) +\partial_{\mu_i} \kformd(u_{h, \mu}, u_{h, \mu}) -\partial_{\mu_i} \kformd(u_{\red, \mu},u_{\red, \mu})  
	+ \partial_{\mu_i}r_\mu^\pr(u_{h,\mu})[p_{h,\mu}] \\  &\qquad - \partial_{\mu_i}r_\mu^\pr(u_{\red,\mu})[p_{\red,\mu}]
    + \underset{=(*)}{\underbrace{r_\mu^\pr(e_{h, \mu}^\pr)[d_{\mu_i}p_{\red, \mu}]}} 
	+ \underset{=(**)}{\underbrace{r_\mu^\du(e_{h, \mu}^\pr,e_{h, \mu}^\du)[d_{\mu_i} u_{\red, \mu}]}}.
	\end{align*}
	For the last two residual terms we have
	\begin{align*}
	(*) &= \lformd(d_{\mu_i}p_{\red, \mu}) - \lformd(d_{\mu_i}p_{\red, \mu}) 
	- \bformd(e_{h, \mu}^\pr,d_{\mu_i}p_{\red, \mu}) \\ 
	&= - \bformd(e_{h, \mu}^\pr,d_{\mu_i}p_{\red, \mu}) 
	+ \partial_{\mu_i} r_\mu^\du(u_{\red, \mu},p_{\red, \mu})[e_{h, \mu}^\pr] 
	+ 2\kformd(d_{\mu_i}u_{\red, \mu}, e_{h, \mu}^\pr) \\ & \qquad
	- \partial_{\mu_i} r_\mu^\du(u_{\red, \mu},p_{\red, \mu})[e_{h, \mu}^\pr] 
	- 2\kformd(d_{\mu_i}u_{\red, \mu}, e_{h, \mu}^\pr) \\
	&= \resd^{\du,d_{\mu_i}}(u_{\red, \mu}, p_{\red, \mu}, d_{\mu_i}u_{\red, \mu}, d_{\mu_i}p_{\red, \mu})[e_{h, \mu}^\pr] 
	- \partial_{\mu_i} r_\mu^\du(u_{\red, \mu},p_{\red, \mu})[e_{h, \mu}^\pr] 
	- 2\kformd(d_{\mu_i}u_{\red, \mu}, e_{h, \mu}^\pr) 
	\end{align*}
	and
	\begin{align*}
	(**) &= \jformd(d_{\mu_i}u_{\red, \mu}) - \jformd(d_{\mu_i}u_{\red, \mu}) 
	+ 2\kformd(d_{\mu_i}u_{\red, \mu}, e_{h, \mu}^\pr) - \bformd(d_{\mu_i}u_{\red, \mu}, e_{h, \mu}^\du). 
	\end{align*}
	Thus, by summing both terms we have
	\begin{align*}
	(*) + (**)= 
	\resd^{\du,d_{\mu_i}}(u_{\red, \mu}, p_{\red, \mu},d_{\mu_i}u_{\red, \mu}, d_{\mu_i}p_{\red, \mu})[e_{h, \mu}^\pr] 
	- \underset{=(***)}{\underbrace{\partial_{\mu_i} r_\mu^\du(u_{\red, \mu},p_{\red, \mu})[e_{h, \mu}^\pr]}} - \bformd(d_{\mu_i}u_{\red, \mu}, e_{h, \mu}^\du) 
	\end{align*}
	and for $(***)$ it holds
	\begin{align*}
	(***) = \partial_{\mu_i}\jformd(e_{h, \mu}^\pr) + 2 \partial_{\mu_i}\kformd(e_{h, \mu}^\pr, u_{\red,\mu}) - \partial_{\mu_i}\bformd(e_{h, \mu}^\pr,p_{\red,\mu}).
	\end{align*}
	Combining $(*)$, $(**)$, and $(***)$ with the previous result, we obtain
	\begin{align*}
	\big(\nabla_\mu \Jhat_h(\mu) &- \nabla_\mu \cJhatn(\mu)\big)_i 
	= \partial_{\mu_i} \jformd(e_{h, \mu}^\pr) + \partial_{\mu_i}\kformd(u_{h, \mu}, u_{h, \mu}) - \partial_{\mu_i} \kformd(u_{\red, \mu},u_{\red, \mu})  \\
	&\qquad - \partial_{\mu_i}\jformd(e_{h, \mu}^\pr) - 2 \partial_{\mu_i}\kformd(e_{h, \mu}^\pr, u_{\red,\mu}) + \partial_{\mu_i}\bformd(e_{h, \mu}^\pr,p_{\red,\mu}) + \partial_{\mu_i}r_\mu^\pr(u_{h,\mu})[p_{h,\mu}]  \\
	&\qquad - \partial_{\mu_i}r_\mu^\pr(u_{\red,\mu})[p_{\red,\mu}] 
	+  \resd^{\du,d_{\mu_i}}(u_{\red, \mu}, p_{\red, \mu},d_{\mu_i}u_{\red, \mu}, d_{\mu_i}p_{\red, \mu})[e_{h, \mu}^\pr] 
	- \bformd(d_{\mu_i}u_{\red, \mu}, e_{h, \mu}^\du)  \\
	& = \partial_{\mu_i} \kformd(e_{h, \mu}^\pr, e_{h, \mu}^\pr) 
	+  \resd^{\du,d_{\mu_i}}(u_{\red, \mu}, p_{\red, \mu},d_{\mu_i}u_{\red, \mu}, d_{\mu_i}p_{\red, \mu})[e_{h, \mu}^\pr] \\
	& \qquad + \underset{=(****)}{\underbrace{\partial_{\mu_i}r_\mu^\pr(u_{h,\mu})[p_{h,\mu}] - \partial_{\mu_i}r_\mu^\pr(u_{\red,\mu})[p_{\red,\mu}] 
			+  \partial_{\mu_i}\bformd(e_{h, \mu}^\pr,p_{\red,\mu}) -  \bformd(d_{\mu_i}u_{\red, \mu}, e_{h, \mu}^\du)}}.
	\end{align*}
	Further, it holds
	\begin{align*}
	\partial_{\mu_i}r_\mu^\pr(u_{h,\mu})[p_{h,\mu}] 
	&-\partial_{\mu_i}r_\mu^\pr(u_{\red,\mu})[p_{\red,\mu}] 
	= \partial_{\mu_i} \lformd(e_{h, \mu}^\du) - 
	\partial_{\mu_i} \bformd(u_{h, \mu}, p_{h, \mu}) + \partial_{\mu_i} \bformd(u_{\red, \mu}, p_{\red, \mu}) \\
	& =  \bformd(d_{\mu_i}u_h, e_{h, \mu}^\du) +\partial_{\mu_i} \bformd(u_h, e_{h, \mu}^\du) 
	- \partial_{\mu_i} \bformd(u_{h, \mu}, p_{h, \mu}) + \partial_{\mu_i} \bformd(u_{\red, \mu}, p_{\red, \mu}),
	\end{align*}
	where we used the discretized version of~\eqref{eq:primal_sens} in the second equality.
	Inserting this into $(*)$ yields
	\begin{align*}
	(***\,*) &=  \bformd(d_{\mu_i}u_h, e_{h, \mu}^\du) -  \bformd(d_{\mu_i}u_{\red, \mu}, e_{h, \mu}^\du) \\ & \qquad +  
	\underset{=0}{\underbrace{\partial_{\mu_i} \bformd(u_h, e_{h, \mu}^\du) - \partial_{\mu_i} \bformd(u_{h, \mu}, p_{h, \mu}) 
			+ \partial_{\mu_i} \bformd(u_{\red, \mu}, p_{\red, \mu}) 
			+\partial_{\mu_i}\bformd(e_{h, \mu}^\pr,p_{\red,\mu}) }} \\ 
		&= \bformd(d_{\mu_i} e_{h, \mu}^\pr, e_{h, \mu}^\du).
	\end{align*}
	In total, we conclude 
	\begin{align*}
	\big(\nabla_\mu \Jhat_h(\mu) - \nabla_\mu \cJhatn(\mu)\big)_i 
	&= \partial_{\mu_i} \kformd(e_{h, \mu}^\pr, e_{h, \mu}^\pr) 
	+ \bformd(d_{\mu_i} e_{h, \mu}^\pr, e_{h, \mu}^\du) \\ & \qquad\qquad
	+  \resd^{\du,d_{\mu_i}}(u_{\red, \mu}, p_{\red, \mu},d_{\mu_i}u_{\red, \mu}, d_{\mu_i}p_{\red, \mu})[e_{h, \mu}^\pr],
	\end{align*}
	which proofs the assertion. \\
	(ii) The estimate follows analogously to (i), by replacing $d_{\mu_i} u_{\red, \mu}$ and $d_{\mu_i} p_{\red, \mu}$ by $\dred{\mu_i} u_{\red, \mu}$ and $\dred{\mu_i} p_{\red, \mu}$, respectively.
\end{proof}

Consequently, $\Delta_{\nabla\Jhat_\red}(\mu)$ and $\Delta_{\widetilde\nabla\Jhat_\red}(\mu)$ both decay with second order as the RB spaces grow; cf.~\cref{sec:TRRB_estimator_study}.
We also point out that $\Delta_{\nabla\Jhat_\red}(\mu)$ is an improved estimator, which can be used to replace the ineffective estimator $\Delta^{\red,*}_{\nabla_\mu \Jhat_\red}(\mu)$.
Certainly, both higher-order estimators $\Delta_{\nabla\Jhat_\red}(\mu)$ and $\Delta_{\tilde{\nabla}\Jhat_\red}(\mu)$ come with the price of computing the dual norm of the sensitivity residuals in~\eqref{sens_res_pr} and~\eqref{sens_res_du} for each direction, which aggravates the computational complexity.
Importantly, however, the gradient estimate is not part of the TR-RB method and therefore, the computational issue is only present in TR-RB algorithms that follow \cite{QGVW2017} or use classical goal-oriented offline-online approaches. 

\subsection{The Hessian of the reduced functional}

We also provide an a posteriori error result for the Hessian of the NCD-corrected reduced functional.
Importantly, just as the sensitivity and gradient estimates from \cref{sec:RB_estimates_for_sens} and \cref{sec:TRRB_RB_estimates_for_grad}, respectively, this error estimator is not part of the TR-RB method and is only stated for completeness.

\begin{proposition}[Upper bound on the model reduction error of the Hessian of the reduced output]
	\label{prop:Hessian_Jhat_error_NCD} 
	For the Hessian $\cHhath(\mu)$ of $\Jhat_h(\mu)$ and the true Hessian 
	$\cHhatn(\mu)$ of the NCD-corrected functional from \cref{prop:true_corrected_reduced_Hessian},
	we have the a posteriori error bound
	\begin{align*}
		\big|\cHhath(\mu) - \cHhatn(\mu)&\big| \leq \Delta_{{\HH}}(\mu) := \Big\| \big(\Delta_{\HH_{i,l}}(\mu)\big)_{i,l} \Big\|_2,
	\end{align*}
	with
	\begin{align*}
		\Delta_{{\HH}_{i,l}}(\mu) :=& \,\, \Delta_\pr(\mu)\Big(
		\cont{\partial_{\mu_i}\partial_{\mu_l} \jformd}
		+ 2\cont{\partial_{\mu_i}\partial_{\mu_l} \kformd} \|u_{\red, \mu}\|
		+\cont{\partial_{\mu_i}\partial_{\mu_l} \bformd}\|p_{\red, \mu}\|
		\\
		&\quad \quad \quad \quad + 2\cont{\partial_{\mu_i} \kformd}\|\dred{\mu_l}u_{\red, \mu}\|
		+ \cont{\partial_{\mu_i} \bformd}\|\dred{\mu_l}p_{\red, \mu}\|
		\Big)
		\\
		&+\Delta_{d_{\mu_l}\pr}(\mu)\Big(
		\cont{\partial_{\mu_i} \jformd} + 2\cont{\partial_{\mu_i} \kformd}\|u_{\red, \mu}\| +\cont{\partial_{\mu_i} \bformd} \|p_{\red, \mu}\|
		\Big)
		\\
		&+\Delta_\du(\mu)\Big(
		\cont{\partial_{\mu_i}\partial_{\mu_l} \lformd}
		+\cont{\partial_{\mu_i}\partial_{\mu_l} \bformd} \|u_{\red, \mu}\|
		+\cont{\partial_{\mu_i} \bformd} \|d_{\mu_l}u_{\red, \mu}\|
		\Big)
		\\
		&+\Delta_{d_{\mu_l}\du}(\mu)\Big(
		\cont{\partial_{\mu_i} \lformd} +\cont{\partial_{\mu_i} \bformd} \|u_{\red, \mu}\|
		\Big) +(\Delta_\pr)^2(\mu)\Big(
		\cont{\partial_{\mu_i}\partial_{\mu_l} \kformd}
		\Big)
		\\
		&+\Delta_\pr(\mu)\Delta_\du(\mu)\Big(
		\cont{\partial_{\mu_i}\partial_{\mu_l} \bformd}
		\Big)+\Delta_\pr(\mu)\Delta_{d_{\mu_l}\pr}(\mu)\Big(
		2\cont{\partial_{\mu_i} \kformd}
		\Big)
		\\
		&+\Delta_\pr(\mu)\Delta_{d_{\mu_l}\du}(\mu)\Big(
		\cont{\partial_{\mu_i} \bformd}
		\Big)+\Delta_{d_{\mu_l}\pr}(\mu)\Delta_\du(\mu)\Big(
		\cont{\partial_{\mu_i} \bformd}
		\Big) \\
		&+ \cont{\partial_{\mu_i}\bformd} \|d_{\mu_l} u_{\red,\mu}\|~\|w_{\red,\mu}\| 
		+ \cont{\partial_{\mu_i}\lformd} \|d_{\mu_l} w_{\red,\mu}\| +  
		\cont{\partial_{\mu_i}\bformd}  \|u_{\red,\mu}\|~\|d_{\mu_l} w_{\red,\mu}\| \\
		&+ 2 \cont{\partial_{\mu_i}\kformd} \|z_{\red,\mu}\|~\|d_{\mu_l} u_{\red,\mu}\| 
		+ \cont{\partial_{\mu_i}a_\mu} \|z_{\red,\mu}\|~\|d_{\mu_l} p_{\red,\mu}\| 
		+ \cont{\partial_{\mu_i}j_\mu} \|d_{\mu_l} z_{\red,\mu}\| \\
		&+ 2\cont{\partial_{\mu_i}k_\mu} \|d_{\mu_l} z_{\red,\mu}\|~\|u_{\red,\mu}\| 
		+ \cont{\partial_{\mu_i}a_\mu} \|d_{\mu_l} z_{\red,\mu}\|~\|p_{\red,\mu}\| \\
		&+ \cont{\partial_{\mu_i}\partial_{\mu_l}\lformd} \|w_{\red,\mu}\| + \cont{\partial_{\mu_i}\partial_{\mu_l}\bformd} \|w_{\red,\mu}\|~\|u_{\red,\mu}\| \\
		&+ \cont{\partial_{\mu_i}\partial_{\mu_l}\jformd}\|z_{\red,\mu}\| 
		+2\cont{\partial_{\mu_i}\partial_{\mu_l}\kformd} \|z_{\red,\mu}\|~\|u_{\red,\mu}\|
		+\cont{\partial_{\mu_i}\partial_{\mu_l}\bformd} \|z_{\red,\mu}\|~\|p_{\red,\mu}\|,
	\end{align*}
	where $\|\cdot\|_2$ denotes the spectral norm for matrices. The norm of the auxiliary functions $\|w_{\red,\mu}\|$, $\|z_{\red,\mu}\|$ and the norm of their sensitivities $\|d_{\mu_l}w_{\red,\mu}\|$ and 
	$\|d_{\mu_l}z_{\red,\mu}\|$ can be estimated by
	\begin{enumerate}
		\item[(i)] $\|z_{\red, \mu}\| \leq \alpha_{\bformd}^{-1} \|r_\mu^\pr(u_{\red,\mu})\|$,
		\item[(ii)] $\|w_{\red, \mu}\| \leq \alpha_{\bformd}^{-1} \left( \|r_\mu^\du(u_{\red,\mu},p_{\red,\mu})\| 
		+ 2\cont{k_\mu} \|z_{\red, \mu}\| \right)$,
		\item[(iii)] 	$\|d_{\mu_l}z_{\red,\mu}\| \leq \alpha_{\bformd}^{-1} \left(\|\resd^{\pr,d_{\mu_i}}(u_{\red, \mu}, d_{\mu_i}u_{\red, \mu})\| + \cont{\partial_{\mu_l} \bformd}\|z_{\red,\mu}\| \right) $,
		\item[(iv)] 	$\|d_{\mu_l}w_{\red,\mu}\|  \leq 
		\alpha_{\bformd}^{-1}\Big(\|\resd^{\du,d_{\mu_i}}(u_{\red, \mu}, p_{\red, \mu}, d_{\mu_i}u_{\red, \mu}, d_{\mu_i}p_{\red, \mu}) \| 
		+ 2\cont{k_\mu}\|d_{\mu_l}z_{\red,\mu}\| \\
		+ 2\cont{\partial_\mu k_\mu}\|z_{\red,\mu}\| + \cont{\partial_\mu \bformd}\|w_{\red,\mu}\| \Big)$.
	\end{enumerate}
\end{proposition}
\begin{proof} 
	To prove the Hessian estimate, we recall that for all $i,l$ we have
	\[
	\begin{aligned}
		\big(\cHhatn(&\mu))_{i, l}  = \partial_\mu \left(j_\mu(d_{\mu_l} u_{\red,\mu})+2k_\mu(u_{\red,\mu},d_{\mu_l} u_{\red,\mu})  - a_\mu(d_{\mu_l} u_{\red,\mu},p_{\red,\mu}+w_{\red,\mu}) \right.\\
		& + r^\pr_\mu(u_{\red,\mu})[d_\nu p_{\red,\mu}+d_{\mu_l} w_{\red,\mu}] - 2 k_\mu (z_{\red,\mu},d_{\mu_l} u_{\mu,\red}) \\
		& + a_\mu(z_{\red,\mu},d_{\mu_l} p_{\red,\mu}) - r^\du_\mu(u_{\red,\mu},p_{\red,\mu})[d_{\mu_l} z_{\red,\mu}]   \\
		& \left.+\partial_\mu (\J(u_{\red,\mu},\mu) + r^\pr_\mu(u_{\red,\mu})[p_{\red,\mu}+w_{\red,\mu}]
		- r_\mu^\du(u_{\red,\mu},p_{\red,\mu})[z_{\red,\mu}])\cdot e_l \right)\cdot e_i,
	\end{aligned}
	\]   
	and
	\begin{align*}
		\big(\cHhath(\mu)\big)_{i, l} \kern -0.1em
		=  \kern -0.1em \partial_\mu\Big(&\kern -0.1em\partial_u \J(u_{h,\mu}, \mu)[d_{\mu_l}u_{h,\mu}] \kern -0.2em+\kern -0.2em r_\mu^\pr(u_{h,\mu})[d_{\mu_l}p_{h,\mu}]\kern -0.2em -\kern -0.1em a_{\mu}(d_{\mu_l}u_{h,\mu}, p_{h,\mu}) \\
		&+\partial_\mu\big(\J(u_{h,\mu}, \mu)+ r_\mu^\pr(u_{h,\mu})[p_{h,\mu}]\big)\cdot e_l\Big)\cdot e_i.
	\end{align*}
	We obtain
	\begin{align*}
		\big|\big(&\cHhath(\mu) - \cHhatn(\mu)\big)_{i, l}\big| \\ &\leq \big| \partial_\mu\Big(\partial_u \J(u_{h,\mu}, \mu)[d_{\mu_l}u_{h,\mu}] - \partial_u \J(u_{\red,\mu}, \mu)[	d_{\mu_l}u_{\red, \mu}] + l_{\mu}(d_{\mu_l}p_{h,\mu}) 
		\\&\quad-  l_{\mu}(d_{\mu_l}p_{\red, \mu}) 
		- a_{\mu}(d_{\mu_l}u_{h,\mu}, p_{h,\mu}) + a_{\mu}(d_{\mu_l}u_{\red, \mu}, p_{\red,\mu}) 
		- a_{\mu}(u_{h,\mu}, d_{\mu_l}p_{h,\mu}) \\
		&\quad+ a_{\mu}(u_{\red,\mu}, d_{\mu_l}p_{\red, \mu}) 
		+ a_\mu(d_{\mu_l} u_{\red,\mu},w_{\red,\mu}) - r^\pr_\mu(u_{\red,\mu})[d_{\mu_l} w_{\red,\mu}]\\ 
		&\quad+ 2 k_\mu (z_{\red,\mu},d_{\mu_l} u_{\mu,\red}) 
		- a_\mu(z_{\red,\mu},d_{\mu_l} p_{\red,\mu}) + r^\du_\mu(u_{\red,\mu},p_{\red,\mu})[d_{\mu_l} z_{\red,\mu}]\\
		&\quad+\partial_\mu\big(\J(u_{h,\mu}, \mu) - \J(u_{\red,\mu}, \mu) 
		+ l_{\mu}(p_{h,\mu}) - l_{\mu}(p_{\red,\mu}) 
		- a_{\mu}(u_{h,\mu}, p_{h,\mu}) \\
		&\qquad\quad+  a_{\mu}(u_{\red,\mu}, p_{\red,\mu}) 
		+ r^\pr_\mu(u_{\red,\mu})[w_{\red,\mu}]- r_\mu^\du(u_{\red,\mu},p_{\red,\mu})[z_{\red,\mu}] \big)\cdot e_l\quad\Big)\cdot e_i \big|. 
	\end{align*}
	For the first terms we see that 
	\begin{align*}
		\partial_\mu\big(\partial_u &\J(u_{h,\mu}, \mu)[d_{\mu_l}u_{h,\mu}] - \partial_u \J(u_{\red,\mu}, \mu)[\dred{\mu_l}u_{\red, \mu}] \big)\cdot e_i 
		\\ & =\partial_\mu \big( j_{\mu}(d_{\mu_l} e_{h, \mu}^\pr) + 2 k_{\mu}(d_{\mu_l} e_{h, \mu}^\pr,u_{h,\mu}) + 2 k_{\mu}(\dred{\mu_l}u_{\red, \mu},e_{h, \mu}^\pr) \big)\cdot e_i. 
	\end{align*}
	Obviously, this equation still incorporates the norm of the FOM solution $u_{h,\mu}$. However, we can simply estimate these FOM quantities by $\|u_{h,\mu}\| \leq \Delta_\pr(\mu) + \|u_{\red,\mu}\|$. 
	For the second terms we have
	\begin{align*}
		\partial_\mu(l_{\mu}(d_{\mu_l}p_{h,\mu}) -  l_{\mu}(\dred{\mu_l}p_{\red, \mu}) )\cdot e_i =  \partial_\mu(l_{\mu}(d_{\mu_l} e_{h, \mu}^\du))\cdot e_i.
	\end{align*}
	The third terms can be determined by
	\begin{align*}
		- \partial_\mu\big( a_{\mu}(d_{\mu_l}u_{h,\mu}, p_{h,\mu}) -& a_{\mu}(\dred{\mu_l}u_{\red, \mu}, p_{\red,\mu})\big)\cdot e_i 
		&= -  \partial_\mu\big( a_{\mu}(d_{\mu_l} e_{h, \mu}^\pr, p_{h,\mu}) - a_{\mu}(\dred{\mu_l}u_{\red, \mu}, e_{h,\mu}^\du)\big)\cdot e_i
	\end{align*}
	and similarly we have for the fourth term
	\begin{align*}
		- \partial_\mu\big( a_{\mu}(u_{h,\mu}, d_{\mu_l}p_{h,\mu}) &- a_{\mu}(u_{\red,\mu}, \dred{\mu_l}p_{\red, \mu})\big)\cdot e_i 
		&= -  \partial_\mu\big( a_{\mu}(e_{h, \mu}^\pr, d_{\mu_l}p_{h,\mu}) - a_{\mu}(u_{\red,\mu},d_{\mu_l} e_{h, \mu}^\du)\big)\cdot e_i.
	\end{align*}
	With the same strategy as above, we have for the first part of the second derivatives that
	\begin{align*}
		\partial_\mu\Big(&\partial_\mu\big(\J(u_{h,\mu}, \mu) - \J(u_{\red,\mu}, \mu) + l_{\mu}(p_{h,\mu}) - l_{\mu}(p_{\red,\mu}) 
		- a_{\mu}(u_{h,\mu}, p_{h,\mu}) +  a_{\mu}(u_{\red,\mu}, p_{\red,\mu}) \big)\cdot e_l\Big)\cdot e_i  \\
		&\hS{40} 
		= \partial_\mu\Big(\partial_\mu\big(j_{\mu}(e_{h, \mu}^\pr) + 2 k_{\mu}(e_{h, \mu}^\pr,u_{h,\mu}) + 2 k_{\mu}(u_{\red, \mu}^\pr,e_{h, \mu}^\pr)  
		+ l_{\mu}(e_{h, \mu}^\du) \\
		&\hS{110}
		- (a_{\mu}(e_{h, \mu}^\pr, p_{h,\mu}) -  a_{\mu}(u_{\red,\mu}, e_{h, \mu}^\du) \big)\cdot e_l\Big)\cdot e_i.
	\end{align*}
	For the rest of the proof we simply use the Cauchy-Schwarz inequality for all terms and sum all pieces together to conclude $\Delta_{{\HH}_{i,l}}(\mu)$.  
	
	For a proof of (i) and (ii), we point to \cref{prop:grad_Jhat_error}.
	For the first sensitivity estimation (iii), we use the equations \eqref{z_eq_sens} and \eqref{sens_res_pr} to obtain
	\begin{align*}
		\alpha_{\bformd} \|d_{\mu_l}z_{\red,\mu}\|^2 &\leq \bformd(d_{\mu_l}z_{\red,\mu},d_{\mu_l}z_{\red,\mu}) \\&= 
		-\partial_\mu(r_\mu^\pr(u_{\red,\mu})[d_{\mu_l}z_{\red,\mu}] + a_\mu(z_{\red,\mu},d_{\mu_l}z_{\red,\mu}))\cdot \nu + a_\mu(d_\nu u_{\red,\mu},d_{\mu_l}z_{\red,\mu}) \\
		&= -\resd^{\pr,d_{\mu_i}}(u_{\red, \mu}, d_{\mu_i}u_{\red, \mu})[d_{\mu_l}z_{\red,\mu}] - \partial_\mu \bformd(z_{\red,\mu},d_{\mu_l}z_{\red,\mu}) \cdot \nu \\
		&\leq \left(\|\resd^{\pr,d_{\mu_i}}(u_{\red, \mu}, d_{\mu_i}u_{\red, \mu})\| + \cont{\partial_{\mu_l} \bformd}\|z_{\red,\mu}\| \right) \|d_{\mu_l}z_{\red,\mu}\|.
	\end{align*}
	For (iv), we instead use \eqref{w_eq_sens} and \eqref{sens_res_du} and yield
	\begin{align*}
		\alpha_{\bformd} \|d_{\mu_l}w_{\red,\mu}\|^2 &\leq \bformd(d_{\mu_l}w_{\red,\mu},d_{\mu_l}w_{\red,\mu}) \\
		&= \partial_\mu (r_\mu^\du(u_{\red,\mu},p_{\red,\mu})[d_{\mu_l}w_{\red,\mu}]-2k_\mu(z_{\red,\mu},d_{\mu_l}w_{\red,\mu})-a_\mu(d_{\mu_l}w_{\red,\mu},w_{\red,\mu}))\cdot\nu \\
		&\quad +2k_\mu(d_{\mu_l}w_{\red,\mu},d_\nu u_{\red,\mu}-d_\nu z_{\red,\mu})-a_\mu(d_{\mu_l}w_{\red,\mu}, d_\nu p_{\red,\mu}) \\
		&= \resd^{\du,d_{\mu_i}}(u_{\red, \mu}, p_{\red, \mu}, d_{\mu_i}u_{\red, \mu}, d_{\mu_i}p_{\red, \mu})[d_{\mu_l}w_{\red,\mu}]  
		- 2k_\mu(d_{\mu_l}w_{\red,\mu},d_\nu z_{\red,\mu}) \\
		&\quad - \partial_\mu (2k_\mu(z_{\red,\mu},d_{\mu_l}w_{\red,\mu})-a_\mu(d_{\mu_l}w_{\red,\mu},w_{\red,\mu}))\cdot\nu \\
		&\leq \Big(\|\resd^{\du,d_{\mu_i}}(u_{\red, \mu}, p_{\red, \mu}, d_{\mu_i}u_{\red, \mu}, d_{\mu_i}p_{\red, \mu}) \| + 2\cont{k_\mu}\|d_{\mu_l}z_{\red,\mu}\| \\
		&\qquad \qquad + 2\cont{\partial_\mu k_\mu}\|z_{\red,\mu}\| + \cont{\partial_\mu \bformd}\|w_{\red,\mu}\| \Big) \|d_{\mu_l}w_{\red,\mu}\|.
	\end{align*}	
\end{proof}

We note that the estimator $\Delta_{{\HH}}(\mu)$ inherits estimation for the primal and dual variable, their sensitivities, and corresponding auxiliary functions.

\subsection{The optimal parameter}
\label{sec:TRRB_RB_estimates_param}

Following ideas from~\cite{dihl15,KTV13}, we derive an error estimation for the optimal parameter consisting of the gradient and Hessian of the FOM cost functional.
This estimator relies on the following second-order condition for a strict local minimum $\bar\mu_h$ of $\Jhat_h$, i.e.
\begin{align}
	\label{coerc-cond}
	\nu \cdot ( \cHhath(\bar \mu_h) \cdot \nu) \geq \lambda_\textnormal{min}\left\|\nu\right\|_2^2 && \textnormal{for all } \nu\in\mathcal{C}(\bar\mu_h)\setminus\left\{0\right\},
\end{align} 
where $\lambda_\textnormal{min}$ is the smallest eigenvalue of $\cHhath(\bar\mu_h)$, since the parameter space is finite-dimensional.
Note that~\eqref{coerc-cond} is equivalent to the second-order sufficient optimality condition from Proposition~\ref{prop:second_order}.
Let~\eqref{coerc-cond} be fulfilled. Then, for any $\tilde{\lambda}$, with $0<\tilde{\lambda}<\lambda_\textnormal{min}$, there exists a radius $r(\tilde{\lambda})>0$, such that for all $\mu\in \mathcal{B}(\bar \mu_h, r(\tilde{\lambda}))$, the closed ball of radius $r(\tilde{\lambda})$ centered in $\bar\mu_h$, it holds:
\begin{align*}
	\nu \cdot ( \cHhath(\mu) \cdot \nu) \geq \tilde{\lambda}\left\|\nu\right\|_2^2 && \textnormal{for all } \nu\in\mathcal{C}(\bar\mu_h)\setminus\left\{0\right\}.
\end{align*} 

\begin{proposition}[Upper bound for optimal parameters with the full-order model] \label{Prop:argmin}
	Let Assumption~{\ref{asmpt:truth}} be satisfied. Moreover, let $\bar\mu_h$ and $\bar \mu_\red$ be strict local minima for the optimization problems~\eqref{Phat_h} and~\eqref{Phat_\red}, respectively. If $\bar\mu_\red \in \mathcal{B}(\bar \mu_h, r(\lambda_\textnormal{min}/2))$, then it holds 
	\begin{equation}
		\label{apo-est-parameters}
		\| \bar \mu_h - \bar \mu_\red \|_2 \leq \Delta_\mu(\bar \mu_\red) := \frac{2}{\lambda_\textnormal{min}} \left\| \zeta \right\|_2,
	\end{equation}
	where $\zeta= (\zeta_i)\in\mathbb{R}^P$ with
	\[
	\zeta_i:= \left\{ \begin{array}{ll} -\min(0,(\nabla\Jhat_h(\bar\mu_\red))_i) & \textnormal{if } \bar\mu_{\red,i}=(\mu_\mathsf{a})_i \\
		-\max(0,(\nabla\Jhat_h(\bar\mu_\red))_i) & \textnormal{if } \bar\mu_{\red,i}=(\mu_\mathsf{b})_i   \\
		-(\nabla\Jhat_h(\bar\mu_\red))_i & \textnormal{otherwise}
	\end{array} \right.
	\]
	for $i=1,\ldots,P$.	
\end{proposition}
\begin{proof}
	Note that Assumption~\ref{asmpt:truth} implies that the distance between $\bar\mu_h$ of~\eqref{Phat_h} and a strict local minimum $\bar\mu$ of $\Jhat$ satisfying~\eqref{coerc-cond} is negligible, thus we can follow the proof of~\cite[Theorem~3.4]{KTV13}.
\end{proof}

\cref{Prop:argmin} requires the strong assumption that the FOM and RB models are accurate enough to have the parameters $\bar\mu_h$ and $\bar\mu_\red$ sufficiently close to a local minimum $\bar\mu$. 
In {\cite{dihl15}}, a sufficient condition based on the FOM gradient and Hessian is given to guarantee this in case $\bar\mu_h,\bar\mu_\red\in \textnormal{int }\Params$.

Due to \cref{Prop:argmin}, we can estimate the distance to the optimal parameter $\bar\mu_h$ without explicitly computing $\bar\mu_h$.
The computation of $\zeta$ is not costly since the FOM adjoint solution is available, cf.~\cref{sec:TR_RB_algorithm}.
On the other hand, the computation of $\lambda_{\textnormal{min}}$ requires the evaluation of the FOM Hessian, which is a costly procedure.
Admittedly, this can be sped up with a cheap estimation of the eigenvalue.
In {\cite[Proposition~6]{dihl15}}, the authors utilize the smallest eigenvalue of the reduced-order Hessian under suitable conditions.
In our numerical tests, these conditions were never confirmed, implying the inapplicability of the mentioned cheap estimate in our case.
For the sake of completeness, let us mention that another technique is to compute $\lambda_{\textnormal{min}}$ in advance on a grid in $\Params\subset\mathbb{R}^P$ when $P$ is sufficiently small.
This approach can even be performed in parallel since each eigenvalue computation is independent; cf.~\cite[Section~6.4.1]{Trenz2017}.

With regard to the fact that the computational cost for evaluating estimate~\eqref{apo-est-parameters} has FOM complexity, we only use it as a post-processing tool for finding sufficiently small termination tolerance of the TR-RB algorithm; cf. \cref{sec:outer_stopping}.
For a numerical experiment of this approach, we refer to \cref{num_test:12params}.

\subsection{The Petrov--Galerkin approach}
\label{sec:TRRB_RB_estimates_pg}

For the Petrov--Galerkin RB approach from \cref{sec:TRRB_pg_approach}, we utilize the same standard residual-based estimation as presented in the former sections with the vital difference that the inf-sup stability w.r.t. the respective test- and ansatz space needs to be used.
The error result is summarized in the following proposition.

\begin{proposition}[Upper error bound for the reduced quantities] \label{prop:error_reduced_quantities}
	For $\mu \in \Params$, let $u_{h, \mu} \in V_h^\pr$ and $p_{h, \mu} \in V_h^\du$ be solutions of~\eqref{eq:state_h} and~\eqref{eq:dual_solution_h}, respectively, and let $u^\textnormal{pg}_{\red, \mu} \in V_\red^\pr$, $p^\textnormal{pg}_{\red, \mu} \in V_\red^\du$ be a solution of the PG-reduced equations~\eqref{eq:state_red_pg} and~\eqref{eq:dual_solution_red_pg}.
	Then it holds
	\begin{enumerate}
		\item[(i)]
		$\|u_{h, \mu} - u^\textnormal{pg}_{\red, \mu}\| \leq \Delta^\textnormal{pg}_\pr(\mu) := (\infsup{\mu})^{-1}\, \|\resd^\pr(u^\textnormal{pg}_{\red, \mu})\|$,
		\item[(ii)] 	$\|p_{h, \mu} - p^\textnormal{pg}_{\red, \mu}\| \leq \Delta^\textnormal{pg}_\du(\mu) := (\infsup{\mu})^{-1}\big(2 \cont{\kformd}\;\Delta^\textnormal{pg}_\pr(\mu) + \|\resd^\du(u^\textnormal{pg}_{\red, \mu}, p^\textnormal{pg}_{\red, \mu})\|\Big)$,
		\item[(iii)]  $|\Jhat_h(\mu) - \cJhatn^\textnormal{pg}(\mu)| \leq \Delta^\textnormal{pg}_{\cJhatn}(\mu)
		:=  \Delta^\textnormal{pg}_\pr(\mu) \|\resd^\du(u^\textnormal{pg}_{\red, \mu}, p^\textnormal{pg}_{\red,\mu})\| + \Delta^\textnormal{pg}_\pr(\mu)^2 \cont{\kformd},$
	\end{enumerate} 
	where $\infsup{\mu}$ and $\cont{\kformd}$ define the inf-sup stability constant of $\bformd$ and the continuity constant of $\kformd$.
\end{proposition}
\begin{proof}
	For $(i)$ and $(ii)$, apart from using the inf-sup stability of $\bformd$, stated in~\eqref{eq:inf_sup_trivial}, in the first step, we proceed along the same lines as in \cref{prop:primal_rom_error} and \cref{prop:dual_rom_error}.
	For $(iii)$, we refer to \cref{prop:Jhat_error}$(ii)$, noting that $r_\mu^\pr(u^\textnormal{pg}_{\red,\mu})[p^\textnormal{pg}_{\red,\mu}] = 0$.
\end{proof}
It is important to mention that the computation of these error estimators includes the computation of the (parameter dependent) inf-sup constant of $\bformd$, which involves an eigenvalue problem on the FOM level.
In practice, cheaper techniques, such as the successive constraint method~\cite{huynh2015methods}, can be used.
For the conforming approach $V_h^\pr = V_h^\du$, the inf-sup constant is equivalent to the coercivity constant $\alpha_{\bformd}$ of $\bformd$, which can be cheaply bounded from below with the help of the min-theta approach; c.f.~\cref{sec:min_theta}.

With the help of a posteriori theory devised for all reduced models from \cref{sec:TRRB_MOR_for_PDEconstr}, we showed that the respective model reduction error of the reduced quantities can be efficiently bounded.
These error estimates are of significant importance for the adaptive TR-RB procedure, introduced in the following section.

\section{Trust-region reduced basis algorithm}
\label{sec:TR_RB_algorithm}

From the former section it is particularly vital that there exists an error estimator for the reduced objective function, such that
\begin{equation} \label{eq:estimation_J}
| \Jhat_h(\mu) - \Jhat_\red(\mu) | \leq \Delta_{\Jhat_\red}(\mu).
\end{equation}
As discussed in \cref{sec:background_PDE-constr}, a standard approach for solving \eqref{Phat_h} is to follow a negligibly long offline phase, performing a goal-oriented greedy-search algorithm constructing a suitable surrogate model as detailed in \cref{sec:background_RB_methods_for_PDEopt}.
However, our work is targeted toward accelerating the overall procedure of a single optimization task, where, as shown in the experiment in \cref{sec:background_RB_methods_for_PDEopt}, the classical offline-online procedure is impractical.

In this section, we present the trust-region reduced basis method for an overall-efficient algorithm to solve Problem~\eqref{Phat_h}.
In particular, we apply a TR method, which iteratively computes a first-order critical point of~\eqref{Phat_h} using a surrogate model as often as possible.
The basic idea of TR methods was presented in \cref{sec:background_TR_methods}, where \cref{alg:BTR} is referred to as the basic trust-region (BTR) algorithm, with the choice of a model function $m^{(k)}$, denoting a cheaply computable approximation of the cost functional $\J$ in a neighborhood of the parameter $\mu^{(k)}$, i.e., the trust-region; cf.~Step 1 of the BTR algorithm.
In Step 2, the TR sub-problem is solved, and, in Step 3, the acceptance of the iterate is verified.
In Step~4, either the TR is shrunk (after rejection) or possibly enlarged (after acceptance).
The error-aware TR-RB version that we present in this section is based on~\cite{QGVW2017} which is based on~\cite{YM2013}, where PDE-equality-constrained parameter optimization problem without bounds on the parameter space is considered.

For PDE-constrained parameter optimization, the TR minimization sub-problem~\eqref{TRsubprob}, for $k\geq 0$, given a TR radius $\delta^{(k)}$, can be rewritten as
\begin{equation}
	\label{eq:TRsubprob}
	\min_{s\in \mathbb{R}^P} m^{(k)}(s) \, \textnormal{ subject to } \|s\|_2 \leq \delta^{(k)},\, \widetilde{\mu}:= \mu^{(k)}+s \in\Paramsad \text{ and } r_{\tilde{\mu}}^\pr(u_{\tilde{\mu}})[v]= 0 \, \textnormal{ for all }  v\in V,
\end{equation}
where the unique solution $\bar s^{(k)}$ is used to compute the next iterate $\mu^{(k+1)} = \mu^{(k)} + \bar s^{(k)}$.
Regarding the acceptance of the TR iterate, we point out that, the quotient  
$
\rho^{(k)} 
$
of true and approximate values, used in the BTR algorithm, are assumed to be unknown in \cite{YM2013,QGVW2017}.
Possible approximate sufficient and necessary conditions for convergence, depending on the approximate generalized Cauchy point (AGC) $\mu^{(k)}_\textnormal{{AGC}}$ (see \cref{def:AGC}), are instead given in~\cite{YM2013}.
Furthermore, in~\cite{QGVW2017}, shrinking but no enlargement of the TR radius is considered.

This chapter is concerned with refactoring the initially proposed TR-RB algorithm of~\cite{QGVW2017} and imposing additional bilateral parameter constraints in~\eqref{TRsubprob}.
The presence of the additional inequality constraints requires a review of the proof of convergence for the TR-RB algorithm, whereas, in~\cite{QGVW2017}, the convergence is based on the results contained in~\cite{YM2013}.

We start by presenting the algorithm in detail, with a strong emphasis on how our method differs from the one in~\cite{QGVW2017}.
We moreover discuss optional features that were particularly added in~\cite{paper2}.
In addition, we present two subsequently deduced convergence results for our method:
The original RB-based convergence study from~\cite{paper1} and the improved version elaborated in~\cite{paper2}.
Subsequently, several basis enrichment technique and many variants are presented.

\subsection{The basic algorithm}
\label{sec:TR_algorithm_details}

We present the TR-RB algorithm with the example of the NCD-corrected reduced functional $\cJhatn$, which has been introduced in \cref{sec:TRRB_ncd_approach}.
The algorithm can also be used by the other reduced formulations that we discussed in \cref{sec:TRRB_MOR_for_PDEconstr} (simply by exchanging the reduced functional, error estimation, and/or derivative information).
For a complete list of variants that we elaborate on, we refer to \cref{sec:TRRB_variants}.
Note that the algorithm discussed in this section contains all features that were independently published in~\cite{paper1} and~\cite{paper2}.

\subsubsection{Inexact error-aware TR sub-problem}
Let the model function be the NCD-corrected RB reduced functional $\cJhatn^{(k)}$ defined in (\ref{eq:Jhat_red_corected}), i.e.~$m^{(k)}(\cdot)= \cJhatn^{(k)}(\mu^{(k)}+\cdot)$ for $k\geq 0$, where the super-index $(k)$ indicates that we use different RB spaces in each iteration. 
We initialize the RB spaces $V^{\pr,(0)}_\red$ and $V^{\du,(0)}_\red$ using the FOM primal- and dual solutions at the initial guess $\mu^{(0)}$.
At every iteration $k$, we enrich the obtained space at the current iterate $\mu^{(k+1)}$ -- for further details on enrichment strategies and optional enrichment, see \cref{sec:TRRB_construct_RB} and \cref{sec:TRRB_optional_enrichment}, respectively.
According to~\cite{QGVW2017}, the inexact RB version of problem~\eqref{eq:TRsubprob} is 
\vspace{-0.3cm}
\begin{equation} 
	\label{TRRBsubprob}
	\min_{\widetilde{\mu}\in\Paramsad} \cJhatn^{(k)}(\widetilde{\mu}) \qquad \textnormal{ such that } \qquad \frac{\Delta_{\Jhat_\red^{(k)}}(\widetilde{\mu})}{\cJhatn^{(k)}(\widetilde{\mu})}\leq \delta^{(k)},
\end{equation}
where $\widetilde{\mu}:= \mu^{(k)}+s$.
The equality constraint $r_{\tilde\mu}^\pr(u_{\tilde{\mu}})[v]= 0$ is hidden in the definition of $\cJhatn$ and the inequality constraints are concealed in the request $\widetilde{\mu}\in \Paramsad$.
Due to the presence of bilateral constraints on the parameters, as defined in \cref{sec:constr_optimization}, we recall the projection operator $\Proj_{\Paramsad}: \R^P \rightarrow {\Paramsad}$ as
\vspace{-0.3cm}
\begin{align*}
	(\Proj_{\Paramsad}(\mu))_i:= \left\{ \begin{array}{ll}
		(\mu_\mathsf{a})_i & \textnormal{if } \mu_i\leq (\mu_\mathsf{a})_i, \\
		(\mu_\mathsf{b})_i & \textnormal{if } \mu_i\geq (\mu_\mathsf{b})_i, \\
		\mu_i & \textnormal{otherwise}
	\end{array} \right.  && \textnormal{for } i=1,\ldots,P.
\end{align*}
The operator $\Proj_{\Paramsad}$ is Lipschitz continuous with Lipschitz constant one; cf.~\cite{kelley}.
The additional TR constraint, instead, is treated with a backtracking technique; cf.~\cite{QGVW2017}.

\subsubsection{Details on the sub-problem}
To solve the TR sub-problem~\eqref{TRRBsubprob}, we use a projected descent optimization method combined with a line-search algorithm as introduced in \cref{sec:constr_optimization}.
We generate a sequence $\{\mu^{(k,\el)}\}_{\el=1}^L$, where $L$ is the last inner iteration, and we set the TR iterate to $\mu^{(k+1)}:=\mu^{(k,L)}$.
The index $k$ refers to the current outer TR iteration, and $\el$ refers to the inner iteration. 
Note that $L$ may be different for each iteration $k$ which we indicate only when strictly necessary. We define
\begin{align}
	\label{eq:General_Opt_Step_point}
	\mu^{(k,\el)}(j):= \Proj_{\Paramsad}(\mu^{(k,\el)} + \kappa^j d^{(k,\el)}) \in{\Paramsad} && \textnormal{for } j\geq 0,
\end{align}
where $\kappa\in(0,1)$, $d^{(k,\el)}\in\mathbb{R}^P$ is the chosen descent direction at the iteration $(k,\el)$, specified further in \cref{sec:sub-problem_solvers}.
Similar to \cref{sec:constr_optimization}, we enforce an Armijo-type condition for inequality constraints 
\begin{subequations}
	\label{Arm_and_TRcond}
	\begin{equation}
		\label{Armijo}\cJhatn^{(k)}(\mu^{(k,\el)}(j)) - \cJhatn^{(k)}(\mu^{(k,\el)}) \leq  -\frac{\kappa_{\mathsf{arm}}}{\kappa^j} \| \mu^{(k,\el)}(j)-\mu^{(k,\el)}\|^2_2,
	\end{equation}
	with $\kappa_{\mathsf{arm}}=10^{-4}$, combined with an additional TR constraint on $\cJhatn^{(k)}$
	\begin{equation}
		\label{TR_radius_condition} q^{(k)}(\mu^{(k,\el)}(j)):= \frac{\Delta_{\Jhat_\red^{(k)}}(\mu^{(k,\el)}(j))}{\cJhatn^{(k)}(\mu^{(k,\el)}(j))} \leq \delta^{(k)}.
	\end{equation}
\end{subequations}
We thus select $\mu^{(k,\el+1)} = \mu^{(k,\el)}(j^{(k,\el)})$ for $\el\geq 1$, where $j^{(k,\el)}<\infty$ is the smallest index for which~\eqref{Arm_and_TRcond} holds.
Moreover, as termination criterion for the optimization sub-problem, we use
\begin{subequations}\label{Termination_crit_sub-problem}
	\begin{equation}
		\label{FOC_sub-problem}
		\big\|\mu^{(k,\el)}-\Proj_{\Paramsad}(\mu^{(k,\el)}-\nabla_\mu \Jhat_\red^{(k)}(\mu^{(k,\el)}))\big\|_2\leq \tau_\textnormal{{sub}}
	\end{equation}
	or
	\begin{equation}
		\label{Cut_of_TR}
		\beta_2\delta^{(k)} \leq \frac{\Delta_{\Jhat_\red^{(k)}}(\mu)}{\Jhat_\red^{(k)}(\mu)} \leq \delta^{(k)},
	\end{equation}
\end{subequations}
\noindent where $\tau_\textnormal{{sub}}\in(0,1)$ is a predefined tolerance and $\beta_2\in(0,1)$, generally close to one.
Condition~\eqref{Cut_of_TR} is used to prevent the optimizer from spending much time close to the boundary of the trust-region, where the model is poor in approximation; cf.~\cite{QGVW2017}.
Note that, without the projection operator $\Proj_{\Paramsad}$, conditions~\eqref{Arm_and_TRcond}-\eqref{Termination_crit_sub-problem} coincide with the ones in~\cite{QGVW2017}, apart from using the NCD-corrected reduced functional.

Regardless of the sub-problem solver, i.e., the method to determine the search direction $d^{(k,\el)}$, we anyway compute the AGC point $\mu^{(k)}_\textnormal{{AGC}}$ in the first iterate, defined in the following:
\begin{definition}[AGC point for simple bounds]
	\label{def:AGC}
	At iteration $k$, we define the AGC point as
	\[
	\mu_\textnormal{{AGC}}^{(k)}:= \mu^{(k,0)}(j^{(k)}_c)=  \Proj_{\Paramsad}(\mu^{(k,0)} + \kappa^{j^{(k)}_c} d^{(k,0)}),
	\]
	where $\mu^{(k,0)}:= \mu^{(k)}$, $d^{(k,0)}:= -\nabla_\mu \cJhatn^{(k)}(\mu^{(k,0)})$ and $j^{(k)}_c$ is the smallest non-negative integer $j$ for which $\mu^{(k,0)}(j)$ satisfies~\eqref{Arm_and_TRcond} for $\el=0$.
\end{definition}
While the AGC point can be computed cheaply, it plays an important role in accepting the sub-problem result and is thus fundamental for the convergence study.
We refer to \cite{YM2013} for further details on the AGC point.

\subsubsection{Sub-problem solvers}
\label{sec:sub-problem_solvers}
For computing the descent direction $d^{(k,\el)}$, we follow either the projected BFGS or the Newton algorithm.
The projected BFGS algorithm is reported in~\cite[Section~5.5.3]{kelley}, and (in the unconstrained version) discussed in \cref{sec:iterative_optimization}.
In both cases, tailored to the control constraints, we use an active set strategy such that only non-boundary components are treated in a standard way.
For the Newton variant, this means that we have
\begin{align*}
	d^{(k,\el)}=-(\mathcal R^{(k)}_\red(\mu^{(k,\el)}))^{-1}\nabla_\mu \cJhatn^{(k)}(\mu^{(k,\el)}) && \textnormal{ for all } k,\el\in \mathbb{N},\, \el\geq 1,
\end{align*}
where 
\vspace{-0.3cm}
\begin{align*}
	\left( \mathcal R^{(k)}_\red(\mu) \right)_{ij} = \left\{ \begin{array}{ll}
		\delta_{ij} & \textnormal{if } i\in \mathcal{A}^\varepsilon(\mu) \textnormal{ or } j\in \mathcal{A}^\varepsilon(\mu)\\
		(\cHhatn(\mu))_{i,j} & \textnormal{otherwise},
	\end{array} \right. && \textnormal{for } \mu\in\Params.
\end{align*}
The function $\delta_{ij}$ indicates the Kronecker delta, and the set $\mathcal{A}^\varepsilon$ is the $\varepsilon$-active set for the parameter constraints, i.e.
\[
\mathcal{A}^\varepsilon(\mu) = \left\{i\in\{1,\ldots,P\}\big| (\mu_\mathsf{b})_i-\mu_i \leq \varepsilon \textnormal{ or } \mu_i-(\mu_\mathsf{a})_i\leq \varepsilon \right\}.
\]
For further details on the projected Newton method, the choice of $\varepsilon$, and its effect on the convergence of the method, we refer to~\cite[Section 5.5]{kelley}.
Note that $\cHhatn(\mu)$ (and thus $\mathcal R^{(k)}_\red(\mu)$) might not be symmetric positive definite for every $\mu\in\Params$.
Therefore, as discussed in \cref{sec:newton_method}, we use a truncated Conjugate Gradient (CG) method to compute $d^{(k,\el)}$, where the CG terminates when a negative curvature condition criterion is triggered.
In such a way, we ensure that $d^{(k,\el)}$ (resulting from the possible premature termination of the CG) is still a descendent direction.
The truncated CG is explained in~\cite[Algorithm 7.1]{Nocedal}.
Both sub-problem solvers use the AGC point, which is not carried out naturally by the projected Newton method.
Although this fact seems disadvantageous compared to the projected BFGS method, where this computation is usually included in the process (cf.~\cite{QGVW2017}), we remark that the search of the AGC point costs only one projected gradient optimization step.
It is used as a warm start for the projected Newton method.
Therefore, the initial cost is justified by the faster locally quadratic convergence of the projected Newton method.
It constitutes an improvement concerning the projected BFGS method, in particular when the optimum is close to the boundary of the parameter set; cf.~\cite{kelley,Nocedal}.

\subsubsection{Acceptance of the TR-iterate}
\label{sec:acceptance_of_TR}

After the inner TR-routine returns a potential iterate $\mu^{(k+1)} = \mu^{(k)} + \bar s^{(k)}$, as usual in TR methods, we require a criterion for acceptance or rejection of the step.
In Step 3 in the BTR algorithm (cf. \cref{alg:BTR}), the quantity $\rho^{(k)}$ is used.
Admittedly, at least before the enrichment, this quantity is computationally costly since it evaluates $\J_h$.
Instead, as pointed out in~\cite{QGVW2017,YM2013}, an error-aware sufficient decrease condition
\begin{align}
	\label{eq:suff_decrease_condition}
	\cJhatn^{(k+1)}(\mu^{(k+1)})\leq \cJhatn^{(k)}(\mu_\textnormal{{AGC}}^{(k)})
\end{align} 
can be used at each iteration $k$ of the TR-RB algorithm.
We further note that, in the context of RB methods, the surrogate model is exact at the enrichment parameters (cf.~\cite{HAA2017}) and thus, computing \eqref{eq:suff_decrease_condition} can be considered equivalent to
\begin{align}
\label{eq:suff_decrease_condition_exact}
\J_h(\mu^{(k+1)})\leq \cJhatn^{(k)}(\mu_\textnormal{{AGC}}^{(k)}).
\end{align}
To avoid confusion, we note that \cite{paper1} has used \eqref{eq:suff_decrease_condition}, whereas \cite{paper2} uses the formulation in~\eqref{eq:suff_decrease_condition_exact}.
Due to the required enrichment of the model in order to have access to $\cJhatn^{(k+1)}$, we require to evaluate the expensive FOM model at $\mu^{(k+1)}$.
Thus, as in~\cite{QGVW2017,YM2013}, we consider cheaply computable sufficient and necessary conditions for~\eqref{eq:suff_decrease_condition}.
We use
\begin{equation}\label{eq:nec_cond}
	\cJhatn^{(k)}(\mu^{(k+1)})+\Delta_{\Jhat_\red^{(k)}}(\mu^{(k+1)})<\cJhatn^{(k)}(\mu^{(k)}_\textnormal{{AGC}})
\end{equation}
as sufficient condition for acceptance and
\begin{equation}\label{eq:suf_cond}
	\cJhatn^{(k)}(\mu^{(k+1)})-\Delta_{\Jhat_\red^{(k)}}(\mu^{(k+1)}) \leq \cJhatn^{(k)}(\mu^{(k)}_\textnormal{{AGC}})
\end{equation}
as necessary condition.
The TR-RB algorithm, then, accepts points that satisfy~\eqref{eq:nec_cond} and rejects any point which does not meet~\eqref{eq:suf_cond} and instead shrinks the TR-radius.
For more details and a derivation of these conditions, we refer to \cite{YM2013}.
If both conditions do not give an immediate decision, we instead check the expensive condition in~\eqref{eq:suff_decrease_condition}.
We enrich the RB spaces and continue with the algorithm if the point is accepted.

\subsubsection{Enlargement of the TR-radius}
\label{sec:enlarging}

To further enhance the algorithm, we consider a condition that allows enlarging the TR radius.
A drawback of the TR algorithm proposed in~\cite{QGVW2017} is that the TR radius may be significantly shrunk at the beginning, i.e.~when the TR model is poor in approximation.
Afterward, even if the RB space is enriched, i.e.~the approximation of the TR model function is improved, the TR radius is kept small.
Thus, one misses the local second-order rate of convergence of the BFGS or Newton method.
More precisely, if $\mu^{(k,\el)}$ is close to the locally optimal solution $\bar \mu^{(k)}$ of the TR sub-problem, we intend to make complete BFGS steps for faster convergence.
The possibility to enlarge the TR radius at each iteration also decreases the number of outer iterations needed to converge.
As a condition for enlarging the radius, we check whether the sufficient reduction predicted by the model function $\cJhatn^{(k)}$ is realized by the objective function, i.e., we check if
\begin{equation}
	\label{TR_act_decrease}
	\varrho^{(k)}:= \frac{\Jhat_h(\mu^{(k)})-\Jhat_h(\mu^{(k+1)})}{\cJhatn^{(k)}(\mu^{(k)})-\cJhatn^{(k)}(\mu^{(k+1)})}  \geq  \eta_\varrho
\end{equation}
for a tolerance $\eta_\varrho \in [3/4,1)$.
Note that $\varrho^{(k)}$ is the TR-RB version of $\rho^{(k)}$ in the BTR algorithm, which is commonly used as acceptance quantity.
In the former section, we explained that $\varrho^{(k)}$ is computationally costly because of the evaluation of the FOM cost functional $\Jhat_h$.
In contrast to the question of acceptance, we request~\eqref{TR_act_decrease} after the iterate is accepted, thus, after the RB space enrichment.
Hence, the quantities in the numerator of~\eqref{TR_act_decrease} are cheaply accessible since we have already solved the FOM to generate the new snapshots.

\subsubsection{Optional basis enrichment}
\label{sec:TRRB_optional_enrichment}

We also introduce the possibility of skipping the basis enrichment of the model if suitable conditions are met.
These conditions can also be used to accept the point $\mu^{(k+1)}$ since they directly imply the error-aware sufficient decrease condition~\eqref{eq:suff_decrease_condition} for the convergence of the method; cf.~\cite{YM2013} and \cref{sec:TRRB_convergence_analysis}.
At first, we define
\begin{equation}
	\label{eq:FOC_def}
	g_h(\mu)  := \|\mu-\Proj_{\Paramsad}(\mu-\nabla_\mu\Jhat_h(\mu))\|_2
\end{equation}
and analogously 
\[
g^{(k)}_\red(\mu) := \|\mu-\Proj_{\Paramsad}(\mu-\nabla_\mu\cJhatn^{(k)}(\mu))\|_2
\]
for all $\mu\in\Params$.
Then, the sufficient condition for skipping the enrichment at iteration $k$ reads as follows:
\begin{equation}
	\label{eq:skip_enrichment_condition}
	\begin{aligned}
		&&\textsf{Skip\_enrichment\_flag}(k) := \left(q^{(k)}(\mu^{(k+1)})\leq \beta_3\delta^{(k+1)}\right)\texttt{ and } \hspace{13mm} \\ &&  \left(\frac{\left|g_h(\mu^{(k+1)})-g_\red^{(k)}(\mu^{(k+1)})\right|}{g_\red^{(k)}(\mu^{(k+1)})}\leq \tau_g\right) \texttt{ and } \\ && \left(\frac{\|\nabla_\mu \Jhat_h(\mu^{(k+1)})-\nabla_\mu \Jhat^{(k)}_\red(\mu^{(k+1)})\|_2}{\|\nabla_\mu \Jhat_h(\mu^{(k+1)})\|_2 } \leq \min\{\tau_\textnormal{{grad}},\beta_3\delta^{(k+1)}\}\right),
	\end{aligned}
\end{equation} 
for given $\tau_g>0$ and $\tau_\textnormal{{grad}},\beta_3\in(0,1)$.
The first part of~\eqref{eq:skip_enrichment_condition} indicates how much the current RB model is trustworthy in the next iteration $k+1$, the second condition is to ensure the convergence of the algorithm (cf. \cref{thm:convergence_of_TR_2}) and the third one is to measure the RB accuracy in reconstructing the FOM gradient of $\Jhat_h$.
Note that these conditions require FOM quantities.
By using an unconditional basis enrichment, as discussed earlier, these quantities are accessible exactly because of the enrichment, therefore it appears contradictory to request them and then skip a basis update.
Here the focus is in fact not to avoid particular FOM solves priorly, but to exploit them in order to keep the dimension of the RB space small.
This is of significant importance when the TR-RB method takes many iterations (as seen in some examples in \cref{sec:TRRB_num_experiments}), where an unconditional enrichment in each iteration would lead to overfitted and too large RB spaces, slowing down the computation in the long run.
Another prominent example, where this feature is relevant, are applications such as PDE-constrained multi-objective optimization by scalarization methods {\cite{BBV2017,Ehrgott2005,IUV2017}}.
In there, many optimization problems have to be solved iteratively so that an efficient algorithm must keep the dimension of the RB space reasonably small. 
For a specific example where the above explained TR-RB method is used for multi-objective optimization problems, we make reference to \cite{banholzer_phd,banholzer2022trust}.

\subsubsection{Outer stopping criterion}
\label{sec:outer_stopping}
As mentioned before, after the (optional) enrichment of the RB space, we have access to the FOM function value $\Jhat_h(\mu^{(k+1)})$ as well as to the FOM gradient $\nabla_\mu \Jhat_h(\mu^{(k+1)})$.
This knowledge is used for the stopping criterion in the outer loop of the algorithm, which is then equivalent to the condition that would be used for a FOM-based algorithm, i.e. for every iterate $\mu^{(k+1)}$, we check
\begin{equation} \label{eq:TR_termination}
	g_h(\mu^{(k+1)}) = \|\mu^{(k+1)}-\Proj_{\Paramsad}(\mu^{(k+1)}-\nabla_\mu\Jhat_h(\mu^{(k+1)}))\|_2
	\leq \tau_{\textnormal{{FOC}}},
\end{equation}
with an appropriate tolerance $\tau_{\textnormal{{FOC}}} > 0$.

As pointed out in \cref{sec:TRRB_RB_estimates_param}, we can also enhance the algorithm with a suitable post-processing after~\eqref{eq:TR_termination} is met.
The reason is that the tolerance $\tau_{\textnormal{{FOC}}}$ is problem-dependent, and it is very likely to happen that the TR-RB terminates due to a wrongly chosen tolerance before the desired convergence is achieved.
Once the TR-RB algorithm has converged, we check if its solution is close enough to $\bar\mu_h$ by evaluating $\Delta_\mu$, given a desired tolerance $\tau_{\mu}>0$, i.e.
\begin{equation}
	\Delta_\mu(\mu^{(k+1)}) \leq \tau_{\mu}.
\end{equation}
If not, we decrease the stopping tolerance $\tau_{\textnormal{{FOC}}}$ and continue with the algorithm until \eqref{eq:TR_termination} is fulfilled again.
We call this approach \emph{parameter control}; cf. \cref{num_test:12params}.

\subsubsection{The algorithm in Pseudo-code}
We summarize the basic concept of the above-explained TR-RB algorithm in \cref{alg:Basic_TR-RBmethod}.
We mention that the originally proposed algorithms in~\cite{paper1} and~\cite{paper2} mainly differ in terms of the sub-problem solver in Line \ref{solve_sub_problem}, the possibility for optional enrichment in Line \ref{enrichment}, and the possibility to perform a parameter control in Line \ref{post_processing}.
For a detailed version of the exact algorithms that are used in~\cite{paper1} and~\cite{paper2}, we refer to \cref{alg:TR-RBmethod_paper1} and \cref{alg:TR-RBmethod_paper2}, respectively.

\begin{algorithm2e}[!h]
	\KwData{Initial TR radius $\delta^{(0)}$, TR shrinking factor $\beta_1\in (0,1)$, tolerance for enlarging the TR radius $\eta_\varrho\in [\frac{3}{4},1)$, initial parameter $\mu^{(0)}$, stopping tolerance for the sub-problem $\tau_\textnormal{{sub}}\ll 1$, stopping tolerance for the first-order critical condition $\tau_\textnormal{{FOC}}$ with $\tau_\textnormal{{sub}}\leq\tau_\textnormal{{FOC}}\ll 1$, safeguard for TR boundary $\beta_2\in(0,1)$, parameter control tolerance $\tau_{\mu}$.}
	Initialize RB model with $\mu^{(0)}$ and set $k=0$\; 
	\While
		{$\|\mu^{(k)}-\Proj_{\Paramsad}(\mu^{(k)}-\nabla_\mu \Jhat_h(\mu^{(k)}))\|_2 > \tau_{\textnormal{{FOC}}}$\label{stopping_condition}}{
		Compute $\mu^{(k+1)}$ as solution of~\eqref{TRRBsubprob} with termination criteria~\eqref{Termination_crit_sub-problem}\label{solve_sub_problem}\;
		\uIf{Sufficient decrease condition~\eqref{eq:suff_decrease_condition} is fulfilled} {
			Accept $\mu^{(k+1)}$, compute $\varrho^{(k)}$ from~\eqref{TR_act_decrease}
			and compute $g_h(\mu^{(k+1)})$\;
			Enrich the RB model at $\mu^{(k+1)}$ (optional enrichment possible)\label{enrichment}\;
			Enlarge the TR-radius if $\varrho^{(k)}\geq\eta_\varrho$\;
		}
		\Else {
			Reject $\mu^{(k+1)}$, shrink the TR radius $\delta^{(k+1)} = \beta_1\delta^{(k)}$ and go to \ref{solve_sub_problem}\;		
			}
		Set $k=k+1$\;
	}
	Optional post-processing: \textbf{if} $\Delta_\mu \geq \tau_{\mu}$ \textbf{then} Decrease $\tau_\textnormal{{FOC}}$ and go to Line \ref{stopping_condition}\label{post_processing}\;
	\caption{Basic TR-RB algorithm} 
	\label{alg:Basic_TR-RBmethod}
\end{algorithm2e}

\begin{algorithm2e}[!h]
	\KwData{Initial TR radius $\delta^{(0)}$, TR shrinking factor $\beta_1\in (0,1)$, tolerance for enlarging the TR radius $\eta_\varrho\in [\frac{3}{4},1)$, initial parameter $\mu^{(0)}$, stopping tolerance for the sub-problem $\tau_\textnormal{{sub}}\ll 1$, stopping tolerance for the first-order critical condition $\tau_\textnormal{{FOC}}$ with $\tau_\textnormal{{sub}}\leq\tau_\textnormal{{FOC}}\ll 1$, safeguard for TR boundary $\beta_2\in(0,1)$.}
	
	Set $k=0$ and \textsf{Loop\_flag}$=$\textsf{True}\; 
	\While{
		{\normalfont\textsf{Loop\_flag}}}{
		Compute $\mu^{(k+1)}$ as solution of \eqref{TRRBsubprob} with termination criteria \eqref{Termination_crit_sub-problem}\label{TRRB-optstep}\;
		\uIf{\label{Suff_condition_TRRB}$\cJhatn^{(k)}(\mu^{(k+1)})+\Delta_{\Jhat_\red^{(k)}}(\mu^{(k+1)})<\cJhatn^{(k)}(\mu^{(k)}_\textnormal{{AGC}})$} {
			Accept $\mu^{(k+1)}$,
			update the RB model at $\mu^{(k+1)}$ and compute $\varrho^{(k)}$ from \eqref{TR_act_decrease}\; 
			\eIf {$\varrho^{(k)}\geq\eta_\varrho$}  {
				Enlarge the TR radius $\delta^{(k+1)} = \beta_1^{-1}\delta^{(k)}$\;	
			}
			{
				Set $\delta^{(k+1)}=\delta^{(k)}$\;
			}
			
		}
		\uElseIf {\label{Nec_condition_TRRB}$\cJhatn^{(k)}(\mu^{(k+1)})-\Delta_{\Jhat_\red^{(k)}}(\mu^{(k+1)})> \cJhatn^{(k)}(\mu^{(k)}_\textnormal{{AGC}})$} {
			Reject $\mu^{(k+1)}$, shrink the TR radius $\delta^{(k+1)} = \beta_1\delta^{(k)}$ and go to \ref{TRRB-optstep}\;
		}
		\Else {
			Update the RB model at $\mu^{(k+1)}$ and compute $\varrho^{(k)}$ from \eqref{TR_act_decrease}\; 
			\eIf {$\cJhatn^{(k+1)}(\mu^{(k+1)})\leq \cJhatn^{(k)}(\mu^{(k)}_\textnormal{{AGC}})$} {
				Accept $\mu^{(k+1)}$\;	
				\eIf {$\varrho^{(k)}\geq\eta_\varrho$} {
					Enlarge the TR radius $\delta^{(k+1)} = \beta_1^{-1}\delta^{(k)}$\;	
				}
				{
					Set $\delta^{(k+1)}=\delta^{(k)}$\;
				}		
			}{
				Reject $\mu^{(k+1)}$, shrink the TR radius $\delta^{(k+1)} = \beta_1\delta^{(k)}$ and go to \ref{TRRB-optstep}\;		
			}
		}
		\If { $\|\mu^{(k+1)}-\Proj_{\Params_{ad}}(\mu^{(k+1)}-\nabla_\mu \Jhat_h(\mu^{(k+1)}))\|_2\leq \tau_{\textnormal{{FOC}}}$ } {
			Set \textsf{Loop\_flag}$=$\textsf{False}\;
		} 
		Set $k=k+1$\;
		
	}
	\caption{TR-RB algorithm as used in \cite{paper1}} 
	\label{alg:TR-RBmethod_paper1}
\end{algorithm2e}

\begin{algorithm2e}
	
	Initialize the ROM at $\mu^{(0)}$, set $k=0$ and \textsf{Loop\_flag}$=$\textsf{True}\; 
	\While{
		{\normalfont\textsf{Loop\_flag}}}{
		Compute the AGC point $\mu^{(k)}_\textnormal{AGC}$\;
		Compute $\mu^{(k+1)}$ as solution of \eqref{TRRBsubprob} with stopping criteria \eqref{Termination_crit_sub-problem}\label{TRRB-optstep_2}\;
		\uIf{\label{Suff_condition_TRRB_2}$\cJhatn^{(k)}(\mu^{(k+1)})+\Delta_{\Jhat_\red^{(k)}}(\mu^{(k+1)})<\cJhatn^{(k)}(\mu^{(k)}_\textnormal{{AGC}})$} {
			Accept $\mu^{(k+1)}$, set $\delta^{(k+1)}=\delta^{(k)}$, compute $\varrho^{(k)}$ and $g_h(\mu^{(k+1)})$\label{AcceptingIterateSufficientCondition}\;
			\eIf { $g_h(\mu^{(k+1)}) \leq \tau_{\textnormal{{FOC}}}$ } {
				Set \textsf{Loop\_flag}$=$\textsf{False}\;
			} 
			{
				\If {$\varrho^{(k)}\geq\eta_\varrho$}  {
					Enlarge the TR radius $\delta^{(k+1)} = \beta_1^{-1}\delta^{(k)}$\;	
				}
				\If {{\normalfont\textbf{not}} {\normalfont\textsf{Skip\_enrichment\_flag}$(k)$}}
				{
					Update the RB model at $\mu^{(k+1)}$ 
					\label{UpdateRBModelSufficientCondition}\;
				}
			}
		}
		\uElseIf {\label{Nec_condition_TRRB_2}$\cJhatn^{(k)}(\mu^{(k+1)})-\Delta_{\Jhat_\red^{(k)}}(\mu^{(k+1)})> \cJhatn^{(k)}(\mu^{(k)}_\textnormal{{AGC}})$} {
			\If { $\beta_1\delta^{(k)} \leq \delta_{\textnormal{{min}}}$ \textbf{ or } \normalfont\textsf{Skip\_enrichment\_flag}$(k-1)$ \label{forced_enrichment}} {
				Update the RB model at $\mu^{(k+1)}$\; 
			}
			Reject $\mu^{(k+1)}$, shrink the radius $\delta^{(k+1)} = \beta_1\delta^{(k)}$ and go to \ref{TRRB-optstep}\; 
		}
		\Else {
			Compute $\Jhat_h(\mu^{(k+1)})$, $g_h(\mu^{(k+1)})$, 
			$\varrho^{(k)}$ 
			and set $\delta^{(k+1)} = \beta_1^{-1}\delta^{(k)}$\; 
			\eIf { $g_h(\mu^{(k+1)})\leq \tau_{\textnormal{{FOC}}}$ } {
				Set \textsf{Loop\_flag}$=$\textsf{False}\;
			} {
				\uIf { \label{AcceptingIterateNoUpdateByExactComputation_condition}{\normalfont\textsf{Skip\_enrichment\_flag}$(k)$} \textbf{ and } $\varrho^{(k)}\geq \eta_\varrho$}
				{ Accept $\mu^{(k+1)}$\label{AcceptingIterateNoUpdateByExactComputation}\;}
				\uElseIf {$\Jhat_h(\mu^{(k+1)}) \leq \cJhatn^{(k)}(\mu^{(k)}_\textnormal{{AGC}})$} {
					Accept $\mu^{(k+1)}$ and update the RB model 
					\label{AcceptingIterateAndUpdateByExactComputation}\;	
					\If {$\varrho^{(k)}<\eta_\varrho$} {
						Set $\delta^{(k+1)}=\delta^{(k)}$\;
					}		
				}
				\Else{
					\If { $\beta_1\delta^{(k)} \leq \delta_{\textnormal{{min}}}$ \textbf{ or } \normalfont\textsf{Skip\_enrichment\_flag}$(k-1)$ \label{forced_enrichment_2}} {
						Update the RB model at $\mu^{(k+1)}$\; 
					}
					Reject $\mu^{(k+1)}$, set $\delta^{(k+1)} = \beta_1\delta^{(k)}$ and go to \ref{TRRB-optstep}\;		
				}
			}
		}
		Set $k=k+1$\;		
	}
	\caption{{TR-RB algorithm as used in \cite{paper2}}}
	\label{alg:TR-RBmethod_paper2}
\end{algorithm2e}

\clearpage

\subsection{Convergence study}
\label{sec:TRRB_convergence_analysis}

To strengthen the practicability of the above-presented TR-RB algorithm, we require a corresponding convergence result that, in contrast to the results in~\cite{YM2013}, accounts for, first, the inequality constraints on the parameter space, and, second, for the optional basis enrichment.
The following convergence study has been carried out in~\cite{paper2,paper1}, where the convergence study in~\cite{paper2} is an improved version of~\cite{paper1}.
We particularly mention that the main contributions in the convergence proof have been conducted by Luca Mechelli (in \cite{paper2,paper1}) and Stefan Banholzer (in \cite{paper2}).
 
We require two basic assumptions for both convergence studies.
In order to guarantee the well-posedness (because of~\eqref{TR_radius_condition}) and the convergence of the method, we enforce:
\begin{assumption}
	\label{asmpt:bound_J}
	The objective functional $\J(u,\mu)$ is strictly positive for all $u\in V$ and all parameters $\mu\in\Params$. 
\end{assumption}
Note that this assumption is not too restrictive since the boundedness from below is a standard assumption in optimization for guaranteeing a solution for the minimization problem.
If a global lower bound for the cost functional is also known, one can add a sufficiently large constant to the objective functional without changing the position of its local minima and maxima.

As explained in \cref{sec:acceptance_of_TR}, outer iterates are only accepted if the error-aware sufficient decrease condition \eqref{eq:suff_decrease_condition} is satisfied.
If the iterate, instead, is rejected, we shrink the TR-radius~$\delta^{(k)}$.
One may be concerned that the TR-RB algorithm may be trapped in an infinite loop where every computed point is rejected, and the TR radius is shrunk all time.
We point out that this never happened in our numerical tests. Anyway, we consider a safety termination criterion, which is triggered when the TR-radius is smaller than the double machine precision.
To prove convergence, in what follows, we then assume that this situation can not occur.
\begin{assumption}
	\label{asmpt:rejection}
	For each $k\geq 0$, there exists a radius $\widetilde{\delta}^{(k)}>\tau_{{\textnormal{mac}}}$ for which a solution of~\eqref{TRRBsubprob} exists such that~\eqref{eq:suff_decrease_condition} is verified, where $\tau_{{\textnormal{mac}}}= 2.22\cdot10^{-16}$ is the double machine precision.
\end{assumption}

We now present the original convergence study that uses properties of the (unconditionally enriched) RB space.
Subsequently, we state the improved version that constitutes a more general result from an infinite-dimensional perspective and tackles optional enrichment.

\subsubsection{RB-based convergence study with unconditional enrichment}
\label{convergence_study_1}

We are concerned with the convergence proof of the TR-RB algorithm as proposed in \cite{paper1}; see \cref{alg:TR-RBmethod_paper1} for a detailed pseudo-code.
In what follows, \cref{lemma:AGC_search} gives a statement about the successful search of the AGC point and \cref{thm:convergence_of_TR} contains the convergence result.
 
\begin{lemma}
	\label{lemma:AGC_search}
	The search of the AGC point defined in \cref{def:AGC} takes finitely many iterations at each step $k$ of the TR-RB Algorithm.
\end{lemma}
\begin{proof}
	We want to prove that there exists an index $j_c^{(k)}<\infty$ for each $k\geq 0$, for which $\mu_\textnormal{{AGC}}^{(k)}= \mu^{(k,0)}(j^{(k)}_c)$ satisfies~\eqref{Arm_and_TRcond} for $\el=0$.
	From~\cite[Theorem~5.4.5]{kelley} (and the subsequent discussion), we conclude that for all $k\in\mathbb{N}$ there exists a strictly positive index $j^{(k)}_1\in\mathbb{N}$, such that $\mu^{(k,0)}(j)$ satisfies~\eqref{Armijo} for $j\geq j^{(k)}_1$ and $\el=0$.
	If $k=0$, by construction, we have that $\Delta_{\Jhat_\red^{(0)}}(\mu^{(0)}) = 0$.
	Therefore, there exists a sufficiently large (but finite) index $j^{(0)}_2\in\mathbb{N}$, such that $\mu^{(0,0)}(j)$ satisfies~\eqref{TR_radius_condition} for all $j\geq j^{(0)}_2$ and $\el=0$.
	This can be concluded from the continuity w.r.t.~$\mu$ of the error estimator $\Delta_{\Jhat_\red^{(k)}}(\mu)$ (cf. \cref{continuity_of_estimator}) and of the cost functional $\cJhatn^{(k)}(\mu)$ for all $k\in\mathbb{N}$.
	We obtain that there exists $j^{(0)}_c=\max(j^{(0)}_1,j^{(0)}_2)<\infty$, for which $\mu^{(0,0)}(j)$ satisfies~\eqref{Arm_and_TRcond} for $\el=0$. If $k\geq 1$, since the model has been enriched, i.e.~$\Delta_{\Jhat_\red^{(k)}}(\mu^{(k)}) = 0$, we can show the claim arguing as we did for $k=0$.
	Note that we increase the iteration counter only when $\mu^{(k)}$ is accepted at iteration $k-1$ and, thus, when the RB model is enriched at this parameter.
\end{proof}

The proof of \cref{lemma:AGC_search} takes significant advantage of the unconditional RB enrichment as well as the continuity of the residual-based error estimator, stated in \cref{continuity_of_estimator}.
Therefore, at least in the first place, the lemma can not be generalized for an arbitrary surrogate model or optional enrichment strategies.
Using \cref{lemma:AGC_search}, we proof the following theorem:

\begin{theorem}
	\label{thm:convergence_of_TR}
	Every accumulation point~$\bar\mu$ of the sequence $\{\mu^{(k)}\}_{k\in\mathbb{N}}\subset \Params$ generated by the TR-RB algorithm is an approximate first-order critical point for $\Jhat_h$ (up to the chosen tolerance $\tau_\textnormal{{sub}}$), i.e., it holds
	\begin{equation}
		\label{eq:TR-RB-FOC}
		\|\bar \mu-\Proj_{\Paramsad}(\bar \mu-\nabla_\mu \Jhat_h(\bar \mu))\|_2 \leq \tau_\textnormal{{sub}}.
	\end{equation}
\end{theorem}
\begin{proof}
	The set $\Paramsad \subseteq \mathbb{R}^P$ is compact.
	Thus, there exists a sequence of indices $\left\{k_i\right\}_{i\in\mathbb{N}}$, such that the sub-sequence $\{\mu^{(k_i)}\}_{i\in\mathbb{N}}$ converges to a point $\bar \mu\in \Params$.
	It remains to show that $\bar \mu$ is an approximate first-order critical point.
	At first, note that once the RB space is enriched at a point $\mu^{(k)}$, we have $\Delta_{\Jhat_\red^{(k)}}(\mu^{(k)})= 0$.
	Hence, also $q^{(k)}(\mu^{(k)}) = 0$ holds with $q^{(k)}$ defined in~\eqref{TR_radius_condition}.
	Note also that $V_h$ is a finite-dimensional space.
	This implies that, at most after  $\dim V_h\leq I<+\infty$ enrichment steps, the RB approximation error and the a posteriori error estimator are zero for each $\mu\in\Params$.
	In particular, it holds $q^{(k_i)}(\mu) = 0$ for all $\mu\in\Params$ and $i\geq I$.
	For this reason, the stopping criterium~\eqref{Cut_of_TR} is not triggered.	
	Hence, we have proved that each $\mu^{(k_{i+1})}$ is an approximate first-order critical point for $\cJhatn^{(k_{i+1}-1)}$ (up to the chosen tolerance $\tau_\textnormal{{sub}}$) for all $i\geq I$, which yields to
	\begin{align*}
		\|\mu^{(k_{i+1})}-\Proj_{\Paramsad}(\mu^{(k_{i+1})}-\nabla_\mu\cJhatn^{(k_{i+1}-1)}(\mu^{(k_{i+1})}))\|_2 \leq \tau_\textnormal{{sub}}, && \textnormal{for all } i\geq I.
	\end{align*}
	Moreover, taking again into account the RB method properties and the fact that $V_h$ is a finite-dimensional space, there exists a constant $I_\nabla>0$ sufficiently large, such that $\nabla_\mu \cJhatn^{(k_i)}(\mu)= \nabla_\mu \Jhat_h(\mu)+\epsilon^{(k_i)}$ for all $\mu$ in a neighborhood of $\bar\mu$ and for $i\geq I_{\nabla}$, with $\epsilon^{(k_i)}\to 0$ as $i\to\infty$. 
	Thus, exploiting the continuity of the projection operator and assuming $i\geq \max(I,I_\nabla)$, we have that
	\[
	\begin{aligned}
		\tau_\textnormal{{sub}} & \geq \|\mu^{(k_{i+1})}-\Proj_{\Paramsad}(\mu^{(k_{i+1})}-\nabla_\mu\cJhatn^{(k_{i+1}-1)}(\mu^{(k_{i+1})}))\|_2 \\
		& = \|\mu^{(k_{i+1})}-\Proj_{\Paramsad}(\mu^{(k_{i+1})}-\nabla_\mu\Jhat_h(\mu^{(k_{i+1})})+\epsilon^{(k_{i+1}-1)})\|_2 \to \|\bar\mu - \Proj_{\Paramsad}(\bar\mu-\nabla_\mu\Jhat_h(\bar \mu))\|_2.
	\end{aligned}
	\]
	Hence, the accumulation point $\bar \mu$ is an approximate first-order critical point (up to the tolerance $\tau_\textnormal{{sub}}$).
\end{proof}

Again, the proof of \cref{thm:convergence_of_TR} relies on the fact that $V_h$ is a finite-dimensional space and that at most after $\dim V_h\leq I<+\infty$ iterations, the RB space is exact.
From a practical point of view, having an RB space of the same dimension as the full-order model will not give any speed-up.
However, with respect to the exponential decay of the Kolmogorov n-width (cf.~\cref{sec:background_basis_gen_RB}), we do not expect such a scenario.
Furthermore, since the optimization sequence will accumulate fast around $\bar\mu$, we expect that the RB model will be accurate enough not to trigger~\eqref{Cut_of_TR} implying the TR method to convergence in a much smaller number of iterations.
The numerical tests in \cref{sec:TRRB_num_experiments} confirm this expectation.

\begin{remark}
	\label{rmk:local_minimum_TR}
	What remains to prove is that $\bar \mu$ is a local minimum of $\Jhat_h$ (or rather a sufficiently close approximation of a local minimum).
	Exploiting the sufficient decrease condition, one can easily show by contradiction that $\bar \mu$ is not a maximum of $\Jhat_h$ but
	it can still be a saddle point.
	In the numerical experiments, to verify that the computed point $\bar \mu$ is a local minimum, we employ the second-order sufficient optimality conditions after the algorithm terminates.
\end{remark}

\subsubsection{Improved convergence study with optional enrichment}
\label{convergence_study_2}

In this section, we improve the convergence analysis of \cref{convergence_study_1} to an infinite-dimensional perspective and enhance it to tackle the discussed optional enrichment strategy.
A detailed description of the TR-RB algorithm that we prove convergence for can be found in \cref{alg:TR-RBmethod_paper2}.

First, to stabilize the algorithm's convergence, we account for an issue that might appear due to skipping an RB basis update in the following situation. Let, at iteration $k-1$, the optimization sub-problem terminate for~\eqref{Cut_of_TR}, the point be accepted, the enrichment be skipped, the radius be enlarged to $\delta^{(k)}=\beta_1^{-1}\delta^{(k-1)}$, but, at iteration $k$, the point $\mu^{(k+1)}$ be rejected, implying to shrink the radius to the old value $\delta^{(k-1)}$.
Suppose the model is also not updated at iteration $k$. In that case, we are solving the same sub-problem again at the next iteration, starting at a point that was already triggering~\eqref{Cut_of_TR}.
Therefore our step would be to compute only the AGC point.
Although the method will converge anyway, this situation might repeat several times before we escape this ``problematic'' region, resulting in a waste of computational time. 
Therefore, we impose an enrichment of the RB model, when the radius is shrunk at iteration $k$ and we skipped the basis update at iteration $k-1$; cf.~Step~\ref{forced_enrichment} and Step~\ref{forced_enrichment_2} of Algorithm~\ref{alg:TR-RBmethod_paper2}.

Further, to improve the convergence results, we required an additional assumption compared to \cref{convergence_study_1}, namely locally Lipschitz-continuous second derivatives from \cref{asmpt:differentiability}.
The assumption is formulated in the following.
\begin{assumption}
	\label{asmpt:Lip_cont}
	The ROM gradient $\nabla_\mu \Jhat_\red^{(k)}$ is uniformly Lipschitz-continuous, i.e. there exists a constant $C_L>0$ independent of $k$, such that
	\[
	\|\nabla_\mu \Jhat_\red^{(k)}(\mu)-\nabla_\mu\Jhat_\red^{(k)}(\nu)\|_2\leq C_L\|\mu-\nu\|_2
	\]
	holds for all $\mu,\nu\in\Params$ and all $k\in\mathbb{N}$.
	Similarly, the second derivatives of $\Jhat_\red^{(k)}$ are locally Lipschitz-continuous.
\end{assumption}
This assumption restricts the set of cost functionals. Nevertheless, it guarantees a locally faster convergence behavior for this class.
\cref{alg:TR-RBmethod_paper2} is still applicable to the general class of quadratic functionals -- and also converges in this case.
We remark that \cref{asmpt:Lip_cont} is needed for proving convergence of the method in an infinite-dimensional perspective; cf. \cref{thm:convergence_of_TR_2}.
In the specific case of the RB model function, it is possible to show that this is satisfied.
The proof follows from the fact that the RB model will exactly approximate the cost functional~$\Jhat_h$ after a finite number of updates.
As a direct consequence of \cref{asmpt:Lip_cont}, we have the following result:
\begin{corollary}
	\label{cor:bound_grad}
	Let Assumption~{\ref{asmpt:Lip_cont}} be satisfied.
	Then there exists a constant $C>0$, such that for any $k\in\mathbb{N}$ it holds
	\[
	\|\nabla_\mu \Jhat_\red^{(k)}(\mu^{(k)})\|_2\leq C.
	\]
\end{corollary}
\noindent Furthermore, the following property of the projection operator $\Proj_{\Paramsad}$ holds:
\begin{lemma} \label{lemma:PropertyOfPParams}
	Let $\mu\in\Params$ and $d \in \R^P$ be arbitrary. Then, it holds
	\begin{equation} 
		\label{eq:PropertyOfPParams}
		\left\Vert \mu - \Proj_{\Paramsad} \left( \mu - t d \right) \right\Vert_2 \geq t \left\Vert \mu - \Proj_{\Paramsad} \left( \mu - d \right) \right\Vert_2
	\end{equation}
	for all $t \in [0,1]$.
\end{lemma}

\begin{proof}
	Let $\mu \in \Params$, $d \in \R^P$ and $t \in [0,1]$ be arbitrary.
	We prove the statement by showing that 
	\begin{equation} \label{eq:PropertyOfPParams:Componentwise}
		\left| \mu_i - \left( \Proj_{\Paramsad} \left( \mu - t d \right) \right)_i \right| \geq t \left|  \mu_i - \left( \Proj_{\Paramsad} \left( \mu - d \right) \right)_i \right|
	\end{equation}
	holds for all components $i= 1,\ldots,P$.
	Let $i \in \{1,\ldots,P\}$ be arbitrary. \\ 
	\emph{Case (1)}: $\left(  \Proj_{\Paramsad} \left( \mu - t d \right) \right)_i \in \{ (\mu_\mathsf{a})_i,(\mu_\mathsf{b})_i \}$. \\
	It clearly holds $\left(  \Proj_{\Paramsad} \left( \mu - t d \right) \right)_i = \left(  \Proj_{\Paramsad} \left( \mu - d \right) \right)_i$, so that
	\[
	\left| \mu_i - \left(  \Proj_{\Paramsad} \left( \mu - t d \right) \right)_i \right| = \left| \mu - \left(  \Proj_{\Paramsad} \left( \mu - d \right) \right)_i \right| \geq t \left| \mu - \left(  \Proj_{\Paramsad} \left( \mu - d \right) \right)_i \right|.
	\]
	\emph{Case (2a)}: $\left(  \Proj_{\Paramsad} \left( \mu - t d \right) \right)_i \in ((\mu_\mathsf{a})_i,(\mu_\mathsf{b})_i)$ and $\left(  \Proj_{\Paramsad} \left( \mu - d \right) \right)_i \in ((\mu_\mathsf{a})_i,(\mu_\mathsf{b})_i)$.
	We can conclude
	\[
	\left| \mu_i - \left(  \Proj_{\Paramsad} \left( \mu - t d \right) \right)_i \right| = t \left| d_i \right| = t \left| \mu_i - \left(  \Proj_{\Paramsad} \left( \mu - d \right) \right)_i \right|,
	\]
	which is what we have to show. \\
	\emph{Case (2b)}: $\left(  \Proj_{\Paramsad} \left( \mu - t d \right) \right)_i \in ((\mu_\mathsf{a})_i,(\mu_\mathsf{b})_i)$ and 
	$\left(  \Proj_{\Paramsad} \left( \mu - d \right) \right)_i \kern-0.1em\in\kern-0.1em \{ (\mu_\mathsf{a})_i,(\mu_\mathsf{b})_i \}$.
	We define $\tilde{t} := \frac{\mu_i -(\mu_{\mathsf{a},\mathsf{b}})_i}{d_i}$.
	Note that $t<\tilde{t}\leq 1$. Then, it holds 
	\begin{align*}
		\left| \mu_i - \left( \Proj_{\Paramsad} \left( \mu - t d \right) \right)_i \right| & = t \left| d_i \right| = \frac{t}{\tilde{t}} \left| \mu_i - \left( \Proj_{\Paramsad} \left( \mu - \tilde{t} d \right) \right)_i \right| \\
		& = \frac{t}{\tilde{t}} \left| \mu_i - \left( \Proj_{\Paramsad} \left( \mu - d \right) \right)_i \right| \geq t \left| \mu_i - \left( \Proj_{\Paramsad} \left( \mu - d \right) \right)_i \right|.
	\end{align*} 
	Consequently, in all cases, for any component $i \in \{1,\ldots,P\}$ the inequality~\eqref{eq:PropertyOfPParams:Componentwise} holds.
	Now, it can be concluded that~\eqref{eq:PropertyOfPParams} holds as well.
\end{proof} 
\noindent The next lemma is also needed to show convergence of the algorithm.
\begin{lemma}
	\label{lemma:BothRBErrorsSmaller2}
	For every iterate $\mu^{(k)}$, $k \in \N$, of \cref{alg:TR-RBmethod_paper2}, it holds
	\begin{equation}
		q^{(k)}(\mu^{(k)}) \leq \beta_3\delta^{(k)} \quad \text{ and } \quad \frac{\left| g_\red^{(k)}(\mu^{(k)}) - g_h(\mu^{(k)}) \right|}{g_\red^{(k)}(\mu^{(k)})} \leq \tau_g. \label{eq:BothRBErrorsSmaller2}
	\end{equation}
\end{lemma}

\begin{proof}
	We show this statement by induction over $k \in \mathbb{N}$.
	For $k = 0$ we trivially have 
	\[
	q^{(0)}(\mu^{(0)}) = 0 \quad \text{ and } \quad \frac{\left| g_\red^{(0)}(\mu^{(0)}) - g_h(\mu^{(0)}) \right|}{g_\red^{(0)}(\mu^{(0)})} = 0
	\]
	since the RB model was constructed at $\mu^{(0)}$. \\
	Now, assume that~\eqref{eq:BothRBErrorsSmaller2} is satisfied for all $1 \leq l \leq k$ for some $k \in \N$ and let $\mu^{(k+1)}$ be the new accepted iterate.
	We consider three cases:
	\begin{enumerate}
		\item $\mu^{(k+1)}$ is accepted in line~\ref{AcceptingIterateSufficientCondition}: \\
		Then, the RB model is updated in line~\ref{UpdateRBModelSufficientCondition}, if
		\[
		\begin{aligned}
			& q^{(k)}(\mu^{(k+1)}) > \beta_3\delta^{(k+1)} \quad \textnormal{ or } \quad \frac{\left| g_\red^{(k)}(\mu^{(k+1)}) - g_h(\mu^{(k+1)}) \right|}{g_\red^{(k)}(\mu^{(k+1)})} > \tau_g \\ & \textnormal{ or } \left(\frac{\|\nabla_\mu \Jhat_h(\mu^{(k+1)})-\nabla_\mu \Jhat^{(k)}_\red(\mu^{(k+1)})\|_2}{\|\nabla_\mu \Jhat_h(\mu^{(k+1)})\|_2 }> \min\{\tau_\textnormal{{grad}},\beta_3\delta^{(k+1)}\}\right).
		\end{aligned}
		\]
		On the one hand, if the RB model is not updated in line~\ref{UpdateRBModelSufficientCondition}, it holds 
		\[
		q^{(k+1)}(\mu^{(k+1)}) = q^{(k)}(\mu^{(k+1)}) \leq \beta_3\delta^{(k+1)}
		\] 
		and
		\[
		\frac{\left| g_\red^{(k+1)}(\mu^{(k+1)}) - g_h(\mu^{(k+1)}) \right|}{g_\red^{(k+1)}(\mu^{(k+1)})} = \frac{\left| g_\red^{(k)}(\mu^{(k+1)}) - g_h(\mu^{(k+1)}) \right|}{g_\red^{(k)}(\mu^{(k+1)})} \leq \tau_g.
		\]
		If instead, the RB model is updated, we have
		\[
		q^{(k+1)}(\mu^{(k+1)}) = 0 \quad \text{ and } \quad \frac{\left| g_\red^{(k+1)}(\mu^{(k+1)}) - g_h(\mu^{(k+1)}) \right|}{g_\red^{(k+1)}(\mu^{(k+1)})} = 0,
		\]
		so that the claim follows in both cases.
		\item $\mu^{(k+1)}$ is accepted in line~\ref{AcceptingIterateNoUpdateByExactComputation}: \\
		In this case the RB model is not updated. Thus, we can directly conclude from the previous if-condition in line~\ref{AcceptingIterateNoUpdateByExactComputation_condition} and the enlarged TR radius that
		\[
		q^{(k+1)}(\mu^{(k+1)}) = q^{(k)}(\mu^{(k+1)}) \leq \beta_3 \beta_1^{-1} \delta^{(k)} = \beta_3\delta^{(k+1)}
		\] 
		and
		\[
		\frac{\left| g_\red^{(k+1)}(\mu^{(k+1)}) - g_h(\mu^{(k+1)}) \right|}{g_\red^{(k+1)}(\mu^{(k+1)})} = \frac{\left| g_\red^{(k)}(\mu^{(k+1)}) - g_h(\mu^{(k+1)}) \right|}{g_\red^{(k)}(\mu^{(k+1)})} \leq \tau_g.
		\]
		Hence, the claim holds also in this case.
		\item $\mu^{(k+1)}$ is accepted in line~\ref{AcceptingIterateAndUpdateByExactComputation}: \\
		Then, the RB model is updated at $\mu^{(k+1)}$, so that we have 
		\[
		q^{(k+1)}(\mu^{(k+1)}) = 0 \quad \text{ and } \quad \frac{\left| g_\red^{(k+1)}(\mu^{(k+1)}) - g_h(\mu^{(k+1)}) \right|}{g_\red^{(k+1)}(\mu^{(k+1)})} = 0.
		\]
		Thus, the claim holds trivially.
		
	\end{enumerate}
	In total, we have shown the claim for every possible case, which concludes the proof.
\end{proof}
We continue the convergence analysis by showing a result about the AGC point $\mu_\textnormal{{AGC}}^{(k)}$, serving as a replacement for \cref{lemma:AGC_search} of the original convergence analysis.
First, we recall the following results from~\cite[Corollary 5.4.4]{kelley}:
\begin{lemma}
	\label{lemma:kelley}
	For all $j,k\in\mathbb{N}$ and $\kappa\in(0,1)$, we have
	\begin{align*}
		\|\mu^{(k)}-\Proj_{\Paramsad}(\mu^{(k)}-&\kappa^j\nabla_\mu\Jhat^{(k)}_\red(\mu^{(k)}))\|^2_2 \\ &\leq \kappa^j \nabla_\mu\Jhat^{(k)}_\red(\mu^{(k)})\cdot(\mu^{(k)}-\Proj_{\Paramsad}(\mu^{(k)}-\kappa^j\nabla_\mu\Jhat^{(k)}_\red(\mu^{(k)}))).
	\end{align*}
\end{lemma}
\noindent To proceed, we assume that the error indicator $q^{(k)}$ in the TR condition~\eqref{TR_radius_condition} is uniformly continuous.
\begin{assumption}
	\label{asmpt:unif_cont_q}
	The function $q^{(k)}:\Params\to \mathbb{R}$ defined in~\eqref{TR_radius_condition} is uniformly continuous in $\Params$, i.e.
	\[
	\textnormal{for all } \varepsilon>0: \exists\eta = \eta(\varepsilon)>0: \textnormal{for all }\mu,\nu\in\Params \quad \|\mu-\nu\|_2<\eta \Rightarrow \left|q^{(k)}(\mu)-q^{(k)}(\nu)\right|<\varepsilon.
	\]
\end{assumption}
\noindent
This assumption is needed for proving the convergence of the method in an infinite-dimensional perspective.
In the case of the RB model, since $\Paramsad$ is compact, one can apply the Heine-Cantor theorem~\cite{Rudin1976} to show that $q^{(k)}$ is uniformly continuous for each $k\in\mathbb{N}$.
Then, the independence from $k$ follows from the fact that the RB model approximation is exact (after a sufficient number of enrichments) and $q^{(k)}=q^{(k+1)}$ when the enrichment is not performed.
Finally, the next result gives a lower and upper bound for the line-search of the AGC point.
This is important because it shows that, at each iteration $k$, the ACG point can be computed in finitely many line-search steps; cf. \cref{lemma:AGC_search}.
\begin{theorem}
	\label{thm:linesearch-agc}
	Let $\mu^{(k)}(j):= \Proj_{\Paramsad}(\mu^{(k)}-\kappa^j\nabla_\mu\Jhat^{(k)}_\red(\mu^{(k)}))$ for $j\in\mathbb{N}$.
	Then, we have that $\mu^{(k)}(j)$ satisfies~\eqref{Arm_and_TRcond} for all
	\begin{equation}
		\label{j_low_bound}
		j \geq \log_\kappa\left(\min\left\{\frac{2(1-\kappa_{\mathsf{arm}})}{C_L},\frac{\eta((1-\beta_3)\tau_{{\textnormal{mac}}})}{C}\right\}\right),
	\end{equation}
	where $\kappa\in(0,1)$ is the backtracking constant introduced in~\eqref{eq:General_Opt_Step_point} and $C_L$, $C$ and $\eta$ are introduced in Assumption~{\ref{asmpt:Lip_cont}}, Corollary~{\ref{cor:bound_grad}} and Assumption~{\ref{asmpt:unif_cont_q}}, respectively.
	Furthermore, for the step-length of the AGC points, it holds
	\begin{equation}
		\label{j_up_bound}
		j^{(k)}_c \leq \log_\kappa\left(\min\left\{\frac{2(1-\kappa_{\mathsf{arm}})\kappa}{C_L},\frac{\eta((1-\beta_3)\tau_{{\textnormal{mac}}})\kappa}{C}\right\}\right).
	\end{equation}
\end{theorem}
\begin{proof}
	We only need to prove~\eqref{j_low_bound} since~\eqref{j_up_bound} is a direct consequence of it.
	Let $j$ satisfying~\eqref{j_low_bound} be arbitrary and consider $y:= \mu^{(k)}-\mu^{(k)}(j)$. Then, it holds
	\begin{equation}
		\label{1step_AGC_thm}
		\begin{aligned}
			\Jhat_\red^{(k)}(\mu^{(k)})-&\Jhat_\red^{(k)}(\mu^{(k)}(j)) = -\int_0^1 \frac{\mathrm d}{\mathrm d s} \Jhat_\red^{(k)}(\mu^{(k)}-sy)\,\mathrm{d} s \\ & = \int_0^1  \nabla_\mu\Jhat^{(k)}_\red(\mu^{(k)}-sy)\cdot y \,\mathrm{d} s \\
			& =  \nabla_\mu\Jhat(\mu^{(k)})\cdot y +\int_0^1\left(\nabla_\mu\Jhat^{(k)}_\red(\mu^{(k)}-sy)-\nabla_\mu\Jhat(\mu^{(k)})\right)\cdot y \,\mathrm{d} s.
		\end{aligned}
	\end{equation}
	Now, the integral term can be estimated, exploiting the Lipschitz continuity of $\nabla_\mu \Jhat_\red^{(k)}$ from \cref{asmpt:Lip_cont}, as
	\begin{equation}
		\label{int_estimate_Lip}
		\begin{split}
			\bigg| \int_0^1 \Big(\nabla_\mu\Jhat^{(k)}_\red(\mu^{(k)}-&sy)-\nabla_\mu\Jhat(\mu^{(k)})\Big)\cdot y \,\mathrm{d} s \bigg| \\&\leq C_L\int_0^1 \|sy\|_2\|y\|_2\,\mathrm{d}s = \frac{C_L}{2}\|\mu^{(k)}(j)-\mu^{(k)}\|^2_2.
		\end{split}
	\end{equation}
	Multiplying~\eqref{1step_AGC_thm} by $\kappa^j$ and using~\eqref{int_estimate_Lip} together with \cref{lemma:kelley}, we obtain
	\[
	\begin{aligned}
		\kappa^j\Big(\Jhat_\red^{(k)}(\mu^{(k)})&-\Jhat_\red^{(k)}(\mu^{(k)}(j))\Big) \\ & \geq \kappa^j\nabla_\mu\Jhat_\red^{(k)}(\mu^{(k)})\cdot(\mu^{(k)}-\mu^{(k)}(j))- \frac{C_L\kappa^j}{2}\|\mu^{(k)}(j)-\mu^{(k)}\|^2_2 \\
		& \geq \left(1-\frac{C_L\kappa^j}{2}\right)\|\mu^{(k)}(j)-\mu^{(k)}\|^2_2.
	\end{aligned}
	\]
	Thus, we have
	\[
	\Jhat_\red^{(k)}(\mu^{(k)}(j))-\Jhat_\red^{(k)}(\mu^{(k)})\leq -\frac{1}{\kappa^j}\left(1-\frac{C_L\kappa^j}{2}\right)\|\mu^{(k)}(j)-\mu^{(k)}\|^2_2.
	\]
	Since $j$ satisfies~\eqref{j_low_bound}, we have that $\kappa^j\leq \frac{2(1-\kappa_{\mathsf{arm}})}{C_L}$.
	Therefore, the Armijo-type condition~\eqref{Armijo} is met.
	It remains to show that~\eqref{TR_radius_condition} holds as well.
	Note that
	\[
	\|\mu^{(k)}(j)-\mu^{(k)}\|_2\leq \|\kappa^j\nabla_\mu\Jhat_\red^{(k)}(\mu^{(k)})\|_2\leq \eta((1-\beta_3)\tau_{{\textnormal{mac}}})
	\]
	by the choice of $j$ and \cref{cor:bound_grad}. Now, \cref{asmpt:rejection}, \cref{lemma:BothRBErrorsSmaller2}, and \cref{asmpt:unif_cont_q} imply that 
	\[
	\begin{aligned}
		q^{(k)}(\mu^{(k)}(j))&\leq|q^{(k)}(\mu^{(k)}(j))-q^{(k)}(\mu^{(k)})|+q^{(k)}(\mu^{(k)}) \\ & < (1-\beta_3)\tau_{{\textnormal{mac}}} +\beta_3\delta^{(k)}\leq \delta^{(k)},
	\end{aligned}\] 
	which completes the proof.
\end{proof}
In the next step of the convergence analysis, we show that \cref{alg:TR-RBmethod_paper2} ensures that the error-aware sufficient decrease condition~\eqref{eq:suff_decrease_condition} is satisfied for every accepted iteration; cf. \cref{sec:acceptance_of_TR}.
\begin{lemma}
	\label{lemma:EASDCFulfilled}
	Let the iterate $\mu^{(k+1)}$ be accepted by \cref{alg:TR-RBmethod_paper2}. Then, the error-aware sufficient decrease condition~\eqref{eq:suff_decrease_condition} is satisfied.
\end{lemma}
\begin{proof}
	When the RB model is updated, we can proceed as in~\cite[Section~4.1]{QGVW2017}. If the model is not enriched we have to distinguish two cases:
	\begin{enumerate}
		\item If $\mu^{(k+1)}$ is accepted in line~\ref{AcceptingIterateSufficientCondition}
		and \textsf{Skip\_enrichment\_flag} is true, we have
		\[
		\cJhatn^{(k+1)}(\mu^{(k+1)}) =\cJhatn^{(k)}(\mu^{(k+1)})\kern-0.1em \leq\kern-0.1em \cJhatn^{(k)}(\mu^{(k+1)}) + \Delta_{\cJhatn^{(k)}}(\mu^{(k+1)}) \kern-0.1em < \kern-0.1em \cJhatn^{(k)}(\mu_\textnormal{{AGC}}^{(k)}).
		\]
		\item $\mu^{(k+1)}$ is accepted in line~\ref{AcceptingIterateNoUpdateByExactComputation}: \\
		Note that it always holds $\cJhatn^{(k)}(\mu^{(k+1)}) \leq \cJhatn^{(k)}(\mu_\textnormal{{AGC}}^{(k)})$, since the used method for solving the TR sub-problem~\eqref{TRRBsubprob} is a descent method and the AGC point $\mu_\textnormal{{AGC}}^{(k)}$ is used as a warm start. Since the RB model is not updated, it holds
		\[
		\cJhatn^{(k+1)}(\mu^{(k+1)}) = \cJhatn^{(k)}(\mu^{(k+1)}) \leq \cJhatn^{(k)}(\mu_\textnormal{{AGC}}^{(k)}),
		\]
		so that~\eqref{eq:suff_decrease_condition} is satisfied.
	\end{enumerate}
	Thus, whenever $\mu^{(k+1)}$ is accepted, regardless updating the RB model or not, the error-aware sufficient decrease condition~\eqref{eq:suff_decrease_condition} is fulfilled.
\end{proof}
We are now able to prove the improved convergence results for superseding \cref{thm:convergence_of_TR}, also taking into consideration the possibility of skipping RB model updates.
\begin{theorem}
	\label{thm:convergence_of_TR_2}
	Every accumulation point $\bar\mu$ of the sequence $\{\mu^{(k)}\}_{k\in\mathbb{N}}\subset \Params$ generated by \cref{alg:Basic_TR-RBmethod} is an approximate first-order critical point for $\Jhat_h$, i.e., it holds
	\begin{equation}
		\label{First-order_critical_condition}
		\|\bar \mu-\Proj_{\Paramsad}(\bar \mu-\nabla_\mu \Jhat_h(\bar \mu))\|_2 = 0.
	\end{equation}
\end{theorem}
\begin{proof}
	Let $k\in\mathbb{N}$ be arbitrary. From \cref{def:AGC},~\eqref{Armijo}, and ~\eqref{eq:suff_decrease_condition} due to \cref{lemma:EASDCFulfilled}, we have
	\[
	\begin{aligned}
		\Jhat_\red^{(k)}(\mu^{(k)})-\Jhat_\red^{(k+1)}(\mu^{(k+1)}) & \geq \Jhat_\red^{(k)}(\mu^{(k)})-\Jhat_\red^{(k)}(\mu^{(k)}_\textnormal{AGC}) \\
		& \geq \frac{\kappa_{\mathsf{arm}}}{\kappa^{j^{(k)}_c}}\|\mu^{(k)}-\mu^{(k)}_\textnormal{AGC}\|^2_2 \\
		& = \frac{\kappa_{\mathsf{arm}}}{\kappa^{j^{(k)}_c}}\|\mu^{(k)}-\Proj_{\Paramsad}(\mu^{(k)}-\kappa^{j^{(k)}_c}\nabla_\mu\Jhat_\red^{(k)}(\mu^{(k)}))\|^2_2.
	\end{aligned}
	\]
	By summing both sides of the previous inequality from $k=0$ to $K$, we obtain
	\[
	\begin{aligned}
		\Jhat_\red^{(0)}(\mu^{(0)})-\Jhat_\red^{(K+1)}(\mu^{(K+1)}) \geq \\ \sum_{k=0}^K \frac{\kappa_{\mathsf{arm}}}{\kappa^{j^{(k)}_c}}\|\mu^{(k)}-&\Proj_{\Paramsad}(\mu^{(k)}-\kappa^{j^{(k)}_c}\nabla_\mu\Jhat_\red^{(k)}(\mu^{(k)}))\|^2_2 \geq 0.
	\end{aligned}
	\]
	For $K\to +\infty$, the term on the left-hand side is bounded from above, due to Assumption~\ref{asmpt:bound_J}. Thus
	\[
	\lim_{k\to\infty} \frac{\kappa_{\mathsf{arm}}}{\kappa^{j^{(k)}_c}}\|\mu^{(k)}-\Proj_{\Paramsad}(\mu^{(k)}-\kappa^{j^{(k)}_c}\nabla_\mu\Jhat_\red^{(k)}(\mu^{(k)}))\|^2_2 = 0.
	\]
	From Theorem~\ref{thm:linesearch-agc}, we have that $\kappa^{j^{(k)}_c}\geq \min\left\{\frac{2(1-\kappa_{\mathsf{arm}})}{C_L},\frac{\eta((1-\beta_3)\tau_{{\textnormal{mac}}})}{C}\right\}:= \widetilde{\kappa}$ for all $k\in\mathbb{N}$. Furthermore, we also have that $\kappa^{j^{(k)}_c}\leq 1$ for all $k\in\mathbb{N}$, because $\kappa\in(0,1)$. Hence,
	\begin{align*}
		\lim_{k\to\infty} \kappa_{\mathsf{arm}}\|\mu^{(k)}-\Proj_{\Paramsad}(\mu^{(k)}-\widetilde{\kappa}\nabla_\mu\Jhat_\red^{(k)}(\mu^{(k)}))\|^2_2 & \leq \\ \lim_{k\to\infty} \frac{\kappa_{\mathsf{arm}}}{\kappa^{j^{(k)}_c}}\|\mu^{(k)}-\Proj_{\Paramsad}(\mu^{(k)}-\kappa^{j^{(k)}_c}&\nabla_\mu\Jhat_\red^{(k)}(\mu^{(k)}))\|^2_2 = 0
	\end{align*}
	which implies
	\[
	\lim_{k\to\infty} \|\mu^{(k)}-\Proj_{\Paramsad}(\mu^{(k)}-\widetilde{\kappa}\nabla_\mu\Jhat_\red^{(k)}(\mu^{(k)}))\|_2 = 0.
	\]
	Lemma~\ref{lemma:PropertyOfPParams} shows that
	\begin{align*}
		\left\| \mu^{(k)} - \Proj_{\Paramsad} ( \mu^{(k)} - \widetilde{\kappa} \nabla_\mu \cJhatn^{(k)}(\mu^{(k)}) ) \right\|_2 \geq \hspace{3em}\\ \widetilde{\kappa} \left\| \mu^{(k)} - \Proj_{\Paramsad} ( \mu^{(k)} - \nabla_\mu \cJhatn^{(k)}(\mu^{(k)}) ) \right\|_2
	\end{align*}
	holds for all $k \in \N$. Thus, we conclude
	\begin{align}
		\label{lim_k_to_zero}
		& \lim_{k \to \infty} \left\| \mu^{(k)} - \Proj_{\Paramsad} ( \mu^{(k)} - \nabla_\mu \cJhatn^{(k)}(\mu^{(k)}) ) \right\|_2 = 0.
	\end{align} 
	By Lemma~\ref{lemma:BothRBErrorsSmaller2}, we have
	\begin{align*}
		& \left| \left\| \mu^{(k)} - \Proj_{\Paramsad} ( \mu^{(k)} - \nabla_\mu \cJhatn^{(k)}(\mu^{(k)}) ) \right\|_2 - \left\| \mu^{(k)} - \Proj_{\Paramsad} ( \mu^{(k)} - \nabla_\mu \Jhat_h(\mu^{(k)}) ) \right\|_2 \right| \\
		& \qquad \leq \left\| \mu^{(k)} - \Proj_{\Paramsad} ( \mu^{(k)} - \nabla_\mu \cJhatn^{(k)}(\mu^{(k)}) ) \right\|_2 \tau_g \to 0 \quad \textnormal{as } k \to \infty,
	\end{align*}
	which, together with~\eqref{lim_k_to_zero}, implies 
	\begin{align}
		\label{lim_k_to_zero_part2}
		\lim_{k \to \infty}  \left\| \mu^{(k)} - \Proj_{\Paramsad} ( \mu^{(k)} - \nabla_\mu \Jhat_h(\mu^{(k)}) ) \right\|_2 = 0.
	\end{align}
	Now let $\bar{\mu}$ be an accumulation point of the sequence $\{\mu^{(k)}\}_{k\in\mathbb{N}}\subset \Params$, i.e., it holds 
	\[
	\mu^{(k_i)} \to \bar{\mu} \quad \textnormal{ as } i \to \infty
	\]
	for some subsequence $\left\{\mu^{(k_i)}\right\}_{i \in \mathbb{N}}$. Using~\eqref{lim_k_to_zero_part2}, we have by the continuity of the gradient $\nabla_\mu \Jhat_h$ and of the projection operator $\Proj_{\Paramsad}$
	\begin{align*}
		\left\| \bar\mu - \Proj_{\Paramsad} ( \bar\mu - \nabla_\mu \Jhat_h(\bar\mu) ) \right\|_2 \kern-0.1em=\kern-0.1em 
		\lim_{i \to \infty} \left\| \mu^{(k_i)} - \Proj_{\Paramsad} ( \mu^{(k_i)} - \nabla_\mu \Jhat_h(\mu^{(k_i)}) ) \right\|_2 = 0,
	\end{align*}
	which concludes the proof.
\end{proof}
The biggest improvement in comparison to the convergence proof from \cref{convergence_study_1} constitutes of the fact that any accumulation point of the sequence $\{\mu^{(k)}\}$ is an actual critical point for $\Jhat_h$, and not an approximated one up to the tolerance $\tau_\textnormal{{sub}}$; cf. \cref{thm:convergence_of_TR}.
Furthermore, there is no direct use of the RB model properties in \cref{thm:convergence_of_TR_2} as it is in \cref{thm:convergence_of_TR}.
This opens the possibility of considering different model functions, similarly to~\cite{YM2013} for unconstrained parameter sets, provided that the requested (and shown) properties hold. 
	%
In addition, in contrast to \cref{convergence_study_1}, we make use of the error-aware sufficient decrease condition~\eqref{eq:suff_decrease_condition} in the proof of convergence and not only to guarantee that the accumulation point is not a local maximum of $\Jhat_h$.

\begin{remark}
	As also highlighted in \cref{rmk:local_minimum_TR}, $\bar\mu$ can still be a saddle point or a local maximum. In the numerical experiments, to verify that the computed point $\bar \mu$ is a local minimum, we check the second-order sufficient optimality conditions (cf. \cref{prop:second_order}) as soon as the algorithm terminates.	
\end{remark}

\subsection{Basis enrichment}
\label{sec:TRRB_construct_RB} 

At each point of the outer iteration of the TR-RB algorithm, presented in \cref{alg:Basic_TR-RBmethod}, we may enrich the model, and multiple enrichment strategies are possible for this.
Importantly, not all of them led to an efficient algorithm in our numerical experiments.
To recall, we are interested in enriching $V^{\pr,k}_{\red}$ and $V^{\du,k}_{\red}$ at a new iterate $\mu^{(k)}$.
At such an enrichment step for $\mu \in \Params$, we assume to have access to the primal and dual solutions $u_{h, \mu}, p_{h, \mu} \in V_h$.
Note that solving for these FOM quantities has high-dimensional complexity.
A straightforward way would be to distribute the primal and dual solutions to their respective separated reduced spaces (which are commonly referred to as Lagrange spaces).
Another choice is to construct only a single (aggregated) RB space for both the primal and dual system.

We introduce a third strategy tailored to the projected Newton method as the sub-problem solver.
Besides the primal and dual solution, we can invest more computational effort to additionally compute their respective sensitivities $d_{\nu} u_{h, \mu}, d_{\nu} p_{h, \mu} \in V_h$ w.r.t. a certain direction $\eta$.
If we use the projected Newton method as the sub-problem solver, we require directional derivatives of the reduced primal and dual solutions for computing the Hessian following \cref{prop:true_corrected_reduced_Hessian}.
However, it can not be guaranteed that a reduced solution $d_{\eta} u_{\red, \mu}$ is a good approximation of $d_{\eta} u_{h, \mu}$; cf. \cref{sec:TRRB_sensitivity_approach}.
As discussed earlier, it is cheap to compute the gradient $\nabla_\mu \Jhat_h(\mu)$ if $u_{h, \mu}, p_{h, \mu} \in V_h$ are already available.
As mentioned above, to proceed with the next sub-problem after enrichment, we need the directional sensitivities of $u_{\red, \mu}$ and $p_{\red, \mu}$ in the direction $\eta=\nabla_\mu \cJhatn(\mu)$.
We can thus add the corresponding FOM directional sensitivities as the third enrichment strategy.

In conclusion, we consider the following three enrichment strategies.
\begin{description}
	\item[(a) Lagrangian RB spaces:] \label{enrich:lag} We add each FOM solution to the RB space that is directly related to its respective reduced formulation, i.e.~for a
	given $\mu \in \Params$, we enrich by
	$
	V^{\pr,k}_{\red} = V^{\pr,k-1}_{\red} \cup \{u_{h,\mu}\},
	V^{\du,k}_{\red} = V^{\du,k-1}_{\red} \cup \{p_{h,\mu}\}.
	$
	\item[(b) Aggregated RB space:] We add all available information into a single RB space,
	i.e. $V^{\pr, k}_{\red} = V^{\du, k}_{\red} = V^{\pr,k-1}_{\red} \cup \{u_{h,\mu}, p_{h,\mu}\}$.
	According to \cref{sec:TRRB_MOR_for_PDEconstr} this results in $\cJhatn(\mu)$ being equal to the standard RB reduced functional from~\eqref{eq:Jhat_red}.
	\label{enrich:single}
	\item[(c) Directional Taylor RB space:] \label{enrich:dir_tay}
	We compute a direction $\eta:=\nabla_\mu \cJhatn(\mu)$ from $u_{h, \mu}$ and $p_{h, \mu}$ and include the directional derivatives to the respective RB space, i.e.
	$
	V^{\pr,k}_{\red} = V^{\pr,k-1}_{\red} \cup \{u_{h,\mu}\} \cup \{d_{\eta} u_{h,\mu}\}, 
	V^{\du,k}_{\red} = V^{\du,k-1}_{\red} \cup \{p_{h,\mu}\} \cup \{d_{\eta} p_{h,\mu}\}.
	$	   
\end{description} 

The above-presented strategies for constructing RB spaces have a significant impact on the performance and accuracy of the TR-RB method.
Note that "$\cup$" also incorporates orthonormalization of the respective RB space w.r.t. the inner product.
The offline-online decomposition for the surrogate model is subsequently carried out for a new basis. 
Importantly, offline computations for the construction of RB models scale quadratically with the number of basis functions in the RB space, cf. \cref{sec:background_off_on_RB}.
Thus, Lagrange RB spaces in (a) are computationally beneficial compared to (b) at a potential loss of accuracy of the corresponding RB models (since less information is added to both RB spaces).
Moreover, Lagrange RB spaces as in (a) destroy the duality of state and adjoint equations, cf. \cref{sec:TRRB_standard_approach}, and hence, the NCD-correction term is required.
We emphasize that the Lagrangian enrichment (a) can be considered the standard procedure for the discussed TR-RB methods, also used in~\cite{QGVW2017}.
As aforementioned, the directional Taylor enrichment (c) is mainly helpful for a Newton-based sub-problem solver.

We also comment on other enrichment strategies and why we do not consider them as helpful for the TR-RB algorithm.
In \cref{sec:TRRB_sensitivity_approach}, we also discussed the possibility to construct separated sensitivity spaces for each canonical direction, i.e. $V_\red^{\pr,d_{\mu_i}}$ and $V_\red^{\du,d_{\mu_i}}$ for all $i = 1, \dots, P$.
In the sense of Lagrangian RB spaces we could thus also compute $d_{\mu_i} u_{h, \mu}$ and $d_{\mu_i} p_{h, \mu}$ component-wise and include it into these spaces.
Alternatively, we could also add them to $V^{\pr}_{\red}$ and $V^{\pr,k}_{\red}$, which constitutes full Taylor RB spaces, or add all bases into a large aggregated space with sensitivities.
From the theory (and also shown by the numerical experiments in \cref{sec:TRRB_estimator_study}), we can expect the highest ROM accuracy and fewer outer iterations in the TR-RB if FOM sensitivities are included in the respective spaces.
However, computing all canonical sensitivities results in a prohibitively sizeable computational effort for high-dimensional parameter spaces since $2 \cdot P$ FOM solves are required for each enrichment.
On top of that, the dimension of the Taylor and aggregated space with sensitivities would grow with a factor of $P$ and $2 \cdot P$, respectively, which is also a computational burden for the offline-online efficiency of the RB scheme.
In contrast to this, strategies (a), (b), and (c) only require $2$, $2$, and $4$ FOM evaluations, respectively, and thus the required FOM evaluations are independent of $P$.
Since our TR-RB is mainly designed for large parameter space cases, we do not follow other strategies than (a)--(c).

Lastly, we mention that enrichment strategies at iteration $k$ are not necessarily restricted to $\mu^{(k)}$.
Instead, it is also feasible to construct a metric box (not dependent on the estimator) around the current iterate with the aim to construct an RB space that is sufficiently accurate in a neighborhood of $\mu^{(k)}$.
This can be referred to as a local greedy-search algorithm.
However, especially for large parameter spaces, we saw in our numerical experiments that the required computational effort for such an intermediate offline-online approach does not pay off for fast convergence of the optimization method.
This can particularly be explained with the same reasoning as we have used for the motivation of adaptively constructed surrogates in \cref{sec:background_RB_methods_for_PDEopt}, namely that it is not feasible to train for a region in the parameter space that is far away from the local optimum.
On the other hand, such an approach may be helpful close enough to convergence, where it can be expected that the optimum lays inside the trained region.
In our experiments, we did not see a major benefit from local greedy strategies which is why we do not follow them in the sequel.

\subsection{TR-RB variants based on adaptive enrichment strategies}
\label{sec:TRRB_variants}
A significant contribution of this thesis is to introduce and analyze variants of adaptive TR-RB methods for efficiently computing a solution of the optimization problem~\eqref{P}.
In terms of performance, we need to account for all computational costs, including the algorithms' offline and online expenses.
Following \cref{sec:TR_algorithm_details}, we propose a TR method that adaptively builds an RB space along the optimization path (see \cref{alg:Basic_TR-RBmethod}).
From a MOR perspective, this diminishes the offline time of the ROM significantly since no global RB space (concerning the parameter domain) has to be built in advance.
We may enrich the model after the sub-problem~\eqref{TRRBsubprob} of the TR method has been solved and if the corresponding iterate is accepted.

The proposed methods differ mainly in the model function, derivative information, sub-problem solvers, and enrichment strategies.

\subsubsection{Choosing the reduced model}
For choosing the ROM in this chapter, we distinguish four different approaches:
\begin{description}
	\item[1. Standard approach:] Following \cref{sec:TRRB_standard_approach} and proposed in \cite{QGVW2017}, the standard approach for the functional is to replace the FOM quantities by their respective ROM counterpart, i.e.~we consider the map $\mu \mapsto \Jnoncor_\red(\mu)$ from~\eqref{eq:Jhat_red}. 
	Gradient information can be computed by $\noncorgrad_\mu \Jnoncor_\red(u_{\red, \mu}, \mu)$ from~\eqref{naive:red_grad}.
	For the corresponding error estimation, $\Delta_{\Jnoncor_\red}(\mu)$ from \cref{prop:Jhat_error}(i) can be used.
	\item[2. Semi NCD-corrected approach:] A first correction strategy is to replace the functional by the NCD-corrected RB reduced functional $\cJhatn$ from~\eqref{eq:Jhat_red_corected} but stick with the inexact gradient of the standard approach.
	This allows using the higher-order estimator for the functional $\Delta_{\cJhatn}(\mu)$ from \cref{prop:Jhat_error}(ii).
	\item[3. NCD-corrected approach:] 
	As explained in \cref{sec:TRRB_ncd_approach}, we propose to consider the NCD-corrected RB reduced functional $\cJhatn$ from~\eqref{eq:Jhat_red_corected} and its actual gradient according to \cref{prop:true_corrected_reduced_gradient_adj}. Again, the estimator $\Delta_{\cJhatn}(\mu)$  \cref{prop:Jhat_error}(ii) is used.
	\item[4. PG approach:] An alternative ROM strategy is explained in \cref{sec:TRRB_pg_approach},
	where we use the Petrov--Galerkin reduced model and use $\cJhatn^{\textnormal{pg}}$ from \cref{eq:Jhat_red_pg} for which the gradient information can be computed following the FOM-type formula in~\eqref{eq:red_gradient}.
	The choice of the functional requires to use $\Delta^{\textnormal{pg}}_{\cJhatn}(\mu)$ from \cref{prop:error_reduced_quantities}(iii) for estimation.
\end{description}

\subsubsection{Choosing the sub-problem solver}

As presented in \cref{sec:sub-problem_solvers}, we choose between the projected BFGS algorithm or the projected Newton algorithm, tailored towards our optimization problem and reduced model.
We refer to these methods by the abbreviations BFGS and Newton.

\subsubsection{Choosing the basis enrichment strategy} \label{sec:choose_enrichment}
For the basis construction, we may use variants (a), (b) or (c) from \cref{sec:TRRB_construct_RB}.
Note, however, that, by using the basis enrichment (b), all ROM approaches $1.$ - $4.$ are equivalent, and (c) is mainly useful for Newton.
On top of that, in \cref{sec:TRRB_optional_enrichment}, we have proposed a notion of skipping basis enrichments and proved that this strategy also leads to convergence, given the circumstances that are devised in \cref{convergence_study_2}.
Thus, we distinguish between the following two approaches:

\begin{description}
	\item[Unconditional Enrichment (UE):] For every new outer iterate $\mu^{(k)}$, we enrich the RB space unconditionally.
	\item[Optional Enrichment (OE):] Before enriching the RB space, we evaluate the criterion in~\eqref{eq:skip_enrichment_condition} to decide for enrichment.
\end{description}

The discussed variants in this section significantly deviate from the original work in~\cite{QGVW2017}.
The presence of inequality constraints, which are missing in~\cite{QGVW2017}, implies a projection-based optimization algorithm.
In addition, differently from~\cite{QGVW2017}, we take advantage of more suitable ROM strategies, enlarging the TR radius, optional RB enrichment, an outer stopping criterion that is independent of the RB a posteriori estimates, and a post-processing strategy for the optimal parameter.

In the numerical experiments in \cref{sec:TRRB_num_experiments}, we extensively study and elaborate our algorithm's different features.
In every numerical experiment, we clarify which of the above-discussed variants of the TR-RB algorithm we use.

\section{Software and implementational design}
\label{sec:software_TRRB}

In this section, we briefly present the numerical software and implementational details that have been used to conduct the numerical experiments for this chapter.
All simulations have been performed in a pure \texttt{Python} implementation and are based on the open-source MOR library \texttt{pyMOR}~\cite{milk2016pymor}, making use of \texttt{pyMOR}'s builtin vectorized \texttt{numpy/scipy}-based discretizer for the FOM and generic MOR algorithms for projection and orthonormalization (such as a stabilized Gram-Schmidt algorithm) to obtain efficient ROMs effortlessly.
In what follows, we review basic concepts of \texttt{pyMOR}'s interfaces and algorithms.
Moreover, we mention the extensions required to implement the TR-RB method.
The source code to reproduce all results (including detailed interactive \texttt{jupyter}-notebooks) is referenced in \cref{apx:code_availability}.

\subsection{Model order reduction with \texttt{pyMOR}}

This section aims to give an (incomplete) overview of \texttt{pyMOR}~\cite{milk2016pymor}.
In particular, we refer to the interfaces and algorithms relevant to our implementation and do not discuss the technical details or other application fields of \texttt{pyMOR}, such as MOR methods in system theory. Elaborated tutorials for getting started with \texttt{pyMOR} are available on the respective web-page and demos can be found in the respective source code.
We emphasize that this section points to the 2020.2 release of \texttt{pyMOR}.

\begin{description}
	\item[Parameters:] The notion of parameters (or inputs) in a parameter space like $\Params \subset \R^P$ and an understanding of how to treat those are essential for MOR methods and is well-implemented in \texttt{pyMOR}.
	Problems and their respective objects always expect to be called for different parameters (or inputs).
	\item[Analytical problems:] There are different ways how to get started with \texttt{pyMOR}.
	Suppose the user only has an analytical problem and does not have any discretization code yet. In that case, it is convenient to use \texttt{pyMOR}'s built-in discretizer that discretizes a predefined analytical problem.
	If, however, the user already has a discretization, the user may start with defining the corresponding \texttt{Model}; see below.
	
	\texttt{pyMOR} has several predefined analytical problems that can be initialized with many different problem parameters.
	In our case, we are interested in stationary elliptic problems, which is why we use the \texttt{StationaryProblem}.
	Respective expression functions (on the domain) and parameter functionals (on $\Params$) for defining the respective data functions can be defined, for instance, by \texttt{ExpressionFunction}
	or \texttt{ExpressionParameterFunctional}.
	To attain a parameter separable function as in \cref{def:parameter_separable}, the \texttt{LincombFunction} can be used.
	Notably, \texttt{pyMOR} is also able to load pictures and converts them to \texttt{BitmapFunctions}, as, for example, used in \cref{sec:TRRB_benchmark_problem}.
	Moreover, \texttt{pyMOR} also provides \texttt{gmsh}-bindings (see~\cite{gmsh}) for a user-defined mesh and boundary values, used in \cref{sec:thermalfin}.
	\item[Discretizers:] Given such an analytical problem, dependent on the problem class, \texttt{pyMOR}'s discretizers are used to discretize the given domain and the problem to obtain the respective \texttt{numpy}-based matrices that are required, for instance, for FEM approximations; cf. \cref{sec:algebraic_fem}.
	For the \texttt{StationaryProblem}, the \texttt{discretize{\textunderscore}stationary{\textunderscore}cg} function constructs such a continuous Galerkin (cg) discretization.
	\item[Models:] The discretizer eventually returns a \texttt{Model}.
	Apart from some technical details, a \texttt{Model} mainly consists of a \texttt{solve(mu)}-method that computes respective solutions, given a parameter $\mu \in \Params$.
	Output functionals, its parameter gradient, and sensitivities of the solution can be obtained by the \texttt{output}-, \texttt{solve{\textunderscore}d{\textunderscore}mu}- and \texttt{output{\textunderscore}d{\textunderscore}mu}-method, respectively.
	Let us also mention that such a \texttt{Model} internally holds \texttt{Operators} for the left-hand side operator, right-hand side, and output functional of the system.
	In our case of a \texttt{StationaryProblem}, we use the corresponding \texttt{StationaryModel} for solving \cref{def:param_elliptic_problem}.
	
	\item[Reductors:] The reduction process in \texttt{pyMOR} is implemented in the \texttt{Reductors}.
	These objects hold information on the FOM, the reduced basis, and other reduction-related objects.
	For the case of RB methods, for the \texttt{StationaryModel}, we may use the \texttt{CoerciveRBReductor}.
	Given a reduced basis, the  \texttt{Reductor} has a \texttt{reduce}-method, which reduces the FOM and,
	given \cref{asmpt:parameter_separable}, prepares the offline-online decomposition as explained in \cref{sec:background_off_on_RB}, including the assembly of the error estimator and the Galerkin projection of the involved matrices and vectors.
	The resulting reduced model is, just as the FOM, a \texttt{StationaryModel} with replaced \texttt{Operator}s and the respective \texttt{error{\textunderscore}estimator}.
	We also mention that \texttt{CoerciveRBReductor} implements the stable and more accurate estimator assembly as is proposed in \cite{buhr14}, whereas \texttt{SimpleCoerciveRBReductor} uses the standard assembly of the error estimator; cf. \cref{sec:background_off_on_error}.
	
	Since the \texttt{Reductor} also holds the RB basis, the corresponding \texttt{reconstruct}-method that can be used to reconstruct the high-fidelity solution from the reduced coefficients.
	
	For an orthonormal reduced basis for the initialization of the respective \texttt{Reductor}, \texttt{pyMOR} also contains a stabilized Gram-Schmidt algorithm.
	If, instead, the \texttt{Reductor} is enriched on the fly, for instance, employed by a weak greedy-search algorithm (which is also available in \texttt{pyMOR} as \texttt{weak{\textunderscore}greedy}), the \texttt{Reductor} contains an \texttt{extend{\textunderscore}basis}-method.
	This method internally and automatically re-orthonormalizes the basis and updates the projected operators and the offline-online decomposition with minimal additional effort.
\end{description}

While \texttt{pyMOR} is already a remarkably flexible MOR library, it has (at least to our knowledge) never before been used for PDE-constrained parameter optimization methods that include explicit derivative information.
Thus, lots of functionalities still had to be added.
Apart from the before-mentioned parameter derivatives, by now, lots of related features have by now entered the releases of \texttt{pyMOR}.
Moreover, a \texttt{pyMOR} tutorial for PDE-constrained parameter optimization problems has been added, which aims to introduce interested readers to the discussed optimization problems.

Nevertheless, for the experiments in this chapter, many additional features that did not yet enter the releases remain, which we elaborate on in the subsequent section.

\subsection{\texttt{pdeopt} extensions to \texttt{pyMOR}}
\label{sec:pdeopt_pymor}
In the following, we explain the \texttt{pdeopt} module that can be used as an extension to \texttt{pyMOR} for PDE-constrained parameter optimization problems.
So far, the output functional in the \texttt{StationaryModel} is expected to be linear, which does not fit the assumption in \eqref{eq:quadratic_J}.
Furthermore, from the explanations in \cref{sec:TRRB_MOR_for_PDEconstr}, we require specialized methods that not only make it possible to compute all quantities but also account for the different variants that have been presented in the respective section.

In order to follow the implementational design of \texttt{pyMOR} as explained in the previous section, we assume to be given a \texttt{StationaryProblem}. Subsequently, we need a specialized \texttt{discretizer}, \texttt{Model}, and \texttt{Reductor}.

\begin{description}
	\item[\texttt{discretizer}:] We implemented a \texttt{discretize{\textunderscore}quadratic{\textunderscore}pdeopt{\textunderscore}stationary{\textunderscore}cg}-function that internally uses the standard \texttt{discretize{\textunderscore}stationary{\textunderscore}cg} and additionally accounts for the quadratic objective functional as it is defined in \eqref{eq:quadratic_J}.
	Respective \texttt{Operators} are built to feed a corresponding model.
	\item[\texttt{Model}:] The main challenge for using \texttt{pyMOR}'s reduction methods is the definition of a respective \texttt{QuadraticPdeoptStationaryModel} that inherits from \texttt{StationaryModel}.
	This model is capable of computing all the mentioned quantities of \cref{sec:TRRB_fom}.
	For instance, to compute the solution of the dual equation, the respective function \texttt{solve{\textunderscore}dual} is available.
	Moreover, sensitivities, gradient, and Hessian information of the quadratic objective functional can be computed.
	We also note that the \texttt{Model} internally uses the original functionality of \texttt{StationaryModel} as often as possible.
	\item[\texttt{Reductor}:] The purpose of the \texttt{QuadraticPdeoptStationaryCoerciveReductor} is to construct the surrogate models that are explained in \cref{sec:TRRB_MOR_for_PDEconstr}, including all \texttt{Operators}, \texttt{estimators} and their offline-online decomposition.
	The main challenge in the implementational design was to make the \texttt{QuadraticPdeoptStationaryCoerciveReductor} flexible w.r.t. their respective reduced scheme.
	This flexibility also has to be prolonged by the \texttt{QuadraticPdeoptStationaryModel} to compute, for instance, the NCD-corrected functional, the adjoint-based gradient, the auxiliary functions, and also the respective Hessian.
	Moreover, all estimators from \cref{sec:TRRB_a_post_error_estimates} have to be computable.
	
	We emphasize that, depending on the ROM, we implemented the \texttt{Reductor} in such a way that it does not invest unnecessary computational time in preparing reduced quantities that are not used by the respective algorithm.
	For instance, if it is priorly known that the optimizer will use no Hessian, the offline-online decomposition of the required parts is not performed at all.
	Furthermore, as they will not be used in the TR-RB algorithm, the assembly of the a posteriori error estimates for the gradient and Hessian are turned off by default.
	
	We would also like to point out that the \texttt{QuadraticPdeoptStationaryCoerciveReductor} is flexible concerning its \texttt{reductor{\textunderscore}type}.
	While the more accurate and stable preassembly of the estimates from~\cite{buhr14} (cf. \cref{sec:background_off_on_RB}) is readily available in the \texttt{CoerciveRBReductor}, the slightly cheaper preassembly of the estimates with the \texttt{SimpleCoerciveRBReductor} was sufficient for our experiments.
	
	For the efficient implementation of the coercivity constant, the $\min$-theta approach as explained in \cref{sec:min_theta} can be used, employing the mesh-independent energy-product $(u, v):= a_{\check{\mu}}(u, v)$ for a fixed parameter $\check{\mu} \in \Params$, which also helps in terms of the effectivity of the estimates.
	For the continuity constants that are needed for the estimate, the $\max$-theta approach from~\cite[Ex.~5.12]{HAA2017} is used to obtain upper bounds on the respective constants.
\end{description}

\subsection{Implementation of the optimization algorithms}

It remains to implement the optimization methods for the TR-RB algorithm. 
We did not use any optimization library for our outer and inner optimization algorithms, mainly because of the direct communication with \texttt{pyMOR}'s objects and the required flexibility in terms of the construction of the surrogate model. 
The \texttt{TR{\textunderscore}algorithm} for the (UE) variant and the \texttt{TR{\textunderscore}algorithm{\textunderscore}with{\textunderscore}optional{\textunderscore}enrichment} for the (OE) variant, cf. \cref{sec:choose_enrichment}, implement the respective TR-RB algorithm.
The other choices of variants, w.r.t. the sub-problem and the reduction scheme, are either hidden in the (predefined) \texttt{Reductor} or given as optional arguments.
We also note that further variants that we tested during the project, some of which are presented in \cref{sec:TRRB_related_approaches}, can be used by the same objects and functions.

We also provide self-written optimization code for Newton, BFGS, and other variants concerning the sub-problem solvers and the FOM methods.

\section{Numerical experiments}

\label{sec:TRRB_num_experiments}
We present numerical experiments to demonstrate the adaptive TR-RB variants from  \cref{sec:TRRB_variants} for quadratic objective functionals~\eqref{eq:quadratic_J} with elliptic PDE constraints as in~\eqref{Phat_second} and compare them to state-of-the-art algorithms from the literature.
Most of the presented experiments are already published in~\cite{paper2,paper1,paper3} which individually are devoted to concentrating on selected variants of the TR-RB method.
All experiments are based on the same implementation, explained in \cref{sec:software_TRRB} (including a re-implementation of~\cite{QGVW2017}) and were performed on the same machine multiple times to avoid caching or multi-query effects.
Timings may thus be used to compare and judge the computational efficiency of the different algorithms.

We start this section by defining a suitable benchmark problem capable of producing differently parameterized model problems.
In \cref{sec:proof_of_concept}, we give an exemplary proof of concept of a basic variant of the presented TR-RB algorithms for a 2-dimensional problem.
The aim is to demonstrate the TR-RB optimization procedure, including extensive visualizations of the performed steps.
The experiments in \cref{sec:exp_paper1} are devoted to analyzing the standard, semi-corrected, and corrected variants and to demonstrating the strength compared to the approach in~\cite{QGVW2017}.
A detailed analysis of the a posteriori error estimates from \cref{sec:TRRB_a_post_error_estimates} is conducted in \cref{sec:TRRB_estimator_study}.
For \cref{sec:exp_paper1}, only the BFGS based sub-problem solver is used, and no optional enrichment is considered.
In contrast, in \cref{sec:exp_paper2}, the main focus lies on the comparison of BFGS and Newton sub-problem solvers as well as the presentation of the optional enrichment and an application of the parameter control post-processing, based on \cref{sec:TRRB_RB_estimates_param}.
Lastly, the Petrov--Galerkin variant is analyzed in \cref{sec:exp_paper3}. 

\subsection{Problem definition and tolerances}

We consider stationary heat transfer in a bounded connected spatial domain $\Omega \subset \R^2$ with polygonal boundary $\partial\Omega$ partitioned into a non-empty Robin boundary $\Gamma_\textnormal{{R}} \subset \partial\Omega$ and possibly empty distinct Neumann boundary $\Gamma_\textnormal{{N}} = \partial\Omega\backslash \Gamma_\textnormal{{R}}$, and unit outer normal $n: \partial\Omega \to \R^2$.
We consider the Hilbert space $V = H^1(\Omega) := \{v \in L^2(\Omega) \,|\, \nabla v \in L^2(\Omega) \}$ of weakly differentiable functions and, for an admissible parameter $\mu \in \Params$, we seek the temperature $u_\mu \in V$ as the solution of
 \begin{equation} \label{eq:heat_equation}
	\begin{split}
		-  \nabla \cdot \left( A_\mu  \nabla u_{\mu} \right) &= f_{\mu} \hS{58} \textnormal{in } \Omega, \\
		c_\mu ( A_\mu  \nabla u_{\mu} \cdot n) &= (u_{\textnormal{out}} - u_{\mu}) \hS{16} \textnormal{on } \Gamma_\textnormal{{R}}, \\
		A_\mu\nabla u_\mu \cdot n &= g_\textnormal{{N}} \hS{58} \textnormal{on }\Gamma_\textnormal{{N}}
	\end{split}
\end{equation}
in the weak sense, with the admissible parameter set, the spatial domain and its boundaries and the data functions $A_\mu \in L^\infty(\Omega)$, $f_\mu \in L^2(\Omega)$, $c_\mu \in L^\infty(\Gamma_\textnormal{{R}})$, and $u_\textnormal{{out}} \in L^2(\Gamma_\textnormal{{R}})$ defined in the respective experiment.
The corresponding bilinear form $a_\mu$ and linear functional $l_\mu$ for \cref{def:param_elliptic_problem}, replicating the primal equation, are thus given for all $\mu \in \Params$ and $v, w \in V$ by
\begin{align}
	a_\mu(v, w) &:= \int_\Omega A_\mu \nabla v \cdot \nabla w \dx + \int_{\Gamma_\textnormal{{R}}}\hspace{-5pt}c_\mu\, vw\ds \qquad\text{and}\\
	l_\mu(v) &:= \int_\Omega f_\mu\,v\dx + \int_{\Gamma_\textnormal{{R}}}\hspace{-5pt}c_\mu\, u_\textnormal{{out}}v\ds + \int_{\Gamma_\textnormal{{N}}}\hspace{-5pt}g_\textnormal{{N}}v\ds.
\end{align}

As detailed in \cref{sec:software_TRRB}, for using the min-theta approach, we define the corresponding fixed energy product
$$
(u, v)_V := a_{\check{\mu}}(u, v),
$$
for a fixed parameter $\check{\mu} \in \Params$.

For the FOM, we fix a fine enough reference simplicial or cubic mesh and define $V_h \subset V$ as the respective space of continuous piecewise linear finite elements; cf. \cref{sec:background_npdgl}.
Note that \cref{asmpt:truth} is fulfilled by this choice of $V_h$, meaning that the FEM approximation can be considered the truth.

For all experiments, we use an initial TR radius of $\delta^{(0)} = 0.1$,
a TR shrinking factor $\beta_1=0.5$, an Armijo step-length $\kappa=0.5$,
a truncation of the TR boundary of $\beta_2 = 0.95$,
a tolerance for enlarging the TR radius of $\eta_\varrho = 0.75$,
a stopping tolerance for the TR sub-problems of $\tau_\textnormal{{sub}} = 10^{-8}$,
a maximum number of TR iteration $K = 40$,
a maximum number of sub-problem iterations $K_{\textnormal{{sub}}}= 400$, and
a maximum number of Armijo iterations of $50$. We also point out that the stopping tolerance for the FOC condition $\tau_\textnormal{{FOC}}$ is specified in each experiment individually.  

\subsection{A flexible benchmark problem for stationary heat}
\label{sec:TRRB_benchmark_problem}

\begin{SCfigure}[][h]
	\centering\footnotesize
	\begin{subfigure}[c]{0.6\textwidth}
		\includegraphics[width=\textwidth]{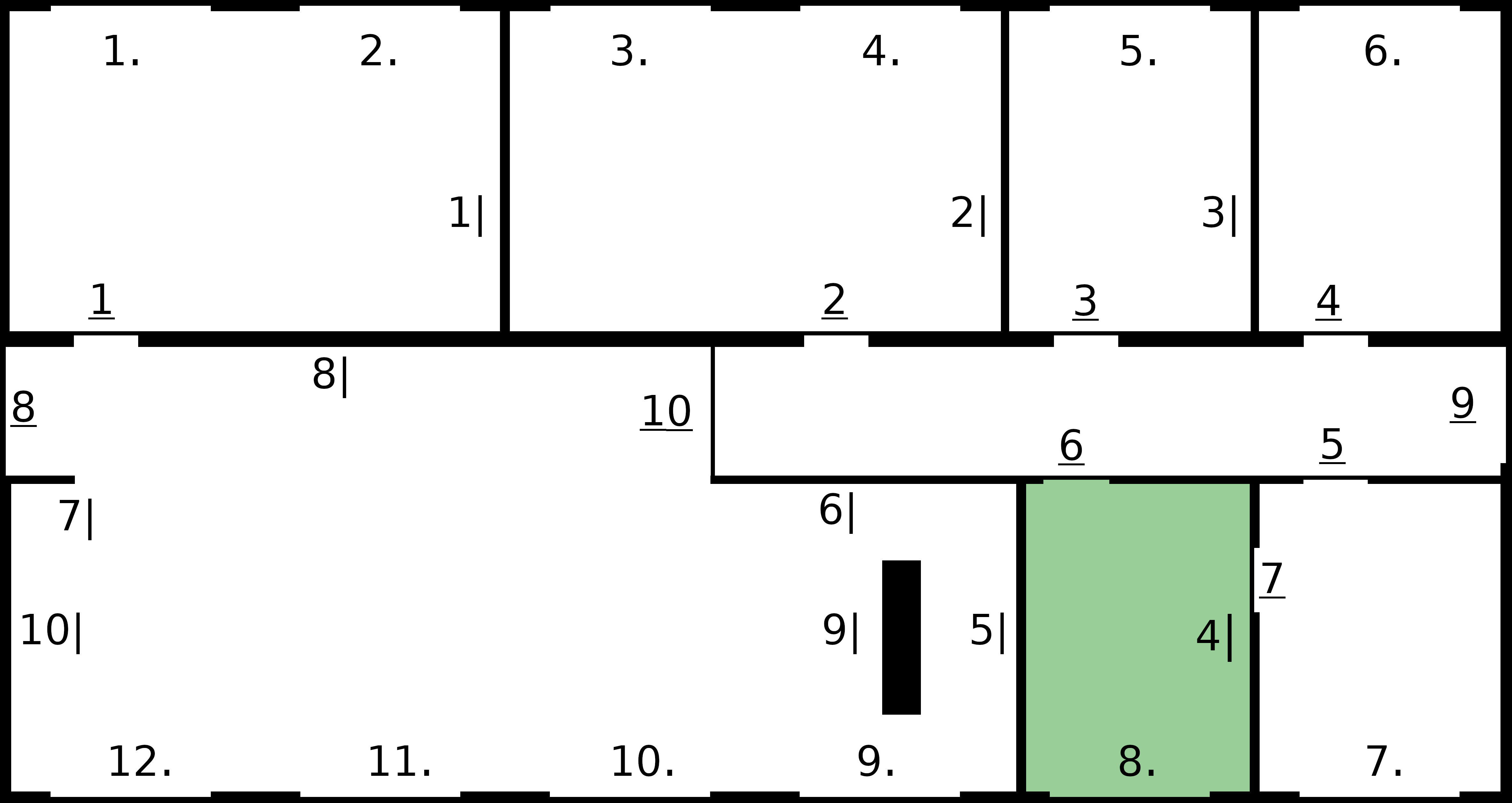}
	\end{subfigure}
	\centering
	\caption[Definition of the benchmark problem]{Benchmark problem with $\Omega := [0,2] \times [0,1] \subset \mathbb{R}^2$.
		Numbers indicate affine components, where $i.$ is a window, $\underbar{i}$ are doors, and $i|$ are walls.
		The $i$-th heater is located under the $i$-th window. With respect to~\eqref{eq:heat_equation}, we consider $\Gamma_\textnormal{{R}}:= \partial \Omega$, where
		$c_{\mu}$ contains outside wall $10|$, outside doors $\underbar{8}$ and $\underbar{9}$ and all windows. All other diffusion components enter $A_\mu$,
		whereas the heaters enter the source term $f_\mu$. We set $u_{\textnormal{out}}=5$ and
		the green region illustrates the domain of interest $D$.}
	\label{ex1:blueprint}
\end{SCfigure}

For all experiments, except the one in \cref{sec:thermalfin}, we consider as objective functional a weighted $L^2$-misfit on a domain of interest $D \subseteq \Omega$ and a weighted Tikhonov-term, i.e.
\begin{align} 
	\mathcal{J}(v, \mu) := \frac{\sigma_D}{2} \int_{D}^{} (v - u^{\textnormal{d}})^2 \dx + \frac{1}{2} \sum^{M}_{i=1} \sigma_i (\mu_i-\mu^{\textnormal{d}}_i)^2 + 1,
	\label{eq:mmexc_example_J}
\end{align}
with given desired state $u^{\textnormal{d}} \in V$, desired parameter $\mu^{\textnormal{d}} \in \Params$, and weights $\sigma_D, \sigma_i$ specified further below. The formulation of $\mathcal{J}(v, \mu)$ is a very general choice. It is applicable to design optimization, optimal control, and inverse problems.
We remark that the constant term $1$ is added to fulfill \cref{asmpt:bound_J} and does not influence the positions of the local minima. With respect to the formulation in~\eqref{eq:quadratic_J}, we have
\begin{align*}
	\Theta(\mu) &= \frac{\sigma_d}{2} \sum^{M}_{i=1} \sigma_i (\mu_i-\mu^{\textnormal{d}}_i)^2 + \frac{\sigma_d}{2} \int_{D}^{} u^{\textnormal{d}} u^{\textnormal{d}} \dx, \\
	j_{\mu}(v) &= -\sigma_d \int_{D}^{} u^{\textnormal{d}}v \dx, \qquad \text{and} \qquad k_{\mu}(v,v) = \frac{\sigma_d}{2} \int_{D}^{} v^2 \dx.
\end{align*}
Motivated by ensuring the desired temperature in a single room of a building,
we consider blueprints with windows, heaters, doors, and walls, yielding parameterized diffusion, forces, and boundary values as sketched in \cref{ex1:blueprint} and refer to the accompanying code for the definition of the data functions.
For simplicity, we omit realistic temperature modeling and restrict ourselves to academic numbers of the diffusion and heat source quantities.
We seek to match a desired temperature $u^\textnormal{d} \equiv 18$ and set $\mu^\textnormal{d}_i = 0$.
For the FOM discretization, we choose a cubic mesh that resolves all features of the data functions derived from \cref{ex1:blueprint}, resulting in $\dim V_h = 80601$ degrees of freedom.

Note that multiple parametric diffusion problems with a parameter dimension of up to $P=44$ can be constructed.
In the sequel, we use $2$-, $10$-, $12$- and $28$-dimensional examples and specify the details in the respective section.
Note that only in \cref{sec:thermalfin} we deviate from this model problem for the sake of a comparison to the case study used in~\cite{QGVW2017}.

\subsection{Experiment 1: Proof of concept}
\label{sec:proof_of_concept}

We start with a simple $2$-dimensional instance of the introduced benchmark problem with the main intention to demonstrate the method with additional visualizations and details.
We also note that the discussed experiment has been used to motivate TR-RB methods in \cref{sec:background_RB_methods_for_PDEopt} and to demonstrate further TR-RB approaches in \cref{sec:TRRB_related_approaches}.

In particular, we optimize over two sets of walls $\{1|,2|,3|,8|\}$ and $\{4|,5|,6|\}$, with $\Paramsad:= [0.005, 0.1]^2$.
Moreover, all doors are opened (meaning that their value is equal to the background), and the other non-parametric quantities are chosen appropriately.
Additionally, we define weights of $\J$ by $\sigma_D=1$ and $\sigma_0 = \sigma_1 = 10$, such that more emphasis lies on minimizing the cost of walls.

\begin{figure}
	\begin{subfigure}{.5\textwidth}
		\centering
		\includegraphics[width=.75\linewidth]{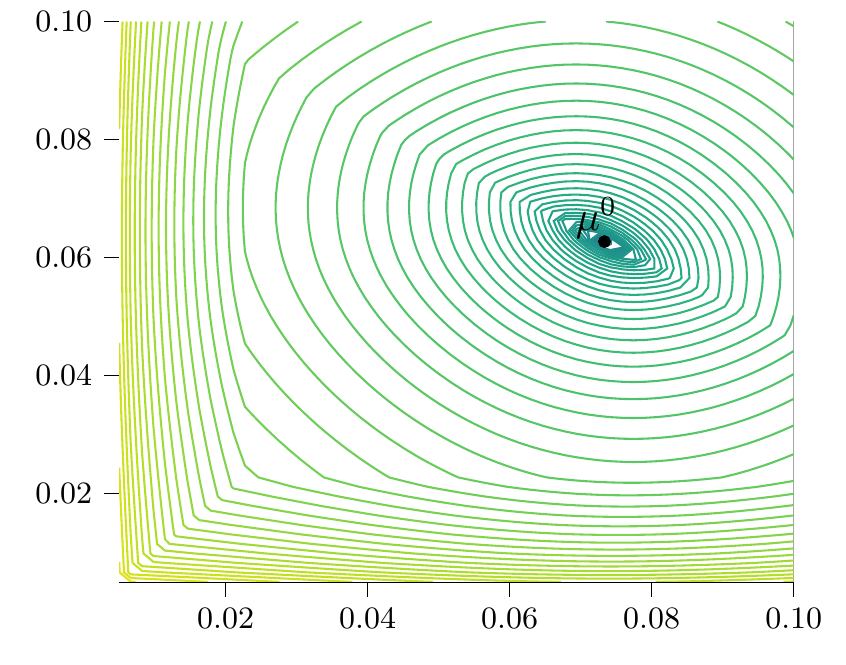}
		\caption{Estimator $\Delta_{\Jhat_\red^{(0)}}$}
		\label{fig:sfig1a}
	\end{subfigure}%
	\begin{subfigure}{.5\textwidth}
		\centering
		\includegraphics[width=.75\linewidth]{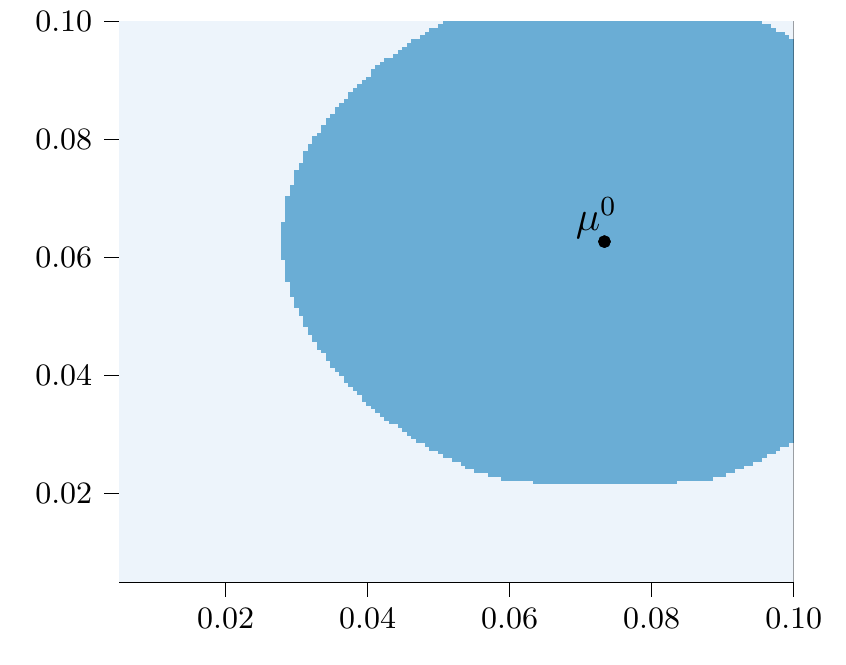}
		\caption{First TR after shrinking TR-radius}
		\label{fig:sfig1b}
	\end{subfigure}

	\begin{subfigure}{.5\textwidth}
		\centering
		\includegraphics[width=.75\linewidth]{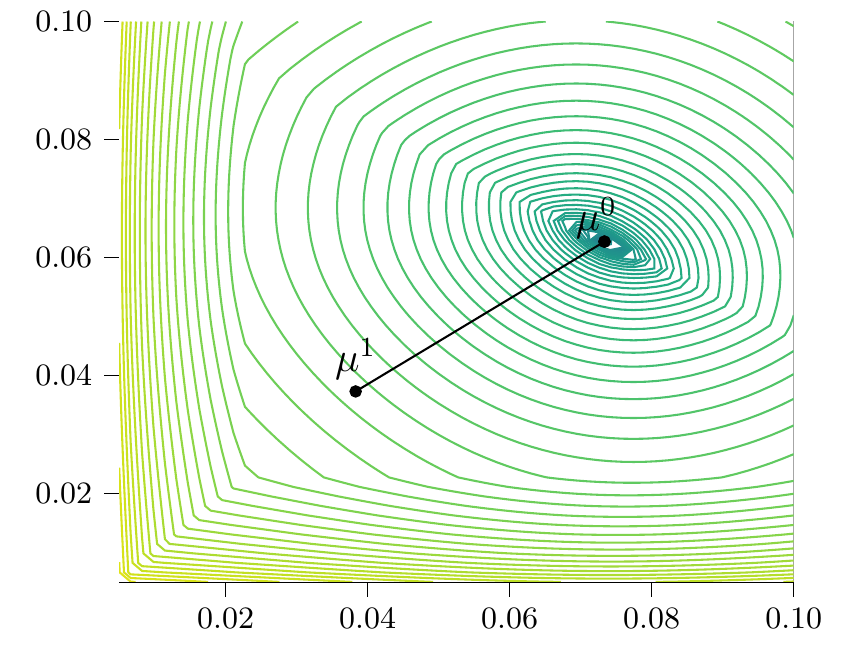}
		\caption{Estimator $\Delta_{\Jhat_\red^{(0)}}$}
		\label{fig:sfig2a}
	\end{subfigure}%
	\begin{subfigure}{.5\textwidth}
		\centering
		\includegraphics[width=.75\linewidth]{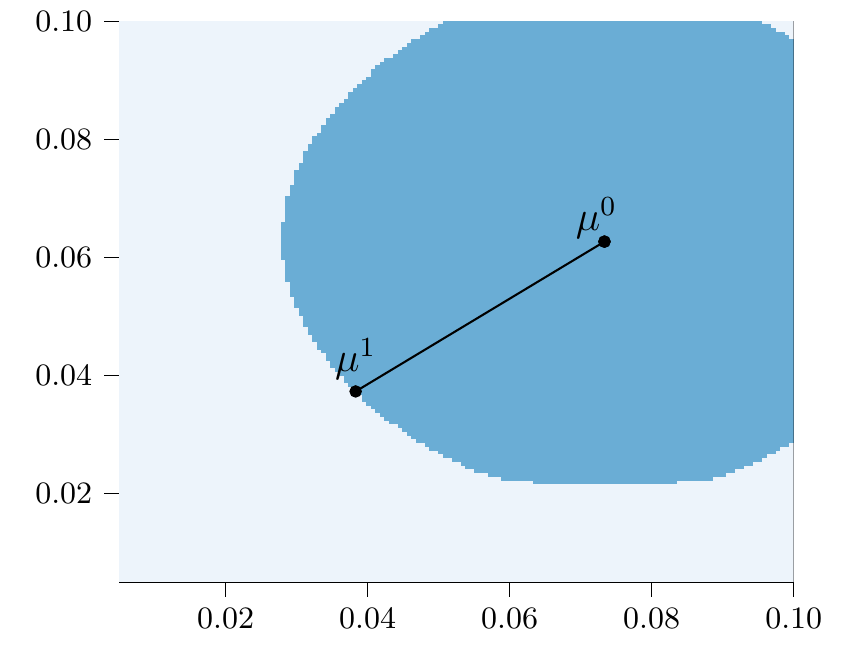}
		\caption{TR before enrichment}
		\label{fig:sfig2b}
	\end{subfigure}

	\begin{subfigure}{.5\textwidth}
	\centering
	\includegraphics[width=.75\linewidth]{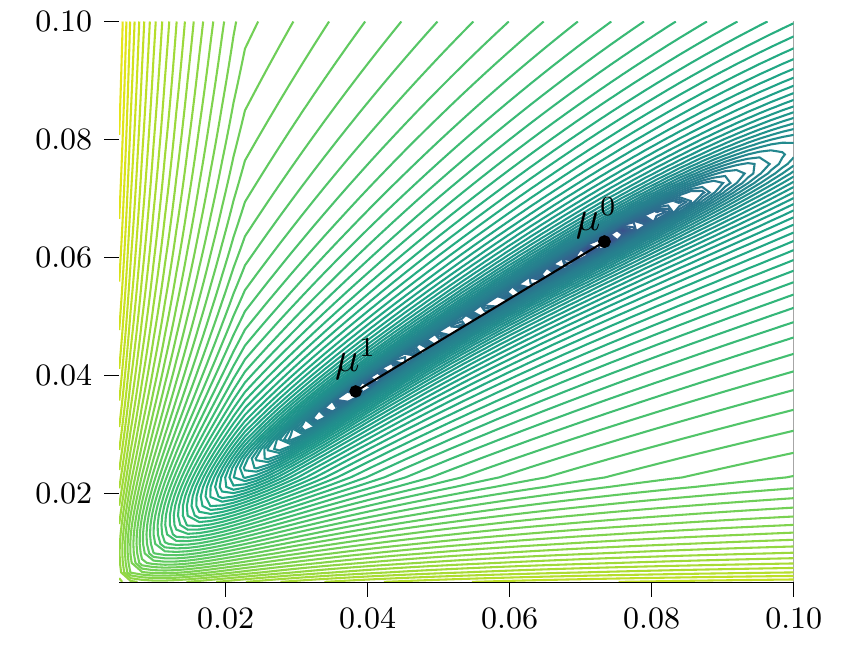}
	\caption{Estimator $\Delta_{\Jhat_\red^{(1)}}$}
	\label{fig:sfig3a}
	\end{subfigure}%
	\begin{subfigure}{.5\textwidth}
	\centering
	\includegraphics[width=.75\linewidth]{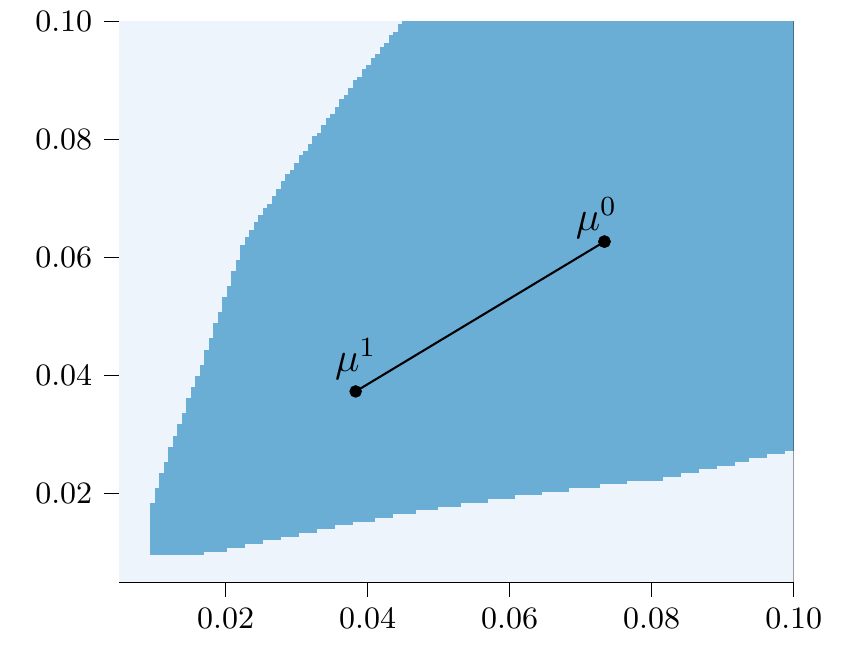}
	\caption{TR after enrichment}
	\label{fig:sfig3b}
	\end{subfigure}

	\begin{subfigure}{.5\textwidth}
	\centering
	\includegraphics[width=.75\linewidth]{Pictures_proof_of_concept/estimator_2.pdf}
	\caption{Estimator $\Delta_{\Jhat_\red^{(1)}}$}
	\label{fig:sfig5a}
\end{subfigure}%
\begin{subfigure}{.5\textwidth}
	\centering
	\includegraphics[width=.75\linewidth]{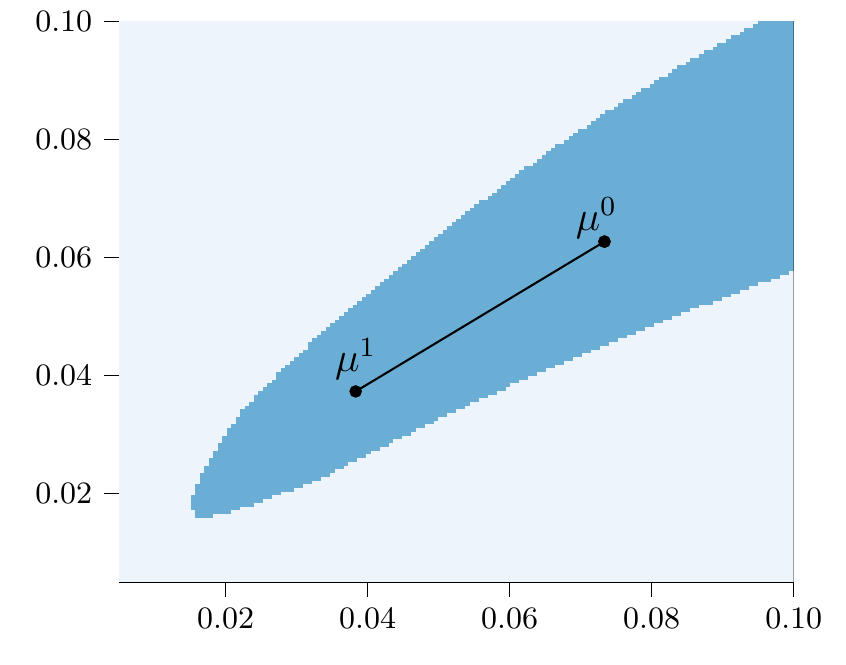}
	\caption{TR after shrinking TR-radius}
	\label{fig:sfig5b}
\end{subfigure}
	\caption{Visualization of iterations $k=0,1$ of the TR-RB procedure with emphasis on the a posteriori error estimator and the resulting Trust-regions.}
\label{fig:fig1}
\end{figure}

\begin{figure}
	\begin{subfigure}{.5\textwidth}
	\centering
	\includegraphics[width=.75\linewidth]{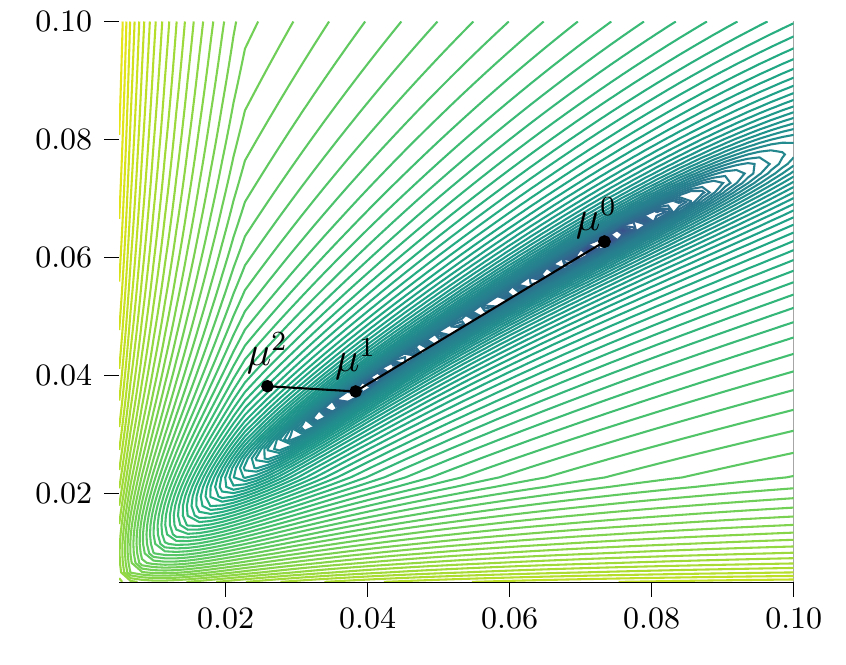}
	\caption{Estimator $\Delta_{\Jhat_\red^{(1)}}$}
	\label{fig:sfig6a}
\end{subfigure}%
\begin{subfigure}{.5\textwidth}
	\centering
	\includegraphics[width=.75\linewidth]{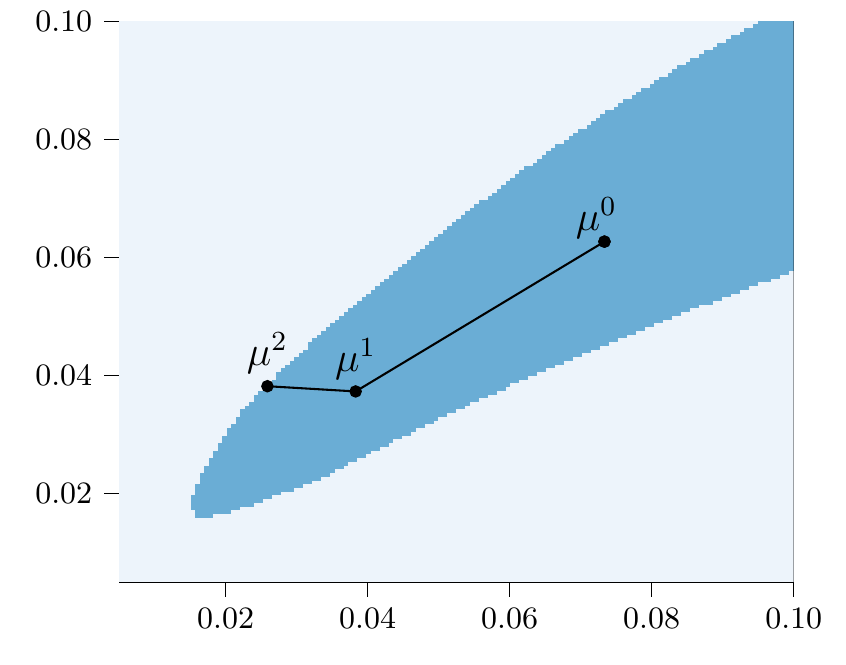}
	\caption{TR before enrichment}
	\label{fig:sfig6b}
\end{subfigure}

	\begin{subfigure}{.5\textwidth}
	\centering
	\includegraphics[width=.75\linewidth]{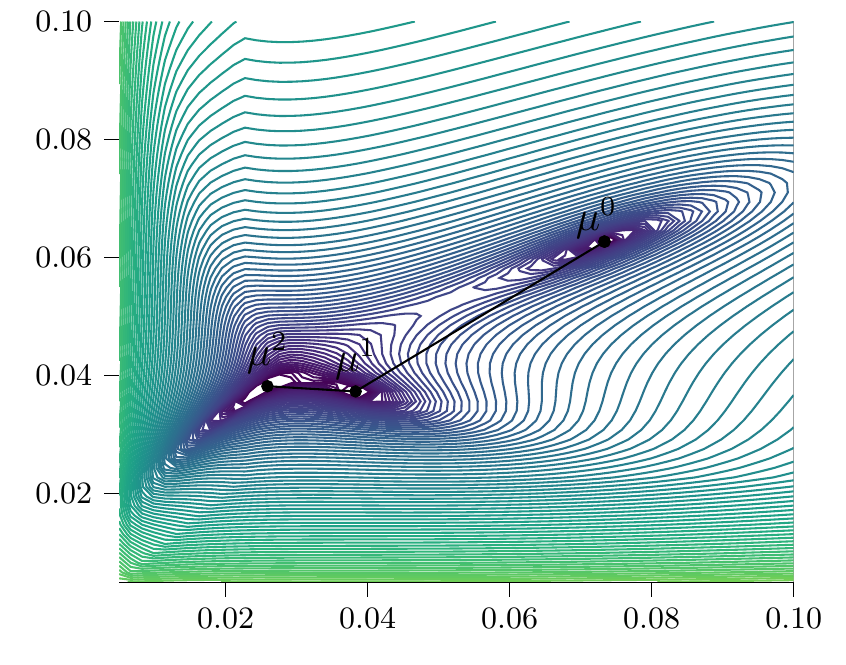}
	\caption{Estimator $\Delta_{\Jhat_\red^{(2)}}$}
	\label{fig:sfig7a}
\end{subfigure}%
\begin{subfigure}{.5\textwidth}
	\centering
	\includegraphics[width=.75\linewidth]{Pictures_proof_of_concept/second_trust_region_2.pdf}
	\caption{TR before enrichment}
	\label{fig:sfig7b}
\end{subfigure}

	\begin{subfigure}{.5\textwidth}
	\centering
	\includegraphics[width=.75\linewidth]{Pictures_proof_of_concept/estimator_3.pdf}
	\caption{Estimator $\Delta_{\Jhat_\red^{(2)}}$}
	\label{fig:sfig8a}
\end{subfigure}%
\begin{subfigure}{.5\textwidth}
	\centering
	\includegraphics[width=.75\linewidth]{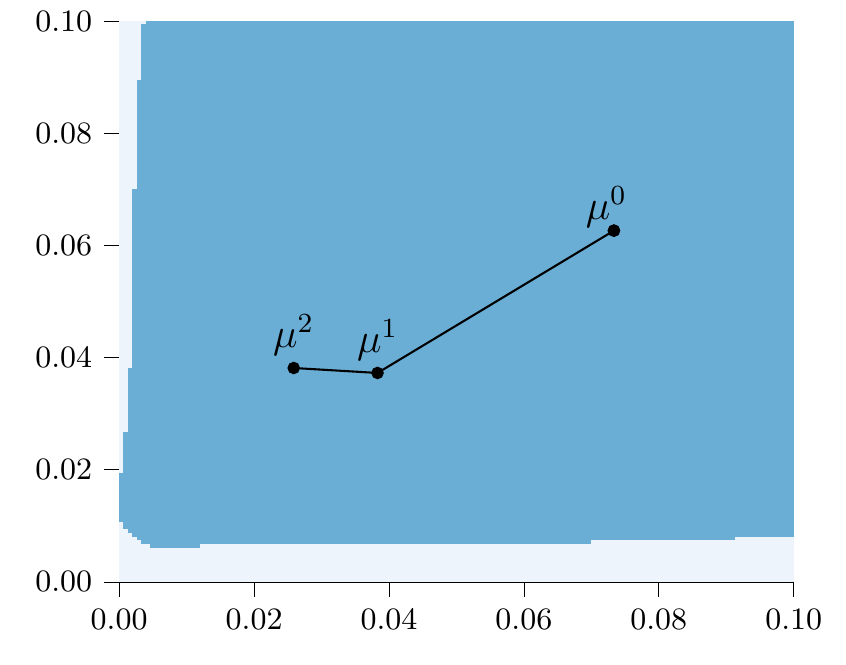}
	\caption{TR after enrichment}
	\label{fig:sfig8b}
\end{subfigure}

	\begin{subfigure}{.5\textwidth}
	\centering
	\includegraphics[width=.75\linewidth]{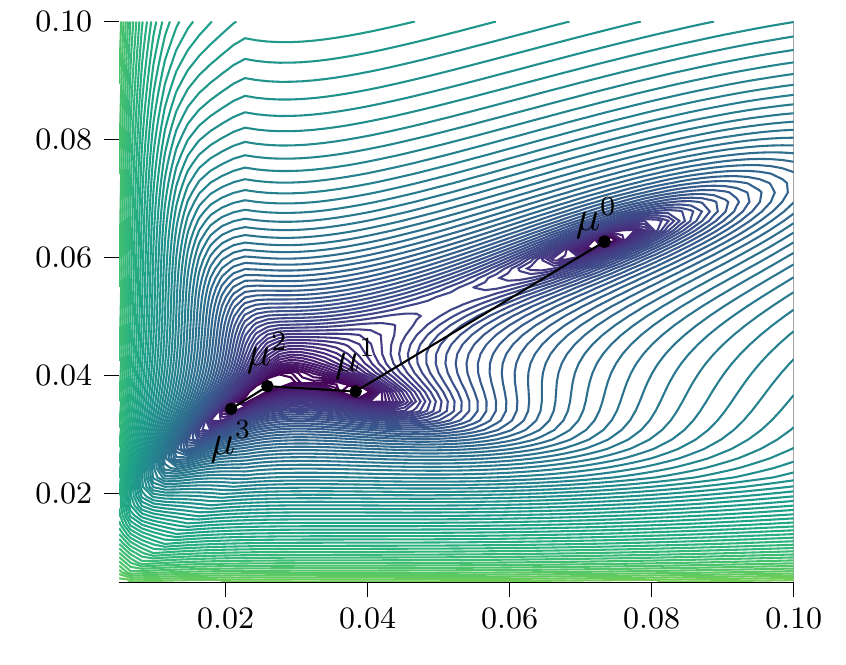}
	\caption{Estimator $\Delta_{\Jhat_\red^{(3)}}$}
	\label{fig:sfig9a}
\end{subfigure}%
\begin{subfigure}{.5\textwidth}
	\centering
	\includegraphics[width=.75\linewidth]{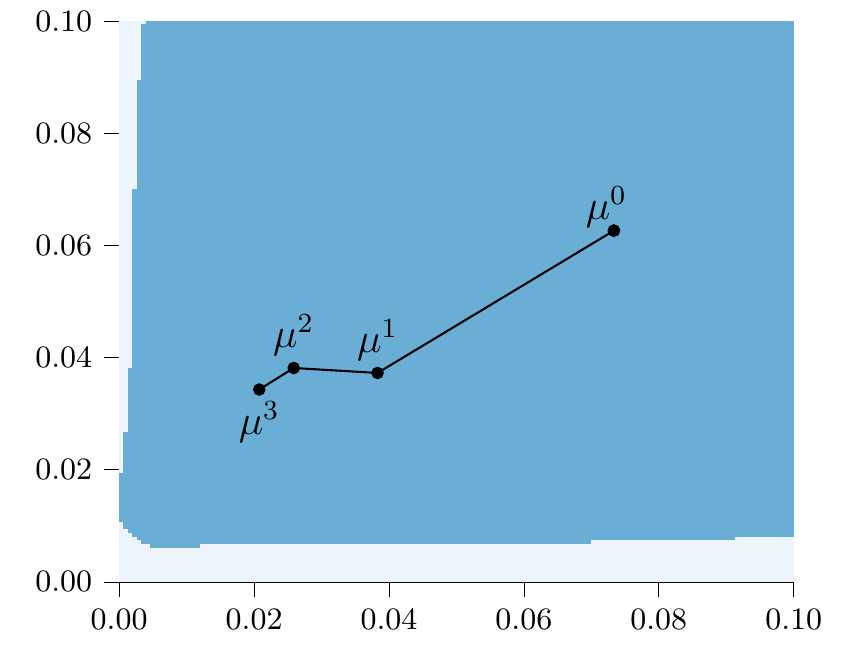}
	\caption{Final TR}
	\label{fig:sfig9b}
\end{subfigure}
	\caption{Visualization of iterations $k=2,3$ of the TR-RB procedure with emphasis on the a posteriori error estimator and the resulting Trust-regions.}
\label{fig:fig2}
\end{figure}

We apply the NCD corrected TR-RB algorithm with unconditional Lagrange enrichment 3(a)(UE) and BFGS as the sub-problem solver.
The detailed steps are illustrated in Figures~\ref{fig:fig1} and~ \ref{fig:fig2}, which we discuss in the following.
Starting with a randomly chosen initial guess $\mu^0 = (0.073, 0.063)$, we begin by initializing the Lagrange RB spaces with $ V_\red^{\pr,0} = \{u_{h, \mu^0}\}$ and $V_\red^{\du,0} = \{p_{h, \mu^0}\}$, modulo the respective orthonormalization.
Then, we solve the first TR-sub-problem \eqref{TRRBsubprob} with conditions~\eqref{Termination_crit_sub-problem}.
In this particular case, the potential iterates are rejected four times according to \cref{sec:acceptance_of_TR} such that the initial TR radius $\delta^{(0)}=0.1$ is subsequently shrunk to $\delta^{(1)}=0.00625$.
In \cref{fig:fig1}(a) and (b), we visualize the initial a posteriori error estimator for $\Jhat_\red^{(0)}$ and the first resulting trust-region for which the resulting iterate is accepted.
It can be seen in \cref{fig:fig1}(d) that the sub-problem, indeed, stops at the boundary of the TR, which indicates that the optimum is located outside of the TR.
Since the iterate $\mu^1$ is accepted, the primal and dual RB spaces are enriched, which results in a new shape of the estimator for $\Jhat_\red^{(1)}$, forming a valley along both iterates.

Again, the trust-region with radius $\delta^{(1)}=0.0063$ (visualized in \cref{fig:fig1}(f)) is too large and hence, for the next iteration, the radius is shrunk another three times, resulting in a radius of $\delta^{(2)}=0.00078$ with corresponding TR in \cref{fig:fig1}(h).
The second iterate $\mu^{(2)}$ lies on the boundary of the TR. 
After $\mu^{(2)}$ has been found and the RB spaces have been enriched, the criterion for enlarging the TR-radius~\eqref{TR_act_decrease} is fulfilled, such that we increase the radius to $\delta^{(2)}=0.0016$.

It can be seen in \cref{fig:fig2}(f) that the resulting TR covers a vast amount of the admissible set.
Consequently, the optimum likely lies inside this region.
Luckily, the TR does not need to be shrunk again and, instead, the final iterate~$\mu^{(3)}$ of the TR-RB is found, which is also the first iterate that has been deduced by triggering the (non-TR-related) stopping criterion \cref{FOC_sub-problem}.
Considering \cref{fig:fig2}(g), we see the valley of the resulting error estimator, which again indicates that our ROM is following the optimization path.

\begin{table}[th] \centering \footnotesize
	\begin{tabular}{|c||c|c|c|c|c|c|} \hline 
		it. $k$ & TR-radius $\delta^{(k)}$ & Inner iterations & TR stopped with & Rejected by & Basis size & $g_h(\mu^{(k)})$\\ \hline
		0 & 0.1 & 4 & \eqref{Cut_of_TR} & \eqref{eq:suff_decrease_condition} & 1 & - \\
		0 & 0.05 & 4 & \eqref{Cut_of_TR} & \eqref{eq:suff_decrease_condition} & 1 & - \\
		0 & 0.025 & 5 & \eqref{Cut_of_TR} & \eqref{eq:suff_decrease_condition} & 1 & - \\		
		0 & 0.00625 & 3 & \eqref{Cut_of_TR} & - & 1 & 4.64e-2 \\ \hline
		1 & 0.00625 & 4 & \eqref{Cut_of_TR} & \eqref{eq:suff_decrease_condition} & 2 & - \\
		1 & 0.00312 & 3 & \eqref{Cut_of_TR} & \eqref{eq:suff_decrease_condition} & 2 & - \\
		1 & 0.00156 & 3 & \eqref{Cut_of_TR} & \eqref{eq:suff_decrease_condition} & 2 & - \\
		1 & 0.00078 & 3 & \eqref{Cut_of_TR} & - & 2 & 3.93e-2\\ \hline
		2 & 0.00078 & 7 & \eqref{FOC_sub-problem} & - & 3 & 1.85e-3 \\ \hline
		3 & 0.00156 & 7 & \eqref{FOC_sub-problem} & - & 4 & 4.18e-7 \\ \hline
	\end{tabular}
	\caption[Experiment 1: Details on the TR-RB iterations]{Details on the TR-RB iterations for the 2-dimensional proof of concept experiment. Depicted is the TR-radius, the triggered local termination and rejection criterion of the sub-problem, the resulting basis sizes and the global termination criterion.}
	\label{tab:poc_details}
\end{table}

We point to \cref{tab:poc_details} for a detailed output of the above-explained TR-RB iterations, where iteration counts, TR-radius sizes, inner iterations, basis sizes, termination criteria, and further details on the rejection of the iterate are stated.
It can be seen that, especially at the beginning of the algorithm, the TR-radius does not fit the approximation quality of the ROM.
Thus, the TR radius is shrunk accordingly, and the first iterates are rejected by the expensive condition~\eqref{eq:suff_decrease_condition}.
As can be seen in \cref{fig:fig1}(d) and \cref{fig:fig2}(b), for $k=0,1$, the sub-problem is stopped by reaching the boundary of the TR (condition \eqref{Cut_of_TR}) and, for $k=2,3$, the solution of the sub-problem terminates by the actual convergence criterion \eqref{FOC_sub-problem}.

We again note that these results are re-discussed for deriving further approaches of the TR-RB algorithm in \cref{sec:TRRB_related_approaches}.

\subsection{BFGS variants and their comparison to existing approaches}
\label{sec:exp_paper1}

This section consists of the experiments that are presented in~\cite{paper1}.
The primary purpose of this section is to analyze the different ROM strategies $1.$-$3.$ concerning their numerical accuracy and to use these ROM strategies in the TR-RB algorithm with basis enrichment strategies (a) and (b); cf. \cref{sec:TRRB_variants}.
Hence, only BFGS variants with unconditional enrichment are used.
Another purpose is to compare our method to the originally proposed TR-RB algorithm from~\cite{QGVW2017} (also by considering the same case study), as well as to existing FOM approaches that do not use adaptive surrogate models at all.

To be precise, we consider the following BFGS-based state-of-the-art methods:
\begin{description}
\item[FOM projected BFGS:]
We consider the standard projected BFGS method, stated in~\cite[Section~5.5.3]{kelley}, which uses FOM evaluations of the forward model for the reduced cost functional and its gradient; cf. \cref{sec:TRRB_fom}.
We restrict the number of iterations by $400$.
\item[TR-RB from~\cite{QGVW2017}:] \hphantom{text}\newline
We consider the same method as in~\cite{QGVW2017}, where the authors used the standard output functional and gradient from \cref{sec:TRRB_standard_approach} (ROM-strategy $1.$) and utilized Lagrangian RB enrichment (a).
For the sub-problems, the authors used the BFGS method.
Furthermore, no enlarging strategy has been used for the TR-radius, and no projection for parameter constraints has been concealed.
Notably, the authors did not take advantage of the fact that, for the outer stopping criterion, the full-order FOC condition~\eqref{eq:TR_termination} is cheaply available after an enrichment step.
Instead, they used the ROM-based FOC condition plus the estimator for the gradient of the
cost functional $\|\noncorgrad_\mu \Jnoncor_\red(\mu^{(k+1)}))\|_2 + \Delta_{\noncorgrad_\mu \Jnoncor_\red}(\mu) \leq \tau_\textnormal{{FOC}}$ as a replacement for~\eqref{eq:TR_termination}.
Note that this approach has multiple drawbacks.
First, the evaluation is more costly due to the estimator.
Second, it is less accurate.
Third, the offline time for the enrichment increases since the estimator needs to be available.
Fourth, it can prevent the TR-RB from converging if the estimator is not small enough (for instance, governed by large overestimation or numerical issues in the estimator).
\end{description}

For the sake of completeness, we compare with the following TR-RB strategies from \cref{sec:TRRB_variants} and introduce corresponding abbreviations that are used in the sequel.
We also refer to \cref{alg:TR-RBmethod_paper1}.

\begin{description}
\item[1(a) TR-RB Lag.] We consider the standard ROM Strategy 1(a) with BFGS as the sub-problem solver and unconditional Lagrangian enrichment (UE).
\item[1(b) TR-RB Agg.] We consider the standard ROM Strategy 1(b) with BFGS as the sub-problem solver and unconditional aggregated enrichment (UE).
\item[2(a) TR-RB Lag.] We consider the semi-corrected ROM Strategy 2(a) with BFGS as the sub-problem solver and unconditional Lagrangian enrichment (UE).
\item[3(a) TR-RB Lag.] We consider the NCD-corrected ROM Strategy 3(a) with BFGS as the sub-problem solver and unconditional Lagrangian enrichment (UE).
\end{description}

\subsubsection{Experiment 2: Elliptic thermal fin model problem}
\label{sec:thermalfin}

\begin{SCfigure}[50][ht]
	\centering\footnotesize
	\begin{tikzpicture}[scale=1.1]
		\draw (-.5,0) -- (-.5,.75);
		\draw (-.5,.75) -- (-3.,.75);
		\node at (-1.8,0.86) {{\tiny $k_1$}};
		\node at (-1.8,1.86) {{\tiny $k_2$}};
		\node at (-1.8,2.86) {{\tiny $k_3$}};
		\node at (-1.8,3.86) {{\tiny $k_4$}};
		\node at (1.8,0.86) {{\tiny $k_1$}};
		\node at (1.8,1.86) {{\tiny $k_2$}};
		\node at (1.8,2.86) {{\tiny $k_3$}};
		\node at (1.8,3.86) {{\tiny $k_4$}};
		\node at (0,2) {{\tiny $k_0$}};
		\draw (-3.,.75) -- (-3.,1) -- (-.5,1);
		\draw (-.5,1) -- (-.5,1.75) -- (-3.,1.75) -- (-3.,2) -- (-.5,2);
		\draw (-.5,2) -- (-.5,2.75) -- (-3.,2.75) -- (-3.,3) -- (-.5,3);
		\draw (-.5,3) -- (-.5,3.75) -- (-3.,3.75) -- (-3.,4) -- (3,4);
		\draw (.5,0) -- (.5,.75) -- (3.,.75) -- (3.,1) -- (.5,1);
		\draw (.5,1) -- (.5,1.75) -- (3.,1.75) -- (3.,2) -- (.5,2);		
		\draw (.5,2) -- (.5,2.75) -- (3.,2.75) -- (3.,3) -- (.5,3);		
		\draw (.5,3) -- (.5,3.75) -- (3.,3.75) -- (3.,4);		
		\draw (-.5,0) -- (.5,0) node [midway,above] {$\Gamma_\textnormal{{N}}$};
		\draw[densely dotted] (-.5,.75) -- (-.5,1);
		\draw[densely dotted] (-.5,1.75) -- (-.5,2);
		\draw[densely dotted] (-.5,2.75) -- (-.5,3);
		\draw[densely dotted] (-.5,3.75) -- (-.5,4);
		\draw[densely dotted] (.5,.75) -- (.5,1);
		\draw[densely dotted] (.5,1.75) -- (.5,2);
		\draw[densely dotted] (.5,2.75) -- (.5,3);
		\draw[densely dotted] (.5,3.75) -- (.5,4);
		\draw[->|] (3.2,0.45) -- (3.2,0.75);
		\draw[->|] (3.2,1.3) -- (3.2,1.);
		\draw[<->|] (0.5,0.45) -- (3,0.45) node[midway, fill=white] {$L$};
		\draw (3,0.75) -- (3,1) node[midway,right] {\small $t$};
	\end{tikzpicture} 
	\caption[Problem definition of the thermal fin example from Section~\ref{sec:thermalfin}]{%
		The thermal fin example from Section~\ref{sec:thermalfin}.
		Depicted is the spatial domain $\Omega$ (with $L=2.5$ and $t=0.25$) with Neumann boundary at the bottom with $|\Gamma_{\textnormal{N}}| = 1$ and Robin boundary $\Gamma_\textnormal{{R}} := \partial\Omega \backslash \Gamma_\textnormal{{N}}$, as well as the values $k_0, \dots, k_4 > 0$ of the diffusion $A_\mu$, which is piecewise constant in the respective parts of the domain.
		The other functions in~\eqref{eq:heat_equation} are given by $f_\mu = 0$, $g_\textnormal{{N}} = -1$, $u_\textnormal{{out}} = 0$ and $c_\mu = \textnormal{Bi} \in \R$, the Biot number.
		We allow to vary the six parameters $(k_0, \dots, k_4, \textnormal{Bi})$ and define the set of admissible parameters as $\Paramsad := [0.1, 10]^5 \times [0.01, 1] \subset \R^P$ with $P = 6$.
		We choose $\check{\mu} = (1, 1, 1, 1, 1, 0.1)$ for the energy product.
	}
	\label{fig:fin}
\end{SCfigure}
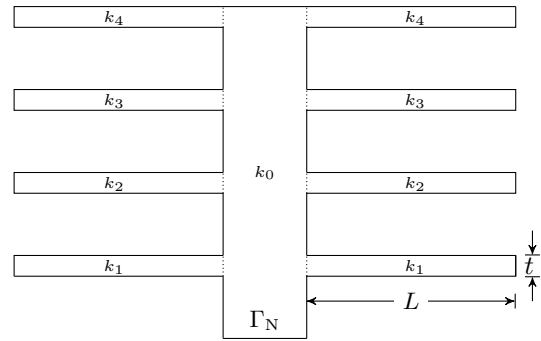

We consider the six-dimensional elliptic thermal fin example from~\cite[Sec.~5.1.1]{QGVW2017} and refer to \cref{fig:fin} for the problem definition.
The purpose of this experiment is to show the applicability of the proposed algorithms and to compare them to the one presented in~\cite{QGVW2017}.
For all optimization runs we prescribe the same desired parameter $\mu^\textnormal{d} \in \Paramsad$ by randomly drawing $k_1, \dots, k_4$ strictly within $\Paramsad$ and by setting $k_0 = 0.1$ and $\textnormal{Bi} = 0.01$, to artificially mimic the situation where parameter constraints have to be tackled.
Defining $T^\textnormal{d} := q(u_{h,\mu^\textnormal{d}})$, where $u_{h,\mu^\textnormal{d}} \in V_h$ is the solution of~\eqref{eq:state_h} associated with the desired parameter and where $q(v) := \int_{\Gamma_\textnormal{{N}}} v \ds$ for $v\in V$ denotes the mean temperature at the root of the fin, we consider a cost functional $\J(u, \mu) := \Theta(\mu) + j_\mu(u) + k_\mu(u, u)$ as in~\eqref{eq:quadratic_J} with $\Theta(\mu) := (\|\mu^\textnormal{d} - \mu\| / \|\mu^\textnormal{d}\|)^2 + {T^\textnormal{d}}^2 + 1$, $j_\mu(v) := - T^\textnormal{d}\,q(v)$ and $k_\mu(v, w) := 1/2\,q(v)\,q(w)$.
Note that this choice slightly deviates from the functional that is explained in \cref{sec:TRRB_benchmark_problem}.
We would like to point out that the authors in~\cite{QGVW2017} dropped the ${T^\textnormal{d}}^2 + 1$ term from the definition of $\Theta$, which we re-add to ensure \cref{asmpt:bound_J}.
This constant term does not change the position of local minima and the derivatives of the cost functional.
However, it makes the trust-region radius shrink, especially initially, slowing down the TR-RB methods.
This does not affect the comparison among the TR-RB methods since all of them suffer from this issue.
Note that for this particular example, the proposed NCD-correction term vanishes; see \cref{rem:fin}.
For the FOM, we generate an unstructured simplicial mesh using \texttt{pyMOR}'s \texttt{gmsh} (see~\cite{gmsh}) bindings, resulting in $\dim V_h = 77537$.

Starting with ten different randomly drawn initial parameters $\mu^{(0)}$, we measure the total computational run time, the number of TR iterations $k$, and the error in the optimal parameter.

\begin{SCfigure}
	\centering\footnotesize
\begin{tikzpicture}[]

\definecolor{color0}{rgb}{0.65,0,0.15}
\definecolor{color1}{rgb}{0.84,0.19,0.15}
\definecolor{color2}{rgb}{0.96,0.43,0.26}
\definecolor{color3}{rgb}{0.99,0.68,0.38}
\definecolor{color4}{rgb}{1,0.88,0.56}
\definecolor{color5}{rgb}{0.67,0.85,0.91}

\begin{axis}[
legend cell align={left},
legend style={fill opacity=0.8, draw opacity=1, text opacity=1, at={(0.99,0.01)}, anchor=south east, draw=white!80!black},
log basis y={10},
tick align=outside,
tick pos=left,
x grid style={white!69.0196078431373!black},
xlabel={time in seconds [s]},
xmajorgrids,
xmin=-5.7143321871758, xmax=236.000975930691,
xtick style={color=black},
y grid style={white!69.0196078431373!black},
ylabel={$\big\| \overline{\mu}-\mu^{(k)} \big\|_2$},
ymajorgrids,
ymin=1.91901353329293e-07, ymax=27.8370378432078,
ymode=log,
ytick style={color=black},
y=0.475cm
]
\addplot [semithick, color0, mark=triangle*, mark size=3, mark options={solid, fill opacity=0.5}]
table {%
0 9.99473140636549
4.34884524345398 7.35446882820324
18.6473112106323 7.3541512685285
30.8121237754822 7.35749530443356
32.9880452156067 6.31413726386274
47.2815909385681 6.31413448410398
58.816871881485 6.31207321123728
61.0466980934143 5.42208140924561
75.8011300563812 5.42204355619545
87.570513010025 5.43839223519159
93.4977016448975 5.12334963173444
108.184760808945 5.12324262569418
118.464211463928 5.12391766323573
122.161113739014 6.18793044571934
136.850122451782 6.18787453046038
139.825917005539 4.10890037737276
142.081712961197 3.43536524159338
145.785378932953 2.91800394897423
160.438888788223 2.9179971551271
162.668388605118 2.07914534719756
164.919844388962 1.947091236104
179.623940706253 1.9470783493143
181.859387397766 0.422463567698884
184.128856897354 0.643786083343712
186.372574567795 0.00702386653971826
201.036522626877 0.00702382140842531
};
\addlegendentry{\hspace{16pt} FOM proj.~BFGS}
\addplot [semithick, color5, mark=diamond*, mark size=3, mark options={solid, fill opacity=0.5}]
table {%
3.60644197463989 9.99473140636549
6.18412208557129 8.91952684680866
10.1918063163757 8.91275320910397
14.3757772445679 8.91255198939559
19.1285471916199 8.9086116684482
24.5155129432678 8.90683550331313
29.9072082042694 8.45403937547982
35.5541481971741 6.77948120396427
41.6052765846252 6.3782435509005
48.167160987854 3.11373071399215
55.3844101428986 1.40622934454799e-05
};
\addlegendentry{\hspace{19pt}TR-RB from \cite{QGVW2017}}
\addplot [semithick, color1, mark=triangle*, mark size=3, mark options={solid,rotate=180, fill opacity=0.5}]
table {%
2.86481094360352 9.99473140636549
5.50435853004456 8.91952684680866
8.40358209609985 8.91275320910397
11.4713573455811 8.91255198939559
15.107709646225 8.90773079070491
19.0532963275909 8.34476797197126
23.4172186851501 7.18301515662573
28.4373028278351 1.71823629414235
34.2188301086426 7.84148139372534e-06
};
\addlegendentry{1(a)\hspace{1pt} TR-RB Lag.}
\addplot [semithick, color2, mark=*, mark size=3, mark options={solid, fill opacity=0.5}]
table {%
2.82310795783997 9.99473140636549
5.37446165084839 8.91952684680867
8.26291704177856 8.91275320910398
11.3642964363098 8.9125519893956
14.8835608959198 8.90773079070492
18.759311914444 8.34476797199037
23.1022231578827 7.18301515637139
28.0834274291992 1.71823629453389
34.0882914066315 7.84148159296135e-06
};
\addlegendentry{1(b) TB-RB Agg.}
\end{axis}

\end{tikzpicture}
	\caption[Experiment 2: Error decay and performance of selected algorithms]{%
		Error decay and performance of selected algorithms for the example from \cref{sec:thermalfin} for a single optimization run with random initial guess $\mu^{(0)}$ for $\tau_\textnormal{{FOC}} = 5\cdot 10^{-4}$: for each algorithm each marker corresponds to one (outer) iteration of the optimization method and indicates the absolute error in the current parameter, measured against the known desired optimum $\bar{\mu} = \mu^\textnormal{d}$.
		In all except the FOM variant, the ROM is enriched in each iteration corresponding to \cref{alg:TR-RBmethod_paper1}, depending on the variant in question.
	}
	\label{fig:fin_timing_mu_d}
\end{SCfigure}
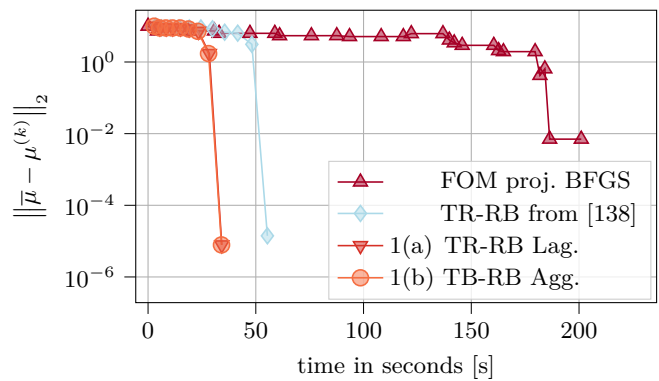

All considered optimization methods converged (up to a tolerance), but we restrict the presentation to the most informative ones (complete results can be found in the accompanying code).
As we observe from \cref{table:fin}, the ROM-based adaptive TR-RB algorithms vastly outperform the FOM variant, noting that the computational time of the ROM variants includes all offline and online computations.
\cref{fig:fin_timing_mu_d} details the decay of the error in the optimal parameter during the optimization for a selected random initial guess.
We observe that, for this example, the choice of the RB enrichment does not impact the performance of the algorithm too much, see \cref{rem:fin}.
Also, Strategies 2(a) and 3(a) show a comparable computational speed (not presented).
We also observe that the method from~\cite{QGVW2017} requires more time and more iterations on average, 
Variants 1(a) and (b) are still faster due to the enlarging of the TR-radius and the use of a termination criterion that does not depend on a posteriori estimates, preventing from additional TR iterations.

\begin{table}
	\centering\footnotesize
	\begin{tabular}{|l|cc|c|cc|} \hline
		& Av.~(min/max) Run time[s]
		& Speedup
		& Av.~(min/max) it. \kern-8pt
		& Rel.~error
		& FOC cond.\\
		\hline
		FOM proj.~BFGS
		& 967.86 (176.69/3401.06)
		& --
		& 111.20 (25/400)
		& $3.13\cdot 10^{-3}$
		& $1.19\cdot 10^{-2}$\\
		TR-RB from~\cite{QGVW2017}
		& 68.06 (43.28/88.21)
		& 10.40
		& 9.20 (8/13)
		& $1.34\cdot 10^{-6}$
		& $4.31\cdot 10^{-5}$\\
		1(a) TR-RB Lag.
		& 44.56 (34.22/74.96)
		& 21.72
		& 8.80 (8/11)
		& $3.08\cdot 10^{-6}$
		& $4.64\cdot 10^{-5}$\\
		1(b) TR-RB Agg.
		& 43.86 (34.09/74.35)
		& 22.07
		& 8.70 (8/10)
		& $3.37\cdot 10^{-6}$
		& $6.40\cdot 10^{-5}$ \\ \hline
	\end{tabular}
	\vspace{0.15cm}
	\caption[Experiment 2: Performance and accuracy of selected BFGS algorithms]{%
		Performance and accuracy of selected algorithms for the example from \cref{sec:thermalfin} for ten optimization runs with random initial guesses $\mu^{(0)}$: average, minimum and maximum total computational time (column 2) and speed-up compared to the FOM variant (column 3); average, minimum and maximum number of iterations $k$ required until convergence (column 4), the average relative error in the parameter (column 5) and average FOC condition at termination(column 6).
	}
	\label{table:fin}
\end{table}

\begin{remark}[Vanishing NCD-correction for the fin problem]
	\label{rem:fin}
	It is important to notice that this model problem is not suitable for fully demonstrating the NCD-corrected approach's capabilities.
	The reason is that the functional choice is a misfit on only the root edge of the thermal fin, plus a Tikhonov regularization term.
	Since the root of the thermal fin is also the source of the primal problem, the dual solutions $p_{\red, \mu}$ of the reduced dual equation \eqref{eq:dual_solution_red} are linearly dependent on the respective primal solutions $u_{\red, \mu}$ and the correction term $r_\mu^\pr(u_{\red,\mu})[p_{\red,\mu}]$ for the NCD-corrected RB reduced functional from~\eqref{eq:Jhat_red_corected} vanishes.
	In general, this is clearly not the case for quadratic objective functionals.
\end{remark}

\subsubsection{Experiment 3: Validation of the a posteriori error estimates}
\label{sec:TRRB_estimator_study}

We now use a ten-dimensional parameter example of the benchmark problem to, first, study the a posteriori error estimates from \cref{sec:TRRB_a_post_error_estimates} for the different choices of ROMs from \cref{sec:TRRB_MOR_for_PDEconstr}, and second, inspect the different TR-RB variants in \cref{sec:mmexc_opt_results}.

We consider the benchmark problem (cf. \cref{sec:TRRB_benchmark_problem}) with three wall sets $\{1|,2|,3|,8|\}$, $\{4|,5|,6|,7|\}$ and $\{9|\}$ and seven heater sets, $\{1,2\}$, $\{3,4\}$, $\{5\}$, $\{6\}$, $\{7\}$, $\{8\}$ and $\{9,10,11,12\}$ (each set governed by a single parameter component).
The set of admissible parameters is given by $\Paramsad= [0.025,0.1]^3\times[0,100]^7$ and we choose $\sigma_D= 100$ and $(\sigma_i)_{1\leq i\leq 10} = (10\sigma_w,5\sigma_w,\sigma_w,2\sigma_h,$ $2\sigma_h,\sigma_h,\sigma_h,\sigma_h,\sigma_h,4\sigma_h)$ in~\eqref{eq:mmexc_example_J}, with $\sigma_w = 0.05$ and $\sigma_h= 0.001$.
The choice of $\sigma_i$ is related to the measure of the walls and how many heaters are considered in each group.
The other components of the data functions are fixed and thus not directly involved in the optimization process.
Briefly, the diffusion coefficient of air and inside doors is set to $0.5$, of the outside wall to $0.001$, of exterior doors $\underline{8}$ and $\underline{9}$ to $0.01$ and of windows to $0.025$.
For the energy product, we choose $\check{\mu}= (0.05,0.05,0.05,10,10,10,10,10,10,10)$.
	
To study the performance of the a posteriori error estimates derived in \cref{sec:TRRB_a_post_error_estimates}, we neglect the outer-loop optimization and simply use a goal-oriented adaptive greedy algorithm~\cite{HDO2011} with Lagrangian basis extension (a) from \cref{sec:TRRB_construct_RB} to generate a ROM, which ensures that the worst relative estimated error for the reduced functional and its gradient over the adaptively generated training set and a randomly chosen validation set is below a prescribed tolerance of $\tau_{\textnormal{{FOC}}}= 5 \cdot 10^{-4}$. 
In particular we first ensure $\Delta_{\hat{J}_\red}(\mu)/\Jnoncor_\red(\mu) < \tau_{{FOC}}$ for $\Delta_{\hat{J}_\red}$ from \cref{prop:Jhat_error}(i) and continue with $\Delta_{\noncorgrad \Jnoncor_\red}(\mu)/\|\noncorgrad \Jnoncor_\red(\mu)\|_2 < \tau_{\noncorgrad \Jnoncor}$ for $\Delta_{\noncorgrad \Jnoncor_\red}$ from \cref{prop:grad_Jhat_error}(i), cf.~\cite[Algorithm 2]{QGVW2017}. 
Let us mention that the goal for $\Delta_{\hat{J}_\red}$ is fulfilled after $24$ basis enrichments.
We have $\Delta_{\noncorgrad \Jnoncor_\red}(\mu)/\|\noncorgrad \Jnoncor_\red(\mu)\| < 4.84$ after $56$ basis enrichments, where we artificially stop the algorithm since the associated computational effort is already roughly 17 hours, demonstrating the need for the proposed adaptive TR-RB algorithm studied in the next section.
	
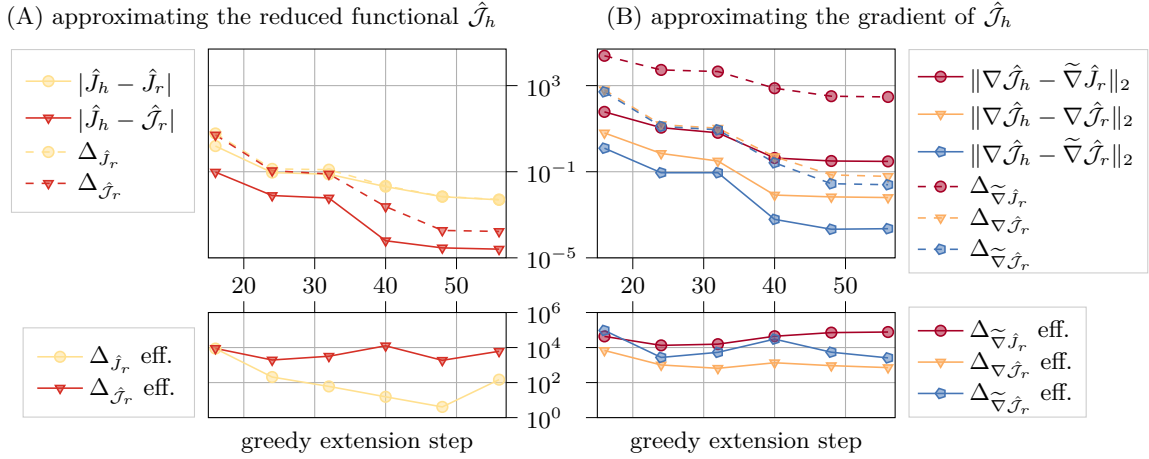
\begin{figure}[t]
	\centering%
	\footnotesize%
\begin{tikzpicture}

\definecolor{color0}{rgb}{0.65,0,0.15}
\definecolor{color1}{rgb}{0.84,0.19,0.15}
\definecolor{color2}{rgb}{0.96,0.43,0.26}
\definecolor{color3}{rgb}{0.99,0.68,0.38}
\definecolor{color4}{rgb}{1,0.88,0.56}
\definecolor{color5}{rgb}{0.67,0.85,0.91}
\definecolor{color6}{rgb}{0.27,0.46,0.71}
\definecolor{color7}{rgb}{0.19,0.21,0.58}

\begin{axis}[
  name=top_left,
width=5.5cm,
height=4.35cm,
legend cell align={left},
legend style={fill opacity=0.8, draw opacity=1, text opacity=1, at={(1.2,0)}, anchor=south, draw=white!80!black},
log basis y={10},
tick align=outside,
ytick pos=right,
xtick pos=bottom,
x grid style={white!69.0196078431373!black},
xmajorgrids,
xmin=15, xmax=57,
xtick style={color=black},
xticklabels={,,,},
y grid style={white!69.0196078431373!black},
ymajorgrids,
ymin=1e-05, ymax=50000,
ymode=log,
ytick style={color=black}
]
\addplot [semithick, color4, mark=*, mark size=2, mark options={solid, fill opacity=0.5}]
table {%
16 1.5602058
24 0.0913711
32 0.0766135
40 0.0203013
48 0.0069036
56 0.0049428
};
\addlegendentry{error J}
\addplot [semithick, color1, mark=triangle*, mark size=2, mark options={solid, rotate=180, fill opacity=0.5}]
table {%
16 0.0960233
24 0.0078391
32 0.0060158
40 0.0000622
48 0.0000290
56 0.0000253
};
\addlegendentry{error J-corr}
\addplot [semithick, color4, mark=*, mark size=2, dashed, mark options={solid, fill opacity=0.5}]
table {%
16 5.9913766
24 0.1355419
32 0.1245623
40 0.0218672
48 0.0070929
56 0.0050998
};
\addlegendentry{estimator J}
\addplot [semithick, color1, mark=triangle*, mark size=2, dashed, mark options={solid, rotate=180, fill opacity=0.5}]
table {%
16 5.1036908
24 0.1120268
32 0.0782738
40 0.0024091
48 0.0001877
56 0.0001684
};
\addlegendentry{estimator J corr}
\legend{};
\end{axis}

\begin{axis}[
  name=top_right,
  at=(top_left.east),
  anchor=west,
  xshift=1.2cm,
width=5.5cm,
height=4.35cm,
legend cell align={left},
legend style={fill opacity=0.8, draw opacity=1, text opacity=1, at={(1.2,0)}, anchor=south, draw=white!80!black},
log basis y={10},
tick align=outside,
x grid style={white!69.0196078431373!black},
xmajorgrids,
xmin=15, xmax=57,
xtick style={color=black},
xtick pos=bottom,
xticklabels={,,,},
y grid style={white!69.0196078431373!black},
ymajorgrids,
ymin=1e-05, ymax=50000,
ymode=log,
yticklabels={,,},
ytick pos=left,
ytick style={color=black}
]
\addplot [semithick, color0, mark=*, mark size=2, mark options={solid, fill opacity=0.5}]
table {%
16 60.3425755
24 11.1470314
32 6.4303011
40 0.4347226
48 0.3154153
56 0.3015742
};
\addlegendentry{error DJ}
\addplot [semithick, color3, mark=triangle*, mark size=2, mark options={solid, rotate=180, fill opacity=0.5}]
table {%
16 6.3335463
24 0.7063809
32 0.3138649
40 0.0082424
48 0.0067837
56 0.0063186
};
\addlegendentry{error DJ-corr}
\addplot [semithick, color6, mark=pentagon*, mark size=2, mark options={solid, rotate=90, fill opacity=0.5}]
table {%
16 1.2382556
24 0.0897214
32 0.0900121
40 0.0006126
48 0.0002131
56 0.0002262
};
\addlegendentry{error DJ-approx}
\addplot [semithick, color0, mark=*, mark size=2, dashed, mark options={solid, fill opacity=0.5}]
table {%
16 24291.8313696
24 5358.8100309
32 4503.6642106
40 756.5755578
48 321.9380287
56 297.4647676
};
\addlegendentry{estimator DJ}
\addplot [semithick, color3, mark=triangle*, mark size=2, dashed, mark options={solid, rotate=180, fill opacity=0.5}]
table {%
16 636.4101662
24 15.3648749
32 10.5050809
40 0.4986274
48 0.0712033
56 0.0608080
};
\addlegendentry{estimator DJ-corr}
\addplot [semithick, color6, mark=pentagon*, mark size=2, dashed, mark options={solid, rotate=90, fill opacity=0.5}]
table {%
16 513.9431314
24 12.7039995
32 8.8771813
40 0.2603803
48 0.0277693
56 0.0249135
};
\addlegendentry{estimator DJ-approx}
\legend{};
\end{axis}

\begin{axis}[
  name=bottom_left,
  at=(top_left.south east),
  anchor=north east,
  yshift=-0.725cm,
width=5.5cm,
height=2.975cm,
legend cell align={left},
legend style={fill opacity=0.8, draw opacity=1, text opacity=1, at={(1.2,0)}, anchor=south, draw=white!80!black},
log basis y={10},
tick align=outside,
xtick pos=top,
ytick pos=right,
x grid style={white!69.0196078431373!black},
xmajorgrids,
xmin=15, xmax=57,
xtick style={color=black},
y grid style={white!69.0196078431373!black},
ymajorgrids,
ymin=1, ymax=1e+6,
ymode=log,
ytick style={color=black}
]
\addplot [semithick, color4, mark=*, mark size=2, mark options={solid, fill opacity=0.5}]
table {%
16 8744.1808226
24 207.2319928
32 60.4071149
40 15.4345741
48 4.0803804
56 146.7465932
};
\addplot [semithick, color1, mark=triangle*, mark size=2, mark options={solid, rotate=180, fill opacity=0.5}]
table {%
16 8930.1670963
24 1948.8936487
32 3188.1462625
40 12586.6862083
48 1865.7813721
56 6113.7695524
};
\legend{};
\end{axis}

\begin{axis}[
  name=bottom_right,
  at=(top_right.south west),
  anchor=north west,
  yshift=-0.725cm,
width=5.5cm,
height=2.975cm,
legend cell align={left},
legend style={fill opacity=0.8, draw opacity=1, text opacity=1, at={(1.2,0)}, anchor=south, draw=white!80!black},
log basis y={10},
tick align=outside,
tick pos=left,
x grid style={white!69.0196078431373!black},
xmajorgrids,
xmin=15, xmax=57,
xtick pos=top,
xtick style={color=black},
y grid style={white!69.0196078431373!black},
ymajorgrids,
yticklabels={,,},
ytick pos=left,
ymin=1, ymax=1e+6,
ymode=log,
ytick style={color=black}
]
\addplot [semithick, color0, mark=*, mark size=2, mark options={solid, fill opacity=0.5}]
table {%
16 43287.5616865
24 13331.8687240
32 15603.3009341
40 44269.8307044
48 71392.3285677
56 77021.5333398
};
\addplot [semithick, color3, mark=triangle*, mark size=2, mark options={solid, rotate=180, fill opacity=0.5}]
table {%
16 6818.7796795
24 1013.1974612
32 653.1983662
40 1352.7644306
48 926.4852844
56 716.3903139
};
\addplot [semithick, color6, mark=pentagon*, mark size=2, mark options={solid, rotate=90, fill opacity=0.5}]
table {%
16 94919.3914497
24 2727.6434365
32 5326.0440363
40 30314.3579579
48 5422.5529477
56 2569.8465150
};
\legend{};
\end{axis}

\begin{customlegend}[legend cell align={left}, legend style={fill opacity=0.8, draw opacity=1, text opacity=1,
    at=(top_left.north west),
    anchor=north east,
    xshift=-5pt,
    inner sep=5pt,
draw=white!80!black},
	legend entries={
    $|\hat{J}_h - \Jnoncor_\red|$,
    $|\hat{J}_h - \cJhatn|$,
    $\Delta_{\Jnoncor_\red}$, 
    $\Delta_{\cJhatn}$, 
  }]
\addlegendimage{semithick, color4, mark=*, mark size=2, mark options={solid, fill opacity=0.5}}
\addlegendimage{semithick, color1, mark=triangle*, mark size=2, mark options={solid, rotate=180, fill opacity=0.5}}
\addlegendimage{semithick, color4, mark=*, mark size=2, dashed, mark options={solid, fill opacity=0.5}}
\addlegendimage{semithick, color1, mark=triangle*, mark size=2, dashed, mark options={solid, rotate=180, fill opacity=0.5}}
\end{customlegend}

\begin{customlegend}[legend cell align={left}, legend style={fill opacity=0.8, draw opacity=1, text opacity=1,
    at=(top_right.north east),
    anchor=north west,
    xshift=5pt,
    inner sep=3pt,
draw=white!80!black},
	legend entries={ 
  $\|\nabla \hat{\mathcal{J}}_h - \noncorgrad\Jnoncor_\red\|_2$,
  $\|\nabla \hat{\mathcal{J}}_h - \nabla\cJhatn\|_2$,
  $\|\nabla \hat{\mathcal{J}}_h - \noncorgrad\cJhatn\|_2$,
	$\Delta_{\noncorgrad\Jnoncor_\red}$,
  $\Delta_{\nabla\cJhatn}$, 
	$\Delta_{\noncorgrad\cJhatn}$,
	}]
\addlegendimage{semithick, color0, mark=*, mark size=2, mark options={solid, fill opacity=0.5}}
\addlegendimage{semithick, color3, mark=triangle*, mark size=2, mark options={solid, rotate=180, fill opacity=0.5}}
\addlegendimage{semithick, color6, mark=pentagon*, mark size=2, mark options={solid, rotate=90, fill opacity=0.5}}
\addlegendimage{semithick, color0, mark=*, mark size=2, dashed, mark options={solid, fill opacity=0.5}}
\addlegendimage{semithick, color3, mark=triangle*, mark size=2, dashed, mark options={solid, rotate=180, fill opacity=0.5}}
\addlegendimage{semithick, color6, mark=pentagon*, mark size=2, dashed, mark options={solid, rotate=90, fill opacity=0.5}}
\end{customlegend}

\begin{customlegend}[legend cell align={left}, legend style={fill opacity=0.8, draw opacity=1, text opacity=1,
    at=(bottom_left.south west),
    anchor=south east,
    xshift=-5pt,
    inner sep=5pt,
draw=white!80!black},
	legend entries={
    $\Delta_{\Jnoncor_\red}$ eff., 
    $\Delta_{\cJhatn}$ eff., 
  }]
\addlegendimage{semithick, color4, mark=*, mark size=2, mark options={solid, fill opacity=0.5}}
\addlegendimage{semithick, color1, mark=triangle*, mark size=2, mark options={solid, rotate=180, fill opacity=0.5}}
\end{customlegend}

\begin{customlegend}[legend cell align={left}, legend style={fill opacity=0.8, draw opacity=1, text opacity=1,
    at=(bottom_right.south east),
    anchor=south west,
    xshift=5pt,
    inner sep=3pt,
draw=white!80!black},
	legend entries={ 
	$\Delta_{\noncorgrad\Jnoncor_\red}$ eff.,
  $\Delta_{\nabla\cJhatn}$ eff., 
	$\Delta_{\noncorgrad\cJhatn}$ eff.,
	}]
\addlegendimage{semithick, color0, mark=*, mark size=2, mark options={solid, fill opacity=0.5}}
\addlegendimage{semithick, color3, mark=triangle*, mark size=2, mark options={solid, rotate=180, fill opacity=0.5}}
\addlegendimage{semithick, color6, mark=pentagon*, mark size=2, mark options={solid, rotate=90, fill opacity=0.5}}
\end{customlegend}

\node[anchor=north, yshift=-2pt] at (bottom_left.south) {greedy extension step};
\node[anchor=north, yshift=-2pt] at (bottom_right.south) {greedy extension step};
\node[anchor=south east, yshift=4pt] at (top_left.north east) {(A) approximating the reduced functional $\Jhat_h$};
\node[anchor=south west, yshift=4pt] at (top_right.north west) {(B) approximating the gradient of $\Jhat_h$};

\end{tikzpicture}
	\caption[Experiment 3: Evolution of the true and estimated model reduction error]{%
		Evolution of the true and estimated model reduction error (top) in the reduced functional and its approximations (A) and the gradient of the reduced functional and its approximations (B), as well as error estimator efficiencies (bottom), during adaptive greedy basis generation for the experiment from Section \ref{sec:TRRB_estimator_study}.
		Top: depicted is the $L^\infty(\Params_\textnormal{val})$-error for a validation set $\Params_\textnormal{val} \subset \Params$ of $100$ randomly selected parameters, i.e.~$|\hat{J}_h - \Jnoncor_\red|$ corresponds to $\max_{\mu \in \Params_\textnormal{val}} |\hat{J}_h(\mu) - \Jnoncor_\red(\mu)|$, $\Delta_{\cJhatn}$ corresponds to $\max_{\mu \in \Params_\textnormal{val}}\Delta_{\cJhatn}(\mu)$, $\|\nabla \hat{\mathcal{J}}_h - \nabla\cJhatn\|_2$ corresponds to $\max_{\mu \in \Params_\textnormal{val}}\|\nabla \hat{\mathcal{J}}_h(\mu) - \nabla\cJhatn(\mu)\|_2$, and so forth.
		Bottom: depicted is the worst effectivity of the respective error estimate (lower: better), i.e.~``$\Delta_{\hat{J}_\red}$ eff.'' corresponds to $\max_{\mu \in \Params_\textnormal{val}} \Delta_{\hat{J}_\red}(\mu) \,/\, |\hat{J}_h(\mu) - \Jnoncor_\red(\mu)|$, and so forth.
	}
	\label{fig:estimator_study_1}
\end{figure}

With the resulting RB spaces, we compare the ROM accuracy of the standard approach from \cref{sec:TRRB_standard_approach}, the NCD-corrected approach from \cref{sec:TRRB_ncd_approach}, as well as the sensitivity based approach from \cref{sec:TRRB_sensitivity_approach}.
Note that, for the sensitivity approach, more computational effort is needed since the FOM sensitivities have to be computed for each parameter dimension, thus increasing the effort to $2 + 2 \cdot 10 = 22$ FOM evaluations, compared to $2$ in the other approaches; cf. the discussion in \cref{sec:TRRB_construct_RB}. 
As we observe from Figure \ref{fig:estimator_study_1}, the error of the NCD-corrected terms is several orders of magnitude smaller than the corresponding terms of the standard approach.
It can also be seen that the (computationally more costly) gradient with approximated sensitivities, i.e.~$\noncorgrad\cJhatn$, from \cref{sec:TRRB_sensitivity_approach} shows the best error.
Certainly, all estimators for the NCD-correction and sensitivity-based quantities show a worse effectivity, hinting that there is still room for improvement in \cref{prop:grad_Jhat_error_sens}.

\subsubsection{Experiment 4: Optimization results for a complex problem}
\label{sec:mmexc_opt_results}

We now use the same parameterized problem as in \cref{sec:TRRB_estimator_study} to demonstrate the TR-RB algorithms on a more suitable model problem, where the NCD-correction term does not vanish, and parameter constraints need to be tackled, cf. \cref{rem:fin}.
Similar to \cref{sec:thermalfin}, starting with ten different randomly drawn initial parameters $\mu^{(0)}$, we measure the total computational run time, the number of TR iterations $k$ and the error in the optimal parameter for selected combinations of adaptive TR algorithms with ROM variants $1.$-$3.$ from  \cref{sec:TRRB_variants}, choice of RB spaces (a) and (b) from \cref{sec:TRRB_construct_RB}, as well as for the mentioned state-of-the-art-methods at the beginning of the section.

All algorithms converged (up to a tolerance) to the same point $\bar\mu$, and it was posteriorly verified that this point is a local minimum of $\Jhat$, i.e.~it satisfies the second-order sufficient optimality conditions. 
The value of $\bar\mu$ to compute the relative error was calculated with the FOM projected Newton method for a FOC condition tolerance of $10^{-12}$ and, thanks to the choice of the cost functional weights, the target $u^\textnormal{d}$ is approximate by $\bar u$ with a relative error of $1.7\cdot 10^{-6}$ in $D$.
We consider the same setup for two different stopping tolerances $\tau_{\textnormal{{FOC}}}= 5\cdot 10^{-4}$ and $\tau_{\textnormal{{FOC}}}= 10^{-6}$ to demonstrate that the performance (both in terms of time and convergence) of the methods vastly depends on the choice of $\tau_{\textnormal{{FOC}}}$.
\begin{table}
	\centering\footnotesize
	\begin{subtable}{\textwidth}
		\centering
		\begin{tabular}{|l|cc|c|cc|} \hline
			\textsc{\textbf{(A) $\tau_{\textnormal{{FOC}}}= 5\mkern-2mu \cdot\mkern-2mu 10^{-4}$} \kern-6pt}
			& Av.~(min/max) Run time[s]
			& \kern-4pt Speedup
			& \kern-2pt Av.~(min/max) it. \kern-8pt
			& Rel.~error
			& FOC cond.\\
			\hline
			FOM proj.~BFGS
			& 332.57 (196.51/591.85)
			& --
			& 44.30 (30/60)
			& $1.40\cdot 10^{-3}$
			& $1.80\cdot 10^{-4}$\\
			TR-RB from~\cite{QGVW2017}
			& 117.87 (70.29/166.31)
			& 2.82
			& 10.10 (6/14)
			& $5.46\cdot 10^{-4}$
			& $1.41\cdot 10^{-4}$\\
			1(a) TR-RB Lag.
			& 91.50 (47.07/230.29)
			& 3.63
			& 8.30 (5/10)
			& $2.01\cdot 10^{-3}$
			& $2.04\cdot 10^{-4}$\\
			1(b) TR-RB Agg.
			& 78.65 (54.69/114.36)
			& 4.23
			& 6.90 (5/9)
			& $2.53\cdot 10^{-4}$
			& $8.23\cdot 10^{-5}$ \\
			2(a) TR-RB Lag.
			& 79.47 (63.38/94.28)
			& 4.18
			& 8.50 (7/10)
			& $5.98\cdot 10^{-5}$
			& $1.02\cdot 10^{-4}$ \\
			3(a) TR-RB Lag.
			& 71.84 (50.38/87.16)
			& 4.63
			& 7.40 (5/9)
			& $1.09\cdot 10^{-3}$
			& $6.12\cdot 10^{-5}$ \\ \hline
		\end{tabular}
		\vspace{0.15cm}
	\end{subtable}
	\vspace{0.15cm}
	\begin{subtable}{\textwidth}
		\centering
		\begin{tabular}{|l|cc|c|cc|} \hline
			\textsc{\textbf{(B) $\tau_{\textnormal{{FOC}}}= 10^{-6}$}\kern+8pt}
			& Av.~(min/max) Run time[s]
			& \kern-4pt Speedup
			& \kern-2pt Av.~(min/max) it. \kern-8pt
			& Rel.~error
			& FOC cond.\\
			\hline
			FOM proj.~BFGS
			& 409.28 (317.25/637.55)
			& --
			& 57.00 (49/71)
			& $2.82\cdot 10^{-6}$
			& $3.35\cdot 10^{-7}$\\
			TR-RB from~\cite{QGVW2017}
			& 614.81 (566.66/671.97)
			& 0.66
			& 40.00 (40/40)
			& $8.46\cdot 10^{-7}$
			& $8.44\cdot 10^{-8}$\\
			1(a) TR-RB Lag.
			& 165.48 (92.26/417.24)
			& 2.47
			& 15.30 (10/40)
			& $3.29\cdot 10^{-6}$
			& $5.43\cdot 10^{-7}$\\
			1(b) TR-RB Agg.
			& 86.39 (62.68/124.43)
			& 4.74
			& 7.80 (6/10)
			& $3.52\cdot 10^{-6}$
			& $3.03\cdot 10^{-7}$ \\
			2(a) TR-RB Lag. 
			& 90.37 (80.97/102.60)
			& 4.53
			& 9.80 (9/11)
			& $8.12\cdot 10^{-7}$
			& $2.26\cdot 10^{-7}$ \\
			3(a) TR-RB Lag.
			& 88.24 (58.18/108.90)
			& 4.64
			& 8.90 (6/10)
			& $2.65\cdot 10^{-6}$
			& $2.73\cdot 10^{-7}$ \\ \hline
		\end{tabular}
		\vspace{0.15cm}
	\end{subtable}
	\caption[Experiment 4: Performance and accuracy of selected BFGS algorithms]{%
		Performance and accuracy of selected algorithms for two choices of $\tau_\textnormal{{FOC}}$ for the example from \cref{sec:mmexc_opt_results} for ten optimization runs with random initial guess, compare \cref{table:fin}.  }
	\label{table:MP2_com}
\end{table}

From \cref{table:MP2_com}, we observe that all proposed TR-RB methods speed up the FOM projected BFGS method with the NCD-corrected approach outperforming the others since the gradient used is the true one of the model function $\Jhat_\red$.
Moreover, independently of the model function, the algorithm from~\cite{QGVW2017} is much slower, demonstrating the positive impact of the suggested improvements on enlarging the TR radius and on the termination criterion based on cheaply available FOM information (instead of relying on a posteriori estimation), also visible in the number of outer TR iterations.
Comparing the TR-RB Variants $1$ and $2$ in terms of iterations, it is more beneficial to consider an aggregated RB space (b), i.e.~$V^\pr_\red = V^\du_\red$.
The enrichment for the aggregated spaces (b) is more costly and the time-to-ROM-solution is slightly larger. On the other hand, the more decadent space seems to allow better approximations of $\Jhat_h$, which can only be achieved by the NCD-corrected variant $3$(a).

All methods approximate the optimal parameter $\bar\mu$ with a small relative error and, at least for  $\tau_\textnormal{{FOC}}= 5\cdot10^{-4}$, reach the desired tolerance for the FOC condition.
However, given the resulting relative error in \cref{table:MP2_com} and \cref{Fig:_FOC_cond}, we observe that the choice $\tau_\textnormal{{FOC}}= 5\cdot10^{-4}$ is not sufficiently small for this model problem.
The methods do not reach an adequately low relative error in approximating $\bar\mu$ which moreover affects the timings by stopping the algorithm too early.
We conclude that the choice $\tau_{\textnormal{{FOC}}}= 10^{-6}$ instead results in a valid optimum of all variants (up to a tolerance of $10^{-6}$).
Note that the issue of choosing a good stopping tolerance can be resolved by utilizing the post-processing as explained in \cref{sec:outer_stopping}, which is demonstrated in \cref{num_test:12params}.
Importantly, for this choice of $\tau_\textnormal{{FOC}}$, we point out that the variant from~\cite{QGVW2017} only stopped because we restricted the maximum number of iterations to $40$. However, the FOM version of the FOC condition (as also used in our variants) dropped under the depicted tolerance of $10^{-6}$.
The reason why the variant from~\cite{QGVW2017} still did not stop is caused by the fact that in~\cite{QGVW2017} the a posteriori estimate, which is summed to the FOM FOC condition, can not become numerically small enough, showing the limit of the proposed stopping criterion in~\cite{QGVW2017}.
\begin{figure}[h]
	\centering\footnotesize
	\input{Pictures_paper1/mu_error_EXC_10_.tex} 	
	\caption[Experiment 4: Error decay and performance of selected algorithms for two choices of $\tau_{\textnormal{{FOC}}}$]{Error decay and performance of selected algorithms for two choices of $\tau_{\textnormal{{FOC}}}$ for the example from Section~\ref{sec:mmexc_opt_results} for a single optimization run with random initial guess, cf.~\cref{fig:fin_timing_mu_d}.
	}
	\label{Fig:_FOC_cond}
\end{figure}
From \cref{Fig:_FOC_cond}(B), we conclude that the NCD-corrected approaches 2(a) and 3(a) outperform the standard ROM variant 1(a), which also
reached the maximum number of iterations for one of the ten samples. 
Consequently, the NCD-correction entirely resolves the issue of the variational crime (introduced by splitting the reduced spaces) since it shows roughly the same performance as
variant 1(b).
Moreover, looking at the minimum and maximum number of computational time in \cref{table:MP2_com}, Variant 3(a) shows a less volatile and more robust
behavior.

In conclusion, in this subsection (\cref{sec:exp_paper1}), we obtained that the NCD-corrected variant with Lagrangian enrichment and BFGS as sub-problem solver $3.$(a) outperforms the other ROM strategies in terms of the best compromise between accuracy and computational efficiency.
Moreover, we demonstrated that our TR-RB variants show a significant improvement of the originally proposed method in~\cite{QGVW2017}.
A drawback that was found in the experiments is the sometimes misleading choice of a too-large $\tau_{\textnormal{FOC}}$ and, related to that, we expect that the ability of the BFGS method to converge for tiny tolerances $\tau_{\textnormal{FOC}}$ can suffer.

While the subsequent subsection demonstrates the parameter control for the choice of $\tau_{\textnormal{FOC}}$, we also show that the BFGS method, indeed, can have problems converging for smaller tolerances and that the Newton method, instead, resolves this issue.

\subsection{Newton approaches with optional enrichment and parameter control}
\label{sec:exp_paper2}

This section is a revised version of the numerical experiments in~\cite{paper2}.
We first demonstrate the TR-RB algorithm's post-processing, second compare the BFGS sub-problem solvers with Newton, and third, we devise the optional enrichment strategies.

To introduce the used abbreviations, we list the selected algorithms; cf. \cref{sec:TRRB_variants}:

\begin{description}	
\item[FOM TR-Newton-CG~\cite{Nocedal}:] Following~\cite[Algorithm~7.2]{Nocedal}, we use a FOM TR method that we presented in \cref{sec:background_TR_methods}. This method considers a standard FOM quadratic approximation for $\Jhat_h$ as the model function. It includes the computation of the Cauchy point as well as a way to handle the box constraints enforced on $\Params$, following~\cite[Section~16.7]{Nocedal}.
\item[BFGS NCD TR-RB (UE) {[Alg.~\ref{alg:TR-RBmethod_paper1}]}:] We use the NCD-corrected method with BFGS sub-problem solver and Lagrange RBs $3$(a) following \cref{alg:TR-RBmethod_paper1} with unconditional enrichment (UE) of the RB spaces, where no reduced Hessian or sensitivities of the primal and dual solutions are required.
\item[Newton NCD TR-RB (UE):] We use the NCD-corrected method with directional Taylor RB spaces $3$(c), a projected Newton method for the TR sub-problems, and unconditional enrichment (UE) of the RB spaces.
\item[Newton NCD TR-RB (OE) {[Alg.~\ref{alg:TR-RBmethod_paper2}]}:]We use the NCD-corrected method with directional Taylor RB spaces $3$(c), a projected Newton method for the TR sub-problems and optional enrichment (OE) of the RB spaces from \cref{alg:TR-RBmethod_paper2}.
\end{description}

\noindent We consider two case studies:
\begin{itemize}
	\item \textbf{Experiment 5:} Optimize 12 Parameters (3 walls, 2 doors, 7 heaters) with $\Paramsad= [0.025,0.1]^5\times[0,100]^7$ to reach the target $u^\textnormal{d}(x) = 18\rchi_{D}(x)$ and $\mu^\textnormal{d}\equiv0$ using also the a posteriori error estimate for the optimal parameter (Proposition~\ref{apo-est-parameters}) as post-processing; cf. \cref{num_test:12params}.
	\item \textbf{Experiment 6:} Optimize 28 Parameters (8 walls, 8 doors, 12 heaters) with target $u^\textnormal{d}= \mathcal{S}_h(\mu^\textnormal{d})$, where $\mu^\textnormal{d}\in\Params$ is given, and $\Paramsad= [0.025,0.1]^{16}\times[0,100]^{12}$; cf. \cref{num_test:28params}. 
\end{itemize}
Just as in the former numerical experiments, the methods are performed with ten different random samples for the starting parameter $\mu^{(0)}$.
For the sake of brevity, we omit additional details on the data functions and again refer to the accompanying code.

\subsubsection{Experiment 5: Parameter control}
\label{num_test:12params}
\begin{figure}
	\footnotesize
	\begin{tikzpicture}
		\definecolor{color0}{rgb}{0.65,0,0.15}
		\definecolor{color1}{rgb}{0.84,0.19,0.15}
		\definecolor{color2}{rgb}{0.96,0.43,0.26}
		\definecolor{color3}{rgb}{0.99,0.68,0.38}
		\definecolor{color4}{rgb}{1,0.88,0.56}
		\definecolor{color5}{rgb}{0.67,0.85,0.91}
		\begin{axis}[
			name=left,
			anchor=west,
			width=10cm,
			scale=0.75,
			log basis y={10},
			tick align=outside,
			tick pos=left,
			x grid style={white!69.0196078431373!black},
			xlabel={time in seconds [s]},
			xmajorgrids,
			xmin=0, xmax=265,
			xtick style={color=black},
			y grid style={white!69.0196078431373!black},
			ymajorgrids,
			ymin=2e-07, ymax=200,
			ymode=log,
			ylabel={\(\displaystyle \| \overline{\mu}_h-\mu^{(k)} \|_2\)},
			yticklabels={,,},
			ytick pos=right,
			yticklabel pos=left,
			ytick style={color=black}
			]
			\addplot [semithick, color0, mark=triangle*, mark size=3, mark options={solid, fill opacity=0.5}]
			table {%
				10.7283747196198 71.4005025909481
				16.8894658088684 70.0913980253683
				23.6503021717072 68.6311266604651
				30.5194776058197 68.6184250989761
				38.0436983108521 65.0536209016139
				46.4489510059357 32.1153217100059
				55.2568151950836 27.0201266807269
				64.8881371021271 15.1115679924662
				80.2292532920837 2.69350898011492
				107.479714155197 0.0120659897452426
			};
			\addplot [semithick, color1, mark=*, mark size=3, mark options={solid, fill opacity=0.5}]
			table {%
				10.7799642086029 71.4005025909481
				22.0994827747345 65.1324982951888
				34.5064990520477 53.9934432453033
				48.954788684845 31.6428672462588
				63.0492160320282 0.22019511096089
			};
			\addplot [semithick, color2, mark=square*, mark size=3, mark options={solid, fill opacity=0.5}]
			table {%
				11.0857253074646 71.4005025909481
				22.2940537929535 65.1324982951888
				35.2455587387085 51.7977039852713
				48.9639213085175 23.3668082033239
				53.176020860672 0.226218016009454
			};
			\addplot [thick, dotted, color0, mark=triangle, mark size=3, mark options={solid, semithick}]
			table {%
				10.7593002319336 71.4005025909481
				16.7505412101746 70.0913980253683
				23.3089096546173 68.6311266604651
				29.9739465713501 68.6184250989761
				37.3876357078552 65.0536209016139
				45.9340124130249 32.1153217100059
				54.7701907157898 27.0201266807269
				64.6972739696503 15.1115679924662
				79.9411873817444 2.69350898011492
				106.865009307861 0.0120659897452426
			};
			\addplot [thick, dashed, black, mark=square, mark size=0, mark options={solid, semithick}]
			table {%
				106.865009307861 0.0120659897452426
				158.72991752624478 0.0120659897452426
			};
			\addplot [thick, dotted, color0, mark=triangle, mark size=3, mark options={solid, semithick}]
			table {%
				158.72991752624478 0.0120659897452426
				174.107760667801 0.000461260770184682
				184.492898464203 0.000456269534405666
				195.489494085312 0.000453628915537476
				207.455162286758 0.000451930571291389
				220.113807439804 0.000438574164874584
				234.120990991592 0.000437426210881451
				251.963258266449 0.000211268220083034
				266.428707361221 0.000211268166754213
			};
			\addplot [thick, dotted, color1, mark=o, mark size=3, mark options={solid, semithick}]
			table {%
				10.7847511768341 71.4005025909481
				22.4626441001892 65.1324982951888
				34.9769263267517 53.9934432453033
				49.4036056995392 31.6428672462588
				63.5665023326874 0.22019511096089
			};
			\addplot [thick, dashed,, black, mark=o, mark size=0, mark options={solid, semithick}]
			table {%
				63.5665023326874 0.22019511096089
				117.5665023326874 0.22019511096089
			};
			\addplot [thick, dotted, color0, mark=o, mark size=3, mark options={solid, semithick}]
			table {%
				117.5665023326874 0.22019511096089
				129.270685434341 0.00317499944918528
			};
			\addplot [thick, dashed,, black, mark=o, mark size=0, mark options={solid, semithick}]
			table {%
				129.270685434341 0.00317499944918528
				178.270685434341 0.00317499944918528
			};
			\addplot [thick, dotted, color0, mark=o, mark size=3, mark options={solid, semithick}]
			table {%
				178.270685434341 0.00317499944918528
				195.078726291656 2.1194544828339e-05
			};
			\addplot [thick, dashed,, black, mark=o, mark size=0, mark options={solid, semithick}]
			table {%
				195.078726291656 2.1194544828339e-05
				246.078726291656 2.1194544828339e-05
			};
			\addplot [thick, dotted, color0, mark=o, mark size=3, mark options={solid, semithick}]
			table {%
				246.078726291656 2.1194544828339e-05
			};
			\addplot [thick, dotted, color2, mark=square, mark size=3, mark options={solid, semithick}]
			table {%
				10.7282385826111 71.4005025909481
				22.0372779369354 65.1324982951888
				34.8145577907562 51.7977039852713
				48.7089738845825 23.3668082033239
				53.176020860672 0.226218016009454 
			};
			\addplot [thick, dashed,, black, mark=square, mark size=0, mark options={solid, semithick}]
			table {%
				53.176020860672 0.226218016009454
				105.33605003356934 0.226218016009454
			};
			\addplot [thick, dotted, color2, mark=square, mark size=3, mark options={solid, semithick}]
			table {%
				105.33605003356934 0.226218016009454 
				116.6974101066592 0.00138082875085338 
			};
			\addplot [thick, dashed,, black, mark=square, mark size=0, mark options={solid, semithick}]
			table {%
				116.6974101066592 0.00138082875085338
				169.05726623535185 0.00138082875085338
			};
			\addplot [thick, dotted, color2, mark=square, mark size=3, mark options={solid, semithick}]
			table {%
				169.05726623535185 0.00138082875085338 
				179.08393168449405 8.66810548955347e-07 
			};
			\addplot [thick, dashed, black, mark=square, mark size=0, mark options={solid, semithick}]
			table {%
				179.08393168449405 8.66810548955347e-07
				231.2080514431 8.66810548955347e-07
			};
			\addplot [thick, dotted, color2, mark=square, mark size=3, mark options={solid, semithick}]
			table {%
				231.2080514431 8.66810548955347e-07
			};
		\end{axis}
		\begin{axis}[
			name=right,
			anchor=west,
			width=10cm,
			at=(left.east),
			xshift=1.3cm,
			scale=0.75,
			legend cell align={left},
			legend style={font=\tiny, fill opacity=0.8, draw opacity=1, text opacity=1, at=(left.north east), anchor=north east, xshift=-7.21cm, yshift=-1.44cm, draw=white!80!black},
			log basis y={10},
			tick align=outside,
			tick pos=left,
			x grid style={white!69.0196078431373!black},
			xlabel={time in seconds [s]},
			xmajorgrids,
			xmin=0, xmax=218,
			xtick style={color=black},
			y grid style={white!69.0196078431373!black},
			ymajorgrids,
			ymin=2e-07, ymax=200,
			ymode=log,
			ytick pos=left,
			ytick style={color=black}
			]
			\addplot [semithick, color0, mark=triangle*, mark size=3, mark options={solid, fill opacity=0.5}]
			table {%
				10.9332320690155 71.4005025909481
				17.0462052822113 70.0913980253683
				23.8179948329926 68.6311266604651
				30.5724768638611 68.6184250989761
				38.0406160354614 65.0536209016139
				46.4146301746368 32.1153217100059
				55.305762052536 27.0201266807269
				65.3911876678467 15.1115679924662
				81.0280475616455 2.69350898011492
				108.373761177063 0.0120659897452426
				126.480749845505 0.000461260770184682
				136.940483808517 0.000456269534405666
				147.813118219376 0.000453628915537476
				160.001117944717 0.000451930571291389
				172.585354804993 0.000438574164874584
				186.605408668518 0.000437426210881451
				204.703733205795 0.000211268220083034
				219.356258392334 0.000211268166754213
			};
			\addlegendentry{BFGS NCD TR-RB (UE)}
			\addplot [semithick, color1, mark=*, mark size=3, mark options={solid, fill opacity=0.5}]
			table {%
				10.9879050254822 71.4005025909481
				22.7747066020966 65.1324982951888
				35.4346406459808 53.9934432453033
				49.5354845523834 31.6428672462588
				63.1189568042755 0.22019511096089
				77.7369198799133 0.00317499944918528
				93.7262444496155 2.1194544828339e-05
			};
			\addlegendentry{Newton NCD TR-RB (UE)}
			\addplot [semithick, color2, mark=square*, mark size=3, mark options={solid, fill opacity=0.5}]
			table {%
				11.1055860519409 71.4005025909481
				22.735223531723 65.1324982951888
				35.5743370056152 51.7977039852713
				49.1154987812042 23.3668082033239
				62.9940128326416 0.226218016009454
				77.7316734790802 0.00138082875085338
				81.8664269447327 8.66810548955347e-07
			};
			\addlegendentry{Newton NCD TR-RB (OE)}
		\end{axis}
		\node[anchor=south, yshift=4pt] at (left.north) {(A) $\tau_{\textnormal{{FOC}}}= 5\times 10^{-4}$ + parameter control};
		\node[anchor=south, yshift=4pt] at (right.north) {(B) $\tau_{\textnormal{{FOC}}}= 1\times 10^{-7}$};
	\end{tikzpicture}
	\captionsetup{width=\textwidth}
	\caption[Experiment 5: Error decay and performance of selected algorithms]{\footnotesize{%
			Error decay and performance of selected algorithms for experiment 5 with unconditional enrichment (UE) vs.~optional enrichment (OE) for a single optimization run with random initial guess $\mu^{(0)}$ for two choices of $\tau_\textnormal{FOC}$ (solid lines) with optional intermediate parameter control according to
			\eqref{apo-est-parameters} (dotted lines): for each algorithm, each marker corresponds to one (outer) iteration of the optimization method and indicates the absolute error in the current parameter, measured against the computed FOM optimum. The dashed black horizontal lines indicate the time taken for the post-processing of the parameter control.
	}}
	\label{Fig:EXC12:mu_error}
\end{figure}
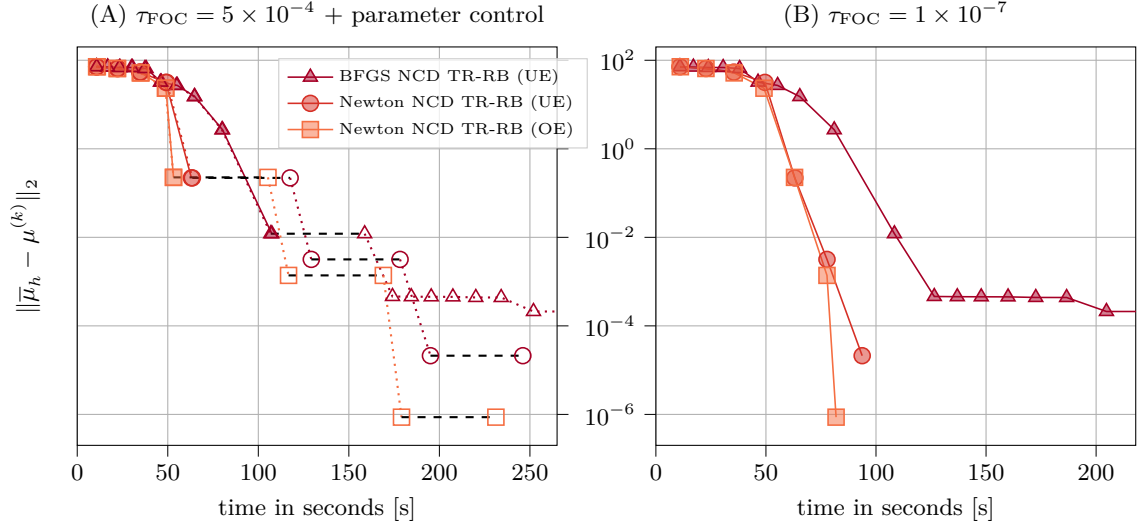

This experiment shows the usability of the a posteriori error estimate~\eqref{apo-est-parameters} and demonstrates the limitations of the projected BFGS method as the sub-problem solver.
We focus, at first, on the behavior of the methods for a given starting parameter $\mu^{(0)}$.
In \cref{Fig:EXC12:mu_error}, the error at each iteration $k$ is reported for the selected TR-RB algorithms, where we omit the FOM TR-Newton-CG due to its comparably large computational time.
We compute the solution with a tolerance $\tau_{\textnormal{{FOC}}}= 5\times 10^{-4}$ in \cref{Fig:EXC12:mu_error}(A).
When the methods reach this tolerance, we evaluate the a posteriori error estimate, and if this is greater than the value $\tau_{\mu}= 10^{-4}$, we decrease the tolerance $\tau_{\textnormal{{FOC}}}$ by two orders of magnitude and repeat the procedure until the a posteriori estimate is below the desired tolerance $\tau_{\mu}$.
In \cref{Fig:EXC12:mu_error}(B), we do not use the estimate and directly compute the solution with a tolerance $\tau_{\textnormal{{FOC}}}= 10^{-7}$ that was picked at last in \cref{Fig:EXC12:mu_error}(A).

We point out that in \cref{Fig:EXC12:mu_error}(A), the cost of computing the a posteriori error estimate~\eqref{apo-est-parameters} is included in the computational time (as a dashed horizontal line), which also includes the costly computation of the smallest eigenvalue of the FOM Hessian affecting the overall performances of the method.
When directly considering a smaller $\tau_{\textnormal{{FOC}}}$, \cref{alg:TR-RBmethod_paper2} is the fastest as visible from \cref{Fig:EXC12:mu_error}(B), which shows that the possibility of skipping enrichments improves the algorithm.

Another important issue that emerges from this numerical test is how the BFGS-based method struggles to reach smaller values of the FOC condition, resulting in a high increase of the computational time and stagnating error, as can be seen in \cref{Fig:EXC12:mu_error}. 
\begin{table}
	\footnotesize
	\centering
	\begin{tabular}{|l|cc|c|cc|} \hline
		& Run time[s]
		& & Iterations $k$ &&\\
		& Avg.~(min/max)
		& Speedup
		& Avg.~(min/max)
		& Rel.~error
		& FOC cond.\\\hline
		FOM   & 1381 (1190/1875) &  -- & 16.8 (14/23) & 2.66e-7 & 9.66e-8 \\
		BFGS (UE) & 818 (722/895) & 1.7 & 60 (60/60) & 2.77e-6 & 4.58e-7 \\
		Newton (UE) & 133 (94/212) & 10.4 & 7.8 (6/10) & 1.70e-7 & 1.73e-8 \\
		Newton (OE)  & 102 (82/141) & 13.6 &  6.9 (6/8) & 1.19e-7 & 1.30e-8 \\ \hline
	\end{tabular}
	\captionsetup{width=\textwidth}
	\caption[Experiment 5: Performance and accuracy of selected BFGS and Newton algorithms]{\footnotesize{%
			Performance and accuracy of the algorithms (abbreviated in order of definition) for Experiment 5 with unconditional enrichment (UE) vs.~optional enrichment (OE) for ten optimization runs with random initial guess $\mu^{(0)}$ and $\tau_{\textnormal{{FOC}}}=10^{-7}$: averaged, minimum and maximum total computational time (column 2) and speed-up compared to the FOM variant (column 3); average, minimum and maximum number of iterations $k$ required until convergence (column 4), the average relative error in the parameter (column 5) and average FOC condition (column 6).
			\label{tab:EXC12}
	}}
\end{table}

In \cref{tab:EXC12}, we report the average computational time and iterations for ten random starting
parameters $\mu^{(0)}$ together with the relative error in reconstructing the local minimizer and the FOC
condition at which the method stops.
Also, here, one can see how the possibility of skipping enrichments and the choice of a projected Newton method improve the results obtained with \cref{alg:TR-RBmethod_paper1}.
In particular, the projected BFGS method struggles to reach the given $\tau_{\textnormal{{FOC}}}$ in all experiments, showing its limitation in that case.
We remark that for larger tolerances, the method from \cref{alg:TR-RBmethod_paper1} is still valid and may converge faster, mainly depending on the given example.

\subsubsection{Experiment 6: Large parameter set}
\label{num_test:28params}
Next, We apply the TR-RB algorithm to a 28-dimensional parameter set.
As shown in the former experiments, the TR-RB algorithm overcomes the issue of a large offline phase.
What might still be problematic is the increase in the number of iterations, which would lead to large RB spaces with unconditional enrichment (UE).
The purpose of this experiment is to demonstrate that skipping enrichment yields similar convergence behavior (in terms of iterations), but at a lower cost.
\begin{table}[h]
	\footnotesize
	\centering
	\begin{tabular}{|l|cc|c|cc|} \hline
		& Run time[s]
		& & Iterations $k$ & & \\
		& Avg.~(min/max)
		& Speedup
		& Avg.~(min/max)
		& Rel.~error
		& FOC cond.\\\hline
		FOM   & 2423 (1962/3006) & -- & 18.5 (16/23) & 5.11e-9 & 4.57e-6 \\
		BFGS (UE) & 197 (156/272) & 12.3 & 12.1 (10/14) & 3.65e-9 & 2.38e-6 \\
		Newton (UE) & 258 (202/387) & 9.4 & 8.1 (7/9) & 5.83e-9 & 1.93e-6 \\
		Newton (OE)  & 168 (145/191) & 14.4 &  8.4 (7/12) & 1.22e-8 & 3.36e-6 \\ \hline
	\end{tabular}
	\captionsetup{width=\textwidth}
	\caption[Experiment 6: Performance and accuracy of the BFGS and Newton algorithms]{\footnotesize{%
			Performance and accuracy of the algorithms for experiment 6 for ten optimization runs with random initial guess $\mu^{(0)}$ and $\tau_{\textnormal{{FOC}}}=10^{-5}$, compare \cref{Tab:EXC12}.
			\label{Tab:EXC28}}}
\end{table}

\cref{Tab:EXC28} reports the average run time and iterations for the tested TR methods together with the relative error in reconstructing $\mu^\textnormal{d}$ and the final FOC value at the termination of the methods.
One can again note that all adaptive TR-RB algorithms speed up the computational time w.r.t.~the FOM TR-New.-CG.
Among all, the best performances are achieved by \cref{alg:TR-RBmethod_paper2}.
The number of outer iterations with or without (UE) is the same, while the computational time decreases.
This is due to two reasons: skipping an enrichment implies no preparation of the new RB space (like preassembling the new a posteriori estimator $\Delta_{\Jhat}$) and faster computations having a smaller RB space.
In \cref{Fig:EXC28} (left), one can see the error between the desired parameter $\mu^\textnormal{d}$ at each iteration of the different adaptive TR-RB methods for the same random starting parameter $\mu^{(0)}$, which confirms what is mentioned above.
Instead, \cref{Fig:EXC28} (right) shows the number of iterations needed to solve each TR sub-problem at the outer iteration $k$ of the method.
One can deduce that the advantages of \cref{alg:TR-RBmethod_paper2}, which is based on Newton, over \cref{alg:TR-RBmethod_paper1}, which is based on BFGS, are not due to a different number of inner iterations, but have to be associated with the reduction of the dimension of the RB space.
From \cref{Fig:EXC28} (right), we also see that the BFGS method might lose its super-linear convergence according to the approximation of the Hessian carried out by the method, which might deteriorate for an increasing number of active components of the parameter $\mu^{(k)}$; see~\cite{kelley}.

\begin{figure}
	\footnotesize\centering
	\begin{tikzpicture}
	\definecolor{color0}{rgb}{0.65,0,0.15}
	\definecolor{color1}{rgb}{0.84,0.19,0.15}
	\definecolor{color2}{rgb}{0.96,0.43,0.26}
	\definecolor{color3}{rgb}{0.99,0.68,0.38}
	\definecolor{color4}{rgb}{1,0.88,0.56}
	\definecolor{color5}{rgb}{0.67,0.85,0.91}
	\begin{axis}[
	name=left,
	anchor=west,
	width=6.5cm,
	height=4.5cm,
	log basis y={10},
	tick align=outside,
	tick pos=left,
	x grid style={white!69.0196078431373!black},
	xlabel={time in seconds [s]},
	xmajorgrids,
	xtick style={color=black},
	y grid style={white!69.0196078431373!black},
	ymajorgrids,
	ymode=log,
	ylabel={\(\displaystyle \| \overline{\mu}-\mu^{(k)} \|_2\)},
	ytick style={color=black}
	]
	\addplot [semithick, color0, mark=triangle*, mark size=3, mark options={solid, fill opacity=0.5}]
	table {%
		11.5934023857117 135.692657995198
		19.2377111911774 90.8446588955084
		27.8750627040863 24.6442190251569
		37.0965969562531 9.75267373981067
		47.1908657550812 4.48305564438341
		58.3313579559326 2.84545730545265
		70.7895290851593 1.68732573675647
		83.6306612491608 1.45223566118595
		97.9136147499084 0.631367397506257
		113.420466423035 0.343342059768175
		128.735586166382 0.307260609200393
		144.771761417389 0.0897779774035193
		167.995193243027 8.72687899118555e-05
		190.073382377625 1.59054173003722e-06
		210.727612018585 7.22008805445801e-07
	};
	\addplot [semithick, color1, mark=*, mark size=3, mark options={solid, fill opacity=0.5}]
	table {%
		11.9374821186066 135.692657995198
		27.3678214550018 90.8498097955523
		44.4709420204163 27.8277653035942
		64.5053310394287 7.95877703598574
		85.2987916469574 3.25287711111348
		111.669808626175 2.16159970164513
		139.39608001709 1.36285213889077
		178.082042694092 9.39198337447456e-05
		210.767727851868 2.56065163956946e-06
	};
	\addplot [semithick, color2, mark=square*, mark size=3, mark options={solid, fill opacity=0.5}]
	table {%
		11.6395268440247 135.692657995198
		27.2452867031097 90.8498097955523
		33.3784108161926 21.2028480910588
		49.9773924350739 13.8959573913167
		70.7744832038879 2.51586455562634
		94.8881878852844 1.59200009037262
		124.319049119949 0.00336062732991893
		156.706687688828 0.000906237762678042
		169.492668151855 2.72660621684241e-06
	};
	\end{axis}
	\begin{axis}[
	name=right,
	anchor=west,
	at=(left.east),
	xshift=0.1cm,
	width=6.5cm,
	height=4.5cm,
	legend cell align={left},
	legend style={nodes={scale=0.7}, fill opacity=0.8, draw opacity=1, text opacity=1, at=(left.north east), anchor=north east, xshift=-2.0cm, yshift=1.0cm, draw=white!80!black},
	tick align=outside,
	tick pos=left,
	x grid style={white!69.0196078431373!black},
	xlabel={outer TR iteration $k$},
	ylabel={sub-problem iterations $L$},
	xmajorgrids,
	xtick style={color=black},
	y grid style={white!69.0196078431373!black},
	ymajorgrids,
	ytick pos=right,
	]
	\addplot [semithick, color0, mark=triangle*, mark size=3, mark options={solid, fill opacity=0.5}]
	table {%
		0 4
		1 4
		2 4
		3 3
		4 3
		5 4
		6 3
		7 6
		8 6
		9 2
		10 5
		11 32
		12 21
		13 7
	};
	\addlegendentry{BFGS NCD TR-RB (UE) [Alg.~\ref{alg:TR-RBmethod_paper1}]}
	\addplot [semithick, color1, mark=*, mark size=3, mark options={solid, fill opacity=0.5}]
	table {%
		0 4
		1 4
		2 4
		3 3
		4 5
		5 4
		6 11
		7 8
	};
	\addlegendentry{Newton NCD TR-RB (UE)}
	\addplot [semithick, color2, mark=square*, mark size=3, mark options={solid, fill opacity=0.5}]
	table {%
		0 4
		1 3
		2 2
		3 4
		4 5
		5 11
		6 11
	};
	\addlegendentry{Newton NCD TR-RB (OE) [Alg.~\ref{alg:TR-RBmethod_paper2}]}
	\end{axis}
	\end{tikzpicture}
	\captionsetup{width=\textwidth}
	\caption[Experiment 6: Error decay and number of inner iterations]{\footnotesize{%
			Error decay w.r.t. the desired parameter $\bar\mu= \mu^\textnormal{d}$ and performance (left) and number of sub-problem iterations in each TR iteration (right) of selected algorithms for experiment 6 with unconditional enrichment (UE) vs.~optional enrichment (OE) for a single optimization run with random initial guess $\mu^{(0)}$ for $\tau_\textnormal{FOC} = 10^{-5}$.}}
	\label{Fig:EXC28}
\end{figure}
To conclude this subsection, we have seen that the Newton sub-problem solver, the optional enrichment, and the parameter post-processing are relevant enhancements for the TR-RB algorithm.
Although the NCD corrected approach shows a robust behavior, we still intend to discuss the PG-based ROM, where the NCD correction term always vanishes in the reduced functional.
However, as we see in the following, this approach is subject to stability issues in the reduced system.

\subsection{Petrov--Galerkin approach}
\label{sec:exp_paper3}

This section analyzes the behavior of the proposed PG-variant of the TR-RB algorithm, introduced by Variant~4 in \cref{sec:TRRB_variants}.
We aim to compare the PG variant's computational time and accuracy to the standard BFGS-based TR-RB variant.
As usual, we state the considered algorithms that we compare in this section:

\begin{description}
\item[FOM projected BFGS:] Equivalently to \cref{sec:exp_paper1}, we consider the standard FOM-based projected BFGS method as high-fidelity reference method, using FOM evaluations for all required quantities.
\item[3(a) TR-RB with NCD:] We consider the standard Galerkin Variant 3, including the NCD-correction term  with BFGS as sub-problem solver and unconditional Lagrangian enrichment (a)(UE).
\item[4(a) TR-RB with PG:] We consider the Petrov--Galerkin Variant 4 with BFGS as sub-problem solver and unconditional Lagrangian enrichment (UE)(a).
\end{description}

As numerical experiment, we again choose the 12-dimensional experiment from \cref{num_test:12params} with a slight modification of the objective functional for the error study, cf. \cref{sec:error_behavior_pg}.
In terms of algorithm details, we only differ in the choice of the stopping tolerance $\tau_{\textnormal{{FOC}}} = 10^{-6}$.

\subsubsection{Experiment 7: Error study for the PG-variant}
\label{sec:error_behavior_pg}
\begin{figure}[ht]
	\centering%
	\footnotesize%
	\begin{tikzpicture}
		
		\definecolor{color0}{rgb}{0.65,0,0.15}
		\definecolor{color1}{rgb}{0.84,0.19,0.15}
		\definecolor{color2}{rgb}{0.96,0.43,0.26}
		\definecolor{color3}{rgb}{0.99,0.68,0.38}
		\definecolor{color4}{rgb}{1,0.88,0.56}
		\definecolor{color5}{rgb}{0.67,0.85,0.91}
		\definecolor{color6}{rgb}{0.27,0.46,0.71}
		\definecolor{color7}{rgb}{0.19,0.21,0.58}
		
		\begin{axis}[
			name=top_middle,
			width=7cm,
			height=5.5cm,
			xshift=2cm,
			log basis y={10},
			tick align=outside,
			ytick pos=right,
			xtick pos=bottom,
			x grid style={white!69.0196078431373!black},
			xmajorgrids,
			xmin=2, xmax=62,
			xtick style={color=black},
			y grid style={white!69.0196078431373!black},
			ymajorgrids,
			ymin=1e-05, ymax=1e03,
			ymode=log,
			ytick style={color=black}
			]
			\addplot [semithick, color6, mark=square*, mark size=2, mark options={solid, rotate=45, fill opacity=0.5}]
			table {%
				4    43.0423546 
				8     4.7622899 
				12     3.8108226 
				16     2.3004248 
				20     0.9208100 
				24     0.1874082 
				28     0.1138518 
				32     0.0583603 
				36     0.0300142 
				40     0.0269577 
				44     0.0338155 
				48     0.0187490 
				52     0.0172461 
				56     0.0165199 
				60     0.0077635 
			};
			\addplot [semithick, color0, mark=diamond*, mark size=2, mark options={solid, rotate=180, fill opacity=0.5}]
			table {%
				4   4.3567802 
				8   0.9740319 
				12   0.7101863 
				16   0.2501180 
				20   0.0288619 
				24   0.0110932 
				28   0.0053377 
				32   0.0025348 
				36   0.0010780 
				40   0.0002019 
				44   0.0001357 
				48   0.0000615 
				52   0.0000577 
				56   0.0000448 
				60   0.0000329 
			};
			\addplot [semithick, color3, mark=triangle*, mark size=2, mark options={solid, fill opacity=0.5}]
			table {%
				4  3.9085827 
				8  1.2450593 
				12  0.7329396 
				16  0.3251995 
				20  0.0247842 
				24  0.0110414 
				28  0.0072402 
				32  0.0094383 
				36  0.0029239 
				40  0.0019344 
				44  0.0011918 
				48  0.0000691 
				52  0.0000788 
				56  0.0000418 
				60  0.0000387 
			};
			\addplot [semithick, color6, mark=square*, mark size=2, mark options={solid, rotate=180, fill opacity=0.5}]
			table {%
				4  113.8903117 
				8   78.8674009 
				12   66.5765970 
				16   41.4825594 
				20    3.3217361 
				24    0.9909078 
				28    0.6019228 
				32    0.5208526 
				36    0.3999896 
				40    0.0916279 
				44    0.0960969 
				48    0.0418367 
				52    0.0430256 
				56    0.0365259 
				60    0.0224140 
			};
			\addplot [semithick, color0, mark=diamond*, mark size=2, mark options={solid, rotate=90, fill opacity=0.5}]
			table {%
				4  116.8902273 
				8   71.2875380 
				12   53.3466076 
				16   21.9862603 
				20    1.4742002 
				24    0.4527781 
				28    0.2880695 
				32    0.1328267 
				36    0.0677626 
				40    0.0183905 
				44    0.0082919 
				48    0.0036965 
				52    0.0053584 
				56    0.0018361 
				60    0.0018599 
			};
			\addplot [semithick, color3, mark=triangle*, mark size=2, mark options={solid, rotate=180, fill opacity=0.5}]
			table {%
				4 116.9413910 
				8  81.3719586 
				12  43.7234391 
				16  37.3951660 
				20   1.3646116 
				24   0.6264351 
				28   0.3770512 
				32   0.4323930 
				36   0.1604367 
				40   0.0830978 
				44   0.1297149 
				48   0.0083390 
				52   0.0044400 
				56   0.0021994 
				60   0.0033864 
			};
			
			\legend{};
		\end{axis}
		
		\begin{axis}[
			name=top_right,
			at=(top_middle.east),
			anchor=west,
			xshift=1.2cm,
			width=7cm,
			height=5.5cm,
			legend cell align={left},
			legend style={fill opacity=0.8, draw opacity=1, text opacity=1, at={(1.2,0)}, anchor=south, draw=white!80!black},
			log basis y={10},
			tick align=outside,
			x grid style={white!69.0196078431373!black},
			xmajorgrids,
			xmin=2, xmax=62,
			xtick style={color=black},
			xtick pos=bottom,
			y grid style={white!69.0196078431373!black},
			ymajorgrids,
			ymin=1e-05, ymax=1e03,
			ymode=log,
			yticklabels={,,},
			ytick pos=left,
			ytick style={color=black}
			]
			\addplot [semithick, color0, mark=*, mark size=2, mark options={solid, fill opacity=0.5}]
			table {%
				4   0.6979035 
				8   0.1683773 
				12   0.0880920 
				16   0.0661280 
				20   0.0349068 
				24   0.0068316 
				28   0.0045575 
				32   0.0036580 
				36   0.0026852 
				40   0.0012673 
				44   0.0009064 
				48   0.0005655 
				52   0.0003778 
				56   0.0003298 
				60   0.0002402 
			};
			\addplot [semithick, color3, mark=triangle*, mark size=2, mark options={solid, rotate=180, fill opacity=0.5}]
			table {%
				4 0.7062068 
				8 0.3404906 
				12 0.1114368 
				16 0.1809833 
				20 0.1392032 
				24 0.0474685 
				28 0.0073499 
				32 0.0077065 
				36 0.0051511 
				40 0.0101835 
				44 0.0090959 
				48 0.0016053 
				52 0.0006615 
				56 0.0005745 
				60 0.0006952 
			};
			\addplot [semithick, color0, mark=pentagon*, mark size=2, mark options={solid, rotate=90, fill opacity=0.5}]
			table {%
				4    40.7460763 
				8    38.8899377 
				12    33.2418317 
				16    23.2686215 
				20     9.4626144 
				24     4.1539872 
				28     3.4412535 
				32     2.8128896 
				36     2.3952862 
				40     1.2969305 
				44     1.1588881 
				48     0.9383411 
				52     0.7765199 
				56     0.6220427 
				60     0.4001464 
			};
			\addplot [semithick, color3, mark=star, mark size=2, mark options={solid, fill opacity=0.5}]
			table {%
				4   40.8146894 
				8   39.3312537 
				12   35.1527006 
				16   35.3298773 
				20   10.6842491 
				24   10.9024490 
				28    8.1679185 
				32   12.9165319 
				36   11.0609593 
				40    8.8701537 
				44   43.5757322 
				48    5.3983823 
				52    7.8105103 
				56    5.2246835 
				60   10.7762531 
			};
			\legend{};
		\end{axis}
		
		\node[anchor=north, yshift=-15pt] at (top_middle.south) {greedy extension step};
		\node[anchor=north, yshift=-15pt] at (top_right.south) {greedy extension step};
		\node[anchor=south, yshift=4pt, xshift=2.6cm] at (top_middle.north west) {(A) functional $\Jhat_h$ and gradient $\nabla \Jhat_h$};
		\node[anchor=south, yshift=4pt, xshift=2.4cm] at (top_right.north west) {(B) primal and dual solution};
		
	\end{tikzpicture}
	\begin{tikzpicture}
		
		\definecolor{color0}{rgb}{0.65,0,0.15}
		\definecolor{color1}{rgb}{0.84,0.19,0.15}
		\definecolor{color2}{rgb}{0.96,0.43,0.26}
		\definecolor{color3}{rgb}{0.99,0.68,0.38}
		\definecolor{color4}{rgb}{1,0.88,0.56}
		\definecolor{color5}{rgb}{0.67,0.85,0.91}
		\definecolor{color6}{rgb}{0.27,0.46,0.71}
		\definecolor{color7}{rgb}{0.19,0.21,0.58}
		
		\begin{customlegend}[legend cell align={left}, legend style={fill opacity=0.8, draw opacity=1, text opacity=1,
				at=(top_right.south),
				anchor=north,
				xshift=-5.5cm,
				yshift=-1.2cm,
				inner sep=5pt,
				/tikz/column 2/.style={
					column sep=4pt},
				/tikz/column 4/.style={
					column sep=4pt},
				/tikz/column 6/.style={
					column sep=4pt},
				/tikz/column 8/.style={
					column sep=4pt},
				draw=white!80!black}, 
			legend columns=4,
			legend entries={
				$|\Jhat_h - \Jnoncor_\red|$,
				$|\nabla\Jhat_h - \tilde{\nabla}\Jnoncor_\red|$,
				$|u_{h,\mu} - u_{\red,\mu}|$,
				$|p_{h,\mu} - p_{\red,\mu}|$, 
				$|\Jhat_h - \cJhatn|$,
				$|\nabla\Jhat_h - \nabla\cJhatn|$,
				$|u_{h,\mu} - u^\textnormal{pg}_{\red,\mu}|$,
				$|p_{h,\mu} - p^\textnormal{pg}_{\red,\mu}|$,
				$|\Jhat_h - \cJhatn^{\mkern-7mu\textnormal{pg}}|$,
				$|\nabla\Jhat_h - \nabla\cJhatn^{\mkern-7mu\textnormal{pg}}|$
			}]
			\addlegendimage{semithick, color6, mark=square*, mark size=2, mark options={solid, rotate=45, fill opacity=0.5}}
			\addlegendimage{semithick, color6, mark=square*, mark size=2, mark options={solid, rotate=180, fill opacity=0.5}}
			\addlegendimage{semithick, color0, mark=*, mark size=2, mark options={solid, fill opacity=0.5}}
			\addlegendimage{semithick, color0, mark=pentagon*, mark size=2, mark options={solid, rotate=90, fill opacity=0.5}}
			\addlegendimage{semithick, color0, mark=diamond*, mark size=2, mark options={solid, rotate=180, fill opacity=0.5}}
			\addlegendimage{semithick, color0, mark=diamond*, mark size=2, mark options={solid, rotate=90, fill opacity=0.5}}
			\addlegendimage{semithick, color3, mark=triangle*, mark size=2, mark options={solid, rotate=180, fill opacity=0.5}}
			\addlegendimage{semithick, color3, mark=star, mark size=2, mark options={solid, fill opacity=0.5}}
			\addlegendimage{semithick, color3, mark=triangle*, mark size=2, mark options={solid, fill opacity=0.5}}
			\addlegendimage{semithick, color3, mark=triangle*, mark size=2, mark options={solid, rotate=180, fill opacity=0.5}}	
		\end{customlegend}
	\end{tikzpicture}
	
	\caption[Experiment 7: Evolution of the true reduction error for PG]{%
		Evolution of the true reduction error during adaptive greedy basis generation. In (a), we visualize the reduced functional, its gradient, and its approximations, and in (b) the primal and dual solutions and their approximations.
		Depicted is the $L^\infty(\Params_\textnormal{val})$-error for a validation set $\Params_\textnormal{val} \subset \Params$ of $100$ randomly selected parameters, i.e.~$|\Jhat_h - \Jnoncor_\red|$ corresponds to $\max_{\mu \in \Params_\textnormal{val}} |\Jhat_h(\mu) - \Jnoncor_\red(\mu)|$, and so forth.
	}
	\label{fig:estimator_study}
\end{figure}
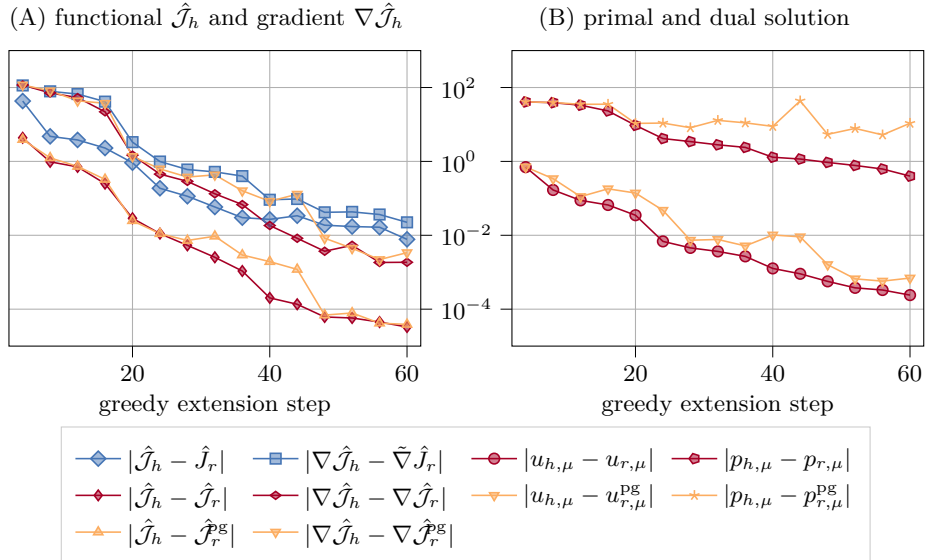

This section aims to show and discuss the MOR for the proposed Petrov--Galerkin approach.
Furthermore, we compare it to the Galerkin strategy from the NCD-corrected approach.
Recall that we denote the PG-solutions of~\eqref{eq:state_red_pg} and~\eqref{eq:dual_solution_red_pg} by $u_{\red}^\textnormal{pg}$, $p_{\red}^\textnormal{pg}$ and the corresponding functional by $\cJhatn^\textnormal{pg}$.
As before, the non-corrected approach from \cref{sec:TRRB_standard_approach} is denoted by $\hat{J}_\red$ and the NCD-corrected functional from \cref{sec:TRRB_ncd_approach} by $\cJhatn$.
To compare the accuracy of the three different functionals, we again employ a standard goal-oriented greedy-search algorithm as also done in \cref{sec:TRRB_estimator_study} with the relative a posteriori error of the objective functional $\Delta_{\cJhatn}(\mu)/\cJhatn(\mu)$.
As pointed out in \cref{sec:TRRB_RB_estimates_pg}, due to $V_h^\pr = V_h^\du$, we can replace the inf-sup constant by a lower bound for the coercivity constant of the conforming approach.
It is also important to mention that, for the following experiment, we have simplified our objective functional $\J$ by setting the domain of interest to the whole domain $D \equiv \Omega$.
As a result, the dual problem is more straightforward, enhancing the PG approach's stability.
The reason for that is further discussed below.
\cref{fig:estimator_study} shows the difference in the error decay and accuracy of the different approaches.
It can be seen that the NCD-corrected approach remains the most accurate approach. At the same time, the objective functional and the gradient of the PG approach show a better approximation than those of the non-corrected version.
The PG approximation of the primal and dual solutions is less accurate, and solely the primal error decays sufficiently stable.

During the shown error study, we experienced instabilities of the PG reduced primal and dual systems.
This is because, for very complex dual problems, the test space of each problem very poorly fits the respective ansatz space.
The decay of the primal error that can be seen in \cref{fig:estimator_study}(B) can not be expected in general and is indeed a consequence of the simplification of the objective functional, where we chose $D \equiv \Omega$.
The stability problems can already be seen in the dual error and partly on the primal error for basis size $20$.
In general, the reduced systems may be unstable for specific parameter values.
Instead, a sophisticated greedy algorithm for deducing appropriate reduced spaces may require constructing a larger dual or primal space by adding stabilizing snapshots, e.g. supremizers~\cite{ballarin2015supremizer}.
Since this section does not aim to provide an appropriate greedy-based algorithm, we instead decided to reduce the complexity of the functional, where stability issues are less present.

\subsubsection{Experiment 8: Optimization results for the PG-variant}

We now compare the PG-variant of the TR-RB with the above-mentioned NCD-corrected TR-RB approach.
Importantly, we use the original version of the model problem, i.e. we pick the domain of interest to be defined as suggested in \cref{ex1:blueprint}.
In fact, we accept high instabilities in the reduced model and show that the method still converges.
We pick ten random starting parameters, perform both algorithms and compare the averaged result.
In \cref{Fig:EXC12} one particular starting parameter is depicted and,
in \cref{Tab:EXC12}, the averaged results for all ten optimization runs are shown.
\begin{figure}
	\footnotesize
	\centering
	\begin{tikzpicture}
		\definecolor{color0}{rgb}{0.65,0,0.15}
		\definecolor{color1}{rgb}{0.84,0.19,0.15}
		\definecolor{color2}{rgb}{0.96,0.43,0.26}
		\definecolor{color3}{rgb}{0.99,0.68,0.38}
		\definecolor{color4}{rgb}{1,0.88,0.56}
		\definecolor{color5}{rgb}{0.67,0.85,0.91}
		\begin{axis}[
			name=left,
			width=8.5cm,
			height=5.5cm,
			log basis y={10},
			tick align=outside,
			tick pos=left,
			legend style={nodes={scale=0.7}, fill opacity=0.8, draw opacity=1, text opacity=1,
				xshift=-1.6cm, yshift=-0.5cm, xshift=4.5cm, draw=white!80!black},
			x grid style={white!69.0196078431373!black},
			xlabel={time in seconds [s]},
			xmajorgrids,
			xtick style={color=black},
			y grid style={white!69.0196078431373!black},
			ymajorgrids,
			ymode=log,
			ylabel={\(\displaystyle \| \overline{\mu}-\mu^{(k)} \|^\textnormal{rel}_2\)},
			ytick style={color=black}
			]
			\addplot [semithick, color0, mark=triangle*, mark size=3, mark options={solid, fill opacity=0.5}]
			table {%
				13.0756642818451 0.574286171259373
				31.8166162967682 0.574017027495371
				38.5897059440613 0.555477788898345
				45.570326089859 0.52429437114602
				52.9812984466553 0.507084209408289
				60.9106016159058 0.41660095554096
				69.3280558586121 0.375510754302739
				78.2450692653656 0.299488619035033
				88.5779669284821 0.113477915062374
				107.952816724777 0.0114236653909785
				135.231210231781 0.000147935221399132
				151.076913118362 2.58254983619181e-05
				164.675134181976 5.26070844357256e-06
			};
			\addlegendentry{TR-RB NCD BFGS}
			\addplot [semithick, color3, mark=*, mark size=3, mark options={solid, fill opacity=0.5}]
			table {%
				12.1607882976532 0.574286171259373
				35.4619424343109 0.574040233206862
				43.4497337341309 0.569126018670776
				52.2917582988739 0.56124503305629
				62.6727526187897 0.548418201394747
				74.8126885890961 0.525321267398303
				88.072719335556 0.524967245607288
				102.938640356064 0.519210519549374
				120.043780326843 0.518929447625113
				138.781098842621 0.504681744974713
				159.799144983292 0.478880041078375
				182.603423595428 0.428336507745106
				208.776380777359 0.405688185871396
				237.186156749725 0.30213730835205
				268.500371694565 0.236613269296937
				302.714802026749 0.202947643541194
				342.248579025269 0.0351113761140485
				384.453255176544 0.0291327949814106
				437.500003099442 3.37918525465338e-05
				486.676156520844 9.02321967828624e-07
			};
			\addlegendentry{TR-RB PG-BFGS }
		\end{axis}
	\end{tikzpicture}
	\captionsetup{width=\textwidth}
	\caption[Experiment 8: Relative error decay for PG]{\footnotesize{%
			Relative error decay w.r.t. the optimal parameter $\bar\mu$ and performance of selected algorithms for a single optimization run with random initial guess $\mu^{(0)}$ for $\tau_\textnormal{FOC} = 10^{-6}$.}}
	\label{Fig:EXC12}
\end{figure}
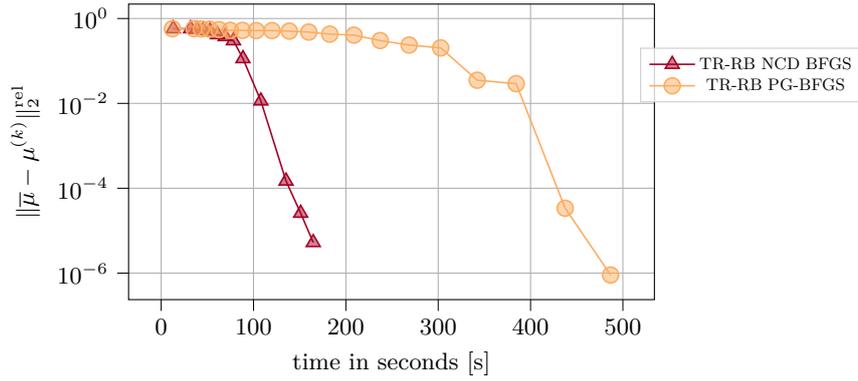
\begin{table}
	\footnotesize
	\centering
	\begin{tabular}{|l|cc|c|cc|} \hline
		& Run time[s]
		& & Iterations $k$ & & \\
		& Avg.~(min/max)
		& Speedup
		& Avg.~(min/max)
		& Rel.~error
		& FOC cond.\\\hline 
		FOM BFGS    & 6955~(4375/15556) & --&471.44~(349/799) & 3.98e-5 & 3.36e-6\\
		TR-RB NCD-BFGS & 171~(135/215)& 40.72&~12.56(11/15) & 4.56e-6& 6.05e-7\\
		TR-RB PG-BFGS  & 424~(183/609)& 16.39&~17.56(11/22) & 4.62e-6& 7.57e-7 \\ \hline
	\end{tabular}
	\captionsetup{width=\textwidth}
	\caption[Experiment 8: Performance and accuracy of PG algorithms]{\footnotesize{%
			Performance and accuracy of the algorithms for ten optimization runs with random initial guess $\mu^{(0)}$ and $\tau_{\textnormal{{FOC}}}=10^{-6}$.
			\label{Tab:EXC12}}}
\end{table}

It can be seen that the PG variant converges sufficiently fast compared to the FOM BFGS method. 
Certainly, it can not be said that the stability issues do not enter the performance of the proposed TR-RB methodology.
Still, regardless of the stability of the reduced system, we note that the convergence result in \cref{thm:convergence_of_TR} holds. 
However, the instability of the reduced equations harms the algorithm from iterating as fast as the NCD-corrected approach.
One reason is that the trust region is larger for the NCD-corrected approach, allowing the method to have faster optimization steps.
We also emphasize that the depicted result in \cref{Fig:EXC12} is neither an instance of the worst nor the best performance of the PG approach but rather an intermediate version.
The comparison highly depends on the starting parameter and the structure of the optimization problem.
As discussed above, the suggested PG approach can benefit from more involved enrichment strategies that account for the mentioned stability issues.

Last but not least, we recall that the experiment above showed weaknesses of the chosen projected BFGS approach as FOM method as well as for the TR-RB sub-problems, which has been extensively studied in \cref{num_test:12params}.
Instead, it is beneficial to choose higher-order optimization methods, such as projected Newton-type methods, as is done in \cref{sec:exp_paper2}, which we omit for the PG approach.

For a final conclusion of the experiments in this section, we refer to \cref{sec:TR_RB_conclusion}.
Beforehand, we present further approaches that are designed to enhance the algorithm with techniques that go beyond the state in \cite{paper2,paper1,paper3}.

\section{Further approaches}
\label{sec:TRRB_related_approaches}

The TR-RB algorithm presented in the former sections showed a robust convergence behavior and proved flexible, incorporating multiple additional features.
While the discussed add-ons like radius enlargement, optional and adapted enrichment, Newton sub-problem solver, and FOM-based stopping criterion proved to enhance or accelerate the algorithm, some ideas that evolved during the project have been left out in order to simplify the presentation of the basic TR-RB algorithm.
This section is devoted to a collection of further approaches that can be used to enhance the TR-RB algorithm concerning particular drawbacks.

Next, we discuss AFEM variants, coarsening of the algorithm's RB spaces, and enrichment of rejected iterates.
Subsequently, we discuss a relaxed version of the TR-RB algorithm, which also plays a vital role in \cref{chap:TR_TSRBLOD} and can further be used as a transition to non-certified algorithms.
Moreover, we show how some of these variants apply to the proof-of-concept experiment from \cref{sec:proof_of_concept}.

\subsection{AFEM variants} \label{sec:afem_approaches}
The presented TR-RB algorithm assumes that the underlying finite-element mesh has been chosen sufficiently small for capturing the data functions entirely.
Such an assumption is standard in the context of RB methods and has been formulated in \cref{asmpt:truth}.
The correct mesh size is indeed known in many applications, but the assumption is not fulfilled in general cases.

As discussed in \cref{sec:background_npdgl}, suitable variants like AFEM adaptively refine the mesh, e.g., by using a posteriori error estimation \cite{Ver1996,Ver2013}, see~\cite{CNSV2016a} for a recent article.
In a general application of MOR methods, it is not immediately apparent why AFEM methods should be used, especially when \cref{asmpt:truth} is enforced.
Certainly, for solving PDE-constrained optimization problems where the only aim is to find the critical point, MOR methods are solely used for accelerating the optimization process.
Thus, it can indeed be feasible to perform a preparatory optimization process on a coarser finite-element mesh and, embedded in the optimization process, verify whether to refine the mesh further w.r.t. certified measures.

As a generic choice of an AFEM variant of the TR-RB algorithm, we may start with a cheap FEM space as the first FOM and perform the TR-RB algorithm.
After convergence, we refine the mesh, e.g., with a fixed contraction factor on $h$, and use the computed critical point as the initial guess for restarting the whole procedure, unless we reach a mesh size where AFEM does not any more call for a refinement.
Many questions are to be answered for this approach, e.g., how information of former mesh sizes can be reused or how the convergence of the optimization routine can be verified.

We point out that adaptive mesh refinement strategies can also be used if $h$ is known according to \cref{asmpt:truth}.
Then, no measures for sufficient mesh refinement need to be enforced.
Instead, computations on coarser FEM mesh sizes can be considered a warm start of the actual TR-RB algorithm, meaning that (with low computational effort) a more appropriate initial point is found in advance.
Critically, such a procedure can still be problematic since the coarser FEM model may lead to an entirely misleading point.

We leave further considerations and details of AFEM variants of the TR-RB algorithm to future research.

\subsection{Coarsening the RB spaces}
\label{sec:coarsening_RB_approach}
The non-conforming choice of the reduction process, i.e.~separating the primal and dual reduced space, was mainly motivated by keeping the reduced spaces as small as possible since the offline time increases significantly with the size of the RB spaces.
Another reason is that the reduced matrices are dense and their (online) solution method is cubic in the respective RB dimension.
Suppose the TR-RB algorithm takes significantly many outer optimization steps.
In that case, the RB spaces also have approximately the same size (minus the cases where the optional enrichment suggested skipping the bases).
Especially for the case of multi-objective optimization problems, where the same TR-RB algorithm is enforced multiple times with accumulated spaces, it can be advantageous to find a notion of how to coarsen the RB spaces by removing some (unused) information from the reduced model.

To motivate that such an approach can be feasible, we again refer to the visualization of the method in \cref{sec:proof_of_concept}; see also \cref{fig:NO_TR_for_poc}.
The initial model is built around the initial guess~$\mu^{(0)}$, which is randomly chosen.
Thus, it can only happen by coincidence that the primal and dual RB basis functions for $\mu^{(0)}$ are related to the parameter region of the optimal value, and the first surrogate models in the algorithm are mainly used for optimization steps towards this region.
It can be expected that every new iterate is (in a way) closer to the optimum, suggesting that the information that has been used to reach the optimal region is not significantly relevant anymore.


Many ideas for removing basis functions can be considered.
For instance, the most unrelated basis functions could be found by the orthogonal projection of the current iterate.
After the algorithm excepts a new iterate, such a coarsening can be performed before or after the enrichment.
Importantly, however, the coarsening strategies should be enforced without excessive computational effort, diminishing the advantage of less dimensional reduced spaces.
We do not detail such approaches and instead refer to \cite{banholzer2022trust}, where first results in this direction have been made in the context of multi-objective optimization problems.

\subsection{FOM-cost-oriented TR-RB algorithm}
\label{sec:FOM_cost_or_TRRB}

The TR-RB algorithm was developed in the spirit of basic TR methods.
This means that, for each TR-RB sub-problem, we do not change the model function $m^{(k)}$ if the iterate is rejected.
For standard TR methods, this strategy is meaningful since updating $m^{(k)}$ in a reasonable way usually can not be done by solely evaluating the high-fidelity functional $\J_h$; cf. \cref{sec:background_TR_methods}.

Note that, in the TR-RB algorithm, the condition for rejecting the iterate \eqref{eq:suff_decrease_condition} can cheaply be verified or contradicted by conditions \eqref{eq:suf_cond} and \eqref{eq:nec_cond}.
On the other hand, if both conditions do not give an immediate decision, \eqref{eq:suff_decrease_condition} needs to be checked itself, which includes evaluating the FOM model at the current iterate.
If the algorithm rejects the iterate because \eqref{eq:suff_decrease_condition} is not fulfilled, we also ignore the FOM information that we needed to check \eqref{eq:suff_decrease_condition}.
Therefore, FOM evaluations are discarded, which may lead to a slower convergence of the overall algorithm.

Instead, we suggest a FOM-cost-oriented TR-RB algorithm where computed FOM snapshots are unconditionally used for updating the model function $m^{(k)}$, regardless of the fact, whether the same sub-problem is solved again (with a shrunk radius).
Importantly, such a strategy does not change the convergence of the TR-RB algorithm.
We omit a formal description of the algorithm but introduce the abbreviation of the FOM-cost-oriented TR-RB algorithm as FCO-TR-RB.
The algorithm is utilized in the numerical experiments below; see \cref{sec:experiments_further_approaches}.

One may be concerned that the FCO-TR-RB contradicts the idea of coarsening the RB spaces since more FOM snapshots are added to the basis than before.
However, a newly added basis function based a rejected parameter can rather be considered a more suitable basis function than the ones before.
The numerical experiment in \cref{sec:experiments_further_approaches} proves this expectation.

The FCO-TR-RB algorithm is mainly concerned with taking maximum advantage of rejections.
In what follows, we follow a completely different strategy and present a variant where early rejections are instead avoided entirely.

\subsection{Relaxed TR-RB algorithm}
\label{sec:Relaxed_TRRB}

We now introduce an adaptive algorithm that can be used as a warm start of the TR-RB algorithm and can formally be defined as a relaxed TR-RB variant (R-TR-RB).
This variant has been published within \cite{paper4}, along with the results that are conducted in \cref{chap:TR_TSRBLOD}.

Let us begin with a further discussion on problematic scenarios that may occur for the basic TR-RB algorithm.
In many cases, the initial TR-radius chosen in the algorithm does not fit the underlying model and optimization problem.
Thus, it is very likely to happen that, especially at the beginning of the TR-RB algorithm, the sufficient decrease condition \eqref{eq:suff_decrease_condition} is not fulfilled after the sub-problem is terminated according to \eqref{Termination_crit_sub-problem}, which causes the algorithm to reject the parameter and shrink the TR-radius.
This situation may repeat multiple times before a suitable TR-radius is found to ensure convergence of the algorithm.

If the TR-radius is shrunk, two occurrences in the algorithm can significantly harm the computational speed of the overall algorithm.
First, as mentioned in \cref{sec:FOM_cost_or_TRRB}, if \eqref{eq:suff_decrease_condition} can not be verified by the cheap conditions, we require FOM evaluations at the current iterate.
In the worst case, the iterate is rejected, and the TR-radius is shrunk even after enrichment.
In \cref{sec:FOM_cost_or_TRRB}, we thus discussed that the FOM information may be used to update the model function, although the iterate is rejected.

A second scenario is that the TR-radius of an accepted sub-problem occurs to be very small, no matter which of the above-described criteria causes a rejection of the iterate.
As a matter of fact, the certified step of the TR-RB algorithm and also later outer iterations potentially are subject to small steps towards the optimum because the cut-off from \eqref{Cut_of_TR} may be tracked too early.
To tackle this, we introduced the possibility of enlarging the TR-radius; cf.~\cref{sec:enlarging}.
However, in an extreme scenario, the computational effort of the TR-RB algorithm may still significantly be harmed.
Importantly, even if the TR-radius is not shrunk initially, the TR-radius can still prevent the algorithm from iterating as far as it potentially could (for instance, if the error estimator for $\J$ is very ineffective).
Furthermore, the choice of a "perfect" TR-radius and shrinking or enlarging factors are problem-dependent and, to our knowledge, always have to be found by trials.

Another important slow-down factor of the TR-RB algorithm is the evaluation of the error estimator $\Delta_{\Jhat_r}$ or the construction of the offline-online decomposed version of it since it involves the computation of Riesz-representatives w.r.t. the inner product of $V_h$.

Following the remarks above, we conclude that the full certification of the model can significantly slow down the convergence of the TR-RB algorithm, although potentially rejected iterates and models would still have been accurate enough to converge faster eventually.
With this in mind, we introduce the relaxed TR-RB algorithm, where certifications of the TR-RB algorithm are essentially ignored for early iterations and the TR methodology is only enforced in later iteration counts.


Let $(\varepsilon_\textnormal{TR}^{(k)})_k$ and $(\varepsilon_\textnormal{cond}^{(k)})_k$ be sequences of relaxation factors for the TR-radius and the sufficient decrease condition, respectively.
Furthermore, we assume that both sequences are decreasing and converge to zero, i.e.
\begin{equation} \label{eq:converge_to_zero}
	\lim_{k \to \infty} \varepsilon_\textnormal{TR}^{(k)} = 0 \qquad \qquad \text{and} \qquad \qquad
	\lim_{k \to \infty} \varepsilon_\textnormal{cond}^{(k)} = 0.
\end{equation} 
After the surrogate model has been initialized with $\mu^{(0)}$, we solve the same problem as \eqref{TRRBsubprob} but with the relaxed trust-region, i.e., for every outer iteration $k$, we solve
the relaxed TR sub-problem
\begin{equation}
	\label{eq:rel_TRRBsubprob}
	\min_{\widetilde{\mu}\in\Paramsad} \cJhatn^{(k)}(\widetilde{\mu}), \qquad \textnormal{ such that } \qquad 
	q^{(k)}(\widetilde{\mu})
	\leq \delta^{(k)} + \varepsilon_\textnormal{TR}^{(k)},
\end{equation}
with $q^{(k)}(\mu) = \frac{\Delta_{\Jhat_\red^{(k)}}({\mu})}{\cJhatn^{(k)}({\mu})}$ as defined in \eqref{TR_radius_condition}. In addition, we relax the TR-related termination criterion \eqref{Cut_of_TR} to
	\begin{equation}
	\label{eq:rel_cut_of_TR}
	\beta_2(\delta^{(k)}+ \varepsilon_\textnormal{TR}^{(k)}) \leq 
	q^{(k)}(\mu)
	\leq \delta^{(k)}+ \varepsilon_\textnormal{TR}^{(k)},
\end{equation}
and again consider
	\begin{equation}
	\label{eq:FOC_sub-problem}
	\big\|\mu^{(k,\el)}-\Proj_{\Paramsad}(\mu^{(k,\el)}-\nabla_\mu \Jhat_\red^{(k)}(\mu^{(k,\el)}))\big\|_2\leq \tau_\textnormal{{sub}},
\end{equation}
equivalent to \eqref{FOC_sub-problem}.
After the sub-problem is terminated, we accept the iterate with the relaxed version of \eqref{eq:suff_decrease_condition}, which reads as 
\begin{align}
	\label{eq:suff_decrease_condition_rel}
	\cJhatn^{(k+1)}(\mu^{(k+1)})\leq \cJhatn^{(k)}(\mu_\textnormal{{AGC}}^{(k)}) +
	\varepsilon_\textnormal{cond}^{(k)} 
\end{align}
and the corresponding relaxed cheap surrogate criteria analogously to \eqref{eq:suf_cond} and \eqref{eq:nec_cond}.
Moreover, we check the same FOM-based termination criterion as in the TR-RB algorithm:
\begin{equation} \label{eq:TR_termination_rel}
	g_h(\mu^{(k+1)}) = \|\mu^{(k+1)}-\Proj_{\Paramsad}(\mu^{(k+1)}-\nabla_\mu\Jhat_h^{(k)}(\mu^{(k+1)}))\|_2
	\leq \tau_{\textnormal{{FOC}}}.
\end{equation}
If the FOM-based criterion is not fulfilled, we enrich the model and again solve \eqref{eq:rel_TRRBsubprob} until convergence.
The rest of the algorithm can be followed analogously to the certified version from \cref{alg:Basic_TR-RBmethod}.

Concerning the convergence, due to conditions \eqref{eq:converge_to_zero} of the relaxation sequences, asymptotic convergence of the R-TR-RB algorithm can be proven along the same lines as in \cref{thm:convergence_of_TR_2}.
We conclude the following theorem:
\begin{theorem}[Convergence of the R-TR-RB algorithm, cf.~\cref{thm:convergence_of_TR_2}]
	\label{Thm:convergence_of_TR}
	Let $(\varepsilon_\text{TR}^{(k)})_k$ and $(\varepsilon_\text{cond}^{(k)})_k$ be relaxation sequences that fulfill \eqref{eq:converge_to_zero} and let sufficient assumptions on the Armijo search to solve \eqref{eq:rel_TRRBsubprob} be given.
	Then, every accumulation point $\bar\mu$ of the sequence $\{\mu^{(k)}\}_{k\in\mathbb{N}}\subset \Params$ generated by the above described R-TR-RB algorithm is an approximate first-order critical point for $\Jhat_h$, i.e., it holds
	\begin{equation}
	\label{First-order_critical_condition__}
	\|\bar \mu-\Proj_\Params(\bar \mu-\nabla_\mu \Jhat_h(\bar \mu))\|_2 = 0.
	\end{equation}
\end{theorem}
Consequently, the first iterations of the R-TR-RB algorithm can be considered a warm start for the algorithm that can certainly lead to early convergence.
In practice, we can expect the same numerical behavior of the two algorithms once $\varepsilon_\text{TR}^{(k)}$ and $\varepsilon_\text{cond}^{(k)}$ are below double machine-precision.
As with all different optimization methods for the discussed optimization problems, it can always happen that the R-TR-RB algorithm finds a different local minimum than the originally proposed TR-RB algorithm.

The R-TR-RB algorithm can be considered a non-certified algorithm in the first iterations, meaning that $\varepsilon_\text{TR}^{(k)}$ and $\varepsilon_\text{cond}^{(k)}$ are chosen large enough such that the error estimator can be ignored.
In conclusion, for these iterations, we do not need to prepare for the efficient computation of the error estimator, which further results in a significant speedup of the algorithm.

\subsection{Non-certified adaptive algorithms}
\label{sec:non-certified_TR}

The relaxed algorithm presented above can further be simplified in terms of certification by removing the convergence and decreasing condition on $(\varepsilon_\textnormal{TR}^{(k)})_k$ and $(\varepsilon_\textnormal{cond}^{(k)})_k$.
Instead, we set $\varepsilon_\textnormal{TR}^{(k)} = \varepsilon_\textnormal{cond}^{(k)} := \infty$ for all $k$.
Thus, the surrogate-based sub-problems are always solved until local convergence and the TR is entirely ignored while every decreasing iterate is accepted.
In fact, the resulting algorithm can not be considered a trust-region algorithm and shall rather be called a non-certified adaptive algorithm.

Due to the missing certification, the algorithm's convergence can not be proven and can only be verified with FOM-based termination criteria.
The sufficient decrease condition \eqref{eq:suff_decrease_condition} can be tracked along with the algorithm (at no cost).

Compared to the TR-RB algorithm, we conclude that the non-certified algorithm does not use any error estimator of the model.
Instead, the surrogate model is trusted unconditionally for every inner iteration.
This also means that the procedure can be used for non-certified surrogate modeling techniques, such as machine learning (ML) based models.

We advised a similar algorithm for a specific application in \cite{AMLenopt}.
In this work, the underlying PDE-constrained optimization problem is a highly complex task, including a three-phase flow system for enhanced oil recovery.
The aim is to maximize the net present value of the recovery process and to choose control points (for instance injection rates) that are dependent on time.
Since the computation of the (vector-valued) objective function is considered a black-box, we were not able to employ RB techniques and, instead, we utilized a deep neural network to approximate the input-output map.
Since the input space depends, e.g., on the time discretization of the flow-model, it cannot be expected that an overall accurate surrogate model can quickly be found.
Similar to what is discussed in this thesis, in such a case, it is not feasible to invest an arbitrary amount of time for generating training data and finding an overall-accurate ML model.

Solving the optimization problem is achieved by an ensemble-based optimization algorithm (with no uncertainty in the problem's geology parameters).
In order to perform a single optimization step in the corresponding FOM method, a significant number of perturbations of the current control point are required to approximate the gradient of the function.
For the adaptive ML-based algorithm, we use these FOM points as training set for constructing the local ML surrogate.
Similar to a non-certified TR-RB algorithm, we then solved a surrogate-based sub-problem until convergence and continued the procedure until the outer optimization algorithm is terminated based on a suitable FOM criterion.
For further information and for suitable references of related fields, we point to \cite{AMLenopt}.

\subsection{Revisit Experiment 1: Further approaches}
\label{sec:experiments_further_approaches}
\begin{figure} \centering
	\input{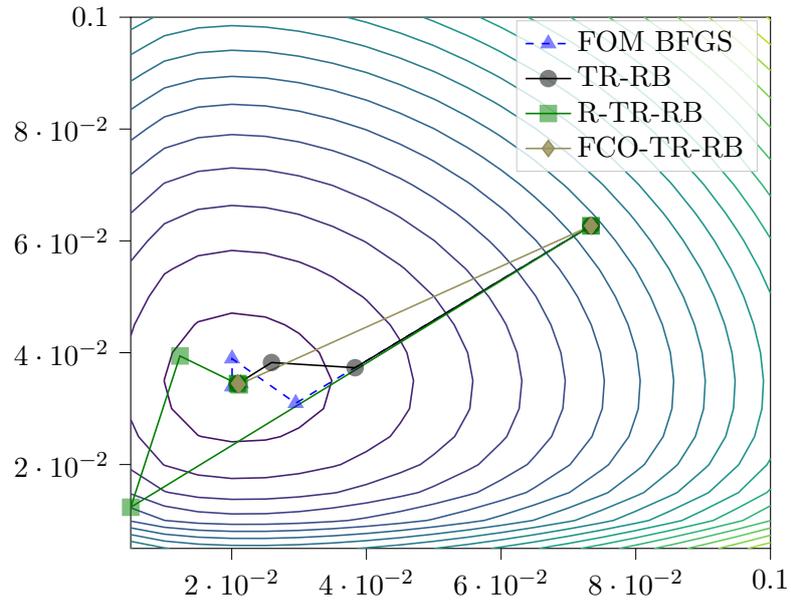}
	\caption{Experiment 1: Optimization path of selected algorithms.}
	\label{fig:NO_TR_for_poc}
\end{figure}

\begin{figure}
	\centering
\begin{tikzpicture}

\definecolor{color1}{rgb}{0.282884,0.13592,0.453427}
\definecolor{color2}{rgb}{0.275191,0.194905,0.496005}
\definecolor{color3}{rgb}{0.258965,0.251537,0.524736}
\definecolor{color4}{rgb}{0.237441,0.305202,0.541921}

\begin{axis}[
legend cell align={left},
legend style={fill opacity=0.8, draw opacity=1, text opacity=1, at={(1,0.5)}, anchor=west,
draw=white!80!black},
width=10cm,
log basis y={10},
tick align=outside,
tick pos=left,
x grid style={white!69.0196078431373!black},
xlabel={time in seconds [s]},
xmajorgrids,
xmin=1.62883254289627, xmax=50.8289027571678,
xtick style={color=black},
y grid style={white!69.0196078431373!black},
ylabel={$\|\muBar - \muBar_h \|_2$},
ymajorgrids,
ymin=2.67990688388142e-10, ymax=0.268952512661232,
ymode=log,
ytick style={color=black}
]
\addplot [semithick, blue, mark=triangle*, mark size=3, mark options={solid, fill opacity=0.5}]
table {%
13.5681953430176 0.0595701054235534
20.2206008434296 0.00454837491910481
22.9678139686584 0.000981823017957724
25.7047007083893 0.000438749415106937
28.6020307540894 2.21513967971858e-05
31.3699910640717 3.11480015640991e-06
34.1037094593048 3.30532942670661e-09
};
\addlegendentry{FOM BFGS}
\addplot [semithick, black, mark=*, mark size=3, mark options={solid, fill opacity=0.5}]
table {%
	24.8370878696442 0.0595701054235534
	38.2350437641144 0.00629077355862692
	41.9597263336182 7.29831583004641e-05
	45.8671171665192 1.71683996399849e-08
};
\addlegendentry{TR-RB BFGS}
\addplot [semithick, green!50!black, mark=square*, mark size=3, mark options={solid, fill opacity=0.5}]
table {%
	4.07095408439636 0.0595701054235534
	7.08349418640137 0.0272133898068323
	10.4665541648865 0.00990211516585644
	13.9668796062469 0.000335775103370097
	17.8843638896942 6.61453159780389e-07
	21.9156887531281 4.97871415991855e-10
};
\addlegendentry{R-TR-RB BFGS}
\addplot [semithick, yellow!50!black, mark=diamond*, mark size=3, mark options={solid, fill opacity=0.5}]
table {%
	14.7062420845032 0.0595701054235534
	18.4671878814697 2.36232028433533e-06
	22.6339585781097 6.15712601700991e-09
};
\addlegendentry{FCO-TR-RB BFGS}
\end{axis}

\end{tikzpicture}
	\caption[Experiment 1: Total wall time comparison of selected algorithms.]{Experiment 1: Total wall time comparison of selected algorithms starting from iteration $k=1$.}
	\label{fig:timings_NO_TR}		
\end{figure}
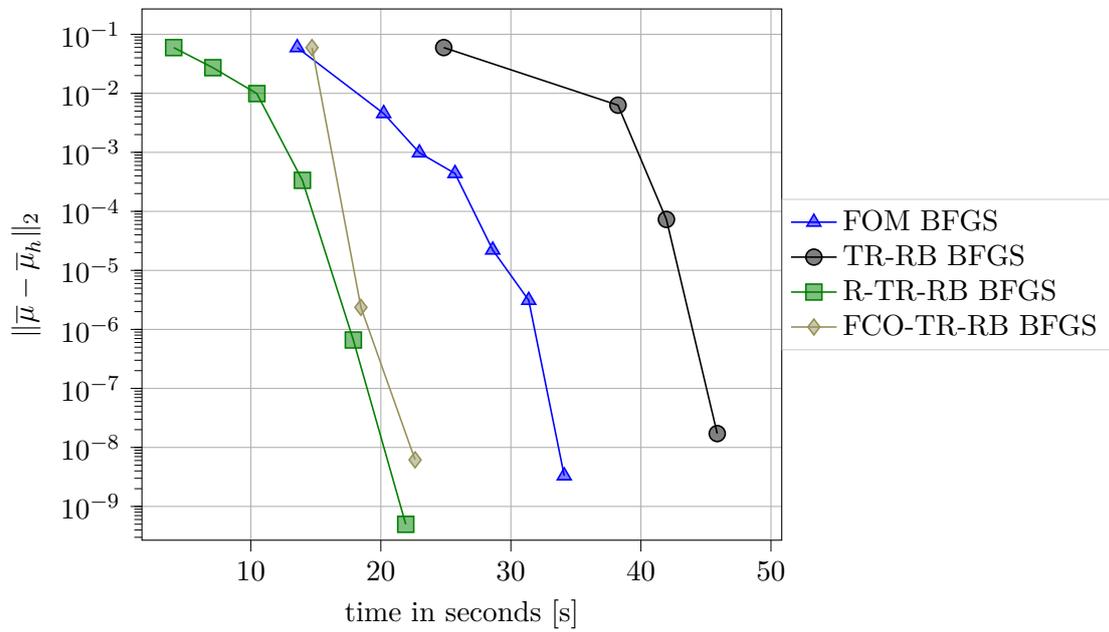

This section is devoted to numerically investigate selected variants that were introduced in this section.
To this end, we apply the presented R-TR-RB algorithm and the FOM-cost-oriented TR-RB to the 2-dimensional example from \cref{sec:proof_of_concept}.
In \cref{fig:NO_TR_for_poc}, we visualize the optimization path of four algorithms:
the TR-RB algorithm used in \cref{sec:proof_of_concept}, the standard FOM-informed BFGS method, the FOM-cost-oriented version (FCO-TR-RB) as explained in \cref{sec:FOM_cost_or_TRRB}, and the R-TR-RB algorithm from~\cref{sec:Relaxed_TRRB}.
For the R-TR-RB, we have chosen $\varepsilon_\textnormal{TR}^{(k)} = \varepsilon_\textnormal{cond}^{(k)} = 10^{(5-k)}$ for all $k$.
We further present additional iteration information in \cref{tab:FCO_TR} and \cref{tab:NO_TR}, and also illustrate the run times of the algorithms in \cref{fig:timings_NO_TR}.

\begin{table}[h] \centering \footnotesize
	\begin{tabular}{|c||c|c|c|c|c|c|} \hline 
		It. $k$ & TR-radius $\delta^{(k)}$ & Inner iterations & TR stopped with & Rejected by & Basis size & $g_h(\mu^{(k)})$\\ \hline \hline
		0 & 0.1 & 4 & \eqref{Cut_of_TR} & \eqref{eq:suff_decrease_condition} & 1 & - \\
		0 & 0.05 & 4 & \eqref{Cut_of_TR} & \eqref{eq:suff_decrease_condition} & 2 & - \\
		0 & 0.025 & 11 & \eqref{FOC_sub-problem} & - & 3 & 6.12e-3 \\ \hline		
		1 & 0.025 & 5 & \eqref{FOC_sub-problem} & - & 4 & 6.08e-5 \\ \hline		
		2 & 0.025 & 6 & \eqref{FOC_sub-problem} & - & 5 & 2.24e-7 \\ \hline					
	\end{tabular}
	\caption[Experiment 1 (revisited): Detailed iteration steps of the FCO-TR-RB algorithm.]{Experiment 1 (revisited): Detailed iteration steps of the FCO-TR-RB algorithm. Depicted are the outer iterations, the TR-radius, inner iterations, details on the TR, basis sizes, and outer stopping.}
	\label{tab:FCO_TR}
\end{table}

\begin{table}[h] \centering \footnotesize
	\begin{tabular}{|c||c|c|c|c|c|} \hline 
		Iteration $k$ & Inner iterations & $q^{(k)}(\mu^{(k)})$ & \eqref{eq:suff_decrease_condition} fulfilled & Basis size
		& Outer stopping $g_h(\mu^{(k)})$
		 \\ \hline \hline
		0 & 3 & 1.34e-02 & yes & 1 & 0.12925 \\ \hline
		1 & 8 & 2.19e-03 & yes & 2 & 0.09420 \\ \hline
		2 & 8 & 2.97e-04 & yes & 3 & 0.00815 \\ \hline
		3 & 5 & 3.09e-09 & yes & 4 & 1.73e-5 \\ \hline
		4 & 4 & 1.06e-15 & yes & 5 & 7.73e-8 \\ \hline
	\end{tabular}
	\caption[Experiment 1 (revisited): Detailed iteration steps of the R-TR-RB algorithm.]{Experiment 1 (revisited): Detailed iteration steps of the R-TR-RB algorithm. Depicted are the outer iterations, inner iterations, further information on the full certification, basis sizes, and outer stopping criterion for \eqref{eq:FOC_sub-problem}.}
	\label{tab:NO_TR}
\end{table}

Although the numerical experiment is relatively simple and only serves as a proof-of-concept, we conclude strong indications that the FCO-TR-RB and, especially, the R-TR-RB approaches show significant improvements to the originally proposed TR-RB from \cref{sec:TR_RB_algorithm}.

\Cref{tab:FCO_TR} can directly be compared to \cref{tab:poc_details}. The FCO-TR-RB algorithm also rejects points in the beginning but, differently from the original TR-RB algorithm, uses the rejected point for adding new basis functions to the ROM.
Consequently, a suitable TR-radius is found earlier, and the outer algorithm only needs two further iterations to converge and has more basis functions in the final ROM.
It can also be seen in \cref{fig:NO_TR_for_poc} that the FCO-TR-RB algorithm immediately steps to the correct region due to the rich RB space from the first sub-problem.

The R-TR-RB method uses an entirely different strategy than the FCO-TR-RB algorithm but shows a comparable computational speedup.
In \cref{tab:NO_TR}, we observe that, due to the relaxation, especially the first iterations indeed converge less sharp to the optimum.
As seen in \cref{fig:NO_TR_for_poc}, the algorithm takes a long way around the optimum.
Importantly, however, this happens at less cost than the TR-RB algorithm since the algorithm rejects no point.
We also note that, for this choice of $\varepsilon_\textnormal{TR}^{(k)}$ and $\varepsilon_\textnormal{cond}^{(k)}$, the certification of the model was never used, which means that a non-certified adaptive algorithm from \cref{sec:non-certified_TR} would have followed the exact same path without using any error estimation (which is even faster).

In \cref{fig:timings_NO_TR}, we see that, just as the greedy-based ROM used as a motivation in \cref{fig:timings_Greedy}, the TR-RB algorithm, introduced in this chapter, also does not give a speed-up w.r.t. the FOM algorithm.
Here, the main reason is the expensive early rejection of the iterates, discussed extensively in the present section.
Instead, the FOM-cost-oriented and relaxed TR-RB algorithm show a significant improvement in terms of computational effort.

For further numerical evidence that the R-TR-RB constitutes a strong improvement to the fully certified variant, we again make reference to \cref{chap:TR_TSRBLOD}.

\section{Summary and outlook}
\label{sec:TR_RB_conclusion}
This chapter introduced recent developments to the TR-RB algorithm for solving PDE-constrained parameter optimization methods with an extensive overview and explanations of variants and additional features.
While the original algorithm was based on \cite{QGVW2017}, we justified the necessity to deviate from the standard ROM choice from~\cite{QGVW2017} and motivated the usage of an NCD-correction term for the reduced model; cf.~\cref{sec:TRRB_MOR_for_PDEconstr}.
Subsequently, we derived a posteriori error estimates for all involved reduced terms and showed that the new ROM techniques indeed enable better estimators for the reduced functional and its gradient; cf.~\cref{sec:TRRB_a_post_error_estimates}.

Apart from the different ROM techniques, the TR-RB algorithm explained in \cref{sec:TR_RB_algorithm} contains significant improvements in terms of parameter constraints, basis enrichment strategies, sub-problem solver, TR-radius enlargement, outer stopping criterion, and parameter post-processing.
The different features and variants were extensively studied in the numerical experiments in \cref{sec:TRRB_num_experiments}, showing that, for the presented problem classes, the suggested TR-RB algorithm, especially with the NCD-corrected approach, is remarkably robust and constitutes a major improvement to~\cite{QGVW2017}.
For our benchmark problems, Newton-type methods with (optional) Taylor- or Lagrange enrichment performed remarkably well.
On the other hand, the Petrov-Galerkin approach showed room for improvements due to the instabilities and was outperformed by the respective Galerkin variant.
As usual for numerical optimization methods, it can not be said which of the TR-RB variants will always converge fastest and which will find a better local minimum.
Certainly, there is a strong indication towards the NCD-corrected approach with Newton-type sub-problem solvers.
In practice, BFGS methods are sometimes still favorable since they do not require explicit Hessian evaluations.

In the last part of this chapter (\cref{sec:TRRB_related_approaches}), we elaborated on still existing potential performance issues of the proposed TR-RB algorithm and we presented ideas to resolve these.
The mentioned AFEM variants tackle the problem of adapting the FEM mesh (even when \cref{asmpt:truth} is fulfilled).
Coarsening of the RB spaces aims to avoid unnecessarily growing RB spaces, and the FOM-cost-oriented variant explicitly uses all FOM information that is gathered during the outer optimization loop.

Importantly, we also introduced the relaxed TR-RB algorithm, where the certification of the TR-RB algorithm is relaxed for the first outer iterations.
Among other advantages, this strategy is beneficial if the error estimator for determining the trust-region has a low effectivity, for instance, caused by extreme overestimation.
In such cases, the R-TR-RB avoids to treat the ROM too conservatively (as done by the fully certified TR-RB algorithm) and allows for maximum convergence of the ROM.
This idea is revisited in \cref{chap:TR_TSRBLOD}, where multiple magnitudes of computational time can be saved.

To motivate the remainder of this thesis, we emphasize that the entire chapter was concerned with global RB methods.
The term "global" refers to the spatial FEM discretization of the system, not to the parameter space, which was clearly localized by the TR method.
The presence of an accurate global discretization internally hides the assumption that a mesh size $h$ in order to fulfill \cref{asmpt:truth} is not arbitrarily small, such that FEM is capable of computing FOM snapshots.
Although higher computational costs for FOM snapshots generally amplify the advantage of TR-RB methods (since these methods avoid expensive FOM computations), global FOM evaluations may as well become prohibitively costly.
This scenario has already been motivated in connection to multiscale problems and the RB-challenge of an inaccessible global discretizations; cf. \cref{sec:multiscale_problems} and \cref{sec:RB_challenges}, respectively.

Transferring the presented TR-RB algorithm to scenarios where local FOM methods are used to attain a solution of the high-fidelity system is relatively straightforward.
To be precise, only the ROM strategy and their related error estimates need to be adjusted.
In \cref{chap:TR_TSRBLOD}, we develop such an algorithm with the help of a localized ROM that is based on a multiscale approach.
The construction of such a localized ROM is one of the primary aims of the subsequent chapter.
\chapter{Two-scale reduced basis localized orthogonal decomposition method}\label{chap:TSRBLOD}

In \cref{chap:background}, we saw that the numerical approximation of parameterized PDEs can raise multiple issues.
Since real-world applications and the involved PDEs become more and more complex, large- and multiscale PDEs have gained a significant research interest.
In cases where standard numerical approximation methods such as FEM fail, localized methods that do not require global solutions can be used.
In the first place, such methods were developed for a (single) deterministic problem without arbitrary or parametric changes in the data.
While the research development of such methods is vast, in this chapter, we are mainly concerned about the class of multiscale problems, already motivated in \cref{sec:multiscale_problems}.

In \cref{sec:RB_challenges}, we formulated the corresponding RB-challenge of an inaccessible global discretization for parameterized systems.
In order to remedy this challenge, we require methods that combine MOR techniques and localized numerical solution methods for parameterized large- or multiscale PDEs.
In this chapter, we use the localized orthogonal decomposition (LOD) method with the primary intention to derive an RB-based MOR approach to solve parameterized multiscale PDEs efficiently.

This chapter is organized as follows: In \cref{sec:localized_MOR}, we start with a general view of localized model reduction for parameterized systems.
In \cref{sec:LOD}, we introduce the LOD as a particular example of a well-established multiscale method with special emphasis on the computational complexity and the related issue of a many-query context.
Besides, we briefly introduce the method from \cite{HKM20} that aims at lowering the computational effort of the LOD for perturbations of a reference problem but does not use an RB-based idea.
In \cref{sec:two_scale}, we provide a two-scale formulation of the LOD, which we analyze in terms of stability and error bounds.
These considerations build the groundwork for an efficient RB-based reduced approach which we present in \cref{sec:TSRBLOD_MOR}.
Afterwards, we demonstrate the new reduction method in dedicated numerical experiments.

\section{Localized model order reduction}
\label{sec:localized_MOR}

This section aims to give a brief introduction to localized model reduction methods for numerically approximating PDEs.
Roughly speaking, a localized model reduction approach is used for cases where the most acceptable scale of the required solution space is not (or only seldom) globally present in the solution method and, instead, low dimensional spaces are constructed.
Thus, in the FE setting, the high-fidelity space $V_h$ with possibly tiny mesh size $h$ is not explicitly used as a global ansatz space.
Instead, the global space is implicitly reconstructed by low-dimensional spaces.
Hence, in what follows, we assume that the global use of $V_h$ results in a system with too many degrees of freedom to fit into the computer resources, and we note that such a mesh size $h$ can always be found.
Let us also mention that the benefit of localized approaches is not limited to such cases, for instance, where the data is only subject to (local) modifications or parameterizations; see, e.g., \cite{BIORSS21} or \cref{chap:TR_TSRBLOD}.

To still accomplish a globally accurate solution, we assume to be given a family of subspaces $\{V^i\}_{i=0, \dots, N_{V}}$, where, for each subdomain, we have $V^i \subset V_h$ and $\dim(V^i) \ll \dim(V_h)$.
Importantly, in order to also capture numerical multiscale approaches with this abstract definition, not all subspaces need to be local.
In particular, the low dimension can also stem from a coarse FE space, as is the case in numerical multiscale methods.
A corresponding direct sum can be build as
\begin{equation}\label{eq:direct_sum}
	X_h = \bigoplus_{i=0}^{N_V} V^i.
\end{equation}
The main idea behind such a space decomposition is that the fine mesh can be decomposed into sufficiently small dimensional subspaces, hoping that the degrees of freedom in the resulting systems are of feasible size for a solution method that fits into the memory or is solvable in a moderate amount of time.
Note that we did not specify that the direct sum space $X_h$ can indeed be interpreted as a representation of a fine-mesh FE space, such that $V_h = X_h$.
Furthermore, $V_h$ does not necessarily have to be the standard FE space.
Also, other choices are possible, for instance, broken Sobolev spaces for discontinuous Galerkin methods.

The variety of localized methods is devastating and can not be tackled in this section.
The main questions for such methods are, first, how the subspaces $V^i$ are chosen, second, what local problems are solved on the subspaces, and third, how the local solutions are coupled to eventually obtain a global solution on an appropriate global approximation space.
These questions are much related and heavily depend on the localized method and the problem to be solved.

In short, we classify localized model reduction into two basic strategies, depending on whether $X_h = V_h$ holds or not.
Examples where $X_h = V_h$ is fulfilled, are, e.g., domain decomposition (DD) methods or discontinuous Galerkin (DG) approximation techniques.
In contrast, in numerical multiscale methods, the space $X_h$ in \eqref{eq:direct_sum} instead serves as a theoretical result and does not directly belong to the solution method.
In most of the cases, however, localized methods are built on a DD-like associated coarse mesh which can be identified by a non-overlapping family of subdomains $\{\Omega_i\}_i$ of $\Omega$.
Importantly, even if the subspaces $V^i$ are directly associated with the respective subdomains $\Omega_i$, it is not necessarily given that also the subspaces are non-overlapping, meaning that $V^i$ is simply the truncated FE mesh of $V_h$ on $\Omega_i$.
Oversampling techniques, e.g., with small fine-mesh layers in DD methods or coarse-element patches as in the LOD, are very widely used.
Let us clarify that such methods can already be phrased as "localized model reduction methods", although they are, in the first place, only formulated for a single deterministic problem.

In a parameterized setting, constructing reduced spaces $V_\red^i$ for the subspaces $V^i$ is useful
and the reduced localized model reduction method can then be interpreted as considering the reduced version of \eqref{eq:direct_sum}, i.e.
\begin{equation}\label{eq:direct_sum_reduced}
X_\red = \bigoplus_{i=0}^{N_V} V^i_\red.
\end{equation}
Such an approach is not limited to the parameterized setting but is also commonly used for accelerating the overall solution procedure of a single non-parameterized problem.
Such a general view on localized model order reduction has been given in \cite{BIORSS21}.

The primary purpose of this chapter is to describe a localized model reduction approach based on a multiscale method that is also capable of handling parameterized multiscale problems.
Let us start by introducing the LOD method.

\section{Localized orthogonal decomposition method}
\label{sec:LOD}

In this section, we review the basic concepts of the LOD utilized for a fixed parameter $\mu \in \Params$.
Moreover, we classify the LOD method w.r.t. \eqref{eq:direct_sum}.
The primary idea of the LOD is to decompose the solution space into a subspace $\Vfh$ of negligible fine-scale variations and an $a_\mu$-orthogonal low-dimensional coarse space of multiscale functions, in which the solution is approximated.
This multiscale space is constructed by computing suitable fine-scale corrections $\QQm(u_H)$ of functions from a given coarse FE space $V_H$.
Due to the dampening of high-frequency oscillations by $a_\mu$, these corrections can then be approximated by the solution of decoupled localized corrector problems.

In recent years, there have been various formulations of the LOD in terms of localization, interpolation, approximation schemes, and applications.
Since we are concerned with the case where storage restrictions may prevent storing explicit solutions of the corrector problems, we focus on the Petrov--Galerkin version of the LOD (PG--LOD)~\cite{elf}, which is favorable in this respect.
For further background, we refer to the LOD-introductory book~\cite{LODbook}.

\subsection{Preliminaries}

To recall from \cref{sec:param_problems}, let $\Params \subset \mathbb{R}^P, P \in \mathbb{N}$, be a parameter space and $\Omega \in \mathbb{R}^d$ a bounded Lipschitz domain.
Similar to \cref{eq:classical_PDE}, we consider the prototypical parameterized elliptic partial differential equation: for a fixed parameter $\mu
\in \Params$,
find $u_{\mu}$ such that
\begin{equation}
\begin{aligned}
\label{eq:problem_classic}
- \nabla \cdot A_{\mu}(x)  \nabla u_{\mu}(x) &= f(x), \qquad x \in \Omega, \\
u_{\mu}(x) & = 0, \quad \;\; \qquad x \in \partial \Omega.
\end{aligned} 
\end{equation}
We assume that the parameter-dependent coefficient field $A_{\mu}$ has a multiscale structure that renders a direct solution of~\eqref{eq:problem_classic} using, e.g., finite elements computationally infeasible due to the high mesh resolution required to resolve all features of $A_\mu$, cf. \cref{sec:multiscale_problems}.
Further, we assume $A_{\mu} \in L^{\infty}(\Omega, \mathbb{R}^{d \times d})$ to be symmetric and uniformly elliptic, such that
\begin{align} 
0 < \alpha &:= \essinf\limits_{x \in \Omega} \inf_{v \in \mathbb{R}^d \setminus \set{0}} \frac{\left( A_\mu(x)v \right) \cdot v}{ v \cdot v}, \label{eq:ellipticity_1} \\
\infty > \beta &:= \esssup\limits_{x \in \Omega} \sup_{v \in \mathbb{R}^d \setminus \set{0}} \frac{\left( A_\mu(x)v \right) \cdot v}{ v \cdot v} \label{eq:ellipticity_2},
\end{align}
and we let $\kappa:=\beta/\alpha$ be the maximum contrast of $A_{\mu}$ for $\mu \in \Params$. Moreover, let $f \in L^2(\Omega)$. 

As explained in \cref{sec:preliminaries}, we consider the corresponding weak formulation of \eqref{eq:problem_classic}: For $\mu \in \Params$, we seek $u_{\mu} \in V:=H^1_0(\Omega)$, such that 
\begin{equation}
\label{eq:problem_weak}
a_{\mu}(u_\mu,v) =  l(v) \qquad \textnormal{for all } v \in V.
\end{equation}
We remark that \eqref{eq:problem_weak} differs from the formulation of a parameterized elliptic problem as introduced in \cref{def:param_elliptic_problem} since the right-hand side is assumed to be non-parametric.
The presented approach in this chapter can be generalized easily to parametric $f$ and other boundary conditions, which is not discussed further as it only aggravates the presentation of the reduced method in \cref{sec:TSRBLOD_MOR}.


Finally, to obtain an online-efficient reduced-order model, we assume parameter-separability of the involved data as formulated in \cref{asmpt:parameter_separable}.
Note that this is trivially the case for $f$, and means that, for $A_\mu$, we have a decomposition $A_\mu = \sum_{\xi=1}^{\Qa} \theta_\xi^a(\mu) A_\xi$
with non-parametric $A_\xi \in L^{\infty}(\Omega, \mathbb{R}^{d \times d})$ and arbitrary $\theta_\xi^a:\Params \to \mathbb{R}$.

\noindent Recall that, for $v\in V = H^1_0(\Omega)$, standard (equivalent) norms are defined by:
\begin{equation*}
\snorm{v}^2       := \int_\Omega \abs{\nabla v(x)}^2\dx, \qquad
\anorm{v}^2       := \int_\Omega \abs{A_\mu^{1/2}\nabla v(x)}^2\dx,
\end{equation*}
where $\snorm{}$ is a norm on $V$ due to Friedrich's inequality; cf. \cref{sec:preliminaries}.

\subsection{Discretization and patches}
Let $V_h \subset V$ be a conforming FE space of dimension $\dimVh$, and let $\grid$ be the corresponding
shape regular mesh over the computational domain $\Omega$; cf. \cref{sec:background_npdgl}. 
We again assume that the mesh size $h$ is chosen such that all features of $A_\mu$ are resolved, making a global solution within $V_h$ infeasible.
Further, we assume to be given a coarse mesh $\Gridh$ with mesh size $H \gg h$ that is aligned with $\grid$ and we let
\begin{equation*}
V_H := V_h \cap \mathcal{P}_1(\Gridh),     
\end{equation*}
where $\mathcal{P}_1(\Gridh)$ denotes $\Gridh$-piecewise affine functions that are continuous on $\Omega$.
We denote the dimension of $V_H$ by $\dimVH$.

\begin{figure}[t]
	\centering
	\includegraphics[scale=0.8]{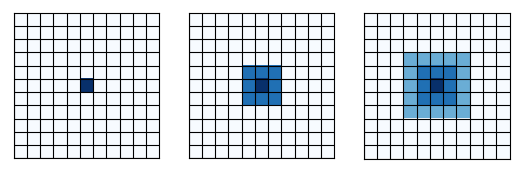}
	\caption[Patches $U_\kl(T)$ for a coarse mesh element $T \in \Gridh$ for $\kl=0,1,2$]{Patches $U_\kl(T)$ for a coarse mesh element $T \in \Gridh$ for $\kl=0,1,2$ (from left to right).}
	\label{patches}
\end{figure}

For an arbitrary set $\omega \subseteq \Omega$, we define coarse grid element patches $U_\kl(\omega) \subset \Omega$ of size $0 \leq \kl \in \mathbb{N}$ by
\begin{equation*}
U_0(\omega) := \omega,  \qquad \text{and}\qquad
U_{\kl+1}(\omega) := \operatorname{Int}\biggl(\,\overline{\bigcup \set[T \in \Gridh]{ \overline{U_\kl(\omega)} \cap \overline{T} \neq \emptyset}}\,\biggr),
\end{equation*}
where $\operatorname{Int}(X)$ is the interior of the set $X$.
Further, for a given patch size $\kl$, let
\begin{equation*}
\CO := \max_{x\in\Omega} \,\#\{T \in \Gridh \,\vert\, x \in U_\kl(T)\}
\end{equation*}
be the maximum number of element patches overlapping in a single point of $\Omega$.
We visualize these coarse element patches in \cref{patches}, where it can easily be seen that, for quadrilateral meshes, we have $\CO = (2 \cdot \kl + 1)^2$.

\subsection{Localized multiscale space}
\label{sec:localized_ms_space}
To define the fine-scale space $\Vfh \subset V_h$, we consider a (quasi-)interpolation operator
$\IH : V_h \to V_H$ mapping a high-fidelity function $v_h \in V_h$ to an element of the coarse FE space $V_H$.
Hence, fine-scale features which the interpolation operator neutralizes can be identified by the respective kernel.
Thus, we define
$$\Vfh := \ker(\IH) \subset V_h.$$
Multiple choices for such an interpolation operator are possible (see~\cite{MH16, PS} for an overview). 
In recent literature, using an operator that is based on local $L^2$-projections~\cite{Pe15} has proven advantageous.
However, for the following error analysis we only require that $\IH$ is linear, $\snorm{}$-continuous and idempotent
on $V_H$, i.e.,
\begin{align*}
\IH(v_H) &= v_H  \qquad\qquad\qquad\quad\textnormal{for all } v_H \in V_H,\\
\snorm{\IH(v_h)} &\leq \CI \snorm{v_h} \qquad\qquad \textnormal{for all } v_h \in V_h.
\end{align*}
Next, we define fine-scale corrections $\QQm(v_h) \in \Vfh$, for a given $v_h \in V_h$, to be the solution of
\[
a_{\mu}(\QQm(v_h),\vf) = a_{\mu}(v_h,\vf) \qquad\qquad \textnormal{for all } \; \vf \in \Vfh.
\]
In conclusion, $\QQm(v_h)$ is the $a_\mu$-orthogonal projection of $v_h$ onto $\Vfh$, and the multiscale space
\begin{equation*}
V_{H, \mu}^{\textnormal{ms}} := (I - \QQm) (V_H)
\end{equation*}
is the $a_\mu$-orthogonal complement of $\Vfh$ in $V_h$, i.e., $V_h = V_{H,\mu}^{\textnormal{ms}} \oplus_{a_{\mu}} \Vfh$.

We need a computable basis to use the multiscale space in a practical implementation.
Since $V_{H, \mu}^{\textnormal{ms}}$ and $V_H$ have equal dimensions, it suffices to apply the fine-scale corrector $\QQm$ on every basis function $\phi_x$ of $V_H$ to obtain a corrected basis, i.e.
\[
\set[\phi_x - \QQm \phi_x]{x \in \mathcal{N}_h},
\]
where $\mathcal{N}_h$ denotes the set of nodes in $\mathcal{T}_h$.
Importantly, even though $\phi_x$ have small local support, $\QQm(v_h)$ will have global support,
and its computation will require the same effort as a global solution of \eqref{eq:problem_weak} in $V_h$.
We further note that the decomposition $V_h = V_{H,\mu}^{\textnormal{ms}} \oplus_{a_{\mu}} \Vfh$ can thus not (yet) be interpreted as a localized approach in the sense of \eqref{eq:direct_sum}.

\begin{figure}[t]
	\begin{subfigure}[b]{0.3\textwidth}
		\includegraphics[width=\textwidth]{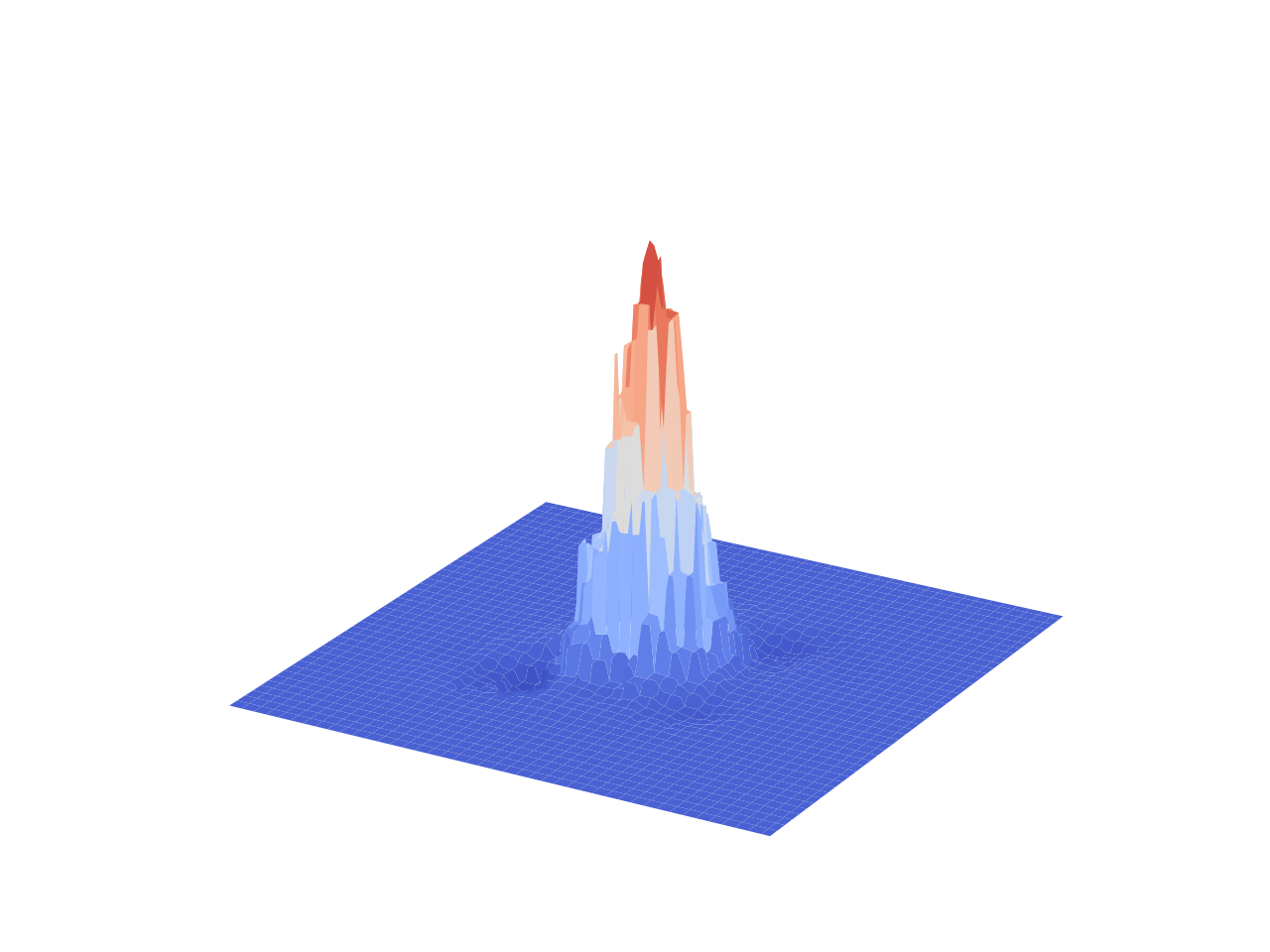}
		\caption{$\phi_x - \QQm \phi_x$.}
	\end{subfigure}
	\begin{subfigure}[b]{0.3\textwidth}
		\includegraphics[width=\textwidth]{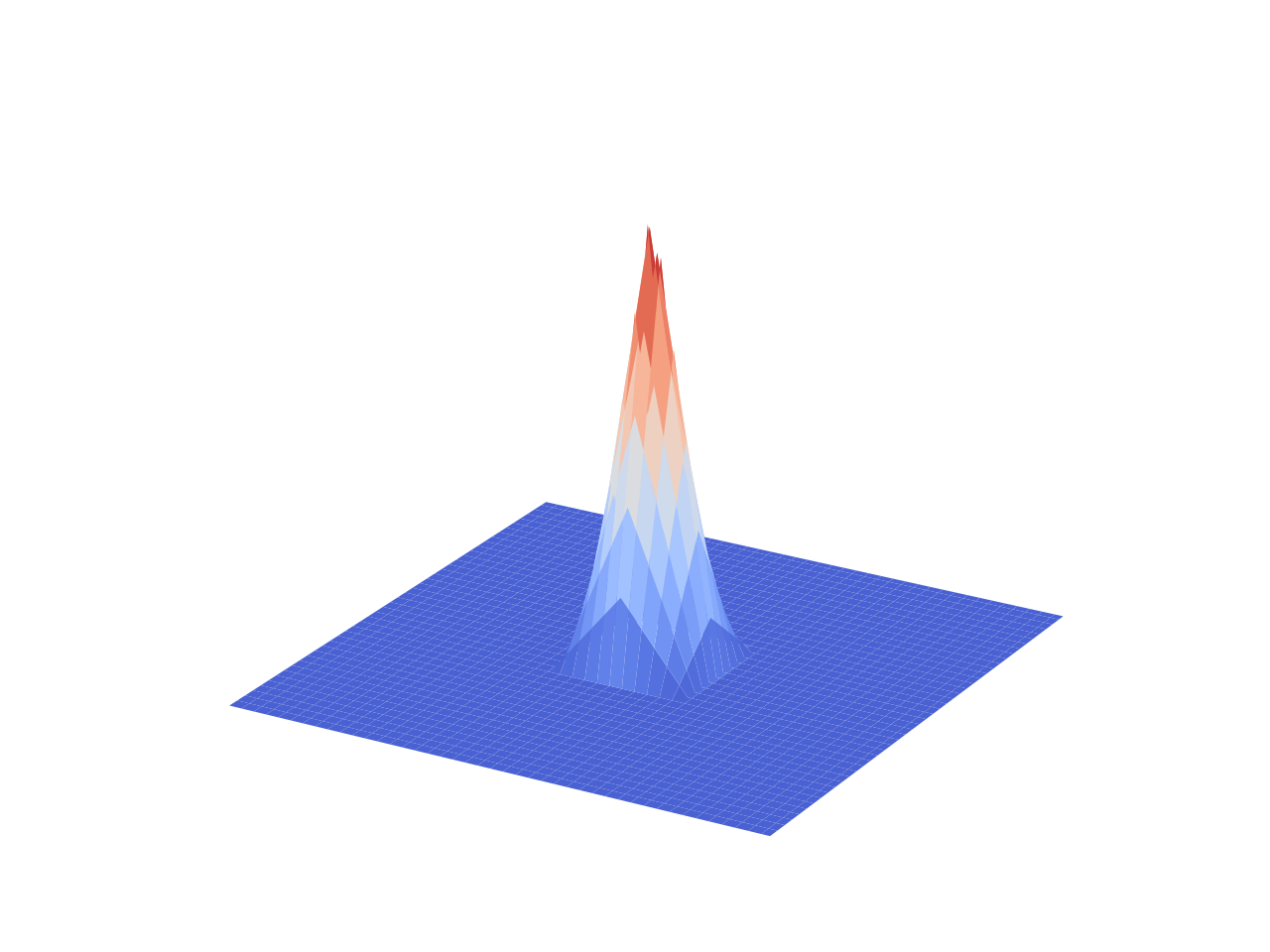}
		\caption{$\phi_x$.}
	\end{subfigure}
	\begin{subfigure}[b]{0.3\textwidth}
		\includegraphics[width=\textwidth]{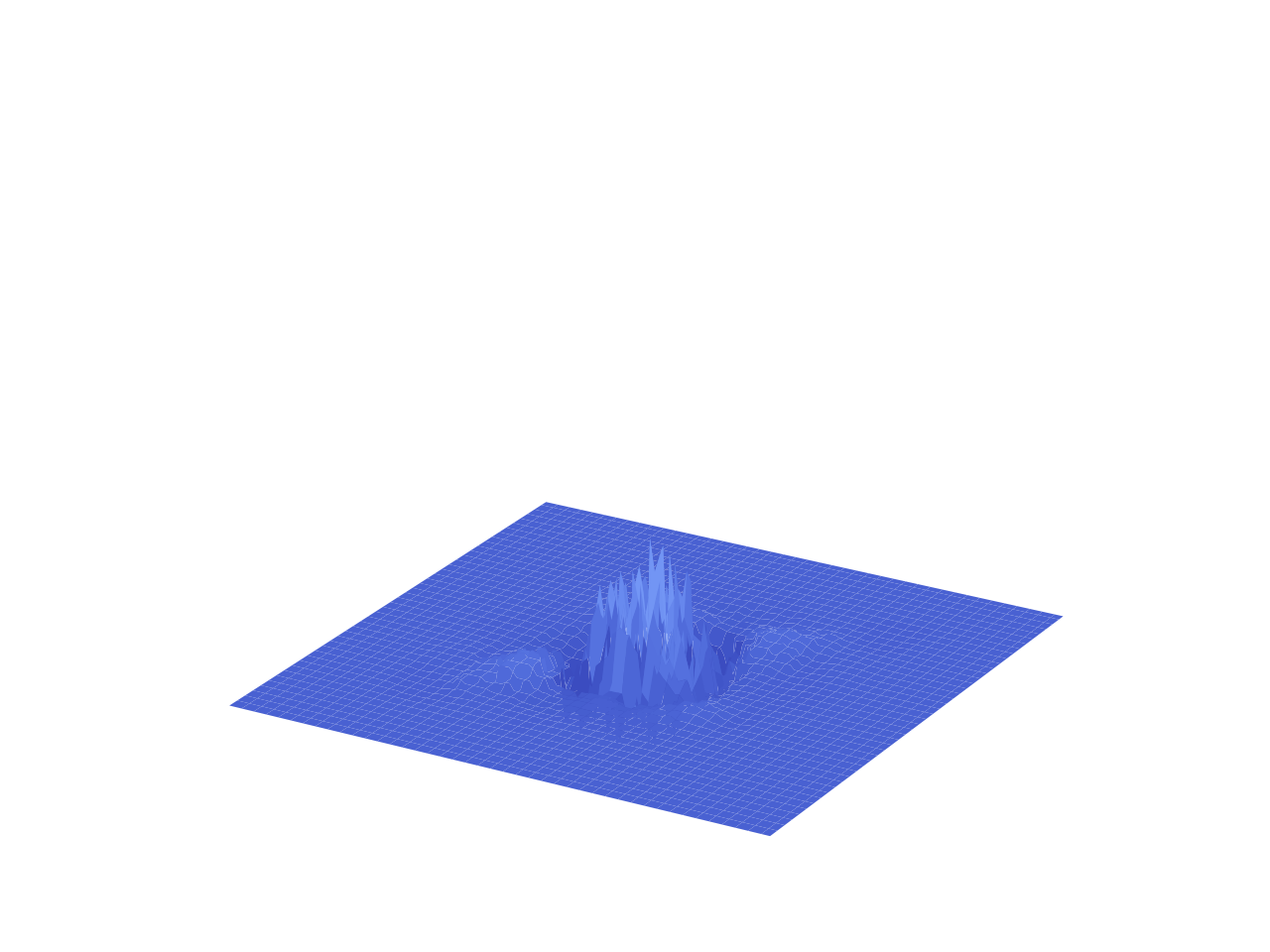}
		\caption{$\QQm \phi_x$.}
	\end{subfigure}
	\caption{Basis function of $V_{H, \mu}^{\textnormal{ms}}$ and its decomposition for $x \in \mathcal{N}_h$.}
	\label{correctorplots}
\end{figure}

Luckily, as illustrated in \cref{correctorplots}, and first proven in \cite{MP14}, the computation can be localized to a small area around the support of each shape function since the associated corrector decays sufficiently fast; see below.

\noindent We approximate $\QQm(v_h)$ using localized correctors $\QQktm(v_h) \in \Vfhkt$ in 
the patch-restricted fine-scale spaces $\Vfhkt: = \Vfh \cap H^1_0(U_\kl(T))$ given by
\begin{align} \label{eq:loc_cor_problems}
a_{\mu}(\QQktm(v_h),\vf) = a^T_{\mu}(v_h, \vf) \qquad \qquad \textnormal{for all } \; \vf \in \Vfhkt,
\end{align}
where $a^T_{\mu}$ denotes the bilinear form obtained by restricting the integration domain in the definition of $a_{\mu}$ to $T \in \Gridh$. 
We then define the localized corrector by 
\begin{equation}\label{eq:loc_corrector}
\QQkm := \sumT\QQktm.
\end{equation}
The exponential decay of the correctors is formulated in the following error bound (see \cite{MP14}):
\begin{equation}\label{eq:decay}
\anorm{(\QQm - \QQkm) v} \leq C \kl^{d/2} \theta^\kl \anorm{v}, 
\end{equation}
with constants $C, \theta$ independent of $H$ and $\kl$, and where $0 < \theta < 1$.
For a detailed discussion on the decay variable $\theta$, we refer to~\cite{MH16, HM19, MP14}.
The fast decay of the truncation error shows that a patch localization parameter of $\kl \approx |\log H|$ suffices in order to consider the localization error negligible.

To conclude, as a replacement of $V_{H, \mu}^{\textnormal{ms}}$, we consider the localized multiscale space $V_{H,\kl,\mu}^{\textnormal{ms}}$ given by
\[
V_{H,\kl,\mu}^{\textnormal{ms}}:= (I - \QQkm)(V_H).
\]

With respect to the originally discussed idea of localized model reduction approaches, concerning the choice of the subspaces $V^i \subset V_h$ in \eqref{eq:direct_sum}, we note that with $N_V = \abs{\mathcal{T}_H} $, we have $V^0 := V_H$, and $V^i := \Vfhkt[T_i]$ for $i=1,\dots,\abs{\mathcal{T}_H}$.
Therefore, the localized multiscale space $V_{H,\kl,\mu}^{\textnormal{ms}}$ can be reinterpreted by the following direct sum
$$ \Vts := V_H \oplus \Vfhkt[T_1] \oplus \cdots \oplus \Vfhkt[\Tlast],$$
in such a way that a function $\uhkmsm \in V_{H,\kl,\mu}^{\textnormal{ms}}$ can uniquely be described by the tuple 
$$\uts_\mu = \left[ \uhkm,\, \QQktm[T_1](\uhkm),\, \ldots,\, \QQktm[\Tlast](\uhkm)\right] \in \Vts.$$
For more details on the coherence of the two spaces, we refer to \cref{prop:two_scale_solution}.
We emphasize that $V^0 = V_H$ is not a localized subspace but still has a low dimension due to the coarse discretization.
In \cref{sec:two_scale}, the corresponding space $\Vts$ is introduced as the two-scale function space for the LOD, which is helpful for the two-scale reduction approach; cf.~\cref{sec:RB_for_LOD}.

\subsection{Petrov--Galerkin projection}
\label{sec:PG--LOD} 
After computing $V_{H,\kl,\mu}^{\textnormal{ms}}$, we determine an approximation of $u_\mu$ in the $N_H$-dimensional multiscale space via Petrov-Galerkin projection, i.e., we let $u_{H,\kl,\mu}^{\textnormal{ms}} \in V_{H,\kl,\mu}^{\textnormal{ms}}$ be the solution of
\begin{equation}
\label{eq:PG_}
a_{\mu}(\uhkmsm,v_H) =  l(v_H) \qquad \textnormal{for all } \; v_H \in V_H.
\end{equation}
We note that the standard Galerkin formulation can be obtained by using $V_{H,\kl,\mu}^{\textnormal{ms}}$ also as the test function space; see \cite{LODbook}.
To ensure that~\eqref{eq:PG_} has a unique solution, inf-sup stability of~$a_\mu$ w.r.t.\ $V_{H,\kl,\mu}^{\textnormal{ms}}$ and $V_H$ is required, which has been shown in~\cite{elf, HM19}, and is conditioned on sufficiently large~$\kl$.
Compared to these references, we use a slightly different definition of the inf-sup stability constant for~\eqref{eq:PG_} by using $\snorm{}$, instead of $\anorm{}$, for the test space:
\begin{equation*}
\CPG := \inf_{0\neq w_H\in V_H}\sup_{0\neq v_H\in V_H}
\frac{a_\mu(w_H-\QQktm(w_H), v_H)}{\anorm{w_H-\QQktm(w_H)}\snorm{v_H}}.
\end{equation*}
For $\kl$ large enough, we have that
\begin{equation}\label{eq:gamma_kl}
\CPG \approx \alpha^{1/2}\CI^{-1},
\end{equation}
which can be proven with a simple modification of the argument in~\cite[Section 4]{HM19}.
In particular, w.l.o.g., we assume that $\CPG \leq \alpha^{1/2}$.

Writing the solution of~\eqref{eq:PG_} as $\uhkmsm=\uhkm - \QQkm(\uhkm)$ with $\uhkm \in V_H$,
we have the following a priori estimate, which was first shown in~\cite{elf}.
\begin{theorem}[A priori convergence result for the PG--LOD]
	\label{thm:pglod_convergence}
	For a fixed parameter $\mu \in \Params$, let 
	$u_{h,\mu} \in V_h$ be the finite-element solution of \eqref{eq:problem_weak} given by
	\begin{equation*} \label{eq:fem_}
	a_\mu(u_{h,\mu}, v_h) = l(v_h) \qquad\textnormal{for all } v_h \in V_h.
	\end{equation*}
	Then, it holds that
	\begin{equation*}
	\norm{u_{h,\mu} - u_{H,\kl,\mu}}_{L^2} + \norm{u_{h,\mu} - u_{H,\kl,\mu}^{\textnormal{ms}}}_{1} \lesssim (H + \theta^\kl
	\kl^{d/2}) \norm{f}_{L^2},
	\end{equation*}
	with $0 < \theta < 1 $ independent of $H$ and $\kl$, but dependent on the contrast $\kappa = \beta/\alpha$.
\end{theorem}
\noindent Although our setting is slightly different from the one in \cite{elf} (in terms of localization and interpolation), the proof can still be followed analogously.
We emphasize that, since $\theta$ depends on the contrast of the problem, the LOD is generally vulnerable to high-contrast problems.
For neglecting the issue of high contrast in the LOD, the interpolation operator $I_H$ has to be adjusted.
Works in this direction have been made in \cite{brown2016multiscale,MH16,PS}, for instance.

\begin{remark}[Right-hand side correction for the PG-variant] \label{rem:rhs_correction}
	As explained in \cite{HKM20,HM19}, a right-hand side correction can be used to improve the accuracy of the PG-variant, which accounts for the fact that the test functions do not inherit the fine-scale features to the right-hand side of \eqref{eq:PG_}.
	For the additional correction terms, we require to solve yet another set of equations like \eqref{eq:loc_cor_problems}, doubling the computational effort on patches.
	While the use of a right-hand side correction is recommended for attaining maximum accuracy for the PG-variant, we note that it is relatively easy to generalize the results of this section to it.
	Moreover, as we shall see later, our approach is anyway mainly interested in the coarse-scale approximations $u_{H,\kl,\mu} \in V_h$ of our system, where storing (local) fine-scale functions to compute $u_{H,\kl,\mu}^{\textnormal{ms}} \in V_{H,\kl,\mu}^{\textnormal{ms}}$ may already be avoided anyway, cf. \cref{rem:coarse_vs_fine} and \cref{sec:TSRBLOD_experiments}.
\end{remark}

\subsection{Computational aspects}
\label{sec:comp_comp_LOD}

We are concerned with the computational effort of the PG--LOD scheme and intend to apply it for multiple samples from $\Params$.
In order to solve \eqref{eq:PG_}, we need to solve multiple corrector problems \eqref{eq:loc_cor_problems} for every $T \in \Gridh$. 
These correctors are then used to assemble a localized multiscale matrix $\Kmu$ given by
\begin{equation} \label{eq:K}
\Kmu:= \sumT \Ktmu, \quad \left(\Ktmu \right)_{ji} := \skal{A_{\mu} (\rchi_T \nabla - \nabla \QQktm) \phi_i}{\nabla \phi_j}_{U_\kl(T)},  
\end{equation}
where $\rchi_T$ denotes the indicator function on $T$ and $\phi_i$ the finite-element basis functions of $V_H$. 
Solving \eqref{eq:PG_} is equivalent to solving the linear system
\begin{equation}\label{eq:PG_in_matrices}
\Kmu \cdot \uhkmsmcoeff = \mathbb{F},
\end{equation}
where $\mathbb{F}_i := l(\phi_i)$, and $\uhkmsmcoeff \in \mathbb{R}^{N_H}$ is the vector of coefficients of $\uhkmsm \in V_{H,\kl,\mu}^{\textnormal{ms}}$ w.r.t.\ the basis of $V_{H,\kl,\mu}^{\textnormal{ms}}$ corresponding to the finite-element basis $\phi_i$; cf. \cref{sec:algebraic_fem} for a similar derivation of FE matrices.

Compared to the Galerkin projection onto $V_{H,\kl,\mu}^{\textnormal{ms}}$, the system matrix $\Kmu$ of the Petrov-Galerkin formulation has a smaller sparsity pattern, and we need less computational work to assemble the matrix.
Every localized corrector and hence each local contribution to $\Kmu$ can be computed in parallel, without any communication and deleted after the contribution $\Ktmu$ has been computed.
In particular, the local contribution matrices $\Ktmu$ only have non-zeros in columns $i$ for which $T \subseteq \supp \phi_i$.

Overall, we summarize the computational procedure of the LOD in the following definition:

\vspace{0.3cm}
\begin{definition}[Computational procedure of the PG--LOD] \label{def:comp_proc_LOD}
	For a new parameter $\mu \in \Params$, solving \eqref{eq:PG_} means to follow
	\begin{description}
		\item[Step 1] For every $T \in \Gridh$: Compute $\QQktm(\phi_i)$ by solving \eqref{eq:loc_cor_problems} for each $i$, s.t. $T \subseteq \supp \phi_i$. Assemble $\Ktmu$ according to \eqref{eq:K}. 
		\item[Step 2] Assemble the localized multiscale stiffness matrix $\Kmu = \sumT \Ktmu$ from the local contributions computed in Step~1.
		\item[Step 3] Solve equation \eqref{eq:PG_in_matrices} to compute $\uhkmsmcoeff$.
	\end{description}
\end{definition}
\vspace{0.3cm}

\begin{figure}[t] \centering \footnotesize
	\begin{tikzpicture}[
	node distance = 6mm and 12mm,
	start chain = A going below,
	base/.style = {draw, minimum width=10mm, minimum height=10mm,
		align=center, on chain=A, rectangle, fill=red!5},
	every edge quotes/.style = {auto=right}]
	\node [base, minimum width=20mm, minimum height=20mm, fill=black!10] {Coarse LOD system \\\eqref{eq:PG_in_matrices} on \textcolor{blue}{$\mathcal{T}_H$} };
	\node [base, above=1cm of A-1, xshift=-2cm] {\footnotesize{$T_0$}\\corrector FOM \\ \eqref{eq:loc_cor_problems} on \textcolor{red}{$U_\kl(T_0)_h$}};
	\node [base, above=1cm of A-1, xshift=1.5cm] {\footnotesize{$T_1$}\\corrector FOM \\\eqref{eq:loc_cor_problems} on \textcolor{red}{$U_\kl(T_1)_h$}};
	\node [base, above=0.1cm of A-1, xshift=4cm, yshift=-0.8cm] {\footnotesize{$T_2$}\\corrector FOM \\\eqref{eq:loc_cor_problems} on \textcolor{red}{$U_\kl(T_2)_h$}};
	\node [base, above=-1.6cm of A-1, xshift=4cm, yshift=-0.8cm] {\footnotesize{$T_3$}\\corrector FOM \\\eqref{eq:loc_cor_problems} on \textcolor{red}{$U_\kl(T_3)_h$}};
	\node [base, above=-3.6cm of A-1, xshift=2.5cm, yshift=-0.5cm] {\footnotesize{$T_4$}\\corrector FOM \\\eqref{eq:loc_cor_problems} on \textcolor{red}{$U_\kl(T_4)_h$}};
	\draw [->] (A-2) -- (A-1) node[midway, left] {$\K_{T_0,\mu}$};
	\draw [->] (A-3) -- (A-1) node[midway, left] {$\K_{T_1,\mu}$};
	\draw [->] (A-4) -- (A-1) node[midway, above] {$\K_{T_2,\mu}$};
	\draw [->] (A-5) -- (A-1) node[midway, above] {$\K_{T_3,\mu}$};
	\draw [->] (A-6) -- (A-1) node[midway, right] {$\K_{T_4,\mu}$};
	
	\draw [->] (-0.7,-1.8) -- (A-1);
	\fill[fill=black] (-0.5,-2.2) circle [radius=0.08];
	\fill[fill=black] (-0.7,-2.2) circle [radius=0.08];
	\fill[fill=black] (-0.9,-2.2) circle [radius=0.08];				
	\end{tikzpicture}
	\caption{Visualization of the PG--LOD FOM procedure.}
	\label{fig:FOM_LOD}
\end{figure}
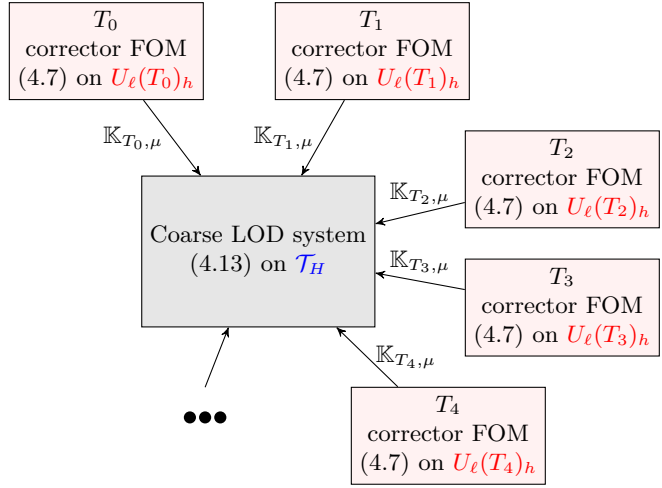

\noindent In general, neither of Steps 1--3 is computationally negligible, and each of these steps must be repeated to obtain a solution for a new $\mu \in \Params$.
Step~1 requires computations on the fine-scale level, whereas Steps 2 and 3 solely depend on the coarse mesh size $H$.
The above-explained procedure can be interpreted as a FOM method of the PG--LOD.
Therefore, we consider the corrector problems in Step~1 as FOM corrector problems.
For a visualization of the above-explained procedure, we refer to \cref{fig:FOM_LOD}.

\begin{remark}[Coarse-scale solution vs. fine-scale solution] \label{rem:coarse_vs_fine}
	In Step~3, we can only compute $\uhkmsm$ if the basis of $V_{H,\kl,\mu}^{\textnormal{ms}}$ is available.
	Again, right-hand side corrections can be used to further improve the approximation quality; cf. \cref{rem:rhs_correction}.
	If, instead, the storage consumption of storing the correctors for $V_{H,\kl,\mu}^{\textnormal{ms}}$ is too high, in Step~3, we can use $\uhkmsmcoeff$ to compute the coarse-scale representation $u_{H,\kl,\mu} \in V_H$.
	As remarked in \cite{elf,EPMH16}, in many applications, using solely $u_{H,\kl,\mu} \in V_H$ is enough.
\end{remark}

\begin{remark}[Periodicity]
	Applications such as composite materials often lead to a periodic structure of the data functions. 
	In this case, correctors may be reusable, which, in the full periodic case, leads to only one corrector problem for full patches and comparably few for the boundary patches.
	This means that memory consumption, as well as complexity, decrease substantially.
	We do not assume the periodicity of the data functions and instead consider a largely non-periodic setting.
\end{remark}

\subsection{Minimizing the cost in a many-query scenario}
\label{sec:minimizing_cost}
The previous section showed that LOD evaluations for multiple changing fine-scale data require recomputing Steps~1--3.
If these configurations are similar or follow a specific parameterization, fine-scale or even coarse-scale functions from the computation steps may be reused.
As discussed in \cref{chap:introduction}, a large number of proposals have been made to speed up the simulation process of solving multiscale problems with similar or parameterized data functions.
For the LOD, in particular, such related approaches can be found in \cite{RBLOD,HKM20,HM19,MV21}.
While the authors in \cite{RBLOD} consider parameterized problems in the sense of \cref{def:param_elliptic_problem} and apply classic MOR to the Galerkin formulation of the LOD, the works \cite{HKM20,HM19,MV21} follow different approaches and use the Petrov--Galerkin formulation of the LOD.

In~\cite{HM19}, similar coefficients over time are considered with the idea that a large number of corrector problems in Step~1 from \cref{def:comp_proc_LOD} are equal or only slightly changed over time.
With the use of a quickly computable error indicator, the correctors that need to be recomputed for achieving the desired accuracy are ascertained and subsequently recomputed in Step~1.
Consequently, at least some of the resulting local matrices $\tilde{\K}_{T,\mu}$ may be inexact.
Luckily, as shown in the error analysis of \cite{HM19}, the effect of this inexactness for Step~3 can be controlled by the error indicator such that reasonable approximations can be computed.
With this strategy, the computational time for Step~1 can be significantly reduced while the effort for Step~2 and Step~3 remain identical.
However, if at least one local problem needs to be recomputed, the resulting computational effort still scales locally with~$h$.

Based on the same idea, in~\cite{HKM20}, the approach from~\cite{HM19} has been generalized to problems with arbitrary perturbations of a reference structure where the fine-scale structure, or at least the resulting local matrices, are readily available; cf.~\cref{fig1}.
Suppose the perturbations are mild or do not reach the entire computational domain. In that case, a large amount of the reference information can be reused, dependent on the error indicator that detects the necessity for recomputation.
In~\cite{HKM20}, a large variety of perturbations, including domain mappings, have been considered, and the approach has shown significant applicability to, e.g., Monte Carlo simulations in material science.
However, just as before, the computational effort for Step~2 and Step~3 is still equal to the original formulation of the PG--LOD, and fine-scale computations may still be required in Step~1.

The approach in~\cite{MV21} follows a more RB-related offline-online-based strategy for Step~1 using structural assumptions on the data.
In particular, the idea is to invest the computational effort that touches the fine-scale data only offline and to rely on the precomputed information in the online phase entirely.
Again, Step~2 and 3 computationally remain the same and, compared to problems with arbitrary parameterization, the approach requires relatively strong assumptions on the data functions.

As already pointed out, the primary purpose of this chapter consists of an efficient RB-based approach for handling parameterized problems such as in~\eqref{eq:problem_weak} for many-query and real-time scenarios.
An outlook w.r.t. \cref{sec:localized_MOR} is given in \cref{sec:RB_for_LOD}.
Beforehand, we resume the main results of the non-RB approach that has been introduced in \cite{HKM20} as one of the first many-query applications for the LOD.

\subsection{Adaptive LOD algorithm for perturbed problems}
\label{sec:adaptive_PGLOD_from_hekema}
This section is devoted to the many-query approach proposed in~\cite{HKM20} as an extension of the algorithm in~\cite{HM19}.
The idea is to handle a given reference structure's arbitrary (but mostly local) perturbations.
Due to the locality of the perturbations, the LOD method allows for recomputing only a fraction of the corrector problems.

\begin{figure}[t]
	\begin{center}
		{	\begin{tikzpicture}[scale=0.7, every node/.style={transform shape}, node distance=1cm]
	 \node[] (market) {};
	 \node[above=of market] (dummy) {};
	 \node[blue-rec, right=of dummy] (t) {};
	 \node[blue-rec, left=of dummy] (g) {};
	 \node[left=of g] (sec_dummy) {};
	 \node[blue-rec, left=of sec_dummy] (a) {};
	 \node[below=of sec_dummy](sec_market) {};
	 \foreach \x in {-11.45,-11.1,-10.75,...,-7.6}{
	 \foreach \y in {-0.45,-0.1,0.25,...,3.4}{
	 		\draw[fill=white] (-0.045+\x,\y) circle [radius=0.11, color=white];
	 	}
	}
	\foreach \x in {1.55,1.9,2.25,...,5.4}{
	\foreach \y in {0.1,0.45,...,3.8}{
			\draw[fill=white] (-0.045+\x+4* 1.7 * \y / 3.85 *\x/3.85-4* 1.7 * \y / 3.85 * 1.55/3.85 -4* 1.7 * \y / 3.85 * \x^2/3.85^2 +4* 1.7 * \y / 3.85 * 2* 1.55 * \x / 3.85^2 -4* 1.7 * \y / 3.85 * 1.55^2/3.85^2 - 4* 1.7 * \y^2 / 3.85^2 *\x/3.85 + 4* 1.7 * \y^2 / 3.85^2 *1.55/3.85 + 4* 1.7 * \y^2 / 3.85^2 *\x^2/3.85^2 -4* 1.7 * \y^2 / 3.85^2 * 2* 1.55 * \x / 3.85^2 + 4* 1.7 * \y^2 / 3.85^2 *1.55^2/3.85^2 ,\y - 0.6)  circle [radius=0.11, color=white];
		}
	}
	
		\foreach \x in {-5.}{
		\foreach \y in {-0.45,-0.1,0.25,...,3.4}{
			\draw[fill=white] (-0.02+\x,\y) circle [radius=0.11, color=white];
		}
	}
	
			\foreach \x in {-4.65}{
	\foreach \y in {-0.45,-0.1,...,2.}{
		\draw[fill=white] (-0.02+\x,\y) circle [radius=0.11, color=white];
	}
}

			\foreach \x in {-4.65}{
	\foreach \y in {2.7,3.05}{
		\draw[fill=white] (-0.02+\x,\y) circle [radius=0.11, color=white];
	}
}

	\foreach \x in {-4.3}{
	\foreach \y in {-0.45,-0.1}{
			\draw[fill=white] (-0.02+\x,\y) circle [radius=0.11, color=white];
		}
	}

	\foreach \x in {-4.3}{
	\foreach \y in {0.6,0.95,...,3.05}{
			\draw[fill=white] (-0.02+\x,\y) circle [radius=0.11, color=white];
		}
	}

	\foreach \x in {-3.95,-3.6,...,-1.4}{
	\foreach \y in {-0.45,-0.1,0.25,0.6,0.95}{
			\draw[fill=white] (-0.02+\x,\y) circle [radius=0.11, color=white];
		}
	}
	
	\foreach \x in {-3.95}{
	\foreach \y in {1.65,2.,...,3.4}{
			\draw[fill=white] (-0.02+\x,\y) circle [radius=0.11, color=white];
		}
	}
		
	\foreach \x in {-3.60}{
	\foreach \y in {1.3,1.65,2.,...,3.4}{
			\draw[fill=white] (-0.02+\x,\y) circle [radius=0.11, color=white];
		}
	}
		
	\foreach \x in {-3.25}{
	\foreach \y in {1.3,1.65,2.,2.7,3.05}{
			\draw[fill=white] (-0.02+\x,\y) circle [radius=0.11, color=white];
		}
	}
	
	\foreach \x in {-2.9}{
	\foreach \y in {1.3,1.65,2.,...,3.4}{
			\draw[fill=white] (-0.02+\x,\y) circle [radius=0.11, color=white];
		}
	}

	\foreach \x in {-2.55}{
	\foreach \y in {1.65,2.,2.35,2.7,3.05}{
			\draw[fill=white] (-0.02+\x,\y) circle [radius=0.11, color=white];
		}
	}
	
	\foreach \x in {-2.2}{
	\foreach \y in {1.3,1.65,2.,2.35,3.05}{
			\draw[fill=white] (-0.02+\x,\y) circle [radius=0.11, color=white];
		}
	}
			
	\foreach \x in {-1.85}{
	\foreach \y in {1.3,1.65,2.,...,3.4}{
			\draw[fill=white] (-0.02+\x,\y) circle [radius=0.11, color=white];
		}
	}
		
	\foreach \x in {-1.5}{
	\foreach \y in {1.3,1.65,2.,...,3.4}{
			\draw[fill=white] (-0.02+\x,\y) circle [radius=0.11, color=white];
		}
	}
			
	\end{tikzpicture}}
	\end{center}
	\caption[Illustration of reference material with defects and domain mappings.]{Illustration of $A_{\textnormal{ref}}(x)$ taking two values in the computational domain (left), random defects (center), and domain mapping of the reference (right).}
	\label{fig1}
\end{figure}
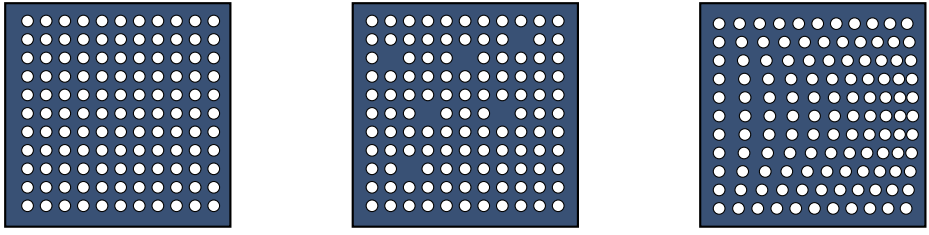

Let $\mu_\textnormal{ref} \in \Params$ be a fixed parameter for which the resulting coefficient $A_\textnormal{ref}$ can be considered a reference coefficient of the perturbed system.
In particular, different from the application of RB methods, the parameter space $\Params$ and the dependency of the data functions can not be stated explicitly to replicate an arbitrary set of perturbations, and, for this section, the standard RB-related \cref{asmpt:parameter_separable} is dropped.
In order to align with the notation of this thesis, we let $A_\mu$ refer to a perturbation sample of $A_\textnormal{ref}$.
In \cref{fig1}, we illustrate two simple examples of perturbations, where a perfect reference diffusion coefficient is subject to either local defects of the particle or to more smooth perturbations that a domain mapping can describe.

To motivate why we cannot simply replace the perturbed coefficient with the reference coefficient, we formulate an artificial problem based on the reference coefficient and right-hand side: find $u_{\textnormal{ref}}\in V$, such that for all $v \in V$,
\begin{equation}
\label{eq:reference_problem}
\int_\Omega A_\textnormal{ref}\nabla u_\textnormal{ref}\cdot\nabla v\,\dx =  \int_\Omega f_{\textnormal{ref}} \, v\dx.
\end{equation}
This equation is equal to \cref{def:param_elliptic_problem} for $\mu = \mu_\textnormal{ref}$, where we note that, different from the formulation in \eqref{eq:problem_classic}, we also consider a perturbed right-hand side to obtain a more general statement in \eqref{eq:apriori_error}.
For a sample $\mu \in \Params$, the error between $u_{\textnormal{ref}}$ and $u_\mu$ can then be bounded in the energy norm by

\begin{equation}
\label{eq:apriori_error}
\begin{split}
\anorm{ u_\mu - u_\textnormal{ref}}^2 &\leq ( A_\mu\nabla  u_\mu- A_\textnormal{ref}\nabla u_\textnormal{ref},\nabla (u_\mu-u_\textnormal{ref}))+( A_\textnormal{ref}\nabla  u_\textnormal{ref}-  A_\mu \nabla u_\textnormal{ref},\nabla ( u_\mu- u_\textnormal{ref})) \\
&= (f_\mu - f_{\textnormal{ref}},  u_\mu- u_\textnormal{ref}) + ((A_\textnormal{ref}- A_\mu)\nabla u_\textnormal{ref},\nabla(u_\mu-u_\textnormal{ref}))\\
&\leq \left( \frac{C_{\textnormal{P}}}{\alpha^{1/2}} \norm{f_\mu - f_{\textnormal{ref}}}_{L^2(\Omega)} + \frac{C_{\textnormal{P}}}{\alpha^{3/2}}\| A_\textnormal{ref}- A_\mu\|_{L^\infty(\Omega)}\|f_{\textnormal{ref}}\|_{L^2(\Omega)} \right) \anorm{ u_\mu - u_\textnormal{ref}},
\end{split}
\end{equation}
where $C_\textnormal{P}$ is the Poincar\'{e} constant for $\Omega$ and $\alpha$ denotes the ellipticity bound from \eqref{eq:ellipticity_1}.
This error bound suggests that even local perturbations in the coefficient structure or the right-hand side, e.g., by a defect or shift, may lead to abysmal accuracy.
This occurs, for example, if we consider a problem with a highly conductive thin channel in the diffusion coefficient and a right-hand side $f$ which has support inside the channel.
If the channel is moved slightly so that the support of the right-hand side is now outside the channel, the solution will behave very differently.
The error concerning perturbations in $f$ is less severe since it is measured in the $L^2$-norm.

Recall that it is unclear how computations for the standard finite element method can be reused reliably.
Instead, the LOD enables such a strategy since the corrector problems are only defined on local subspaces.

\subsubsection{Perturbations}
\label{sec:perturbations}
To simplify the presentation, we consider perturbations of coefficients that only take the two values $1$ and $0 < \alpha < 1$.
We emphasize that this assumption is unnecessary for the proposed method to work or for the theory to hold.
However, it highlights the application to composite materials, which has inspired this work and simplifies the presentation.
Let $\Omega_{1}, \Omega_{\alpha} \subseteq \Omega$ be two disjoint subdomains of $\Omega$ with $\Omega_{1} \cup \Omega_{\alpha} = \Omega$.
Let $A_{\textnormal{ref}}$ be defined by
\begin{equation}
 \label{assumption_on_A_ref}
A_{\textnormal{ref}} = \rchi_{\Omega_{1}} + \alpha \rchi_{\Omega_{\alpha}},  
\end{equation}
where $\rchi$ is the indicator function.

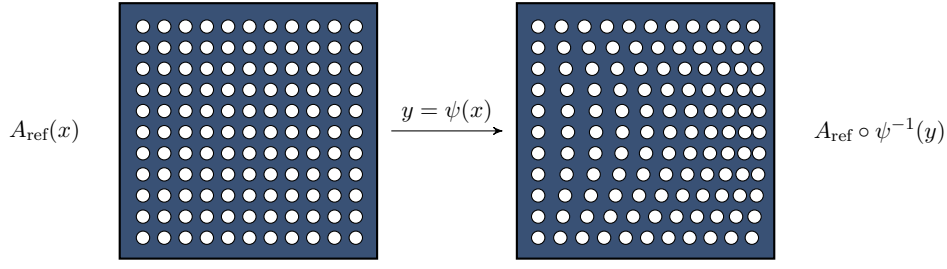
\begin{figure}
	\begin{center}
		\begin{tikzpicture}[scale=0.8, every node/.style={transform shape}, node distance=1cm]
		\node[] (market) {};
		\node[above=of market] (dummy) {};
		\node[blue-rec, right=of dummy] (t) {};
		\node[blue-rec, left=of dummy] (g) {}
		edge[pil,->, bend left=0] node[auto] {$y = \psi(x)$} (t);
		\node [left= .5cm of g]{$A_\textnormal{ref}(x)$};
		\node [right= .5cm of t]{$A_\textnormal{ref} \circ \psi^{-1}(y)$ };
		\foreach \x in {-11.45,-11.1,-10.75,...,-7.6}{
			\foreach \y in {-0.45,-0.1,0.25,...,3.4}{
				\draw[fill=white] (\x+6.45,\y-0.05) circle [radius=0.11, color=white];
			}
		}
		\foreach \x in {1.55,1.9,2.25,...,5.4}{
			\foreach \y in {0.1,0.45,...,3.8}{
				\draw[fill=white] (\x-0.05+4* 1.7 * \y / 3.85 *\x/3.85-4* 1.7 * \y / 3.85 * 1.55/3.85 -4* 1.7 * \y / 3.85 * \x^2/3.85^2 +4* 1.7 * \y / 3.85 * 2* 1.55 * \x / 3.85^2 -4* 1.7 * \y / 3.85 * 1.55^2/3.85^2 - 4* 1.7 * \y^2 / 3.85^2 *\x/3.85 + 4* 1.7 * \y^2 / 3.85^2 *1.55/3.85 + 4* 1.7 * \y^2 / 3.85^2 *\x^2/3.85^2 -4* 1.7 * \y^2 / 3.85^2 * 2* 1.55 * \x / 3.85^2 + 4* 1.7 * \y^2 / 3.85^2 *1.55^2/3.85^2 ,\y - 0.6)  circle [radius=0.11, color=white];
			}
		}
		\end{tikzpicture}
		\caption{Illustration of a domain mapping $\psi$ with $D \equiv 0$.}
		\label{neu}
	\end{center}
\end{figure}

We formalize the two types of perturbations that we consider (see Figure \ref{fig1}) by introducing a defect perturbation $D$ and a domain mapping $\psi$.
A perturbation from a defect can be expressed by $D = (1 - \alpha) \rchi_{\omega}$ where $\omega \subseteq \Omega_{1}$. Then, $A_\mu = A_{\textnormal{ref}} - D$ can be considered the perturbed coefficient.
For shift perturbations, we assume that the domain mapping perturbation can be described as a variable transformation with a perturbation function $\psi:\Omega\rightarrow\Omega$ which maps the reference coefficient (expressed in $x$-coordinates) to a mapped coefficient (expressed in $y$-coordinates).
We assume that $\psi$ maps the boundary to itself (i.e.\ $\Gamma = \{\psi(x)\,:\,x\in\Gamma\}$) and that it is a one-to-one mapping in $\Omega$. Figure \ref{neu} provides an example of a variable transformation.
We denote the corresponding Jacobi matrix
\begin{equation*}
(\Jac)_{ij}(x)=\left[\frac{\partial \psi_i}{\partial x_j}(x)\right],
\end{equation*}
for $1\leq i,j\leq d$, and assume it to be bounded with bounded inverse for a.e.\ $x\in \Omega$.
The two perturbation types can either be combined, or be considered individually by letting $D \equiv 0$ or $\psi = \id$.
Note that, apart from the formalized perturbations, also arbitrary perturbations can be handled.
However, the approach only shows its strength if these perturbations occur locally.

Using $A_{\textnormal{ref}}$, $D$ and $\psi$, we formulate the \emph{mapped problem} in the $y$-variable with $y = \psi(x)$ by
\begin{equation} \label{eq:domain_mapping_problem}
\begin{split}
- \nabla_y \cdot A_y \nabla_y u_y & = f_y, \qquad \textnormal{ in }\Omega, \\
u_y & = g_y, \qquad \textnormal{ on } \Gamma,
\end{split}
\end{equation}
where the coefficient is defined by $A_y = (A_\textnormal{ref}-D) \circ \psi^{-1}$, and the derivatives have been distorted accordingly.
The $y$-variable corresponds to the physical spatial variable in a typical situation.
Depending on the physics being modeled, either the mapped right-hand side $f_y \in L^2(\Omega)$ or the perturbed $f_\mu$ (below) can be considered given.
It makes no difference for the development of the numerical method but may affect the choice of $f_{\textnormal{ref}}$. 
We also note that the boundary $\Gamma$ is mapped to itself.
The solution in the perturbed domain is denoted $u_y(y)$.

Next, we use the mapped problem to define the perturbed problem that serves as the parameterized problem in equation \eqref{eq:param_elliptic_problem}.
The gradient operator $\nabla_y$ in the mapped domain can be expressed in terms of $\nabla_x = \nabla$ by
\begin{equation}
\nabla_y v(x) := \left[ \diff{v}{y_i}(x)  \right]_i = \left[ \sum_{j}^{}\diff{v}{x_j}(x) \diff{x_j}{y_i}(x) \right]_i = \Jac^{-T}(x)  \nabla_x v(x).
\end{equation}
Based on the elliptic operator in equation \eqref{eq:domain_mapping_problem} we define the perturbed bilinear form for the mapped problem as 
\begin{equation*}
\begin{split}
a_\mu(v,w) &= \int_{{\Omega}}^{} \left( (A_\textnormal{ref}-D) \circ \psi^{-1} \right )  \nabla_y \left(v \circ \psi^{-1} \right) \cdot  \nabla_y \left(w \circ \psi^{-1} \right) \, \dd y\\
&= \int_{\Omega}^{} \det(\Jac) \Jac^{-1} (A_\textnormal{ref}-D) \Jac^{-T}  \nabla_x v \cdot  \nabla_x w \, \dd x
\end{split}
\end{equation*}
and the corresponding linear functional
\begin{equation*}
l_\mu(w) = \int_{ \Omega} f_y \left(w \circ \psi^{-1}\right) \,\dd y = \int_{\Omega} \det(\Jac) \left(f_y\circ \psi\right)  w\, \dd x.
\end{equation*}
We see that this now fits the formulation of the parameterized (perturbed) problem \eqref{eq:param_elliptic_problem} with
\begin{equation*}
A_\mu=\det(\Jac) \Jac^{-1} (A_\textnormal{ref}-D) \Jac^{-T}, \qquad f_\mu = \det(\Jac) \left(f_y\circ \psi\right),
\end{equation*}
and the mapped solution $u_y = u_\mu \circ \psi^{-1}$.
For the problem to be well-posed, we assume $\Jac$ and $\Jac^{-1}$ to be bounded almost everywhere and $A_\mu$ to be symmetric positive definite.
We note that the perturbed coefficient $A_\mu$ can be computed from the reference coefficient using the Jacobian matrix.
We emphasize that the domain mapping transforms a shift defect into a change-in-value perturbation.
For many coefficients, this is advantageous as seen in equation \eqref{eq:apriori_error}, where now the $L^{\infty}$-norm can be expressed entirely in terms of how much $\Jac$ differs from the identity.
The domain mapping covers continuous (possibly global) perturbations, while defects cover discontinuous (often local) perturbations.
With respect to Figure \ref{fig1} we see that the middle picture corresponds to $\psi = \id$ and the right picture to $D \equiv 0$. 

\subsubsection{Error indicator and adaptive method}
The primary idea of the approach of the following adaptive procedure is that the local stiffness matrix contributions $\Ktmuref$ defined in \eqref{eq:K} are readily available from the reference coefficient $A_{\textnormal{ref}}$ and can partly be reused for a new sample $\mu \in \Params$.
Concerning \cref{def:comp_proc_LOD}, this means that in Step 1, we aim at replacing $\Ktmu$ by $\Ktmuref$ when the resulting approximation error is below a specific tolerance.
In order to decide where to recompute the respective local stiffness matrix contribution $\Ktmu$, we the following error indicator:
\begin{definition}[Error indicator] \label{def:indicators}
	For each $T \in \mathcal{T}_H$, we define
	\begin{equation}  \label{indidef}
	\begin{split}
	E^2_{\QQ V_H,T, \mu} &:= 
	\max_{\substack{ w |_{T}, \,w \in V_H 
		}
	} \frac{\norm{(A_\mu-A_{\textnormal{ref}}) A_\mu^{-1/2}(\rchi_T  \nabla w -  \nabla {\QQ ^{\textnormal{ref}}_{\kl,T}}w)  }^2_{L^2(U_\kl(T))}}{\|A_\mu^{1/2}\nabla w\|^2_{L^2(T)}},\\
	\end{split}
	\end{equation}
	where $\rchi_T$ denotes the indicator function for an element $T \in \mathcal{T}_H$ and $\QQ ^{\textnormal{ref}}_{\kl,T}$ denotes the corrector for $\mu_\textnormal{ref}$.
\end{definition}

It is important to notice that $E^2_{\QQ V_H,T, \mu}$ can be computed quickly for a new sample $\mu$.
The other required quantities can be assembled once $\QQ ^{\textnormal{ref}}_{\kl,T}$ are available.
Using the error indicators, we present an adaptive method that decides where to use $A_\mu$ and where to use $A_{\textnormal{ref}}$ for the local stiffness matrix contribution at $T$.
\begin{definition}[PG--LOD with adaptively updated correctors] \label{def:adaptive_LOD}
	The proposed method follows five steps, where the first step can be considered a preparatory step:
	\begin{enumerate}
		\item
		Follow Step 1 as explained in \cref{def:comp_proc_LOD} for $A_\textnormal{ref}$, i.e. compute (for all $T\in\mathcal{T}_H$) reference correctors $\QQ^{\textnormal{ref}}_{\kl,T}$ and compute the corresponding stiffness matrix contribution $\Ktmuref$.
	\end{enumerate}
	For every new perturbation $\mu \in \Params$:
	\begin{enumerate} \setcounter{enumi}{1}
		\item Compute (for all $T\in\mathcal{T}_H$) the error indicator $E_{\QQ V_H,T,\mu}$ and mark the elements $T$ for which the following inequality holds:
		$$
		E_{\QQ V_H,T,\mu} \le \textnormal{TOL}.
		$$
		Denote the set of marked elements by $\mathcal{T}^{\textnormal{ref}}_H\subset \mathcal{T}_H$.
		\item
		Compute (for all $T \in\mathcal{T}_H\setminus \mathcal{T}^{\textnormal{ref}}_H$) the mixed correctors $\tilde{\QQ }^T_{\kl,\mu}$, based on the following definitions of the mixed right-hand side and correctors:
		\begin{equation*}
		\tilde{\QQ }_{\kl,\mu}^T = \left\{
		\begin{array}{ll}
		{\QQ }^{\textnormal{ref}}_{\kl,T}, \\
		\QQktm.
		\end{array}\right.
		\end{equation*}
		Further, let $\tilde{\QQ }_{\kl, \mu}=\sum_{T\in\mathcal{T}_H}\tilde{\QQ }^T_{\kl,\mu}$.
		\item Assemble the adaptively updated mixed LOD stiffness matrix by
		\begin{equation*} \label{eq:K_tilde}
		\tilde{\K}_\mu = \sum_{ T \in \mathcal{T}^{\textnormal{ref}}_H}^{} \Ktmuref + \sum_{T \in {\mathcal{T}}_H \setminus \mathcal{T}^{\textnormal{ref}}_H}^{} \Ktmu.
		\end{equation*}
		\item Solve for the component vector of $\tilde u_{H,\kl,\mu} \in  V_H$ in
		\begin{equation}\label{proposedmethod}
		\begin{aligned}
		\tilde{\K}_\mu \underline{\tilde u}_{H,\kl,\mu} = \mathbb{F}
		\end{aligned}
		\end{equation}
		and compute the solution in the mixed LOD space as
		\begin{equation*}\label{eq:approximate_tildeuk}
		\tilde u^{\textnormal{ms}}_{H,\kl,\mu}= \tilde{u}_{H,\kl,\mu} -\tilde{\QQ }_{\kl,\mu}^T\tilde{u}_{H,\kl,\mu}.
		\end{equation*}
	\end{enumerate}
\end{definition}

\begin{figure}[t] \centering \footnotesize
	\begin{tikzpicture}[
		node distance = 6mm and 12mm,
		start chain = A going below,
		base-red/.style = {draw, minimum width=10mm, minimum height=10mm,
			align=center, on chain=A, rectangle, fill=red!5},
		base-blue/.style = {draw, minimum width=10mm, minimum height=10mm,
			align=center, on chain=A, rectangle, fill=blue!5},
		every edge quotes/.style = {auto=right}]
		\node [base-blue, minimum width=20mm, minimum height=20mm, fill=black!10] {mixed coarse LOD system \\ \eqref{proposedmethod} on \textcolor{blue}{$\mathcal{T}_H$} };
		\node [base-blue, above=1cm of A-1, xshift=-2cm] {\footnotesize{$T_0$}\\$E_{\QQ V_H,T_0,\mu} \le \textnormal{TOL}$ \\ reuse reference};
		\node [base-blue, above=1cm of A-1, xshift=1.5cm] {\footnotesize{$T_1$}\\$E_{\QQ V_H,T_1,\mu} \le \textnormal{TOL}$ \\ reuse reference};
		\node [base-red, above=0.1cm of A-1, xshift=5cm, yshift=-0.8cm] {\footnotesize{$T_2$}\\$E_{\QQ V_H,T_2,\mu} > \textnormal{TOL}$ \\ recompute corrector};
		\node [base-blue, above=-1.6cm of A-1, xshift=5cm, yshift=-0.8cm] {\footnotesize{$T_3$}\\$E_{\QQ V_H,T_3,\mu} \le \textnormal{TOL}$ \\ reuse reference};
		\node [base-red, above=-3.6cm of A-1, xshift=2.5cm, yshift=-0.5cm] {\footnotesize{$T_4$}\\$E_{\QQ V_H,T_4,\mu} > \textnormal{TOL}$ \\ recompute corrector};
		\draw [->] (A-2) -- (A-1) node[midway, left] {$\K_{T_0,\mu_\textnormal{ref}}$};
		\draw [->] (A-3) -- (A-1) node[midway, left] {$\K_{T_1,\mu_\textnormal{ref}}$};
		\draw [->] (A-4) -- (A-1) node[midway, above] {$\K_{T_2,\mu}$};
		\draw [->] (A-5) -- (A-1) node[midway, above] {$\K_{T_3,\mu_\textnormal{ref}}$};
		\draw [->] (A-6) -- (A-1) node[midway, right] {$\K_{T_4,\mu}$};
		
		\draw [->] (-0.7,-1.8) -- (A-1);
		\fill[fill=black] (-0.5,-2.2) circle [radius=0.08];
		\fill[fill=black] (-0.7,-2.2) circle [radius=0.08];
		\fill[fill=black] (-0.9,-2.2) circle [radius=0.08];				
	\end{tikzpicture}
	\caption[Illustration of the adaptive PG--LOD procedure for perturbations.]{Illustration of the adaptive PG--LOD procedure for perturbations.
		Compared to \cref{fig:FOM_LOD}, the corrector problems of $A_\textnormal{ref}$ may be reused.}
	\label{fig:Hekema}
\end{figure}
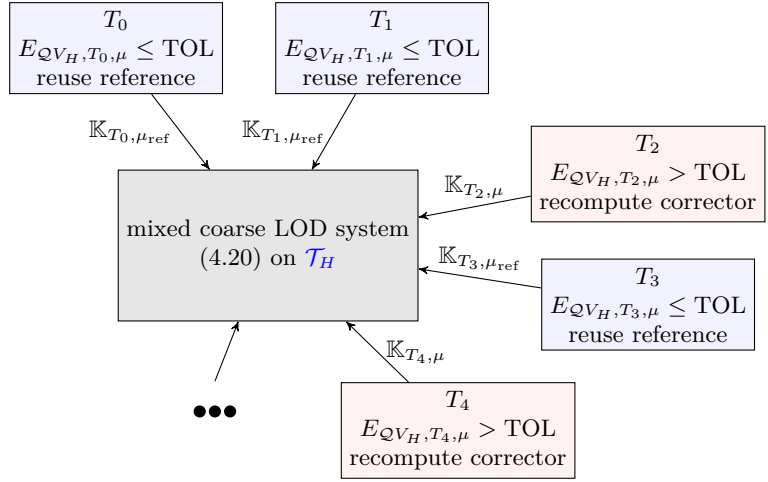

We note that, different from the source \cite{HKM20}, the algorithm in \cref{def:adaptive_LOD} has been simplified to only treat the standard correctors and to be aligned with the notation of parameterized problems.
Instead, in \cite{HKM20}, right-hand side corrections and non-zero boundary values are considered, which require introducing two additional error indicators.
However, for simplicity, we omit a detailed description of these features; cf.~\cref{rem:rhs_correction}.
In \cref{fig:Hekema} we illustrate the basic procedure of the adaptive approach compared to the classical PG--LOD, visualized in \cref{fig:FOM_LOD}.
While the computational time for Step~1 is reduced, Step~2 and Step~3 computationally remain the same.

The corresponding approximation result can be stated as follows. 

\begin{theorem}[Error bound for the PG--LOD with adaptively updated correctors] \label{theorem}
	If $$
	\max_{T\in \mathcal{T}^{\textnormal{ref}}_H}(E_{\QQ V_H,T})\leq \textnormal{TOL},
	$$
	there exist $\kl_0 > 0 $ and $\tau_0 > 0$, such that for all $\kl > \kl_0$ and $0 < \tau < \tau_0$, with $\textnormal{TOL} = \tau \kl^{-d/2}$, the error bound
	\begin{equation*} \label{eq:theorem_4.2}
	\anorm{u - \tilde u^{\textnormal{ms}}_{H,\kl,\mu}} \lesssim (H+\kl^{d/2} ( \theta^\kl  + \textnormal{TOL})) \norm{f}_{L^2}
	\end{equation*}
	is satisfied. Here $0 < \theta < 1$ is independent of $H$, $\kl$, $\tau$ and TOL.
\end{theorem}
\noindent For a proof of this result, we refer to \cite{HM19}.
It means that by guaranteeing that the error indicators are less than TOL we obtain an approximate solution, which is arbitrarily close to the PG--LOD approximation.
If the perturbations are local in space, the reference coefficient can be used in a large part of the computational domain. 
We again mention that \cref{theorem} has been generalized for right-hand side corrections and boundary value problems in \cite{HKM20}, where the difference is, that the $H$ dependency is diminished thanks to the more accurate approximation strategy of the right-hand side correction.
With respect to the fact that storage requirements may make it infeasible to store correctors, we mention that the coarse-space representation $\tilde{u}_{H,\kl,\mu} \in V_H$ can be considered sufficient in many applications; cf.~\cite{elf}.
In this case, Step~3 and the computation of $\tilde u^{\textnormal{ms}}_{H,\kl,\mu}$ may be left out.
We finish this section by demonstrating defects and domain mappings with two numerical examples.

\subsubsection{Experiment 9: Local defects and local domain mappings}
\label{sec:hekema_experiments}

\label{clean up this section to the fundamentals !!}
The following experiments are a small selection of \cite{HKM20} and have been performed using the \texttt{Python} PG--LOD implementation \texttt{gridlod}~\cite{gridlod}.
The source code is referenced in \cref{apx:code_availability}, where also details on how the perturbed diffusion coefficients for this section have been constructed can be found.
More details on \texttt{gridlod} are furthermore discussed in \cref{sec:gridlod}.

Our experiments are performed on a $2$d-quadrilateral mesh on $\Omega = [0,1]^2$. We let $H=2^{-5}$, $h=2^{-8}$, $\kl=4$, and use a reference coefficient that is piece-wise constant on every fine mesh element.
In the first experiment, $A_{{\textnormal{ref}}}$ can be expressed as stated in \eqref{assumption_on_A_ref}, i.e.\ it takes two values $1$ and $\alpha$ and only defect-perturbations are present, whereas, in the second experiment, perturbations that can be described by domain mappings are demonstrated.
For the sake of convenience, we neglect an explicit definition of $A_{\textnormal{ref}}$ and $A_\mu$ as we visualize them in the figures.

\begin{figure}[h]
	\begin{subfigure}[b]{0.3\textwidth}
		\includegraphics[width=\textwidth]{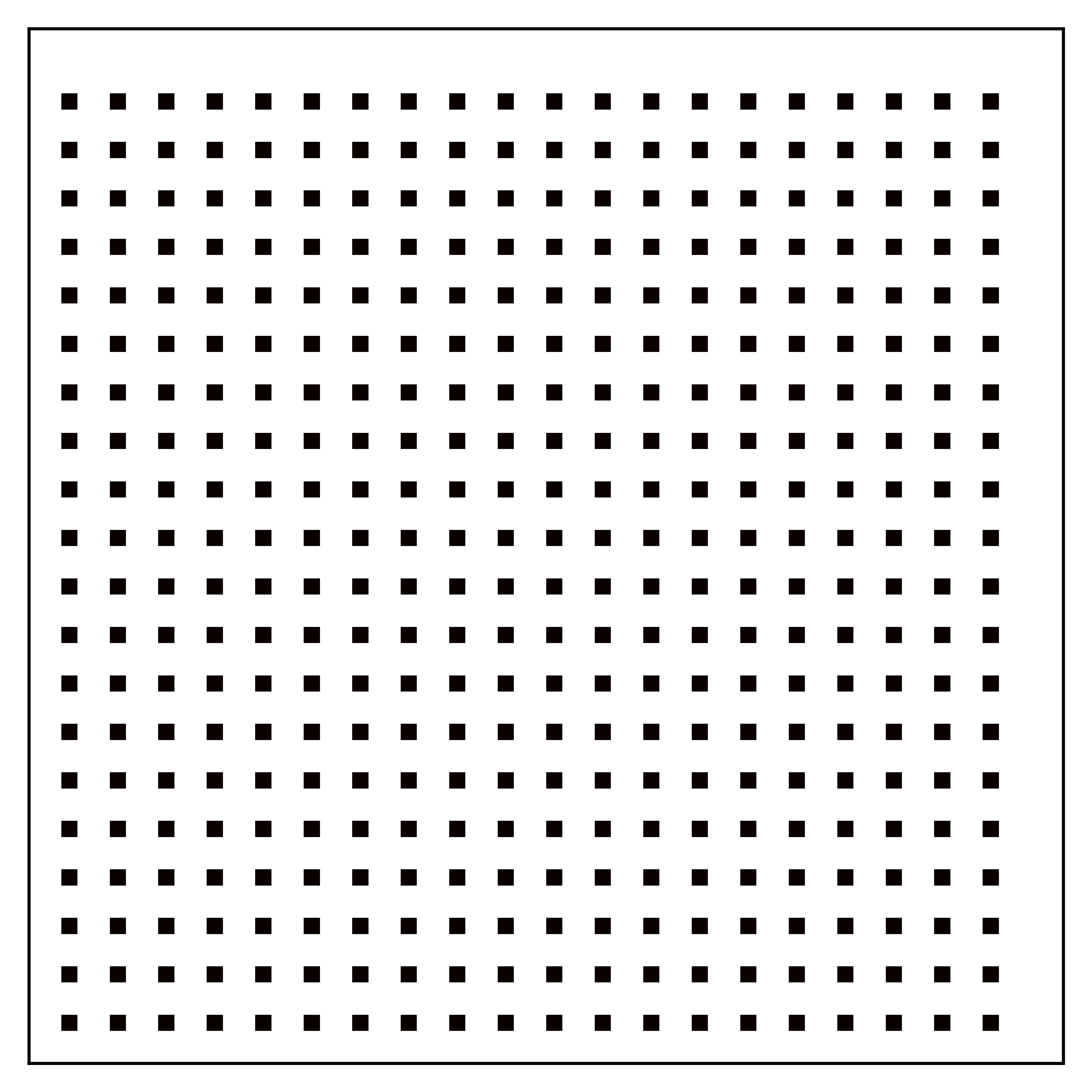}
	\end{subfigure}
	\begin{subfigure}[b]{0.3\textwidth}
		\includegraphics[width=\textwidth]{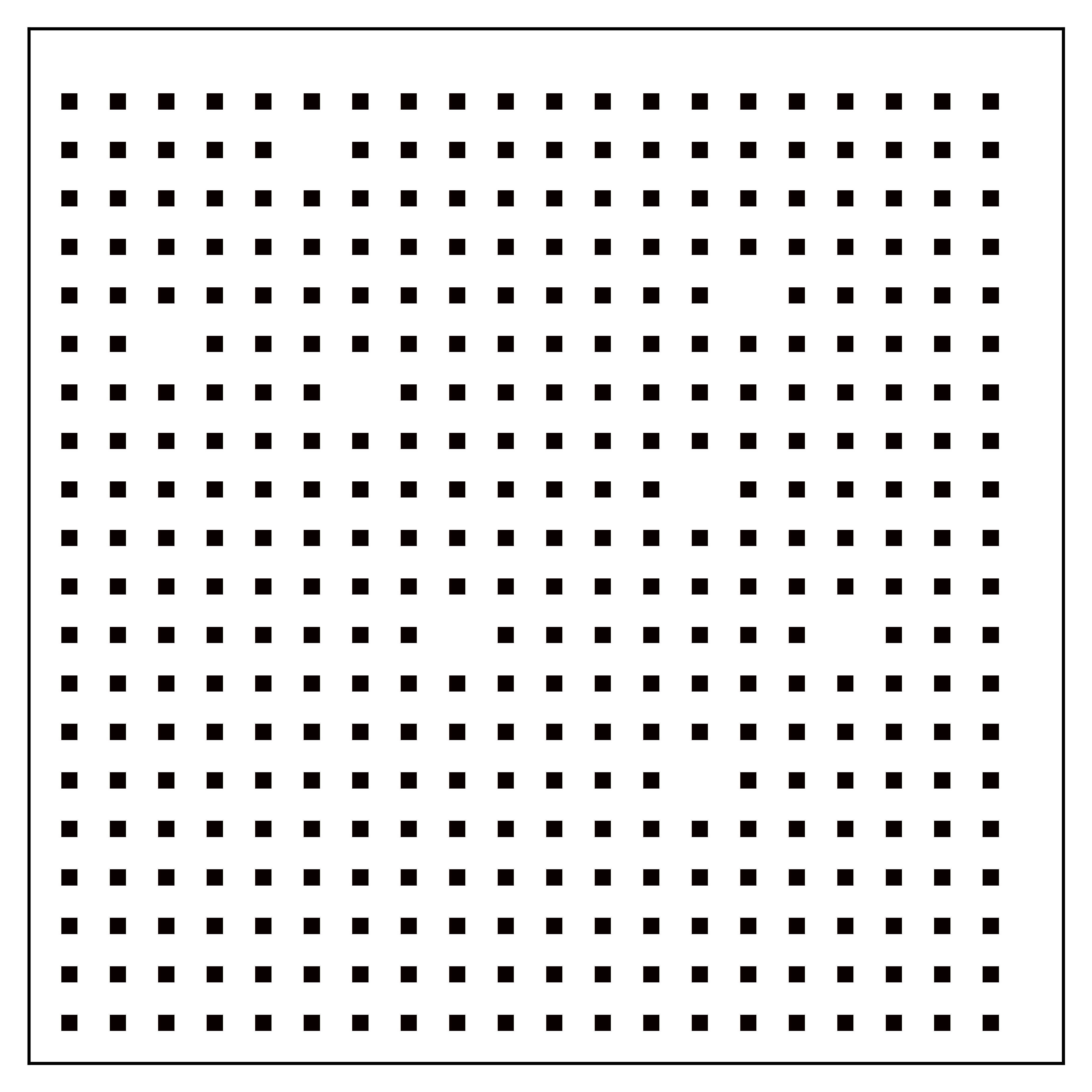}
	\end{subfigure}
	\begin{subfigure}[b]{0.37\textwidth}
	\includegraphics[width=\textwidth]{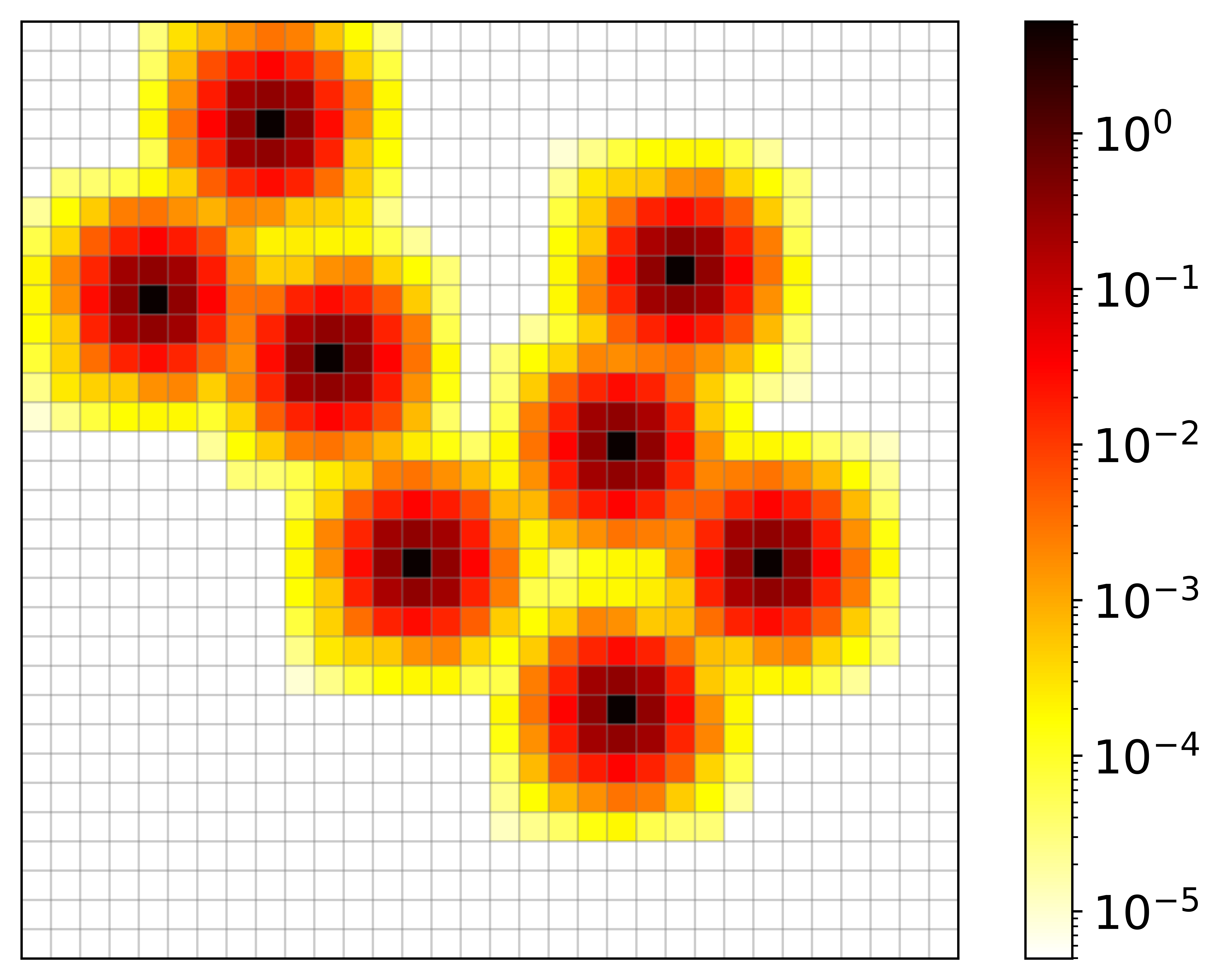}
	\end{subfigure}
	\centering
	\caption[Experiment 9: Reference coefficient and perturbations]{Reference coefficient $A_{\textnormal{ref}}$ (left) and defect perturbation $A$ (middle). Black is $1$, white is $0.1$. Furthermore, $E_{\QQ V_H,T}$ is shown (right).}
	\label{ex1:coefficients}
\end{figure}

In the experiments, a defect in the material means that a particle is equalized to the background (compare Figure \ref{ex1:coefficients}).
These defects occur with a probability of $2\%$. The right-hand side $f_y$ is defined by $f_y(y) = \rchi_{[1/8,7/8]^2}(y)$. 

We consider the relative error
$$
\mathcal{E}_{\textnormal{rel}}(\tilde u^{\textnormal{ms}}_{H,\kl,\mu},u^{\textnormal{ms}}_{H,\kl,\mu}) = \frac{\anorm{\tilde u^{\textnormal{ms}}_{H,\kl,\mu} - u^{\textnormal{ms}}_{H,\kl,\mu}}}{\anorm{\tilde u^{\textnormal{ms}}_{H,\kl,\mu}}},
$$ 
where $u^{\textnormal{ms}}_{H,\kl,\mu}$ is the best PG--LOD solution and $\tilde u^{\textnormal{ms}}_{H,\kl,\mu}$ is the solution of \eqref{eq:approximate_tildeuk} for a specific of the tolerance TOL in the algorithm of \cref{def:adaptive_LOD}.
We note that, for this experiment, the right-hand side correction and corresponding correctors and indicators as proposed in \cite{HKM20} are used, which has the only effect that the PG--LOD error is generally smaller and that the adaptive method requires additional error indicators, which has not been detailed in this section; cf.~\cref{rem:rhs_correction}.
For $\textnormal{TOL}=\infty$, we clearly have $0\%$ updates of the correctors whereas $\textnormal{TOL}= 0$ corresponds to $100\%$ updates. 
For $100\%$ updates, we then end up with the standard Petrov--Galerkin LOD error dependent on our data and discretizations. 
In order to observe the complete behavior of $\mathcal{E}_{\textnormal{rel}}(u_\kl,\tilde{u}_\kl)$, 
we compute $\mathcal{E}_{\textnormal{rel}}(u_\kl,\tilde{u}_\kl)$ for every possible choice of TOL (and thus for every percentage of updates).
The relative best PG--LOD error $\mathcal{E}_{\textnormal{rel}}(u_{h,\mu},\tilde u^{\textnormal{ms}}_{H,\kl,\mu})$ is always around $10^{-3}$ which means that we are comparing to a sufficiently accurate solution.

\begin{figure}[h]
	\centering
	\includegraphics[width=0.6\textwidth]{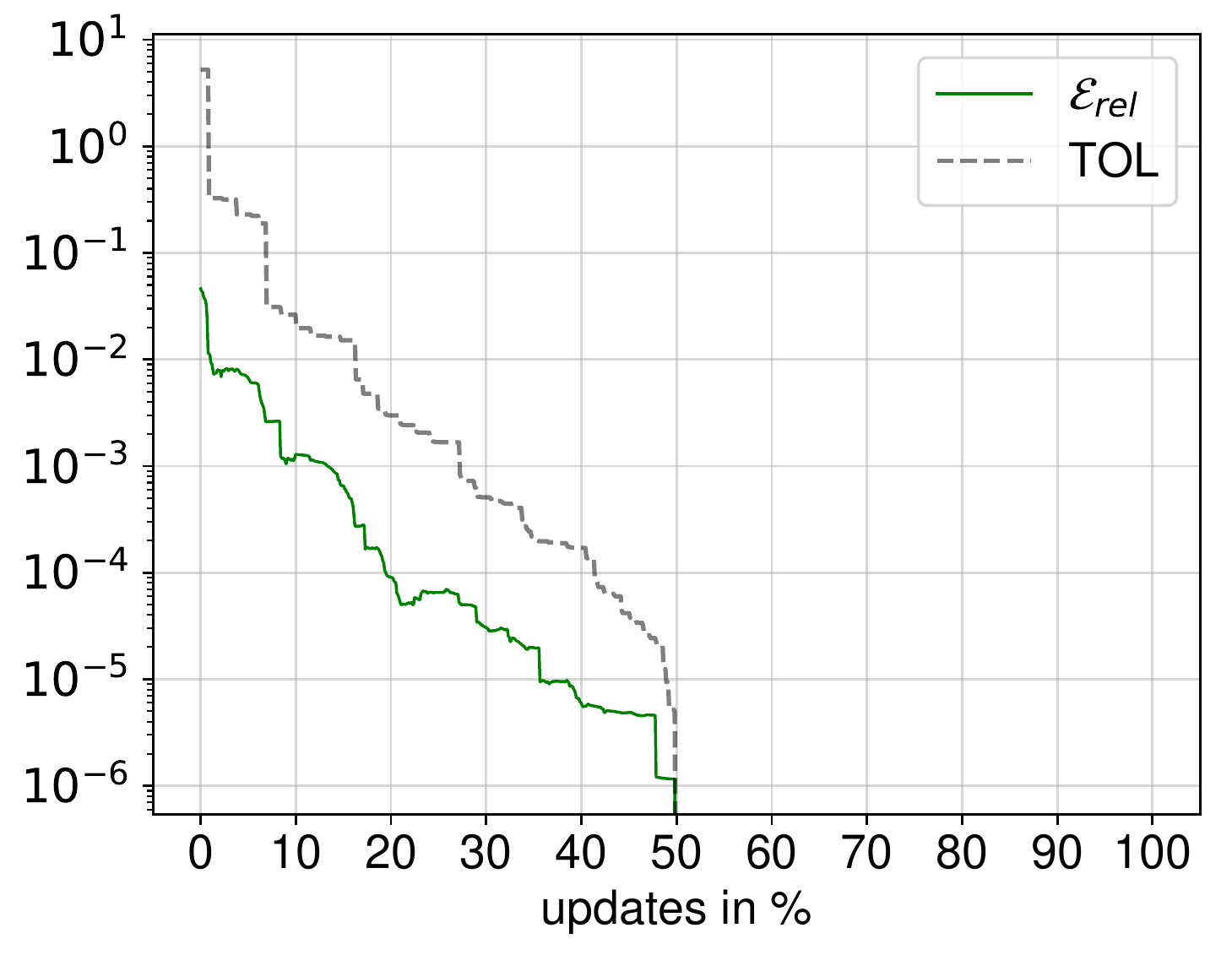}
	\caption{Experiment 9: Relative error improvement for defects.}
	\label{ex1:errors}
\end{figure}

\paragraph{Local defects}
In the first experiment, we let $\psi=\id$, which means we only consider defects, and for the background $\Omega_{\alpha}$ we choose $\alpha=0.1$.
Figure \ref{ex1:coefficients} displays the coefficient and its perturbation. Also, the error indicator $E_{\QQ V_H,T}$ is plotted for each $T$, where the coarse mesh is visible in the background.
We can see that the error indicator $E_{\QQ V_H,T}$ detects the defects in the coefficient correctly. 
Furthermore, $E_{\QQ V_H,T}$ is exponentially decaying away from each defect. 
From $\mathcal{E}_{\textnormal{rel}}(u_\kl,\tilde{u}_\kl)$ in \cref{ex1:errors} we see a big improvement for few updates of the correctors and a sufficiently fast convergence to the best PG--LOD solution.
Notably, after the correctors that attain the most significant error (black elements in \cref{ex1:coefficients}(right)) have been updated, the error in \cref{ex1:errors} already gains multiple magnitudes.
Furthermore, since only approximately 50\% of the correctors are affected by the perturbations, not all correctors have to be updated to achieve the possible PG--LOD error.
This experiment concludes that the method can efficiently be used for local defects.

\begin{figure}[h!]
	\begin{subfigure}[b]{0.3\textwidth}
		\includegraphics[width=\textwidth]{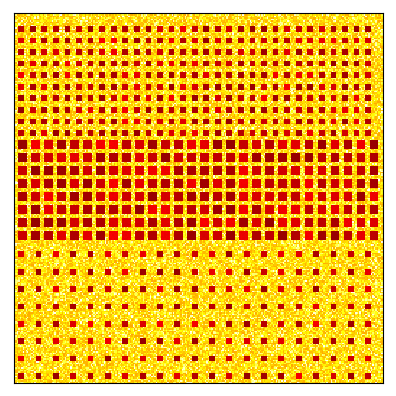}
	\end{subfigure}
	\begin{subfigure}[b]{0.3\textwidth}
		\includegraphics[width=\textwidth]{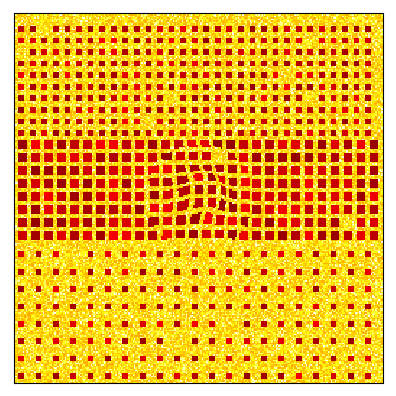}
	\end{subfigure}
	\begin{subfigure}[b]{0.375\textwidth}
		\hbox{
			\includegraphics[width=\textwidth]{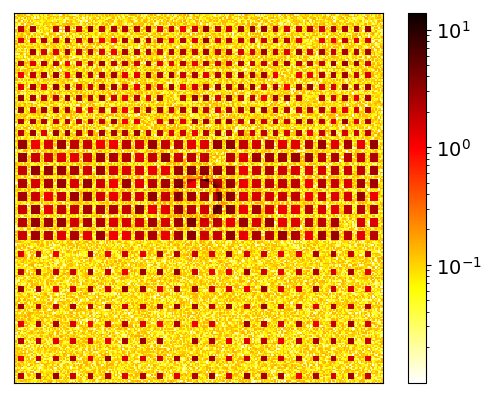}}
	\end{subfigure}
	\centering
	\caption[Experiment 9: Reference coefficient and perturbations with domain mappings]{Reference coefficient $A_{\textnormal{ref}}$ (left), perturbation in the physical domain $A_y$ (center) and corresponding change in value perturbation $A$ (right).}
	\label{ex2:coefficients}
\end{figure}

\paragraph{Local defects and local domain mappings}
In the second experiment, we add domain mappings and neglect the assumption in \eqref{assumption_on_A_ref} on $A_{\textnormal{ref}}$ to demonstrate the widened applicability of our adaptive method.
We choose $A_{\textnormal{ref}}$ to consist of differently sized and valued particles with values in $[1,5]$ (to replace $\Omega_1$) and a noisy background with values in $[10^{-2},5 \cdot 10^{-1}]$ (instead of one value $\alpha$). See Figure \ref{ex2:coefficients}(left) for a visualization of $A_{\textnormal{ref}}$.
Furthermore, we choose $\psi$ to be a local distortion in the middle of the domain, which can be seen in Figure \ref{ex2:coefficients}(center).
With the help of domain mappings, the reference coefficient $A_{\textnormal{ref}}$ is subjected to a simple change in value which means that $\Omega_1$ does not change its position. 
This is visualized in Figure \ref{ex2:coefficients}(right).

\begin{figure}[h!]
	\begin{minipage}{0.45\textwidth}
		\includegraphics[width=\textwidth]{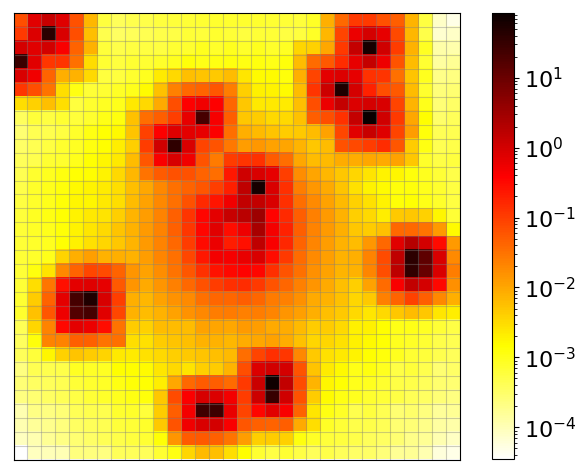}
	\end{minipage}\hfill
	\begin{minipage}{0.5\textwidth}
		\includegraphics[width=\textwidth]{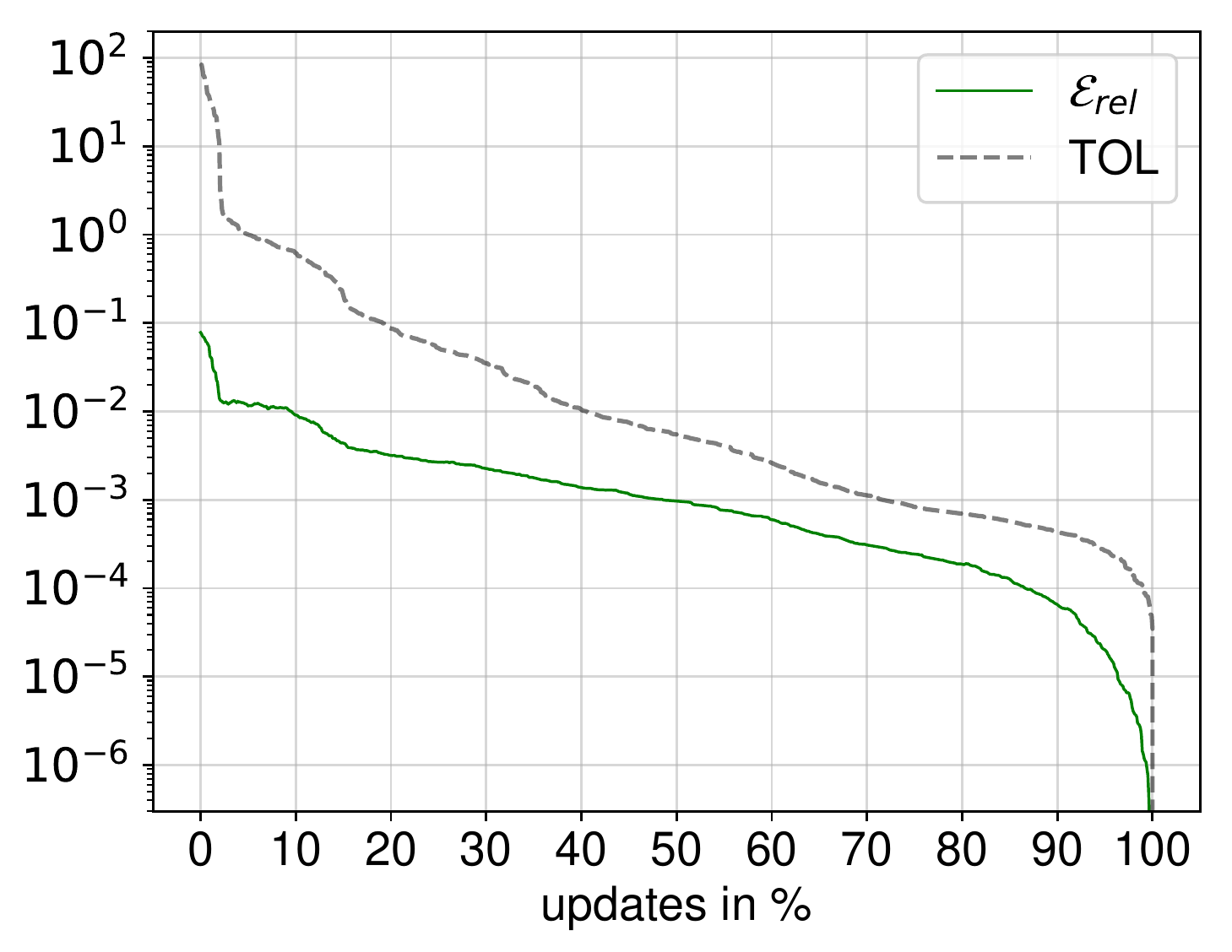}
	\end{minipage}
	\centering
	\caption[Experiment 9: Relative error improvement for domain mappings]{Error indicator $E_{\QQ V_H,T}$ (left) and relative error improvement for the local domain mapping (right).}
	\label{ex2:errors}
\end{figure}

As illustrated in \cref{ex2:errors}(left), the domain mapping, as well as the defects, can be seen in $E_{\QQ V_H,T}$, whereas the defects stick out compared to the domain mapping. 
Because of the relative error behavior in \cref{ex2:errors}(right), we observe a similar effect as in the first experiment. Again, the defects are treated first, giving a rapid error decay for a few updates.
However, due to the domain mapping, it takes comparatively longer to converge to the optimal PG--LOD solution. 

\subsubsection{Concluding remarks}

As discussed in \cref{sec:minimizing_cost}, the presented adaptive LOD approach has the weakness that, concerning computational effort, it only diminishes some of the computational requirements of Step 1.
On the other hand, the construction of the adaptive method is quite flexible because it does not assume a specific structure of the perturbations.
Note that \eqref{assumption_on_A_ref} was only required to simplify the theory and has been relaxed in the numerical experiments.
Standard assumptions like \cref{asmpt:parameter_separable} for RB methods are not required.
As particularly shown in the experiments in \cite{HKM20}, however, it can be expected that for large or global perturbations, the number of correctors that have to be recomputed is very high.
In the worst case, the adaptive method recomputes all correctors, which does not give any advantage for the many-query example.

\subsection{RB approaches for the LOD}
\label{sec:RB_for_LOD}

As we have already seen in this thesis, RB methods are exceptionally competitive in many-query and real-time scenarios for parameterized problems.
In contrast to the previously discussed adaptive LOD method~\cite{HKM20}, the idea of parameterized systems is not to consider only a perturbed system from a reference configuration but, instead, more general parameterized problems where \cref{asmpt:parameter_separable} can be exploited.
As this thesis is highly concerned with RB methods, in what follows, we aim at reviewing and developing RB approaches for the LOD.

\begin{figure}[t] \centering \footnotesize
	\begin{tikzpicture}[
	node distance = 6mm and 12mm,
	start chain = A going below,
	base/.style = {draw, minimum width=10mm, minimum height=10mm,
		align=center, on chain=A, rectangle, fill=green!5},
	every edge quotes/.style = {auto=right}]
	\node [base, minimum width=20mm, minimum height=20mm, fill=black!10] {Coarse LOD system \\ on \textcolor{blue}{$\mathcal{T}_H$} };
	\node [base, above=1cm of A-1, xshift=-2cm] {\footnotesize{$T_0$}\\corrector ROM \\ with size \textcolor{olive}{$N_\textnormal{rb}^{T_0}$}};
	\node [base, above=1cm of A-1, xshift=1.5cm] {\footnotesize{$T_1$}\\corrector ROM \\ with size \textcolor{olive}{$N_\textnormal{rb}^{T_0}$}};
	\node [base, above=0.1cm of A-1, xshift=4cm, yshift=-0.8cm] {\footnotesize{$T_2$}\\corrector ROM \\ with size \textcolor{olive}{$N_\textnormal{rb}^{T_0}$}};
	\node [base, above=-1.6cm of A-1, xshift=4cm, yshift=-0.8cm] {\footnotesize{$T_3$}\\corrector ROM \\ with size \textcolor{olive}{$N_\textnormal{rb}^{T_0}$}};
	\node [base, above=-3.6cm of A-1, xshift=2.5cm, yshift=-0.5cm] {\footnotesize{$T_4$}\\corrector ROM \\ with size \textcolor{olive}{$N_\textnormal{rb}^{T_0}$}};
	\draw [->] (A-2) -- (A-1) node[midway, left] {$\K^{\textnormal{rb}}_{T_0,\mu}$};
	\draw [->] (A-3) -- (A-1) node[midway, left] {$\K^{\textnormal{rb}}_{T_1,\mu}$};
	\draw [->] (A-4) -- (A-1) node[midway, above] {$\K^{\textnormal{rb}}_{T_2,\mu}$};
	\draw [->] (A-5) -- (A-1) node[midway, above] {$\K^{\textnormal{rb}}_{T_3,\mu}$};
	\draw [->] (A-6) -- (A-1) node[midway, right] {$\K^{\textnormal{rb}}_{T_4,\mu}$};
	
	\draw [->] (-0.7,-1.8) -- (A-1);
	\fill[fill=black] (-0.5,-2.2) circle [radius=0.08];
	\fill[fill=black] (-0.7,-2.2) circle [radius=0.08];
	\fill[fill=black] (-0.9,-2.2) circle [radius=0.08];				
	\end{tikzpicture}
	\caption[Illustration of the PG--RBLOD procedure.]{Illustration of the PG--RBLOD procedure. Compared to \cref{fig:FOM_LOD}, the corrector problems are solved with ROMs.}
	\label{fig:RBLOD}
\end{figure}
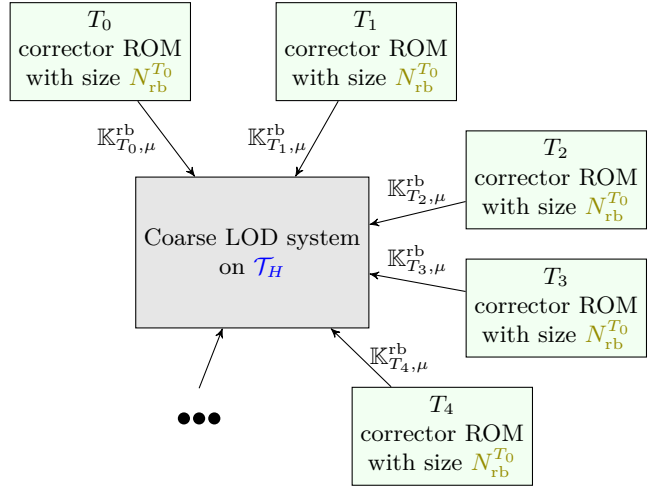

In~\cite{RBLOD}, RB approximations of the local corrector problems were introduced to obtain a reduced model for Step~1 that is independent of the size of $\grid$.
Since the corresponding RB systems are constructed offline, just as in \cite{MV21}, the approach gains significant improvements in terms of online time, clearly outperforming the time-to-solution for a single sample of the fine-scale data compared to \cite{HKM20,HM19}, if offline times are neglected.
The resulting approach is called the RBLOD and was developed based on the original Galerkin formulation of the LOD.
Since, as justified earlier, in this thesis, we are solely concerned about the Petrov--Galerkin variant of the LOD, we consider the methodology of the RBLOD in a Petrov--Galerkin variant.
The resulting PG-based RBLOD is visualized in \cref{fig:RBLOD}, where we mention that, due to the PG variant, the local matrices $\K^{rb}_{T,\mu}$ can be computed independently.
Concerning the overall computational cost, we note that for large $\Gridh$, the costs of solving the reduced corrector problems in Step~1 and the further computations in Steps~2 and~3 still can be high.
Thus, the online efficiency of the original RBLOD is bounded by the size of $\Gridh$.

Next, we classify the RBLOD concerning the abstract definition of a localized reduced approach with the direct sum of local spaces $V^i \subset V_h$ in \eqref{eq:direct_sum}.
In \cref{sec:localized_ms_space}, we already mentioned that the construction of the localized multiscale space $V_{H,\kl,\mu}^{\textnormal{ms}}$ constitutes the global approximation space for the PG--LOD; cf.~\eqref{eq:PG_}.
The space can also be interpreted by the space
\begin{equation*}
\Vts = V_H \oplus \Vfhkt[T_1] \oplus \cdots \oplus \Vfhkt[\Tlast],
\end{equation*}
in such a way that a function in $V_{H,\kl,\mu}^{\textnormal{ms}}$ has a unique representation in $\Vts$.
This formulation is particularly suitable for explaining the key idea of the RBLOD for parameterized problems.
As discussed in \cref{sec:localized_MOR}, reducing localized approaches for the sake of efficiently solving parameterized problems consists of finding reduced spaces~$V^i_\red$ for the subspaces~$V^i$ to construct~\eqref{eq:direct_sum_reduced}.
Since, in the LOD, only the local corrector spaces $\Vfhkt$ inherit parts of the fine-scale mesh $\mathcal{T}_h$, the reduction strategy of the RBLOD is to find reduced spaces $\Vfrbkt$ for these corrector spaces, i.e., to set $V^i_\red \coloneqq \Vfrbkt[T_i]$ for $i= 1, \dots, \abs{\mathcal{T}_H}$.
Importantly, the coarse space $V^0 = V_H$ is not changed, such that the RBLOD considers
\begin{equation*}
	 \Vtsrblod := V_H \oplus \Vfrbkt[T_1] \oplus \cdots \oplus \Vfrbkt[\Tlast] \subset \Vts,
\end{equation*}
as a subspace of $\Vts$.
Notably, the function space $\Vtsrblod$ might still have a high dimension and, instead of seeking the two-scale approximation in this space, meaning to solve all corrector problems and the coarse system at once, the RBLOD method is computationally much more feasible if the reduced version of Step~1 and the classical procedure of Steps~2 and~3 are followed.
We also refer to the discussion related to \cref{fig:TSRBLOD} for the impracticability of computing solutions directly on $\Vtsrblod$.
We conclude that the resulting coarse approximation procedure that we summarized in Steps 2 and 3 has not changed and only the local corrector problems can now be assessed quickly by their reduced approximations.
Certainly, the $H$ dependence is still present in the respective space.

In the remainder of this chapter, we are interested in the overall online efficiency of a reduced model for the PG--LOD, taking all computational steps into account and yielding a reduced-order model independent of the sizes of both $\grid$ and $\Gridh$.
In other words, with regard to the above-discussed function spaces, this means that we are aiming at finding a low dimensional subspace $\Vtsrb$ of $\Vtsrblod$ that makes it possible to reduce the entire two-scale space $\Vts$.
This results in a quasi-optimal efficient reduced approach for the PG--LOD.
The preliminaries of this approach are presented in the subsequent section.

\section{Two-scale formulation of the PG--LOD}
\label{sec:two_scale}
In the sequel, we formulate the PG--LOD method in a two-scale formulation.
In particular, we aim to consider the PG--LOD solution as the solution of one single system where the coarse system \eqref{eq:PG_} and all fine-scale corrections~\eqref{eq:loc_cor_problems} are solved at the same time.
This formulation will be the basis for the Stage~2 ROM, constructed in \cref{sec:stage_2_red}.

\subsection{The two-scale bilinear form}
Let, again, $\Vts$ denote the two-scale function space given by the direct sum of Hilbert spaces
\begin{equation*}
\Vts:= V_H \oplus \Vfhkt[T_1] \oplus \cdots \oplus \Vfhkt[\Tlast],
\end{equation*}
such that for $\uts = (u_H, \uft[T_1], \dots, \uft[\Tlast]) \in \Vts$ we define the two-scale $H^1$-norm of $\uts$ by
\begin{equation*}
\Snorm{\uts}^2 := \snorm{u_H}^2 + \sumT \snorm*{\uft}^2.
\end{equation*}
On this space, we define the two-scale bilinear form $\Btsm \in \textnormal{Bil}(\Vts)$ given by
\begin{multline*}
\Btsm\left((u_H, \ufti{1}, \dots, \ufti{{\abs{\mathcal{T}_H}}}),  (v_H, \vfti{1}, \dots, \vfti{{\abs{\mathcal{T}_H}}})\right):=
\\ a_{\mu}(u_H - \sumT \uft, v_H) + \rho^{1/2}\sumT a_{\mu}(\uft, \vft) - a_{\mu}^T(u_H, \vft),
\end{multline*}
with a stabilization parameter $\rho \geq 1$ that will be chosen later.
Further, let $\Fts \in \Vts^\prime$ be defined by
\begin{align*}
\Fts\left((v_H, \vfti{1}, \dots, \vfti{{\abs{\mathcal{T}_H}}})\right) &:= l(v_H),
\end{align*}
and let $\utsm \in \Vts$ of the two-scale solution of the PG--LOD by the variational problem
\begin{equation}
\Btsm\left(\utsm,  \vts \right) = \Fts(\vts) \qquad \textnormal{for all } \vts \in \mathcal{V}. \label{eq:two_scale_PGLOD}
\end{equation}
We show, that \eqref{eq:two_scale_PGLOD} is equivalent to the original PG--LOD formulation
\eqref{eq:loc_cor_problems},~\eqref{eq:PG_}:
\begin{proposition} \label{prop:two_scale_solution}
	The two-scale solution $\utsm \in \Vts$ of \cref{eq:two_scale_PGLOD} is uniquely determined and given by
	\begin{equation}\label{eq:two_scale_solution}
	\utsm = \left[\uhkm,\, \QQktm[T_1](\uhkm),\, \ldots,\, \QQktm[\Tlast](\uhkm)\right].
	\end{equation}
\end{proposition}
\begin{proof}
	With $\utsm$ as in~\eqref{eq:two_scale_solution} we have for any $\vts = (v_H, \vft[T_1], \ldots, \vft[\Tlast])$
	\begin{align*}
	\Btsm (\utsm, \vts) &=
	a_\mu(\uhkm - \sumT \QQktm(\uhkm), v_H)\\
	&\hS{50}+ \rho^{1/2} \sumT a_\mu(\QQktm(\uhkm), \vft) - a_\mu^T(\uhkm, \vft) \\
	&= l(v_H) = \Fts(\vts),
	\end{align*}
	where we have used the definition of $\QQktm(\uhkm)$ to eliminate the sum over $\Gridh$, the definition of $\QQkm$ in \cref{eq:loc_corrector} and the definition of $\uhkm$ in \cref{eq:PG_}.
	
	To show that $\utsm$ is the only solution of~\eqref{eq:two_scale_solution}, it suffices to show that $\Btsm(\uts, \vts) = 0$ for all $\vts\in \Vts$ implies $\uts = (u_H, \uft[T_1], \ldots, \uft[\Tlast]) = 0$.
	For that, first note that for $1 \leq i \leq \Tlast$ and each $\vft[T_i] \in \Vfhkt$ we have:
	\begin{equation*}
	a_\mu(\uft[T_i], \vft[T_i]) - a_\mu^T(u_H, \vft[T_i]) =
	\rho^{-1/2}\cdot \Btsm(\uts, (0, \ldots, 0, \vft[T_i], 0 \ldots, 0)) = 0,
	\end{equation*}
	hence, $\uft[T_i] = \QQktm(u_H)$.
	This implies $a_\mu(u_H - \QQkm(u_H), v_H) = \Btsm(\uts, (v_H, 0, \ldots, 0)) = 0$ for all $v_H \in V_H$, which means $u_H = 0$ due to the inf-sup stability of the PG--LOD bilinear form.
	Certainly, $u_H = 0$ implies $\QQktm(u_H) = 0$ for all $T \in \Gridh$.
\end{proof}

Concerning the computational complexity, we consider \eqref{eq:two_scale_PGLOD} as the FOM version of the two-scale formulation of the PG--LOD that performs all Steps~1-3 from \cref{def:comp_proc_LOD} at the same time.
We emphasize that, computationally, such an approach is highly prohibitive in the FOM version since the system matrix for \eqref{eq:two_scale_PGLOD} for $\kl > 1$ is even more significant than the system matrix for a classical FEM solution.
With \cref{prop:two_scale_solution}, however, a solution of \eqref{eq:two_scale_PGLOD} can be found with the original PG--LOD \eqref{eq:loc_cor_problems},~\eqref{eq:PG_}.

\subsection{Analysis of the two-scale bilinear form}\label{sec:two_scale_form}

We introduce two weighted norms $\Anorm{}$, $\SMnorm{}$ on $\Vts$, w.r.t.\ which we show the inf-sup stability of $\Btsm$ and derive approximation error bounds.
For arbitrary $\uts = (u_H, \uft[T_1], \dots, \uft[\Tlast]) \in \Vts$ these norms are given by
\begin{align*}
\Anorm{\uts}^2 &:= 
\anorm{u_H - \sumT \uft}^2
+ \rho\sumT \anorm{\QQktm(u_H) - \uft}^2,\\
\SMnorm{\uts}^2 &:= 
\snorm{u_H}^2
+ \rho\sumT \snorm{\QQktm(u_H) - \uft}^2.
\end{align*}

\begin{proposition}
	$\Anorm{}$ and $\SMnorm{}$ are norms on $\Vts$ for all $\mu \in \mathcal{P}$.
\end{proposition}
\begin{proof}
	Since $\uts \mapsto u_H - \sumT \uft$ and $\uts \mapsto \QQktm(u_H) - \uft$ are linear, the pull-back norms $\anorm{u_H - \sumT \uft}$ and $\anorm{\QQktm(u_H) - \uft}$ are semi-norms on $\Vts$.
	Hence, $\Anorm{}$ is a semi-norm on $\Vts$ as well.
	
	Further, we have $$\snorm{u_H} = \snorm{\IH(u_H - \sumT \hS{-2} \uft)} \leq \CI \snorm{u_H - \sumT \hS{-2} \uft} \leq \CI \alpha^{-1/2}\anorm{u_H - \sumT \hS{-2} \uft}.$$
	So, $\Anorm{u} = 0$ implies $u_H = 0$.
	This, in turn, implies $\anorm{\uft} = \anorm{\QQktm(u_H) - \uft} = 0$ for all $T \in \Gridh$, so $\utsm = 0$.
	Hence, $\Anorm{}$ is indeed a norm on $\Vts$.
	The argument for $\SMnorm{}$ is similar.
\end{proof}

We intend to show that $\Anorm{}$ and $\SMnorm{}$ are equivalent norms, for which we require some technical results.

\begin{lemma} \label{lemma:norm_prop_1}
	Let $\vft\in\Vfhkt$ for each $T \in \Gridh$ be given. Then, we have:
	\begin{equation*}
	\snorm{\sumT\vft}^2 \leq \CO \sumT\snorm{\vft}^2 \quad\text{and}\quad
	\anorm{\sumT\vft}^2 \leq \CO \sumT\anorm{\vft}^2.
	\end{equation*}
\end{lemma}
\begin{proof}
	Using Jensen's inequality, we have:
	\begin{equation*}
	\begin{aligned}
	\anorm{\sumT\vft}^2 &= \int_\Omega \abs{\sumT A_\mu^{1/2}(x)\nabla\vft(x)}^2\dx \\
	&\leq \int_\Omega \CO\cdot \sumT\abs{A_\mu^{1/2}(x)\nabla\vft(x)}^2\dx
	= \CO \sumT\anorm{\vft}^2.
	\end{aligned}
	\end{equation*}
	The proof for $\snorm{}$ is the same.
\end{proof}

\begin{lemma}\label{thm:qqktm_bound}
	For arbitrary $u_H \in V_H$, we have
	\begin{equation*}
	\left(\sumT \anorm{\QQktm (u_H)}^2\right)^{1/2} \leq \anorm{u_H}.
	\end{equation*}
\end{lemma}
\begin{proof}
	By definition of $\QQktm$, we have
	\begin{equation*}
	\begin{split}
	\sumT &\anorm{\QQktm (u_H)}^2 
	= \sumT a_\mu(\QQktm (u_H), \QQktm (u_H)) 
	= \sumT a_\mu^T(u_H, \QQktm (u_H)) \\
	&\!\!\!\!\!\!\!\leq \sumT \biggl(\int_T\abs{A_\mu^{1/2}(x)\nabla u_H(x)}^2\dx\biggr)^{1/2}
	\cdot \biggl(\int_T\abs{A_\mu^{1/2}(x)\nabla\QQktm(u_H)(x)}^2\dx\biggr)^{1/2} \\
	&\!\!\!\!\!\!\!\leq \biggl(\sumT \int_T\abs{A_\mu^{1/2}(x)\nabla u_H(x)}^2\dx\biggr)^{1/2}
	\cdot \biggl(\sumT \int_T\abs{A_\mu^{1/2}(x)\nabla\QQktm(u_H)(x)}^2\dx\biggr)^{1/2} \\
	&\!\!\!\!\!\!\!\leq \anorm{u_H}^2 \cdot \biggl( \sumT \anorm{\QQktm u_H}^2 \biggr)^{1/2}.
	\end{split}
	\end{equation*}
	Dividing by the second factor yields the claim.
\end{proof}

\noindent Now, we are prepared to show the equivalence of both norms.
\begin{proposition}\label{thm:norm_equiv}
	$\Anorm{}$ and $\SMnorm{}$ are equivalent norms on $\Vts$ with the following bounds for every $\uts \in \Vts$:
	\begin{equation*}
	\CI^{-1}\alpha^{1/2}\SMnorm{\uts}
	\leq \Anorm{\uts}
	\leq \sqrt{3}(1+\CO)^{1/2}\beta^{1/2}\SMnorm{\uts}.
	\label{eq:SMnormAnorm}
	\end{equation*}
\end{proposition}
\begin{proof}
	Let $\uts = (u_H, \uft[T_1], \ldots, \uft[\Tlast])$.
	To bound $\SMnorm{\uts}$ by $\Anorm{\uts}$, note that
	\begin{equation}\label{eq:two_scale_norm_equiv_proof_sa_bound}
	\snorm{u_H}^2 \leq \CI^2 \snorm{u_H - \sumT\uft}^2
	\leq \alpha^{-1}\CI^2\anorm{u_H - \sumT\uft}^2,
	\end{equation}
	and, using $\CI \geq 1$, it holds
	\begin{equation*}
	\rho\sumT\snorm{\QQktm(u_H) - \uft}^2 \leq
	\alpha^{-1}\CI^2\rho\sumT\anorm{\QQktm(u_H) - \uft}^2.
	\end{equation*}
	To bound $\Anorm{\uts}$ by $\SMnorm{\uts}$, we use \cref{lemma:norm_prop_1}, \cref{thm:qqktm_bound} and $\rho \geq 1$ to obtain
	\begin{align*}
	&\anorm{u_H - \sumT\uft}^2\\
	&\qquad\leq 3\anorm{u_H}^2 + 3\CO\sumT\anorm{\QQktm(u_H) - \uft}^2 + 3\CO\sumT\anorm{\QQktm(u_H)}^2 \\
	&\qquad\leq 3(1 + \CO)\anorm{u_H}^2 + 3\CO\rho\sumT\anorm{\QQktm(u_H) - \uft}^2\\
	&\qquad\leq 3(1 + \CO)\beta\snorm{u_H}^2 + 3\CO\beta\rho\sumT\snorm{\QQktm(u_H) - \uft}^2.
	\end{align*}
	Adding
	\begin{equation*}
	\rho\sumT\anorm{\QQktm(u_H) - \uft}^2 \leq
	\beta\rho\sumT\snorm{\QQktm(u_H) - \uft}^2
	\end{equation*}
	to both sides yields the claim.
\end{proof}
\noindent Finally, we show that $\Btsm$ is $\Anorm{}$-$\Snorm{}$ inf-sup stable with controllable constants.
\begin{proposition}\label{thm:two_scale_inf_sup}
	Let
	\begin{equation*}
	\rho := \CO\cdot\kappa,
	\end{equation*}
	then, $\Btsm$ is $\Anorm{}$-$\Snorm{}$-continuous and inf-sup stable with the following bounds on the respective constants:
	\begin{equation*}
	\sup_{0 \neq \uts \in \Vts}
	\sup_{0 \neq \vts \in \Vts}
	\frac{\Btsm(\uts,\vts)}{\Anorm{\uts}\cdot\Snorm{\vts}} 
	\leq \beta^{1/2} \quad \text{and}\quad
	\inf_{0 \neq \uts \in \Vts}
	\sup_{0 \neq \vts \in \Vts}
	\frac{\Btsm(\uts,\vts)}{\Anorm{\uts}\cdot\Snorm{\vts}} 
	\geq \CPG/\sqrt{5}.
	\end{equation*}
\end{proposition}
\begin{proof}
	We first bound the continuity constant of $\Btsm$.
	Let $\uts = (u_H, \uft[T_1], \ldots, \uft[\Tlast]) \in \Vts$ and $\vts = (v_H, \vft[T_1], \ldots, \vft[\Tlast]) \in \Vts$ be arbitrary.
	Then, we have
	\begin{equation*}
	\begin{aligned}
	\Btsm(\uts, \vts)
	& = a_\mu(u_H - \sumT\uft, v_H)
	- \rho^{1/2}\sumT\left(a_\mu(\uft,\vft) - a_\mu^T(u_H,\vft)\right)\\
	& \leq \anorm{u_H-\sumT\uft}\anorm{v_H}
	+ \sumT\rho^{1/2}\anorm{\uft-\QQktm(u_H)}\anorm{\vft}\\
	& \leq \left[\anorm{u_H - \sumT\uft}^2
	+ \sumT\rho\anorm{\uft-\QQktm(u_H)}^2 \right]^{1/2}\\
	& \qquad\qquad\cdot\left[\anorm{v_H}^2 + \sumT\anorm{\vft}^2\right]^{1/2}\\
	& \leq \Anorm{\uts} \cdot \beta^{1/2} \cdot \Snorm{\vts}.
	\end{aligned}
	\end{equation*}
	To prove inf-sup stability, first note that
	\begin{align*}
	\anorm{u_H-\sumT\uft}\hspace{-6em}& \\
	& \leq \anorm{u_H-\sumT\QQktm(u_H)}
	+ \anorm{\sumT(\QQktm(u_H) - \uft)} \\
	& \leq (\CPG)^{-1} \sup_{0\neq v_H\in V_H} \frac{a_\mu(u_H-\sumT\QQktm(u_H),v_H)}{\snorm{v_H}}
	+ \anorm{\sumT(\QQktm(u_H) - \uft)} \\
	& \leq (\CPG)^{-1}\underbrace{\sup_{0\neq v_H\in V_H}\hspace{-2pt}
		\frac{\Btsm(\uts,(v_H,0,\ldots,0))}{\Snorm{(v_H,0,\ldots,0)}}}_{A}
	\\&\quad+ (\CPG)^{-1}\hspace{-2pt}\sup_{0\neq v_h\in V_H}\hspace{-2pt}
	\frac{a_\mu(\sumT(\QQktm(u_H)-\uft),v_H)}{\snorm{v_H}} + \anorm{\sumT\hspace{-2pt}(\QQktm(u_H) - \uft)} \\
	& \leq (\CPG)^{-1}A
	+ \left((\CPG)^{-1}
	\sup_{0\neq v_H\in V_H}\frac{\anorm{v_H}}{\snorm{v_H}}
	+ 1\right) \anorm{\sumT(\QQktm(u_H) - \uft)} \\
	& \leq (\CPG)^{-1}A 
	+ ((\CPG)^{-1}\beta^{1/2} + 1)\CO^{1/2}\left(\sumT\anorm{\QQktm(u_H) - \uft}^2\right)^{1/2},
	\end{align*}
	hence:
	\begin{equation*}
	\begin{aligned}
	\Anorm{\uts}^2
	& = \anorm{u_H-\sumT\uft}^2 + \rho\sumT\anorm{\QQktm(u_H)-\uft}^2\\
	& \leq 2(\CPG)^{-2}A^2
	+ (2(\CPG)^{-2}\beta\CO +2\CO + \rho)\sumT\anorm{\QQktm(u_H) - \uft}^2.
	\end{aligned}
	\end{equation*} 
	Further,
	\begin{equation*}
	\begin{aligned}
	\left(\sumT\anorm{\uft - \QQktm(u_H)}^2\right)^{1/2}
	&= \frac{\sumT a_\mu(\uft-\QQktm(u_H), \uft-\QQktm(u_H))}{(\sumT\anorm{\uft-\QQktm(u_H)}^2)^{1/2}}\\
	& \leq \sup_{\vft\in\Vfhkt}\frac{\sumT a_\mu(\uft-\QQktm(u_H), \vft)}{(\sumT\anorm{\vft}^2)^{1/2}}\\
	& \leq \alpha^{-1/2}\sup_{\vft\in\Vfhkt}\frac{\sumT a_\mu(\uft, \vft) - a_\mu^T(u_H, \vft)}{(\sumT\snorm{\vft}^2)^{1/2}}\\
	& = \alpha^{-1/2}\rho^{-1/2}\underbrace{\sup_{\vft\in\Vfhkt}
		\frac{\Btsm(\uts, (0,\vft[T1],\ldots,\vft[\Tlast]))}
		{\Snorm{(0,\vft[T1],\ldots,\vft[\Tlast])}}}_{B}.
	\end{aligned}
	\end{equation*}
	Combining both estimates yields
	\begin{equation*}
	\begin{aligned}
	\Anorm{\uts}^2
	&\leq 2(\CPG)^{-2}A^2 + (2(\CPG)^{-2}\alpha^{-1}\beta\CO\rho^{-1} +
	2\alpha^{-1}\CO\rho^{-1}+\alpha^{-1})B^2\\
	&= 2(\CPG)^{-2}A^2 + \underbrace{(2(\CPG)^{-2} +
		2\alpha^{-1}\kappa^{-1}+\alpha^{-1})}_{\leq 5\cdot(\CPG)^{-2}}B^2.\\
	&\leq 5\cdot(\CPG)^{-2}(A^2 + B^2) \\
	&=5\cdot(\CPG)^{-2}\left(\sup_{\vts\in\Vts}\frac{\Btsm(\uts,\vts)}{\Snorm{\vts}}\right)^2,
	\end{aligned}
	\end{equation*}
	where we have used $\CPG \leq \alpha^{1/2}$ and $\kappa \geq 1$ in the first equality.
	In the last equality we have used the fact that the square norm of a linear functional on a direct sum of Hilbert spaces is the sum of the square norms of the functional restricted to the respective subspaces.
\end{proof}

\subsection{Error bounds}
\label{sec:TSRBLOD_a-posteriori_estimators}
Exploiting the inf-sup stability of the two-scale bilinear form $\Btsm$, we now quickly obtain error bounds w.r.t.\ the PG--LOD solution.
We start with an a posteriori bound and define for arbitrary $\uts \in \Vts$ the residual-based error indicators
\begin{align}\label{eq:stage2_est_a}
\eta_{a,\mu}(\uts) &:= 
\sqrt{5}(\CPG)^{-1}
\sup_{v\in\Vts}
\frac{\Fts(\vts) - \Btsm(\uts,\vts)}
{\Snorm{\vts}},\\\label{eq:stage2_est_s}
\eta_{1,\mu}(\uts) &:= 
\sqrt{5}\CI\alpha^{-1/2}(\CPG)^{-1}
\sup_{v\in\Vts}
\frac{\Fts(\vts) - \Btsm(\uts,\vts)}
{\Snorm{\vts}}.
\end{align}
These error indicators provide strict and efficient upper bounds on the error between $\uts$ and the  two-scale PG--LOD solution:
\begin{theorem}[A-posteriori Bound]\label{thm:a_posteriori}
	Let $\uts= (u_H, \uft[T_1], \ldots, \uft[\Tlast])\in\Vts$ be an arbitrary two-scale function, denote by $\utsm$ the solution of the two-scale solution as in \cref{eq:two_scale_solution} for a given parameter $\mu$, and let $\rho$ be given as in \cref{thm:two_scale_inf_sup}.
	Then, the following energy error bounds hold:
	\begin{equation}\label{eq:a_posteriori_anorm}
	\Anorm{\utsm - \uts}
	\leq \eta_{a,\mu}(\uts)
	\leq \sqrt{5}(\CPG)^{-1}\beta^{1/2} \Anorm{\utsm - \uts}.
	\end{equation}
	Further, we have 
	\begin{equation}\label{eq:a_posteriori_snorm}
	\left(\snorm{\uhkm - u_H}^2 + \rho\sumT\snorm{\QQktm(u_H) - \uft}^2\right)^{1/2}
	\leq \eta_{1,\mu}(\uts),
	\end{equation}
	and
	\begin{equation}\label{eq:a_posteriori_snorm_eff}
	\begin{split}
	\eta_{1,\mu}(\uts) &\leq \sqrt{15}\CI(\CO+1)^{1/2}\kappa^{1/2}(\CPG)^{-1}\beta^{1/2}\\&\qquad\qquad\cdot
	\left(\snorm{\uhkm - u_H}^2 + \rho\sumT\snorm{\QQktm(u_H) - \uft}^2\right)^{1/2}.
	\end{split}
	\end{equation}
\end{theorem}
\begin{proof}
	Since $\utsm$ is a solution of \cref{eq:two_scale_PGLOD}, we have
	\begin{equation*}
	\Btsm(\utsm - \uts,\vts) =
	\Fts(\vts) - \Btsm(\uts, \vts).
	\end{equation*}
	Hence, \cref{eq:a_posteriori_anorm} directly follows from \cref{thm:two_scale_inf_sup}.
	\Cref{eq:a_posteriori_snorm,eq:a_posteriori_snorm_eff} follow from \cref{eq:a_posteriori_anorm} using \cref{thm:norm_equiv} and noting that for each $T \in \Gridh$ we have 
	\begin{equation*}
	\snorm{\QQktm(\uhkm-u_H) - (\QQktm(\uhkm)-\uft)}=\snorm{\QQktm(u_H)-\uft}.
	\end{equation*}~
\end{proof}
\noindent
Finally, we also show a corresponding a priori result:
\begin{theorem}[A priori bound]\label{thm:a_priori}
	Let $\overline{\Vts}$ be an arbitrary linear subspace of $\Vts$ and let $\overline{\uts}$ be the solution of the residual-minimization problem
	\begin{equation}\label{eq:red_solution_abstract}
	\overline{\uts}_\mu:= \argmin_{\uts\in\overline{\Vts}}
	\sup_{\vts\in \Vts}\frac{\Fts(v) - \Btsm(\uts, \vts)}
	{\Snorm{\vts}},
	\end{equation}
	then we have
	\begin{equation*}
	\Anorm{\utsm - \overline{\uts}_\mu}
	\leq \sqrt{5}(\CPG)^{-1}\beta^{1/2}
	\min_{\overline{\vts}\in\overline{\Vts}}\Anorm{\utsm - \overline{\vts}},
	\end{equation*}
	and
	\begin{equation*}
	\begin{split}
	\left(\snorm{\uhkm - u_H}^2 + \rho\sumT\snorm{\QQktm(u_H) - \uft}^2\right)^{1/2}\hspace{-20em}&\\
	&\leq \sqrt{15}\CI(\CO+1)^{1/2}\kappa^{1/2}(\CPG)^{-1}\beta^{1/2} \\&\qquad\qquad\cdot
	\min_{\overline{\vts}\in\overline{\Vts}}
	\left(\snorm{\uhkm - \overline{v}_H}^2 + \rho\sumT\snorm{\QQktm(u_H) - \overline{v}_T^\textnormal{f}}^2\right)^{1/2}.
	\end{split}
	\end{equation*}
\end{theorem}
\begin{proof}
	This follows directly from \cref{thm:a_posteriori} and the definition of $\overline{\uts}_\mu$.
\end{proof}

Note that for $\kl$ large enough, we have $\CPG \approx \alpha^{1/2}\CI^{-1}$, such that the a priori and a posteriori bounds have efficiencies of that scale with $\kappa^{1/2}$ in the energy norm and with $\kappa$ in the 1-norm.
This agrees with what is to be expected for these error bounds in a standard finite element setting; cf.~\cref{sec:background_npdgl}.

\section{Two-scale reduced basis approach}
\label{sec:TSRBLOD_MOR} 
This section describes the construction of the two-scale RB approach for the LOD (TSRBLOD).
As with all RB methods, the ROM is defined by a projection of the original model equations, in our case, the two-scale formulation \eqref{eq:two_scale_PGLOD}, onto a reduced approximation space (\cref{sec:stage_2_rom}).
Rigorous upper and lower bounds for the MOR error are given by a residual-based a posteriori error estimator (\cref{sec:stage_2_err}).
To be able to assemble the ROM for a new parameter $\mu$ quickly and to evaluate the error estimator efficiently, an offline-online decomposition of the ROM must be performed (\cref{sec:stage_2_offon}).
As part of the offline phase, the reduced space is constructed as the linear span of FOM solutions $\mathfrak{u}_{\mu^*}$, where the snapshot parameters $\mu^*$ are selected via an iterative greedy-search over $\Params$ (\cref{sec:stage_2_bas_gen}).

The computation of the solution snapshots $\mathfrak{u}_{\mu^*}$ for each new $\mu^*$ via \eqref{eq:two_scale_solution}, however, requires the recomputation of all corrector problems.
On top of that, the dimension of the system matrix of \cref{eq:two_scale_solution} and thus, also the dimension of the solution snapshots scale with the number of coarse-mesh elements times the number of fine-mesh elements in each patch.
For large problems, this may be computationally infeasible.
To remedy this, we combine our approach (Stage~2) with a preceding preparatory step similar to~\cite{RBLOD}, where each corrector problem is replaced by an efficient ROM surrogate (Stage~1); cf. \cref{fig:RBLOD}.
These ROMs are then used in the offline phase of Stage~2 to compute approximate solution snapshots $\mathfrak{u}_{\mu^*}$.
Again, Stage~1 is divided into ROM construction (\cref{sec:stage_1_rom}), error estimation (\cref{sec:stage_1_err}), offline-online decomposition (\cref{sec:stage_1_of/on}) and construction of the reduced spaces (\cref{sec:stage1_basisgen}).
After a Stage~1 ROM is constructed, all associated fine-mesh data can be deleted.
In particular, all computations in Stage~2 are independent of $\grid$. We describe the two stages for a fixed localization parameter $\kl \in \mathbb{N}$.

\subsection{Stage~1: RB approximations of the fine-scale correctors}
\label{sec:stage1}
\subsubsection{Definition of the reduced-order model}\label{sec:stage_1_rom}%
Let $T \in \Gridh$ be fixed, and let $\Vfrbkt \subset \Vfhkt$ be an approximation space for the correctors $\QQktm(v_H)$ of $v_H$ for arbitrary $v_H \in V_H$ and $\mu \in \mathcal{P}$.
Then, for given $v_H \in V_H$ and $\mu \in \mathcal{P}$, we determine an approximate corrector $\QQktmrb(v_H) \in \Vfrbkt$ via Galerkin projection onto $\Vfrbkt$ as the solution of
\begin{equation}\label{eq:stage1_rom}
a_\mu(\QQktmrb(v_H), \vft) = a_\mu^T(v_H, \vft) \qquad \textnormal{for all } \vft \in
\Vfrbkt.
\end{equation}
Note that $\QQktmrb(v_H)$ is well-defined since $a_\mu$ is a coercive bilinear form.
Using these reduced correction operators we can define an approximate localized multiscale matrix $\Kmurb$ given by
\begin{equation*}
\Kmurb:= \sumT \Ktmurb, \quad \left(\Ktmurb \right)_{ji} := \skal{A_{\mu} (\chi_T \nabla - \nabla \QQktmrb)
	\phi_i}{\nabla \phi_j}_{U_\kl(T)}.
\end{equation*}

\begin{remark}[Flexibility in the reduction process] \label{rem:reduction_flexibility}
	Although it is not directly visible in \eqref{eq:stage1_rom}, the technical reduction process for the corrector problems can be realized differently, which can also affect the approximation behavior.
	Recall that we require a corrector function for all shape functions on $T$, where the different systems differ in terms of their right-hand side since different shape functions are inserted (for details, see \cref{sec:stage_1_of/on}).
	In a nutshell, our reduction approach treats these different right-hand sides as an additional parameter and the reduced basis is built for all right-hand sides at once.
	This has the effect that the two-scale matrix in Stage~2 is smaller.
	As we briefly discuss later, indeed, the reduction approach that has been chosen for the RBLOD~\cite{RBLOD} constructs a reduced basis for all shape functions.
	Another approach that we do not detail further is the simultaneous RB construction of a block system with a non-parameterized right-hand side, similar to the two-scale matrix.
	To the best of our knowledge, such a reduction process has not yet been considered in a general setting and is left to future research.
\end{remark}

\subsubsection{Error estimation}\label{sec:stage_1_err}
We employ standard RB tools for the a posteriori error estimation of the reduced system \eqref{eq:stage1_rom}.
In detail, we use the residual-norm based estimate:
\begin{equation}\label{eq:stage1_error_estimate}
\anorm{\QQktm(v_H) - \QQktmrb(v_H)} 
\leq \eta_{T,\mu}(\QQktmrb(v_H))
\leq \kappa^{1/2}\anorm{\QQktm(v_H) - \QQktmrb(v_H)},
\end{equation}
where
\begin{equation} \label{eq:stage_1_estimator}
\eta_{T,\mu}(\QQktmrb(v_H)) := \alpha^{-1/2} \sup_{\vft\in\Vfhkt}
\frac{a_\mu^T(v_H, \vft) - a_\mu(\QQktmrb(v_H), \vft)}{\snorm{\vft}}.
\end{equation}
The bounds in \eqref{eq:stage1_error_estimate} easily follow from the definition of $\QQktm(v_H)$ and the equivalence of $\anorm{}$ and $\snorm{}$.

\subsubsection{Offline-online decomposition}\label{sec:stage_1_of/on}
In order to compute $\QQktmrb(v_H)$ and subsequently $\Ktmurb$, let $\dimVfrbkt := \dim \Vfrbkt$, and choose a basis $\rbt{n}$, $1 \leq n \leq \dimVfrbkt$, of $\Vfrbkt$.
Expanding $\QQktmrb(v_H)$ w.r.t.\ this basis as
\begin{equation*}
\QQktmrb(v_H) := \sum_{n=1}^{\dimVfrbkt} c_n\cdot\rbt{n},
\end{equation*}
the coefficient vector $c \in \mathbb{R}^{\dimVfrbkt}$ is given as the solution of the $\dimVfrbkt \times \dimVfrbkt$-dimensional linear system
\begin{equation}\label{eq:stage1_rom_in_matrices}
\mathbb{A}^T_\mu \cdot c = \mathbb{G}^T_\mu(v_H),
\end{equation}
where
\begin{equation*}
(\mathbb{A}^T_{\mu})_{m,n} := a_\mu(\rbt{n}, \rbt{m}) \qquad\text{and}\qquad
\mathbb{G}^T_\mu(v_H)_m := a^T_\mu(v_H, \rbt{m}).
\end{equation*}
While~\eqref{eq:stage1_rom_in_matrices} can be solved quickly when $\dimVfrbkt$ is sufficiently small, we still need to reassemble this equation system for each new parameter $\mu$ and coarse-scale function $v_H$.
This requires time and memory that scales with $\dim \Vfhkt$.
To avoid these high-dimensional computations we exploit \cref{asmpt:parameter_separable} and pre-assemble matrices and vectors
\begin{equation*}
(\mathbb{A}^T_{\xi})_{m,n} := a_\xi(\rbt{n}, \rbt{m}) \qquad\text{and}\qquad
(\mathbb{G}^T_{\xi,j})_{m} := a^T_\xi(\phi_{i_{T,j}}, \rbt{m}),
\end{equation*}
for $1 \leq q \leq \Qa$ and $1 \leq j \leq \KT$, where $\KT$ is the number of finite-element basis functions $\phi_{i}$ of $V_H$ with support containing $T$ and $i_{T,1}, \ldots, i_{T,\KT}$ is an enumeration of these basis functions.
Then, $\mathbb{A}^T_\mu$ and $\mathbb{G}^T_\mu(v_H)$ can be determined as
\begin{align*}
\mathbb{A}^T_{\mu} &:= \sum_{\xi=1}^{\Qa}\theta_\xi^a(\mu) \mathbb{A}^T_{\xi} \qquad\text{and}\qquad \\
\mathbb{G}^T_\mu(v_H) &:= \sum_{\xi=1}^{\Qa}\sum_{j=1}^{\KT}\theta_\xi^a(\mu)\lambda_{i_{T,j}}(v_H) \mathbb{G}^T_{\xi,j},
\end{align*}
where by $\lambda_{i}\in V_H^\prime$ we denote the dual basis of $\phi_i$.

To compute $\Kmurb$, we further store the matrices 
\begin{equation*} 
\left(\Ktqo\right)_{j,i} := \skal{A_\xi \chi_T \nabla \phi_i}{\nabla \phi_j}_{U_\kl(T)}
\qquad\text{and}\qquad
\left(\Ktqrb\right)_{j,n} := \skal{A_{\xi} \nabla \rbt{n}}{\nabla \phi_j}_{U_\kl(T)}.
\end{equation*}
Then we have:
\begin{equation*} 
\left(\Ktmurb\right)_{j,i} = 
\sum_{\xi=1}^{\Qa}\theta_\xi^a(\mu)\left[\left(\Ktqo\right)_{j,i} -
\sum_{n=0}^{N_T} c^i_n(\mu) \left(\Ktqrb\right)_{j,n}\right],
\end{equation*}
where $c^i(\mu) \in \mathbb{R}^{N_T}$ is given as the solution of
\begin{equation*}
\mathbb{A}^T_\mu \cdot c^i(\mu) = \mathbb{G}^T_\mu(\phi_i).
\end{equation*}
Note that $(\K^T_{\xi,0})_{j,i}$, $(\K^T_\xi)_{j,n}$ and $c^i(\mu)$ are zero unless the support of $\phi_i$ is non-disjoint from $T$ and the support of $\phi_j$ is non-disjoint from $U_\kl(T)$.
In particular, only $\KT$ reduced problems have to be solved to determine $\Ktmurb$.
The total computational effort for solving these problems is of order $\mathcal{O}(\Qa N_T^2 + N_T^3 + \KT N_T^2)$, where the first term corresponds to the assembly of $\mathbb{A}_\mu^T$, the second term to its LU decomposition and the third term to the solution of the $\KT$ linear systems using forward/backward substitution.

Finally, to efficiently evaluate $\eta_{T,\mu}(\QQktmrb(v_H))$, first note that
\begin{equation*}
\eta_{T,\mu}(\QQktmrb(v_H)) = \alpha^{-1/2} \norm{\mathcal{R}_{\Vfhkt}(a_\mu^T(v_H, \cdot) - a_\mu(\QQktmrb(v_H),
	\cdot))}_1,
\end{equation*}
where $\mathcal{R}_{\Vfhkt}: (\Vfhkt)^\prime \to \Vfhkt$ denotes the Riesz-isomorphism.

Following~\cite{buhr14}, let $\Wfrbkt$ denote the $M_T$-dimensional linear subspace of $\Vfhkt$ that is spanned by the vectors
\begin{align*}
\set[\mathcal{R}_{\Vfhkt}(a_\xi(\rbt{n}, \cdot))]{1\leq q \leq \Qa, 1\leq n \leq N_T} \,\cup\,\hspace{-20em}& \\&
\qquad\qquad\set[\mathcal{R}_{\Vfhkt}(a_\xi^T(\phi_{i_{T,j}}, \cdot))]{1\leq q \leq \Qa, 1\leq j \leq \KT}.
\end{align*}
Choose an $H^1$-orthonormal basis $\rbestt{m}$ of $\Wfrbkt$.
Since $\mathcal{R}_{\Vfhkt}(a_\mu^T(v_H, \cdot) - a_\mu(u_N, \cdot))$ lies in $\Wfrbkt$ for each $u_N \in \Vfrbkt$, we obtain
\begin{equation*}
\begin{aligned}
\eta_{T,\mu}(\QQktmrb(v_H)) 
&= \alpha^{-1/2} \norm*{\left[\left(\mathcal{R}_{\Vfhkt}(a_\mu^T(v_H, \cdot) - a_\mu(\QQktmrb(v_H), \cdot)),
	\rbestt{m}\right)_1\right]^{M_T}_{m=1}} \\
&=\alpha^{-1/2} \norm*{\left[a_\mu^T(v_H, \rbestt{m}) - a_\mu(\QQktmrb(v_H), \rbestt{m})\right]^{M_T}_{m=1}}.
\end{aligned}
\end{equation*}
Hence, defining the $M_T$-dimensional vectors and $M_T \times N_T$ matrices
\begin{equation}
(\hat{\mathbb{G}}^T_{\xi,j})_{m} := a^T_\xi(\phi_{i_{T,j}}, \rbestt{m})\qquad\text{and}\qquad
(\hat{\mathbb{A}}^T_{\xi})_{m,n} := a_\xi(\rbt{n}, \rbestt{m}),
\end{equation}
we obtain
\begin{equation}\label{eq:stage1_error_estimate_in_matrices}
\eta_{T,\mu}(\QQktmrb(v_H)) = \alpha^{-1/2}
\norm*{
	\sum_{\xi=1}^{\Qa}\theta_\xi^a(\mu)\left(\sum_{j=1}^{\KT}\lambda_{i_{T,j}}(v_H) \hat{\mathbb{G}}^T_{\xi,j}
	- \hat{\mathbb{A}}^T_{\xi} \cdot c\right)
},
\end{equation}
where $c$ is again given by \eqref{eq:stage1_rom_in_matrices}.
We remark that $M_T$ can be bounded by $\Qa(N_T + \KT)$.
Hence, the cost of evaluating \eqref{eq:stage1_error_estimate_in_matrices} is of order  $\mathcal{O}(\Qa^2(N_T + \KT)^2)$.
 
\subsubsection{Basis generation}\label{sec:stage1_basisgen}

To build the reduced spaces $\Vfrbkt$, we use a standard weak greedy approach, as explained in \cref{sec:background_basis_gen_RB}, in order to minimize the model order reduction error $\QQktm(v_H) - \QQktmrb(v_H)$ for all $\mu \in \Params$ and $v_H \in V_H$.
To this end, we choose a target error tolerance $\varepsilon_1$ and an appropriate training set of parameters $\ParamsTrain$, over which we estimate the maximum reduction error.
The initial reduced space $\Vfrbkt$ is chosen as the zero-dimensional space.
Then, in each iteration, the reduction error is estimated with $\eta_{T,\mu}(\QQktmrb(v_H))$ for all $\mu \in \ParamsTrain$ and $\phi_{i_{T,j}}$, and a pair $\mu^*$, $\phi_{i_{T,j^*}}$ maximizing the estimate is selected.
Thanks to the offline-online decomposition of $\eta_{T,\mu}$ this step does not involve any high-dimensional computations, so $\ParamsTrain$ can be chosen large.
Since~\eqref{eq:stage1_rom} is linear, it suffices to consider the basis functions $\phi_{i_{T,j}}$ as $v_H$.
After $\mu^*$, $j^*$ have been found, $\QQktmrb[T][\mu^*](\phi_{i_{T,j^*}})$ is computed. 
The reduced space $\Vfrbkt$ is extended with this solution-\emph{snapshot}, and the offline-online decomposition for this expanded reduced space is computed.
The iteration ends when the maximum estimated error drops below $\varepsilon_1$.
A formal definition of the procedure is given in Algorithm~\ref{alg:stage1_greedy}.

\begin{algorithm2e}
	\KwData{Coarse element $T$, training set $\ParamsTrain$, tolerance $\varepsilon_1$}
	\KwResult{$\Vfrbkt$}
	$\Vfrbkt \leftarrow \{0\}$\;
	\While{$\max_{\mu\in\ParamsTrain}\max_{1\leq j\leq \KT}\eta_{T,\mu}(\QQktmrb(\phi_{i_T,j})) > \varepsilon_1$}{
		$(\mu^*, j^*) \leftarrow \argmax_{(\mu^*, j^*) \in \ParamsTrain \times \{1,\ldots,\KT\}}\eta_{T,\mu}(\QQktmrb(\phi_{i_T,j}))$\;
		$\Vfrbkt \leftarrow \Span(\Vfrbkt \cup \{\QQktm[T][\mu^*](\phi_{i_{T,j^*}})\})$\;
	}
	\caption{Weak greedy algorithm for the generation of $\Vfrbkt$.}\label{alg:stage1_greedy}
\end{algorithm2e}

\subsection{Stage~2: RB approximations of the two-scale solutions}
\label{sec:stage_2_red}

\subsubsection{Definition of the reduced-order model}\label{sec:stage_2_rom}
To find an approximate solution of $\utsm$, we assume to be given an appropriate reduced subspace $\Vtsrb$ of $\Vts$. 
As we have proven inf-sup stability of $\Btsm$ in \cref{thm:two_scale_inf_sup}, and since the inf-sup stability is preserved by restricting $\Btsm$ to a linear subspace, we can define the reduced two-scale solution $\utsmrb$ as the unique solution of the residual minimization problem
\begin{equation}\label{eq:stage2_rom}
\utsmrb := \argmin_{\uts \in \Vtsrb} \sup_{\vts \in \Vts} \frac{\Fts(\vts) - \Btsm(\uts, \vts)}{\Snorm{\vts}}.
\end{equation}
As a direct consequence of \cref{thm:a_priori}, $\utsmrb$ is a quasi best-approximation of $\utsm$ within $\Vtsrb$. 

\subsubsection{Error estimation}\label{sec:stage_2_err}
We have already defined a posteriori error estimators for approximations of the two-scale solution~$\utsm$ in \cref{sec:TSRBLOD_a-posteriori_estimators}.
In \cref{thm:a_posteriori}, we have shown that these estimators yield efficient upper bounds for the approximation errors in the two-scale energy norm as well as in the Sobolev 1-norm.
Note that, even though we will use the Stage~1 approximations $\Ktmurb$ of~$\Ktmu$ to build the reduced space $\Vtsrb$, the derived error estimates are with respect to the true LOD solution and take these approximation errors into account.
Thus, in contrast to~\cite{RBLOD}, the derived error estimator entirely takes the effect of the errors of the reduced corrector problems on the resulting global solution into account and, thus, rigorously bounds the error of the ROM w.r.t.\ the LOD solution.

\subsubsection{Offline-Online Decomposition}~\label{sec:stage_2_offon}
For the offline-online decomposition of \eqref{eq:stage2_rom}, we proceed similar to the decomposition of the Stage~1 error estimator $\eta_{T,\mu}$.
Denote by $\mathcal{R}_{\Vts}: \Vts^\prime \to \Vts$ the Riesz-isomorphism for $\Vts$.
Then~\eqref{eq:stage2_rom} is equivalent to solving
\begin{equation}\label{eq:stage2_rom_riesz}
\utsmrb := \argmin_{\uts \in \Vtsrb} \Snorm{\mathcal{R}_{\Vts}(\Fts) - \mathcal{R}_{\Vts}(\Btsm(\uts, \cdot))}^2.
\end{equation}
Let $N := \dim \Vtsrb$, and let $\rbts{n}$, $1 \leq n \leq N$ be a basis of $\Vtsrb$. 
We again construct an $\Snorm{}$-orthonormal basis $\rbestts{m}$ for the $M$-dimensional subspace $\Wtsrb$ of $\Vts$ spanned by the vectors
\begin{equation}\label{eq:two_scale_estimator_generators}
\{\mathcal{R}_{\Vts}(\Fts)\} \cup \{\mathcal{R}_{\Vts}(\Btsq(\rbts{n}, \cdot)) \,|\, 1\leq n \leq N, 1 \leq q
\leq \Qa \},
\end{equation}
where
\begin{equation*}
\begin{split}
\Btsq\left((u_H, \ufti{1}, \dots, \ufti{{\abs{\mathcal{T}_H}}}),  (v_H, \vfti{1}, \dots, \vfti{{\abs{\mathcal{T}_H}}})\right) &:=
\\ a_{\xi}(u_H - \sumT \uft, v_H) + &\rho^{1/2}\sumT a_{\xi}(\uft, \vft) - a_{\xi}^T(u_H, \vft),
\end{split}
\end{equation*}
such that $\Btsm$ has the decomposition:
$
\Btsm = \sum_{\xi=1}^{\Qa} \theta_\xi^a(\mu) \Btsq.
$
Using these bases, we define matrices $\hat{\mathbb{A}}_\xi \in \mathbb{R}^{M\times N}$ and the vector $\hat{\mathbb{F}} \in \mathbb{R}^M$ by
\begin{equation*}
(\hat{\mathbb{A}}_\xi)_{m,n} := \Btsq(\rbts{n}, \rbestts{m}) \qquad\text{and}\qquad
\hat{\mathbb{F}}_m := \Fts(\rbestts{m}).
\end{equation*}
Then, with $\hat{\mathbb{A}}_\mu := \sum_{\xi=1}^{\Qa} \theta_\xi^a(\mu)\hat{\mathbb{A}}_\xi,$ solving~\eqref{eq:stage2_rom_riesz} is equivalent to solving the least-squares problem
\begin{equation}\label{eq:stage2_rom_in_matrices}
c(\mu) := \argmin_{c \in \mathbb{R}^N} \norm{\hat{\mathbb{F}} - \hat{\mathbb{A}}_\mu\cdot c}^2,
\end{equation}
where $\utsm = \sum_{n=1}^N c_n(\mu) \rbts{n}$.
In the same way we can evaluate the error bounds $\eta_{a,\mu}(\utsm)$ and~$\eta_{1,\mu}(\uts)$ as~
\begin{align*}
\eta_{a,\mu}(\utsm) &= 
\sqrt{5}(\CPG)^{-1}
\norm{\hat{\mathbb{F}} - \hat{\mathbb{A}}_\mu\cdot c(\mu)},\\ 
\eta_{1,\mu}(\utsm) &= 
\sqrt{5}\CI\alpha^{-1/2}(\CPG)^{-1}
\norm{\hat{\mathbb{F}} - \hat{\mathbb{A}}_\mu\cdot c(\mu)}.
\end{align*}
Since $M \leq \Qa N + 1$, the computational effort for assembling the least-squares system is of order $\mathcal{O}(\Qa^2N^2)$.
Solving the system requires $\mathcal{O}(\Qa N^3)$ operations, and evaluating the estimators requires $\mathcal{O}(\Qa^2 N^2)$ operations.
In particular, the computational effort is entirely independent of~$h$ and $H$. 

We still need to show how the matrices $\hat{\mathbb{A}}_\xi$ and the vector $\hat{\mathbb{F}}$ can be computed after Stage~1 without using any data or operations associated with the fine mesh $\grid$.
To this end, we assume that $\Vtsrb \subseteq V_H \oplus \Vfrbkt[T_1] \oplus \dots \oplus \Vfrbkt[T_{\Tlast}] \subset \Vts$.
By construction of the $\Wfrbkt$, we see that for such a $\Vtsrb$, $\Wtsrb$ is a linear subspace of $V_H \oplus \Wfrbkt[T_1] \oplus \dots \Wfrbkt[\Tlast]$.
Choose a basis $\rbts{n}$ of $\Vtsrb$ with coefficient vectors
\begin{equation*}
\rbtscoeff{n} = (\rbtscoeff{n,H}, \rbtscoeff{n,T_1}, \dots, \rbtscoeff{n,\Tlast}) \in
\mathbb{R}^{\dimVH} \oplus \mathbb{R}^{N_{T_1}} \oplus \dots \oplus \mathbb{R}^{N_{\Tlast}}
\end{equation*}
w.r.t.\ the finite element basis $\phi_i$ of $V_H$ and the reduced bases $\rbt{n}$ of $\Vfrbkt$.
Denote by $\mathbb{S}$ the $\dimVH \times \dimVH$ matrix of the Sobolev 1-inner product on $V_H$ given by
\begin{equation*}
\mathbb{S}_{j,i} := \int_\Omega \nabla \phi_i \cdot \nabla \phi_j \,\mathrm{dx},
\end{equation*}
and use the $H^1$-orthonormal bases $\rbestt{m}$ to isometrically represent the vectors~\eqref{eq:two_scale_estimator_generators} as coefficient vectors in the direct sum Hilbert space $\underline{\mathfrak{W}}^{rb} := \mathbb{R}^{\dimVH} \oplus \mathbb{R}^{M_{T_1}} \oplus \dots \mathbb{R}^{M_{\Tlast}}$ equipped with the $\mathbb{S}$-inner product in the first and with the Euclidean inner products in the remaining components.
Checking the definitions of $\Btsq$, $\Fts$, $\hat{\mathbb{A}}^T_\xi$, $\hat{\mathbb{G}}^T_{\xi,j}$, $\Ktqo$ and $\Ktqrb$, we obtain the vectors
\begin{equation}\label{eq:two_scale_est_gen_in_coord_f}
(\mathbb{S}^{-1}\cdot[l(\phi_i)]_i,0\dots,0),
\end{equation}
and
\begin{equation}\label{eq:two_scale_est_gen_in_coord_b}
\begin{pmatrix}
\mathbb{S}^{-1}\cdot\sumT \left(\Ktqo \cdot \rbtscoeff{n,H} - \Ktqrb \cdot \rbtscoeff{n,T}\right) \\
\rho^{1/2}\hat{\mathbb{A}}^{T_1}_\xi\rbtscoeff{n,T_1}-
\rho^{1/2}\sum_{j=1}^{J_{T_1}}\rbtscoeff{n,H,i_{T_1,j}}\hat{\mathbb{G}}^{T_1}_{\xi,j} \\
\vdots\\
\rho^{1/2}\hat{\mathbb{A}}^{\Tlast}_\xi\rbtscoeff{n,\Tlast}-
\rho^{1/2}\sum_{j=1}^{J_{\Tlast}}\rbtscoeff{n,H,i_{\Tlast,j}}\hat{\mathbb{G}}^{\Tlast}_{\xi,j}
\end{pmatrix},
\end{equation}
for $1 \leq n \leq N$.
After having computed~\eqref{eq:two_scale_est_gen_in_coord_f} and~\eqref{eq:two_scale_est_gen_in_coord_b}, we compute a $\underline{\mathfrak{W}}^{rb}$-orthonormal basis $\rbesttscoeff{m}$ for these vectors with coefficients
\begin{equation*}
\rbesttscoeff{m} = (\rbesttscoeff{m,H}, \rbesttscoeff{m,T_1}, \dots, \rbesttscoeff{m,\Tlast}),
\end{equation*}
such that 
\begin{equation*}
\rbestts{m} :=
(\sum_{i=1}^{\dimVH}\rbesttscoeff{m,H,i}\cdot\phi_i,
\sum_{l=1}^{M_{T_1}}\rbesttscoeff{m,T_1,l}\cdot\rbestt{l},
\dots,
\sum_{l=1}^{M_{\Tlast}}\rbesttscoeff{m,\Tlast,l}\cdot\rbestt{l}
)
\end{equation*}
is an $\Snorm{}$-orthonormal orthonormal basis for $\Wtsrb$.

Finally, following the definitions again, we see that $\hat{\mathbb{A}}_\xi$ and $\hat{\mathbb{F}}$ can be computed as
\begin{equation*}
\hat{\mathbb{F}}_m = \Fts(\rbestts{m}) = \sum_{i=1}^{\dimVH}\rbesttscoeff{m,H,i}\cdot l(\phi_i),
\end{equation*}
and
\begin{align*}
(\hat{\mathbb{A}}_\xi)_{m,n} 
&= \Btsq(\rbts{n}, \rbestts{m}) \\
&=
\begin{multlined}[t]
\sumT  \rbesttscoeff{m,H}^T \cdot (\Ktqo \cdot \rbtscoeff{n,H} - \Ktqrb \cdot \rbtscoeff{n,T}) \\
+ \rho^{1/2}\sumT\left( (\rbesttscoeff{m,T})^T\cdot\hat{\mathbb{A}}^T_\xi\cdot\rbtscoeff{n,T}
- \sum_{j=1}^{J_{\Tlast}}(\rbesttscoeff{m,T})^T\cdot\hat{\mathbb{G}}^{T}_{\xi,j}\cdot\rbtscoeff{n,H,i_{T,j}}\right).
\end{multlined}
\end{align*}

\subsubsection{Basis generation}\label{sec:stage_2_bas_gen}
To build the reduced Stage~2 space $\Vtsrb$, we follow the same methodology as in \cref{sec:stage1_basisgen} and use a greedy-search procedure to iteratively extend $\Vtsrb$ until a given error tolerance $\varepsilon_2$ for the MOR error estimate $\eta_{a,\mu}(\utsmrb)$ is reached.
However, in contrast to \cref{sec:stage1_basisgen}, we will not use the full-order model~\eqref{eq:two_scale_PGLOD}, or equivalently~\eqref{eq:PG_in_matrices}, to compute solution snapshots, but rather its Stage~1 approximation, i.e., we solve
\begin{equation}\label{eq:stage_2_fom}
\Kmurb \cdot \uhkmcoeff = \mathbb{F}
\end{equation}
to determine the $V_H$-component of the solution, followed by solving the Stage~1 corrector ROMs~\eqref{eq:stage1_rom} for each $T \in \Gridh$ to determine the fine-scale components $\QQktmrb(\uhkm)$ of the two-scale solution snapshot.
The complete algorithm is given by \cref{alg:stage2_greedy}. 

Since we only extend $\Vtsrb$ with approximations of the true solution snapshots of the full-order model, note that \cref{alg:stage2_greedy} is no longer a weak greedy algorithm in the sense of~\cite{binev2011convergence}.
In particular, the model reduction error is generally non-zero even for parameters $\mu^*$ for which the corresponding Stage~1 snapshot has been added to $\Vtsrb$.
Thus, when $\varepsilon_1$ is chosen too large in comparison to $\varepsilon_2$, a single $\mu^*$ might be selected twice, causing Algorithm~\ref{alg:stage2_greedy} to fail.
In such a case, the individual Stage~1 errors for each $T \in \Gridh$ can be estimated to enrich further the Stage~1 spaces for which the error is too large.
We will not discuss such an approach in more detail here and instead note that, in practice, it is feasible to choose $\varepsilon_1$ small enough to avoid such issues (cf. \cref{sec:TSRBLOD_experiments}).
More so, for sufficiently small $\varepsilon_1$, we can expect the convergence rates of a weak, greedy algorithm with exact solution snapshots to be preserved by Algorithm~\ref{alg:stage2_greedy} up to the given target tolerance $\varepsilon_2$.

\begin{algorithm2e}
	\KwData{Training set $\ParamsTrain$, tolerance $\varepsilon_2$}
	\KwResult{$\Vtsrb$}
	$\Vtsrb \leftarrow \{0\}$\;
	\While{$\max_{\mu\in\ParamsTrain} \eta_{a,\mu}(\utsmrb) > \varepsilon_2$}{
		$\mu^* \leftarrow \argmax_{\mu^* \in \ParamsTrain}\eta_{a,\mu}(\utsmrb)$\;
		$\uhkm[\mu^*] \leftarrow \textnormal{solution of~\eqref{eq:stage_2_fom}}$\;
		$\Vtsrb \leftarrow \Span(\Vtsrb \cup \{(\uhkm[\mu^*], \QQktmrb[T_1][\mu^*](\uhkm[\mu^*]), \dots, \QQktmrb[\Tlast][\mu^*](\uhkm[\mu^*]))\})$\;
	}
	\caption{Weak greedy algorithm for the generation of $\Vtsrb$.}\label{alg:stage2_greedy}
\end{algorithm2e}

\subsection{Computational complexity}

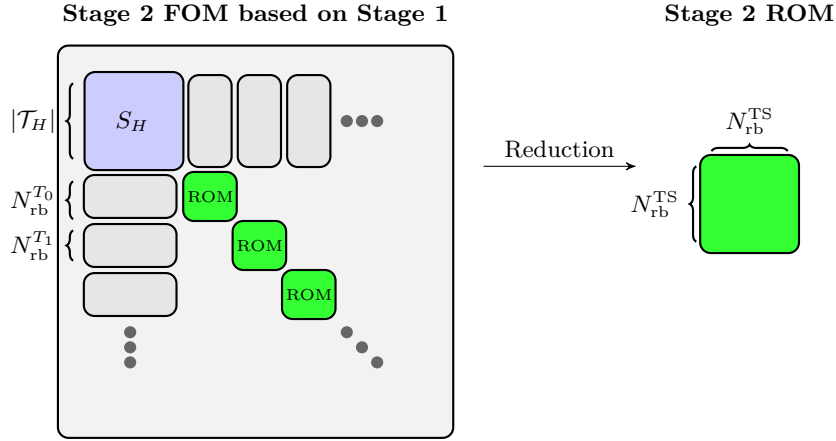
\begin{figure}[h]
	\centering \footnotesize
	\begin{tikzpicture}
	\node[Smatrixhuge] at (0,0) {};
	\node[Smatrixblue] at (-1.6,1.6) {$S_H$};
	
	\node[Smatrixgreen2] at (-0.6,0.6) {\tiny{\kern-0.2emROM\kern-0.2em}};
	\node[Smatrixgreen2] at (0.05,-0.05) {\tiny{\kern-0.2emROM\kern-0.2em}};
	\node[Smatrixgreen2] at (0.7,-0.7) {\tiny{\kern-0.2emROM\kern-0.2em}};
	
	\node[Smatrixblack] at (-0.6,1.6) {};
	\node[Smatrixblack] at (0.05,1.6) {};
	\node[Smatrixblack] at (0.7,1.6) {};
	
	\node[Smatrixblack2] at (-1.65, 0.6) {};
	\node[Smatrixblack2] at (-1.65,-0.05) {};
	\node[Smatrixblack2] at (-1.65,-0.7) {};
	
	\fill[fill=black!60] (1.2,-1.2) circle [radius=0.08];
	\fill[fill=black!60] (1.4,-1.4) circle [radius=0.08];
	\fill[fill=black!60] (1.6,-1.6) circle [radius=0.08];
	
	\fill[fill=black!60] (-1.65,-1.2) circle [radius=0.08];
	\fill[fill=black!60] (-1.65,-1.4) circle [radius=0.08];
	\fill[fill=black!60] (-1.65,-1.6) circle [radius=0.08];
	
	\fill[fill=black!60] (1.2,1.6) circle [radius=0.08];
	\fill[fill=black!60] (1.4,1.6) circle [radius=0.08];
	\fill[fill=black!60] (1.6,1.6) circle [radius=0.08];
	
	\draw [->] (3,1) -- (5,1) node[midway, above, thick] {Reduction};
	\node at (0,3) {\textbf{Stage 2 FOM based on Stage 1}};
	\node at (6.5,3) {\textbf{Stage 2 ROM}};
	
	\node[Smatrixgreen] at (6.5,0.5) {};
	\draw [decorate, decoration = {brace}, thick] (5.8,0) --  (5.8,1)
	node[pos=0.5,left=3pt,black]{$N^{\textnormal{TS}}_\textnormal{rb}$};
	\draw [decorate, decoration = {brace}, thick] (6,1.2) --  (7,1.2)
	node[pos=0.5,above=3pt,black]{$N^{\textnormal{TS}}_\textnormal{rb}$};
	
	\draw [decorate, decoration = {brace}, thick] (-2.4,1.1) --  (-2.4,2.1)
	node[pos=0.5,left=3pt,black]{$\abs{\mathcal{T}_H}$};
	
	\draw [decorate, decoration = {brace}, thick] (-2.4,0.3) --  (-2.4,0.8)
	node[pos=0.5,left=3pt,black]{$N^{T_0}_{\textnormal{rb}}$};
	
	\draw [decorate, decoration = {brace}, thick] (-2.4,-0.25) --  (-2.4,0.15)
	node[pos=0.5,left=3pt,black]{$N^{T_1}_{\textnormal{rb}}$};
	
	\end{tikzpicture}
	\caption{Visualization of the Stage~2 reduction process of the TSRBLOD.}
	\label{fig:TSRBLOD}
\end{figure}

Compared to~\cite{RBLOD}, the Stage~2 TSRBLOD ROM is not only independent from $\grid$ but also from the number of coarse-mesh elements in $\Gridh$.
Its size only depends on the number of selected basis vectors in \cref{sec:stage_2_bas_gen}.
The last reduction process in Stage~2 can be visualized via \cref{fig:TSRBLOD}.
The left matrix of \cref{fig:TSRBLOD} works as a certified replacement for the prohibitive system matrix of \eqref{eq:two_scale_PGLOD}.
Importantly, due to the dependency on $\Gridh$ and the local ROM sizes, also the left system matrix in \cref{fig:TSRBLOD} does not have to be assembled at all.
Instead, we only need the corresponding vectorized structure to attain the vector-matrix multiplication with the snapshots, which can be efficiently computed with the RBLOD methodology built from Stage~1.

\section{Software and implementational design}
\label{sec:software_TSRBLOD}

All experiments have been implemented in \texttt{Python} using \texttt{gridlod}~\cite{gridlod} for the PG--LOD discretization and \texttt{pyMOR}~\cite{milk2016pymor} for the model order reduction.
The complete source code, including setup instructions, is referenced in \cref{apx:code_availability}.
While \texttt{pyMOR} has already been detailed in \cref{sec:software_TRRB}, in this section, we first present the basics of \texttt{gridlod} and subsequently explain extensions that were necessary to conduct the numerical experiments of this section.

\subsection{PG--LOD with \texttt{gridlod}}
\label{sec:gridlod}

The \texttt{Python} module \texttt{gridlod}~\cite{gridlod} implements the Petrov-Galerkin version of the LOD for a quadrilateral mesh on the hypercube $[0,1]^d$.
In what follows, we briefly describe the different modules.
The library does not have a release system yet. Hence, we point to the respective version as it is referenced in \cite{tsrblod}.
For more information on the computational details of the LOD, we again refer to \cref{sec:comp_comp_LOD} and, in particular, to \cref{def:comp_proc_LOD}.

{\setlist[description]{labelindent=10pt,style=multiline,leftmargin=2.7cm}
\begin{description}
	\item[\texttt{world}] Contains the \texttt{world}-class for global information on the fine- and coarse grid, and the boundary conditions.
	\item[\texttt{coef}] Contains the code for localizing data on fine-scale patches.
	\item[\texttt{util}] Contains technical mapping functions for the required index arithmetic on element patches.
	\item[\texttt{interp}] Contains the interpolation operator for the corrector saddle point problem.
	\item[\texttt{linalg}] Provides code for solving linear systems appearing in Galerkin methods.
	\item[\texttt{fem}] Contains code for the finite element matrices.
	\item[\texttt{femsolver}] Code for computing FEM reference solutions.
	\item[\texttt{lod}] Contains the code for the corrector problems and the computation of the respective element correctors and the error indicators (used in \cref{sec:adaptive_PGLOD_from_hekema}).
	\item[\texttt{pglod}] Contains the code for the final PG--LOD, e.g., the multiscale stiffness matrix~$\Kmu$. 
\end{description}
}

Detailed examples for using \texttt{gridlod} are available in~\cite{gridlod}.
In particular, we note that the resulting user-code for computing the LOD is based on functions that are separately called in a loop over the coarse elements $T$.
Thus, given a suitable parallel system, Step~1 in \cref{fig:FOM_LOD} can be computed in parallel.
It is also optional whether only the coarse matrices $\Ktmu$ are returned or if the corrector solutions are also desired (as long as they fit into the memory).
For further details, we refer to~\cite{gridlod} and the code that we used for \cite{HKM20,tsrblod}, referenced in \cref{apx:code_availability}.

\subsection{Using \texttt{gridlod} in \texttt{pyMOR}}

To carry out the experiments for \cite{tsrblod}, we developed bindings-code for \texttt{gridlod} in order to use it within \texttt{pyMOR}.
In contrast to the TR-RB implementation that we explained in \cref{sec:software_TRRB}, we started with an existing \texttt{gridlod} discretization and wrote suitable \texttt{Model}s for the parts of the LOD where we intend to employ model order reduction for.
We call the resulting implementation the \texttt{rblod} module.
The respective objects are explained in the sequel.
We again refer to \cref{def:comp_proc_LOD} for the computational details (Steps~1--3) of the LOD:

\begin{description}
	\item[Stage~1-\texttt{Model}:] For the computation of Step~1 and the corresponding reduction process as detailed in \cref{sec:stage1}, for any $T \in \mathcal{T}_H$, we use a constructor that builds a respective \texttt{Model} to implement the corrector problem.
	The resulting \texttt{CorrectorProblem{\textunderscore}for{\textunderscore}all{\textunderscore}rhs} inherits from \texttt{pyMOR}'s \texttt{StationaryModel} and is mainly initialized by the respective \texttt{patch} object.
	Let us also mention that the fine-scale data of, for instance, $A_\mu$ can entirely be constructed locally; cf.~\cref{sec:exp_RBLOD_2}.
	The corresponding output functional of the \texttt{Model} is the contribution $\Ktmu$ to the coarse system matrix $\Kmu$.
	With the \texttt{Model} at hand, we can use \texttt{pyMOR}'s standard \texttt{CoerciveRBReductor} to reduce the corrector problems.
	Moreover, the \texttt{weak{\textunderscore}greedy} algorithm can be used to construct the respective basis.
	
	We note that, as discussed in \cref{rem:reduction_flexibility}, the models for the reduction scheme that follows the original RBLOD is instead implemented, for every $T$ and shape function, in \texttt{SeparatedCorrectorProblem} and can also be reduced by \texttt{CoerciveRBReductor}.
	
	\item[Stage~2-\texttt{Model}:] For the Stage~2 FOM, as it is explained in \cref{sec:stage_2_red}, we require a respective \texttt{Model}, which we called the \texttt{Two{\textunderscore}Scale{\textunderscore}Problem}.
	For this \texttt{StationaryModel} it is crucial to keep the idea that the two-scale system as it is illustrated in \cref{fig:TSRBLOD}(left) is never used to compute the FOM solutions.
	Instead, the original idea of the RBLOD, cf.~\cref{fig:RBLOD} is followed.
	To keep the memory requirement of Stage~2 low, a specialized \texttt{CoerciveRBReductorForTwoScale} was needed, internally accounting for the block structure in \cref{fig:TSRBLOD}(left) and assembling the error estimator.
	Again, the \texttt{weak{\textunderscore}greedy} can be used to construct the reduced basis.
	\item[Parallelization:] We emphasize that the initialization and reduction process for the Stage~1 problems can be parallelized by construction, inheriting from the implementational design in \texttt{gridlod}.
In our case we made use of \texttt{pyMOR}'s builtin \texttt{MPI}-parallelization routines.
	\item[Optimized \texttt{numpy}-based models with few memory requirements:]
	Apart from the relatively optimized offline time, we stress that \texttt{pyMOR} lacks computational efficiency in the online phase when it comes to dedicated online timings, mainly because of the technical nested scopes and assertions that are hidden in \texttt{pyMOR}.
	Since, in the numerical experiments, we aimed at ruling out such effects, we implemented the \texttt{OptimizedNumpyModelStage1} class for Stage~1 and the \texttt{OptimizedTwoScaleNumpyModel} for Stage~2 that are no longer a \texttt{Model} in the \texttt{pyMOR}-sense.
	Certainly, methods like \texttt{solve} and \texttt{output} can still be used the same way.
	For calling these methods for a new parameter, these models assemble the left-hand side and right-hand side with a \texttt{numpy-einsum}-call and call a direct \texttt{numpy}-based solver.

	Moreover, both optimized models contain a \texttt{minimal{\textunderscore}object}-method, making it possible to store only above-mentioned \texttt{numpy-arrays} on disc for minimal access to the surrogate models.
\end{description}

For an elaborated showcase of the above-mentioned functionality, we again refer to \cref{apx:code_availability}.
We also note that the explained \texttt{rblod}-code has been further extended for carrying out the experiments in \cref{chap:TR_TSRBLOD}, where also TR-RB related code fragments and a dedicated \texttt{gridlod-discretizer} for \texttt{pyMOR}'s \texttt{StationaryProblem} were needed.

\section{Numerical experiments}\label{sec:TSRBLOD_experiments}
This section applies the new TSRBLOD reduction method to two test cases and evaluates its efficiency.
For the first smaller problem, we mainly investigate the MOR error and the performance of the certified error estimator.
The second large-scale problem assesses the computational speedup achieved by the TSRBLOD.
As mentioned above, we use an \texttt{MPI} distributed implementation to benefit from parallelization of the localized corrector problems, which gives significant speed-ups for all considered methods.
Our computations have been performed on an HPC cluster with 1024 parallel processes.

In both cases, we use structured 2D-grids $\grid$, $\Gridh$ with quadrilateral elements on the domain $\Omega=[0,1]^2$.
We specify the number of elements by $n_h \times n_h$ and $n_H \times n_H$ respectively. 
We use the interpolation operator from \cite[Example 3.1]{LODbook} and choose the oversampling parameter $\kl$ as the first integer to satisfy $\kl > |\log(H)|$.
For the evaluation of the error estimator $\eta_{a, \mu}$~\eqref{eq:stage2_est_a} we approximate $\CPG \approx \alpha^{1/2}\CI^{-1}$ and assume $C_{\mathcal{I}_H} \approx 1$.
As aforementioned, the overlapping constant can be explicitly computed for quadrilateral elements by $\CO = (2\kl + 1)^2$.
Moreover, we approximate $\alpha$ and the contrast $\kappa$ by replacing $\mathcal{P}$ in \eqref{eq:ellipticity_1} and \eqref{eq:ellipticity_2} by the training set $\ParamsTrain$.
Thus, while it is not guaranteed that our approximations yield strict upper bounds on the model order reduction error, the decay rate and choice of snapshot parameters will not be affected.

\begin{description}
\item[Petrov--Galerkin variant of the RBLOD method]
In order to compare our method to the RBLOD approach introduced in~\cite{RBLOD}, we have implemented a corresponding version of the RBLOD that applies to our LOD formulation.
In particular, we use the interpolation operator from \cite[Example 3.1]{LODbook} and Petrov-Galerkin projection~\eqref{eq:PG_} in contrast to the Cl\'ement interpolation and Galerkin projection used in~\cite{RBLOD}.
Further, in~\cite{RBLOD}, individual ROMs for the local correctors $\QQktm(\phi_i)$ are constructed, where the snapshot parameters are chosen identically among all $\QQktm[T^\prime](\phi_i)$, $T^\prime \subset \operatorname{supp} \psi_z$.
In the RBLOD variant implemented by us, we independently train the corrector ROMS for each $T \in \Gridh$ as is done in Stage~1 of the TSRBLOD.
We note that for both RBLOD variants, the number of ROMs constructed for each coarse element $T$ is equal to the number of coarse-mesh basis functions $\phi_i$ supported on this element, whereas Stage~1 of the TSRBLOD constructs a single ROM to approximate all local correctors; cf.~\cref{rem:reduction_flexibility}
\end{description}

\paragraph{Error measures}
To quantify the accuracy of the TSRBLOD, we use a validation set $\ParamsVer \subset \mathcal{P}$ of $10$ random parameters and compute the maximum relative approximation errors w.r.t.\ the PG--LOD solutions for this set, i.e.,
\begin{align*}
e^{*, \textnormal{rel}}_{\textnormal{LOD}} := \max_{\mu\in\ParamsVer} \frac{\| \uhkm - \tilde{u}_{\mu} \|_{*} }{\|\uhkm\|_{*}}, 
\end{align*}
where $\tilde{u}_{\mu}$ denotes the coarse-scale component of either the RBLOD or TSRBLOD solution, $\uhkm$ denotes the coarse component of the PG--LOD solution~\eqref{eq:PG_}, and $^*$ stands for either the $H^1$- or $L^2$-norm.
With the solution of the FEM approximation $u_{h,\mu}$ w.r.t.~$\grid$, we let
\begin{align*}
e^{L^2, \textnormal{rel}}_{\textnormal{FEM}} &:=  \max_{\mu\in\ParamsVer} \frac{\| u_{h,\mu} - \tilde{u}_{\mu}\|_{L^2} }{\|u_{h,\mu}\|_{L^2}}, \qquad
&e^{L^2, \textnormal{rel}}_{\textnormal{LOD-FEM}} &:=  \max_{\mu\in\ParamsVer} \frac{\| u_{h,\mu} - \uhkm \|_{L^2} }{\|u_{h,\mu}\|_{L^2}}.
\end{align*}

\paragraph{Time measures}
For an analysis of the computational wall times, we consider the following quantities:
\begin{itemize}[leftmargin=*]
	\item $t^\textnormal{offline}_{1,\textnormal{av}}(T)$: Average (arithmetic mean) time for creating a ROM for the localized corrector problem(s) corresponding to a single coarse element $T \in \Gridh$ as discussed in \cref{sec:stage1}.
	For the RBLOD, this involves all individual corrector problems for the four basis functions that are supported on $T$.
	\item $t^\textnormal{offline}_{1}$: Total time required to build the corrector ROMs.
	In the case of full parallelization, this is equal to the maximum of the Stage~1 times over all $T \in \Gridh$.
	\item $t^\textnormal{offline}_2$: Time for creating the Stage~2 TSRBLOD ROM as discussed in \cref{sec:stage_2_red}.
	\item $t^{\textnormal{offline}}$: Total offline time for building the final reduced model.
	For the RBLOD this is equal to $t^\textnormal{offline}_{1}$. For the TSRBLOD this additionally includes $t^\textnormal{offline}_2$.	
	\item $t^{\textnormal{online}}$: Average time for solving the obtained ROM for a single new $\mu \in \ParamsVer$.
	In the case of the RBLOD, the reduced corrector problems are solved sequentially on a single compute node.
	\item $t^{\textnormal{LOD}}$: Average time needed to compute the PG--LOD with parallelization of the corrector problems for a single new $\mu \in \ParamsVer$.
\end{itemize}

\subsection{Experiment 10: Test case of Section 4.1 in \cite{RBLOD}}
\label{sec:exp_RBLOD_1}
\begin{figure}[h]
	\centering
	\includegraphics[width=2.in]{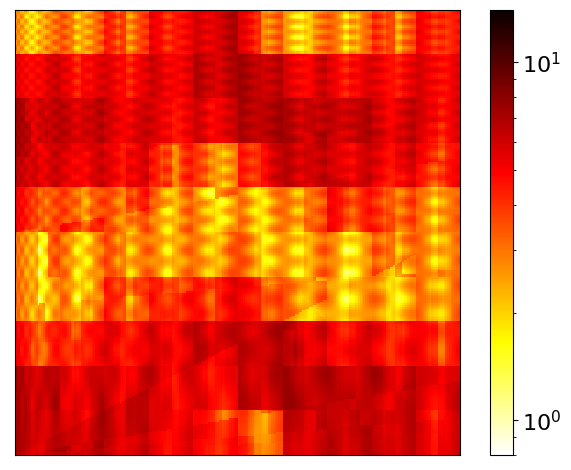}
	\includegraphics[width=2.in]{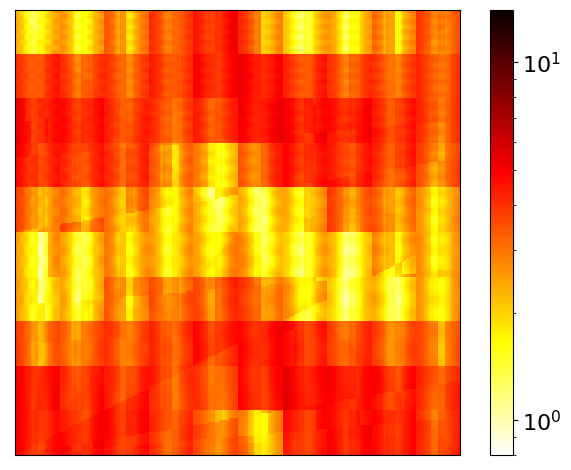}
	\includegraphics[width=2.in]{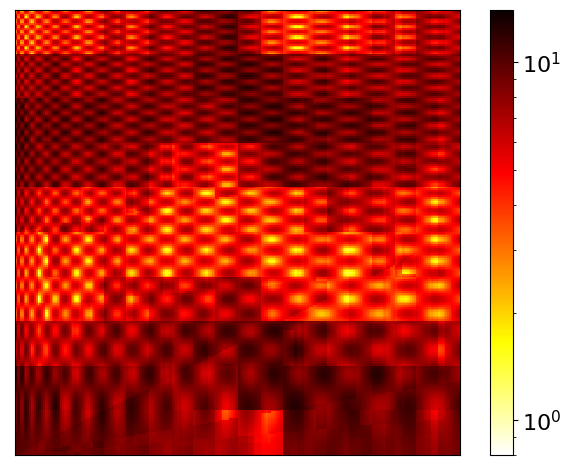}
	\captionof{figure}[Diffusion coefficient $A_\mu$ for different $\mu$.]{Diffusion coefficient $A_\mu$ for $\mu_1 = 1.8727$ (left), $\mu_2 = 2.9040$ (middle), and $\mu_3 = 4.7536$ (right) for Experiment 10.}
	\label{mp1:Pictures_TSRBLOD}
	\centering
\end{figure}

The first test case is taken from~\cite[Section 4.1]{RBLOD} and has a one-dimensional parameter space $\Params:= [0,5]$. The coefficient $A_\mu$ is visualized in \cref{mp1:Pictures_TSRBLOD}, and we set $f \equiv 1$.
For the sake of brevity, we refer to~\cite{RBLOD} for an exact definition of $A_\mu$.
In order to resolve the microstructure of the problem, we choose $n_h = 2^{8}$ which results in $65,536$ fine-scale elements.
The approximated maximum contrast of $A_\mu$ is $\kappa \approx 13$.

\begin{table}
	\footnotesize
	\centering
	{\def\arraystretch{1.7}\tabcolsep=3pt
		\begin{tabular}{|l|c|c|c|c|c|c|} \hline
			Mesh size $n_H$& \multicolumn{2}{c|}{$2^{3}$} & \multicolumn{2}{c|}{$2^{4}$} &\multicolumn{2}{c|}{$2^{5}$} \\ \hline
			Method  &\hphantom{T}RBLOD\hphantom{T} & TSRBLOD&\hphantom{T}RBLOD\hphantom{T} & TSRBLOD 
			&\hphantom{T}RBLOD\hphantom{T} & TSRBLOD \\ \hline  \hline
			$t^\textnormal{offline}_{1, \textnormal{av}}(T)$ 	&	41		&  	61		&  	39		& 	61		&   33		&	55		\\
			$t^\textnormal{offline}_{1}$   			 	&	71		&  	106		&  	67		& 	102		&   63 		&	98		\\
			$t^\textnormal{offline}_2$   					&	- 		& 	8 		&  	-		& 	56  	&	- 		&	472	\\
			$t^\textnormal{offline}$    					&	71		& 	114 	&  	67		& 	158		&   63		&	570	\\ \hline
			Cum. size Stage 1 							&	2346	&  	1670	&	8718	& 	6134	&	31810 	&	22189	\\
			Av. size Stage 1           				&	9.16	&  	26.09	& 	8.51	& 	23.96	&	7.77 	&	21.67	\\
			Size Stage 2 								&	-		&  	8		& 	-		& 		9	&		 	&	9		\\
			\hline \hline
			$t^{\textnormal{LOD}}$ 
			& 	\multicolumn{2}{c|}{0.69} 	
			& 	\multicolumn{2}{c|}{0.49}	
			& 	\multicolumn{2}{c|}{0.90}	 \\ \hline
			$t^\textnormal{online}$ 		 				&	0.0610	&  	0.0003	&  	0.2272	& 	0.0003	&   1.0462 	&	0.0003	 \\ \hline \hline
			Speedup LOD 	 						&	11		&  	2506	& 	2 		& 	1536 	&   1		& 	2714	 \\ \hline \hline
			$e^{H^1, \textnormal{rel}}_{\textnormal{LOD}}$ 	 	&	1.97e-5	&	7.30e-4	&	5.08e-5 &	2.94e-4 &	1.11e-4 &	4.21e-4	 \\
			$e^{L^2, \textnormal{rel}}_{\textnormal{LOD}}$ 	 	&	4.89e-6 &	2.71e-4 &	6.77e-6 &	1.03e-4 &	7.70e-6	&	1.32e-4  \\
			\hline \hline
			$e^{L^2, \textnormal{rel}}_{\textnormal{FEM}}$ 		&	2.46e-2 &	2.46e-2	&	9.05e-3	&	9.05e-3 &	3.98e-3	&	3.98e-3 \\ \hline
			$e^{L^2, \textnormal{rel}}_{\textnormal{LOD-FEM}}$	 
			& 	\multicolumn{2}{c|}{2.46e-2} 
			& 	\multicolumn{2}{c|}{9.05e-3} 
			& 	\multicolumn{2}{c|}{3.98e-3}\\ \hline
	\end{tabular}}
	\captionsetup{width=\textwidth}
	\caption[Experiment 10: Performance, ROM sizes and accuracy of the methods]{\footnotesize{%
			Performance, ROM sizes and accuracy of the methods for test case 1 with tolerances $\varepsilon_1 = 0.001$ and $\varepsilon_2 = 0.01$ and varying coarse-mesh sizes. All times are given in seconds.	\label{Tab:mp1:1}}}
\end{table}

\subsubsection{Performance and error comparison}

We vary the coarse mesh size, $n_H \in \{2^{3},2^{4},2^{5}\}$, and show results for $\varepsilon_1 = 0.001$, $\varepsilon_2=0.01$ and a training set $\ParamsTrain \subset \Params$ of $50$ equidistant parameters in \Cref{Tab:mp1:1}.

We observe that this choice of tolerances suffices for both the RBLOD and TSRBLOD to match the error of the PG--LOD solution w.r.t.\ the FEM reference solution.
The TSRBLOD produces a ROM with 8 or 9 basis functions without losing accuracy for all coarse mesh sizes.
The offline phase of Stage~1 for the TSRBLOD is longer and yields larger ROMS per element than
the RBLOD.
However, for the RBLOD four ROMS per element are required, such that the total number of basis vectors is smaller for the TSRBLOD.
The offline time for Stage~2 increases for larger $n_H$ due to the increasing complexity of assembling and solving the coarse-scale problem, even with an RB approximation of the corrector problems.
However, this only affects the offline phase for the TSRBLOD, and the online times are at a constant level since the dimension of the Stage~2 ROM is largely unaffected by the number of coarse elements.
Considering the speedups achieved by both methods, it can be seen that the RBLOD has just a slight benefit over the (parallelized) PG--LOD. In contrast, the TSRBLOD shows a significant speedup for all $n_H$ without requiring parallelization.

\subsubsection{Error and estimator decay}
\label{sec:RBLOD_estimator_study}

\begin{figure}[t]
	\footnotesize
	\centering
	\begin{tikzpicture}
	\definecolor{color0}{rgb}{0.65,0,0.15}
	\definecolor{color1}{rgb}{0.84,0.19,0.15}
	\definecolor{color2}{rgb}{0.96,0.43,0.26}
	\definecolor{color3}{rgb}{0.99,0.68,0.38}
	\definecolor{color4}{rgb}{1,0.88,0.56}
	\definecolor{color5}{rgb}{0.67,0.85,0.91}
	\begin{axis}[
	name=left,
	anchor=west,
	width=8.5cm,
	height=5.5cm,
	log basis y={10},
	tick align=outside,
	tick pos=left,
	legend style={nodes={scale=0.7}, fill opacity=0.8, draw opacity=1, text opacity=1, 
	},
	x grid style={white!69.0196078431373!black},
	xlabel={basis enrichments},
	ylabel={max estimator $\eta_{a, \mu}$},
	xmajorgrids,
	xtick style={color=black},
	y grid style={white!69.0196078431373!black},
	ymajorgrids,
	ymode=log,
	ymin=1e-04, ymax=2,
	ytick style={color=black},
	ytick pos=left,
	]
	\addplot [semithick, color0, mark=triangle*, mark size=3, mark options={solid, fill opacity=0.5}]
	table {%
		1 1.302208490599 
		2 0.542298633478 
		3 0.169259772169
	};
	\addlegendentry{$\varepsilon_1 = 0.1\hphantom{000}$}
	\addplot [semithick, color3, mark=*, mark size=3, mark options={solid, fill opacity=0.5}]
	table {%
		2 0.548045241003 
		3 0.140901484479 
		4 0.104570132829 
		5 0.058864107218 
		6 0.034589798089 
		7 0.023123372984
	};
	\addlegendentry{$\varepsilon_1 = 0.01\hphantom{00}$}
	\addplot [semithick, color5, mark=square*, mark size=3, mark options={solid, fill opacity=0.5}]
	table {%
		6 0.032599625814 
		7 0.014522966903 
		8 0.010109952521 
		9 0.003802794378 
		10 0.003003471958 
		11 0.001834614254 
		12 0.001739345745 
		13 0.001618788041
	};
	\addlegendentry{$\varepsilon_1 = 0.001\hphantom{0}$}
	\addplot [semithick, color1, mark=diamond*, mark size=3, mark options={solid, fill opacity=0.5}]
	table {%
		10 0.002914232833 
		11 0.001344315689 
		12 0.000991445067 
		13 0.000905171861 
		14 0.000392978852 
		15 0.000316746807 
		16 0.000214242687
	};
	\addlegendentry{$\varepsilon_1 = 0.0001$}
	\end{axis}
	\end{tikzpicture}
	\captionsetup{width=\textwidth}
	\caption[Experiment 10: Evolution of the maximum estimated error $\eta_{a, \mu}$, $\mu \in \ParamsTrain$]{\footnotesize{Evolution of the maximum estimated error $\eta_{a, \mu}$, $\mu \in \ParamsTrain$ for different Stage~1 tolerances $\varepsilon_1$ during the greedy algorithm of Stage~2 ($n_H = 2^4$). The greedy algorithm has been continued until the enrichment failed. The first values for smaller tolerances are left out to improve readability.}}
	\label{Fig:estimator_study}
\end{figure}
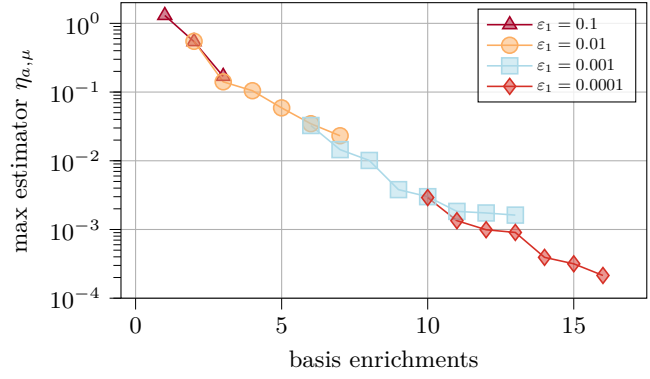

\begin{figure}[t]
	\footnotesize
	\centering
	\begin{tikzpicture}
	\definecolor{color0}{rgb}{0.65,0,0.15}
	\definecolor{color1}{rgb}{0.84,0.19,0.15}
	\definecolor{color2}{rgb}{0.96,0.43,0.26}
	\definecolor{color3}{rgb}{0.99,0.68,0.38}
	\definecolor{color4}{rgb}{1,0.88,0.56}
	\definecolor{color5}{rgb}{0.67,0.85,0.91}
	\begin{axis}[
	name=left,
	anchor=west,
	width=7.5cm,
	height=5.5cm,
	log basis y={10},
	tick align=outside,
	legend cell align={left},
	legend style={nodes={scale=0.7}, fill opacity=0.8, draw opacity=1, text opacity=1, at=(left.west), anchor=north east, xshift=-0.2cm, yshift=0.7cm, draw=white!80!black},
	x grid style={white!69.0196078431373!black},
	xlabel={basis enrichments},
	xmajorgrids,
	xtick style={color=black},
	y grid style={white!69.0196078431373!black},
	ymajorgrids,
	ymode=log,
	ymin=5e-08, ymax=15,
	ytick style={color=black},
	ytick pos=right
	]
	\addplot [semithick, color0, mark=x, mark size=4, mark options={solid, fill opacity=0.5}]
	table {%
		1 0.20583663158433416
		2 0.11887146368778598
		3 0.029783870160454723
		4 0.01489840725378891
		5 0.014067080101043957
		6 0.007252894748508427
		7 0.003190748490468686
		8 0.0016313851696358674
		9 0.0007856473067503389
		10 0.0007520128559338216
		11 0.00038889201711461214
		12 0.00038219438646071996
		13 0.00034000389283784436
	};
	\addplot [semithick, black, mark=diamond*, mark size=3, mark options={solid, fill opacity=0.2}]
	table {%
		1 0.013304715193246983
		2 0.009712437622371539
		3 0.002816468173892014
		4 0.0011318376672740094
		5 0.0010401116607180129
		6 0.0006290529115014175
		7 0.00014645370454404016
		8 0.00008648122221639462
		9 0.000045722930894291935
		10 0.00004620880598336765
		11 0.000021310199120887216
		12 0.00001932049892755662
		13 0.00001260466561162649
	};
	\addplot [semithick, color1, mark=square*, mark size=2, mark options={solid, fill opacity=0.2}]
	table {%
		1 0.01167264689554406
		2 0.008934578133502196
		3 0.002677895817283881
		4 0.000982896580291272
		5 0.0008893416590385863
		6 0.0005798216284228025
		7 6.493345512553352e-05
		8 5.704378512652952e-05
		9 3.59582145370613e-05
		10 3.677364143091256e-05
		11 1.7109809385720184e-05
		12 1.2277935969088866e-05
		13 2.4157504039053084e-06
	};
	\addplot [semithick, black, mark=square*, opacity=0.4, mark size=1.5, mark options={solid, fill opacity=0.5}]
	table {%
		1 0.20556287671660609
		2 0.11853521920869325
		3 0.029663239804960596
		4 0.014895622925657302
		5 0.014041884984912538
		6 0.00722968112105605
		7 0.0031902180521013682
		8 0.0016303875546281278
		9 0.0007852344218980461
		10 0.0007511131970519845
		11 0.00038867189820279414
		12 0.00038202110052247513
		13 0.00034000203935049633
	};
	
	\addplot [semithick, color0, mark=*, mark size=3, dashed, mark options={solid, fill opacity=0.5}]
	table {%
		1 1.309575235
		2 0.548060703
		3 0.141635678
		4 0.104847183
		5 0.059237218
		6 0.032599626
		7 0.014522967
		8 0.010109953
		9 0.003802794
		10 0.003003472
		11 0.001834614
		12 0.001739346
		13 0.001618788
	};
	
	\addplot [semithick, color5, mark=square*, mark size=3, mark options={solid, fill opacity=0.5}]
	table {%
		1 7.12906094171107
		2 7.3606022714764565
		3 6.969506743881207
		4 7.084436819956133
		5 6.662971068069086
		6 6.681572675067291
		7 6.359053731851468
		8 6.567192608633759
		9 6.530567472242883
		10 6.6733601290328055
		11 6.676811460567216
		12 6.59853821648667
		13 6.544061489519691
	};
	\addplot [semithick, color3, mark=triangle*, dashed, mark size=3, mark options={solid, fill opacity=0.5}]
	table {
		1 0.064086520
		2 0.033603723
		3 0.012175896
		4 0.004836665
		5 0.004318154
		6 0.002758986
		7 0.000401127
		8 0.000339618
		9 0.000153555
		10 0.000157558
		11 0.000087057
		12 0.000054408
		13 0.000013076
	};
	\end{axis}
	\begin{axis}[
	name=right,
	anchor=west,
	at=(left.east),
	xshift=1.1cm,
	width=7.5cm,
	height=5.5cm,
	log basis y={10},
	legend cell align={left},
	legend style={nodes={scale=0.7}, fill opacity=0.8, draw opacity=1, text opacity=1, at=(left.east), anchor=north east, xshift=-13.5cm, yshift=-1.4cm, draw=white!80!black},
	tick align=outside,
	tick pos=left,
	x grid style={white!69.0196078431373!black},
	xlabel={basis enrichments},
	xmajorgrids,
	ymode=log,
	ymin=5e-08, ymax=15,	
	xtick style={color=black},
	y grid style={white!69.0196078431373!black},
	ymajorgrids,
	yticklabels={,,},
	ytick pos=left,
	]
	\addplot [semithick, color0, mark=x, mark size=4, mark options={solid, fill opacity=0.5}]
	table {
		1 0.20583663158433418
		2 0.11694206395586734
		3 0.02973104282919254
		4 0.0187861354798649
		5 0.013858503372350388
		6 0.012416870025239479
		7 0.0026686216441154298
		8 0.002635502193860233
		9 0.0007400980127743919
		10 0.0005443247667907771
		11 0.0004620402045593204
		12 0.0004127701959231501
		13 0.00034019560933311744
		14 0.00034018968282338157
		15 0.00034006575626278755
		16 0.00034005631077660875
		17 0.00034000677770389305
		18 0.00034000614761857314
		19 0.0003400047029756307
		20 0.00034000469392983946
		21 0.0003400046813423982
		22 0.0003400046754256699
		23 0.0003400046507017614
		24 0.00034000451666054874
		25 0.00034000442076877424
	};
	\addlegendentry{two-scale error}
	\addplot [semithick, black, mark=diamond*, mark size=3, mark options={solid, fill opacity=0.2}]
	table {
		1 0.013304715193246948
		2 0.009511015595946179
		3 0.0028113477767885784
		4 0.0010418776425366626
		5 0.000695950849759363
		6 0.0007067760148770652
		7 0.0001288287554918048
		8 0.00012289411897222592
		9 4.5151058847703486e-05
		10 2.9558799782869537e-05
		11 2.237760999975541e-05
		12 1.9830356206882944e-05
		13 1.2628238317123184e-05
		14 1.2625312273148166e-05
		15 1.2624222828842822e-05
		16 1.2622672901671814e-05
		17 1.2608270221566213e-05
		18 1.2607216790604492e-05
		19 1.2605206742112442e-05
		20 1.2605264311413092e-05
		21 1.2605117349898283e-05
		22 1.260517373978866e-05
		23 1.2605059145931392e-05
		24 1.2604687600214781e-05
		25 1.2605079504893575e-05		
	};
	\addlegendentry{LOD error}
	\addplot [semithick, color1, mark=square*, mark size=2, mark options={solid, fill opacity=0.2}]
	table {
		1 0.01167264689554405
		2 0.008727740382171471
		3 0.0026728177591846536
		4 0.0007560533227358547
		5 0.0005811156126797668
		6 0.0004809482874223695
		7 7.457221511493489e-05
		8 5.032735352846313e-05
		9 3.5721462437374986e-05
		10 2.3204288829405183e-05
		11 1.492326673635207e-05
		12 1.0182680403255702e-05
		13 2.0210071303287623e-06
		14 2.027257098624611e-06
		15 1.6766888163742802e-06
		16 1.1714618468972557e-06
		17 1.4449245167868556e-06
		18 1.12423583570912e-06
		19 1.1229860536982722e-06
		20 1.1230851612871708e-06
		21 1.1228418969885156e-06
		22 1.1229347412743485e-06
		23 1.122478451417174e-06
		24 1.1225804318428979e-06
		25 1.122916984414993e-06
	};
	\addlegendentry{$V_H$-error}
	\addplot [semithick, black, mark=square*, mark size=1.5, opacity=0.4, mark options={solid, fill opacity=0.5}]
	table {
		1 0.20556287671660609
		2 0.11661592031142051
		3 0.02961065607070313
		4 0.01877402382243087
		5 0.013853923304489067
		6 0.012407552134426684
		7 0.002667964814102145
		8 0.0026350286117542217
		9 0.0007392354466838959
		10 0.0005439894059769174
		11 0.00046195181513729793
		12 0.0004126528759191995
		13 0.0003401935686116118
		14 0.00034018763966547634
		15 0.00034006375702974264
		16 0.0003400542929829571
		17 0.00034000491152417437
		18 0.00034000428896149035
		19 0.0003400028484408183
		20 0.00034000283906762354
		21 0.00034000282728356853
		22 0.00034000282106018176
		23 0.0003400027978428308
		24 0.00034000266346419613
		25 0.00034000256646054275
	};
	\addlegendentry{corrector error}
	
	\addplot [semithick, color0, mark=*, mark size=3, dashed, mark options={solid, fill opacity=0.5}]
	table {
		1 1.309575235
		2 0.538279469
		3 0.143527566
		4 0.129353986
		5 0.096122047
		6 0.091427914
		7 0.013125357
		8 0.012644495
		9 0.002991285
		10 0.002984935
		11 0.002913277
		12 0.002555585
		13 0.001619900
		14 0.001619876
		15 0.001619153
		16 0.001619044
		17 0.001618832
		18 0.001618806
		19 0.001618797
		20 0.001618799
		21 0.001618796
		22 0.001618799
		23 0.001618800
		24 0.001618798
		25 0.001618796
	};
	\addlegendentry{$\eta_{a, \mu}$}		
	
	\addplot [semithick, color3, mark=triangle*, mark size=3, dashed, mark options={solid, fill opacity=0.5}]
	table {
		1 0.064086520
		2 0.032843009
		3 0.012154224
		4 0.005207466
		5 0.002802488
		6 0.003292533
		7 0.000361567
		8 0.000210742
		9 0.000153361
		10 0.000119555
		11 0.000072419
		12 0.000062428
		13 0.000008881
		14 0.000008243
		15 0.000007220
		16 0.000004308
		17 0.000007923
		18 0.000001845
		19 0.000001591
		20 0.000000897
		21 0.000000279
		22 0.000000216
		23 0.000000209
		24 0.000000103
		25 0.000000080
	};
	\addlegendentry{$\eta_{a, \mu}$, $\rho = 0$}
	\addplot [semithick, color5, mark=square*, mark size=3, mark options={solid, fill opacity=0.5}]
	table {
		1 7.129060941711069
		2 7.3656688584410945
		3 6.9824580069300275
		4 6.885609149899999
		5 7.333056748247623
		6 7.3632013268725025
		7 6.471238841821415
		8 6.566019043110956
		9 6.626142811888608
		10 6.624813706577451
		11 6.658917543462653
		12 6.556560885917927
		13 6.475931924186771
		14 6.483625325095644
		15 6.4830024301230464
		16 6.477027190303722
		17 6.482799619476211
		18 6.481894302317812
		19 6.481905020327957
		20 6.478519232309123
		21 6.4822406600512705
		22 6.482246722522054
		23 6.482690110890613
		24 6.4811915922350805
		25 6.481407203351824
	};
	\addlegendentry{effectivity $\eta_{a, \mu}$}
	
	\end{axis}
	\node[anchor=south east, yshift=4pt, xshift=-30pt] at (left.north east) {(a) train with $\eta_{a, \mu}$};
	\node[anchor=south west, yshift=4pt, xshift=25pt] at (right.north west) {(b) train with $\eta_{a, \mu}$, $\rho = 0$};
	\end{tikzpicture}
	\captionsetup{width=\textwidth}
	\caption[Experiment 10: Comparison of the Stage~2 training error decay using different error estimators]{\footnotesize{Comparison of the Stage~2 training error decay using different error estimators ($\varepsilon_1=0.001$).
	The greedy algorithm is continued until the enrichment fails or $25$ enrichments are reached.
	Depicted is the maximum value of the full two-scale error $\Anorm{\utsm - \uts^{rb}_\mu}$, the LOD error given by $\Anorm{\utsm - \uts^{rb}_\mu}$ with $\rho=0$, the $V_H$-error $\anorm{\uhkm - \uhkm^\textnormal{rb}}$, the fine-scale corrector error $\rho^{1/2} \cdot (\sumT \anorm{\QQktm(\uhkm^\textnormal{rb}) -  u^\textnormal{rb,f}_T}^2)^{1/2}$, the estimator $\eta_{a, \mu}$, its effectivity $\eta_{a, \mu}/\Anorm{\utsm - \uts^{rb}_\mu}$, and the part of $\eta_{a, \mu}$ corresponding to the LOD residual (obtained by setting $\rho=0$).
	The maximum is computed over the training set $\ParamsTrain$.
	The error estimator used in the greedy algorithm is either (a) the two-scale error estimator $\eta_{a, \mu}$ or (b) only its LOD-residual part ($\rho=0$).}}
	\label{Fig:estimator_study_2}
\end{figure}
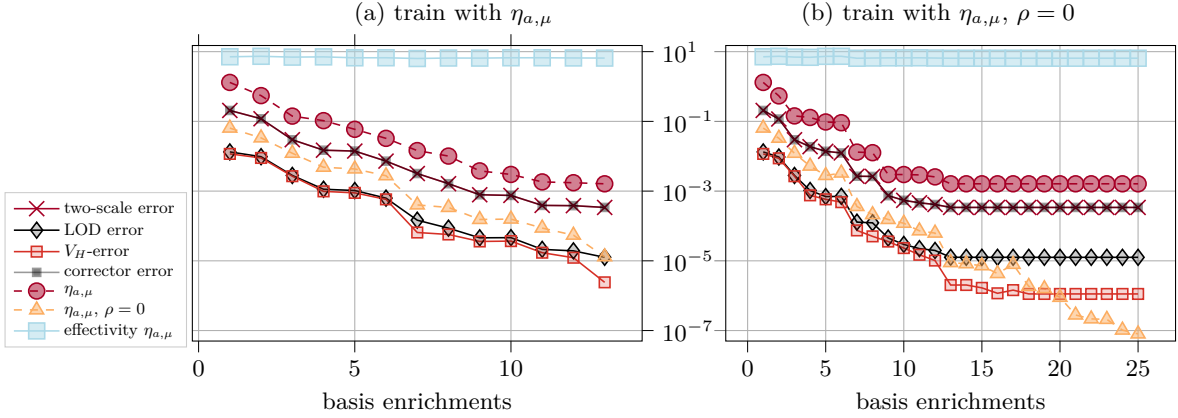

Next, we study the influence of the Stage~1 tolerance $\varepsilon_1$ on the training of the Stage~2 ROM.
To this end, we fix the number of coarse-mesh elements $n_H = 2^4$ and depict in \cref{Fig:estimator_study} for different $\varepsilon_1$ the maximum estimated training error w.r.t.\ the number of basis functions of the Stage~2 ROM.
For all $\varepsilon_1$, the Stage~2 training was continued until enrichment failed due to repeated selection of the same snapshot parameter (cf. \cref{sec:stage_2_bas_gen}).
As expected, a sufficiently small $\varepsilon_1$ is required to achieve small errors for the Stage~2 ROM.
Note, however, that choosing a smaller $\varepsilon_1$ for the same $\varepsilon_2$ only affects the offline time of the Stage~2 training, but not the efficiency of the resulting Stage~2 ROM.

In \cref{Fig:estimator_study_2}(a), we study in more detail how the two-scale error $\Anorm{\utsm - \utsm^\textnormal{rb}}$ and its estimator $\eta_{a, \mu}$ are affected by the coarse- and fine-scale errors in the two-scale system.
We observe that the greedy algorithm aborts when the error in the fine-scale correctors stagnates at the lower bound determined by the fixed Stage~1 ROMs, which are used to generate the Stage~2 solution snapshots.\@
While the LOD error ($\Anorm{}$ with $\rho = 0$) and the corresponding residual norms ($\eta_{a,\mu}$ with $\rho = 0$) decrease over all iterations, the corrector residuals dominate~$\eta_{a,\mu}$, finally causing the same snapshot parameter to be selected twice.

To verify that the enrichment procedure should, indeed, be stopped at this point, we perform another experiment where we neglect the corrector residuals and use $\eta_{a, \mu}$ with $\rho = 0$ as the error surrogate in the greedy algorithm.
This corresponds to treating the RBLOD coarse system as the FOM w.r.t.\ which the MOR error is estimated.
The result is visualized in \cref{Fig:estimator_study_2}(b).
The estimator rapidly decays over all 25 enrichments but both the two-scale and LOD errors do not decrease any further.
This is expected as both error measures involve the error in the fine-scale correctors.
Admittedly, also the $V_H$-error $\anorm{\uhkm - \uhkm^\textnormal{rb}}$ stagnates after one further iteration, and the estimator eventually underestimates all three errors.
This underlines that the fine-scale corrector errors need to be considered, even if one is only interested in a coarse-scale approximation in $V_H$.

\subsection{Experiment 11: Large-scale example}
\label{sec:exp_RBLOD_2}
\begin{figure}[h]
	\begin{center}
		\begin{subfigure}[b]{0.38\textwidth}
			\includegraphics[width=\textwidth]{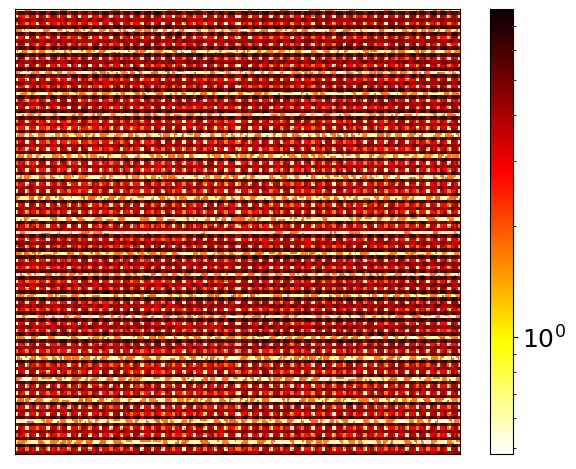}
		\end{subfigure}
	\end{center}
	\begin{subfigure}[b]{0.26\textwidth}
		\includegraphics[width=\textwidth]{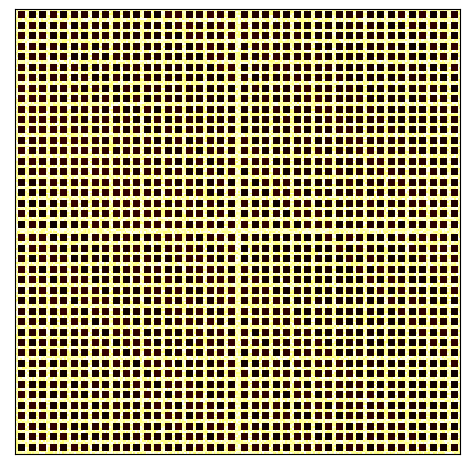}
	\end{subfigure}
	\begin{subfigure}[b]{0.26\textwidth}
		\includegraphics[width=\textwidth]{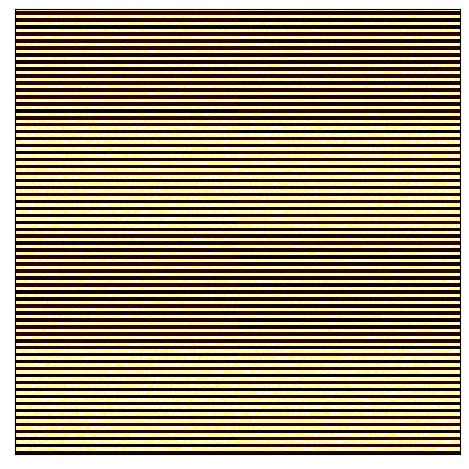}
	\end{subfigure}
	\begin{subfigure}[b]{0.328\textwidth}
		\includegraphics[width=\textwidth]{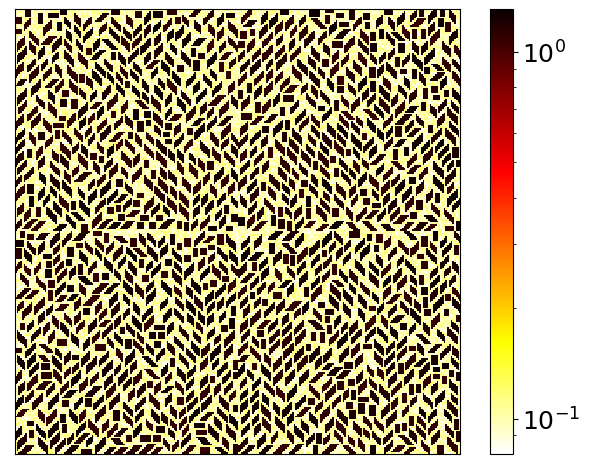}
	\end{subfigure}
	\centering
	\caption[Experiment 11: Coefficient $A_\mu$ on $4$ of $4096$ coarse-mesh elements]{\footnotesize{Coefficient $A_\mu$ on $4$ of $4096$ coarse-mesh elements for $\mu=(1,2,3)^T$ (top center) and $A_\xi$ for all $\xi=1,2,3$ (bottom from left to right).
			All coefficients $A_\xi$ are subjected to normally distributed noise in the interval $[1,1.2]$ for the particles (black) and in the interval $[0.03,0.11]$ for the background (yellow). To ensure reproducibility, for every coarse-mesh element $T$, we choose the random seed for the noise and for the distribution of the randomly shaped particles in $A_3$ as the global index of $T$ in $\Gridh$.}}
	\label{mp2:Pictures_TSRBLOD}
\end{figure}

We now consider a more complex test case with a significantly larger fine-scale mesh $\grid$.
The three-dimensional parameter space is given by $\Params := [1,5]^3$, and we let $A_\mu:=\sum_{\xi=1}^{3}\mu_\xi A_\xi$, where for a representative patch of four coarse-mesh elements the randomly generated functions $A_\xi$ are given according to \cref{mp2:Pictures_TSRBLOD}, and we again set $f \equiv 1$.
The exact definition of the functions $A_\xi$ is again left out for brevity and can be found in the accompanying code.
The approximate maximum contrast is $\kappa \approx 16$.

To fully resolve the microstructure on all $4096$ ($n_H = 2^{6}$) coarse elements, we need to choose $n_h = 2^{13}$ which results in about $67.1$ million degrees of freedom.
The local corrector problems have roughly $1.3$ million degrees of freedom.
With 1024 available parallel processes, each must train ROMs for $4$ coarse elements.

We choose a training set $\ParamsTrain \subset \Params$ of $4^3$ equidistant parameters and show results for $\varepsilon_1 = 0.01$ and $\varepsilon_2=0.02$ in \cref{Tab:mp2:1}.
Again, the TSRBLOD shows high online efficiency with speed-ups of up to $10^6$, clearly outperforming the RBLOD.
We note that the reported online time of $t^{\textnormal{online}}=4.39\,\mathrm{s}$ for the RBLOD roughly splits into $2.2\,\mathrm{s}$ required for sequentially solving the reduced corrector problems (Step~1 in \cref{def:comp_proc_LOD}), $1.2\,\mathrm{s}$ for assembling the coarse system (Step~2) and
$1\,\mathrm{s}$ for solving the coarse system (Step~3).
In particular, the further speed-up that could be achieved for the RBLOD by parallelizing the corrector problems is bounded by a factor of approximately 2.

Also, the storage requirements are noteworthy: the reduced data required for evaluating the TSRBLOD ROM and $\eta_{a,\mu}$ is only 28~KB in size, whereas the RBLOD requires 409~MB.

\begin{table}
	\footnotesize
	\centering
	{\def\arraystretch{1.7}\tabcolsep=3pt
		\begin{tabular}{|l|c|c|} \hline
			Method  &\hphantom{T}RBLOD\hphantom{T} & TSRBLOD 
			\\ \hline  \hline
			$t^\textnormal{offline}_{1}(T)$   &   10278	& 	11289	\\
			$t^\textnormal{offline}_{1}$   	&   49436 	& 	54837	\\
			$t^\textnormal{offline}_2$   		& 	- 		&	9206 	\\
			$t^\textnormal{offline}$    		& 	49436 	&	64043	\\ \hline
			Cum. size Stage 1 				&	278528	&	193289 	\\
			Av. size Stage 1           	&	17.00  	&	47.19	\\
			Size Stage 2 					& 	- 		& 	16 		\\
			\hline \hline
			$t^{\textnormal{LOD}}$ &  	\multicolumn{2}{c|}{515} \\ \hline
			$t^\textnormal{online}$ 		 	            & 	4.39 	& 	0.0005	\\ \hline \hline
			Speedup w.r.t LOD 	 		            & 	117  	&	9.57e5	\\ \hline \hline
			$e^{H^1, \textnormal{rel}}_{\textnormal{LOD}}$ 	 	&	1.95e-5 &	4.43e-4 \\
			$e^{L^2, \textnormal{rel}}_{\textnormal{LOD}}$ 	 	&	2.36e-5 &	4.49e-4 \\ 
			\hline
	\end{tabular}}
	\captionsetup{width=\textwidth}
	\caption[Experiment 11: Performance, ROM sizes and accuracy of the methods]{\footnotesize{%
			Performance, ROM sizes and accuracy of the methods for test case 2 with $\varepsilon_1 = 0.01$ and $\varepsilon_2=0.02$. All times are given in seconds.
			\label{Tab:mp2:1}}}
\end{table}

\section{Summary and outlook}\label{sec:TSRBLOD_conclusion}
In this chapter, we started with a general view on localized model reduction and explained the PG--LOD method for solving (parameterized) multiscale problems.
In particular, we assessed the applicability of the PG--LOD to many-query and real-time applications.
In \cref{sec:adaptive_PGLOD_from_hekema}, the adaptive PG--LOD for perturbed problems showed the re-usability of previous (related) corrector problems based on a reference structure.

For tackling parameterized problems as they are usually defined in the context of MOR methods, we derived a new two-scale RB scheme for efficient real-time evaluations of the PG--LOD.
Due to the two-stage reduction process and the independence of the resulting ROMs from the size of both the fine-scale and the coarse-scale meshes, our approach is practical even for large-scale problems and, given \cref{asmpt:parameter_separable}, has extremely low storage requirements.
For ease of presentation, we have assumed non-parametric right-hand sides and did not consider output functionals.
Incorporating both into our approach is straightforward (and also used in \cref{chap:TR_TSRBLOD}).
Furthermore, for large coarse meshes, additional intermediate reduction stages can be added to further reduce the needed computational effort in Stage~2.
Instead of a fixed a priori choice of the Stage~1 tolerance $\varepsilon_1$, the Stage~1 ROMs can be adaptively enriched during Stage~2 when an insufficient approximation quality of some of the Stage~1 ROMs is detected.
The presented methodology could also be applied to other problem classes and multiscale methods.

To conclude, the derived TSRBLOD and also the RBLOD constitute instances of localized model order reduction techniques for tackling parameterized multiscale problems in the sense of \cref{sec:localized_MOR}.
As the last step of this thesis, the subsequent chapter is devoted to using these LMOR methods to develop localized versions of the TR-RB algorithm from \cref{chap:TR_RB} for tackling the RB challenge of an inaccessible global discretization for PDE-constrained parameter optimization problems; cf. \cref{sec:RB_challenges}.

\chapter{Localized trust-region reduced basis methods}\label{chap:TR_TSRBLOD}
\noindent 
The purpose of this chapter is to combine localized RB techniques with the adaptive trust-region method for PDE-constrained parameter optimization problems that involve large-scale and multiscale PDEs where FEM evaluations are considered prohibitive.
The resulting trust-region localized RB method (TR-LRB) adaptively constructs local RB models and deviates from a classical global FEM approximation for the FOM state.

While combining error-aware TR methods with localized RB techniques is new, localized MOR methods are an extensively studied field which we detailed in \cref{chap:TSRBLOD}.
The primary focus of an efficient adaptive reduced method is to derive localized error bounds for detecting where the model requires local basis updates and to find efficient global coupling techniques.
As an explicit example, we use the localized RB technique for the LOD multiscale method from \cref{chap:TSRBLOD} that already internally inherits all requirements for extending the TR-RB algorithm from \cref{chap:TR_RB} to a TR-LRB approach.

Again, we emphasize that, as far as this thesis is concerned, using (localized) MOR for solving a single PDE-constrained parameter optimization problem aims at the overall-efficiency of the solution method.
Thus, neither offline computation nor online computations of the TR-LRB approach can be considered negligible.
The results of this chapter have partly been published in~\cite{paper4}.

The chapter is organized as follows: In \cref{sec:TR-LRB}, we discuss the general formulation of a TR-LRB algorithm.
Subsequently, we show that the RBLOD and TSRBLOD can be used as an instance of a TR-LRB algorithm.
Lastly, we present numerical experiments based on a multiscale benchmark problem that demonstrate the benefit of localized techniques as well as the relaxed TR-RB algorithm, introduced in \cref{sec:Relaxed_TRRB}.

\section{Trust-region localized reduced basis algorithm}
\label{sec:TR-LRB}

This section is devoted to a general view of how localized model reduction techniques can extend the TR-RB method.
Concerning the basic description of the TR-RB algorithm, as explained in \cref{chap:TR_RB} (and summarized in \cref{alg:Basic_TR-RBmethod}), we can speak about a TR-LRB method if a localized model replaces the FOM and the reduced model, in turn, is a localized reduced RB model; cf. \cref{sec:localized_MOR}.

To this end, let $\Jhat_h^{\textnormal{loc}}$ be the output functional of a model that uses localized computations w.r.t. $\grid$ without having to assemble system matrices on the high-fidelity space $V_h$ (for instance, a multiscale, or domain decomposition method in the sense of \cref{sec:localized_MOR}).
Thus, in $\Jhat_h^{\textnormal{loc}}$, the $h$ index indicates that local computations w.r.t. the fine mesh $\grid$ are used.
Furthermore, let $\Jhat_\red^{\textnormal{loc}}$ be the output functional of the corresponding localized RB approach, which, for now, is not specified further. For theoretical purposes, the replacement of $\Jhat_h$ by $\Jhat_h^{\textnormal{loc}}$ is only valid by assuming that the discretization error of the localized FOM can be neglected.

\noindent
Analog to \cref{asmpt:truth}, we make the following assumption:

\begin{assumption}[The localized FOM is the ``truth'']
	\label{asmpt:loc_truth}
	We assume that the localized discretization error $| \Jhat_h^{\textnormal{loc}}(\mu) - \Jhat(\mu)|$ can be neglected, where $\Jhat_h^{\textnormal{loc}}$ defines the output functional based on an approach that does not require such a global approximation. This also translates to the primal and dual approximation error being negligible.
\end{assumption}

Let $\Jhat^{\textnormal{loc}}_\red$ denote an RB reduced functional of $\Jhat_h^{\textnormal{loc}}$.
Just as in \cref{chap:TR_RB}, for characterizing the trust-region, we assume that $\Jhat^{\textnormal{loc}}_\red$ admits an a posteriori error result, such that
\begin{equation} \label{eq:local_loc_red_estimator}
| \Jhat^{\textnormal{loc}}_h(\mu) - \Jhat^{\textnormal{loc}}_\red(\mu) | \leq \Delta_{\Jhat^{\textnormal{loc}}_r}(\mu),
\end{equation}
where $\Delta_{\Jhat^{\textnormal{loc}}_r}(\mu)$ can be computed without explicitly evaluating $\Jhat^{\textnormal{loc}}_h$ at $\mu \in \Params$.
We emphasize that, in the original formulation in \cref{chap:TR_RB}, the reduced functional is based on a global RB method, and the respective error result has been presented in \cref{sec:TRRB_a_post_error_estimates}.
If the reduced scheme instead stems from a localized model reduction approach, we abbreviate the resulting method as the TR-LRB method.

Importantly, if \cref{asmpt:loc_truth} is not fulfilled, the localized FOM can not be trusted (or only to a certain extent), and thus, the FOM-based stopping criterion can not be reliably used if convergence w.r.t. FEM is desired.
Moreover, in such a case, $\Delta_{\Jhat^{\textnormal{loc}}_r}$ loses its validity if the TR-LRB is employed for finding the true critical point of the underlying FEM (or infinite-dimensional) problem.
If, instead, the distance to the true (or FEM) solution $| \Jhat_\red^{\textnormal{loc}}(\mu) - \Jhat_h(\mu)|$ can be bounded, the algorithm can still reliably be followed but requires FEM evaluations for the stopping criterion, suitable error estimation on the reduced gradient w.r.t. FEM (similar to what is used in \cite{QGVW2017}), or an alternative FOM-based termination criterion.
For the remainder of this chapter, we do not consider such cases and, for simplicity, assume that \cref{asmpt:loc_truth} is always fulfilled and a suitable error result is given.

In our numerical experiments, we also consider a case where \cref{asmpt:loc_truth} is slightly relaxed, emphasizing the broad applicability of the TR-LRB method.
It is important to mention that the fulfillment of \cref{asmpt:loc_truth} heavily depends on the localized approach and on the specific problem.
Also, it is not always priorly known that a localized approach fulfills \cref{asmpt:loc_truth}.

\begin{remark}[Convergence of the TR-LRB method] \label{rmk:unconditional_convergence}
	The convergence of TR-LRB does not require \cref{asmpt:loc_truth} to be fulfilled.
	As the convergence study from \cref{sec:TRRB_convergence_analysis} suggests, the TR-LRB method at least converges to a critical point of $\Jhat_h^{\textnormal{loc}}$ even if this point is not a critical point of $\Jhat_h$.
	In this case, we can not rigorously say much about the critical point of the TR-LRB method unless FEM evaluations or reliable error estimators w.r.t. the true solution are available.
\end{remark}

To conclude, all features of the TR-RB that we detailed in \cref{sec:TR_RB_algorithm} can naturally be used for the localized case, whereas exceptional circumstances may occur depending on the localized reduction approach at hand.
Again, as usual for RB methods, \cref{asmpt:parameter_separable} plays a vital role in obtaining a fast surrogate model and respective estimation.
We particularly emphasize that the TR-LRB method can also be used within the relaxed variant (R-TR-LRB), which uses the concepts from \cref{sec:Relaxed_TRRB}. 

In the following, we refer to the TR-RB summarized in \cref{alg:Basic_TR-RBmethod}.
In a dedicated TR-LRB algorithm, Line~1 must contain the initialization of the localized RB model, which, for solving the TR sub-problem in Line~4, needs to have an online-efficient estimator for the output functional readily available.
Moreover, for the localized RB model, we require a notion of adaptive enrichment in Line~7.
Concerning the FOM-based convergence criterion, the quantity~$\Jhat_h$ that is to be computed in Line~3 can be replaced by $\Jhat_h^{\textnormal{loc}}$ and the same applies to $\varrho^{(k)}$ as it also contains $\Jhat_h$ terms.

While the required changes that come with TR-LRB methods sound reasonably straightforward, the absence of a dedicated FEM-based FOM model raises new challenges that need to be specified depending on the localized approach.
Moreover, optional local RB basis enrichment opens the possibility of new adaptive strategies, where computational effort and basis sizes are particularly targeted.
Due to the wide variety of localized reduction approaches, we can not discuss all features and challenges in TR-LRB methods.
Instead, we discuss two specific examples of a TR-LRB algorithm, using the reduced RBLOD and TSRBLOD multiscale approach presented in \cref{chap:TSRBLOD}.

\section{Trust-region (TS)RBLOD method}
\label{sec:TR-RBLOD_methods}

We develop a TR-LRB method based on the RB approaches for the LOD from \cref{chap:TSRBLOD}.
To this end, we give a precise definition of the localized FOM functional $\Jhat_h^{\textnormal{loc}}$, the localized ROM functional $\Jhat_\red^{\textnormal{loc}}$, and its gradient information $\nabla\Jhat_h^{\textnormal{loc}}$, and $\nabla\Jhat_\red^{\textnormal{loc}}$, respectively.
Furthermore, we elaborate on the error estimator~$\Delta_{\textnormal{loc}}$ and give details on how the localized RB enrichment is performed.

We note that, for applying the PG--LOD, we use the same assumptions as in \cref{chap:TSRBLOD}.
In particular, $A_\mu \in L^\infty(\Omega, \R^{d \times d})$ is a symmetric, uniform elliptic, and parameter separable multiscale coefficient.
Furthermore, we again assume that the objective functional $\J$ is linear-quadratic, such that it fulfills \cref{asmpt:quadratic_J}.

\subsection{Localized full-order model}
\label{sec:localized_FOM}
As detailed in \cref{sec:LOD}, we choose a coarse-mesh space $V_H$ and compute a corresponding localized LOD space by 
\[
V_{H,\kl,\mu}^{\textnormal{ms}}= (I - \QQkmpr)(V_H),
\]
where $\QQkmpr$ contains localized corrector functions $\QQktmpr$ computed on the localized fine-scale spaces $\Vfhkt: = \Vfh \cap H^1_0(U_\kl(T))$.
Therefore, the localized primal state $\uhkmsm \in V_{H,\kl,\mu}^{\textnormal{ms}}$ can be computed via the Petrov--Galerkin version of the LOD in Equation \eqref{eq:PG_}, i.e.
\begin{equation}
\label{eq:PG_LOD}
	a_{\mu}(\uhkmsm,v_H) =  l(v_H) \qquad \textnormal{for all } \; v_H \in V_H.
\end{equation}
In order to relate to the two-scale based view on the PG--LOD, we further note that there exists a uniquely defined two-scale representation $\utsm \in \Vts$ of $\uhkmsm \in V_{H,\kl,\mu}^{\textnormal{ms}}$, where
$$ \Vts := V_H \oplus \Vfhkt[T_1] \oplus \cdots \oplus \Vfhkt[\Tlast],$$
and $\utsm \in \Vts$ is the solution of 
\begin{equation}
\Btsm\left(\utsm,  \vts \right) = \Fts(\vts) \qquad \textnormal{for all } \vts \in \mathcal{V}. \label{eq:two_scale_TRLRB}
\end{equation}
As proven in \cref{prop:two_scale_solution}, the two-scale solution $\utsm \in \Vts$, can always be constructed from~\eqref{eq:PG_LOD}.

We we did not consider any objective functional in \cref{chap:TSRBLOD}, the corresponding primal state can easily be used to compute the corresponding localized FOM objective functional, i.e.
\begin{equation}
	\Jhat_h^{\textnormal{loc}}(\mu) := \J(\uhkm, \mu).
\end{equation}
Again, the subindex $h$ in $\Jhat_h^{\text{loc}}$ refers to the fact that the construction of the solution space $V_{H,\kl,\mu}^{\text{ms}}$ for solving \eqref{eq:PG_LOD} internally requires the computation of the correctors that resolve the fine-scale mesh, which can then be discarded immediately.
Note that we do not plugin $\uhkmsm \in V_{H,\kl,\mu}^{\text{ms}}$ into $\J$ since the basis of $V_{H,\kl,\mu}^{\text{ms}}$ may not be available.
This is the case if the corrector problems can not be stored, which, in many cases, can be sufficient since the coarse-scale behavior is anyway captured by $\uhkm \in V_H$; cf. \cite{elf}.
To align with this, we employ the following structural assumption on $\J$, which agrees with what is usually expected in a multiscale setting.
\begin{assumption}[$\J$ is a coarse functional] \label{asmpt:coarse_J}
	We assume that the objective functional is a coarse functional, i.e.~for all $u_H \in V_H$ and $u^\text{f} \in \Vfh$, we have $$\J(u_H + u^\text{f}, \mu) = \J(u_H, \mu).$$ 
\end{assumption}
We conclude that we aim to find a local optimum of the following reduced PDE-constrained parameter optimization problem
\begin{align}
\min_{\mu \in \Paramsad} \Jhat_h^{\textnormal{loc}}(\mu).
\tag{$\hat{\textnormal{P}}^{\textnormal{loc}}_h$}\label{Phat_h_loc}
\end{align}
To solve \eqref{Phat_h_loc} with the TR-LRB algorithm, we require the gradient information of $\Jhat_h^{\text{loc}}$.
In our approach, a dual model is utilized, which we also solve with the PG--LOD.
While \eqref{eq:PG_LOD} works as a replacement for \eqref{eq:state_h}, we formulate a corresponding PG--LOD version of the dual problem for \eqref{eq:dual_solution_h}:
Seek a function $\phkmsm \in V_{H,\kl,\mu}^{\text{ms}}$, such that
\begin{equation}
\label{eq:PG_LOD_dual}
a_{\mu}(v_H, \phkmsm) =  \partial_u \J(\uhkm, \mu)[v_H] \qquad \text{for all } \; v_H \in V_H.
\end{equation}
Note that \cref{asmpt:coarse_J} justifies that $\uhkm$ is used for the right-hand side of \eqref{eq:PG_LOD_dual}.
Looking at \eqref{eq:PG_LOD_dual}, we conclude that, just as the FOM in \cite{paper1}, the localized FOM is a conforming choice in the sense that $\uhkmsm$ and $\phkmsm$ belong to the same space~$V_{H,\kl,\mu}^{\text{ms}}$.
Note that this choice only makes sense if the given multiscale coefficient $A_\mu$ is symmetric, as we have assumed throughout this article.
In that case, the recaptured multiscale effects for the primal and dual operator are the same.
If instead, $A_\mu$ is not symmetric, different LOD spaces must be constructed, which we do not consider.

With 
\begin{align} \label{eq:dual_two_scale_rhs}
\Fts^\du_{\mu,\uhkm}\left((v_H, \vfti{1}, \dots, \vfti{{\abs{\mathcal{T}_H}}})\right) := \partial_u \J(\uhkm, \mu)[v_H],
\end{align}
we can formulate the two-scale dual solution $\ptsm \in \Vts$ of
\begin{equation}
\Btsm\left(\ptsm,  \vts \right) = \Fts^\du_{\mu,\uhkm}(\vts) \qquad \text{for all } \vts \in \mathcal{V}, \label{eq:two_scale_dual_TRLRB}
\end{equation}
where we note that we did not flip the arguments, due to the symmetry of $a_\mu$.

Concerning the definition of $\Jhat_h^{\textnormal{loc}}$, one may be concerned about the fact that the corresponding residual term associated to~\eqref{eq:PG_} with test function $\phkmsm \in V_{H,\kl,\mu}^{\textnormal{ms}}$, i.e., $r_{\mu}^\pr(\uhkmsm)[\phkmsm]$ has to be added to use the Lagrangian functional, analogously to the term that has been used as NCD correction in \cref{sec:TRRB_ncd_approach}.
With respect to the discussion from above that a basis of $V_{H,\kl,\mu}^{\textnormal{ms}}$ can not be stored, the term is not computable.
If we only insert the coarse dual solution $\phkm \in V_H$ to the residual, the resulting term is zero since it vanishes on $V_H$.

Finally, we compute the gradient information with the following FOM-based formula:
\begin{equation} \label{eq:LOD_gradient}
\nabla_\mu \Jhat_h^{\text{loc}}(\mu) = \partial_\mu \J(\uhkm, \mu) +  \partial_\mu r_{\mu}^\pr(\uhkm)[\phkm].
\end{equation}
Again, we do not use $\uhkmsm \in V_{H,\kl,\mu}^{\text{ms}}$ and $\phkmsm \in V_{H,\kl,\mu}^{\text{ms}}$ but instead their coarse-scale representations $\uhkm \in V_H$ and $\phkm \in V_H$ to be able to discard corrector information directly after their computation.
Instead, we need to sacrifice the potential gain of accuracy.
Certainly, we could still plugin $\uhkmsm  \in V_{H,\kl,\mu}^{\text{ms}}$, e.g., for the linear terms of~$\J$ since the related terms can be prepared simultaneously to the assembly of $\Ktmu$.

For keeping the theory short, we do not consider Hessian information in the localized approach but note that using Newton's method is relatively straightforward.

\begin{remark}[Fulfillment of \cref{asmpt:loc_truth}]
	As discussed earlier, multiscale approaches like the LOD are known to struggle, for instance, with high-contrast problems or rapid coarse-scale changes induced by high conductivity channels.
	With \cref{asmpt:coarse_J}, we rule out the concern that the dual problem may not fit to the right-hand side structure of \eqref{eq:PG_}.
	To fulfill \cref{asmpt:loc_truth}, we may thus assume that a good choice of the coarse mesh size $H$, localization parameter $\kl$, and the fine mesh size $h$ is made to cope with the underlying problem.
	For the LOD, this aligns with the aim to fulfill the a priori result of \cref{thm:pglod_convergence} sufficiently well, i.e.~that the errors $\norm{\uhm - \uhkm}_{L^2}$, $\norm{\phm - \phkm}_{L^2}$, and the corresponding $H^1$-errors are sufficiently small.
\end{remark}

\subsection{Localized reduced-order models}

It remains to explain how the reduced objective functional and its gradient are computed with the localized reduced LOD model.
Suitable localized RB models have been investigated in \cref{chap:TSRBLOD}.
The Petrov--Galerkin version of the RBLOD constructs RB models for the corrector problems and gathers the local reduced solutions in a global Petrov--Galerkin LOD scheme.
Furthermore, the new TSRBLOD also uses these RB correctors and uses the two-scale formulation to reduce the global scheme.

The TSRBLOD showed to be more beneficial in terms of online efficiency, whereas the additional coarse-scale reduction introduces an additional approximation error.
According to \cref{sec:RBLOD_estimator_study}, the additional error of the TSRBLOD can rigorously be controlled. 
As seen in \cref{sec:exp_RBLOD_2}, for large coarse systems, the online-acceleration w.r.t. the RBLOD can be of multiple orders.
On the other hand, the offline construction of the TSRBLOD is higher than the RBLOD since an additional offline-online decomposition is to be performed; cf.~\cref{sec:stage_2_offon}.

In what follows, we consider both localized RB techniques as potential models for the TR-LRB.
We introduce the following abbreviations as a subclass of TR-LRB methods and refer to \cref{sec:TSRBLOD_MOR} for the respective details of the approaches.
We again note that all additional features of the TR-RB from \cref{sec:TR_RB_algorithm} may be incorporated into the TR-RBLOD and TR-TSRBLOD.

\subsubsection{Reduced model based on the RBLOD}\label{sec:TRLRB_with_RBLOD}
We use a Petrov--Galerkin version of the RBLOD, which, in its Galerkin version, has been initially proposed in \cite{RBLOD}.
Different from our version of the RBLOD in \cref{sec:TSRBLOD_experiments}, which has mainly been constructed for ease of comparison to~\cite{RBLOD}, in what follows, we employ the same Stage~1 reduction process as used for the TSRBLOD; see also the discussion in \cref{rem:reduction_flexibility}.
To recall, we construct only a single reduced corrector for each $T$ equivalently to \cref{sec:stage1} and form a reduced LOD space $V_{H,\kl,\mu}^{\textnormal{rb},\textnormal{ms}}$ and a corresponding reduced two-scale space
$$\Vtsrblod := V_H \oplus \Vfrbkt[T_1] \oplus \dots \oplus \Vfrbkt[T_{\Tlast}] \subset \Vts.$$
Note that we do not need to store a basis of these spaces.
In the online-phase, we compute the local contributions $\Ktmurb$ to the global approximation of the PG--LOD stiffness matrix $\Kmurb$, which is then used for a replacement of~\eqref{eq:PG_in_matrices}, such that we instead solve \eqref{eq:stage_2_fom}, i.e.
\begin{equation}\label{eq:stage_2_fom_2}
\Kmurb \cdot \uhkmrblodcoeff = \mathbb{F},
\end{equation}
where $\uhkmrblodcoeff$ denotes the coefficient vector of $\uhkmrblod \in V_H$ and $\uhkmsmrblod \in V_{H,\kl,\mu}^{\textnormal{rb},\textnormal{ms}}$.
We again refer to \cref{fig:RBLOD} for a visualization of the computational procedure and mention that, concerning \cref{def:comp_proc_LOD}, Step~1 is locally reduced.
At the same time, Step~2 and Step~3 remain the same in terms of computational effort.
To conclude, the reduced localized objective functional for the RBLOD is defined as
\begin{equation}
	\Jhat_\red^{\textnormal{rblod}}(\mu) := \J(\uhkmrblod, \mu).
\end{equation}
Due to \cref{asmpt:coarse_J}, we do not insert $\uhkmsmrblod \in V_{H,\kl,\mu}^{\textnormal{rb},\textnormal{ms}}$, which constitutes to the fact that we can not store fine-scale correctors.
Similarly, as for the primal equation, reduced corrector spaces can be built for the dual equation if $a_\mu$ is not symmetric, which we do not detail further due to the symmetry assumption on $a_\mu$.
To approximate~\eqref{eq:PG_LOD_dual}, the dual version of~\eqref{eq:stage_2_fom_2} can be stated by
\begin{equation}\label{eq:stage_2_fom_dual}
\Kmurb \cdot \phkmrblodcoeff = \mathbb{F}^\du,
\end{equation} 
where $(\mathbb{F}^\du)_i = \partial_u \J(\uhkmrblod, \mu)[\phi_i]$.
Finally, using the corresponding solution $\phkmrblod \in V_H$, gradient information can be computed by the standard reduced approach, i.e.
\begin{equation} \label{eq:RBLOD_gradient}
	\nabla_\mu \Jhat_\red^{\textnormal{rblod}}(\mu) = \partial_\mu \J(\uhkmrblod, \mu) +  \partial_\mu r_{\mu}^\pr(\uhkmrblod)[\phkmrblod].
\end{equation}

\subsubsection{Reduced model based on the TSRBLOD} 
\label{sec:TRLRB_with_TSRBLOD}
To also reduce the primal and dual PG--LOD, we use the TSRBLOD approach that has been described in \cref{sec:TSRBLOD_MOR}.
Since the primal and dual equations have different right-hand sides, we consider two two-scale reduced spaces $\Vtsrbpr, \Vtsrbdu \subset \Vtsrblod$.
Reducing the primal equation \eqref{eq:two_scale_TRLRB} with the TSRBLOD, given $\Vtsrbpr$, is completely analog to \cref{sec:stage_2_red} and means to further reduce~\eqref{eq:stage_2_fom_2}.
We compute the two-scale reduced primal solution $\utsmrb \in \Vtsrbpr$ by
\begin{equation}\label{eq:tsrblod_primal}
\utsmrb := \argmin_{\uts \in \Vtsrbpr}\, \sup_{\vts \in \Vts} \frac{\Fts(\vts) - \Btsm(\uts, \vts)}{\Snorm{\vts}},
\end{equation}
which is equal to~\eqref{eq:stage2_rom}.
Further, let $\uhkmrb \in V_H$ denote the resulting TSRBLOD coarse-scale approximation, which can be reconstructed from $\utsmrb$, just by using the $V_H$ part of the respective basis of $\Vtsrbpr$.
Then, we can define the corresponding reduced functional by
\begin{equation}\label{eq:TSRBLOD_func}
	\Jhat_\red^{\textnormal{rb}}(\mu) := \J(\uhkmrb, \mu).
\end{equation}
Given the dual two-scale reduced space $\Vtsrbdu$, the dual problem can be defined analogously to the primal TSRBLOD reduction, with the vital difference that the right-hand side of the Stage~2 FOM system needs to be replaced by the dual right-hand side.
To be precise, we define 
\begin{align} \label{eq:dual_two_scale_rhs_rb}
\Fts^\du_{\mu,\uhkmrb}\left((v_H, \vfti{1}, \dots, \vfti{{\abs{\mathcal{T}_H}}})\right) := \partial_u \J(\uhkmrb, \mu)[v_H].
\end{align}
By replacing $\Fts$ by $\Fts^\du_{\uhkmrb}$ from~\eqref{eq:dual_two_scale_rhs} in~\eqref{eq:tsrblod_primal} and using $\Vtsrbdu$ instead, we obtain the two-scale dual solution $\ptsmrb \in \Vtsrbdu$ by
\begin{equation}\label{eq:tsrblod_dual}
\ptsmrb := \argmin_{\pts \in \Vtsrbdu} \,\sup_{\vts \in \Vts} \frac{\Fts^\du_{\uhkmrb}(\vts) - \Btsm(\pts, \vts)}{\Snorm{\vts}},
\end{equation}
which again uses the fact that $a_\mu$ is symmetric.
With the resulting coarse approximate $\phkmrb \in V_H$ reconstructed from $\ptsmrb$, we can compute the reduced gradient by
\begin{equation} \label{eq:TSRBLOD_gradient}
	\nabla_\mu \Jhat_\red^{\textnormal{rb}}(\mu) = \partial_\mu \J(\uhkmrb, \mu) +  \partial_\mu r_{\mu}^\pr(\uhkmrb)[\phkmrb].
\end{equation}

\begin{remark}[Generalization of the TSRBLOD approach to parameterized right-hand sides]\label{rem:TSRBLOD_has_mu_rhs}
	For the reduction process of \cref{eq:dual_two_scale_rhs}, we emphasize that this requires generalizing the TSRBLOD approach to parameterized right-hand sides. Thus, the offline-online decomposition as explained in \cref{sec:stage_2_offon} changes slightly.
	We still omit a further technical description for brevity, noting that the online efficiency remains the same.
\end{remark}

\subsection{A posteriori error estimate for the reduced functional}

Having set up the localized FOM and ROM approximation schemes, we aim to derive the error estimator $\Delta_{\Jhat_r^{\text{loc}}}$ of the reduced functional from \eqref{eq:local_loc_red_estimator}, which is needed for characterizing the TR in~\eqref{TRRBsubprob}.
Luckily, for both approaches, \cref{thm:a_posteriori} can be used with simple modifications for both the primal and the dual problems of the RBLOD and TSRBLOD.
On top of that, similar to  \cref{prop:Jhat_error}, we combine the results to obtain an estimator for the reduced functional.
We show that the resulting estimators admit an estimator such that both reduced models and can be used to follow a TR-LRB algorithm.
To recall from \cref{thm:a_posteriori}, we define the following two-scale estimators:
\begin{align}\label{eq:stage2_est_a_pr}
	\eta^\pr_{a,\mu}(\uts) &:= 
	\sqrt{5}(\CPG)^{-1}
	\sup_{v\in\Vts}
	\frac{\Fts(\vts) - \Btsm(\uts,\vts)}
	{\Snorm{\vts}},\\\label{eq:stage2_est_a_du}
	\eta^\du_{a,\mu}(\pts) &:= 
	\sqrt{5}(\CPG)^{-1}
	\sup_{v\in\Vts}
	\frac{\Fts^\du_{\mu,\uhkm}(\vts) - \Btsm(\pts,\vts)}
	{\Snorm{\vts}}.
\end{align}

\begin{proposition}[Upper bound on the local primal model reduction error] \label{prop:primal_rom_error_LOD}
	For $\mu \in \Params$ let $\utsm \in \Vts$ be the solution of~\eqref{eq:two_scale_TRLRB}.
	\begin{enumerate}
		\item [\emph{(i)}] 
	Let $\uhkm^{\textnormal{rblod}} \in V_H$ be the primal RBLOD solution of \eqref{eq:stage_2_fom_2} and let $\utsmrblod \in \Vtsrblod$ be the corresponding two-scale solution that additionally includes the reduced solutions of all corrector problems needed to assemble \eqref{eq:stage_2_fom_2}.
	Then, it holds
	\begin{align}
	\Anorm{\utsm - \utsmrblod} \leq \Delta_{\pr}^{\textnormal{rblod}}(\mu) := \eta_{a,\mu}(\utsmrblod),
	\end{align}
	with $\eta_{a,\mu}$ defined in \cref{thm:a_posteriori}.
	\item [\emph{(ii)}]
	Let $\utsmrb \in \Vtsrb$ be the two-scale reduced solution of~\eqref{eq:tsrblod_primal}
	.
	Then, it holds
	\begin{align}
	\Anorm{\utsm - \utsmrb} \leq \Delta_{\pr}^{\textnormal{rb}}(\mu) := \eta_{a,\mu}(\utsmrb).
	\end{align}
	\end{enumerate}
\end{proposition}
\begin{proof}
	The assertions follow directly from \cref{thm:a_posteriori}.
\end{proof}

We note that the estimator $\Delta_{\pr}^{\textnormal{rblod}}$ is based on the two-scale formulation and its a posteriori error result. Still, the RBLOD does not use a two-scale reduction.
A posteriori error analysis w.r.t. the true LOD solution has not been proposed in \cite{RBLOD} and \cref{prop:primal_rom_error_LOD}(i), instead, enables such a result.

\begin{remark}[Equivalence of the two-scale norms] \label{rmk:norm_equivalence_LOD}
	Due to the definitions of $\Anorm{\cdot}$, $\Snorm{\cdot}$, $\anorm{\cdot}$, and $\snorm{\cdot}$ and the  equivalences of $\Anorm{\cdot}$ and $\Snorm{\cdot}$, as well as $\anorm{\cdot}$, and $\norm{\cdot}$, respectively, we note that the fine-scale errors $\anorm{\uhkmsm - \uhkmsmrb}$, $\norm{\uhkmsm - \uhkmsmrb}$, and the coarse-scale errors $\anorm{\uhkm - \uhkmrb}$, and $\norm{\uhkm - \uhkmrb}$ can be bounded by $\Delta_{\pr}^{\text{rb}}$ with the respective equivalence constants, cf.~\cite{tsrblod}.
\end{remark}

\noindent A similar result is also available for the equivalent $\mu$-dependent $\Snorm{\cdot}$-norm.

The corresponding dual estimates account for the fact that the right-hand side contains the reduced primal solution instead of the true LOD solution.

\begin{proposition}[Upper bound on the local dual model reduction error]
\label{prop:dual_rom_error_LOD}
	For $\mu \in \Params$ let $\ptsm \in \Vts$ be the solution of the two-scale dual equation \eqref{eq:two_scale_dual_TRLRB}.
\begin{enumerate}
	\item [\emph{(i)}]
	Let $\phkm^{\textnormal{rblod}} \in V_H$ be the dual RBLOD solution of \eqref{eq:stage_2_fom_dual} and let $\ptsmrblod \in \Vtsrblod$ be the corresponding two-scale solution that additionally includes the reduced solutions of all corrector problems needed to assemble \eqref{eq:stage_2_fom_dual}.
	Then, it holds
	\begin{align}
	\Anorm{\ptsm - \ptsmrblod} \leq \Delta_{\du}^{\textnormal{rblod}}(\mu) \mkern-3mu := \mkern-3mu \frac{\sqrt{5}}{\CPG} \left(2 \cont{k_\mu}\Delta_{\pr}^{\textnormal{rblod}}(\mu) + \eta^\du_{a,\mu}(\ptsmrblod)\right).
	\end{align}
	\item [\emph{(ii)}]
	Let $\ptsmrb \in \Vtsrb$ be the two-scale reduced dual solution of~\eqref{eq:tsrblod_dual}.
	Then, it holds
	\begin{align}
	\Anorm{\ptsm - \ptsmrb} \leq \Delta_{\du}^{\textnormal{rb}}(\mu) := \frac{\sqrt{5}}{\CPG} \left(2 \cont{k_\mu}\Delta_{\pr}^{\textnormal{rb}}(\mu)   + \eta^\du_{a,\mu}(\ptsmrb)\right).
	\end{align}
\end{enumerate}
\end{proposition}
\begin{proof}
	Again, the proof of (i) and (ii) are the same. For proving (ii), we use the shorthands $\mathfrak{e}_{\mu}^{\du} := \ptsm - \ptsmrb \in \Vtsrb$ and $e_{H, \mu}^{\pr} := \uhkm - \uhkmrb$, where $\uhkm \in V_H$ and $\uhkmrb$ are the $V_H$ parts of $\utsm$ and $\utsmrb$, respectively.
	With the inf-sup stability from \cref{thm:two_scale_inf_sup}, we have
	\begin{align*}
	\CPG/\sqrt{5} \, \Anorm{\mathfrak{e}_{\mu}^{\du}} &\leq \sup_{0 \neq \vts \in \Vts}
	\frac{\Btsm(\mathfrak{e}_{\mu}^{\du},\vts)}{\Snorm{\vts}}  = 
	 \sup_{0 \neq \vts \in \Vts} \left(
	\frac{\Fts^\du_{\uhkm}(\vts)}{\Snorm{\vts}} - \frac{\Btsm(\ptsmrb,\vts)}{\Snorm{\vts}}\right) \\
	&= 	 \sup_{0 \neq \vts \in \Vts} \left(
	\frac{\Fts^\du_{\uhkm}(\vts)}{\Snorm{\vts}} - \frac{\Fts^\du_{\uhkmrb}(\vts)}{\Snorm{\vts}} + \frac{\Fts^\du_{\uhkmrb}(\vts)}{\Snorm{\vts}} - \frac{\Btsm(\ptsmrb,\vts)}{\Snorm{\vts}}\right) \\	
	&\leq 2 \|k_\mu\|\; \|e_{H, \mu}^\pr\|\; + \eta^\du_{a,\mu}(\ptsmrb),
	\end{align*}
	with $\Fts^\du_{*}$ defined in~\eqref{eq:dual_two_scale_rhs} and~\eqref{eq:dual_two_scale_rhs_rb}, which is linear in its sub-index argument due to~\eqref{eq:quadratic_J}.
	We attain the desired result utilizing Proposition \ref{prop:primal_rom_error_LOD} and \cref{rmk:norm_equivalence_LOD}.
\end{proof}

Similar to \cref{rmk:norm_equivalence_LOD}, we note that the respective dual estimators can also bound the corresponding dual norms from \cref{prop:dual_rom_error_LOD}.
Finally, we derive the a posteriori error result for the reduced objective functional.

\begin{proposition}[Upper bound for the reduced functionals]
	\label{prop:LOD_reduced_functional}
	For $\mu \in \Params$ let $\utsm \in \Vts$ be the two-scale solution of \eqref{eq:two_scale_TRLRB} with coarse part $\uhkm \in V_H$ and LOD-space representation $\uhkmsm \in V_{H,\kl,\mu}^{\textnormal{ms}}$.
	Further, let $\ptsm \in \Vts$ be the two-scale solution of \eqref{eq:two_scale_dual_TRLRB} with coarse part $\phkm \in V_H$ and LOD-space representation $\phkmsm \in V_{H,\kl,\mu}^{\textnormal{ms}}$.
	\begin{enumerate}
		\item[\emph{(i)}] We have for the RBLOD reduced cost functional
		\begin{align*}
		\hspace{-1cm}|\Jhat_h^{\textnormal{loc}}(\mu) - \Jhat_\red^{\textnormal{rblod}}(\mu)| \lesssim \Delta_{\Jhat_\red^{\textnormal{rblod}}}(\mu) :=  \Delta_{\pr}^{\textnormal{rblod}}(\mu) &\eta^\du_{a,\mu}(\ptsmrblod) 
		+ (\Delta_{\pr}^{\textnormal{rblod}}(\mu))^2 \cont{\kformd} + \Delta_{\textnormal{trunc}}^{\textnormal{rblod}}(\mu),
 		\end{align*}
		where $\ptsmrblod \in \Vtsrb$ denotes the two-scale representation of the RBLOD solution and 
		$\Delta_{\textnormal{trunc}}^{\textnormal{rblod}}(\mu)$ is a truncation-reduction-based homogenization term which is specified below. 
		\item [\emph{(ii)}] Furthermore, we have for the TSRBLOD reduced cost functional
		\begin{align*}
		|\Jhat_h^{\textnormal{loc}}(\mu) - \Jhat_\red^{\textnormal{rb}}(\mu)| \lesssim \Delta_{\Jhat_\red^{\textnormal{rb}}}(\mu) :=  \Delta_{\pr}^{\textnormal{rb}}(\mu) &\eta^\du_{a,\mu}(\ptsmrb) 
		+ (\Delta_{\pr}^{\textnormal{rb}}(\mu))^2 \cont{\kformd}  + \Delta_{\textnormal{trunc}}^{\textnormal{rb}}(\mu),
		\end{align*}
		where $\ptsmrb \in \Vtsrbdu$ denotes the two-scale reduced dual equation. 
		\item[(iii)] For the TSRBLOD, the truncation-reduction-based homogenization term $\Delta_{\textnormal{trunc}}^{\textnormal{rb}}(\mu)$ is defined as 
		\begin{equation} \label{eq:hom_red_term}
		 \Delta_{\textnormal{trunc}}^{\textnormal{rb}}(\mu) := a_\mu(\ehms,\QQkmrb(\phkmrb)) 
		\end{equation}
		and can be estimated by
		\begin{equation} \label{eq:estimation_hom_red_term}
\hspace{-0.5cm}\abs{\Delta_{\text{trunc}}^{\textnormal{rb}}(\mu)} \leq \Delta_{\pr}^{\textnormal{rb}}(\mu) \left(2 c \, \kl^{d/2} \theta^k \snorm{\phkmrb} + \eta^\du_{a, \mu}(\ptsmrb)\right) + \alpha^{-1/2}\snorm{\phkmrb} \eta^\pr_{a, \mu}(\utsmrb),
\end{equation}
		with respective coarse- and two-scale-space primal and dual solutions and constant $c>0$. The RBLOD truncation term $\Delta_{\textnormal{trunc}}^{\textnormal{rblod}}(\mu)$ is defined analogously.
	\end{enumerate}
\end{proposition}
\begin{proof}
	The proof for (i) and (ii) is the same.
	To proof (ii), we utilize \cref{asmpt:coarse_J} to incorporate the estimates of \cref{prop:primal_rom_error_LOD} and \cref{prop:dual_rom_error_LOD}.
	The rest of the proof is similar to \cref{prop:Jhat_error}(i).
	By using the shorthands $\ehms := \uhkmsm - \uhkmsmrb$ and $\eh := \uhkm - \uhkmrb$, we have 
	\begin{align*}
		|\Jhat_h^{\textnormal{loc}}&(\mu) - \Jhat_\red^{\textnormal{rb}}(\mu)| 
		= |\J(\uhkm, \mu) - \J(\uhkmrb, \mu)| 
		\\
		&= |j_\mu(\ehms) + \kformd(\uhkm, \uhkm) - \kformd(\uhkmrb, \uhkmrb) 
		- a_\mu(\ehms,\phkmsmrb) + a_\mu(\ehms,\phkmsmrb)| \\
		&= |r_{\mu}^\du(\uhkmrb, \phkmsmrb)[\ehms] + \kformd(\ehms, \ehms) + a_\mu(\ehms,\phkmsmrb)| \\
		&\leq \eta^\du_{a,\mu}(\ptsmrb)\; \|\ehms\| +\cont{\kformd} \; \|\ehms\|^2 
		+ |a_\mu(\ehms,\QQkmrb(\phkmrb))|,
	\end{align*}
	where we have used that
	\begin{align*}
	a_\mu(\ehms,\phkmsmrb) &= - a_\mu(\ehms,\QQkmrb(\phkmrb)).
	\end{align*}
	For (iii), we further note that	
	\begin{align*}
		a_\mu(\ehms,\QQkmrb(\phkmrb))  &= a_\mu(\ehms,\QQm(\phkmrb))  \mkern-3mu
		- a_\mu(\ehms,(\QQm \mkern-3mu - \mkern-3mu \QQkm \mkern-3mu + \mkern-3mu\QQkm \mkern-3mu - \mkern-3mu \QQkmrb)(\phkmrb))  
	\end{align*}
	and
	\begin{align*}
		a_\mu(\ehms,\QQm(\phkmrb)) 
		&= a_\mu((\QQkm \mkern-3mu - \mkern-3mu \QQm)(\uhkm - \uhkmrb) + (\QQkm \mkern-3mu - \mkern-3mu \QQkmrb)(\uhkmrb),\QQm(\phkmrb)).
	\end{align*}
	So,
	\begin{align*}
		|a_\mu(\ehms,\QQkmrb(\phkmrb))| &\leq \anorm{\ehms} (\anorm{(\QQm \mkern-3mu - \mkern-3mu \QQkm)(\phkmrb)}
		+ \anorm{(\QQkm \mkern-3mu - \mkern-3mu \QQkmrb)(\phkmrb)}) \\
		& \quad + \anorm{\QQm(\phkmrb)} (\anorm{(\QQm \mkern-3mu - \mkern-3mu \QQkm)(\eh)} +
		\anorm{(\QQkm \mkern-3mu - \mkern-3mu \QQkmrb)(\uhkmrb)}) \\
		& \lesssim \Delta_{\pr}^{\textnormal{rb}}(\mu) ( k^{d/2} \theta^k \anorm{\phkmrb}
		+ \eta_{a, \mu}(\ptsmrb)) \\
		& \quad + \anorm{\phkmrb} ( k^{d/2} \theta^k \anorm{\eh} + \eta_{a, \mu}(\utsmrb)) \\
		& \leq \Delta_{\pr}^{\textnormal{rb}}(\mu) \left(2 k^{d/2} \theta^k \anorm{\phkmrb}
		+ \eta^\du_{a, \mu}(\ptsmrb)\right) + \anorm{\phkmrb} \eta^\pr_{a, \mu}(\utsmrb),
	\end{align*}
	where we have used the a priori result on the corrector decay (cf. \cref{thm:pglod_convergence}) and \cref{rmk:norm_equivalence_LOD}. Using the equivalence of $\anorm{\cdot}$ and $\snorm{\cdot}$ attains the assertion.
\end{proof}
\noindent In conclusion, for the RBLOD as well as for the TSRBLOD, \cref{prop:LOD_reduced_functional} ensures that~\eqref{eq:local_loc_red_estimator} is fulfilled.

\begin{remark}[Truncation-reduction-based homogenization term] \label{rem:trunc-red-hom-term}
	In \cref{prop:LOD_reduced_functional}, we intentionally separated the error estimation from the homogenization term $\Delta_{\text{trunc}}^{rb}(\mu)$, and presented a rather naive estimation of it.
	The reason is that the term can be interpreted as a truncation term that (without reduction) vanishes for  true LOD-space functions, i.e.
	\begin{equation}
	a_\mu(u_{H,\mu}^{\text{ms}},\QQm(\phkmrb)) = 0,
	\end{equation}
	for all $u_{H,\mu}^{\text{ms}} \in V_{H,\mu}^{\text{ms}}$, since $\QQm(\phkmrb) \in \Vfh$ and $V_h = V_{H,\mu}^{\text{ms}} \oplus_{a_{\mu}} \Vfh$.
	Due to \cref{asmpt:truth}, the a priori term can be neglected from \eqref{eq:estimation_hom_red_term}, such that $\Delta_{\text{trunc}}^{rb}$ can be computed efficiently.
\end{remark}

The computation of the above-explained TSRBLOD estimator can be offline-online decomposed with a numerically stable procedure, which has carefully been explained in~\cite{tsrblod}.
Although the Stage~1 models are also offline-online efficient, it is important to mention the RBLOD still contains loops over all $T$ to compute the two-scale estimates. 

The additional orthonormalization of the residual terms of Stage~1 is necessary for the RBLOD as well as for the Stage~2 residual.
In contrast, the additional expenses for stabilizing the Stage 2 residual are not strictly needed in our approach.
Indeed, concerning the overall cost of the TR-LRB algorithm, we omit the offline-online decomposition of Stage 2 entirely and instead compute the residual and its Riesz-representative whenever needed, cf.~\cref{sec:TR_TS_pseudo-code}.

\subsection{Local basis enrichment}
\label{sec:local_basis_enrichment}
It remains to explain the adaptive localized enrichment strategy for a parameter $\mu \in \Params$, for instance, a newly accepted iterate $\mu^{(k)} \in \Params$ of the TR-LRB algorithm.
In \cref{chap:TR_RB}, it is advised to either update the RB space unconditionally or optionally w.r.t. some enrichment flag.
Since we deal with local RB models, the situation is more complex.
For obtaining certified convergence in the sense of \cref{Thm:convergence_of_TR}, we always perform an enrichment.
However, some local models may reject or dismiss the snapshots if, e.g., the selected parameter does not influence the local model.
For this reason, a localization of the error result is crucial since we may employ the local error estimators for the local models to decide whether an update is required. Such a strategy is also commonly known as adaptive online enrichment.

First of all, we note that the estimates $\eta_{a, \mu}^\pr$ and $\eta_{a, \mu}^\du$ can indeed be boiled down to their respective local reduction errors, namely the standard RB estimation of Stage~1 of the reduction process for both the RBLOD and TSRBLOD.
To recall from \cref{sec:stage_1_err}, for each $T \in \Gridh$, we use the residual-norm based estimate
\begin{equation}\label{eq:stage1_error_estimate_rephrased}
\anorm{\QQktm(v_H) - \QQktmrb(v_H)} 
\leq \eta_{T,\mu}(\QQktmrb(v_H)),
\end{equation}
where
\begin{equation} \label{eq:stage_1_estimator_rephrased}
\eta_{T,\mu}(\QQktmrb(v_H)) := \alpha^{-1/2} \sup_{\vft\in\Vfhkt}
\frac{a_\mu^T(v_H, \vft) - a_\mu(\QQktmrb(v_H), \vft)}{\snorm{\vft}}.
\end{equation}

For the TR-LRB scheme, at an enrichment point, for every $T$, we use the Stage~1 estimator $\Delta^T_{\text{loc}}(\mu) := \eta_{T,\mu}(\QQktmrb(v_H))$ to decide for a local enrichment.
If the estimator is below a specific tolerance $\tau_{{\text{loc}}} > 0$, i.e.
\begin{equation}\label{eq:optional_enrichment}
\Delta^T_{\text{loc}}(\mu) \leq \tau_{{\text{loc}}},
\end{equation}
the enrichment is skipped.
For a sufficiently small $\tau_{{\text{loc}}}$, the enrichment strategy can be considered unconditionally.
However, the computational effort can still decrease significantly, e.g., if a local modal is not associated with the parameter.
We also note that the online adaptive approach is also motivated by the numerical experiments in~\cref{sec:TSRBLOD_experiments} and \cite{RBLOD},
where it was demonstrated that moderate choices of $\tau_{{\text{loc}}}$ already produce acceptable reduced models.
Admittedly, the choice of the tolerance $\tau_{{\text{loc}}}$ is highly problem-dependent.
If the tolerance is chosen too large, the method could be stagnant (due to the missing local basis quality).
In such cases, it is recommend to refine the tolerance adaptively.
For simplicity we omit such a strategy.

While the local enrichment is helpful for the Stage~1 reduction of the RBLOD and TSRBLOD, another crucial question is how the Stage~2 reduction is performed for the TSRBLOD.
Note that the TS\-RB\-LOD model is based on the reduced models from Stage~1, meaning that, whenever the Stage~1 models are enriched, the Stage~2 reduction needs to be restarted from scratch.
Thus, there is more freedom in choosing the enrichment parameters for the TSRBLOD model.
With respect to the fact that, at iteration $k$, the Stage~1 models are exact (up to the tolerance $\tau_{{\text{loc}}}$), we propose to enrich the TSRBLOD model for the same sequence of TR iterates $\mu^{(i)}$, for $i = 0, \dots, k$.
Greedy-based enrichments of the TSRBLOD are also possible, mainly because the snapshot generation with Stage~1 is fast.
Certainly, our experiments suggested that greedy-search algorithms do not provide a significant update to the accuracy of the TSRBLOD model.

Let us also mention that we did not use optional local enrichment strategies for the R-TR-LRB algorithm since we do not have access to the estimators in the first iterations.

\subsection{Algorithm in Pseudo-code} \label{sec:TR_TS_pseudo-code}

For clarity, we summarize the (R)-TR-TSRBLOD algorithms for parameter optimization of multiscale problems in the following.
We emphasize that the respective version of the (R)-TR-RBLOD method are left out since they are equivalent to the (R)-TR-TSRBLOD algorithm with the only difference that Stage~2 is left out.

\begin{algorithm2e}[!h]
	\KwData{$\delta^{(0)}$, $\beta_1\in (0,1)$, $\eta_\varrho\in [\frac{3}{4},1)$, $\mu^{(0)}$, $\tau_\textnormal{{sub}}$, $\tau_\textnormal{{FOC}}$, and $\beta_2\in(0,1)$ analog to \cref{alg:Basic_TR-RBmethod}, tolerance for online enrichment $\tau_{{\text{loc}}}$}
	Initialize TSRBLOD model with $\mu^{(0)}$\;
	Set $k=0$\; 
	\While
	{$\|\mu^{(k)}-\Proj_{\Paramsad}(\mu^{(k)}-\nabla_\mu \Jhat^{\textnormal{loc}}_h(\mu^{(k)}))\|_2 > \tau_{\textnormal{{FOC}}}$\label{stopping_condition_L}}{
		Compute $\mu^{(k+1)}$ from~\eqref{TRRBsubprob} with termination criteria~\eqref{Termination_crit_sub-problem}\label{solve_sub_problem_L}\;
		\uIf{Sufficient decrease condition~\eqref{eq:suff_decrease_condition} is fulfilled (cf.~\cite{paper1})\label{suff_dec}} {
			Accept $\mu^{(k+1)}$ and possibly enlarge the TR-radius (cf.~\cite{paper1})\label{accept}\;
			Before enrichment: check \eqref{eq:TR_termination} for early termination\label{stop_early}\;
			\tcp{TSRBLOD enrichment}
			\textbf{Stage 1:} enrich the local RB corrector models at $\mu^{(k+1)}$ if \eqref{eq:optional_enrichment} including pre-assembly of the local estimators $\Delta^T_{\text{loc}}$\label{enrichment_L}\;
			\textbf{Stage 2:} construct the primal and dual two-scale models and enrich for all $\mu^{(k')}$, $k'=0,\dots,k+1$, and do not pre-assemble $\Delta_{\Jhat_\red^{\text{loc}}}$\;
		}
		\Else {
			Reject $\mu^{(k+1)}$, shrink the TR radius $\delta^{(k+1)} = \beta_1\delta^{(k)}$ and go to \ref{solve_sub_problem}\;		
		}
		Set $k=k+1$\;
	}
	\caption{TR-TSRBLOD algorithm}
	\label{alg:TR-TSRBLOD}
\end{algorithm2e}

We emphasize that the TR-TSRBLOD procedure in \cref{alg:TR-TSRBLOD} is analog to the algorithm presented in \cite{paper1} but with the major difference that the localized LOD-based FOM and the localized TSRBLOD reduced model, including its respective estimator is used.
The details of \cref{alg:R-TR-TSRBLOD} are explained in \cref{sec:Relaxed_TRRB}.
As stated in Line \ref{stop_early} of \cref{alg:TR-TSRBLOD}, we check the FOM termination criterion prior to the enrichment.
This is because online enrichment, including the assembly of the respective estimators, is relatively more expensive than the pure computation of the termination criterion.
A similar strategy is used in Line \ref{stop_early_rel} of \cref{alg:R-TR-TSRBLOD}.
However, since Stage 1 is generally less expensive for the relaxed variant, we check the FOM-based termination criterion between Stage~1 and Stage~2 of the TSRBLOD.

We further note that in Lines \ref{suff_dec} and \ref{accept} of both algorithms, we have neglected detailed information on the exact computational procedure concerning the cheap conditions for the sufficient decrease conditions according to \cite{YM2013} and the enlarging of the TR-radius with a suitable accessible condition.

\begin{algorithm2e}[!h]
	\KwData{$\delta^{(0)}$, $\beta_1\in (0,1)$, $\eta_\varrho\in [\frac{3}{4},1)$, $\mu^{(0)}$, $\tau_\textnormal{{sub}}$, $\tau_\textnormal{{FOC}}$, and $\beta_2\in(0,1)$ analog to \cref{alg:Basic_TR-RBmethod}, relaxation sequences $(\varepsilon_\text{TR}^{(k)})_k$ and $(\varepsilon_\text{cond}^{(k)})_k$.}
	Initialize TSRBLOD model with $\mu^{(0)}$ without constructing error estimates\;
	Set $k=0$\; 
	\While
	{$\|\mu^{(k)}-\Proj_{\Paramsad}(\mu^{(k)}-\nabla_\mu \Jhat^{\textnormal{loc}}_h(\mu^{(k)}))\|_2 > \tau_{\textnormal{{FOC}}}$\label{stopping_condition_R}}{
		Compute $\mu^{(k+1)}$ from~\eqref{eq:rel_TRRBsubprob} with relaxed termination~\eqref{FOC_sub-problem} and \eqref{eq:rel_cut_of_TR}\label{solve_sub_problem_R}\;
		\uIf{Relaxed sufficient decrease condition~\eqref{eq:suff_decrease_condition_rel} is fulfilled (cf.~\cite{paper1})} {
			Accept $\mu^{(k+1)}$ and possibly enlarge the TR-radius (cf.~\cite{paper1})\;
			\tcp{TSRBLOD enrichment}
			\textbf{Stage 1:} enrich the local RB corrector models at $\mu^{(k+1)}$ and skip offline assembly of the local estimators if $\varepsilon_\text{TR}^{(k)}$ is large enough\;
			Before Stage 2 enrichment: check \eqref{eq:TR_termination} for early termination\label{stop_early_rel}\;
			\textbf{Stage 2:} construct the primal and dual two-scale models and enrich for all $\mu^{(k')}$, $k'=0,\dots,k+1$, and do not assemble estimator\;
		}
		\Else {
			Reject $\mu^{(k+1)}$, shrink the TR radius $\delta^{(k+1)} = \beta_1\delta^{(k)}$ and go to \ref{solve_sub_problem}\;		
		}
		Set $k=k+1$\;
	}
	\caption{Relaxed TR-TSRBLOD algorithm}
	\label{alg:R-TR-TSRBLOD}
\end{algorithm2e}

\subsection{A note on computational efficiency}

As usual for localized model order reduction approaches, the computational procedures for the discussed TR-LRB are complex and depend on many criteria, especially considering the exact implementational design and the available computer power.
The theory of the previous sections suggests that the RBLOD and TSRBLOD both have computational advantages and disadvantages.

An advantage of the RBLOD is that it does not use an additional reduction process after the Stage~1 models are enriched, and thus, the enrichment time compared to the TSRBLOD can be much lower.
However, the computation of the reduced solutions, as well as the computation of the error estimators, is always proportional to $\abs{\mathcal{T}_H}$, which may harm the computational time that is required for the TR-sub-problem.
In contrast, the TSRBLOD naturally enables a fast computation of the reduced solutions and estimates, which is especially helpful for the fast computation of the TR sub-problem solution.
Admittedly, as already said, the online speed comes with a more expensive enrichment phase that is needed to perform the Stage~2 reduction process.

In conclusion, it is not trivial to declare the most efficient reduction process among the RBLOD and TSRBLOD.
Despite the above-mentioned issues, we again emphasize that \cref{asmpt:parameter_separable} ensures the most achievable offline-online efficiency of the respective reduced approaches as explained in \cref{sec:TSRBLOD_MOR}.

\section{Software and implementational design}
\label{sec:software_TR_TSRBLOD}

The main concepts of the implementational design are already discussed concerning the TR-RB algorithm, in \cref{sec:software_TRRB}, and the (TS)RBLOD, in \cref{sec:software_TSRBLOD}.
For the numerical experiments in this chapter, we combine both implementations.
While the basic \texttt{TR{\textunderscore}algorithm} and  \texttt{Relaxed{\textunderscore}TR{\textunderscore}algorithm} that has been used in \cref{chap:TR_RB} can directly be reused for the TR-LRB, for the localized reduction schemes, we need a \texttt{gridlod}-based \texttt{pdeopt}-\texttt{discretizer}, \texttt{Model} and \texttt{Reductor} that are similar to those discussed in \cref{sec:pdeopt_pymor}.
We note that for the code in \cref{sec:software_TSRBLOD}, we did not define a \texttt{Model} for the original PG--LOD, and, instead, computed these solutions in a \texttt{gridlod} way.
In particular, we did not use a \texttt{StationaryProblem} for the benchmark problems and directly defined a corresponding \texttt{Model}.
Since, in this chapter, we aim at comparing global methods, such as the ones in \cref{chap:TR_RB}, to localized methods, such as the ones in \cref{chap:TSRBLOD}, we required a mutual understanding of the underlying analytical problem.
For this reason, we implemented a way to discretize a \texttt{StationaryProblem} given that the resulting \texttt{Model} is based on \texttt{gridlod}.
Therefore, we enable a closer linkup of \texttt{gridlod} to \texttt{pyMOR}'s builtin discretizer used for the global methods.
Selected implemented objects and methods are mentioned in the following.

\begin{description}
	\item[\texttt{gridlod}-\texttt{discretizer}:] The function \texttt{discretize{\textunderscore}gridlod} interprets the given data functions for treating them in \texttt{gridlod}, prepares the localized patches and sets additional LOD-related parameters.
	For comparison purposes, the coarse LOD system $\mathcal{T}_H$ is discretized with \texttt{pyMOR}s standard discretizer, where we note that the DoF mappings coincide.
	\item[\texttt{GridlodModel}:] The discretizer returns a \texttt{GridlodModel} that inherits from \texttt{StationaryModel} and internally uses the computational procedure that is summarized in \cref{def:comp_proc_LOD}.
\end{description}

\noindent With the \texttt{GridlodModel} at hand, we can now proceed an define the \texttt{gridlod}-specific counterparts for the \texttt{pdeopt} objects:

\begin{description}
	\item[\texttt{discretizer}:] The \texttt{discretize{\textunderscore}quadratic{\textunderscore}pdeopt{\textunderscore}with{\textunderscore}gridlod}-function internally calls the \texttt{discretize{\textunderscore}gridlod}-function and prepares the \texttt{gridlod}-specific view on the quadratic objective functional $\Jhat$.
	\item[\texttt{Model}:] The \texttt{QuadraticPdeoptStationaryModel} from the \texttt{pdeopt} can be reused and only small adjustments are made. Based on the given FOM, the \texttt{Model} automatically detects, whether FEM or PG--LOD shall be used as FOM.
	\item[\texttt{Reductor}:] Since the enrichment process and the estimators substantially differ from the global RB case, a new \texttt{QuadraticPdeoptStationaryCoerciveLODReductor} is implemented.
	In there, the entire reduction routine that we used in \cref{chap:TSRBLOD}, as well as further specifics, for instance, w.r.t. the dual problem, are hidden.
\end{description}

\noindent For more details, we again refer to \cref{apx:code_availability}.
Let us also mention that the LOD-based method is written in such a way that it can profit from parallelization, which is again implemented with \texttt{pyMOR}'s \texttt{MPI}-parallelization tools.

\section{Numerical experiments}
\label{sec:loc_experiments}

We analyze the presented TR-LRB approaches with three experiments with the same problem description, only differing in their respective multiscale complexity. We define the fine-mesh by $n_h \times n_h$ and the coarse-mesh by $n_H \times n_H$ quadrilateral grid-blocks of $\Omega:=[0,1]^2$, used to determine the standard FE mesh $\mathcal{T}_h$ and $\mathcal{T}_H$, respectively, with traditional $\mathcal{P}^1$-FE spaces $V_h$ and $V_H$.
The mesh-sizes $h$ and $H$ can be computed from $n_h$ and $n_H$.
In the first small experiment, we compare the localized methods to FEM-based (TR-RB) methods, whereas, in the second large experiment, we neglect FEM entirely, as it is computationally infeasible.
The third experiment serves the only purpose to further showcase the computational difference between the RBLOD and TSRBLOD.

In the experiments, we focus on the number of evaluations relative to their complexity to the fine FEM mesh-size~$h$, the coarse LOD mesh-size~$H$, or the respective low dimensions of the reduced models.
Moreover, we provide run time comparisons that hint at the computational efficiency observed with our specific implementation.

We emphasize that our computations were performed on an HPC cluster with $400$ parallel processes.
Nevertheless, the observed run times can not be interpreted as minimal computational times of the localized algorithms.
More HPC-based and non-\texttt{Python}-based implementations can strengthen the localized approaches even more.
We further emphasize that the Stage~2 reduction has been implemented as a serialized process that neither gathers the observed data from Stage~1 efficiently nor uses the sparsity pattern of the two-scale system matrix to a minimal extent.

\begin{figure}
	\begin{subfigure}[b]{0.49\textwidth}
		\includegraphics[width=\textwidth]{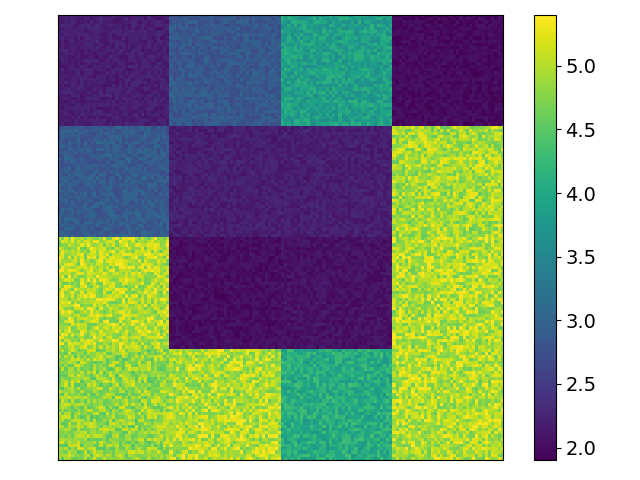}
	\end{subfigure}
	\begin{subfigure}[b]{0.49\textwidth}
		\includegraphics[width=\textwidth]{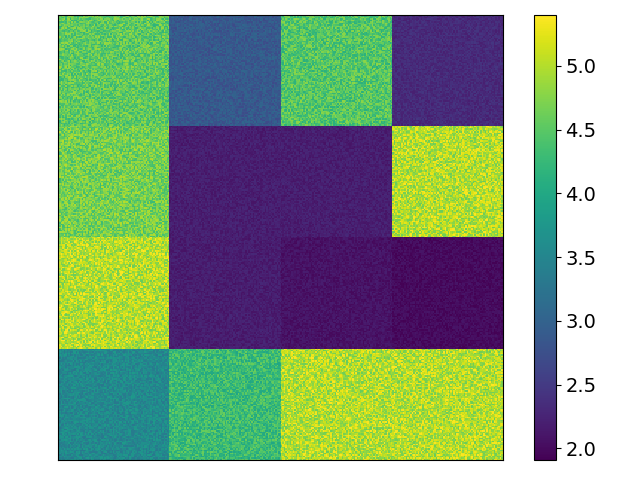}
	\end{subfigure}
	\centering
	\caption[Coefficient $A^1_\mu$ and $A^2_\mu$.]{{Coefficient $A^1_\mu$ with $N_1 = 150$ (left) and $A^2_\mu$ with $N_2 = 300$ (right) for the desired state of $\mu^{\text{d}} \in \Params$.}}
	\label{fig:thermal_block}
\end{figure}

As usual, we use the $L^2$-misfit objective functional with a Tikhonov-regularization term, i.e.
\begin{align*}
\J(v, \mu) = \frac{\sigma_d}{2} \int_{D}^{} (I_H(v) - u^{\text{d}})^2 \dx + \frac{1}{2} \sum^{P}_{i=1} \sigma_i (\mu_i-\mu^{\text{d}}_i)^2 + 1.
\end{align*} 
Here, $\mu^{\text{d}} \in \Params$ is the desired parameter and $u^{\text{d}} = I_H(u_{\mu^{\text{d}}})$ the corresponding desired solution specified in each experiment.
The use of the interpolation operator $I_H$ ensures \cref{asmpt:coarse_J}.
Moreover, using an actual solution as desired temperature and the respective desired parameter in the objective functional ensures that the optimization problem is sufficiently regular, such that all optimization methods converge to the same point.
As detailed earlier, $\J$ can easily be written in the linear-quadratic form as in \cref{asmpt:quadratic_J}.
For the diffusion or conductivity coefficient $A_\mu$ in the symmetric bilinear form $a_\mu$, we consider a 4x4-thermal block problem with two different thermal block multiscale coefficients $A_\mu^1$ and $A_\mu^2$, i.e.
$$
A_\mu = \sum_{\xi = 1}^{16} \mu_\xi A_\mu^{1,\xi} +  \sum_{\xi =17}^{32} \mu_\xi A_\mu^{2,\xi}.
$$
Each of the 4x4 blocks has its linear parameter value, i.e., $\Params \subseteq \mathbb{R}^{32}$.
The respective parameterized multiscale blocks are given by $A_\mu^{1,\xi} = A_\mu^{1}\big|_{\Omega_{i,j}}$ and $A_\mu^{2,\xi} = A_\mu^{2}\big|_{\Omega_{i,j}}$, where $\Omega_{i,j}$ denotes the $(i,j)$-th thermal block for $i,j=1,2,3,4$. 
The multiscale features are randomly constructed with normally distributed values in $\mathcal{N}([0.9,1.1])$ on a $N_1 \times N_1$ (for $A_\mu^{1}$) and $N_2 \times N_2$ (for $A_\mu^{2}$) quadrilateral grid.
The specific values for $N_1$ and $N_2$ are given for each experiment.
Due to the affine decomposition of $A_\mu$, the multiscale features are linearly scaled by $\mu_\xi$ within each thermal block.
We would like to point out that the multiscale data does not admit periodicity or any other additional structure; see \cref{fig:thermal_block} for a visualization.
Moreover, both coefficients $A_\mu^{1}$ and $A_\mu^{2}$ have low-conductivity blocks in the middle of the domain, i.e., for $\Omega_{i,j}$, $i,j=2,3$.
The low conductivity is enforced by a restriction on the parameter space, i.e. we consider the admissible parameter set $\Params = [1,4]^{24} \times [1,1.2]^8$.
We choose the non-parameterized constant function $f_\mu \equiv 10$ as right-hand side function.
As also used in \cref{chap:TR_RB}, for the inner product of $V_h$, we use the energy norm $\norm{\cdot} := \norm{\cdot}_{a,\check{\mu}}$ for a fixed parameter $\check{\mu} \in \Params$ in the middle of the parameter space. Thus, constants in the estimators can easily be deduced by the min/max-theta approach.
For the two-scale estimators, we further approximate the maximum contrast $\kappa$ and the respective constants $\alpha$ and $\beta$ accordingly.

Details on the fine- and coarse-mesh are given in the respective experiments.
The desired parameter $\mu^{\text{d}} \in \Params$ is equal for both experiments, and mimics the case where boundary constraints are active, i.e. we set $\mu^{\text{d}}_i = 4$ for $i=3,4,6,7,8,9,11,14$ and $\mu^{\text{d}}_i = 1.2$ for $i=28,29,30,31$.
The remaining values of $\mu^{\text{d}}$ are chosen randomly, see \cref{fig:thermal_block}.
The initial guess $\mu^{(0)}$ is also chosen randomly, where we note that the choice of the initial guess is not relevant for the shown results, which is why we only show the experiments for a single initial guess.
Furthermore, the weights for the objective functional are chosen as $\sigma_d=100$ and $\sigma_i = 0.001$ for each $i = 1, \dots, P$.

Similar to the experiments in \cref{chap:TR_RB}, we choose an initial TR radius of $\delta^{(0)} = 0.1$,
a TR shrinking factor $\beta_1=0.5$, an Armijo step-length $\kappa=0.5$,
a truncation of the TR boundary of $\beta_2 = 0.95$,
a tolerance for enlarging the TR radius of $\eta_\varrho = 0.75$,
a stopping tolerance for the TR sub-problems of $\tau_\text{{sub}} = 10^{-8}$,
a maximum number of TR iteration $K = 40$,
a maximum number of sub-problem iterations $K_{\text{{sub}}}= 400$,
a maximum number of Armijo iteration of $50$, and
a stopping tolerance for the FOC condition $\tau_\text{{FOC}}= 10^{-6}$.

\paragraph{State-of-the-art methods}
The following algorithms are used to compare our algorithms to the literature.
\setlist[description]{font=\normalfont\bfseries\space,labelindent=10pt}
\begin{description}
	\item[1. FEM BFGS:] Analogously to what was used in~\cref{sec:exp_paper1}, we perform a standard projected BFGS method that is solely based on classical FEM evaluations. This means that the high-fidelity space $V_h$ is used as the only approximation space, and no reduced approach is used. \vspace{0.4cm}
	\item[2. TR-RB BFGS \cref{alg:TR-RBmethod_paper1}:] Given the FEM discretization used for Method 1, we use a trust-region reduced basis algorithm with full certification and global RB evaluations based on FEM enrichments.
	As the reduced model, we choose the NCD-corrected approach with Lagrangian enrichment.
	Moreover, we use the projected BFGS as the ROM-based TR sub-problem.
	After each sub-problem, the global RB model is enriched after acceptance of the iterate, and the algorithm is terminated with a FEM-based FOC-type criterion. We also point to Variant 3(a) in \cref{sec:exp_paper1}.
\end{description}

\paragraph{Selected new methods}

This chapter uses the relaxation of the TR-(L)RB method and further proposes a localized approach tailored toward multiscale problems.
Consequently, we consider the following relaxed variant of Method 2:
\begin{description}
	\item[2.r R-TR-RB BFGS:] In this method, we use the relaxed trust-region reduced basis variant for the FEM-based TR-RB algorithm, as explained in \cref{sec:TRRB_related_approaches}.
	We use global NCD-corrected RB evaluations with relaxed certification.
	For the relaxation, we choose the relaxation sequences $\varepsilon_\text{TR}^{(k)} = \varepsilon_\text{cond}^{(k)}:= 10^{10-k}$.
	Therefore, the first iterations can be considered certification-free, and, as discussed in  \cref{sec:Relaxed_TRRB}, we do not pre-assemble the TR-estimator for $k \leq 8$.
	Moreover, from iteration count $k > 27$, we relax the TR-RB method below double machine-precision and follow the TR-RB method in its original form.
\end{description}
As detailed in \cref{sec:TR-RBLOD_methods}, given a respectively accurate LOD discretization, we formulate the following localized methods.
\begin{description}
	\item[3. PG--LOD BFGS:] For comparison, we utilize the standard projected BFGS method with PG--LOD evaluations without using reduced models.
	The PG--LOD system is always constructed from scratch and does not use any prior knowledge from previous parameters.
	In particular, if FEM is not accessible and \cref{asmpt:truth} is given, this method is considered the FOM method.
	\item[4. TR-TSRBLOD BFGS (TR-TS):] We use the localized TSRBLOD reduction process as detailed in \cref{sec:TRLRB_with_TSRBLOD} for the PG--LOD and use the TR-LRB method with PG--LOD enrichments and local RB evaluations based on the TS\-RB\-LOD.
	The procedure is summarized in \cref{alg:TR-TSRBLOD}.
	The sub-problems are again solved with the BFGS method, and the outer algorithm is fully certified.
	Moroever, we use a local enrichment tolerance $\tau_{{\text{loc}}}=10^{-3}$ which has proven to be sufficient for our experiment, cf.~\cref{sec:local_basis_enrichment}.
	\item[4.r R-TR-TSRBLOD BFGS (R-TR-TS):] Just as explained in Method~2.r, we devise the relaxed version of Method~4., by choosing $\varepsilon_\text{TR}^{(k)} = \varepsilon_\text{cond}^{(k)}:= 10^{10-k}$.
	Again, this means that the first iterations are certification-free, and no error estimators need to be prepared. The procedure is summarized in \cref{alg:R-TR-TSRBLOD}.
	\item[5. TR-RBLOD BFGS (TR-RBLOD):] We use the localized RBLOD reduction process as detailed in \cref{sec:TRLRB_with_RBLOD} for the PG--LOD and use the TR-LRB method with PG--LOD enrichments and local RB evaluations based on the RB\-LOD.
	\item[5.r R-TR-RBLOD BFGS (R-TR-RBLOD):] As explained above, we devise the relaxed version of Method~5., analog to Method~4.r. 	
\end{description}

\paragraph{Complexity measures}
To assess the presented methods w.r.t. their computational demands, we count the evaluations of the respective systems.
To be precise, we deviate between the following complexities:
\begin{description}
	\item[FEM:] FEM evaluations, proportional to $\mathcal{T}_h$, which are needed for approximating \eqref{eq:state_h} or \eqref{eq:dual_solution_h} with FEM or for enriching the respective global RB model for Methods 2.r and 2.
	\item[RB:] Global RB evaluations for approximating the FEM system, proportional to the global basis size, used in Methods 2.r and 2.
	\item[LOD coarse:] Coarse PG--LOD system evaluations with exact corrector data, meaning to solve \eqref{eq:PG_LOD} or \eqref{eq:PG_LOD_dual}, proportional to the coarse mesh $\mathcal{T}_H$.
	These are only required in Method 3 and for the FOM-based termination criterion in Methods 4.r and 4.
	\item[LOD local:] Local evaluations of all FOM corrector problems that are required for assembling the multiscale stiffness matrix of \eqref{eq:PG_LOD} and \eqref{eq:PG_LOD_dual}, locally proportional to $U_\kl(T)_h$.
	\item[RBLOD coarse:] Coarse PG--LOD system evaluations with RB-based corrector data for the multiscale stiffness matrix, required for the online phase of the RBLOD, as well as for the snapshots generation in Stage~2 of the TSRBLOD, proportional to the coarse mesh $\mathcal{T}_H$.
	\item[RBLOD local:] Evaluations of RB corrector problems, proportional to the local RB sizes.
	\item[TSRBLOD:] Evaluations of the TSRBLOD system, proportional to the two-scale RB size.
\end{description}

\paragraph{Error measures}
As the optimization target, we validate the respective accuracy of the methods by considering the relative error in the optimal value of $\Jhat$, i.e., we consider
$$
e^{\Jhat,\text{rel}}(\bar{\mu}) := |\Jhat(\mu^{\text{d}})  - \Jhat(\bar{\mu})| / \Jhat(\mu^{\text{d}}),
$$
where $\bar{\mu}$ is the respective convergence point of the optimization methods and $\Jhat$ is either the FEM-based objective functional~$\Jhat_h$ or the LOD-based objective functional~$\Jhat^\textnormal{loc}_h$.

\subsection{Experiment 12: Comparison with FEM-based methods}

\begin{figure}
	\centering
	\begin{subfigure}{.48\textwidth}
		\centering
		\includegraphics[width=.8\linewidth]{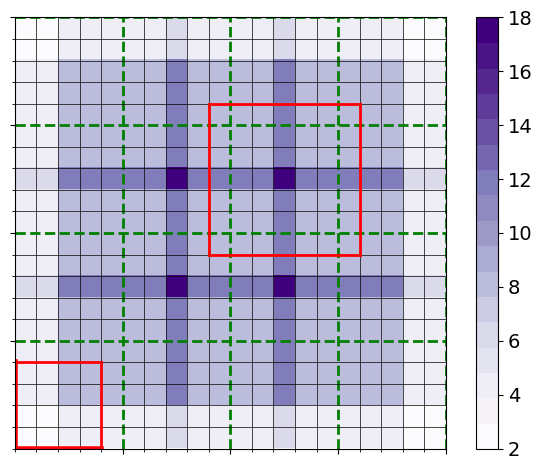}  
	\end{subfigure}
	\begin{subfigure}{.48\textwidth}
		\centering
		\includegraphics[width=.8\linewidth]{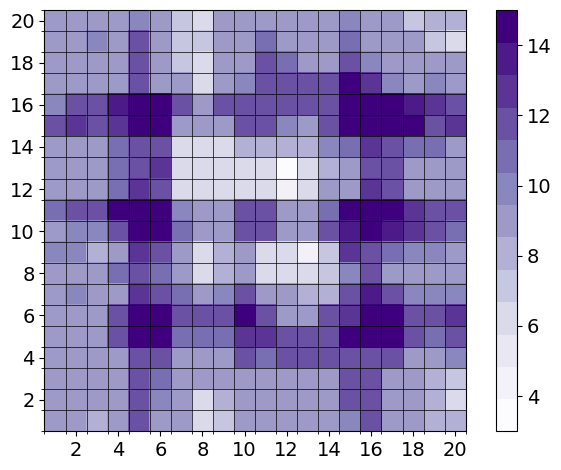}  
	\end{subfigure}
	\caption[Experiment 12: Number of affine coefficients and local RB sizes]{Left: Number of affine coefficients in each patch problem $T$ for $n_H=20$ and $\kl=3$.
		The thermal block structure of $A^{1,\xi}_\mu$ and $A^{2,\xi}_\mu$ is highlighted in green, two patch instances are highlighted in red. Right: local corrector RB sizes of the Stage~1 models after successful termination of Method 4.}
	\label{fig:aff_coefs}
\end{figure}
In what follows, we consider an experiment where FEM solves are computationally affordable.
To this end, we set the resolution of the multiscale coefficients to $N_1 = 150$ and $N_2 = 300$.
For the fine-mesh, we thus choose $n_h = 1200$ to ensure at least $4$ quadrilateral grid cells in each of the rapidly varying multiscale features.
Therefore, the FEM mesh has $1.4$ Mio degrees of freedom.
For the coarse-grid, we choose $n_H=20$, which results in only $400$ coarse grid cells and, in particular, $\kl=3$ and $176.400$ fine-mesh elements for full patches $U_\kl(T)$.
Concerning, the objective functional, we compute $u_{\mu^{\text{d}}}$ as the FEM solution of \eqref{eq:state_h} for $\mu^{\text{d}}$.

We emphasize that for this experiment, \cref{asmpt:truth} is not entirely fulfilled.
Instead, an approximation error of the PG--LOD in $\Jhat$ at the desired parameter $\mu^{\text{d}}$ is still observable, and we have
$$
| \Jhat_h^{\text{loc}}(\mu^{\text{d}}) - \Jhat_h(\mu^{\text{d}})| = 8.25 \cdot 10^{-6}.
$$
Although this violates \cref{asmpt:truth}, we can expect that all methods converge up to the LOD discretization error, which is sufficiently close for this experiment.

In \cref{fig:aff_coefs}(left), we visualize the number of affine components of the local corrector models that are directly associated with the number of affine components in $A_\mu$.
Thus, the number of affine components can be determined by the number of thermal blocks that lie in the patch.
The thermal blocks are highlighted in green, and since $\kl=3$, the resulting affine components can be counted in the plot.
For instance, the lower-left element's patch only reaches the lower-left thermal block (resulting in $2$ affine coefficients). Moreover, the elements that directly lie inside the inner thermal blocks have a patch that reaches up until all $9$ neighboring blocks (resulting in $18$ affine components each).
The discussed patch instances are highlighted in red in \cref{fig:aff_coefs}.
We conclude that the corrector problems have a more minor parameter dependence than globalized RB methods.
In turn, the local RB models can be expected to require less basis functions.

In \cref{fig:aff_coefs}(right), the final local RB size of the Stage~1 models is depicted. It can be seen that the model requires a relatively rich space at the coarse elements that are close to the "jumps" in the desired parameter, cf~\cref{fig:thermal_block}.
As expected, the low conductivity blocks in the middle of the domain do not require many RB enrichments since the optimization problem in these blocks is less demanding.
In addition, from solely looking at \cref{fig:aff_coefs}(left), one would guess that the local patch problems that admit the highest number of affine components require the most basis functions.
The fact that this expectation is not valid further proves that the optional enrichment plays a significant role in the algorithm.

We note that all compared methods indeed converged up to the chosen tolerance $\tau_\text{{FOC}}$ to the same point, and it was posteriorly verified that the point is indeed a local optimum.
In \cref{tab:TRRBLOD_1} and \cref{tab:TRRBLOD_2}, we report relevant information on the evaluation counts, the iteration, and the observed run times.

\begin{table}[h] \centering \footnotesize
	\begin{tabular}{|l|c|c|c|c|c|c|c|c|c|} \hline
		& \multicolumn{2}{c|}{} & \multicolumn{2}{c|}{LOD} 	& \multicolumn{2}{c|}{Stage 1} & TS &
		\multicolumn{2}{c|}{} \\ \hline
		Evaluations 		& FEM 		   & RB  & Coa.		& Local 	  & Coa. 	& Local & & O.it. & Time 
		\\ \hline  \hline
		Cost factor   & $h$     & $N_{\text{RB}}$ & $H$ & $U(T_{H})_{h}$ & $H$  & $N_{\text{RB}}$  & $N_{\text{RB}}$ & 
		\multicolumn{2}{c|}{} \\ \hline
		1.\hphantom{a} FEM 	    		
		& \textbf{280} & -   & -   & -  			& -    	& -   	 & -  & 92 & 11402s
		\\
		2.\hphantom{c} TR-RB   		
		& \textbf{12}  & 1546 & -   & -  			& -    	& -   	 & -  & 5  & \hphantom{0}1226s
		\\ 
		2.r R-TR-RB   		
		& \textbf{10}  & 862 & -   & -  			& -    	& -   	 & -  & 4  & \hphantom{00}410s
		\\ \hline\hline
		3.\hphantom{a} PG--LOD 			
		& -  		   & -   & 244 & \textbf{131200}& -   	& -   	 & -  & 80 & \hphantom{00}620s
		\\
		4.\hphantom{c} TR-TS  		
		& -  		   & -   & 10   & \textbf{8000}	& {30}  	& 48000  & \textbf{910}	& 6  & 		\hphantom{00}579s
		\\
		4.r R-TR-TS 		
		& -  		   & -   & 8   & \textbf{6400}	& {30}  	& 48000  & \textbf{316}	& 5  & 		\hphantom{00}272s 
		\\
		5.\hphantom{r} TR-RBLOD 		
		& -  		   & -   & 8   & \textbf{6400}  & \textbf{682} 	& 208000 & -  & 5  &
		\hphantom{00}862s
		\\
		5.r R-TR-RBLOD 		
		& -  		   & -   & 4   & \textbf{4800}  & \textbf{347} 	& 208000 & -  & 3  & 		\hphantom{00}345s
		\\		
		\hline
	\end{tabular}
	\caption{Experiment 12: Evaluations and timings of selected methods.}
	\label{tab:TRRBLOD_1}
\end{table}

\begin{table}[h] \centering \footnotesize
	\begin{tabular}{|l|c|c|c|c|c|c|c|c|} \hline
		& 		& 		  & \multicolumn{2}{c|}{Online} 	& \multicolumn{3}{c|}{Offline} & \\ \hline
		Method 				  & Total 			  & Speedup & Outer & Inner 					
		& FEM & Stage 1	 &  Stage 2  & $e^{\Jhat_h,\text{rel}}$ 
		\\ \hline  \hline
		1.\hphantom{a} FEM 	  & 11402s 			  & -  &  11402s			& -  			  
		& - & - & - & 4.18e-10 	\\
		2.\hphantom{c} TR-RB 			  & \hphantom{0}1226s & 9  &  \hphantom{000}41s	& \hphantom{00}3s
		& 1165s & - & - & 5.37e-11 	\\
		2.r R-TR-RB			  & \hphantom{00}410s & 28 &  \hphantom{000}34s	& \hphantom{00}4s
		& \hphantom{0}348s & - & - & 7.75e-10 	\\ \hline \hline
		3.\hphantom{a} PG--LOD & \hphantom{00}620s & 18 &  \hphantom{00}620s	& - 			 
		& - & - & - & 4.22e-06 	\\
		4.\hphantom{c} TR-TS			  & \hphantom{00}579s & 20 &  \hphantom{000}41s	& 100s 			  
		& - & 311s & 127s & 4.22e-06 	\\
		4.r R-TR-TS			  & \hphantom{00}272s & 42 &  \hphantom{0000}4s	& \hphantom{00}4s 
		& - & 193s & \hphantom{0}71s & 4.22e-06 	\\
		5.\hphantom{r} TR-RBLOD	& \hphantom{00}862s	& 13 &  \hphantom{000}55s & 520s 	& -		& 287s 	& -   & 4.22e-06	  \\
		5.r R-TR-RBLOD	& \hphantom{00}345s	& 33 &  \hphantom{0000}4s & 236s  & -			& 105s 	& -   	& 4.22e-06 		  \\
		\hline
	\end{tabular}
	\caption{Experiment 12: More details on run times and accuracy of selected methods.}
	\label{tab:TRRBLOD_2}
\end{table}

From \cref{tab:TRRBLOD_1}, we conclude that all Methods~2-5 give a significant speedup to the standard FEM Method~1.
Although $92$ iterations of Method~1 and the corresponding $280$ FEM evaluations are relatively few for a $32$-dimensional optimization problem, the computational effort required to perform a FEM solution with $1.4$ Mio. DoFs harm the speed of the method.
As shown in \cref{chap:TR_RB}, the TR-RB Method~2 is mainly designed to avoid expensive FEM evaluations and converges already after $4$ and $5$ outer iterations, which only requires $12$ and $10$ FEM-based enrichments of the reduced spaces.
On the other hand, the $1546$ and $862$ inner RB evaluations are cheap and do not harm the computational speed of the method.

As expected, the localized methods only converge until the priorly known approximation error of the PG--LOD is reached.
However, it can be seen that the TR-LRB methods find the same point and are not subject to severe approximation issues.

A significant reason why the TR-TSRBLOD method is particularly suitable for the TR procedure work is its very efficient online phase.
This result can be obtained by \cref{tab:TRRBLOD_2}, where extended timings are given for the TR-TSRBLOD method compared to the TR-RBLOD.
Just as it is the case for the global TR-RB methods, only a few seconds are required to solve the sub-problems in the relaxed variant, independent of the coarse LOD mesh.
Instead, the sub-problem is more demanding for the fully certified variant, which goes back to the fact that we did not afford the offline time to reduce the respective error estimator in Stage~2.
Moreover, we see that the TR-RBLOD generally has a more demanding online time but saves time due to fewer iterations and less offline time thanks to the missing Stage~2 reduction.

It can further be noticed that the localized methods show a comparably good convergence speed w.r.t. the FEM-based methods, although FEM is still comparably fast.
We also see that Method~3, considered the localized FOM, shows a strong convergence speed.
This is due to the relatively small patch problems such that the localized corrector problems and the corresponding Stage~2 reduction do not pay off massively.

The relaxed versions of the TR-(L)RB methods show a remarkably fast convergence behavior in this experiment.
One reason for this is that the fully enforced certification cannot detect the full benefit from the surrogate model and truncates the sub-problems too early.
With the specific relaxation sequences, the R-TR methods unconditionally trust the first surrogate models, allowing the overall algorithm to converge with fewer outer iterations.
On top of that, the enrichment time is significantly lower, which goes back to the left-out pre-assembly of the error estimates.

In conclusion, FEM based-methods can reliably be replaced by localized methods already for moderately small fine-mesh sizes.
Certainly, the full benefit of the TR-LRB approaches can only be deduced for scenarios where the PG--LOD is costly in itself, shown in the next experiment.

\subsection{Experiment 13: Large scale example}

\label{sec:exp_13}

We consider a large-scale example where the global FEM mesh does not fit into the machine's memory.
We set the multiscale resolution to $N_1 = 1.000$ and $N_2 = 250$.
For the fine-mesh, we again aim to have at least $4$ fine mesh entities in each multiscale cell and choose $n_h = 4000$.
Therefore, the FEM mesh would have $16$ Mio degrees of freedom, which can be considered prohibitively large.
Hence, we do not utilize FEM-based methods and only compare Methods 3 and 4.
For the coarse-grid, we choose $n_H=40$, which results in $1600$ coarse-grid cells, $\kl=4$, and $810.000$ fine-mesh elements for full patches $U_\kl(T)$.
Since FEM evaluations are not available, the desired solution $u_{\mu^{\text{d}}}$ is computed with the PG--LOD, i.e. we solve \eqref{eq:PG_LOD} for $\mu^{\text{d}}$.

\begin{figure}
	\centering
	\begin{subfigure}{.48\textwidth}
		\centering
		\includegraphics[width=.8\linewidth]{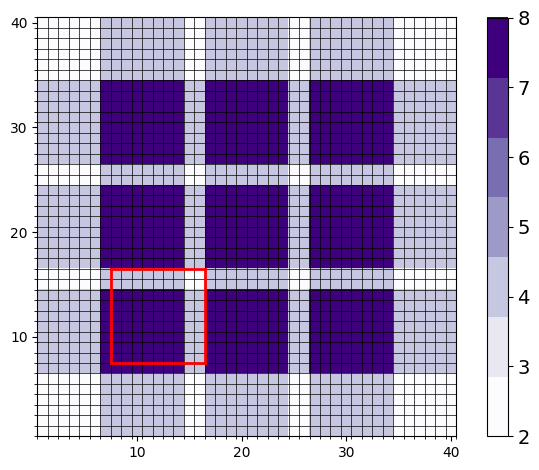}  
	\end{subfigure}
	\begin{subfigure}{.48\textwidth}
		\centering
		\includegraphics[width=.8\linewidth]{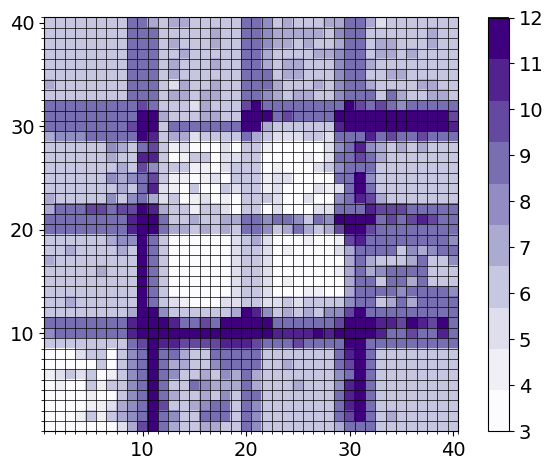}  
	\end{subfigure}
	\caption[Experiment 13: Number of affine coefficients and local RB sizes]{Left: Number of affine coefficients in each patch problem $T$ for $n_H=40$ and $\kl=4$.
		One patch instance is highlighted in red. Right: local corrector RB sizes of the Stage~1 models after successful termination of Method 4.}
	\label{fig:aff_coefs_2}
\end{figure}

Similar to the above illustrations, in \cref{fig:aff_coefs_2}, we report the respective number of affine components of the patch problems as well as the final local RB sizes of the certified TR-TSRBLOD method with optional enrichment (Method 4).
In particular, \cref{fig:aff_coefs_2} can be interpreted as the refined version of \cref{fig:aff_coefs}, where it is even more visible that the local corrector problems have fewer affine components and need more basis functions for the corrector problems that are largely affected by the "jumps" in the desired thermal block state, depicted in \cref{fig:thermal_block}.
Just as before, it can be seen that the amount of basis functions is also associated with the intensity of the respective "jumps", and the low conductivity in the middle of the domain is well visible.

In \cref{tab:TRRBLOD_3} and \cref{tab:TRRBLOD_4}, we again provide an extensive comparison concerning evaluations, run time, and iteration counts of the methods.
It can be seen that the (TS)RBLOD-based methods successfully reduce the computational effort of Method 3, which is mainly due to the increasing number of fine-mesh DoFs in the patches.
With increasing complexity of the multiscale problem, we thus expect even more speedups.
We also note that the speedup w.r.t. the FOM method is also dependent on the outer iteration counts, cf. \cite{paper1,paper2}.
It can be expected that the benefit of reduced models is even more present for increasing complexity of the optimization problem.

At the same time, as mentioned before, in \cref{tab:TRRBLOD_4} the results confirm that our specific implementation leaves room for improvements concerning the construction of the Stage~2 model, which has an additional slow-down effect on the computational time.

\begin{table}[h!] \centering \footnotesize
	\begin{tabular}{|l|c|c|c|c|c|c|c|c|c|} \hline
		& \multicolumn{1}{c|}{} & \multicolumn{2}{c|}{LOD} 	& \multicolumn{2}{c|}{St.1} & St.2 &
		\multicolumn{3}{c|}{} \\ \hline
		Evaluations 		& FEM 		   & Coa.		& Local 	  & Coa. 	& Local & & O.it. & Time 
		& $e^{\Jhat_h^\textnormal{loc},\text{rel}}$
		\\ \hline  \hline
		Cost factor   & $h$     & $H$ & $U(T_{H})_{h}$ & $H$  & $N_{\text{RB}}$  & $N_{\text{RB}}$ &
		\multicolumn{3}{c|}{} \\ \hline
		3.\hphantom{a} PG--LOD
		& -  		     & 307 & \textbf{665600}& -   	& -   	 & -  & 100 & 11317s
		& 2.07e-10
		\\
		4.\hphantom{c} TR-TS
		& -  		     & 10   & \textbf{32000}	& \textbf{20}  	& 128000  & \textbf{938} & 5  & \hphantom{0}5393s
		& 1.57e-10
		\\
		4.r R-TR-TS
		& -  		     & 8   & \textbf{51200}	& \textbf{30}  	& 128000  & \textbf{422} & 4  & \hphantom{0}2998s
		& 1.32e-11
		\\
		5.\hphantom{c} TR-RBLOD
		& -  		     & 10   & \textbf{32000}	& \textbf{676}  	& 793600  & - & 5  & \hphantom{0}6423s
		& 6.11e-12
		\\
		5.r R-TR-RBLOD
		& -  		     & 6   & \textbf{25600}	& \textbf{402}  	& 851200  & - & 3  & \hphantom{0}2260s
		& 1.68e-12
		\\
		
		\hline
	\end{tabular}
	\caption{Experiment 13: Evaluations and accuracy of selected methods.}
	\label{tab:TRRBLOD_3}
\end{table}

\begin{table}[h!] \centering \footnotesize
	\begin{tabular}{|l|c|c|c|c|c|c|c|} \hline
		& 		& 		  & \multicolumn{2}{c|}{Online} 	& \multicolumn{3}{c|}{Offline}  \\ \hline
		Method 				  & Total 			  & Speedup & Outer & Inner 					
		& FEM & Stage 1	 &  Stage 2 
		\\ \hline  \hline
		3.\hphantom{a} PG--LOD & 11317s & - &  11317s	& - 			 
		& - & - & - 
		\\
		4.\hphantom{c} TR-TS & \hphantom{0}5393s & 2 &  \hphantom{00}486s	& \hphantom{0}360s  & - & 4042s & 505s
		\\
		4.r R-TR-TS			 & \hphantom{0}2998s & 4 &  \hphantom{0000}5s	& \hphantom{000}5s & - & 2690s & 336s
		\\
		5.\hphantom{c} TR-RBLOD & \hphantom{0}6423s & 2 &  \hphantom{00}552s	& 1524s  & - & 4347s & -
		\\
		5.r R-TR-RBLOD			 & \hphantom{0}2260s & 5 &  \hphantom{000}7s	& \hphantom{0}689s & - & 1564s & -
		\\
		\hline
	\end{tabular}
	\caption{Experiment 13: More details on run times of selected methods.}
	\label{tab:TRRBLOD_4}
\end{table}

\subsection{Experiment 14: RBLOD vs. TSRBLOD}

As the last experiment in this chapter, we showcase that the RBLOD indeed struggles for increasing number of coarse elements in the LOD discretization and that the TSRBLOD still lacks a more efficient implementation of Stage~2.
To this end, we perform the exact same experiment as in Experiment 12 but with $n_h=600$ for the fine mesh and $n_H = 60$ for the coarse mesh, such that $3600$ corrector problems need to be constructed and the two-scale matrix is growing accordingly.
Note that this discretization choice only serves as a proof-of-expectation.
It is clear that the FEM-based method will be accessible and fast (since the FE mesh has $361.201$ DoFs).
For this reason, we omit an elaborated comparison between all methods in this section and only compare the TR-(TS)RBLOD algorithm.
More so, we neglect the fully certified algorithms and only consider the relaxed variants in order to concentrate on the offline- and online times that can be obtained without the construction of error estimates.
\vspace{0.4cm}
\begin{table}[h!] \centering \footnotesize
	\begin{tabular}{|l|c|c|c|c|c|c|c|} \hline
						& 			  & 		  &\multicolumn{2}{c|}{Online time} 	& \multicolumn{2}{c|}{Offline time} &  \\ \hline
		Method 			& Outer it.	  & Total 			   & Outer & Inner 					
		& Stage 1	 &  Stage 2  & $e^{\Jhat_h^\textnormal{loc},\text{rel}}$ 
		\\ \hline  \hline
		4.r R-TR-TS & 4 & 1245s  &  13s	& \hphantom{000}5s   & 389s & 838s & 5.13e-11
		\\
		5.r R-TR-RBLOD & 3 &2742s &   25s	& \textbf{2743s}  & 269s & - & 3.64e-14
		\\
		\hline
	\end{tabular}
	\caption{Experiment 14: Run times and accuracy of R-TR-TSRBLOD and R-TR-RBLOD method.}
	\label{tab:TRRBLOD_5}
\end{table}
\vspace{0.1cm}

In \cref{tab:TRRBLOD_5}, we illustrate the performance of the variants.
We observe that the missing online efficiency of the RBLOD-based method for increasing number of coarse elements.
This complexity is not visible in the online time of the TSRBLOD-based method.
However, the Stage~2 reduction remains vulnerable to the growing size of the two-scale system matrix.
Although it can already be seen that the slow-down effect of the TSRBLOD is less severe, the observed run time still indicated room for improvements in the respective implementation; cf. the discussion~\cref{sec:exp_13}.

\section{Summary and outlook} \label{sec:TR_TS_conclusion}

In this chapter, we combined localized reduced basis methods for efficiently solving parameterized multiscale problems with optimization methods that adaptively construct such localized reduced methods in the context of an iterative error-aware trust-region algorithm for accelerating PDE-constrained optimization.
We also devised the relaxed version of the TR-RB algorithm from \cref{sec:Relaxed_TRRB} to neglect the strong certification in the first iterations.

To this end, we discretized the optimality system of the PDE-constrained optimization problem \eqref{P} with a localized ansatz based on the Petrov-Galerkin version of the localized orthogonal decomposition method.
For an online efficient reduced model with optional local basis enrichment, such that the sub-problems of the TR algorithm can be solved fast, we used the RBLOD and TSRBLOD reduction approach for the LOD from \cref{chap:TSRBLOD}.

The resulting TR-LRB method has proven advantageous both in terms of computational effort and adaptivity concerning the localized RB models.
In the numerical experiments, we observed that these localized RB approaches can efficiently replace FEM-based techniques, especially for growing complexity of the multiscale system.

Although the underlying multiscale data is indeed highly heterogeneous and non-periodic, and the LOD approach showed good approximation properties w.r.t. FEM, it is commonly known that the LOD struggles, e.g., for high-contrast problems or complex coarse data such as thin channels.
For using the TR-TSRBLOD, it has to be verified priorly that \cref{asmpt:truth} is given up to an acceptable tolerance.
To remedy this, the discussed concepts can be generalized to other multiscale methods, always dependent on the respective multiscale task.
On the other hand, it would be desirable to derive a posteriori error theory for the LOD such that the homogenization term from \eqref{eq:hom_red_term} can be used to validate the approximation properties of the LOD, cf.~\cref{rem:trunc-red-hom-term}.
Let us also mention that, in the work at hand, we have enforced several problem assumptions, e.g., ellipticity, symmetry, and homogeneous boundary conditions, to simplify the presentation of the algorithm.
However, it is straightforward to generalize the methodology to more challenging problem classes.

Concerning the specific instance of the TR-TSRBLOD, the numerical experiments showed a significant overall speedup w.r.t. FEM and the PG--LOD with our implementation, and further improvements are possible.
Our theoretical findings expect better run times with an even more HPC-oriented implementational design.
Moreover, an intermediate preparatory reduction of the two-scale system can be used to decrease further offline expenses of Stage~2.

Lastly, the described TR-LRB can be enhanced in terms of the choice of the local enrichment tolerance $\tau_{{\text{loc}}}$, such that an appropriate choice for the respective optimization problem or model can efficiently be found with an adaptive refinement, cf.~\cref{sec:local_basis_enrichment}.
\chapter{Conclusion and Outlook}\label{chap:outlook}

\noindent The last part of this work is devoted to concluding the results obtained in this thesis with particular emphasis on the prescribed research goals and the connection of the individual content.
Lastly, we elaborate on different research perspectives directly associated with the presented concepts or generalizations to further problem classes.
Note that some specialized summaries and outlooks were already given in the respective chapters.

\section{Conclusion}
In this work, we presented and rigorously analyzed significant advances in model order reduction used in the context of parameterized multiscale problems as well as for large-scale PDE-constrained parameter optimization.

In \cref{chap:background}, we elaborated on the state-of-the-art concerning reduced basis methods and their application to PDE-constrained parameter optimization problems.
In addition, we motivated future research challenges concerning the curse of dimensionality and the inaccessibility of global discretization schemes, such as the finite element method.
In the context of PDE-constrained optimization problems, one of the primary goals of this thesis was to derive overall-efficient methods for finding a local solution to a single optimization problem, meaning that the offline time for constructing a surrogate model cannot be ignored.

As shown in \cref{chap:TR_RB}, the trust-region reduced basis (TR-RB) algorithm is an overall-efficient algorithm since a dedicated surrogate model is built progressively along the optimization path, whereas FOM evaluations are avoided as often as possible.
In detail, the initial surrogate model is only constructed with the initial guess of the iterative optimization method.
The TR method then internally solves an inner sub-problem that only accepts iterates in a local area of the parameter set, in which the RB model is trustable.
This region can be detected quickly by using the a posteriori error estimator of the RB model.
If the iterate is accepted with respect to certified conditions, the surrogate model is enriched further and the next inner sub-problem is started unless the TR-RB algorithm is reliably terminated.

While such TR-RB methods were already proposed in \cite{QGVW2017}, in \cref{chap:TR_RB}, we explained and demonstrated significant enhancements of the TR-RB algorithm, also by generalizing the concept to additional parameter constraints.
These include, for instance, a non-conforming dual (NCD) correction term for the reduced model, higher-order sub-problem solvers, enlargement of the TR radius, FOM-based stopping criteria, parameter control, and optional basis enrichment; cf. \cref{sec:TR_algorithm_details}.
An appropriate convergence analysis was presented in \cref{sec:TRRB_convergence_analysis}.

For the corresponding numerical experiments in \cref{sec:TRRB_num_experiments}, we introduced a flexible, possibly high-dimensional, benchmark optimization problem for with an underlying elliptic PDE.
With these benchmark problems, it was possible to demonstrate the remarkably robust behavior of the proposed TR-RB algorithm concerning many aspects.
As thoroughly summarized in \cref{sec:TR_RB_conclusion}, the experiments indicated that the suggested algorithms cover many different cases of optimization problems and that significant improvement in comparison to \cite{QGVW2017} were achieved.

Apart from the basic TR-RB algorithm that was used in \cite{paper2,paper1,paper3}, in \cref{sec:TRRB_related_approaches}, we also discussed several more advances concerning the convergence speed of the outer loop of the algorithm.
As examples, we discussed mesh adaptivity for dropping \cref{asmpt:truth}, coarsening of the RB spaces to avoid large RB spaces, and an updated surrogate model after parameter rejection to take maximum advantage of FOM evaluations.
Notably, we presented the relaxed TR-RB (R-TR-RB) approach that can be used for a computationally faster procedure since it relaxes the certification of the TR method in the first iterations but admits the same convergence result asymptotically.
Despite applying this approach in the context of \cref{chap:TR_TSRBLOD}, we already conducted a proof-of-concept experiment demonstrating that such further variants can lead to a significant speedup of the algorithm.

This thesis also concerns the fact that spatially global schemes like the finite element method can become infeasible, which was thoroughly motivated by \cref{sec:multiscale_problems} and further formalized in \cref{sec:RB_challenges}.
To this end, in \cref{chap:TSRBLOD} and \cref{chap:TR_TSRBLOD}, we considered localized model reduction approaches.

In \cref{chap:TSRBLOD}, we explained recent developments for the localized orthogonal decomposition (LOD) multiscale method applied to many-query and real-time scenarios of parameterized multiscale problems.
We introduced the two-scale reduced basis localized orthogonal decomposition (TSRBLOD) method, which internally benefits from the reduced basis localized orthogonal decomposition (RBLOD) method from \cite{RBLOD}.
In the TSRBLOD, the primary motivation was to develop an approach for parameterized multiscale problems that, in the online phase, is independent of the fine mesh and, additionally, of the coarse mesh.
Hence, no loops over coarse-mesh elements are required anymore.
On top of that, in the TSRBLOD, we were able to efficiently control the full reduction error with respect to the true LOD, which was not yet discovered by \cite{RBLOD}.
Numerical experiments indeed showed remarkable applicability of the TSRBLOD to large-scale problems with massive online speedups w.r.t.~both the LOD and RBLOD; see also the discussion in \cref{sec:TSRBLOD_conclusion}.

The TSRBLOD has proven to be an efficient localized model order reduction method for multiscale problems. Furthermore, the a posteriori error estimation in \cref{chap:TSRBLOD} showed that it is also possible to construct the local RB problems for both the TSRBLOD as well as the RBLOD in an online adaptive manner.
Therefore, using the respective reduced schemes in a trust-region localized reduced basis TR-LRB algorithm for parameter optimization problems constrained by multiscale PDEs seemed natural.
Importantly, such an algorithm resolves the RB challenge of an inaccessible global discretization for the TR-RB algorithm devised in \cref{chap:TR_RB}.

For the successful application of the TR-LRB algorithm, we used suitable assumptions concerning the approximation quality of the localized FOM approach in \cref{asmpt:loc_truth}, the error estimation for the objective functional, and the parameter separability from \cref{asmpt:parameter_separable}.
In particular, the affine decomposition of the data functions was essential for achieving low storage requirements of the reduced method and played a significant role in the efficiency.
Moreover, as usual for multiscale problems, we assumed that the underlying optimization problem is focused on a coarse-scale behavior of the objective functional; cf. \cref{asmpt:coarse_J}.
The a posteriori error analysis of \cref{chap:TSRBLOD} enabled us to efficiently characterize the trust-region for both the RBLOD and the TSRBLOD.

In the numerical experiments in \cref{sec:loc_experiments}, we introduced a benchmark problem for multiscale scenarios and started with a case where the finite element method was still accessible.
While the relaxed TR method already showed its strength in the 2-dimensional small experiment from \cref{sec:TRRB_related_approaches}, it also proved computationally highly advantageous for the depicted 32-dimensional multiscale example, both for the localized and global discretization schemes.
The relaxed TR method is especially well-suited for the localized models, where the actual approximation error of the surrogate model can be overestimated, harming the fully certified methods in their convergence speed.
More so, the left out offline-preparation of the estimators in the relaxed TR method also plays a significant role in less computational effort in the (localized) RB scheme.
The localized approaches showed a promising and robust convergence behavior.
Moreover, our implementation sufficiently took advantage of parallelization in the localized schemes.
In conclusion, the results highlighted that these methods could efficiently replace the FEM-based methods, already in scenarios where FEM is accessible.
Last but not least, we demonstrated the effect of online enrichment, where some local models only needed a few basis functions, especially in the regions where the optimizer was less active.

The classical PG--LOD, which was used as the localized FOM method, benefited quite heavily from the parallelization and was thus already very fast without reduction.
We also elaborated on the case where the patch problems incorporate many DoFs, comparable to standard FEM meshes.
A first hint towards a substantial speedup of the localized reduced schemes in these cases was given in the second multiscale experiment, where FEM was considered inaccessible.
In the last experiment, we saw that the missing online efficiency of the RBLOD for large coarse meshes harms the computational speed of the associated TR-LRB method.
Instead, as expected from the experiments in \cref{chap:TSRBLOD}, the TSRBLOD resolves this problem, despite a relatively large additional offline serialized construction.
We also point to \cref{sec:TR_TS_conclusion} for more details.

\section{Outlook}
\noindent Many research questions are left for the future, both in optimization methods for different applications and in (localized) model order reduction schemes and their combination.

The presented trust-region methods are very flexible concerning the choice of the full-order model, its surrogate, or other algorithm-related features.
Concerning the continued development of the specific TR-RB algorithm, as also concluded above, further hints towards more involved algorithms were given in \cref{sec:TRRB_related_approaches}, e.g., coarsening of the RB spaces or more FOM-cost oriented procedures.

Importantly, further research perspectives of MOR-informed TR methods are usually unaffected by the choice of the (localized) full- and reduced-order model since these methods can be interpreted as a black-box reduced approach as long as the models can be trusted concerning their approximation qualities.
Indeed, TR methods are also applicable to many other problem classes or other types of parameter optimization.

Concerning the variational formulation in \cref{def:elliptic_problem}, which we restricted to the elliptic case, several different applications are possible.
Model order reduction techniques and a posteriori error estimates are, for instance, also well-established for parameterized parabolic problems \cite{grepl2005posteriori, karcher2014posteriori,Ohlberger2017169,QGVW2017}, and hence, it is relatively straightforward to transfer the TR-RB algorithm to such scenarios.

For illustration, this thesis strongly utilized the reduced formulation \eqref{Phat} of the optimization problem \eqref{P}, which is internally based on the fact that the primal equation is uniquely solvable and the parameter-to-state mapping $\mathcal{S}(\mu)$ exists.
However, in a more general setting, the optimal solution pair is not uniquely determined by $\bar{\mu}$.
On the contrary, so-called \emph{all-at-once} approaches assume the control and state variables to be independent \cite{Clever2012,ganzler2006sqp,HV06}.
A popular method is the sequential quadratic programming (SQP) method, for which also TR approaches can be used for a globally convergent scheme.
Again, the involved model functions can be reduced by RB techniques similar to this thesis, with the aim to reduce the overall computational cost of the algorithm.

The algorithms in this thesis were used for finding a local minimum of an optimization task, where the existence of an optimum was verified, but no uniqueness was given.
If one is instead interested in a global minimum, an algorithm with multiple different starting parameters with interactions of different ROMs can be employed.
In such cases, FOM information and reduced models may also be reused.
A related problem class where several optimization loops are required are multi-objective PDE-constrained optimization problems; see~\cite{banholzer_phd, banholzer2022trust} and the references therein.
In such cases, coarsening of RB spaces has already proved to be advantageous; see also the discussion in \cref{sec:coarsening_RB_approach}.

Another prevalent class directly related to parameter optimization are inverse problems or parameter estimation, where the forward mapping is given by a nonlinear ill-posed operator that maps from a parameter set to the state space.
For such problems, we particularly mention iteratively regularized Gauss-Newton methods (IRGNM) \cite{kaltenbacher2018convergence}.
The involved forward operators and suitable regularization terms can once again be modeled by respective surrogate models, making it possible to transfer ideas from this thesis.  

A central assumption enforced in our work is the accuracy of the respective full-order model, and a posteriori error control of the respective FOM has not been taken into account; cf. Assumptions~\ref{asmpt:truth} and~\ref{asmpt:loc_truth}.
In \cref{sec:afem_approaches}, we discussed the possibility to drop such an assumption and elaborated on future perspectives of algorithms that do not trust the FOM unconditionally.
Instead, such a method starts with an intentionally coarse discretization and refines it adaptively.
For FEM, it is relatively straightforward to devise such an AFEM-based algorithm since the respective theory already exists.
In this context, we mention the recently developed $hp$-AFEM, advised in \cite{CNSV2016a}.
Open questions concerning the prolongation of snapshots from coarse-mesh reduced models during the optimization routine are still to be answered.
At the same time, the choice of a sufficiently good mesh size $h$ for FEM is usually easy to find, and hence, \cref{asmpt:truth} is not very restrictive, as long as the respective number of DoFs fits into the memory of the machine and the computational effort is acceptable.

For localized FOMs, however, the situation is more challenging.
In the absence of a reliable a posteriori theory for the particular approach, \cref{asmpt:loc_truth} can only be verified by prior knowledge of the system or by a priori results.
Admittedly, these a priori results are usually less sharp, and the individual choice of mesh sizes and localization parameters can be difficult.

For the specific case of the PG--LOD, a posteriori results for the FOM are left to further research, and adaptive mesh refinements were also not considered.
Furthermore, we have only shown a multiscale setup with standard assumptions.
The generalization of the devised methods to more challenging multiscale tasks is thus also on further due; cf. \cref{sec:TSRBLOD_conclusion} and \cref{sec:TR_TS_conclusion}.

From a priori theory of the PG--LOD, we know that there exist problem classes where the LOD has severe approximability issues, such as thin high-conductivity channels with possibly high contrast.
Although future work is devoted to fixing these problems in the LOD, future research must also be concerned with combining the TR-LRB ideas with more flexible localized reduced schemes.
This is especially of interest for non-elliptic or nonlinear problems.
In a nonlinear setting, a key ingredient is the empirical interpolation method, already mentioned for treating non-parameter separable data functions.
Suitable online adaptive reduced-order models can also be derived for such cases.
For a recent result in the context of time-dependent problems, we refer to \cite{sleeman2022goal}.

A particular promising and flexible approach to the problems considered (and beyond) is the localized reduced basis multiscale (LRBMS) method \cite{AHKO2012,OS2017} which relies on a DD-based discontinuous Galerkin scheme with the idea of constructing a reduced basis for each element of the domain decomposition.
The great advantage of this approach is that both the domain decomposition itself as well as the corresponding local basis functions can be tailored towards the problem characteristics.
For instance, the local basis functions of multiscale methods such as the LOD or GFEM can be used to perform the reduced method's online enrichment.
Promising work of quasi-optimal local spaces has been discussed, e.g., in \cite{BIORSS21}.
However, many questions regarding the efficient application of the LRBMS in TR methods are yet to be addressed, such as the choice of the optimization approach, the construction of the right-hand side of the dual model (if applicable), and the efficient derivation of a suitable error estimator.
The latter has already been achieved in the context of the ArbiLoMod \cite{BEOR2017}.
If employed for overall-efficient optimization methods, it is also of interest to assess the offline and online computational expenses of the LRBMS.
Moreover, it holds further interest whether the approach requires extensive offline training to achieve an appropriate (local) accuracy, which is not finally answered for the LRBMS.
A suitable online enrichment has already been presented in \cite{OSS2018,OS2017}.

Last but not least, future work is also devoted to more flexible, HPC-oriented, and, most importantly, sustainable software concepts that allow for a fast transition to more problem classes, applicable for more researchers that work on related topics.
\backmatter 
{\small
	\bibliographystyle{abbrv}

}
\cleardoublepage 


\appendix
\chapter{Appendix}\label{chap:appendix}
\section{Code availability}
\label{apx:code_availability}
Most of the program code for the numerical experiments of this thesis has been published along with the respective articles.
We also refer to \cref{sec:software_TRRB}, \cref{sec:software_TSRBLOD}, and \cref{sec:software_TR_TSRBLOD} for more information on the used software packages and the implementational design.
To avoid an overload of referenced code fragments, we only cite the supplementary code for this thesis in~\cite{code:thesis}.
In the respective README, we have listed other related repositories and provided very simple setup instructions.
For brevity, we point to the respective code reference for each numerical experiment in this thesis.

\begin{description}
	\item[\cref{chap:background}] The only numerical experiment in this chapter was used as a motivation for adaptive RB models in the context of PDE-constrained parameter optimization in \cref{sec:background_RB_methods_for_PDEopt}.
	The 2-dimensional experiment was depicted from \cref{sec:proof_of_concept}, where a proof-of-concept experiment is presented.
	The entire experiment has not been published along with any of the articles but can be found in the supplementary code to this thesis~\cite{code:thesis}.
	\item[\cref{chap:TR_RB}:] As explained above, the proof-of-concept experiment that is used in \cref{sec:proof_of_concept} can be found in~\cite{code:thesis}.
	The same code has been used for the numerical experiments for the R-TR-RB and the FCO-TR-RB algorithms in \cref{sec:experiments_further_approaches}.
	
	The rest of the experiments has been published within the cited articles.
	To recall, the code reference for \cref{sec:exp_paper1}, \cref{sec:exp_paper2}, and \cref{sec:exp_paper3} can be found in \cite{paper1}, \cite{paper2}, and \cite{paper3}, respectively.
	\item[\cref{chap:TSRBLOD}:] The adaptive LOD algorithm for perturbed problems, which we summarized in \cref{sec:adaptive_PGLOD_from_hekema} can be found in the respective paper~\cite{HKM20}.
	
	The numerical experiments for the TSRBLOD in \cref{sec:TSRBLOD_experiments} have been cited in \cite{tsrblod}.
	\item[\cref{chap:TR_TSRBLOD}:] Most of the program code that was used for the experiments in \cref{sec:loc_experiments} is available within the respective preprint \cite{paper4}.
	The results that were not published within \cite{paper4} can be found in the supplementary code in~\cite{code:thesis}.
\end{description}


\end{document}